\documentclass[b5paper,twoside]{book}
\newcommand{\myTitle}{Large Deviations of Irreversible Processes\xspace}
\title{\myTitle}
\newcommand{\myName}{Mikola Christoph Schlottke\xspace}
\author{\myName}
\date{}

\usepackage{mathpazo} 
\usepackage[scaled]{helvet} 
\usepackage{courier} 
\usepackage[T1]{fontenc}
\usepackage[utf8]{inputenc}
\usepackage[font=small,labelfont={bf}]{caption}
\usepackage{epigraph}
\usepackage[explicit]{titlesec}
\usepackage{fancyhdr}
\renewcommand{\chaptermark}[1]{\markboth{#1}{}}
\renewcommand{\sectionmark}[1]{\markright{#1}}
\usepackage{afterpage}
\usepackage{xspace}	
\usepackage{makeidx}
\usepackage{amsmath}
\usepackage{amssymb}
\usepackage{mathtools}
\usepackage{amsthm}
\usepackage{enumerate}
\usepackage[shortlabels]{enumitem}
\usepackage{pinlabel}
\usepackage{subfig}
\usepackage{booktabs}
\usepackage{tikz-cd}
\usepackage{tikz}
\usetikzlibrary{matrix}
\pgfdeclarelayer{background}
\pgfsetlayers{background,main}
\usepackage{verbatim}
\usepackage[normalem]{ulem}
\usepackage{bbm}
\newtheorem{theorem}{Theorem}[section]
\newtheorem{lemma}[theorem]{Lemma}
\newtheorem{proposition}[theorem]{Proposition}

\theoremstyle{definition}
\newtheorem{definition}[theorem]{Definition}
\newtheorem{condition}[theorem]{Condition}	
\newtheorem{assumption}[theorem]{Assumption}
\newtheorem{example}[theorem]{Example}
\newtheorem{example*}{Example}
\newtheorem{remark}[theorem]{Remark}
\newtheorem*{running_example*}{Running Example}
\usepackage{xcolor}

\definecolor{thm_color_red}{RGB}{245, 188, 174}
\definecolor{thm_color}{RGB}{238, 187, 240}
\definecolor{dark_green}{RGB}{0, 123, 2}
\definecolor{very_light_blue}{RGB}{211, 255, 255}
\definecolor{dark_blue}{RGB}{0, 80, 255}
\definecolor{light_orange}{RGB}{255, 165, 0}
\definecolor{red_one}{RGB}{183, 0, 0}
\definecolor{gold}{RGB}{255, 255, 240}
\definecolor{gradient_blue_red}{RGB}{91, 40, 128}
\definecolor{turquoise}{RGB}{29, 169, 255}
\definecolor{turquoise_darker}{RGB}{0, 174, 255}
\definecolor{block_bgcolor1}{RGB}{29, 169, 255}
\definecolor{red_alert}{RGB}{222, 20, 20}

\newcommand{\myeqdef}{\stackrel{\mathclap{\text{def}}}{=}}
\newcommand{\dd}{ \mathrm{d}}
\newcommand{\PP}{\mathrm{P}}
\newcommand{\cA}{\mathcal{A}}

\newcommand{\cC}{\mathcal{C}}
\newcommand{\cD}{\mathcal{D}}
\newcommand{\cE}{\mathcal{E}}
\newcommand{\cF}{\mathcal{F}}
\newcommand{\cG}{\mathcal{G}}
\newcommand{\cH}{\mathcal{H}}
\newcommand{\cI}{\mathcal{I}}

\newcommand{\cK}{\mathcal{K}}
\newcommand{\cL}{\mathcal{L}}

\newcommand{\cP}{\mathcal{P}}

\newcommand{\cX}{\mathcal{X}}
\newcommand{\cY}{\mathcal{Y}}

\newcommand{\bE}{\mathbb{E}}

\newcommand{\bN}{\mathbb{N}}

\newcommand{\bP}{\mathbb{P}}

\newcommand{\bR}{\mathbb{R}}

\newcommand{\bT}{\mathbb{T}}

\newcommand{\bfH}{\mathbf{H}}

\newcommand{\bfR}{\mathbf{R}}

\newcommand{\bfV}{\mathbf{V}}

\def\calB{\mathcal B}
\def\calC{\mathcal C}
\def\calD{\mathcal D}

\def\calI{\mathcal I}
\def\calJ{\mathcal J}

\def\calL{\mathcal L}
\def\calM{\mathcal M}

\def\calO{\mathcal O}
\def\calP{\mathcal P}
\def\calQ{\mathcal Q}

\newcommand{\R}{\mathbb R}
\renewcommand{\P}{\mathbb P}
\newcommand{\arcsinh}{\operatorname{arcsinh}}
\newcommand{\unif}{\mathrm{Unif}}
\newcommand{\Leb}{\mathrm{Leb}}

\DeclareMathOperator*{\LIM}{LIM}

\newcommand{\vn}[1]{\left| \! \left| #1\right| \! \right|} 

\newcommand{\ip}[2]{\langle #1,#2\rangle}
\newcommand{\PR}{\mathbb{P}}
\newcommand{\bONE}{\mathbbm{1}}

\renewcommand{\u}{\cup}

\newcommand{\norm}[1]{\| #1 \| }

\usepackage{geometry}
\begin{document}
\frontmatter
\maketitle 
\thispagestyle{empty}
\noindent

\vfill

\begin{tabular}{p{109mm}}
%
%
%
%

\end{tabular}

\vfill

\noindent
\begin{tabular}{p{110mm}}
\\

Printed by ProefschriftMaken\\[6mm]

Cover: Formulas and pictures that evoke happy memories. Design by Mercedes Benjaminse, ProefschriftMaken and myself, with a blackboard background taken from www.freepik.com \\[6mm]

\noindent
A catalogue record is available from the Eindhoven University of Technology Library\\[2mm]
ISBN: 978-90-386-5072-2\\[6mm]
\\[6mm]

Copyright \copyright{} 2020 by \myName. All Rights Reserved. No part of this publication may be reproduced, stored in a retrieval system, or transmitted, in any form or by any means, electronic, mechanical, photocopying, recording or otherwise, without prior permission of the author.\\[6mm]
\end{tabular}

\cleardoublepage
\thispagestyle{empty}
\noindent

\begin{center}

\huge 


Large Deviations of \\ Irreversible Processes

\vspace{3em}

\Large
PROEFSCHRIFT

\vspace{3em}

ter verkrijging van de graad van doctor aan de Technische Universiteit Eindhoven, op gezag van de rector magnificus prof.dr.ir. F.P.T. Baaijens,
voor een commissie aangewezen door het College voor Promoties, in het openbaar te verdedigen op \\woensdag 8 juli 2020 om 16:00 uur

\vspace{3em}

door

\vspace{3em}

Mikola Christoph Schlottke

\vspace{3em}

geboren te Erlangen, Duitsland

\end{center}

\newpage
\thispagestyle{empty}
\normalsize

\noindent
Dit proefschrift is goedgekeurd door de promotoren en de samenstelling van de promotiecommissie is als volgt:

\vspace{1em}
\noindent
\begin{tabular}{lp{9cm}}
voorzitter:     & prof.dr. J.J. Lukkien \\[5pt]
1e promotor:    & prof.dr. M.A. Peletier \\[5pt]
2e promotor:    & prof.dr. F.H.J. Redig (Technische Universiteit Delft) \\[5pt]
leden:  		& prof.dr. J. Feng (University of Kansas) \\[5pt]
				& prof.dr. M.G. Westdickenberg (RWTH Aachen University) \\[5pt]
				& dr. O.T.C. Tse	 \\[5pt]
				& prof.dr. A.P. Zwart \\[5pt]
				& prof.dr. G. A. Pavliotis (Imperial College London)
\end{tabular}
\vspace*{\fill}

\noindent
Het onderzoek dat in dit proefschrift wordt beschreven is uitgevoerd in overeenstemming met de TU/e Gedragscode Wetenschapsbeoefening.
\cleardoublepage

\
\thispagestyle{empty}
\newpage
\
\thispagestyle{empty}

\thispagestyle{empty}
\clearpage
$\,$

%
%
%

\bigskip

$\,$
\bigskip

$\,$
\bigskip

$\,$
\bigskip

$\,$
\bigskip

$\,$
\bigskip

$\,$
\bigskip

$\,$
\bigskip

$\,$
\bigskip

$\,$
\bigskip

$\,$
\bigskip

$\,$
\bigskip

$\,$
\bigskip
$\,$
\bigskip

$\,$
\bigskip

\qquad\qquad\qquad\qquad\qquad\qquad\qquad\qquad
\qquad\qquad\qquad\qquad
\emph{To Mia}
\thispagestyle{empty}
\thispagestyle{empty}
\cleardoublepage
\chapter*{Abstract}
Time-irreversible stochastic processes are frequently used in natural sciences to explain non-equilibrium phenomena and to design efficient stochastic algorithms.
Our main goal in this thesis is to analyse their dynamics 
by means of large deviation theory. 
\smallskip

We focus on processes that become deterministic in a certain limit, and characterize their fluctuations around that deterministic limit by Lagrangian rate functions. Our main techniques for establishing these characterizations rely on the connection between large deviations and Hamilton-Jacobi equations. We sketch this connection with examples in the introductory parts of this thesis.
\smallskip

The second part of the thesis is devoted to irreversible processes that are motivated from molecular motors, Markov chain Monte Carlo (MCMC) methods and stochastic slow-fast systems. 
We characterize the asymptotic dynamics of molecular motors by Hamiltonians defined in terms of principal-eigenvalue problems. 
From our results about the zig-zag sampler used in MCMCs, we learn that maximal irreversibility corresponds to an optimal rate of convergence. 
In stochastic slow-fast systems, our main theoretical contributions are techniques to work with the variational formulas of Hamiltonians that one encounters in mean-field systems coupled to fast diffusions.
\smallskip

In the final part of the thesis, we study a family of Fokker-Planck equations whose solutions become singular in a certain limit.
The associated gradient-flow structures do not converge since the relative entropies diverge in the limit. To remedy this, we propose to work with a different variational formulation that takes fluxes into account, which 
is motivated by density-flux large deviations.
\smallskip

\paragraph{Keywords.} Large deviations, partial differential equations, viscosity solutions, comparison principle, variational techniques,~$\Gamma$-convergence, gradient flows.
\tableofcontents
\mainmatter
\chapter{Introduction}
\epigraph{\emph{There are three rules for writing the novel. Unfortunately, no one knows what they are.}}{W. Somerset Maugham.}
\section{Irreversible stochastic processes}
Many phenomena in natural sciences such as biology, chemistry and physics, are modelled by stochastic processes. In this thesis, we encounter for instance stochastic models of molecular motors~\cite{JulicherAjdariProst1997,KolomeiskyFisher07,Kolomeisky13}. Various other examples may be found in the monograph of Risken on Fokker-Planck equations~\cite[Chapters~1,~3~and~12]{Risken1996}.
The stochasticity is usually introduced in order to model the effect of noise in the dynamical systems. Our general objective in the works presented in this thesis is to analyse the dynamics of several examples of stochastic processes. 
\smallskip

Frequently, the dynamics simplifies in a certain limit where it becomes predictable. An example of such a simplification is the transition from microscopic to macroscopic scales. To illustrate this transition, imagine we would see the world only through a strong microscope. Then a familiar phenomenon such as a glas of water 
would all of a sudden appear complicated. Peering into the glas with our microscope, we observe the particles erratically moving back and forth, bouncing off and chasing each other in an unpredictable way. 
However, the moment we lay aside the microscope, this microscopic chaos disappears from our view; on the macroscopic scale, the density of particles does not evolve randomly, but becomes predictable.
When describing the particle density as a stochastic process, we should find that this stochastic process becomes deterministic in the limit of infinitely many particles.
\smallskip

There is a vast activity in probability theory and analysis to investigate mathematical theories of both microscopic and macroscopic dynamics. In particular, the focus lies on deriving a relationship between the dynamics at micro- and macroscales. Liggett~\cite{Liggett2004} as well as Kipnis and Landim~\cite{KipnisLandim1998} review and summarize works on interacting particle systems. Typically, the stochastic dynamics on the microscale incorporates basic features such as repulsion or attraction between particles~(for instance, the exclusion or the inclusion processes).
A common characteristic of the stochastic models is that in the limit of infinitely many particles, the particle density evolves deterministically according to a partial differential equation, such as the diffusion equation. 
\smallskip

Phenomena on the macroscopic scale such as first-order phase transitions originate from their underlying microscopic dynamics and may be explained using such micro-macro connections. More background and examples on this matter may be found in the books of Berglund and Gentz about noise-induced phenomena~\cite{BerglundGentz2005} and of Bovier and den Hollander on metastability~\cite{BovierDenHollander2016}. We remark that the randomness in the microscopic stochastic models is often rather
put in by hand than derived from first principles. 
This point of view builds up on two aspects coming together. First, the modelled system is chaotic in the sense of being highly sensitive towards the initial condition. Second, we have only partial knowledge about the initial condition. The system's behaviour appears to be random if both aspects, chaos and ignorance, come together. In this sense, the stochastic system may be seen as the approximation of a chaotic deterministic system. We refer to Bricmont~\cite{Bricmont1996} for more background on chaos.
\smallskip

In this thesis, we analyse irreversible stochastic processes by means of large deviation theory. Our central goals are to derive their limiting dynamics, and to characterize their fluctuations around this limiting dynamics by means of Lagrangian rate functions. As we shall further discuss at the end of Section~\ref{intro:sec:large-deviatoin-theory}, we are motivated by the fact that while reversible processes lead via large deviation theory to gradient flows, it is an open question of which variational formulations can, in principle, be derived for irreversible processes.
\smallskip

In Section~\ref{intro:sec:large-deviatoin-theory}, we introduce our main tool for the analysis of irreversible processes, large deviation theory. Then we give examples that clarify the concepts of pathwise large deviation principles and Lagrangian rate functions. In Section~\ref{intro:sec:overview-of-the-thesis}, we give a more detailed overview of the thesis. In Chapter~2, we provide an introduction to our main method for proving large deviation principles, the Feng-Kurtz method~\cite{FengKurtz2006}. In Chapter~3, we consider stochastic models of walking molecular motors. In Chapter~4, we analyse a Markov chain Monte Carlo method based on the irreversible zig-zag sampler. In Chapters~5 and~6, our interest lies in deriving---by means of large deviation principles---limiting evolution equations of mean-field interacting particles that are coupled to fast external processes. In Chapter~7, we consider a limit problem of variational structures of certain PDEs. Finally, we discuss our results in Chapter~\ref{chapter:discussion}.
\section{Large deviation theory}
\label{intro:sec:large-deviatoin-theory}
The first unified treatment of large deviation theory in the sense of an abstract framework is attributed to Srinivasa Varadhan, who laid the ground for decades of active mathematical research by his landmark paper~\cite{Varadhan1966}. Varadhan was honored in 2007 with the Abel Prize "for his fundamental contributions to probability theory and in particular for creating a unified theory of large deviations".
Numerous works have enriched the scope of large deviation theory by connecting it to other mathematical fields and applications in natural sciences. The most commonly used techniques for studying large deviations are summarized in a number of different books and papers; we only give an incomplete list here.
Varadhan relates among other things function space integrals with large deviations in his lectures~\cite{Varadhan1984}.
Freidlin and Wentzell were the first to explore pathwise large deviations of stochastic processes~\cite{FreidlinWentzell1998}. Ellis shows the relation of large deviations and statistical mechanics~\cite{Ellis1985}. 
Deuschel and Strook introduced the term \emph{exponential tightness}~\cite{DeuschelStroock1989}.
Numerous abstract techniques that are frequently used in large deviation theory are presented by Dembo and Zeitouni~\cite{DemboZeitouni1998}. A concise overview of large deviations with many examples may be found in the lectures of den Hollander~\cite{denHollander2000}. Bovier and den Hollander also give a brief overview in their book on metastability~\cite[Chapter~6]{BovierDenHollander2016}. Recent monographs focusing on stochastic processes are the semigroup approach of Feng and Kurtz~\cite{FengKurtz2006}, and the weak-convergence approach initiated by Dupuis and Ellis~\cite{DupuisEllis1997}, which Budhiraja and Dupuis extend in~\cite{BudhirajaDupuis2019}.
\smallskip

In this section, we first exemplify the general definition of a large deviation principle. The first example is a simple observation of \emph{exponential decay} of probabilities. With the second example, we illustrate a \emph{concentration effect} that occurs \emph{exponentially}, and furthermore motivate the notion of \emph{large deviations}. The examples provide a useful mental image for interpreting the general definition. For further reading and examples suitable for familiarization, we refer to Richard Ellis' beautiful note on Boltzmann's discoveries~\cite[Section~3]{Ellis1999}, where he illustrates how relative entropies arise naturally from Stirling's formula. Further illustrating examples may also be found in Ellis' lectures on large deviations~\cite{Ellis1995}, and in particular in Hugo Touchette's review~\cite[Section~2]{Touchette2009}.
We also refer to Terence Tao's note~\cite{Tao2015_275A}
for a short introduction to the mathematical notions from probability theory we use below.
\smallskip

After having introduced the concept of a large deviation principle, we will specialize further to the setting of this thesis: pathwise large deviations for stochastic processes.  We illustrate by means of a classical example some interesting aspects of a pathwise large deviation principle, with a focus on the so-called \emph{action-integral representation} of the rate function. 
\newpage

\begin{example*}
	If we toss a fair coin~$n$ times,
	the probability of observing 
	"only heads"
	is
	\begin{equation*}
	\mathbb{P}\left(\text{only heads}\right) = \left(\frac{1}{2}\right)^n = e^{-n \log 2}.
	\end{equation*}
	Let us point out the following observations:
	\begin{enumerate}[$\bullet$]
		\item If $n$ is large, the event "only heads" is \emph{unlikely} or \emph{improbable}.
		\item As we let $n$ grow, the event "only heads" becomes \emph{increasingly} unlikely.
		\item The probability of observing "only heads" is \emph{exponentially small} with respect to~$n$. The event "only heads" decays exponentially with rate~$\log 2$. 
	\end{enumerate}
\end{example*}
\begin{example*}\label{ex:Gaussian-rv}
	Let~$X_1,X_2,\dots$ be a sequence of~i.i.d.
	real-valued random variables. Suppose each~$X_i$ is normally distributed with mean~$\mu\in\mathbb{R}$ and variance one,
	\begin{equation*}
	\mathbb{P}\left(X_i \in A\right) = \int_A \rho(x)\,\dd x,\quad \rho(x)= \frac{1}{\sqrt{2\pi}}e^{-(x-\mu)^2/2}.
	\end{equation*}
	Let us focus on the behaviour of their partial sums $S_n:=\sum_{i=1}^nX_i$ for large~$n$. 
	The probability distribution~$\rho_n$ of the averages~$\frac{1}{n}S_n$ is depicted in~Figure~\ref{fig:gaussians_concentrating}.
	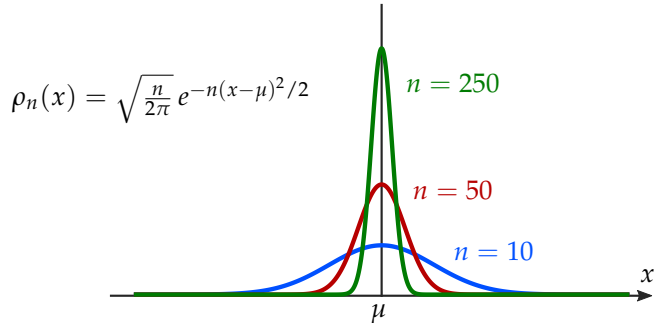
\begin{figure}[!htbp]
	\centering
	{\labellist
	\pinlabel $x$ at 1050 50
	\pinlabel $\rho_n(x)=\sqrt{\frac{n}{2\pi}}\,e^{-n(x-\mu)^2/2}$ at 100 400
	\pinlabel $\mu$ at 525 -25
	\pinlabel {\color{dark_green}{$n=250$}} at 670 430
	\pinlabel {\color{red_one}{$n=50$}} at 670 220
	\pinlabel {\color{dark_blue}{$n=10$}} at 750 100
	\endlabellist
	\centering
	\includegraphics[scale=.2]{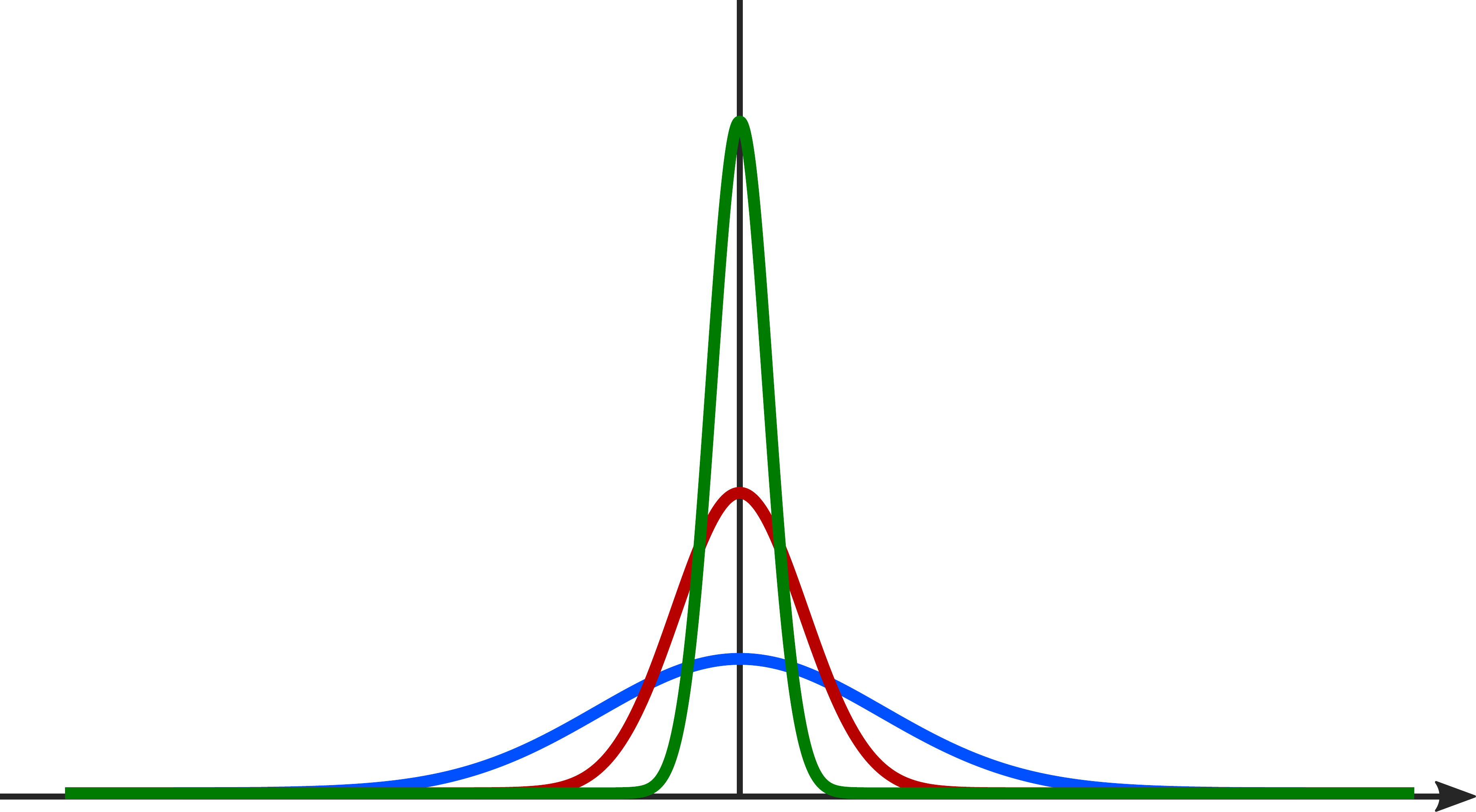}
	}
	\caption{
	The probability distribution~$\rho_n(x)$ of the averages~$\frac{1}{n}S_n$. As~$n$ increases, 
	the probability distribution concentrates around the mean~$\mu$.
	}
	\label{fig:gaussians_concentrating}
	\end{figure}
	We observe a concentration effect of the distribution around the mean as~$n$ increases. That means for large~$n$, we are likely to observe~$\frac{1}{n}S_n\approx \mu$.
	To summarize this concentration effect, let~$\varepsilon>0$, and write~$B_\varepsilon(\mu):=(\mu-\varepsilon,\mu+\varepsilon)$ for the small interval around~$\mu$. In accordance with the weak law of large numbers, we have
	\begin{equation}
	\label{intro:ex2:LLN}
	\mathbb{P}\left(\frac{1}{n}S_n\in B_\varepsilon(\mu)\right) \xrightarrow{n\to\infty} 1.
	\end{equation} 
	A natural question is: how fast does the distribution of the averages concentrate around the mean as~$n$ tends to infinity? 
	Let us show that there is a rate with which it concentrates \emph{exponentially}.
	We abbreviate the quadratic function in the exponent of~$\rho_n$ by~$\mathcal{I}(x):=(x-\mu)^2/2$. For~$\varepsilon>0$, we want to estimate~$\mathbb{P}\left((1/n)S_n \notin B_\varepsilon(\mu)\right)$.
	Using the formula of the probability density~$\rho_n$ and exploiting its symmetry, we find
	\begin{align}\label{intro:eq:Ex2:log-of-prob}
	\frac{1}{n}\log \mathbb{P}\left(\frac{1}{n}S_n\notin B_\varepsilon(\mu)\right) = \frac{1}{n}\log\left(2\, \sqrt{\frac{n}{2\pi}}\right) + \frac{1}{n}\log \int_{\mu+\varepsilon}^{\infty} e^{-n \,\mathcal{I}(x)}\,\dd x.
	\end{align}
	As~$n\to\infty$, the first term vanishes. In the second term, the lowest value of the exponent 
	dominates the integral. The precise statement is the Laplace principle; for a set~$A$ and a function~$g$ bounded from below,
	\begin{equation*}
	\frac{1}{n}\log\int_A e^{-n\,g(x)}\,\dd x \xrightarrow{n\to\infty}\sup_{x\in A} \left(-g(x)\right) = -\inf_{x\in A}g(x).
	\end{equation*}
	With these remarks, we find
	\begin{align*}
	\lim_{n\to\infty}\frac{1}{n}\log \mathbb{P}\left(\frac{1}{n}S_n\notin B_\varepsilon(\mu)\right) = - \inf_{x\in(\mu+\varepsilon,\infty)}\mathcal{I}(x) = -\frac{\varepsilon^2}{2}.
	\end{align*} 
	Hence for any~$\delta>0$ (smaller than~$\varepsilon^2/2$), we find for all~$n$ large enough that
	\begin{equation*}
	\mathbb{P}\left(\frac{1}{n}S_n\notin B_\varepsilon(\mu)\right) \leq \exp\left\{-n \left[\frac{1}{2}\varepsilon^2 - \delta\right]\right\}.
	\end{equation*}
	In this sense the concentration effect~\eqref{intro:ex2:LLN} occurs exponentially.
	With similar reasoning, we find for any closed~$A\subseteq \mathbb{R}$ not containing~$\mu$,
	\begin{equation}\label{intro:eq:ex2-log-P-converges}
	\frac{1}{n}\log\mathbb{P}\left(\frac{1}{n}S_n\in A\right) \xrightarrow{n\to\infty} - \inf_{x\in A}\mathcal{I}(x).
	\end{equation}
	Therefore, with~$\mathrm{r}(A):=\inf_A\mathcal{I}(\cdot)$, we find for~$\delta>0$ that for~$n$ sufficiently large,
	\begin{equation*}
	\mathbb{P}\left(\frac{1}{n}S_n\in A\right)\leq e^{-n \left(\mathrm{r}(A)-\delta\right)}.
	\end{equation*}
	Let us summarize: in regions~$A$ away from the mean, the probability mass is exponentially small with respect to~$n$, and the exponential decay rate~$r(A)$ is the minimum of the quadratic function~$\mathcal{I}(\cdot)$ evaluated over~$A$.
	\smallskip
	
	We close this example by pointing out in what sense the above considerations are related to large deviations. The random variable~$\sqrt{n}\left((1/n)S_n-\mu\right)$ is normally distributed around zero with variance one. This means that observations of the type~$S_n \approx n\mu + \sqrt{n}\,x$ are 
	normally distributed for large~$n$. This is a deviation from what we expect by~$\sqrt{n}\,x$, and in that sense, fluctuations of order~$\sqrt{n}$ are "normal" (the generalization of this statement is the central limit theorem). For any~$\varepsilon$, the event~$\frac{1}{n}S_n\notin B_\varepsilon(\mu)$ corresponds to observing events of the type~$S_n\approx n\mu + n\,\varepsilon$. This is a deviation of order~$n$, which is no longer captured by the central limit theorem. Therefore, these type of events are called large deviations.
	The generalization of these observations beyond this example ($X_i$ that are not normally distributed) is known as Cramér's theorem~\cite{Cramer1938,CramerTouchette2018}, and we refer to~\cite[Chapter~I and Theorem~I.4]{denHollander2000} for more details.\qed
\end{example*}
\subsubsection{Precise formulation of a large deviation principle.}
We typically consider sequences of probability measures~$\PP_n$ on a state space~$\mathcal{X}$, concentrating at a single element~$x\in\mathcal{X}$. Above in Example~2, the measures~$\PP_n$ correspond to the distribution of the averages~$\frac{1}{n}S_n$ with state space~$\mathcal{X}=\mathbb{R}$, that means~$\PP_n=\mathbb{P}((1/n)S_n\in \cdot)$. The single element is the mean value~$x=\mu$, and the concentration effect can be formulated as a weak law of large numbers; denoting by~$B_\varepsilon(x)$ the ball of radius~$\varepsilon>0$ around~$x$,
\begin{equation*}
	\PP_n\left(B_\varepsilon(x)\right)\xrightarrow{n\to\infty} 1,
\end{equation*}
and for any Borel set~$A\subseteq\cX $ whose closure does not contain~$x$,
\begin{equation*}
	\PP_n(A) \xrightarrow{n\to\infty} 0.
\end{equation*}
Frequently, we can observe an exponential decay of these probabilities; at least intuitively, we find a rate~$\mathrm{r}(A)$ depending in the set~$A$ with which for large~$n$,
\begin{equation*}
	\PP_n(A) \approx e^{-n \cdot \mathrm{r}(A)}.
\end{equation*}
One attempt of making this rigorous would be to say: a sequence of probability measures~$\PP_n$ satisfies a large deviation principle if there is a rate~$\mathrm{r}:\mathcal{B}(\mathcal{X})\to[0,\infty)$ with which for any Borel subset~$A\subseteq\mathcal{X}$,
\begin{equation}\label{eq:intro:LDP:log-conv-to-r-A}
	\frac{1}{n}\log\PP_n(A) \xrightarrow{n\to\infty} -\mathrm{r}(A).
\end{equation}
Furthermore, the example from above suggests that this rate can be characterized by a so-called rate function~$\mathcal{I}:\mathcal{X}\to[0,\infty]$ as
\begin{equation}\label{eq:intro:LDP:rate-inf-RF}
	\mathrm{r}(A) = \inf_{x\in A} \mathcal{I}(x).
\end{equation}
Varadhan's definition is a suitable more general form of~\eqref{eq:intro:LDP:log-conv-to-r-A}. We first give his definition here. A complete separable metric space~$\mathcal{X}$ is called a \emph{Polish space}. We call a map~$\mathcal{I}:\mathcal{X}\to[0,\infty]$ a \emph{rate function} if the sublevel sets~$\{x\in\mathcal{X}\,:\,\mathcal{I}(x)\leq C\}$ are compact for all~$C\geq 0$. In the literature, such rate functions are called good rate functions---since all rate functions we encounter in this thesis are good, we adopt the convention of~\cite{BudhirajaDupuis2019} and omit the adjective "good". For a Borel subset~$A\subseteq\mathcal{X}$, we let~$\mathrm{int}(A)$ be its interior and~$\mathrm{clos}(A)$ be its closure.
\begin{definition}[Large Deviation Principle]\label{def:LDP}
	For~$n = 1,2,\dots,$ let~$\PP_n$ be a probability measure on a Polish space~$\mathcal{X}$. We say the family of measures $\{\PP_n\}_{n\in \mathbb{N}}$ satisfies a \emph{large deviation principle} with rate function $\mathcal{I}:\mathcal{X} \to [0,\infty]$ if for any Borel subset~$A\subseteq \mathcal{X}$,
	\begin{align*}
		-\inf_{x\in \mathrm{int}(A)}\mathcal{I}(x) 
		&\leq \liminf_{n\to\infty}\frac{1}{n}\log\PP_n(A)\\
		&\leq\limsup_{n\to\infty}\frac{1}{n}\log\PP_n(A) 
		\leq -\inf_{x\in\mathrm{clos}(A)}\mathcal{I}(x).
		\tag*\qed
	\end{align*}
\end{definition}
Let~$X^n$ be a random variable with law~$\PP_n\in\mathcal{P}(\mathcal{X})$. We say that the sequence of random variables~$\{X^n\}_{n\in\mathbb{N}}$ satisfies a large deviation principle if the sequence of their laws~$\{\PP_n\}_{n\in\mathbb{N}}$ does. In this case, we write~$\PP_n=\mathbb{P}\left(X^n\in\cdot\right)$ for the law. Furthermore, we abbreviate the statement that~$\{X^n\}_{n\in\mathbb{N}}$ satisfies a large deviation prinicple with rate function~$\mathcal{I}:\mathcal{X}\to[0,\infty]$ as follows:
\begin{equation}\label{intro:LDP-written-in-local-form}
\mathbb{P}\left(X^n\approx x\right) \sim e^{-n\,\mathcal{I}(x)},\quad n\to\infty.
\end{equation}
We alert the reader that the tilde has no mathematical precise meaning. Equation~\eqref{intro:LDP-written-in-local-form} should rather be read as a total statement summarizing all the essential information; the probability that~\emph{$X^n$ is close to~$x$} (~$X^n\approx x$) decays \emph{exponentially} as~$n$ tends to infinity (~$\sim e^{-n\,\mathcal{I}(x)}$).
The notation is motivated by the fact that
\begin{equation*}
\lim_{\varepsilon\to 0} \limsup_{n\to\infty}\frac{1}{n}\log \mathbb{P}\left(X^n\in B_\varepsilon(x)\right) = -\mathcal{I}(x).
\end{equation*}
Let us mention how a large deviation principle really corresponds to the exponential decay of probabilities. If the rate function is continuous, then we recover~\eqref{eq:intro:LDP:log-conv-to-r-A} for Borel subsets~$A\subseteq\mathcal{X}$ satisfying~$\mathrm{clos}(\mathrm{int}(A))=\mathrm{clos}(A)$~\cite[Corollary~1]{Ellis1999}, and in particular for any~$\varepsilon>0$, if~$n$ is large enough,
\begin{equation*}
	e^{-n(\mathrm{r}(A)+\varepsilon)} \leq \PP_n(A) \leq e^{-n(\mathrm{r}(A)-\varepsilon)}.
\end{equation*}
In general, if a non-trivial rate function has a unique minimizer~$x$, then for a Borel set~$A$ whose closure does not contain~$x$, we have~$\mathrm{r}(\mathrm{clos}(A))>0$, and the limsup bound implies exponential decay of~$\PP_n(A)$~\cite[Corollary~2]{Ellis1999}.
\smallskip

We may motivate Varadhan's definition of a large deviation principle in terms of the liminf- and limsup bounds by analogy to weak convergence of probability measures. To that end, consider~$\PP_n$ and~$\PP\in\mathcal{P}(\mathcal{X})$. The measures~$\PP_n$ are said to converge weakly to~$\PP$ if for any Borel set~$A\subseteq\mathcal{X}$,
\begin{equation}\label{eq:intro:LDP:liminf-limsup-weak-convergence}
	\PP(\mathrm{int}(A))\leq \liminf_{n\to\infty}\PP_n(A) \leq \limsup_{n\to\infty}\PP_n(A) \leq \PP(\mathrm{clos}(A)).
\end{equation} 
Demanding "pointwise" convergence~$\PP_n(A)\to \PP(A)$ for all~$A$ would exclude examples such as~$\PP_n = \delta_{1/n}$ and~$\PP=\delta_0$ with~$\mathcal{X}=\mathbb{R}$. Also in Example~2 from above with~$\PP_n$ the law of~$(1/n)S_n$ and~$\PP=\delta_\mu$, the singleton set~$A=\{\mu\}$ violates this strong convergence condition. The notion of weak convergence applies to many interesting examples while still providing useful information. By the Portmanteau Theorem~\cite[Theorem~2.1]{Billingsley1999}, weak convergence is equivalent to the convergence of expectations;~$\int f\,\dd\PP_n\to\int f\,\dd\PP$ for any function~$f\in C_b(\mathcal{X})$. An equivalent formulation is to demand the liminf-bound for all open sets and the limsup-bound for all closed sets~\cite[Theorem~2.1]{Billingsley1999}.
%
\smallskip

Next, we motivate the fact that the exponential decay rates~$\mathrm{r}(A)$ are characterized by minimizing a rate function~$\mathcal{I}(\cdot)$ over the region~$A$. For two real-valued positive sequences $a_n,b_n$, suppose $a_n>b_n$ for all $n$ sufficiently large. Then
\begin{equation*}
\left|\frac{1}{n}\log\left(a_n+b_n\right) - \frac{1}{n}\log (a_n)\right| = \frac{1}{n}\log\left(1 + b_n/a_n\right) \leq \frac{1}{n}\log 2\xrightarrow{n\to\infty} 0.
\end{equation*}
Hence the maximal value~$a_n=\max(a_n,b_n)$ dominates the sum on the logarithmic scale. This fact is known as the so-called \emph{the-winner-takes-it-all} principle. Now suppose a set~$A$ satisfies~\eqref{eq:intro:LDP:log-conv-to-r-A} with some rate~$\mathrm{r}(A)>0$, and suppose we can decompose~$A=A_1\cup A_2$ into disjoint sets~$A_1,A_2$ satisfying~\eqref{eq:intro:LDP:log-conv-to-r-A} as well. Then using additivity,~$\PP_n(A_1\cup A_2)=\PP_n(A_1)+\PP_n(A_2)$, we find by the winner-takes-it-all principle
\begin{equation*}
	\frac{1}{n}\log\PP_n(A_1\cup A_2) \xrightarrow{n\to\infty} -\min\left\{\mathrm{r}(A_1),\mathrm{r}(A_2)\right\}.
\end{equation*}
Therefore, we may expect the exponential rates to be given by~\eqref{eq:intro:LDP:rate-inf-RF}. Similar to the definition of weak convergence~\eqref{eq:intro:LDP:liminf-limsup-weak-convergence}, passing to the interior and closure in Definition~\ref{def:LDP} is necessary in order for the limits to hold for any Borel set~$A$. 
\smallskip

A large deviation principle is a type of concentration inequality, and therefore implies a strong type of convergence of random variables. The minimizers of the rate function are the elements corresponding to the strong law of large numbers, as demonstrated by the following theorem. For a rate function~$\mathcal{I}:\mathcal{X}\to[0,\infty]$, we denote by $\{\mathcal{I} = 0 \}$ the set of its global minimizers.
\begin{theorem}
	\label{thm:math-formulation-LDP:LDP-implies-as}
	For $n=1,2\dots$, let $X^n$ be a random variable taking values in a Polish space~$(\mathcal{X},d)$. Suppose that~$\{X^n\}_{n \in \mathbb{N}}$ satisfies a large deviation principle with rate function $\mathcal{I}$. Then $d(X^n,\{\mathcal{I} = 0\}) \to 0$ almost surely as $n \to \infty$.
\end{theorem}
This theorem can be proven via the limsup-bound of the large deviation principle, and applying the Borel-Cantelli Lemma. In many examples, we can verify uniqueness of the minimizer~$x_0$ of a rate function. Then by Theorem~\ref{thm:math-formulation-LDP:LDP-implies-as}, a large deviation principle implies~$X^n\to x_0$ almost surely. We point out that the rate function in Theorem~\ref{thm:math-formulation-LDP:LDP-implies-as} is assumed to have compact sub-level sets.
\smallskip

Next to the law of large numbers, the central limit theorem can as well be understood from a large deviation principle.
Specialising to $E=\mathbb{R}$, a formal Taylor expansion around a minimizer $x_0$ of the rate function yields for $x\approx x_0$,
\begin{equation*}
\mathbb{P}\left(X^n\approx x\right) \sim e^{-n\left[\mathcal{I}(x_0) + (x-x_0)\mathcal{I}'(x_0) + \frac{1}{2}(x-x_0)^2\mathcal{I}''(x_0)\right]} = e^{-n \mathcal{I}''(x_0)(x-x_0)^2/2}.
\end{equation*}
In that sense, fluctuations around minimizers of the rate function are normally distributed. The curvature of the rate function is inverse proportional to the variance: if the rate function is rapidly growing near the minimizer, then the variance is small, and vice versa. Bryc makes this connection precise in~\cite{Bryc1993}.
\subsubsection{Pathwise large deviations in stochastic systems.}
In this thesis, we will mostly focus our attention on stochastic processes~$X^n$ that become deterministic in the limit of a parameter~$n$ tending to infinity. In particular, we are interested in situations in which this transition to a deterministic limit occurs exponentially in the sense of a large deviation principle. In this context, we speak of \emph{pathwise} large deviations, because we make statements about the paths of~$X^n$.
Here, we illustrate with a classical example what makes a pathwise large deviation principle interesting. In the example, we will point out the following two central features. First, the \emph{typical behaviour}: the expected trajectory of~$X^n$, corresponding to the law of large numbers limit, 
is recovered from the minimizer of the rate function. Second, the \emph{least-action principle}: if the stochastic process realizes an event far away from this expected trajectory, the most likely way in which this event occurs can be determined by minimizing the rate function.
%
%
\begin{example*}
	Let $E=\mathbb{R}$,~$x_0\in E$. For $n\in\mathbb{N}$, consider the process~$X^n$ solving
	\[
	\dd X^n_t = \frac{1}{\sqrt{n}}\,\dd B_t,\quad X^n(0) =x_0,
	\]
	where~$B_t$ is the standard Brownian motion. For large~$n$, the process~$X^n$ corresponds to a small-diffusion regime. 
%
	The transition probabilities $P_n(t,x,\dd y)$ of~$X^n$ are normal distributions,
	\[
	\mathbb{P}\left(X^n(t)\in \dd y\,|\,X^n(0)=x\right) \myeqdef P_n(t,x,\dd y) = \sqrt{\frac{n}{2\pi t}}\,e^{-n(y-x)^2/2t}\,\dd y.
	\]
	We fix a time interval~$[0,T]$. Let~$\mathcal{X}=C_\mathbb{R}[0,T]$ be the set of continuous maps~$x:[0,T]\to\mathbb{R}$, equipped with the uniform norm. We consider the~$X^n$ as random variables in~$\mathcal{X}$, and are interested in the behaviour of~$X^n$ in the limit~$n\to\infty$.
	\smallskip
	
	For large values of $n$, typical realizations of~$X^n$ are shown in Figure~\ref{fig:small-diffusion-10-realizations}. 
	\begin{figure}[!htbp]
		\centering
		{\labellist
			\pinlabel $t$ at 1150 50
			\pinlabel $x_0$ at -40 250
			\pinlabel $0$ at -40 50
			\pinlabel $X^n_t(\omega)$ at 100 500
			\endlabellist
			\centering			\includegraphics[scale=.2]{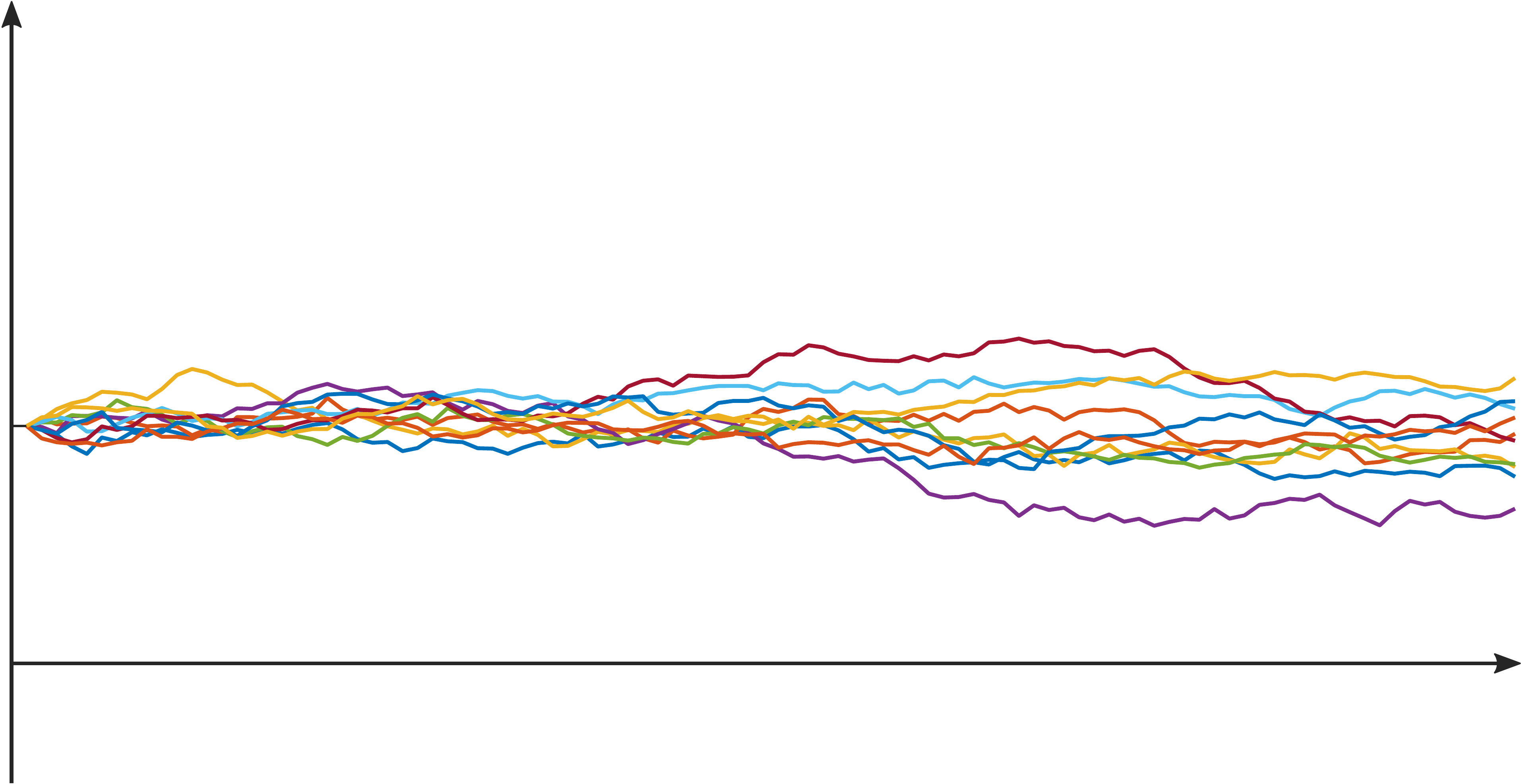}
		}
		\caption{
			Ten realizations~$X^n(\omega)$ of the stochastic process~$X^n$ with deterministic initial condition~$X^n_0=x_0$.
		}
		\label{fig:small-diffusion-10-realizations}
	\end{figure}
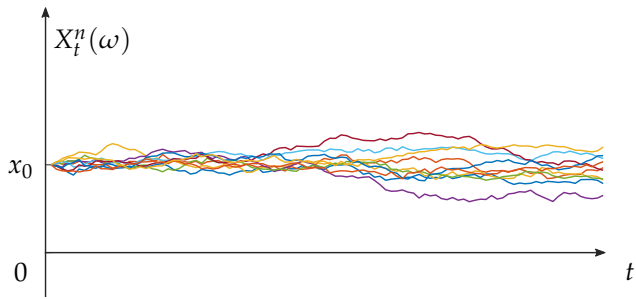 
	Judging by eye, most realizations are close to the constant path $\overline{x}\equiv x_0$ determined by the initial starting point~$x_0$. Indeed, for any~$t\in[0,T$], the one-dimensional time marginals~$X^n(t)$ are converging to~$x_0$, as can be seen from the transition probabilities.
	In fact, the probability of observing realizations of~$X^n$ that deviate from the constant path~$\overline{x}(t) := x_0$ vanishes exponentially fast as~$n\to\infty$: the path measures~$\PP_n := \mathbb{P}\left(X^n\in \cdot\right)\in\mathcal{P}(\mathcal{X})$ satisfy a large deviation principle with rate function ~$\mathcal{I}:\mathcal{X}\to[0,\infty]$ given by
	\begin{equation}\label{intro:eq:RF-Schilders-Theorem}
	\mathcal{I}(x) = \int_0^T \frac{1}{2}|\partial_t x(t)|^2\,\dd t.
	\end{equation}
	This fact is known as Schilder's theorem (e.g.~\cite[Theorem~5.2.3]{DemboZeitouni1998}), which is a special case of the Freidlin-Wentzell theorem (e.g.~\cite[Theorem~5.6.3]{DemboZeitouni1998}). If a trajectory~$x$ is not absolutely continuous or~$x(0)\neq x_0$, then~$\mathcal{I}(x)=\infty$. 
	As we discussed below the definition of a large deviation principle, an informal but useful interpretation is to say that for a path~$x$ satisfying $x(0)=x_0$, we have
	\begin{equation}\label{intro:eq:LDP-Schilders-Theorem}
	\mathbb{P}\left(X^n \approx x\right) \sim e^{-n \,\mathcal{I}(x)}, \quad \text{as}\;n\to\infty.
	\end{equation}
	Alternatively, 
	let~$B_\varepsilon(x)$ be the ball of radius~$\varepsilon$ around~$x$ with respect to the uniform norm in~$\mathcal{X}$. Then with~$\mathcal{I}(B_\varepsilon(x)):=\inf_{B_\varepsilon(x)}\mathcal{I}$,
	\begin{equation*}
		\mathbb{P}\left(X^n\in B_\varepsilon(x)\right)\sim \exp\{-n\, \mathcal{I}(B_\varepsilon(x))\},\quad \text{as}\;n\to\infty.
	\end{equation*}
	In terms of the topology on~$\mathcal{X}$, this means the probability of~$X^n$ being inside an~$\varepsilon$-tube around~$x$ decays exponentially with respect to~$n$.
\smallskip

	Let us point out two interesting conclusions from the large deviation principle~\eqref{intro:eq:LDP-Schilders-Theorem}.
	First, suppose that~$\mathcal{I}(x)>0$. Then the probability of~$X^n$ being close to~$x$ with respect to the uniform norm decays exponentially with increasing~$n$. Since~$\mathcal{I}(x)>0$ whenever~$x$ has a non-zero velocity, we conclude that realizations of~$X^n$ are with high probability close to the minimizer~$\overline{x}$ of the rate function~\eqref{intro:eq:RF-Schilders-Theorem}. The minimizer is unique and given by the constant path~$\overline{x} \equiv x_0$. This identifies~$\overline{x}$ as the law of large number limit of~$X^n$, by Theorem~\ref{thm:math-formulation-LDP:LDP-implies-as} from above.
	\smallskip
	
	Second, we illustrate the least-action principle. Consider a closed subset of trajectories~$A\subseteq\mathcal{X}$ not containing~$\overline{x}$. This set represents an atypical event. Suppose~$x$ is the unique trajectory minimizing the rate function evaluated over~$A$,
	\begin{equation*}
		\mathcal{I}(x) = \min_{y\in A} \mathcal{I}(y).
	\end{equation*}
	Then if the event~$A$ occurs, it will most likely be realized as~$X^n\approx x$. More precisely, for any~$\varepsilon>0$, we have by~\cite[Theorem~1.4]{BudhirajaDupuis2019} that
	\begin{equation*}
		\mathbb{P}\left(X^n\in B_\varepsilon(x)\,|\,X^n\in A\right) \xrightarrow{n\to\infty} 1.
	\end{equation*}
	For instance, fix~$C>0$ and consider~$A=\{x\in\mathcal{X}\,:\,x(0)=x_0,\,x(T)\geq x_0+C\}$. This event corresponds to~$X^n$ exceeding the threshold~$x_0+C$ at final time~$T$. To determine the most likely way in which this rare event occurs, we have to solve the corresponding minimization problem with rate function~\eqref{intro:eq:RF-Schilders-Theorem},
	\begin{equation*}
		\mathcal{I}(x) = \min_{y\in A} \int_0^T \frac{1}{2}|\partial_t y(t)|^2\,\dd t.
	\end{equation*}
	Solving the Euler-Lagrange equation with the boundary conditions~$y(0)=x_0$ and~$y(T) = \lambda$ (any~$\lambda\geq C$), we find that the minimizing trajectory~$x$ is the path with constant velocity~$C/T$, that is~$x(t) = x_0 + t\,C/T$.
	\qed
\end{example*}
The example illustrates in what sense a large deviation principle for stochastic processes contains more information than the law of large numbers. We have an exponential estimate on the probabilities of deviating from the law of large number limit, and the rate function contains information about the rare-event behaviour. In the example, the Brownian motion exceeds the threshold $x_0+C$ most likely by following the path with a constant slope. Determining the rare-event behaviour for more involved examples is an interesting topic, but we will not study it in this thesis. For more background on the least-action principle, we refer to the following papers and the references therein. Weinan, Ren and Vanden-Eijnden use Freidlin-Wentzell theory to study rare events in a couple of perturbed dynamical systems, including for instance the one-dimensional Ginzburg-Landau model~\cite{WeinanRenVanden-Eijnden2004}, and introduced the string method~\cite{WeinanRenVanden-Eijnden2002}. Metzner, Sch{\"u}tte and Vanden-Eijnden provide an overview of illustrating examples~\cite{MetznerSchutteVanden-Eijnden2006}, and Grafke and Vanden-Eijnden explore numerical methods for various rare-event algorithms~\cite{GrafkeVanden-Eijnden2019}.
\smallskip

The form of the rate function~\eqref{intro:eq:RF-Schilders-Theorem} is a special case of a more general principle. For many examples, we can derive rate functions of the form
\begin{equation}\label{intro:eq:RF-action-integral}
	\mathcal{I}(x) = \int_0^T \mathcal{L}(x(t),\partial_t x(t))\,\dd t.
\end{equation}
The map~$\mathcal{L}:\mathbb{R}\times\mathbb{R}\to[0,\infty]$ appearing in the rate function is called the \emph{Lagrangian}. In the above example,~$\mathcal{L}(x,v) = v^2/2$ is independent of~$x$. We call~\eqref{intro:eq:RF-action-integral} an \emph{action-integral representation} or~\emph{Lagrangian rate function}, which is motivated from the least-action principle that we discussed above. If a process~$X^n$ satisfies a large deviation principle with a Lagrangian rate function, then its limiting dynamics~$\overline{x}=\overline{x}(t)$ can be determined by solving~$\mathcal{L}(\overline{x}(t),\partial_t\overline{x}(t))=0$. 
\smallskip

For more involved stochastic processes, it is often difficult to derive an exact characterization of their limiting dynamics in the first place. In these situations, deriving the Lagrangian provides one way of finding a good characterization. This is what we do in the first part of this thesis, where we are interested in two main aspects: deriving Lagrangians and extracting useful information from them. A common feature making the stochastic processes that we study interesting is their irreversibility with respect to time. We close this chapter by pointing out our motivation for considering irreversible processes.
\subsubsection{The role of irreversibility---an open question.}
Jordan, Kinderlehrer and Otto demonstrated that the solution of the diffusion equation is the steepest descent of the relative entropy~\cite[Theorem~5.1]{JordanKinderlehrerOtto1998}. Their variational formulation is motivated by the backward Euler approximation scheme, and represents an example of a \emph{gradient flow}---we introduce these concepts in more detail in Section~\ref{GF_NGF:sec:GF} of Chapter~\ref{chapter:GF-to-NGF}. A special role in the gradient flow is played by the Wasserstein distance between probability measures, which serves as the metric in the gradient flow. Therefore, this variational formulation is called the Wasserstein gradient flow. Such a variational structure involving the Wasserstein distance can be recognized in many other PDEs, e.g.~\cite{AmbrosioGigliSavare2008,BlanchetCalvezCarrillo2008,CarrilloDiFrancescoFigalliLaurentSlepcev2011,CarrilloYoungPilTse2019,CarlenGangbo2004,Gigli2010,GianazzaSavareToscani2009,MatthesMcCannSavare2009,Savare2007,Lisini2009}. Many of these PDEs arise from stochastic particle systems, but it is a priori not clear how to find a corresponding gradient flow. Therefore, it is helpful to know how to derive the corresponding Wasserstein gradient flows from the microscopic dynamics. A recent example of such a derivation is the study of Gavish, Nyquist and Peletier~\cite{GavishNyquistPeletier2019} about hard-rod systems.
\smallskip

Adams, Dirr, Peletier and Zimmer derived the Wasserstein gradient flow for the diffusion equation by means of large deviation theory~\cite[Theorem~3]{AdamsDirrPeletierZimmer2011}, by considering the empirical density of independent Brownian motions and sending the number of particles to infinity. Soon after, Mielke, Peletier and Renger revealed 
that the
gradient flow is a consequence of microscopic \emph{reversibility} of the Brownian motions~\cite{MielkePeletierRenger2014}. The argument exploits an action-integral form of the rate function. Let us briefly state in what sense.
For~$n$ independent Brownian motions~$\{B^i\}_{i=1,\dots,n}$, the empirical particle density defined by~$\rho^n = \frac{1}{n}\sum_{i=1}^n \delta_{B^i}$
is a measure-valued process that converges in the narrow topology to the solution of the diffusion equation as~$n\to\infty$. That means~$\rho^n\rightharpoonup\rho$, where~$\partial_t\rho = \Delta \rho$. The sequence also satisfies a large deviation principle in~$\mathcal{X}=C_{\mathcal{P}(\mathbb{R})}[0,\infty)$ with a rate function given by
\begin{equation}\label{intro:eq:action-RF-measure-valued-process}
\mathcal{I}(\mu) = \int_0^\infty\mathcal{L}(\mu(t),\partial_t\mu(t))\,\dd t.
\end{equation}
The rate function satisfies~$\mathcal{I}(\rho)=0$. We ignore here the initial conditions and do not go into details, but refer to~\cite[Theorem~13.3]{FengKurtz2006} for the precise statement. Mielke, Peletier and Renger show in~\cite[Section~4.2]{MielkePeletierRenger2014} how to decompose the Lagrangian in~\eqref{intro:eq:action-RF-measure-valued-process} in such a way that one can recognize the Wasserstein gradient-flow structure in the rate function~\eqref{intro:eq:action-RF-measure-valued-process}.
\smallskip

The argument that connects the rate function to a gradient flow is based on reversibility. This argument extends from the abovementioned example to a wider class of reversible Markov processes, which by large deviations give rise to so-called \emph{generalized gradient flows}~\cite{MielkePeletierRenger2014} (see also Section~\ref{GF_NGF:sec:GF} of Chapter~\ref{chapter:GF-to-NGF}).
Triggered by this connection between variational structures of PDEs and large deviations, a natural question we can ask is: \emph{which variational structures can we derive from irreversible processes?} While this question is still open, we remark that the starting point for this connection in the reversible case is a Lagrangian rate function. The main questions we ask in the first part of this thesis are thus: \emph{how can we prove large deviation principles for irreversible dynamics and obtain action-integral representations of the rate functions?} \emph{How can we establish useful characterizations of the Lagrangians?} We hope that the techniques we develop by answering such questions can contribute to extending the abovementioned connection to a suitable class of irreversible processes. The study of irreversible processes is also of independent interest, since various non-equilibrium phenomena are modelled by irreversible processes; we refer to the note of Harris and Touchette~\cite[Section~1.2]{TouchetteHarris2011} for more. A broader overview on irreversibility may be found in Bricmont's note~\cite[Section~3]{Bricmont1996}.
\smallskip

Since we will come back to reversibility, let us close this section by formulating this property here.
%
For a state space~$E$ and a trajectory~$\gamma\in\mathcal{X}=C_E[0,T]$, let~$\mathrm{rev}(\gamma)$ be the trajectory defined by
\begin{equation*}
\mathrm{rev}(\gamma)(t):= \gamma(T-t),\quad t\in [0,T].
\end{equation*}
That means~$\mathrm{rev}(\gamma)$ is the time-reversed trajectory of~$\gamma$. For a Borel set of trajectories~$A\subseteq \mathcal{X}$, let~$\mathrm{rev}(A) := \{\mathrm{rev}(\gamma)\,:\,\gamma\in A\}$. For~$\mu\in\mathcal{P}(E)$, we write~$\mathbb{P}_\mu$ for the path distribution of a process~$X$ with initial distribution~$X(0)\sim \mu$. A measure~$\pi\in\mathcal{P}(E)$ is \emph{stationary} if~$\mathbb{E}_\pi f(X(t))$ is constant in time for any observable~$f=f(x)$.
\begin{definition}[Reversibility]\label{intro:def:reversibility}
	Let~$X$ be a Markov process with path distributions~$\mathbb{P}_\mu$ and stationary measure~$\pi\in\mathcal{P}(E)$. We say~$X$ is \emph{reversible} with respect to~$\pi$ if for any Borel subset~$A\subseteq\mathcal{X}$,
	\begin{align}
	\mathbb{P}_\pi\left(X\in A\right) = \mathbb{P}_\pi\left(X\in \mathrm{rev}(A)\right).
	\tag*\qed
	\end{align}
\end{definition}
For illustration, an example of a reversible process is a jump process on~$\{1,2,3\}$ with uniform nearest-neighbor jump rates; its stationary measure is the uniform measure. A counterexample is a jump process on~$\{1,2,3\}$ with jumps only clockwise,~$r(1\to 2)=r(2\to 3)=r(3\to 1)>0$, and all other jump rates equal to zero. The stationary measure is also the uniform measure, but for the set~$A_\circlearrowright$ containing all trajectories only going clockwise,~$\mathrm{rev}(A_\circlearrowright)=A_\circlearrowleft$, and therefore
\begin{equation*}
\mathbb{P}_\pi\left(X\in A_\circlearrowright\right) = 1\quad \text{and}\quad \mathbb{P}_\pi\left(X\in \mathrm{rev}(A_\circlearrowright)\right) = 0.
\end{equation*}
The notion of reversibility of Definition~\ref{intro:def:reversibility} is sometimes also refered to as \emph{microscopic reversibility} or \emph{time reversibility}. For Markov processes, there are several useful equivalent characterizations of reversibility that we will work with. For instance, reversibility is equivalent to symmetry of the infinitesimal generator or the semigroup, in the sense made precise in~\cite[Proposition~5.3]{Liggett2004}.
\section{Overview of the thesis}
\label{intro:sec:overview-of-the-thesis}
Here we outline the content of the subsequent chapters. We further detail the relation of our results to the literature in the introductory parts of the chapters.
\paragraph{Chapter~\ref{chapter:LDP-via-HJ}: Large Deviations via Hamilton-Jacobi Equations.} In this chapter we demonstrate how to prove pathwise large deviation principles by exploiting the connection to Hamilton-Jacobi equations~\cite{FengKurtz2006}. The gist of this connection is that solving certain PDEs of Hamilton-Jacobi type allows us to prove an action-integral representation of the rate function involving the so-called Lagrangian. 
The crucial insight we take from this chapter is an algorithm that allows us to rigorously derive the Lagrangian starting from microscopic dynamics. 
\smallskip

While the results in this chapter are not novel, some proofs simplify because we choose to illustrate all concepts in a simpler setting. The extension to theorems including the general settings are presented in the monograph of Jin Feng and Thomas Kurtz~\cite{FengKurtz2006}. We close the chapter by outlining the relation of our presentation to such general settings. 
\paragraph{Chapter~\ref{chapter:LDP-for-switching-processes}: Large Deviations of Switching Processes.}
This chapter is based on a joint work with Mark Peletier~\cite{PeletierSchlottke2019}. Our work is inspired by a series of papers by Mirrahimi, Perthame and Souganidis about PDEs describing molecular motors~\cite{PerthameSouganidis09a, PerthameSouganidis2009Asymmetric, Mirrahimi2013}. We consider a general class of switching Markov processes that comprise the PDE models as a special case, and prove pathwise large deviation principles. The large-deviation theorems extend and generalize the results of~\cite{PerthameSouganidis09a, PerthameSouganidis2009Asymmetric, Mirrahimi2013}. The main tool we work with is the connection of large deviations to Hamilton-Jacobi equations. In particular, this connection allows us to study within the same framework multiple limit regimes as well as continuous and discrete models of molecular motors.
\smallskip

As an application, we show how macroscopic transport properties of molecular motors can be deduced from associated principal-eigenvalue problems. We work with variational formulas of principal eigenvalues to demonstrate that breaking detailed balance is necessary for obtaining transport. In Section~\ref{subsec:model-of-molecular-motor} we discuss an example of a continuous molecular-motor model that illustrates our more general results.
\paragraph{Chapter~\ref{chapter:LDP-of-empirical-measures}: Large Deviations of Empirical Measures.}
This chapter is based on a joint work with Joris Bierkens and Pierre Nyquist~\cite{BierkensNyquistSchlottke2019}.
Joris Bierkens and Gareth Roberts discovered the zig-zag process as a scaling limit of the Lifted Metropolis-Hastings~\cite{BierkensRoberts2017}. The zig-zag process is an example of a piecewise deterministic Markov process in position and velocity space.
The process can be designed to have an arbitrary Gibbs-type marginal probability density for its position coordinate, which makes it suitable for Monte Carlo simulation of continuous probability distributions. An important question in assessing the efficiency of this method is how fast the empirical measure converges to the stationary distribution of the process. We provide a partial answer to this question by characterizing the large deviations of the empirical measure from the stationary distribution. Based on the Feng-Kurtz approach to large deviations~\cite{FengKurtz2006}, we develop an abstract framework aimed at encompassing piecewise deterministic Markov processes in position-velocity space. We derive explicit conditions for the zig-zag process to allow the Donsker-Varadhan variational formulation of the rate function, both for a compact setting (the torus) and one-dimensional Euclidean space. 
\smallskip

For reversible processes, Donsker and Varadhan offer an exact formula of the rate function involving the stationary measure. There is no generic formula for irreversible processes, which makes it generally harder to draw conclusions from the rate function.
For the zig-zag process however, we derive an explicit expression for the Donsker-Varadhan functional for the case of a compact state space. We use this form of the rate function to address a key question concerning the optimal choice of the switching rate of the zig-zag process. We show that maximal irreversibility corresponds to the fastest possible convergence to the stationary distribution.
\paragraph{Chapter~\ref{chapter:LDP-in-slow-fast-systems}: Large Deviations in Stochastic Slow-Fast Systems.} 
This chapter is based on a work in progress with Richard Kraaij.
We give conditions for proving pathwise large deviations in stochastic slow-fast systems in the limit of time-scale separation tending to infinity. The conditions are imposed in order to solve the corresponding Hamilton-Jacobi equations. In the limit regime we consider, the convergence of the slow variable to its deterministic limit and the convergence of the fast variable to equilibrium are competing at the same scale. We cast the rate functions in action-integral form and interpret the Lagrangians in two ways: in terms of a double-optimization problem of the slow variable's velocity and the fast variable's distribution, and in terms of a principal-eigenvalue problem associated to the slow-fast system.
\smallskip

As an application, we provide a large-deviation theorem for the empirical density-flux pair of mean-field interacting particles coupled to fast diffusion. This system cannot be treated with classical methods. We further show how the Lagrangian can be used to derive an averaging principle from the large deviation principle.
\paragraph{Chapter~\ref{chapter:CP-for-two-scale-H}: Comparison Principle for Two-Scale Hamiltonians.} This chapter is based on a joint work with Richard Kraaij~\cite{KraaijSchlottke2019}. We study the well-posedness of Hamilton-Jacobi-Bellman equations on subsets of~$\mathbb{R}^d$. The Hamiltonian consists of two parts: an internal Hamiltonian depending on an external control variable and a cost function penalizing the control. We show under suitable assumptions that if a comparison principle holds for the Hamilton-Jacobi equation involving only the internal Hamiltonian, then the comparison principle holds for the Hamilton-Jacobi-Bellman equation involving the full Hamiltonian. In addition to establishing uniqueness, we give sufficient conditions for existence of solutions. 
Our key features are that the internal Hamiltonian is allowed to be non-Lipschitz and non-coercive in the momentum variable, and that we allow for discontinuous cost functions. To compensate for the greater generality of our approach, we assume sufficient regularity of the cost function on its sub-level sets and that the internal Hamiltonian satisfies a comparison principle uniformly in the control variable on compact sets.
As an application, we show our established result to cover interesting examples that were posed as open problems in the literature as well as mean-field Hamiltonians that cannot be treated with standard methods.
\paragraph{Chapter~\ref{chapter:GF-to-NGF}: Gradient Flow to Non-Gradient-Flow.} This chapter is based on a work in progress with Mario Maurelli and Mark Peletier.
We study a singular limit problem arising in modelling chemical reactions. At finite~$\varepsilon>0$, the model is a Fokker-Planck equation corresponding to a particle diffusing in a double-well potential. In the limit~$\varepsilon=0$, the solution concentrates at the two potential wells. Arnrich, Mielke, Peletier, Savaré and Veneroni~\cite{ArnrichMielkePeletierSavareVeneroni2012} considered a \emph{symmetric} double-well potential and proved Gamma convergence of the associated Wasserstein gradient-flow structures. We take the double-well potential to be \emph{asymmetric}. In that case, the Wasserstein gradient flows do no longer converge. This is because the relative entropies diverge in the limit. To obtain a meaningful limit of a variational structure associated to the family of equations, we  consider density-flux functionals rather than density functionals. The Wasserstein gradient flow is obtained from the density-flux functional by contraction.
\paragraph{Chapter~\ref{chapter:discussion}: Discussion and Future Questions.} In this final chapter, we first summarize the results presented in this thesis. Then we discuss their limitations and point out questions that we could not answer so far.
\chapter{Introduction to Large Deviations via Hamilton-Jacobi Equations}
\label{chapter:LDP-via-HJ}
\chaptermark{LDP via Hamilton-Jacobi Equations}
\section{A general strategy of proof}
\label{BG:sec:math-gem}
This chapter is an introduction to a connection between two mathematical subjects: pathwise large deviations of stochastic processes on the one hand, and Hamilton-Jacobi equations on the other hand. Jin Feng and Thomas Kurtz show in their monograph~\cite{FengKurtz2006} how to rigorously connect these subjects by means of mathematical theorems.
The scope of the approach is demonstrated by the examples given in~\cite[Section~I.1.4]{FengKurtz2006}. 
\smallskip

When I first tried to work with the theory, I had difficulties to get started. This was mainly because 
the general conditions are involved, which can make it difficult for a newcomer to grasp the essence. I write this chapter with the intention to facilitate for other newcomers the process of getting started. To do so, I sacrifice generality for clarity, and answer three straightforward questions I initially struggled to answer for myself, and to which I could not find straight answers in the literature. Before we get to the questions, let us first have a look at the gist of the connection.
\subsubsection{The connection in a nutshell.} For $E:=\mathbb{R}^d$ and a finite~$T>0$, let~$\mathcal{X}:= C_E[0,T]$ be the set of~$E$-valued continuous trajectories~$x:[0,T]\to E$, equipped with the supremum norm. 
Consider a sequence of Markov processes~$\{X^n\}_{n=1,2\dots}$, where each~$X^n$ is regarded as a random variable in~$\mathcal{X}$, with deterministic initial conditions~$X^n(0)=x_0$. 
\smallskip

We will typically consider~$X^n$ that become deterministic in the limit~$n\to\infty$: frequently we expect by the law of large numbers that there exists a trajectory~$\overline{x}\in\mathcal{X}$ such that~$X^n\to \overline{x}$ almost surely as~$n\to\infty$. Then for any closed set of trajectories~$A\subseteq \mathcal{X}$ not containing~$\overline{x}$, we have~$\mathbb{P}\left(X^n\in A\right)\to 0$ as~$n\to\infty$.
We say~$X^n$ satisfies a \emph{pathwise large deviation principle} if these probabilities are exponentially small with respect to~$n$ in the sense of Definition~\ref{def:LDP}. Our goal is both to prove a large deviation principle and to find a useful formula of the rate function.
\smallskip

Let us state the connection to Hamilton-Jacobi equations. We denote the transition probabilities of~$X^n$ by~$P_n(t,x,\dd y)$. Define for a bounded measurable function~$f\in B(E)$ and~$t\geq 0$ the function~$V_n(t)f$ by
\begin{equation}\label{eq:BG:def-nonlinear-semigroup-Vn}
	V_n(t) f(x) := \frac{1}{n}\log\int_{E}e^{n f(y)}\,P_n(t,x,\dd y).
\end{equation}
For each~$n=1,2,\dots$, the family~$\{V_n(t)\}_{t\geq 0}$ forms a one-parameter semigroup of maps acting on~$B(E)$. 
Further below, we prove that 
the convergence of these semigroups~$V_n$ to a limiting semigroup~$\{V(t)\}_{t\geq 0}$ implies a pathwise large deviation principle of~$X^n$. This limiting semigroup can be regarded as the \emph{semigroup flow} of a Hamilton-Jacobi equation; 
there is a 
map~$\mathcal{H}:E\times\mathbb{R}^d\to\mathbb{R}$ called the \emph{Hamiltonian} with which the function~$u(t,x) := V(t)f(x)$
is the solution to
\begin{equation*}
	\begin{cases}
		\partial_t u(t,x) = \mathcal{H}(x,\nabla_x u(t,x)),\\
		u(0,x) = f(x).
	\end{cases}
\end{equation*}
In which precise sense~$u$ solves this equation is not important here. The Hamiltonian fully characterizes the large-deviation fluctuations via its Legendre dual defined as the map~$\mathcal{L}(x,v) := \sup_{p}[p\cdot v-\mathcal{H}(x,p)]$, which we call the \emph{Lagrangian}. Indeed, frequently the large-deviation rate function~$\mathcal{I}:\mathcal{X}\to[0,\infty]$ satisfies
\begin{equation}\label{BG:eq:intro-RF-action-integral}
	\mathcal{I}(x) = \int_0^T\mathcal{L}(x(t),\partial_t x(t))\,\dd t.
\end{equation}
This is a useful formula which allows us to determine the law of large number limit for complicated processes~$X^n$, namely as the path~$\overline{x}$ satisfying~$\mathcal{I}(\overline{x})=0$, which solves the equation~$\mathcal{L}(\overline{x}(t),\partial_t\overline{x}(t))=0$. We call~\eqref{BG:eq:intro-RF-action-integral} an \emph{action-integral representation}.
\smallskip

The Hamiltonian can be derived by taking the limit of the so-called \emph{nonlinear generators}~$H_n$ of the semigroups~$V_n(t)$, which are formally determined by~$H_n := \frac{\dd}{\dd t}|_{t=0}V_n(t)$.
These nonlinear generators converge in a suitable sense to a limiting operator~$H$ acting on functions as~$Hf(x)=\mathcal{H}(x,\nabla f(x))$, where~$\mathcal{H}$ is the Hamiltonian from above. This derivation will provide us with a recipe for three aspects at once: finding the Hamiltonian, giving a rigorous proof of large deviations, and proving the action-integral formula~\eqref{BG:eq:intro-RF-action-integral}. The goal of this chapter is to prove a rigorous version of this recipe in a simplified setting.\qed
\subsubsection{Three questions that we answer in this chapter.} 
Our first question is:
\begin{enumerate}[label=(\arabic*)]
	\item \emph{Why is verifying the convergence of the nonlinear semigroups~$V_n(t)$ to a limiting semigroup~$V(t)$ sufficient for proving pathwise large deviation principles?}
\end{enumerate}
We answer this question in Section~\ref{BG:sec:LDP-from-convergence-of-semigroups} by means of Theorem~\ref{thm:LDP-via-semigroup-convergence:compact}---the additional assumption of exponential tightness appearing therein is not important for now. 
\smallskip

In practice, verifying convergence of the nonlinear semigroups~$V_n(t)$ from scratch is hard. In that sense, the result formulated in Theorem~\ref{thm:LDP-via-semigroup-convergence:compact} really only serves as a stepping stone to obtain useful and applicable results. 
The bulk of the general functional analytical work in~\cite{FengKurtz2006} lies in detecting useful conditions to verify the convergence of nonlinear semigroups from the convergence of their \emph{generators}. Let us briefly sketch the idea.
For a bird's-eye view on semigroups, we refer to Chapters~I and~VII of Engel's and Nagel's monograph~\cite{EngelNagel1999}. 
\begin{example}\label{ex:semigroups-are-exponentials}
	Let $T:[0,\infty)\to\mathbb{C}$ be a continuous map forming a semigroup, that means~$T(t+s)=T(t)T(s)$ and~$T(0)=1$. Cauchy and Abel proved the existence of a unique scalar $g\in \mathbb{C}$ with which the semigroup is given by $T(t)=e^{tg}$ (\cite[Theorem~1.4]{EngelNagel1999}). We call $g$ the \emph{generator} of the semigroup~$T(t)$. Let us point out two aspects about this result:
	\begin{enumerate}[$\bullet$]
		\item The whole semigroup~$T$ is uniquely identified by its generator~$g$.
		\item While the map $T$ is only assumed to be continuous, its semigroup property $T(t+s)=T(t)T(s)$ actually enforces differentiability. Its generator is uniquely determined by $g = \frac{d}{dt}T(0)$.
	\end{enumerate}
	Based on this result, we can prove the following recipe for convergence of a sequence of semigroups $\{T_n\}_{n\in\mathbb{N}}$. First, identify their generators by computing $g_n = \frac{d}{dt}T_n(0)$. Second, identify the limit $g:=\lim_n g_n$. Then this limit generates a semigroup by $T(t) := e^{tg}$, and the semigroups~$T_n$ converge to~$T$ uniformly over compact time intervals.
	\qed
\end{example}
In the spirit of this example, the natural question we can ask is whether there exists a similar recipe for proving convergence of the semigroups~$V_n(t)$.
That means first identifying generators~$H_n$ by making sense of $H_n = \frac{d}{dt}V_n(0)$, and then secondly identifying a suitable limit $H:=\lim_n H_n$. 
In the above example the semigroups are complex scalars, and the fact that the limit~$g$ is a complex scalar is sufficient to generate a semigroup by means of the formula~$T(t):=e^{tg}$. 
Since the semigroups~$V_n(t)$ are nonlinear maps defined on~$B(E)$, the conditions on a limit~$H$ are more involved. 
\smallskip

Therefore, our second question is:
\begin{enumerate}[label=(\arabic*)]
	\setcounter{enumi}{1}
	\item \emph{How does the recipe from Example~\ref{ex:semigroups-are-exponentials} for verifying convergence of semigroups carry over to the nonlinear semigroups~$V_n(t)$?}
\end{enumerate}
The answer we give in Section~\ref{BG:sec:LDP-from-convergence-of-semigroups} identifies the generators~$H_n$ as certain nonlinear operators and establishes conditions on a limit operator~$H$ to generate a nonlinear semigroup~$V(t)$. We find that the convergence of generators~$H_n\to H$ indeed implies the desired convergence of semigroups~$V_n\to V$.
The conditions on the limit operator~$H$ are imposed in order to make sense of the formula~$V(t)=e^{tH}$. In Theorem~\ref{thm:LDP-classical-sol} in Section~\ref{subsec:LDP-via-classical-sol}, we first see how this program leads to the problem of solving PDEs of the form
\begin{equation*}
	(1-\tau H) f = h,
\end{equation*}
where~$\tau>0$ and the function~$h=h(x)$ are given. 
In a running example, by which we illustrate intermediate results, this PDE is
\[
f(x) - \tau \frac{1}{2}|\nabla f(x)|^2 = h(x),\quad x\in\mathbb{R}.
\]
We show in Section~\ref{subsec:LDP-via-visc-sol} why the notion of viscosity solutions provides the right tools to solve these type of PDEs.
The recipe we obtain for proving pathwise large deviation principles is summarized in Theorem~\ref{thm:LDP-visc-sol}, and we apply this theorem to the running example. The drawback of this theorem is that the rate function is still intricate.
Therefore, our third question is:
\begin{enumerate}[label=(\arabic*)]
	\setcounter{enumi}{2}
	\item \emph{How can we prove an action-integral representation~\eqref{BG:eq:intro-RF-action-integral} of the rate function?}
\end{enumerate}
We provide an answer in Section~\ref{BG:sec:semigroup-flow-HJ-eq}. The required Hamiltonian~$\mathcal{H}$ is identified from the limit operator~$H$, by recognizing the latter to act on functions by
\begin{equation*}
	Hf(x) = \mathcal{H}(x,\nabla f(x)).
\end{equation*}
Let us summarize where we stand after we will have answered the above three questions. We find an algorithm that provides us with a convenient method for guessing the form of a rate function.
For a sequence of~$\mathbb{R}^d$-valued Markov processes~$X^n$, the algorithm can be divided into five steps. First, start from the generators~$L_n$ of~$X^n$. Second, compute the nonlinear generators defined by acting on functions as~$H_n f := \frac{1}{n}e^{-nf}L_n e^{nf}$. Third, identify the limit operator~$Hf=\lim_n H_nf$. Fourth, identify the Hamiltonian~$\mathcal{H}:\mathbb{R}^d\times\mathbb{R}^d\to\mathbb{R}$ as the map satisfying~$Hf(x)=\mathcal{H}(x,\nabla f(x))$ for all~$f$ in the domain of~$H$. Finally, define the Lagrangian as the Legendre-Fenchel dual~$\mathcal{L}(x,v):=\sup_{p\in\mathbb{R}^d}\left[p\cdot v - \mathcal{H}(x,p)\right]$. Now details aside, the rate function satisfies~\eqref{intro:eq:RF-action-integral} with this Lagrangian.

In Section~\ref{BG:sec:setting}, we briefly summarize some aspects about Markov processes that we will work with in the subsequent sections. Each subsequent section is devoted to answer one of the three questions posed above. Section~\ref{BG:sec:LDP-from-convergence-of-semigroups} answers the first question about semigroup convergence, Section~\ref{sec:LDP-via-generators} the second question about deriving semigroup convergence from generator convergence, and Section~\ref{BG:sec:semigroup-flow-HJ-eq} the third question about the
action-integral representation.
\section{Setting: Markov process in compact state space}
\label{BG:sec:setting}
We denote by~$E$ a Polish space, that is a complete separable metric space. We will assume~$E$ to be compact. For~$T>0$, let~$\mathcal{X}:=C_E[0,T]$ be the set of continuous maps~$\gamma:[0,T]\to E$, equipped with the supremum norm. 
We consider a set of transition probabilities $\{P(t,x,\dd y)\}_{t\geq 0}$ such that:
\begin{enumerate}[label=(\roman*)]
	\item For any~$x\in E$,~$P(t,x,\cdot)$ is a probability measure on~$E$, and~$P(0,x,\cdot)=\delta_x$.
	\item For any Borel subset $A\subseteq E$, the map $x\mapsto P(t,x,A)$ is measurable on~$E$, and for any $s\leq t$, we have~$P(s+t,x,A) = \int_E P(s,y,A) \, P(t,x,\dd y)$.
\end{enumerate}
By Theorem~1.1 in Chapter~IV of~\cite{EthierKurtz1986}, such a collection of transition probabilities gives rise to a corresponding Markov process~$X(t)|_{t\geq 0}$; for~$t\in[0,T]$, we have a random variable~$X(t)$ in~$E$, and~$X$ is a random variable in~$\mathcal{X}$. The Markov process is identified with the path distributions~$\{\mathbb{P}_x\}_{x\in E}$, where each~$\mathbb{P}_x$ is a probability measure on~$\mathcal{X}$ describing the law of the process when starting at~$x$.
\smallskip
 
If we think of the process as describing a particle that moves in~$E$, then the value~$P(t,x,A)$ corresponds to the probability that starting from~$x$, the particle propagates in time~$t$ into the region~$A$. 
It is the conditional probability
\[
P(t,x,A) = \mathbb{P}\left(X(t)\in A\,|\,X(0)=x\right) = \mathbb{P}_x\left(X(t)\in A\right).
\]

Let~$B(E)$ be the set of bounded and measurable functions on~$E$. We call the family of maps~$\{S(t)\}_{t\geq 0}$, with $S(t):B(E)\to B(E)$ given by
\begin{equation}\label{BG:eq:def:Markov-semigroup}
S(t)f(x) := \int_E f(y)\,P(t,x,\dd y),
\end{equation}
the \emph{semigroup} associated to the Markov process $X$. Its semigroup property, that is~$S(t+s)=S(t)S(s)$, is inherited from the transition probabilities.
%
\section{Large deviations via convergence of semigroups}
\label{BG:sec:LDP-from-convergence-of-semigroups}
The main point of this section is to answer our first key question: why is the convergence of nonlinear semigroups sufficient for proving pathwise large deviation principles? We answer it by proving Theorem~\ref{thm:LDP-via-semigroup-convergence:compact} below, which is a simplification of~\cite[Theorem~5.15, Corollary~5.17]{FengKurtz2006}. 
\begin{definition}[Nonlinear Semigroup associated to Markov process]\label{def:nonlinear-semigroups}
	Let $E$ be a Polish space.
	For a Markov process~$X^n$ with transition probabilities $P_n(t,x,\dd y)$, define the map $V_n(t):B(E)\to B(E)$ 
	by~\eqref{eq:BG:def-nonlinear-semigroup-Vn}; that is for~$f\in B(E)$,
	\begin{equation}\label{eq:def:nonlinear-semigroup}
	V_n(t)f(x) := \frac{1}{n}\log\int_E e^{nf(y)}\,P_n(t,x,\dd y).
	\end{equation}
	We call the family~$\{V_n(t)\}_{t\geq 0}$ the \emph{nonlinear semigroup} associated to the process~$X^n$.\qed
\end{definition}
We will see further below how Bryc's formula leads us directly to consider these nonlinear semigroups. The family $\{V_n(t)\}_{t \geq 0}$ inherits its semigroup property from the semigroup $\{S_n(t)\}_{t\geq 0}$ of the Markov process~$X^n$, since
\begin{align*}
	V_n(t+s)f(y) &\myeqdef \frac{1}{n}\log S_n(t+s) e^{nf}(y)\\
	&= \frac{1}{n}\log S_n(t)\left[S_n(s)e^{nf}\right](y)\\
	&= \frac{1}{n}\log S_n(t)\left[e^{n V_n(s)f}\right](y) 
	\myeqdef V_n(t) \left[V_n(s) f\right](y).
\end{align*}
For the theorem, we also need the following condition. 
\begin{definition}[Exponential tightness]
	Let $\{\PP^n\}_{n\in\mathbb{N}}$ be a sequence of probability measures on a Polish space~$\mathcal{X}$. The sequence $\{\PP^n\}_{n\in\mathbb{N}}$ is \emph{exponentially tight} if for any $\ell > 0$, there exists a compact set $K_\ell \subseteq E$ such that
	\begin{align*}
		\limsup_{n\to\infty} \frac{1}{n} \log \PP^n\left(\mathcal{X}\setminus K_\ell\right) \leq -\ell.
		\tag*\qed
	\end{align*}
\end{definition}
	Exponential tightness means the mass of the probability measures $\PP^n$ concentrates exponentially fast on compact sets: given an arbitrary rate $\ell > 0$, for any $\varepsilon > 0$ there exists a compact set $K_{\ell,\varepsilon}\subseteq \mathcal{X}$ such that for $n$ large enough,
	\[
	\PP^n\left(\mathcal{X}\setminus K_{\ell,\varepsilon}\right)\leq e^{-n(\ell-\varepsilon)}.
	\]
	We comment further below on the role of exponential tightness. Let us first formulate the theorem. For a function~$g\in B(E)$, we write~$\|g\|_E:=\sup_E |g|$.
	\begin{theorem}[Large deviations via convergence of nonlinear semigroups]
		\label{thm:LDP-via-semigroup-convergence:compact}
		For~$n=1,2,\dots$, let~$X^n$ be a Markov process in~$\mathcal{X}=C_E[0,T]$ with path distribution denoted by~$\PP^n := \mathbb{P}(X^n\in \cdot)\in\mathcal{P}(\mathcal{X})$, and with the corresponding nonlinear semigroup~$V_n$ from Definition~\ref{def:nonlinear-semigroups}. Assume the following:
		\begin{enumerate}[label=(\roman*)]
			\item The sequence $\{\PP^n\}_{n\in\mathbb{N}}$ is exponentially tight in $\mathcal{X}$.
			\item\label{BG:item:LDP-via-semigroup:conv-Vn} There are maps $V(t):C(E)\to C(E)$, $t\geq 0$, such that for any sequence of functions $f_n\in B(E)$ and $f\in C(E)$, 
			\begin{equation*}
			\text{if}\quad \|f-f_n\|_E\xrightarrow{n\to\infty}0,\qquad\text{then}\quad
			\|V(t)f-V_n(t)f_n\|_E\xrightarrow{n\to\infty}0.
			\end{equation*}
		\end{enumerate} 
		Suppose furthermore that the initial conditions~$\{X^n(0)\}_{n\in\mathbb{N}}$
		satisfy a large deviation principle in~$E$ with rate function $\mathcal{I}_0:E\to[0,\infty]$. Then the sequence $\{X^n\}_{n\in \mathbb{N}}$ satisfies a large deviation principle in $\mathcal{X}$ with rate function $\mathcal{I}:\mathcal{X}\to[0,\infty]$ given by~\eqref{eq:RF-with-semigroup-approach} below.
	\end{theorem}
For any~$n\in\mathbb{N}$, the sequence~$X^n$ with initial conditions~$X^n(0)\sim \nu_n$ has a path distribution~$\PP^n = \mathbb{P}_{\nu_n}(X^n\in \cdot)\in\mathcal{P}(\mathcal{X})$. The theorem gives two conditions under which a large deviation principle for the initial condition bootstraps to these path distributions.
The rate function in Theorem~\ref{thm:LDP-via-semigroup-convergence:compact} is determined by the limit~$V(t)$,
\begin{equation}\label{eq:RF-with-semigroup-approach}
\mathcal{I}(x) = \mathcal{I}_0(x(0)) + \sup_{k\in\mathbb{N}}\sup_{(t_1,\dots,t_k)} \sum_{i=1}^k \mathcal{I}_{t_i-t_{i-1}}(x(t_i)|x(t_{i-1})),
\end{equation}
where $\mathcal{I}_t(z|y)$ given by
\begin{equation}\label{BG:eq:RF-1d-marginals}
\mathcal{I}_t(z|y) = \sup_{f\in C(E)}\left[f(z)-V(t)f(y)\right].
\end{equation}
We will encounter the functions~$\mathcal{I}_t(\cdot|y)$ as the rate functions for the one-dimensional marginals. Before we give an overview of the proof of Theorem~\ref{thm:LDP-via-semigroup-convergence:compact}, a remark on exponential tightness. This property is always expected: if a sequence of probability measures on a Polish space satisfies the large-deviation upper bound, then the sequence is exponentially tight~\cite[Exercise~4.1.10]{DemboZeitouni1998}. 
In general, the role of exponential tightness is to bootstrap the large-deviation upper bound from compact to closed sets~\cite[Lemma~1.2.18]{DemboZeitouni1998}. 
\smallskip

In our context, it allows us to deduce pathwise large deviations from large deviations of the finite-dimensional marginals. If for each tuple $0\leq t_1<t_2<\dots<t_k$, the marginals $\left\{(X^n(t_1),\dots,X^n(t_k))\right\}_{n\in\mathbb{N}}$ satisfy large deviations in $E^k$ with rate function $\mathcal{I}_{t_1\dots t_k}:E^k\to[0,\infty]$, then the sequence of processes $\{X^n\}_{n\in\mathbb{N}}$ satisfies a large deviation principle in $\mathcal{X}$ with rate function $\mathcal{I}:\mathcal{X}\to[0,\infty]$ given by
\begin{equation}\label{eq:RF-sup-over-finite-dim-distr}
\mathcal{I}(x) := \sup_{k\in\mathbb{N}}\sup_{(t_1,\dots,t_k)}\mathcal{I}_{t_1\dots t_k}\left(x_1(t_1),\dots,x_k(t_k)\right),
\end{equation}
where the supremum is taken over all finite tuples $t_1<t_2<\dots<t_k$.
A proof of this fact can be found in~\cite[Theorem~4.28]{FengKurtz2006}. The rate function~\eqref{eq:RF-sup-over-finite-dim-distr} is an example of bootstrapping large deviations from lower to higher-dimensional spaces, known as the Dawson-G\"artner Theorem~\cite[Theorem~4.6.1]{DemboZeitouni1998}.
We postpone the problem of how to obtain exponential tightness to Section~\ref{sec:LDP-via-generators}.
	\begin{proof}[Overview of the proof of Theorem~\ref{thm:LDP-via-semigroup-convergence:compact}]
	The pathwise large deviation principle follows from the large deviation principles of finite-dimensional distributions by exponential tightness~\cite[Theorem~4.28]{FengKurtz2006}. We first prove in Section~\ref{BG:sec:LDP-1d} the large deviation principle for the one-dimensional time marginals~$X^n(t)$.
	Then we see how the argument iterates to finite-dimensional distributions in Proposition~\ref{prop:LDP-multidim-compact:semigroup}. That gives the pathwise large deviation principle of~$\{X^n\}_{n\in\mathbb{N}}$ with rate function given by the formula~\eqref{eq:RF-sup-over-finite-dim-distr}. Finally, we prove in Proposition~\ref{prop:LDP-semigroups:conditional-RF} the formula~\eqref{eq:RF-with-semigroup-approach} of the rate function.
	\end{proof}
\subsection{Varadhan's Lemma and Bryc's Formula}
The main point of this section is to formulate the equivalence of large deviations and asymptotic evaluation of integrals of continuous functions, since this equivalence will be our starting point for proving large deviations of finite-dimensional distributions. For the following theorems, we consider a sequence of probability measures~$\mathrm{Q}_n$ on a compact Polish space~$\mathcal{S}$.
\begin{theorem}[Varadhan's Lemma]\label{thm:Varadhans-lemma}
	Suppose that the sequence $\{\mathrm{Q}_n\}_{n\in\mathbb{N}}$ satisfies a large deviation principle with rate function~$\mathcal{I}:\mathcal{S}\to[0,\infty]$. Then for any bounded and continuous function $f:\mathcal{S}\to\mathbb{R}$,
	\begin{equation*}
	\lim_{n\to\infty}\frac{1}{n}\log\int_\mathcal{S} e^{nf(x)}\,\mathrm{Q}_n\left(\dd x\right) = \sup_{x\in \mathcal{S}} \left[f(x)-\mathcal{I}(x)\right].
	\end{equation*}
\end{theorem}
\begin{theorem}[Bryc's Formula]\label{thm:Brycs-formula}
	Suppose that
	for any~$f\in C(\mathcal{S})$, the limit
	\begin{equation}\label{BG:eq:rate-transform}
	\Lambda(f) := \lim_{n\to\infty}\frac{1}{n}\log\int_\mathcal{S} e^{nf(x)}\,\mathrm{Q}_n\left(\dd x\right).
	\end{equation}
	exists. Then the sequence~$\{\mathrm{Q}_n\}_{n\in\mathbb{N}}$ satisfies a large deviation principle with rate function~$\mathcal{I}:\mathcal{S}\to[0,\infty]$ given by
	\begin{equation}\label{BG:eq:Brycs-formula}
	\mathcal{I}(x) = \sup_{f\in C(\mathcal{S})}\left[f(x)-\Lambda(f)\right].
	\end{equation}
\end{theorem}
Varadhan's Lemma is a generalization of winner-takes-it-all principle. For a continuous function $g:[a,b]\to\mathbb{R}$ on a closed interval $[a,b]$, we have
\[
\frac{1}{n}\log \int_a^b e^{ng(x)}\,\dd x \xrightarrow{n\to\infty} \sup_{x\in[a,b]} g(x) =:\overline{g}.
\]
This follows from the fact that~$(g-\overline{g})\leq 0$ on $[a,b]$ and
\[
\frac{1}{n}\log \int_a^b e^{ng(x)}\,\dd x = \overline{g}+\frac{1}{n}\log\int_a^be^{n(g(x)-\overline{g})}\,\dd x.
\]
Consider a sequence of probability measures $\mathrm{Q}_n\in\mathcal{P}(\mathbb{R})$ satisfying large deviations. At least intuitively, this means an approximation of the type
\[
\mathrm{Q}_n\left(\dd x\right) \approx e^{-n\mathcal{I}(x)}\,\dd x
\] 
is valid for large~$n$.
Then for a bounded and continuous function $f$ on $\mathbb{R}$, 
\[
\int_\mathbb{R} e^{nf(x)}\,\mathrm{Q}_n\left(\dd x\right) \approx \int_\mathbb{R} e^{nf(x)} e^{-n\mathcal{I}(x)}\,\dd x,\quad n\to\infty
\]
Hence on the logarithmic scale, we expect the dominant contribution to come from the maximal value of $g:=f-\mathcal{I}$. Varadhan's Lemma states both the existence of the limit and that it equals to what we expect from the winner-takes-it-all principle. 
Proofs closely following the above sketch are given for instance by Budhiraja and Dupuis~\cite[Theorem~1.5]{BudhirajaDupuis2019} and Frank den Hollander~\cite[Theorem~III.13]{denHollander2000}. A proof based on the exponential Chebyshev inequality is given by Feng and Kurtz in~\cite[Proposition~3.8]{FengKurtz2006}. Dembo and Zeitouni prove it in regular topological spaces under an additional tail bound~\cite[Theorem~4.3.1]{DemboZeitouni1998}.
\smallskip

Bryc proved the inverse to Varadhan's Lemma in~\cite{Bryc1990}.
The point of Bryc's formula is: if we can compute the rate transforms, then we obtain a rate function. Bryc's formula focuses on the existence of the limit~$\Lambda(f)$ and does not require to identify a rate function beforehand. We refer to the map~$\Lambda:C(\mathcal{S})\to\mathbb{R}$ as the \emph{rate transform} associated to $\{\mathrm{Q}_n\}_{n\in\mathbb{N}}$. By Varadhan's Lemma, we have
\begin{equation*}
\Lambda(f) = \sup_{x\in \mathcal{S}}\left[f(x)-\mathcal{I}(x)\right].
\end{equation*}
\subsection{Large deviations of one-dimensional marginals}
\label{BG:sec:LDP-1d}
Here we show how Bryc's formula directly leads to a proof of large deviations of the one-dimensional time marginales. To recall the setting of Theorem~\ref{thm:LDP-via-semigroup-convergence:compact}, we consider a sequence of Markov processes~$X^n$ with paths in~$\mathcal{X}=C_E[0,T]$, where~$E$ is a compact Polish space. For any~$t\in[0,T]$, the time marginal~$X^n(t)$ is a random variable in~$E$. We denote its distribution by~$\mathbb{P}^n_t\in\mathcal{P}(E)$. 
\smallskip

Fix~$t\in[0,T]$. By Bryc's Formula, the sequence $\{X^n(t)\}_{n \in \mathbb{N}}$ satisfies a large deviation principle if for any $f\in C(E)$, the sequence
\begin{equation*}
	\Lambda^n_t(f) := \frac{1}{n}\log \int_E e^{nf(y)}\,\mathbb{P}_t^n\left(\dd y\right)
\end{equation*}
converges as $n$ tends to infinity. First, let us suppose that the initial condition is deterministic, that means $X^n(0) \sim \mathbb{P}^n_0 := \delta_{x_0}$ for some fixed~$x_0\in E$. 
\smallskip

Let~$S_n$ be the semigroup corresponding to~$X^n$. For every $f\in B(E)$, we have by conditioning (e.g.~\cite[Proposition~4.1.6]{EthierKurtz1986} or~\cite[Definition~1.6]{Liggett2004})
\begin{equation}\label{BG:eq:setting:1d-distr-by-conditioning}
\int_E f(y)\,\mathbb{P}_t^n\left(\dd y\right)
= \int_E S_n(t)f(x)\,\mathbb{P}_0^n\left(\dd x\right).
\end{equation}
Therefore
\begin{equation}\label{BG:eq:Vn-condition-to-P0}
\int_E e^{n f(y)}\,\mathbb{P}^n_t\left(\dd y\right)
\overset{\eqref{BG:eq:setting:1d-distr-by-conditioning}}{=} 
\int_E S_n(t) e^{nf(x)}\,\mathbb{P}^n_0(\dd x) 
\myeqdef 
\int_E e^{n V_n(t)f(x)}\,\mathbb{P}^n_0\left(\dd x\right).
\end{equation}
Hence using~$\mathbb{P}^n_0=\delta_{x_0}$, we find
\begin{align*}
	\Lambda_t^n(f)&\myeqdef \frac{1}{n}\log \int_E e^{nf(y)}\,\mathbb{P}_t^n\left(\dd y\right)
	\overset{\eqref{BG:eq:Vn-condition-to-P0}}{=} \frac{1}{n}\log \int_E e^{nV_n(t)f(x)}\,\mathbb{P}_0^n(\dd x) \\
	&=\frac{1}{n}\log e^{nV_n(t)f(x_0)}= V_n(t)f(x_0).
\end{align*}
This is how the semigroups~$V_n(t)$ arise directly from Bryc's formula. Recall that we assume the convergence~$V_n(t)\to V(t)$ as~$n\to\infty$. Hence with the special initial condition~$\mathbb{P}_0^n=\delta_{x_0}$, the rate transform~$\Lambda(f)$ from~\eqref{BG:eq:rate-transform} is
\begin{equation*}
	\Lambda(f) = V(t)f(x_0).
\end{equation*}
By Bryc's formula~\eqref{BG:eq:Brycs-formula}, the rate function $\mathcal{I}_t(\cdot|x_0) : E \to [0,\infty]$ takes the form
\begin{equation*}
	\mathcal{I}_t(x|x_0) = \sup_{f\in C(E)} \left[f(x) - V(t)f(x_0)\right].
\end{equation*}
This is the conditional rate function introduced in~\eqref{BG:eq:RF-1d-marginals}. We just proved that the conditional probability measures
\begin{equation*}
	A\mapsto P_n(t,x_0,A) = \mathbb{P}\left(X^n(t)\in A\,|\,X^n(0)=x_0\right)
\end{equation*}
satisfy a large deviation principle with rate function~$\mathcal{I}_t(\cdot|x_0)$.
\smallskip

Both the fact that we only need convergence of $V_n(t)f$ at the point $x_0$ and that the limit $\Lambda(f)$ depends on $x_0$ are an artefact of the special form of the initial distribution, $X^n(0) \sim \delta_{x_0}$.
%
\begin{proposition}\label{prop:LDP-1d-compact:semigroup}
	Let~$t\in[0,T]$. Under the conditions of Theorem~\ref{thm:LDP-via-semigroup-convergence:compact}, the sequence of one-dimensional time marginals $\{X^n(t)\}_{n\in\mathbb{N}}$ satisfies a large deviation principle in~$E$ with rate function $\mathcal{I}_t:E\to[0,\infty]$ given by
	\begin{equation*}
	\mathcal{I}_t(x) = \sup_{f\in C(E)}\left\{f(x)-\Lambda_0\left[V(t)f\right]\right\},
	\end{equation*}
	where $\Lambda_0$ is the rate transform~\eqref{BG:eq:rate-transform} associated to the initial conditions~$\{X^n(0)\}$.
\end{proposition}
\begin{proof}[Proof of Proposition~\ref{prop:LDP-1d-compact:semigroup}]
	Let~$\mathbb{P}_t^n\in\mathcal{P}(E)$ be the distribution of~$X^n(t)$. If for any function~$f\in C(E)$ the rate transform
\begin{equation*}
	\Lambda_t(f) := \lim_{n\to\infty}\frac{1}{n}\log\int_E e^{n f(y)}\,\mathbb{P}^n_t\left(\dd y\right)
\end{equation*}
exists, then by Bryc's formula,~$\{X^n(t)\}$ satisfies a large deviation principle with rate function~$\mathcal{I}_t:E\to[0,\infty]$ given by
\begin{equation*}
	\mathcal{I}_t(x) = \sup_{f\in C(E)}\left[f(x)-\Lambda_t(f)\right].
\end{equation*}
Since the initial conditions~$\{X^n(0)\}$ satisfy large deviations, the rate transform
\begin{equation*}
	\Lambda_0(g):= \lim_{n\to\infty} \frac{1}{n}\log\int_E e^{ng(y)}\,\mathbb{P}^n_0\left(\dd y\right),\quad g\in C(E),
\end{equation*}
exists by Varadhan's Lemma.
Hence we can prove the Proposition by showing that~$\Lambda_t(f)$ is equal to~$\Lambda_0\left[V(t)f\right]$. Let~$f\in C(E)$. Then
\begin{align*}
\int_E e^{n f(y)}\,\mathbb{P}^n_t\left(\dd y\right) 
\overset{\eqref{BG:eq:Vn-condition-to-P0}}{=} \int_Ee^{n V_n(t)f(x)}\,\mathbb{P}^n_0\left(\dd x\right).
\end{align*}
The functions~$h_n:=V_n(t)f\in B(E)$ converge by assumption uniformly to the function~$h:=V(t)f\in C(E)$. The map $g\mapsto \Lambda_0^n(g) := \frac{1}{n}\log\int_E e^{ng(x)}\,\mathbb{P}_0^n\left(\dd x\right)$ is well-defined on~$B(E)$ and satisfies the bounds
\begin{equation*}
	-\|h-h_n\|_E + \Lambda_0^n(h)  \leq \Lambda_0^n(h_n) \leq \Lambda_0^n(h) + \|h-h_n\|_E.
\end{equation*}
Now the equality $\Lambda_t(f) = \Lambda_0\left[V(t)f\right]$ follows by taking the limit~$n\to\infty$. This last step requires the limit $V(t)f$ to be a continuous function on~$E$, since Varadhan's Lemma a priori only guarantees the rate transform~$\Lambda_0$ on continuous functions.
\end{proof}
\begin{running_example*}[Small diffusion]
	We illustrate the above result for the process $\dd X^n_t=n^{-1/2}\dd B_t$ on~$E=\mathbb{R}$, ignoring for the moment the fact that~$\mathbb{R}$ is not compact. The transition probabilities~$P_n(t,x,\dd y)$ are explictly known,
	\begin{equation*}
		P_n(t,x,\dd y) = \sqrt{\frac{n}{2\pi t}} \exp\{-n(y-x)^2/2t\}\,\dd y.
	\end{equation*}
	Hence the nonlinear semigroups are
	\[
	V_n(t)f(x) = \frac{1}{n}\log\left(\int_\mathbb{R}\exp\left\{n\left[f(y)-\frac{1}{2t}(y-x)^2\right]\right\}\,\dd y\right) + \frac{1}{n}\log\left(\sqrt{\frac{n}{2\pi t}}\,\right).
	\]
	The second term vanishes in the limit $n\to\infty$. In the integral term, the highest value dominates in the limit, which gives
	\[
	V(t)f(x) \myeqdef \lim_{n\to\infty} V_n(t)f(x) = \sup_{z\in\mathbb{R}}\left[f(z)-\frac{1}{2t}(z-x)^2\right].
	\]
	This expression is the well-known Hopf-Lax formula. Thus~$X^n(t)$ conditioned to~$X^n(0)=x$ satisfies a large deviation principle with rate function
	\begin{align*}
		\mathcal{I}_t(y|x) &= \sup_{f\in C_b(\mathbb{R})}\left[f(y) - V(t)f(x)\right]
	\end{align*}
	Inserting~$V(t)$ and evaluating, we find by proving two inequalities that
	\begin{align*}
	\mathcal{I}_t(y|x) &= \sup_{f\in C_b(\mathbb{R})}\inf_{z\in\mathbb{R}}\left[f(y)-f(z)+\frac{1}{2t}(z-x)^2\right] = \frac{1}{2t}(y-x)^2.
	\end{align*}
	This confirms what we can readily see from the transition probabilites.
\end{running_example*}
\subsection{Large deviations of finite-dimensional marginals}
We first consider two-dimensional time marginals $\{(X^n(t_1),X^n(t_2))\}$ for some fixed $t_2>t_1\geq 0$. To that end, denote the distribution of~$(X^n(t_1),X^n(t_2))$ by~$\mathbb{P}^n_{t_1t_2}\in\mathcal{P}(E\times E)$. We copy the strategy of one-dimensional marginals based on Bryc's formula: for proving large deviations of~$\{(X^n(t_1),X^n(t_2))\}_{n\in\mathbb{N}}$, it is sufficient to prove for any $f\in C(E\times E)$ the existence of the following limit:
\begin{equation*}
\Lambda_{t_1t_2}\left(f\right) := \lim_{n\to\infty}\frac{1}{n}\log\int_{E\times E}e^{n f(x,y)}\,\mathbb{P}^n_{t_1t_2}\left(\dd x\dd y\right).
\end{equation*}
We would like to use conditioning in order to reduce this convergence problem to convergence of the nonlinear semigroups~$V_n(t)$, as in the proof regarding one-dimensional marginals. To that end, we would like to consider only functions of the form $f_{12}(y,z) = f_1(y)+f_2(z)$, with~$f_1,f_2\in C(E)$. The fact that proving convergence for functions of this form is sufficient is the content of the following Lemma. It can be seen as strengthening Bryc's formula for Cartesian products. To shorten the presentation, we just indicate below where to find the proof.
\begin{lemma}[Considering sums is sufficient]\label{lemma:LDP-multidim-distr:sums-are-good}
	Let~$\mathcal{S}_1,\mathcal{S}_2$ be a compact Polish spaces and $\{\PP_n\}_{n\in\mathbb{N}}$ be a sequence of probability measures on~$\mathcal{S}_1\times \mathcal{S}_2$. For~$f_1\in C(\mathcal{S}_1)$ and $f_2\in C(\mathcal{S}_2)$, we write~$f_{12}\in C\left(\mathcal{S}_1\times \mathcal{S}_2\right)$ for the function~$f_{12}(y,z):=f_1(y)+f_2(z)$. Suppose that for any~$f_1,f_2$, the rate transform
	\[
	\Lambda(f_{12}) := \lim_{n\to\infty}\frac{1}{n}\log\int_{\mathcal{S}_1\times\mathcal{S}_2}e^{nf_1(y) + nf_2(z)}\,\PP_n\left(\dd y\dd z\right) 
	\]
	exists. Then the family $\left\{\PP_n\right\}_{n\in\mathbb{N}}$ satisfies a large deviation principle with rate function $\mathcal{I}:\mathcal{S}_1\times \mathcal{S}_2\to[0,\infty]$ given by
	\[
	\mathcal{I}(y,z) = \sup_{\substack{f_1\in C(\mathcal{S}_1)\\ f_2\in C(\mathcal{S}_2)}}\left[f_1(y) + f_2(z) - \Lambda(f_{12})\right].
	\]
\end{lemma}
\begin{proof}[Sketch of proof of Lemma~\ref{lemma:LDP-multidim-distr:sums-are-good}]
	If two subsets of functions $D_1\subseteq C(\mathcal{S}_1)$ and~$D_2\subseteq C(\mathcal{S}_2)$ are bounded above and isolates points, then the set of functions on $\mathcal{S}_1\times \mathcal{S}_2$ defined by $D_{12}:=\{f_1+f_2\,|f_1\in D_1,f_2\in D_2\}$ is bounded above and isolates points~\cite[Lemma~3.22]{FengKurtz2006}. Hence by Proposition~3.20 of~\cite{FengKurtz2006}, the set
	$F_{12}:=\{f_1+f_2\,:\,f_1\in C(\mathcal{S}_1),\,f_2\in C(\mathcal{S}_2)\}\subseteq C(\mathcal{S}_1\times\mathcal{S}_2)$ contains a set that is bounded above and isolates points. Therefore~$F_{12}$ 
	is rate-function determining in the sense of Definition~3.15 of~\cite{FengKurtz2006}.
\end{proof}
\begin{proposition}\label{prop:LDP-multidim-compact:semigroup}
	Let~$0<t_1<t_2\leq T$. Under the conditions of Theorem~\ref{thm:LDP-via-semigroup-convergence:compact}, the sequence of two-dimensional time marginals $\{\left(X^n(t_1),X^n(t_2)\right)\}_{n\in\mathbb{N}}$ satisfies a large deviation principle with rate function $\mathcal{I}_{t_1t_2}:E\times E\to[0,\infty]$ given by
	\begin{equation*}
		\mathcal{I}_{t_1t_2}(y,z) = \sup_{f_1,f_2\in C(E)}\left\{f_1(y)+f_2(z)-\Lambda_0\left[V(t_1)\left(f_1 + V(t_2-t_1)f_2\right)\right]\right\}.
	\end{equation*}
\end{proposition}
\begin{proof}[Proof of Proposition~\ref{prop:LDP-multidim-compact:semigroup}]
	Let~$\mathbb{P}_{t_1t_2}^n$ be the distribution of~$(X^n(t_1),X^n(t_2))$. We know by Lemma~\ref{lemma:LDP-multidim-distr:sums-are-good}, if for any function of the form~$f_{12}(y,z) := f_1(y) + f_2(z)$ with functions~$f_1,f_2\in C(E)$ the rate transform
	\begin{equation*}
		\Lambda_{t_1t_2}\left(f_{12}\right) := \lim_{n\to\infty}\frac{1}{n}\log\int_{E\times E}e^{nf_{12}(y,z)}\,\mathbb{P}^n_{t_1t_2}\left(\dd y\dd z\right)
	\end{equation*}
	exists, then the large deviation principle holds with rate function
	\begin{equation*}
		\mathcal{I}_{t_1t_2}(y,z) = \sup_{f_1, f_2 \in C(E)}\left\{f_1(y) + f_2(z) - \Lambda_{t_1t_2}\left[f_{12}\right]\right\}.
	\end{equation*}
	The initial conditions~$\{X^n(0)\}$ satisfy a large deviation principle by assumption. Hence by Varadhan's Lemma, for any~$g\in C(E)$ the rate transform
	\begin{equation*}
		\Lambda_0(g):= \lim_{n\to\infty} \frac{1}{n}\log\int_E e^{ng(y)}\,\mathbb{P}^n_0\left(\dd y\right)
	\end{equation*}
	exists. We have~$V(t_1)(f_1+V(t_2-t_1)f_2)\in C(E)$ since~$V(t)$ is a map from~$C(E)$ to~$C(E)$.
	Therefore, the proposition follows if we prove
	\begin{equation*}
		\Lambda_{t_1t_2}\left[f_{12}\right]=\Lambda_0\left[V(t_1)\left(f_1+V(t_2-t_1)f_2\right)\right].
	\end{equation*} 
	As in the proof of one-dimensional distributions, we condition to earlier times (\cite[Proposition~4.1.6]{EthierKurtz1986}), and find
	\begin{align*}
	\int_{E\times E}e^{n (f_1(y)+f_2(z))}\,\mathbb{P}_{t_1t_2}^n\left(\dd y\dd z\right)
	&=
	\int_E e^{n (f_1(y) + V_n(t_2-t_1)f_2(y))}\,\mathbb{P}_{t_1}^n(\dd y)\\
	&=
	\int_E e^{n V_n(t_1)\left[f_1 + V_n(t_2-t_1)f_2\right](x)}\,\mathbb{P}_0^n\left(\dd x\right).
	\end{align*}
	By the convergence assumption on the nonlinear semigroups $V_n(t)$,
	\begin{equation*}
	f_1 + V_n(t_2-t_1)f_2 \xrightarrow{n\to\infty} f_1 + V(t_2-t_1)f_2
	\end{equation*}
	uniformly on $E$. Hence again by the convergence assumption,
	\begin{equation*}
	h_n:= V_n(t_1)\left[f_1 + V_n(t_2-t_1)f_2\right] \xrightarrow{n\to\infty} h:=V(t_1)\left[f_1 + V(t_2-t_1)f_2\right].
	\end{equation*}
	The map $g\mapsto \Lambda_0^n(g) := \frac{1}{n}\log\int_E e^{ng(x)}\,\mathbb{P}_0^n\left(\dd x\right)$ is well-defined on~$B(E)$ and satisfies the bounds
	\begin{equation*}
		-\|h-h_n\|_E + \Lambda_0^n(h)  \leq \Lambda_0^n(h_n) \leq \Lambda_0^n(h) + \|h-h_n\|_E,
	\end{equation*}
	and the desired equality follows by taking the limit~$n\to\infty$.
\end{proof}
The convergence condition on the nonlinear semigroups $V_n(t)$ is sufficient for iterating to finite-dimensional marginals $\left(X^n(t_1),\dots,X^n(t_k)\right)$. The rate function is then given by
\begin{equation}
\label{eq:LDP-semigroup-multidim:RF-multidim}
\mathcal{I}_{t_1\dots t_k}(x_1,\dots,x_k) = \sup_{f_1\dots f_k \in C(E)}\left\{\sum_i f_i(x_i)-\Lambda_{t_1\dots t_k}\left[f_1,\dots,f_k\right]\right\},
\end{equation}
where the rate transform $\Lambda_{t_1\dots t_k}$ includes concatinations of the limiting map~$V(t)$,
\[
\Lambda_{t_1\dots t_k}\left[f_1,\dots,f_k\right] = \Lambda_0\left[V(t_1)\left(f_1 + V(t_2-t_1)(f_2+\dots+V(t_k-t_{k-1})f_k)\dots))\right)\right].
\]
\subsection{Conditional structure of the rate function}
In this section, we show how to cast the rate function for finite-dimensional distributions from~\eqref{eq:LDP-semigroup-multidim:RF-multidim} into the more convenient form~\eqref{eq:prop:condtional-RF-finite-dim-marginals} given below. 
\begin{proposition}[]\label{prop:LDP-semigroups:conditional-RF}
	For $t_k>t_{k-1}>\dots>t_1> = 0$, consider the rate function $\mathcal{I}_{t_1\dots t_k}:E^k\to[0,\infty]$ of finite-dimensional time-marginals $\left\{X^n(t_1),\dots,X^n(t_k)\right\}_{n\in\mathbb{N}}$ given by~\eqref{eq:LDP-semigroup-multidim:RF-multidim}.
	Then
	\begin{equation}\label{eq:prop:condtional-RF-finite-dim-marginals}
	\mathcal{I}_{t_1\dots t_k}(x_1,\dots,x_k) = \mathcal{I}_{t_1}(x_1) + \mathcal{I}_{t_2-t_1}(x_2|x_1) + \dots + \mathcal{I}_{t_k-t_{k-1}}(x_{k}|x_{k-1}),
	\end{equation}
	where $\mathcal{I}_{t}$ is the rate function for $X^n(t_1)$ and the conditional rate functions $\mathcal{I}_{t}(z|y)$ are
	\begin{equation}\label{eq:LDP-semigroups:conditional-RF}
	\mathcal{I}_t(z|y) = \sup_{f\in C(E)}\left[f(z)-V(t)f(y)\right].
	\end{equation}
\end{proposition}
\begin{proof}[Proof of Proposition~\ref{prop:LDP-semigroups:conditional-RF}]
	We consider the case $k=2$. Then
	\[
	\mathcal{I}_{t_1t_2}(x_1,x_2) = \sup_{f_1,f_2 \in C(E)}\left\{\,f_1(x_1)+f_2(x_2)-\Lambda_0\left[V(t_1)\left(f_1+V(t_2-t_1)f_2\right)\right]\,\right\}.
	\]
	Concatinating the supremum and adding zero, we obtain
	\begin{multline*}
		\mathcal{I}_{t_1t_2}(x_1,x_2) = \sup_{f_2\in C(E)}\bigg[ f_2(x_2) - V(t_2-t_1)f_2(x_1) \\
		+\sup_{f_1\in C(E)}f_1(x_1) + V(t_2-t_1)f_2(x_1) - \Lambda_0\left[V(t_1)\left(f_1+V(t_2-t_1)f_2\right)\right] \bigg].
	\end{multline*}
	Since $V(t):C(E)\to C(E)$, we may shift in the second supremum to functions of the form $f_1 = g_1 - V(t_2-t_1)f_2$, with $g_1\in C(E)$, to obtain
	\begin{align*}
	\mathcal{I}_{t_1t_2}(x_1,x_2) &= \sup_{f_2}\left[ f_2(x_2) - V(t_2-t_1)f_2(x_1)\right] + \sup_{g_1}\left[g_1(x_1)-\Lambda_0\left[V(t_1)g_1\right]\right]\\
	&\myeqdef \mathcal{I}_{t_2-t_1}\left(x_2|x_1\right) + \mathcal{I}_{t_1}(x_1).
	\end{align*}
	This finishes the proof for $k=2$. Similarly, we obtain for $k=3$
	\begin{align*}
		\mathcal{I}_{t_1t_2t_3}(x_1,x_2,x_3) &= \mathcal{I}_{t_3-t_2}(x_3|x_2) + \mathcal{I}_{t_1t_2}(x_1,x_2).
	\end{align*}
	The general case follows by induction.
\end{proof}
Proposition~\ref{prop:LDP-semigroups:conditional-RF} represents the fact that for a Markov process $X^n$, the time marginals such as~$X^n(t_1)$ and~$X^n(t_2)$ for~$t_1<t_2$ are in general not independent, but correlated. From the large deviation principles
\[
\mathbb{P}\left[X^n(t_1)\approx x_1\right]\sim e^{-n \mathcal{I}_{t_1}(x_1)}\quad\text{and}\quad
\mathbb{P}\left[X^n(t_2)\approx x_2\right]\sim e^{-n \mathcal{I}_{t_2}(x_2)},
\]
we can not conclude the large deviation principle of the pair~$(X^n(t_1),X^n(t_2)$ as
\[
\mathbb{P}\left[X^n(t_1)\approx x_1,X^n(t_2)\approx x_2\right] \overset{?}{\sim} e^{-n \left[\mathcal{I}_{t_1}(x_1)+\mathcal{I}_{t_2}(x_2)\right]}.
\]
Rather, the rate functions reflect the fact that the event $X^n(t_1)\approx x_1$ takes place before the event $X^n(t_2)\approx x_2$. This condition appears in the rate function~$\mathcal{I}_{t_1t_2}$ of the joint event:
\[
\mathbb{P}\left[X^n(t_1)\approx x_1,X^n(t_2)\approx x_2\right]\sim
e^{-n\left[\mathcal{I}_{t_1}(x_1) + \mathcal{I}_{t_2-t_1}(x_2|x_1)\right]}.
\]
\begin{running_example*}[Small diffusion]
	Let $E=\mathbb{R}$ and $\dd X^n_t = n^{-1/2}\dd B_t$, and suppose $X^n(0)=x_0$. Again we ignore that~$\mathbb{R}$ is not compact. We find
	\begin{equation*}
	\mathbb{P}\left[X^n(t_1)\approx x_1,X^n(t_2)\approx x_2\right] \sim e^{-n \left[\mathcal{I}_{t_1}(x_1) + \mathcal{I}_{t_2-t_1}(x_2|x_1)\right]},
	\end{equation*}
	where we already computed $\mathcal{I}_t(x_2|x_1) = (x_2-x_1)^2/2t$. Since the process starts at~$x_0$, we have $\mathcal{I}_{t_1}(x_1)=(x_1-x_0)^2/2t$. For a partition~$0=t_0<t_1<\dots t_k= T$ of the time interval~$[0,T]$,
	\begin{equation*}
		\mathbb{P}\left[X^n(t_1)\approx x_1,\dots,X^n(t_k)\approx x_k\right]\sim \exp\{-n \cdot \mathcal{I}_{t_1\dots t_k}(x_1,\dots,x_k)\}, \quad n\to\infty.
	\end{equation*}
	Suppose~$t_i-t_{i-1}\approx \Delta t>0$ is small. Then massaging the rate function a bit,
	\begin{align*}
		\mathcal{I}_{t_1\dots t_k}(x_1,\dots,x_k) &\overset{\eqref{eq:prop:condtional-RF-finite-dim-marginals}}{=} \mathcal{I}_{t_1}(x_1) + \mathcal{I}_{t_2-t_1}(x_2|x_1) + \dots + \mathcal{I}_{t_k-t_{k-1}}(x_{k}|x_{k-1})\\
		&= \frac{1}{2}\frac{(x_1-x_0)^2}{(t_1-t_0)} + \frac{1}{2} \frac{(x_2-x_1)^2}{(t_2-t_1)} + \dots + \frac{1}{2} \frac{(x_k-x_{k-1})^2}{(t_k-t_{k-1})}\\
		&\approx \frac{1}{2} \sum_{i=1}^k \left(\frac{x_{i}-x_{i-1}}{t_i-t_{i-1}}\right)^2 \Delta t.
	\end{align*}
	Hence with a fine partition and regarding the points~$x_i$ as corresponding to a path~$x:[0,T]\to E$ via~$x_i=x(t_i)$, we expect
	\begin{equation}
	\mathcal{I}_{t_1\dots t_k}(x_1,\dots,x_k) \overset{k\gg 1}{\approx} \int_0^T \frac{1}{2} (\partial_t x(t))^2\,\dd t.
	\end{equation}
	The rigorous version of this derivation is Schilder's Theorem, and we will prove the corresponding rigorous statement further below. This small calculation is based on the explicit formula for the conditional rate functions. In general, we will be able to obtain something like
	\begin{equation*}
	\mathcal{I}_{t_1\dots t_k}(x_1,\dots,x_k) \overset{k\gg 1}{\approx} \int_0^T \mathcal{L}(\partial_t x(t))\,\dd t,\quad \mathcal{L}(\cdot)\;\text{convex}.
	\end{equation*}
	We provide more details in Section~\ref{BG:sec:semigroup-flow-HJ-eq} below.
\end{running_example*}
\section{Large deviations via convergence of generators}
\label{sec:LDP-via-generators}
In the previous section, we introduced in Definition~\ref{def:nonlinear-semigroups} the nonlinear semigroups $V_n(t)$ associated to a Markov process $X^n$ with paths in~$\mathcal{X}=C_E[0,T]$. We summarized the main preliminary result in Theorem~\ref{thm:LDP-via-semigroup-convergence:compact}, which identifies two conditions for proving pathwise large deviation principles:
\begin{enumerate}[label=(\roman*)]
	\item The sequence~$\{X^n\}_{n\in\mathbb{N}}$ is exponentially tight.
	\item The semigroups $V_n(t)|_{t\geq 0}$ converge to a semigroup $V(t)|_{t\geq 0}$.
\end{enumerate}
We say that Theorem~\ref{thm:LDP-via-semigroup-convergence:compact} is preliminary for a couple of reasons: 
\begin{enumerate}[$\bullet$]
	\item Verifying exponential tightness is a nasty and unfortunate task that we would like to avoid carrying out on a case-by-case analysis.
	\item Typically, the nonlinear semigroups $V_n(t)$ are not computable, and it is hard to even \emph{identify} a possible limit candidate $V(t)$ in the first place, yet proving convergence.
	\item The formula~\eqref{eq:RF-with-semigroup-approach} for the rate function is complicated. Even a simple question like "what is its minimizer?" is hard to answer. 
\end{enumerate}
In this section, we answer the second key question from Section~\ref{BG:sec:math-gem}: how can we verify the convergence of semigroups from convergence of generators? To that end, let us turn to the recipe of semigroup convergence as outlined in Example~\ref{ex:semigroups-are-exponentials} above: we want to identify the generator~$H_n=(d/dt)V_n(0)$ and then identify a suitable limit~$H=\lim_n H_n$. Finally, we hope to conclude the semigroup convergence~$V_n\to V$. We start with deriving~$H_n$.
\smallskip

For a Markov process with semigroup~$S(t)$, the generator~$L$ is a linear operator characterizing the infinitesimal time evolution by
\begin{equation*}
S(t+\Delta t)f(x) = \mathbb{E}\left[f(X(t+\Delta t))|X(t)=x\right] = f(x) + Lf(x) \Delta t + \mathcal{O}(\Delta t^2).
\end{equation*}
\begin{definition}[Infinitesimal generator]
	\label{BG:def:infinitesimal-generator}
	Consider a strongly continuous contraction semigroup $S(t):C(E)\to C(E)$. Its corresponding \emph{infinitesimal generator} $L$ is a linear operator $L:\mathcal{D}(L)\subseteq C(E)\to C(E)$, where for any~$f\in C(E)$, if there is some~$g\in C(E)$ such that uniformly on~$E$,
	\[
	g = \lim_{t\to 0}\frac{1}{t}\left(S(t)f-f\right),
	\]
	then $f\in \mathcal{D}(L)$ and $Lf:= g$.
	\qed
\end{definition}
We consider a Markov process~$X^n\in C_E[0,T]$ with corresponding transition probabilities~$P_n(t,x,\dd y)$ and linear semigroup $S_n(t)f(x) := \int_Ef(y)P_n(t,x,\dd y)$. Let~$L_n$ be its infinitesimal generator. For any function~$f\in\mathcal{D}(L)$,
\begin{equation*}
	L_n f = \frac{\dd}{\dd t}\bigg|_0 S_n(t)f.
\end{equation*}
The nonlinear semigroups $V_n(t)$ from Definition~\ref{def:nonlinear-semigroups} are given by
\begin{equation*}
V_n(t)f(x) = \frac{1}{n}\log S_n(t)e^{nf(\cdot)}|_x.
\end{equation*}
Taking the time derivative and evaluating at zero, the chainrule formally yields
\begin{align*}
	\frac{\dd}{\dd t}\bigg|_{t=0} V_n(t)f(x) &=  \frac{1}{n}\frac{1}{S_n(0)e^{nf(\cdot)}|_x}\frac{\dd}{\dd t}\bigg|_{t=0}S_n(t)e^{nf(\cdot)}|_x\\
	&= \frac{1}{n}e^{-nf(x)}L_n e^{nf(x)}.
\end{align*}
This suggests the operators~$H_nf := n^{-1}e^{-nf}L_n e^{nf}$ are the generators of~$V_n(t)$.
\begin{running_example*}
	Let~$E=\mathbb{R}$ and~$\dd X^n_t = n^{-1/2}\dd B_t$. The linear generator is~$L_nf=(2n)^{-1}\Delta f$ with domain~$\mathcal{D}(L_n)=C_b^2(\mathbb{R})$. We find
	\begin{equation*}
	H_nf \myeqdef \frac{1}{n}e^{-nf}L_n e^{nf} = \frac{1}{2}\frac{1}{n}\Delta f + \frac{1}{2}|\nabla f|^2,
	\end{equation*}
	writing $\nabla f = f'$ and $\Delta f = f''$.
	\qed
\end{running_example*}
\begin{definition}[Nonlinear generators]\label{def:nonlinear-generators}
	Let~$E$ be a compact Polish space and let a linear operator~$L_n:\mathcal{D}(L_n)\subseteq C(E)\to C(E)$ be the generator of an $E$-valued Markov process. The corresponding \emph{nonlinear generator} $H_n$ is defined as the map
	\begin{equation}\label{BG:eq:nonlinear-generators}
	H_nf (x):= \frac{1}{n}e^{-nf(x)}L_ne^{nf(\cdot)}(x),
	\end{equation}
	defined on the domain $\mathcal{D}(H_n):=\{f\,|\,e^{nf}\in\mathcal{D}(L_n)\}$. 
	\qed
\end{definition}
Here, the operators~$H_n$ have to be understood as \emph{formal} generators of~$V_n(t)$. We only took the above calculation as a motivation, but do not claim the nonlinear generator to be a generator in the mathematically precise sense as for instance in the Hille-Yosida Theorem. Also in~\cite{FengKurtz2006} it is never claimed that we can make precise sense of $\frac{d}{dt}V_n(0)=H_n$. The formal calculations merely suggest that the limiting behaviour of~$H_n$ is closely related to the limiting behaviour of~$V_n(t)$. 
Jump processes form an important exception, where we will indeed find the relation~$\frac{d}{dt}V_n(0)=H_n$.
\smallskip

Equipped with Definition~\ref{def:nonlinear-generators}, we can tackle the task of deriving semigroup convergence from generator convergence.
	In Section~\ref{subsec:LDP-via-classical-sol}, we find conditions under which convergence of the nonlinear generators~$H_n$ to a limiting operator~$H$ implies large deviations (Theorem~\ref{thm:LDP-classical-sol}). 
	The main ingredient of the proof of Theorem~\ref{thm:LDP-classical-sol} is a convergence statement that translates the Trotter-Kato approximation theorem for linear semigroups to the nonlinear setting. We call this convergence statement the Feng-Kurtz approximation theorem. In addition, we have to pose conditions on the limit~$H$ in order to construct a semigroup~$V(t)$ from it. Below the proof of Theorem~\ref{thm:LDP-classical-sol}, we illustrate with the running example which condition is hard to verify.
	\smallskip
	
	In Section~\ref{subsec:LDP-via-visc-sol}, we motivate the consideration of viscosity solutions. With this type of weak solutions, the conditions on the limit~$H$ are verifiable. The summarize this main result in Theorem~\ref{thm:LDP-visc-sol}, and we verify its conditions for the running example.
\subsection{Using classical solutions}\label{subsec:LDP-via-classical-sol}
The main point of this section is to show under which conditions convergence of generators implies the large deviation principle (Theorem~\ref{thm:LDP-classical-sol}).
In the following definitions, nonlinear operators~$H$ acting on Banach spaces~$B$ are regarded as subsets of~$B\times B$. We denote by~$\overline{H}$ the closure of~$H$ with respect to the graph norm. For a Polish space~$E$, we will consider the Banach space~$B(E)$ of measurable bounded functions on~$E$, equipped with the supremum norm denoted by~$\|\cdot\|$. The following two properties are posing solvability conditions on an equation of the type~$(1-\tau H)f=h$ for a nonlinear operator~$H$, where~$\tau>0$ and~$h(x)$ are given and a solution~$f(x)$ in the domain of~$H$ is sought. We say~$f$ is a \emph{classical solution} if~$f\in\mathcal{D}(H)$ and~$(1-\tau H)f=h$.
\begin{definition}[Dissipative operator]
	For a Polish space~$E$, a nonlinear operator $H\subseteq B(E)\times B(E)$ with domain $\mathcal{D}(H)$ is called \emph{dissipative} if for all $\tau>0$ and any $f_1,f_2\in\mathcal{D}(H)$, the following estimate is satisfied:
	\begin{align*}
	\|f_1-f_2\| \leq \|(f_1-\tau Hf_1) - (f_2-\tau Hf_2)\|.
	\tag*\qed
	\end{align*}
\end{definition}
Dissipativity corresponds to uniqueness of classical solutions.
For $\tau>0$ and $h\in B(E)$, suppose two functions~$f_1,f_2\in\mathcal{D}(H)$ satisfy $(1-\tau H)f_1=h$ and $(1-\tau H)f_2=h$. If $H$ is a dissipative operator, then $\|f_1-f_2\| \leq 0$. 
\begin{definition}[Range condition]\label{def:range-condition}
	Let $E$ be a Polish space and let $H$ be a nonlinear operator $H\subseteq B(E)\times B(E)$ with domain $\mathcal{D}(H)$. We say that $H$ satisfies the \emph{range condition} if there exists a $\tau_0>0$ such that for all $0<\tau <\tau_0$, we have
	\begin{align*}
	\mathcal{D}(H) \subseteq \overline{\text{range}\left(1-\tau H\right)}.
	\tag*\qed
	\end{align*}
\end{definition}
The range condition corresponds to the existence of classical solutions.
	For dissipative operators, we have $\overline{\text{range}(1-\tau H)}=\text{range}(1-\tau \overline{H})$. If a dissipative operator $H$ satisfies the range condition, then for any $h\in\mathcal{D}(H)$ and $\tau>0$ sufficiently small, there exists a function $f\in\mathcal{D}(\overline{H})$ such that $(1-\tau \overline{H})f=h$.
\begin{theorem}[{Crandall-Liggett,~\cite{CrandallLiggett1971}}]
	Let $E$ be a Polish space and let $H$ be a nonlinear operator $H\subseteq B(E)\times B(E)$ with domain $\mathcal{D}(H)$. Suppose that $H$ is dissipative and satisfies the range condition. Then for each $f\in \overline{\mathcal{D}(H)}$, the map
	\[
	V(t)f := \lim_{k\to\infty}\left(1-\frac{t}{k}\overline{H}\right)^{-k}f 
	\]
	exists and $\{V(t)\}_{t\geq 0}$ defines a contraction semigroup $V(t):\overline{\mathcal{D}(H)}\to \overline{\mathcal{D}(H)}$.
\end{theorem}
In Example~\ref{ex:semigroups-are-exponentials}, we considered semigroups taking values in~$\mathbb{C}$. In that context, given a generator~$g\in\mathbb{C}$, we can use different equivalent formulas of the exponential map to associate a continuous semigroup $T(t)=e^{tg}$ to it:
\begin{equation*}
T(t) := \sum_{k=0}^\infty \frac{1}{k!}(tg)^k\quad \text{or}\quad T(t) := \lim_{k\to\infty}\left(1-\frac{t}{k}\,g\right)^{-k}.
\end{equation*}
The first formula can be used in the context of linear semigroups whose generators are bounded operators (Hille-Yosida Theorem). The Crandall-Liggett Theorem is based on the second formula. For an operator~$H$ and~$\tau>0$, it uses the resolvent~$R(\tau) := \left(1-\tau H\right)^{-1}$ defined by finding a unique solution~$f$ to~$(1-\tau H)f=h$ for each given~$\tau>0$ and~$h\in C(E)$. Then~$V(t)=\lim_k R(t/k)^k$ serves as the rigorous version of~$V(t)=e^{tH}$.
\smallskip

If $E$ is compact and we work with $\mathcal{D}(H)\subseteq C(E)$ dense, then the associated semigroup consists of maps $V(t):C(E)\to C(E)$. 
%
The following Theorem is a simplification of Proposition~5.5 in~\cite{FengKurtz2006}. 
\begin{theorem}[Feng-Kurtz approximation]\label{thm:feng-kurtz-approx}
	Let~$E$ be a compact Polish space and let $G_n:B(E)\to B(E)$ and $H:\mathcal{D}(H)\subseteq C(E)\to C(E)$ be two dissipative operators that both satisfy the range condition with the same~$\tau_0$. Let~$V_n(t)$ and~$V(t)$ be the corresponding generated semigroups in the Crandall-Liggett sense. Suppose the following:
	\begin{enumerate}[label=(\roman*)]
		\item For each~$f\in\mathcal{D}(H)$, there exist~$f_n \in B(E)$ such that
		\[
		\|f-f_n\|_E\xrightarrow{n\to\infty}0\quad\text{and}\quad
		\|Hf-G_nf_n\|_E\xrightarrow{n\to\infty}0.
		\]
	\end{enumerate}
	Then for any~$f\in\overline{\mathcal{D}(H)}$ and $f_n\in B(E)$ such that $\|f-f_n\|_E\to 0$, we have
	\begin{align*}
	\|V(t)f-V_n(t)f_n\|_E \xrightarrow{n\to\infty}0.
	\tag*\qed
	\end{align*}
\end{theorem}
We will apply the Feng-Kurtz approximation to operators~$G_n$ that are Hille-Yosida approximations of Markov generators. These Hille-Yosida approximations are generators of jump processes and satisfy the conditions of the Feng-Kurtz approximation theorem. A detailed discussion on their construction is given in~\cite[Section~IV.2]{EthierKurtz1986}.
The fact that they are dissipative and satisfy the range condition is proven in~\cite[Lemma~5.7]{FengKurtz2006}. We now use Lemmas of~\cite[Chapter~5]{FengKurtz2006} to prove the following simplification of~\cite[Corollary~5.19]{FengKurtz2006}.
%
\begin{theorem}[Large deviations via classical solutions]\label{thm:LDP-classical-sol}
	Let $E$ be a compact Polish space and for $n=1,2,\dots$, let $L_n:\mathcal{D}(L_n)\subseteq C(E)\to C(E)$ be the generator of an $E$-valued Markov process $X^n_t|_{t\geq 0}$ with continuous sample paths in~$\mathcal{X}=C_E[0,T]$. Let $H_n$ be the nonlinear generators~\eqref{BG:eq:nonlinear-generators}. Suppose the following:
	\begin{enumerate}[label=(\roman*)]
		\item \label{item:LDP-classical-sol:H-to-Hn}There exists a densely defined operator $H:\mathcal{D}(H)\subseteq C(E)\to C(E)$ such that for any $f\in \mathcal{D}(H)$, there are $f_n\in\mathcal{D}(H_n)$ satisfying
		\[
		\|f-f_n\|_E\xrightarrow{n\to\infty}0 \quad \text{and}\quad \|Hf-H_nf_n\|_E \xrightarrow{n\to\infty}0.
		\]
		\item \label{item:LDP-classical-sol:range} The operator $H$ satisfies the range condition (Definition~\ref{def:range-condition}).
	\end{enumerate}
	Suppose furthermore that $X^n(0)$ satisfies the large deviation principle in~$E$ with rate function $\mathcal{I}_0:E\to[0,\infty]$.
%
	Then~$H$ satisfies the conditions of the Crandall-Liggett Theorem and hence generates a semigroup~$V(t)$, and the sequence $\{X^n\}_{n\in\mathbb{N}}$ satisfies the large deviation principle in $\mathcal{X}$ with a rate function $\mathcal{I}:\mathcal{X}\to[0,\infty]$ given by~\eqref{eq:RF-with-semigroup-approach}.
\end{theorem}
\begin{proof}[Sketch of proof of Theorem~\ref{thm:LDP-classical-sol}]
	We verify the conditions of Theorem~\ref{thm:LDP-via-semigroup-convergence:compact}, according to which a large deviation principle of~$\{X^n\}$ follows from two conditions: exponential tightness and convergence of the nonlinear semigroups $V_n(t)$ to some limiting semigroup~$V(t):C(E)\to C(E)$.
	\smallskip
	
	Under the above convergence condition on the nonlinear generators~$H_n$, exponential tightness of~$\{X^n\}_{n\in\mathbb{N}}$ follows by~\cite[Corollary~4.17]{FengKurtz2006}. We do not give the details here, but comment briefly on why: (a) the exponential compact containment condition is always satisfied for compact spaces, (b) we can take $F=C(E)$ since $\mathcal{D}(H)$ is dense in~$C(E)$, and (c) exploits the fact that by the convergence condition $H_n\to H$, the sequences $H_nf_n$ are uniformly bounded.
	\smallskip
	
	We are left with showing that 1) we can define~$V(t)$ in terms of the limit operator~$H$ by means of the Crandall-Liggett Theorem and 2) that we obtain the semigroup-convergence $V_n(t)\to V(t)$ as specified in Theorem~\ref{thm:LDP-via-semigroup-convergence:compact}. 
	\smallskip
	
	1) By assumption, the operator~$H$ satisfies the range condition, and we only need to verify dissipativity. To that end, we henceforth only work with the full generator of~$X^n$, the graph in~$B(E)\times B(E)$ defined as (see \cite[Section~1.1.5]{EthierKurtz1986})
		\begin{equation*}
		\left\{\left(f,g\right)\in B(E)\times B(E)\,:\,\;\forall \,t,\, S_n(t)f-f=\int_0^tS_n(s)g\,\dd s\right\}.
		\end{equation*}
	We will denote them as well by $L_n$, and their associated nonlinear generators as well by~$H_n$. The reason for considering the full generator is that by Proposition~5.1 in~\cite{EthierKurtz1986}, it is a linear dissipative operator with resolvent
	\begin{equation*}
	\left(\lambda - L_n\right)^{-1}h = \int_0^\infty e^{-\lambda t}S_n(t)h\,\dd t.
	\end{equation*} 
	%
	Consider for $\varepsilon_n:=\exp\{-n^2\}$ the Hille-Yosida approximations $L_n^{\varepsilon_n}$ defined by
	\[
	L_n^{\varepsilon_n} := L_n\left(1-\varepsilon_n L_n\right)^{-1}.
	\]
	The map~$L_n^{\varepsilon_n}:B(E)\to B(E)$ defines a bounded, linear and dissipative operator (\cite[Lemma~1.2.4]{EthierKurtz1986}) that generates a Markov jump process on~$E$. Define the associated nonlinear generators $H_n^{\varepsilon_n}:B(E)\to B(E)$ by
	\[
	H_n^{\varepsilon_n}f:= \frac{1}{n}e^{-nf}L_n^{\varepsilon_n}e^{nf}.
	\]
	Then~$H_n^{\varepsilon_n}$ is dissipative~\cite[Lemma~5.7]{FengKurtz2006}. We prove below that our Assumption~\ref{item:LDP-classical-sol:H-to-Hn} on the convergence $H_n\to H$ implies that $H_n^{\varepsilon_n}\to H$ in the same sense. That establishes dissipativity of~$H$ as the limit of the dissipative operators~$H_n^{\varepsilon_n}$; for any $f_1,f_2\in \mathcal{D}(H)$, let $f_1^n,f_2^n\in B(E)$ be such that $f_1^n\to f_1$ and $f_2^n\to f_2$ uniformly on~$E$. Then using dissipativity of $H_n^{\varepsilon_n}$ and that the corresponding images converge uniformly, we find that
	\begin{align*}
		\|f_1-f_2\|_E &\leq \|f_1^n-f_2^n\|_E + o(1)_{n\to\infty}\\
		&\leq \|\left(f_1^n-\tau H_n^{\varepsilon_n}f_1^n\right)-\left(f_2^n-\tau H_n^{\varepsilon_n}f_2^n\right)\|_E +o(1)_{n\to\infty}\\
		&\leq \|\left(f_1-\tau Hf_1\right) - \left(f_2-\tau Hf_2\right)\|_E + o(1)_{n\to\infty}.
	\end{align*}
	Now taking the limit~$n\to\infty$ shows that~$H$ is dissipative.
	\smallskip
	
	We are left with verifying $H_n^{\varepsilon_n}\to H$. For~$f\in \mathcal{D}(H)$, let~$f_n\in\mathcal{D}(H_n)$ be such that $f_n\to f$ and $H_nf_n\to Hf$, both uniformly on~$E$. Then since~$H_nf_n$ is bounded, $n\varepsilon_nH_nf_n\to 0$ as~$n\to\infty$. Hence $e^{nf_n}(1-n\varepsilon_n \, H_nf_n) > 0$ eventually. We show that the functions~$f_n^{\varepsilon_n}$ defined by
	\[
	e^{nf_n^{\varepsilon_n}} := e^{nf_n}\left(1-n\varepsilon_n H_n f_n\right) =\left(1-\varepsilon_n L_n\right) e^{nf_n}
	\]
	satisfy $f_n^{\varepsilon_n}\to f$ and $H_n^{\varepsilon_n}f_n^{\varepsilon_n}\to Hf$. The first convergence follows from the fact that~$f_n\to f$ and~$e^{n(f_n^{\varepsilon_n}-f_n)}\to 1$. We find by the definition of the Hille-Yosida approximants~$L_n^{\varepsilon_n}$ that
	\begin{align*}
	H_n^{\varepsilon_n}f_n^{\varepsilon_n} &= \frac{1}{n} e^{-nf_n^{\varepsilon_n}} L_n^{\varepsilon_n} (1-\varepsilon_nL_n) e^{nf_n}\\ &= e^{-nf_n^{\varepsilon_n}} \frac{1}{n}L_n e^{nf_n} = e^{n(f_n-f_n^{\varepsilon_n})} H_nf_n.
	\end{align*}
	Hence $H_n^{\varepsilon_n}f_n^{\varepsilon_n}\to Hf$ is implied by $H_nf_n\to Hf$. That finishes the proof of~1): the operator~$H$ is dissipative and satisfies the range condition, and hence generates a semigroup~$V(t)$.
	\smallskip
	
	2) Since the operators $H_n^{\varepsilon_n}$ defined above are dissipative and satisfy the range condition (\cite[Lemma~5.7]{FengKurtz2006}), they generate a semigroup~$V_n^{\varepsilon_n}(t)$ acting on~$B(E)$. We showed above the convergence $H_n^{\varepsilon_n}\to H$. Hence by the Feng-Kurtz approximation (Theorem~\ref{thm:feng-kurtz-approx}) applied to $G_n=H_n^{\varepsilon_n}$, we obtain $V_n^{\varepsilon_n}(t)\to V(t)$: for any function $f\in C(E)$ and functions $f_n\in B(E)$ such that $\|f-f_n\|_E\to 0$,
	\[
	\|V(t)f-V_n^{\varepsilon_n}(t)f_n\|_E \xrightarrow{n\to\infty}0.
	\]
	Furthermore, the semigroup~$V_n^{\varepsilon_n}(t)$ approximates~$V_n(t)$, in the sense that for any function $f_n\in\mathcal{D}(H_n)$,
	\begin{equation}\label{eq:Vn-close-to-Vn_varep}
	\|V_n^{\varepsilon_n}(t)f_n-V_n(t)f_n\| \leq \sqrt{2\varepsilon_n t}\, e^{2n\|f_n\|}\|H_nf_n\|,
	\end{equation}
	which is proven in~\cite[Lemma~5.11]{FengKurtz2006}. The choice $\varepsilon_n=\exp\{-n^2\}$ implies that the difference vanishes in the limit $n\to\infty$. With that,
	\begin{align*}
		\|V(t)f-V_n(t)f_n\| &\leq \|V(t)f-V_n^{\varepsilon_n}(t)f_n\| + \|V_n^{\varepsilon_n}(t)f_n-V_n(t)f_n\| \to 0,
	\end{align*}
	which finishes the proof.
\end{proof}
\begin{running_example*}
	We illustrate on the small-diffusion process which condition of Theorem~\ref{thm:LDP-classical-sol} is difficult to verify in practice. We consider the small-diffusion process on the flat torus~$E=\mathbb{T}=\mathbb{R}/\mathbb{Z}$;
	that means the infinitesimal generator is the map $L_n:C^2(\mathbb{T})\to C(\mathbb{T})$ given by
	\[
	L_nf(x) = \frac{1}{2} \frac{1}{n} \Delta f(x).
	\]
	Therefore, the nonlinear generators~$H_n:C^2(\mathbb{T})\to C(\mathbb{T})$ read
	\[
	H_nf(x) = \frac{1}{n} e^{-nf(x)}L_n e^{nf(x)} = \frac{1}{2}\frac{1}{n}\Delta f(x) + \frac{1}{2}|\nabla f(x)|^2.
	\]
	They converge to~$Hf(x) := \frac{1}{2}|\nabla f(x)|^2$. Indeed, if we take for instance the domain~$\mathcal{D}(H) := C^2(\mathbb{T})$, then for any~$f\in\mathcal{D}(H)$, the constant sequence~$f_n:=f$ satisfies
	\[
	\|Hf-H_n f_n\|_\mathbb{T} =\frac{1}{2}\frac{1}{n} \|\Delta f\|_\mathbb{T} \to 0,
	\]
	with $\|\cdot\|_\mathbb{T}$ the supremum norm. We are only left with verifying the range condition for~$H$ in order to apply Theorem~\ref{thm:LDP-classical-sol}. The definition translates to the following PDE-problem: for a $C^2$ function $h:\mathbb{T}\to\mathbb{R}$ and for $\tau > 0$, find $u:\mathbb{T}\to\mathbb{R}$ in the domain of $H$ such that for any $x\in\mathbb{T}$, we have
	\[
	u(x) - \tau \frac{1}{2}|\nabla u(x)|^2 = h(x).
	\]
	There is no general theory available providing the existence of such a solution. The problem lies in the differentiability that solutions have to satisfy.
	\qed
\end{running_example*}
Using viscosity solutions makes the semigroup approach to large deviations applicable. In the words of Jin Feng and Thomas Kurtz~\cite[Preface]{FengKurtz2006}:
\smallskip

\emph{"This work began as a research paper intended to show how the convergence of nonlinear semigroups associated with a sequence of Markov processes implied the large deviation principle for the sequence. We expected the result to be of little utility for specific applications, since classical convergence results for nonlinear semigroups involve hypotheses that are very difficult to verify, at least using classical methods. We should have recognized at the beginning that the modern theory of viscosity solutions provides the tools needed to overcome the classical difficulties."}
\smallskip

We sketch in the next section how the approach using viscosity solutions works out in the compact setting.
\subsection{Using viscosity solutions}
\label{subsec:LDP-via-visc-sol}
In the previous section, we discussed how to verify the convergence of nonlinear semigroups $V_n$ from the convergence of associated formal nonlinear generators~$H_n$. The Feng-Kurtz approximation (Theorem~\ref{thm:feng-kurtz-approx}) was the key to obtain semigroup convergence from generator convergence. The example above illustrates that while finding a candidate limit~$H$ of the~$H_n$ is often straightforward, verifying the range condition for~$H$ is hard. We required the range condition to generate a semigroup by the Crandall-Liggett theorem. Here, we discuss why viscosity solutions are well suited for generating the desired limiting semigroup.
\smallskip

The basic idea is to use weak solutions~$u$ of~$(1-\tau H)u=h$ that are not required to be in the domain of~$H$. Then we define an auxiliary operator~$\widehat{H}$ by adding the weak solutions to the domain of~$H$ and the corresponding ranges to the image of~$H$.
If the requirement on a solution is weak enough, we can find enough solutions until the domain~$\mathcal{D}(\widehat{H})$ is dense in~$C(E)$, such that the operator~$\widehat{H}$ automatically satisfies the range condition. However, we also want~$\widehat{H}$ to be a dissipative operator in order to use the Crandall-Liggett Theorem. The limit operators~$H$ that we start from are dissipative, which follows from the convergence~$H_n\to H$ (see the proof of Theorem~\ref{thm:LDP-classical-sol}). Therefore, we are searching for weak solutions such that we keep dissipativity while enlargening~$H$.
\smallskip

We now motivate why viscosity solutions are suited for that purpose. A generator~$L$ of a Markov process satisfies the positive maximum principle; for a function~$f$ in the domain of~$L$, if $x$ is a local maximum of~$f$, then $Lf(x)\leq 0$. 
This propery carries over to their nonlinear generators~$Hf=e^{-f}Le^f$, where we obtain that if $(f_1-f_2)(x)=\sup(f_1-f_2)$, then $Hf_1(x)- Hf_2(x)\leq 0$.
In general, operators satisfying the positive maximum principle are dissipative.
\smallskip

Hence adding weak solutions such that the extended operator~$\widehat{H}$ still satisfies the positive maximum principle suffices for our purposes. Now given a "weak solution" $u$ to the equation $(1-\tau H)u=h$, consider the extended operator $\widehat{H} := H\cup \left(u,(u-h)/\tau\right)$; that is we added the weak solution and its corresponding image "$Hu$"$=(u-h)/\tau$ to the graph of~$H$. Let us see how the newly added elements affect the maximum principle. If $u$ is a weak solution and $(u-f)(x)=\sup(u-f)$, with~$f\in\mathcal{D}(H)$ such that~$\widehat{H}f=Hf$, then 
\begin{align*}
	(\widehat{H}u-\widehat{H}f)(x) &= \frac{1}{\tau}(u-h)(x) - Hf(x)\\ &= \frac{1}{\tau} \left[u(x)-\tau Hf(x)-h(x)\right] \stackrel{!}{\leq}0.
\end{align*}
When considering $(f-u)$, then $(f-u)(x)=\sup(f-u)$ should imply
\begin{align*}
	(\widehat{H}f-\widehat{H}u)(x) &= Hf(x) -\frac{1}{\tau}(u-h)(x)\\ &= -\frac{1}{\tau} \left[u(x)-\tau Hf(x)-h(x)\right] \stackrel{!}{\leq}0.
\end{align*}
This motivates the following definition.
\begin{definition}[Viscosity solutions]\label{def:viscosity-sol:beginner-guide}
	For a compact Polish space~$E$, let~$\tau > 0$ and~$h\in C(E)$. For an operator $H:\mathcal{D}(H)\subseteq C(E)\to C(E)$ with domain~$\mathcal{D}(H)$, consider the equation~$(1-\tau H)u=h$.
	\begin{enumerate}[label=(\roman*)]
		\item We say that a function~$u_1:E\to\mathbb{R}$ is a \emph{viscosity subsolution} if it is bounded, upper semicontinuous and for any function~$f\in\mathcal{D}(H)$, if a point~$x\in E$ is such that $(u_1-f)(x)=\sup_E(u_1-f)$, then
		\[
		u_1(x)-\tau Hf(x)-h(x)\leq 0.
		\]
		\item We call a function~$u_2:E\to\mathbb{R}$ a \emph{viscosity supersolution} if it is bounded, lower semicontinuous and for any function~$f\in\mathcal{D}(H)$, if a point~$x\in E$ is such that $(f-u_2)(x)=\sup_E(f-u_2)$, then
		\[
		u_2(x)-\tau Hf(x)-h(x)\geq 0.
		\]
		\item A function $u:E\to\mathbb{R}$ is a \emph{viscosity solution} if it is both a viscosity subsolution and a viscosity supersolution. \qed
	\end{enumerate}
\end{definition}
A viscosity solution in the sense of Definition~\ref{def:viscosity-sol:beginner-guide} is both upper- and lower semicontinuous, and hence continuous.
\begin{definition}[Comparison principle]
	We say that $(1-\tau H)u=h$ as in Definition~\ref{def:viscosity-sol:beginner-guide} satisfies the \emph{comparison principle} if for any viscosity subsolution~$u_1$ and viscosity supersolution~$u_2$, the inequality $u_1\leq u_2$ holds on~$E$.\qed
\end{definition}
If the comparison principle holds, then any two viscosity solutions~$u,v$ are equal: since $u$ is a viscosity subsolution and~$v$ a viscosity supersolution,~$u\leq v$. Reversing the roles, we find $v\leq u$. Hence~$u=v$, and thus the comparison principle implies uniqueness of viscosity solutions.
\smallskip

We now formulate the viscosity-analogue of Theorem~\ref{thm:LDP-classical-sol}. For that purpose, we define for an operator $H:\mathcal{D}(H)\subseteq C(E)\to C(E)$ its extension~$\widehat{H}$ as follows. If for any~$\tau > 0$ and $h\in C(E)$ there exists a unique viscosity solution~$u$ of the equation~$(1-\tau H)u=h$, then we denote it by $R(\tau)h:=u$. The map~$R(\tau)$ is called the resolvent. We denote by~$\widehat{H}\subseteq C(E)\times C(E)$ the operator defined as the graph
\[
\widehat{H} := \bigcup_{\tau > 0}\left\{\left(R(\tau)h,\frac{1}{\tau}(R(\tau)h-h)\right)\,:\,h \in C(E)\right\}.
\]
\begin{theorem}[Theorem~6.14 in~\cite{FengKurtz2006}, Large Deviations via Viscosity Solutions]\label{thm:LDP-visc-sol}
	Let~$E$ be a compact Polish space and~$\{X^n\}$ be a sequence of Markov processes in~$\mathcal{X}=C_E[0,T]$, with generators~$L_n$ and associated nonlinear generators~$H_n$ from Definition~\ref{def:nonlinear-generators}. Assume the following:
	\begin{enumerate}[label=(\roman*)]
		\item There is a densely defined operator $H:\mathcal{D}(H)\subseteq C(E) \to C(E)$ such that $H_n$ converges to~$H$; for every~$f\in\mathcal{D}(H)$, there are functions~$f_n\in\mathcal{D}(H_n)$ such that 
		\[
		\|f-f_n\|_E\xrightarrow{n\to\infty}0 \quad \text{and}\quad \|Hf-H_nf_n\|_E \xrightarrow{n\to\infty}0.
		\]
		\item For~$\tau>0$ and~$h\in C(E)$, the comparison principle holds for~$(1-\tau H)u=h$.
	\end{enumerate}
	Suppose furthermore that the initial conditions~$X^n(0)$ satisfy a large deviation principle with rate function~$\mathcal{I}_0$.
	
	Then the sequence~$\{X^n\}$ satsifies a large deviation principle in~$\mathcal{X}$ with a rate function~$\mathcal{I}$ given by~\eqref{eq:RF-with-semigroup-approach}, where the semigroup~$V(t)$ is generated by the operator~$\widehat{H}$: for every~$f\in C(E)$, we have~$V(t)f = \lim_{k\to\infty}[R(t/k)]^k f$.
\end{theorem}
\begin{proof}[Sketch of proof of Theorem~\ref{thm:LDP-visc-sol}]
	Just as in the previous section, we want to verify exponential tightness and convergence of the nonlinear semigroups~$V_n(t)$ to some limit~$V(t)$. Then the large-deviation statement follows from Theorem~\ref{thm:LDP-via-semigroup-convergence:compact}. Exponential tightness follows from the convergence condition~$H_n\to H$, just as we indicated in the proof of Theorem~\ref{thm:LDP-classical-sol} in the previous section.
	\smallskip
	
	We are left with showing that the operator~$\widehat{H}$ satisfies the conditions of the Crandall-Liggett Theorem (dissipativity and the range condition) with dense domain, so that it generates a semigroup~$V(t)$ acting on~$C(E)$, and that we have convergence~$V_n(t)\to V(t)$ as in Theorem~\ref{thm:LDP-via-semigroup-convergence:compact}. The argument is based on the same techinque as in Theorem~\ref{thm:LDP-classical-sol}: we use the Hille-Yosida approximations~$L_n^{\varepsilon_n}$ of~$L_n$ and their corresponding nonlinear generators~$H_n^{\varepsilon_n}$. 
	\smallskip
	
	The operator~$\widehat{H}$ is defined via the existence of unique viscosity solutions. We first show that its domain is dense in~$C(E)$. To that end, fix~$\tau>0$ and a function~$h\in C(E)$. Since~$\text{range}(1-\tau H_n^{\varepsilon_n})=C(E)$ (\cite[Lemma~5.7]{FengKurtz2006}), there exists a classical solution~$f_n\in C(E)$ to~$(1-\tau H_n^{\varepsilon_n})f_n = h$. In particular,~$f_n$ is a viscosity solution. Define the functions~$u_1,u_2:E\to\mathbb{R}$ by
	\begin{align*}
		u_1(x)&:=\lim_{k\to\infty}\sup\left\{f_n(z)\,|\,n\geq k,\, d(x,z)\leq \frac{1}{k}\right\}\\
		u_2(x)&:=\lim_{k\to\infty}\inf\left\{f_n(z)\,|\,n\geq k,\, d(x,z)\leq \frac{1}{k}\right\} 
	\end{align*}
	It is shown in~\cite[Lemma~6.9]{FengKurtz2006} that~$u_1$ is a viscosity subsolution and~$u_2$ is a viscosity supersolution of~$(1-\tau H)u=h$. By construction, $u_1\geq u_2$. By assumption, the comparison principle holds, which gives~$u_1\leq u_2$. Hence the function $u:=u_1=u_2$ is the unique viscosity solution to~$(1-\tau H)u=h$. Define the resolvent map $R(\tau):C(E)\to C(E)$ by setting $R(\tau)h:=u$. Lemma~6.9 also establishes the estimate
	\[
	\|u-f\|\leq \|h-(f-\tau Hf)\|
	\]
	for any~$f\in\mathcal{D}(H)$. Specializing to~$h$ in the domain of~$\mathcal{D}(H)$ and choosing in the estimate $f=h$, this implies
	\[
	\|R(\tau)h-h\| \leq \tau \|Hh\| \xrightarrow{\tau \to 0} 0.
	\]
	That demonstrates $\mathcal{D}(H)\subseteq \overline{\mathcal{D}(\widehat{H})}$, and we conclude that~$\mathcal{D}(\widehat{H})$ is dense in~$C(E)$ since~$\mathcal{D}(H)$ is dense in~$C(E)$.
	\smallskip
	
	We showed in the proof of Theorem~\ref{thm:LDP-classical-sol} that the assumed convergence condition~$H_n\to H$ implies~$H_n^{\varepsilon_n}\to H$, and that the operator~$H$ is dissipative as the limit of the dissipative operators~$H_n^{\varepsilon_n}$. The fact that dissipativity transfers further to~$\widehat{H}$ is proven in~\cite[Theorem~6.13]{FengKurtz2006}, part~(c); the operators~$H_n$ appearing therein are the Hille-Yosida approximants~$H_n^{\varepsilon_n}$. The range condition on~$\widehat{H}$ is satisfied by construction.
	Now the fact that~$V(t)=\lim_k R(t/k)^k$ follows from the Feng-Kurtz approximation theorem.
\end{proof}
In summary, with introducing viscosity solutions, we weakened the requirement on a function being a solution. The existence of sub- and supersolutions is guaranteed as a consequence of the convergence $H_n\to H$. However, dissipativity of~$H$ is no longer sufficient for uniqueness of viscosity solutions. That is because viscosity solutions are in general not in the domain of~$H$, which breaks the argument shown below Definition~\ref{def:range-condition}. This contrasts the classical approach, where uniqueness of solutions is for free while existence of solutions remains open.
We close this section by illustrating Theorem~\ref{thm:LDP-visc-sol}.
\begin{running_example*}
	Consider the small-diffusion process $\dd X^n_t=n^{-1/2}\dd B_t$ on the torus~$\mathbb{T}$. The linear generators are $L_nf = (2n)^{-1}\Delta f$, and
	\[
	H_nf = \frac{1}{n}e^{-nf}L_n e^{nf} = \frac{1}{2}\frac{1}{n}\Delta f + \frac{1}{2}|\nabla f|^2.
	\]
	We already checked in the previous section that they converge to $Hf = \frac{1}{2}|\nabla f|^2$. Here, we can take for instance the domain~$\mathcal{D}(H)=C^{17}(\mathbb{T})$.
	\smallskip
	
	For applying Theorem~\ref{thm:LDP-visc-sol}, we must verify the comparison principle. To that end, fix~$\tau >0$ and~$h\in C(\mathbb{T})$, and let~$u_1$ be a vioscosity subsolution and~$u_2$ be a viscosity supersolution of $(1-\tau H)u=h$. We want to prove that $u_1\leq u_2$. 
	
	For illustration, suppose first that they are classical sub- and supersolutions; they are in the domain of~$\mathcal{D}(H)$ and for any~$x\in \mathbb{T}$,
	\begin{align*}
		u_1(x)-\tau Hu_1(x) - h(x)\leq 0 \quad \text{and}\quad u_2(x)-\tau Hu_2(x)-h(x)\geq 0.
	\end{align*}
	Let $x_m$ be a point such that $(u_1-u_2)(x_m) = \sup_\mathbb{T}(u_1-u_2)$. Then we have $\nabla u_1(x_m)=\nabla u_2(x_m)$, and by the sub-and supersolution inequalities, we obtain 
	\begin{align*}
		(u_1-u_2)(x) \leq (u_1-u_2)(x_m) &\leq Hu_1(x_m) - Hu_2(x_m) 
		\\&= \frac{1}{2}\left(|\nabla u_1(x_m)|^2-|\nabla u_2(x_m)|^2\right) = 0.
	\end{align*}
	That shows uniqueness of classical solutions to~$(1-\tau H)u=h$.
	\smallskip
	
	For viscosity sub- and supersolutions, we can not rely on~$u_1,u_2$ being in the domain of~$H$. The classical trick is to use distance-like functions that are in the domain in order to approximate $\sup_\mathbb{T}(u_1-u_2)$. Here, we can take
	\[
	\Psi(x,y) := \sin^2(\pi(x-y)).
	\]
	Then $\Psi(\cdot,y)$ and $\Psi(x,\cdot)$ are smooth on~$\mathbb{T}$, and hence in the domain~$\mathcal{D}(H)=C^{17}(\mathbb{T})$. For $\alpha>0$, define
	\[
	\Phi_\alpha(x,y) := u_1(x)-u_2(y) - \alpha \Psi(x,y).
	\]
	By the semi-continuity properties of~$u_1,u_2$, for every~$\alpha>0$ there are~$x_\alpha,y_\alpha$ such that
	\[
	\Phi_\alpha(x_\alpha,y_\alpha) = \sup_{x,y}\Phi_\alpha(x,y).
	\]
	The point is that since~$u_1,u_2$ are bounded, these maximizing points~$x_\alpha,y_\alpha$ will converge to each other as~$\alpha\to\infty$; indeed, observing that $\Phi_\alpha(x_\alpha,x_\alpha)\leq \Phi(x_\alpha,y_\alpha)$, we obtain
	\[
	\Psi(x_\alpha,y_\alpha) \leq \frac{2}{\alpha}(u_2(x_\alpha)-u_2(y_\alpha)) \leq \frac{4}{\alpha}\|u_2\| \xrightarrow{\alpha\to\infty} 0.
	\]
	Since $\Psi(x,y)\geq 0$ and~$\Psi(x,y)=0$ if and only if $x=y$, we approximate the supremum of $(u_1-u_2)$ in the sense that $x_\alpha \approx y_\alpha$ and
	\begin{align*}
		\sup_{\mathbb{T}}(u_1-u_2) &= \sup_{x\in\mathbb{T}}\left(u_1(x)-u_2(x)-\alpha \Psi(x,x)\right)\\
		&\leq \sup_{x,y}\Phi_\alpha(x,y) = \Phi_\alpha(x_\alpha,y_\alpha)\\
		&\leq u_1(x_\alpha)-u_2(y_\alpha).
	\end{align*}
	Now we can use the sub- and supersolution inequalities. The test functions defined by
	\begin{align*}
		f_1^\alpha(x) := u_2(y_\alpha) + \alpha\Psi(x,y_\alpha)\quad\text{and}\quad f_2^\alpha(y):=u_1(x_\alpha) - \alpha\Psi(x_\alpha,y)
	\end{align*}
	are smooth, and hence are both in the domain of~$H$. By construction,
	\[
	(u_1-f_1^\alpha)(x_\alpha) = \sup_\mathbb{T}(u_1-f_1^\alpha)\quad\text{and}\quad (f_2^\alpha-u_2)(y_\alpha) = \sup_{\mathbb{T}}(f_2^\alpha-u_2),
	\]
	so that with the sub- and supersolution inequalities, 
	\begin{align*}
		u_1(x_\alpha)-\tau Hf_1^\alpha(x_\alpha)-h(x_\alpha)\leq 0\quad\text{and}\quad 
		u_2(y_\alpha)-\tau Hf_2^\alpha(y_\alpha)-h(y_\alpha)\geq 0.
	\end{align*}
	With these, we can further estimate $u_1(x_\alpha)-u_2(y_\alpha)$ to arrive at
	\[
	\sup_{\mathbb{T}}(u_1-u_2) \leq \tau \left[Hf_1^\alpha(x_\alpha)-Hf_2^\alpha(y_\alpha)\right] + h(x_\alpha)-h(y_\alpha).
	\]
	Since $Hf=(1/2) |\nabla f|^2$ depends only on the gradient and $\nabla f_1^\alpha(x_\alpha)=\nabla f_2^\alpha(y_\alpha)$, the difference of the Hamiltonians is zero. The function~$h$ is uniformly continuous on the compact space~$\mathbb{T}$. Then since $\Psi(x_\alpha,y_\alpha)\to 0$, we obtain finally
	\begin{align*}
		\sup_\mathbb{T}(u_1-u_2) \leq \liminf_{\alpha\to 0}|h(x_\alpha)-h(y_\alpha)| = 0,
	\end{align*}
	which finishes the verification of the comparison principle.
	\qed
\end{running_example*}
The running example also illustrates a principle that applies more generally. We can choose the domain of the limiting operator~$H$ as small as we want, provided that it contains sufficient functions to allow for verifying the comparison principle. In the example, merely using smooth functions was sufficient. As a rule of thumb, in compact spaces one wants to make sure that distance functions are in the domain of~$H$.
\section{Action-integral representation of rate functions}
\label{BG:sec:semigroup-flow-HJ-eq}
Let us summarize where we stand after the previous section. We considered a sequence of Markov processes $X^n$ in~$\mathcal{X}=C_E[0,T]$ and established that the following two conditions imply a pathwise large deviation principle:
\begin{enumerate}[label=(\roman*)]
	\item The nonlinear generators converge $H_n\to H$.
	\item The comparison principle holds for $(1-\tau H)u=h$. 
\end{enumerate}
We illustrated on the example of small diffusion how one can verify these conditions in practice. The rate function is given via the limiting semigroup~$V(t)$ generated by~$H$, based on finding unique viscosity solutions of~$(1-\tau H)u=h$.
\smallskip

In this section, we focus on this rate function, which is given by
\begin{equation}\label{eq:RF-with-semigroup-approach:last-sec}
\mathcal{I}(x) = \mathcal{I}_0(x(0)) + \sup_{k\in\mathbb{N}}\sup_{(t_1,\dots,t_k)} \sum_{i=1}^k \mathcal{I}_{t_i-t_{i-1}}(x(t_i)|x(t_{i-1})),
\end{equation}
where the conditional rate functions~$\mathcal{I}_t(z|y)$ are
\begin{equation}\label{BG:eq:RF-1d-marginals:last-sec}
\mathcal{I}_t(z|y) = \sup_{f\in C(E)}\left[f(z)-V(t)f(y)\right].
\end{equation}
We specialize henceforth to the state space~$E=\mathbb{T}$, the one-dimensional flat torus. Denote by~$\cA\cC_E[0,T]$ the set of absolutely continuous curves in~$E$. Our aim in this section is to find conditions under which the rate function~$\mathcal{I}:\mathcal{X}\to[0,\infty]$ is given by a Lagrangian~$\mathcal{L}:\mathbb{R}\to[0,\infty]$ via the formula
\begin{equation}
\label{BG:eq:action-integral-formula}
\mathcal{I}(x) = 
\begin{cases}
\mathcal{I}_0(x(0))+\int_0^T \mathcal{L}(\partial_t x(t))\,\dd t, &\qquad x\in\cA\cC_E[0,T],\\
+\infty, &\qquad \text{otherwise}.
\end{cases}
\end{equation}
We first indicate how to obtain~\eqref{BG:eq:action-integral-formula} from~\eqref{eq:RF-with-semigroup-approach:last-sec} via an informal calculation. Then we show how this can be obtained rigorously based on identifying the semigroup~$V(t)$ at least formally as a Hamilton-Jacobi semigroup~$V_\mathcal{H}(t)$---we give details below by Proposition~\ref{BG:prop:good-Lagrangians-give-action}. Finally, we show in what sense the required equality~$V(t)=V_\mathcal{H}(t)$ follows from solving a Hamilton-Jacobi equation.
\subsubsection{Action-integral via an informal calcuation.}
	Here we consider the 
	operator $Hf = \mathcal{H}(\nabla f)$ with the Hamiltonian $\mathcal{H}(p)=\frac{1}{2}p^2$. Let~$x\in\cA\cC_E[0,T]$. We want to show~\eqref{BG:eq:action-integral-formula} starting from~\eqref{eq:RF-with-semigroup-approach:last-sec}. To that end, we  compute~$\mathcal{I}_{t_2-t_1}\left(z|y\right)$ 
	for $y$ close to $z$ and $t=t_2-t_1>0$ small, having in mind that $y=x(t_1)$ and $z=x(t_2)$ are close to each other.
	With the formal expansions
	\[
	e^{tH}\approx 1+tH\quad\text{and}\quad f(z)-f(y) \approx (z-y)\cdot\nabla f(y),
	\]
	and thinking of~$V(t)=e^{tH}$ (the generator of~$V(t)$ is~$H$), we obtain from~\eqref{BG:eq:RF-1d-marginals:last-sec}
	\begin{align*}
		\mathcal{I}_t(z|y) &\approx \sup_f \left[f(z)-f(y) - tHf(y)\right]\\
		&\approx t \cdot \sup_f \left[\nabla f(y) \cdot \frac{z-y}{t} - \mathcal{H}(\nabla f(y))\right]
		= t\cdot \sup_p\left[p \cdot \frac{z-y}{t} - \mathcal{H}(p)\right].
	\end{align*}
	Hence with $\mathcal{L}(v)=\sup_p\left[pv-\mathcal{H}(p)\right]$, which here is equal to $v^2/2$, we find 
	\begin{align*}
		\mathcal{I}_{t_2-t_1}(z|y) &\approx (t_2-t_1)\mathcal{L}\left(\frac{z-y}{t_2-t_1}\right) = \int_{t_1}^{t_2} \mathcal{L}(\partial_s \gamma_{t_1t_2}(s))\,\dd s,	
	\end{align*}
	where $\gamma_{t_1t_2}:[t_1,t_2]\to\mathbb{R}$ is the linear path connecting $y$ and $z$. Now starting from~\eqref{eq:RF-with-semigroup-approach:last-sec}, the action-integral formula follows from summing up all the conditional rate functions, since the linear paths~$\gamma_{t_kt_{k+1}}$ approximate~$x$ in~$[t_k,t_{k+1}]$.
\subsubsection{Action-integral via rewriting the semigroup.}
Here we indicate how to make the above informal calculation rigorous. We start from an operator~$H$ acting on functions as~$Hf(x)=\mathcal{H}(\nabla f(x))$, with a convex Hamiltonian~$\mathcal{H}:\mathbb{R}\to\mathbb{R}$ satisfying~$\mathcal{H}(0)=0$. Define the Lagrangian~$\mathcal{L}$ as the Legendre dual~$\mathcal{L}(v):=\sup_{p\in\mathbb{R}}\left[pv-\mathcal{H}(p)\right]$, and the semigroup~$V_\mathcal{H}(t)$ by
\begin{equation}\label{BG:eq:control-semigroup}
V_\mathcal{H}(t)f(x) := \sup_{\substack{\gamma\\\gamma(0)=x}} \left[f(\gamma(t))-\int_0^t\mathcal{L}(\partial_s\gamma(s))\,\dd s\right],
\end{equation}
where the supremum is taken over absolutely continuous paths $\gamma:[0,t]\to E$.
Formally taking the time derivative, exchanging limit and supremum, we obtain
\begin{align*}
	\frac{\dd}{\dd t}\bigg|_{t=0} V_\mathcal{H}(t)f(x)&= \sup_{\gamma(0)=x} \left[\nabla f(\gamma(0)) \cdot\partial_t\gamma(0) - \mathcal{L}(\partial_t\gamma(0))\right]\\
	&= \sup_{v\in\mathbb{R}}\left[\nabla f(x)\cdot v - \mathcal{L}(v)\right] = \mathcal{H}(\nabla f(x)).
\end{align*} 
That is why we indeed expect the the operator~$H$ to be the generator of~$V_\mathcal{H}(t)$.
\begin{proposition}
	\label{BG:prop:good-Lagrangians-give-action}
	Suppose that~$V(t)=V_\mathcal{H}(t)$. Then the rate function~\eqref{eq:RF-with-semigroup-approach:last-sec} satisfies the action-integral form~\eqref{BG:eq:action-integral-formula}.
\end{proposition}
\begin{proof}[Sketch of proof of Proposition~\ref{BG:prop:good-Lagrangians-give-action}]
	We first show that the Lagrangian is superlinear, that means~$(\mathcal{L}(v)/|v|) \to \infty$ as~$|v|\to\infty$. Then for any~$t,M\geq0$, the sub-level sets $\{\gamma\in\mathcal{X}\,|\,\int_0^t\mathcal{L}\left(\partial_s\gamma(s)\right)\,\dd s\leq M\}$ are compact in~$\mathcal{X}$---we do not prove this compactness statement here, but comment on it in Section~\ref{BG:sec:bibliographical-notes}.
	Regarding superlinearity, for any~$c>0$ we have
	\begin{align*}
	\frac{\mathcal{L}(v)}{|v|} &=\sup_{p\in\mathbb{R}}\left[p\cdot \frac{v}{|v|}-\frac{\mathcal{H}(p)}{|v|}\right]\\
	&\geq \sup_{|p|=c}\left[p\cdot \frac{v}{|v|}-\frac{\mathcal{H}(p)}{|v|}\right]
	\geq c - \frac{1}{|v|}\sup_{|p|=c}\mathcal{H}(p).
	\end{align*}
	The convex Hamiltonian is continuous, and therefore~$\sup_{|p|=c}\mathcal{H}(p)$ is finite. Hence for arbitrary~$c>0$, we have~$\mathcal{L}(v)/|v| > c/2$ for all~$|v|$ large enough.
	\smallskip
	
	Let $x:[0,T]\to E$ be absolutely continuous and take two arbitrary~$t_1<t_2$. We show for $y:=x(t_1)$ and $z:=x(t_{2})$ that
	\begin{equation}\label{eq:action-int:two-point-RF-is-L}
	\mathcal{I}_{t_2-t_1}\left(z|y\right) =  \inf_{\substack{\gamma(t_1)=y\\\gamma(t_2)=z}}\int_{t_1}^{t_2} \mathcal{L}\left(\partial_s\gamma\right)\,\dd s,
	\end{equation}
	where the infimum is taken over absolutely continuous paths $\gamma:[t_1,t_2]\to E$. Once we have this equality established, we obtain for arbitrary~$k\in\mathbb{N}$ and points in time $t_1,\dots,t_k=T$ the esimate
	\begin{align*}
		\mathcal{I}_{t_1}(x_1|x_0) + \mathcal{I}_{t_2-t_1}(x_2|x_1) + \dots + \mathcal{I}_{t_{k}-t_{k-1}}(x_k|x_{k-1})
		\leq \int_0^{T} \mathcal{L}(\partial_s x)\,\dd s,
	\end{align*}
	since~$x$ satisfies the begin- and endpoint contraints. 
	For the reverse inequality, we note that adding time points increases the two-point rate functions since we add a condition on the paths; for $t_1<t_2<t_3$,
	\begin{align*}
	\mathcal{I}_{t_3-t_1}(x_3|x_1) &= \inf_{\substack{\gamma(t_1)=x_1\\\gamma(t_3)=x_3}}\left[\int_{t_1}^{t_2} \mathcal{L}(\partial_t\gamma)\,\dd t + \int_{t_2}^{t_3}\mathcal{L}(\partial_t\gamma)\,\dd t\right] \\
	&\leq 
	\inf_{\substack{\gamma(t_1)=x_1\\\gamma(t_2)=x_2}}\left[\int_{t_1}^{t_2} \mathcal{L}(\partial_t\gamma)\,\dd t\right]  +
	\inf_{\substack{\gamma(t_2)=x_2\\\gamma(t_3)=x_3}}\left[ \int_{t_2}^{t_3}\mathcal{L}(\partial_t\gamma)\,\dd t\right] \\
	&= \mathcal{I}_{t_3-t_2}(x_3|x_2) + \mathcal{I}_{t_2-t_1}(x_2|x_1).
	\end{align*}
	The partitions of a time interval~$[0,T]$ give rise to a monotonically increasing sequence. In the limit, we obtain
	\begin{align*}
	\sup_k \sup_{t_i}\sum_{i=1}^k \mathcal{I}_{t_i-t_{i-1}}(x(t_i)|x(t_{i-1}))
	= \int_0^T \mathcal{L}\left(\partial_s x(s)\right)\,\dd s.
	\end{align*}
	We do not show that here, but refer to~\cite[Definition~7.11,~Example~7.12]{Villani2008}. We now show how~\eqref{eq:action-int:two-point-RF-is-L} follows from the compact sub-level sets. Starting from the assumption~$V(t)=V_\mathcal{H}(t)$, we have
	\begin{align*}
		\mathcal{I}_t\left(z|y\right) &\myeqdef \sup_f\left[f(z)-V(t)f(y)\right]
		= \sup_f\left[f(z)-V_\mathcal{H}(t)f(y)\right]\\
		&= \sup_{f}\inf_{\substack{\gamma(0)=y}}\left[f(z)-f(\gamma(t))+\int_0^t\mathcal{L}\left(\partial_s \gamma\right)\,\dd s\right].
	\end{align*}
	For any $f\in C(E)$, 
	\[
	\inf_{\gamma(0)=y}\left[f(z)-f(\gamma(t))+\int_0^t\mathcal{L}(\partial_s\gamma)\,\dd s\right] \leq \inf_{\substack{\gamma(0)=y\\\gamma(t)=z}} \int_0^t\mathcal{L}(\partial_s\gamma)\,\dd s,
	\]
	since $\{\gamma:\gamma(0)=y\}$ contains $\{\gamma:\gamma(0)=y,\gamma(t)=z\}$. Taking the supremum over all~$f$ shows the inequality "$\leq$".
	\smallskip
	
	For the reverse, let~$f\in C(E)$. There are curves $\gamma_m$ satisfying~$\gamma_m(0)=y$ and
	\begin{align*}
		\inf_{\gamma(0)=y}\left[f(z)-f(\gamma(t))+\int_0^t\mathcal{L}(\partial_s\gamma)\,\dd s\right] + \frac{1}{m} \geq f(z)-f(\gamma_m(t))+\int_0^t\mathcal{L}(\partial_s\gamma_m)\,\dd s.
	\end{align*}   
	Since $f$ is bounded, this implies~$\limsup_{m\to\infty}\int_0^t\mathcal{L}(\partial_s\gamma_m)\,\dd s < \infty$.
	By compactness of sublevel sets, we can pass to a converging subsequence (denoted as well by~$\gamma_m$).
	If $\gamma_m(t)\not\to z$, then~$\mathcal{I}_t(z|y)=\infty$,
	and the desired estimate holds. If $\gamma_m(t)\to z$, then by lower semicontinuity of $\gamma\mapsto\int_0^t\mathcal{L}(\partial_s\gamma)\,\dd s$,
	\begin{align*}
		\inf_{\gamma(0)=y}\left[f(z)-f(\gamma(t))+\int_0^t\mathcal{L}(\partial_s\gamma)\,\dd s\right]&\geq \liminf_{m\to\infty}f(z)-f(\gamma_m(t))+\int_0^t\mathcal{L}(\partial_s\gamma_m)\,\dd s\\
		&\geq \int_0^t\mathcal{L}(\partial_s\gamma)\,\dd s \geq \inf_{\substack{\gamma(0)=y\\\gamma(t)=z}}\int_0^t\mathcal{L}(\partial_s\gamma)\,\dd s,
	\end{align*}
	and the reverse inequality follows. 
\end{proof}
\subsubsection{Rewriting the semigroup via solving the Hamilton-Jacobi equation.}
We saw above that if~$V(t)=V_\mathcal{H}(t)$, then the action-integral form of the rate function follows. In this section we illustrate how to verify this equality. The semigroup~$V(t)$ is defined by the resolvent map $h\mapsto R(\tau)h:=u$, where~$u$ is the unique viscosity solution of~$(1-\tau H)u=h$; for~$f\in C(E)$, we have uniformly
\[
V(t)f = \lim_{k\to\infty}[R(t/k)]^k f.
\]
For a Hamiltonian~$\mathcal{H}$ with corresponding Lagrangian~$\mathcal{L}$, define~$R_\mathcal{H}(\tau)$ by
\begin{equation}\label{BG:HJ-resolvent}
R_\mathcal{H}(\tau)h(x) := \sup_{\substack{\gamma\in\cA\cC_E[0,\infty)\\\gamma(0)=x}}\int_0^\infty \frac{1}{\tau}e^{-s/\tau} \left[h(\gamma(s)) - \tau \mathcal{L}\left(\partial_s \gamma(s)\right)\right]\,\dd s.
\end{equation}
One can show, under suitable conditions on the Lagrangian, that also
\begin{equation}\label{BG:eq:H-semigroup-from-resolvent}
V_\mathcal{H}(t)f = \lim_{k\to\infty}[R_\mathcal{H}(t/k)]^k f.
\end{equation}
Therefore the desired equality~$V(t)=V_\mathcal{H}(t)$ follows if we prove~$R(\tau)=R_\mathcal{H}(\tau)$ for all~$\tau>0$. Let us first focus on establishing~$R(\tau)=R_\mathcal{H}(\tau)$, and defer the problem of obtaining~\eqref{BG:eq:H-semigroup-from-resolvent}. We will show that~$R_\mathcal{H}(\tau)$ gives viscosity solutions to~$(1-\tau H)u=h$. Then~$R_\mathcal{H}(\tau)=R(\tau)$ follows by definition of~$R(\tau)$. The following definition summarizes the key properties to look after.
\begin{definition}\label{def:visc-sol-machine}
	For~$\tau>0$, let~$\mathrm{R}(\tau)$ be a map~$\mathrm{R}(\tau):C(E)\to C(E)$. We call the family~$\{\mathrm{R}(\tau)\}_{\tau >0}$ a \emph{contractive pseudo-resolvent} if:
	\begin{enumerate}[label=(\roman*)]
		\item \label{item:visc-machine:resolv-id} For any $0<\tau_1<\tau_2$, we have
		\[
		\mathrm{R}(\tau_2) = \mathrm{R}(\tau_1)\left[\mathrm{R}(\tau_2)-\frac{\tau_1}{\tau_2}\left(\mathrm{R}(\tau_2)-\mathbf{1}\right)\right].
		\]
		\item \label{item:visc-machine:contr} The map $\mathrm{R}(\tau)$ is contractive: for any two functions $h_1,h_2\in C(E)$, we have the estimate $\sup_E\left(\mathrm{R}(\tau)h_1-\mathrm{R}(\tau)h_2\right)\leq \sup_E \left(h_1-h_2\right)$.
	\end{enumerate}
\end{definition}
\begin{theorem}[]\label{thm:viscosity-sol-machine}
	Let~$\mathrm{R}(\tau)$ be a contractive pseudo-resolvent and suppose that for any~$f\in\mathcal{D}(H)$, we have~$f = \mathrm{R}(\tau) (\mathbf{1}-\tau H)f$ on~$E$.
	Then for any~$\tau>0$ and~$h\in C(E)$, the function~$\mathrm{R}(\tau)h$ is a viscosity solution of $(1-\tau H)u=h$.
\end{theorem}
For the proof Theorem~\ref{thm:viscosity-sol-machine}, we will use the following simplification of~\cite[Lemma~7.8]{FengKurtz2006} (the proof in there is incorrect---see~\cite[Lemma~3.5]{Kraaij2019GenConv}). 
\begin{lemma}\label{lemma:visc-sol-machine}
	Let~$f,g:E\to\mathbb{R}$ be two continuous functions on a compact Polish space~$E$. Suppose that for any~$\varepsilon>0$, the inequality $\sup_E f \leq \sup_E (f-\varepsilon g)$ holds true. Then there is a point~$x\in E$ such that both $f(x)=\sup_E f$ and $g(x)\leq 0$. Similarly, if $\inf_E f\geq \inf_E (f-\varepsilon g)$, then $f(x)=\inf_E f$ and $g(x)\geq 0$ for some point~$x\in E$. 
\end{lemma}
\begin{proof}[Proof of Theorem~\ref{thm:viscosity-sol-machine}]
	Fix~$\tau>0$ and~$h\in C(E)$. Let~$f\in\mathcal{D}(H)$. For every~$\varepsilon>0$, we show below the estimate
	\begin{equation}\label{BG:eq:estimate-to-get-subsolution}
	\sup_E \left[\mathrm{R}(\tau)h-f\right] \leq \sup_E\left\{\mathrm{R}(\tau)h-f -\varepsilon\left[ \frac{1}{\tau}\left(\mathrm{R}(\tau)h-h\right)- Hf\right]\right\}.
	\end{equation}
	Then by Lemma~\ref{lemma:visc-sol-machine}, there is a point~$x\in E$ such that both
	\[
	(\mathrm{R}(\tau)h-f)(x)=\sup_E \left(\mathrm{R}(\tau)h-f\right) \quad\text{and}\quad \left[ \frac{1}{\tau}\left(\mathrm{R}(\tau)h-h\right)- Hf\right](x)\leq 0,
	\]
	which establishes that~$\mathrm{R}(\tau)h$ is a viscosity subsolution. The argument for the supersolution case is similar. We now prove the estimate~\eqref{BG:eq:estimate-to-get-subsolution}. We use the resolvent identity~\ref{item:visc-machine:resolv-id} to rewrite $\mathrm{R}(\tau)h$, and the fact that~$f=\mathrm{R}(\varepsilon) (\mathbf{1}-\varepsilon H)f$, to find
	\begin{equation*}
	\mathrm{R}(\tau)h-f = \mathrm{R}(\varepsilon)\left[\mathrm{R}(\tau)h-\frac{\varepsilon}{\tau}\left(\mathrm{R}(\tau)h-h\right)\right]-\mathrm{R}(\varepsilon)\left(f-\varepsilon Hf\right).
	\end{equation*}
	Since~$\mathrm{R}(\tau)$ is contractive~\ref{item:visc-machine:contr},
	\begin{align*}
		\sup_E\left[\mathrm{R}(\tau)h-f\right]
		\leq \sup_E \left\{\mathrm{R}(\tau)h-\frac{\varepsilon}{\tau}\left(\mathrm{R}(\tau)h-h\right)-\left(f-\varepsilon Hf\right)\right\},
	\end{align*}
	which establishes the desired estimate.
\end{proof}
To summarize where we are, Theorem~\ref{thm:viscosity-sol-machine} tells us that if~$R_\mathcal{H}(\tau)$ is a contractive pseudo-resolvent satisfying~$f = R_\mathcal{H}(\tau) (\mathbf{1}-\tau H)f$ for all~$f\in\mathcal{D}(H)$, then we have proven our desired equality~$R(\tau)=R_\mathcal{H}(\tau)$. We finish this section by showing that reasonable Hamiltonians give indeed rise to pseudo-resolvents. 
In the following theorem, we consider an operator~$H:\mathcal{D}(H)\subseteq C^1(E)\to C(E)$ acting functions by~$Hf(x)=\mathcal{H}(\nabla f(x))$, with a dense domain~$\mathcal{D}(H)\subseteq C(E)$. We associate the Lagrangian~$\mathcal{L}(v)=\sup_p[pv-\mathcal{H}(p)]$.
\begin{theorem}\label{BG:thm:resonable-H-gives-pseudo-resolvent}
	Suppose~$\mathcal{H}:\mathbb{R}\to\mathbb{R}$ is convex, continuously differentiable, and that~$\mathcal{H}(0)=0$. Then~$R_\mathcal{H}(\tau)$ defined by~\eqref{BG:HJ-resolvent} is a contractive pseudo-resolvent such that for all functions~$f\in\mathcal{D}(H)$,~$f = R_\mathcal{H}(\tau) (\mathbf{1}-\tau H)f$, and~\eqref{BG:eq:H-semigroup-from-resolvent} holds true.
\end{theorem}
\begin{proof}[Sketch of proof of Theorem~\ref{BG:thm:resonable-H-gives-pseudo-resolvent}]
	We first verify~$f=R_\mathcal{H}(1-\tau H)f$. For~$f\in\mathcal{D}(H)$,
	\begin{align*}
	R_\mathcal{H}(\tau)(f-\tau Hf)(x) \overset{\mathrm{def}}{=} \sup_{\gamma(0)=x}\int_0^\infty e^{-s/\tau}\left[\frac{1}{\tau}f(\gamma(s))- Hf(\gamma(s)) - \mathcal{L}(\partial_s\gamma(s)) \right]\dd s.
	\end{align*}
	Since~$Hf=\mathcal{H}(\nabla f)$ and~$\mathcal{L}(v)+\mathcal{H}(p)\geq pv$ for any~$p,v\in\mathbb{R}$, we have
	\begin{align*}
	R_\mathcal{H}(\tau)(f-\tau Hf)(x) \leq \sup_{\gamma(0)=x}\int_0^\infty e^{-s/\tau}\left[\frac{1}{\tau}f(\gamma(s))-\nabla f(\gamma(s))\cdot \partial_s\gamma(s)\right]\dd s.
	\end{align*}
	Using~$\nabla f(\gamma)\cdot\partial_s\gamma=(d/ds) f(\gamma)$ and integration by parts, we find the estimate
	\begin{equation*}
	R_\mathcal{H}(\tau)(f-\tau Hf)(x) \leq f(x).
	\end{equation*}
	For the reverse inequality, we find a path~$\gamma$ such that~$\gamma(0)=x$ and
	\begin{equation*}
	\int_0^\infty e^{-s/\tau}\left[\frac{1}{\tau}f(\gamma(s))- Hf(\gamma(s)) - \mathcal{L}(\partial_s\gamma(s)) \right]\dd s \geq f(x).
	\end{equation*}
	We will in fact prove equality. Let~$\gamma$ be the path solving
	\begin{align}\label{BG:eq:proof-R-is-pseudores:eq-in-Young}
	\partial_t \gamma(t) = \partial_p\mathcal{H}(\nabla f(\gamma(t))),\qquad t\geq 0.
	\end{align}
	Such a path exists since the vector field~$F(x)=\partial_p\mathcal{H}(\nabla f(x))$ is continuous and bounded---continuity follows by~$\mathcal{H}\in C^1(\mathbb{R})$ and~$f\in C^1(\mathbb{T})$, and boundedness from compactness of~$E=\mathbb{T}$. The precise argument for the existence is given for instance in~\cite[Lemma~3.4]{Kraaij2016}, which is~based on~\cite{Crandall1972}. With this path~$\gamma$,
	\begin{equation}\label{BG:eq:exact-eq-in-Young-inequ}
	Hf(\gamma) + \mathcal{L}(\partial_s\gamma) = \nabla f(\gamma)\cdot\partial_s\gamma,
	\end{equation}
	and~\eqref{BG:eq:proof-R-is-pseudores:eq-in-Young} follows from~$\nabla f(\gamma)\cdot\partial_s\gamma=(d/ds)f(\gamma)$ and integration by parts.
	\smallskip
	
	Since~$\mathcal{H}(0)=0$, the Lagrangian is non-negative. With that, the properties of a contractive pseudo-resolvent are verified by writing out the definitions. For instance, for~$f\in C(E)$ and and any path~$\gamma$,
	\begin{equation*}
	\int_0^\infty\frac{1}{\tau}e^{-s/\tau}\left[f(\gamma(s))-\tau \mathcal{L}(\partial_s\gamma(s))\right]\,\dd s \leq \|f\| \int_0^\infty \frac{1}{\tau}e^{-s/\tau}\,\dd s = \|f\|.
	\end{equation*} 
	Therefore~$\|R(\tau)f\|\leq \|f\|$. With a similar estimate, taking arbitrary~$h_1,h_2$, we find contractivity of~$R(\tau)$ in the sense of~\ref{item:visc-machine:contr}.
	The resolvent identity~\ref{item:visc-machine:resolv-id} follows from rearrangements involving integration by parts. 
	\smallskip
	
	The argument for proving~\eqref{BG:eq:H-semigroup-from-resolvent} is given in~\cite[Lemma~8.18]{FengKurtz2006}. One exploits the fact that with a unit exponential random variable~$\Delta$,
	\begin{equation*}
	R_\mathcal{H}(\tau)h(x) = \sup_{\gamma(0)=x} \mathbb{E}\left[h(\gamma(\tau\Delta)) - \int_0^{\tau\Delta}\mathcal{L}(\partial_s\gamma(s))\,\dd s\right].
	\end{equation*}
	The path~$x^t$ in the proof of Lemma~8.18, cited from~Lemma~8.16, is the one satisfying
	\begin{equation*}
	V_\mathcal{H}(t)f(x) = f(x^t(t))-\int_0^t\mathcal{L}(\partial_sx^t(s))\,\dd s.
	\end{equation*}
	The proof of~\cite[Lemma~8.16]{FengKurtz2006} carries over verbatim; the fact that for every~$x_0$, there exists a path~$\gamma$ satisfying~$\gamma(t_0)=x_0$ and
	\begin{equation}\label{BG:eq:zero-cost-flow-of-L}
	\int_{t_0}^\infty \mathcal{L}(\partial_t\gamma(t))\,\dd t=0
	\end{equation}
	follows in our case from~\eqref{BG:eq:exact-eq-in-Young-inequ} specialized to~$f\equiv 1$.
\end{proof}
%
\section{Bibliographical notes}
\label{BG:sec:bibliographical-notes}
Outlines of the idea of~\cite{FengKurtz2006} are also offered for instance in Feng's paper~\cite{Feng2006} and the introduction of Kraaij's PhD thesis~\cite{Kraaij2016PhDThesis}. The focus in~\cite{FengKurtz2006} lies on conditions for proving large deviation principles for particle systems that lead to Hamilton-Jacobi equations in the space of probability measures, which requires to solve various functional-analytic problems in non-locally compact spaces. Here we comment on some relations to our simplified treatment.
\paragraph{Markov processes via solution to martingale problems.} We specified Markov processes from the existence of transition probabilities. Since these are generally unknown, Markov processes are in practice not obtained by writing down an explicit family of transition probabilities. A common strategy is to find a semigroup satisfying the conditions of~\cite[Proposition~1.3]{Liggett2004}, which by~\cite[Theorem~1.5]{Liggett2004} gives a Markov process defined by means of a family of path distributions. In general, the semigroup determines all finite-dimensional distributions~\cite[Proposition~4.1.6]{EthierKurtz1986}, which induces the path distribution of a stochastic process by the Daniell-Kolmogorov extension Theorem~\cite[Theorem~3.38]{BovierDenHollander2016}. However, frequently we only have an idea about the infinitesimal time evolution of the process. Hence we would like to construct the Markov process by specifying its generator. This point of view is explored for instance in~\cite[Chapter~4]{EthierKurtz1986} and~\cite[Sections~5.3]{BovierDenHollander2016}. The starting assumption in many theorems in~\cite{FengKurtz2006} is the well-posedness of the so-called martingale problem, which associated a path distribution~$\mathbb{P}_\nu$ to a generator~$A$ and an initial distribution~$\nu$. Overviews of the martingale approach can be found for instance in~\cite[Sections~4.3-4.5 and~8]{EthierKurtz1986},~\cite[Section~5.4]{BovierDenHollander2016} and~\cite[Section~1.5]{Liggett2004}.
\paragraph{Large deviations via convergence of semigroups.} Theorem~\ref{thm:LDP-via-semigroup-convergence:compact} is a special case of~\cite[Theorem~5.15]{FengKurtz2006}. The proofs we outlined in Section~\ref{BG:sec:LDP-from-convergence-of-semigroups} carry over to non-compact Polish spaces. The only adaptions are: replace~$C(E)$ with~$C_b(E)$ (continuous and bounded), demand convergence of semigroups bounded and uniformly on compact sets (buc-convergence), and exploit exponential tightness of the initial conditions to conclude as in the proof of Proposition~\ref{prop:LDP-1d-compact:semigroup}. The idea behind the proof of Proposition~\ref{prop:LDP-semigroups:conditional-RF} has been communicated to me by Richard Kraaij. This proof carries over verbatim to the non-compact setting. Finally, a generalization of a collection of compact subsets that is required for non-locally compact spaces is given by~\cite[Definition~2.5]{FengKurtz2006}.
\paragraph{Large deviations via convergence of generators.} In Section~\ref{subsec:LDP-via-classical-sol}, we indicated how Theorem~\ref{thm:LDP-classical-sol} (which is a simplification of~\cite[Corollary~5.19]{FengKurtz2006}) can be proven with Lemmas from~\cite[Section~5]{FengKurtz2006}. The extension to locally-compact state spaces such as~$\mathbb{R}^d$ can be executed by replacing the convergence conditions on the~$H_n$ with uniform convergence on compact subsets---this means specialising as shown in~\cite[Example~2.6]{FengKurtz2006}, which leads to the notion of buc-convergence. The convergence of corresponding semigroups carries over as shown in~\cite[Lemma~5.13~(b)]{FengKurtz2006}. 
\smallskip

A notable difference is that in contrast to the compact case, exponential tightness is no longer a direct consequence of the convergence of~$H_n$. Instead, exponential tightness follows if in addition one can verify the exponential containment condition~\cite[Condition~2.8]{FengKurtz2006}. This condition corresponds to controlling the probability of the process escaping compact sets. There is a convenient way of verifying this condition from the limit operator~$H$; finding a so-called good containment function is sufficient. A detailed account on this is offered in the appendix of Collet's and Kraaij's paper~\cite{CoKr17}, in particular~Proposition~A.15 therein. In our running example on~$\mathbb{R}$, the function~$\Upsilon(x)=\log(1+x^2)$ is a good containment function. More general conditions based on Lyapunov function techniques are given in~\cite[Section~4]{FengKurtz2006}.
\paragraph{Comparison principle in non-compact state space.} There is an extensive literature on comparison principles. The concept of viscosity solutions is outlined in the user's guide of Crandall, Ishii and Lions~\cite{CIL92}. Another introduction to techniques for verifying the comparison principle in~$\mathbb{R}^d$ can be found in Bardi's and Cappuzzo-Dolcetta's monograph~\cite[Chapter~2]{BardiDolcetta1997}. The proof of comparison principle follows the same idea as in the running example at the end of Section~\ref{subsec:LDP-via-visc-sol}; one only has to use a good containment function~$\Upsilon$ to reduce the analysis to compact sets. This point of view is further explained in Chapter~6 and in~\cite[Appendix~A]{CoKr17}. The analysis is more involved in infinite dimensions; see for instance the works of Tataru~\cite{Tataru1992,Tataru1994} and Feng~\cite{Feng2006} in linear spaces, and the recent paper by Feng, Mikami and Zimmer~\cite{FengMikamiZimmer2019} in the space of probability measures, where methods from~\cite[Chapter~13]{FengKurtz2006} are extended.
\paragraph{Action-integral representation.} In Section~\ref{BG:sec:semigroup-flow-HJ-eq}, we worked with~$E=\mathbb{T}$ and proper convex Hamiltonians of the form~$\mathcal{H}=\mathcal{H}(p)$, which allowed us to exploit superlinearity of the Lagrangian defined as the Legendre dual. 
The information contained in the superlinearity is generalized by~\cite[Condition~8.9]{FengKurtz2006}; the compactness of sub-level sets we used in Proposition~\ref{BG:prop:good-Lagrangians-give-action}---for which I could not find a simple proof---is proven under the more general Condition~8.9 in~\cite[Proposition~8.13]{FengKurtz2006}. Furthermore, we assumed differentiability of the Hamiltonian in in Theorem~\ref{BG:thm:resonable-H-gives-pseudo-resolvent} in order to find a path~$\gamma$ satisfying~\eqref{BG:eq:exact-eq-in-Young-inequ}. We used this path to prove~$R_\mathcal{H}(\tau) (\mathbf{1}-\tau H)f\geq f$, and to find the zero-cost flow~\eqref{BG:eq:zero-cost-flow-of-L} required for proving~\eqref{BG:eq:H-semigroup-from-resolvent}. The condition of finding a path~$\gamma$ is generalized by~\cite[Condition~8.11]{FengKurtz2006}, and the existence of a zero-cost flow is generalized by~\cite[Condition~8.10]{FengKurtz2006}. The generalization of Proposition~\ref{BG:prop:good-Lagrangians-give-action} is given by~\cite[Theorem~8.14]{FengKurtz2006}, and the generalizations of verifying equality of semigroups from equality of resolvents are found in~\cite[Corollaries~8.28, 8.29]{FengKurtz2006}.
\paragraph{Pseudo-resolvents.}
Richard Kraaij proofs large deviation principles by generalizing the concept of pseudo-resolvents~\cite{Kraaij2019ExpResolv,Kraaij2019GenConv}. Instead of working with the Hille-Yosida approximants as in the proof of Theorem~\ref{thm:LDP-classical-sol}, Kraaij shows in~\cite{Kraaij2019ExpResolv} how to rigorously obtain the semigroup~$V_n(t)$ from the nonlinear generators~$H_n$ via resolvents~$R_n$. The resolvents are defined by finding viscosity solutions to~$(1-\tau H_n)u=h$ via a control problem involving the relative entropy as a cost---this step replaces the argument in~\cite{FengKurtz2006} of passing to Hille-Yosida approximants. The existence of viscosity solutions follows from generalizing the concept outlined in Theorem~\ref{thm:viscosity-sol-machine}. Our case is a simplified version of the strategy carried out in~\cite[Section~8.4]{FengKurtz2006}. In the generalizations, one proves the fact that the images~$R(\tau)h$ are continuous functions by exploiting the comparison principle, and passing to lower- and upper semicontinuous regularizations first. We avoided these details to clarify the idea. 
\chapter{Large Deviations of Switching Processes}
\label{chapter:LDP-for-switching-processes}
\section{Introduction---molecular motors}
\label{section:intro}
In this chapter we focus on switching Markov processes motivated by stochastic models of walking molecular motors. Broadly speaking, \emph{molecular motors} are proteins that are capable of binding on and moving on filaments in a living cell. Molecular motors such as kinesin and dynein drag vesicles along while moving, and thereby they transport them within the cell. The motors achieve their directed mechanical motion by converting chemical energy of surrounding ATP molecules. In that sense, molecular motors enable living cells to organise directed transport of vesicles.
Jonathan Howard provides an overview of the phenomenon of molecular motors in~\cite{Howard2001}.
\smallskip

There are several mathematical models of molecular motors describing the motor's movement on a filament. J{\"u}licher, Ajardi and Prost review the most common approaches in the Physics literature, with a focus on the relation between models and numerous experimental results~\cite{JulicherAjdariProst1997}. Recent overviews of mathematical models are offered for instance by Anatoly Kolomeisky and Michael Fisher~\cite{Kolomeisky13,KolomeiskyFisher07}. 
\smallskip

Mathematical models can help us to answer questions about the mechanism behind the transport phenomenon based on molecular motors. For instance, is there an underlying common working principle? How do transport properties such as the effective velocity, energy efficiency, stability with respect to perturbations, and response to external forces depend on physical characteristics such as the involved chemical reaction times or the structure of the polymeric filaments? 
\smallskip

Peskin, Dwight and Elston show that the \emph{Brownian ratchet model}
predicts a decrease of the motor's speed when increasing stiffness of the string connecting motor and cargo~\cite{PeskinElston2000}. On the other hand, the \emph{correlation ratchet model} introduced by Peskin, Ermentrout and Oster~\cite{PeskinErmentroutOster1995} predicts an increase of the speed when increasing stiffness~\cite{PeskinDwightElston2000}. Deville and Vanden-Eijnden consider models of multiple motors pulling the same cargo to investigate synchronization effects~\cite{DevilleVandenEijnden2008}. In the models, the gait of a motor regularizes when pulling a cargo, and multiple motors synchronize their actions when pulling together. The authors reveal a similar effect for Brownian ratchet models~\cite{DevilleVandenEijnden2007}.
\smallskip

We focus on correlation ratchet models for a single motor. In Section~\ref{subsec:model-of-molecular-motor} below, we describe a stochastic version of these models. What makes them interesting is the fact that they do not prescribe a directional movement by introducing a uniform bias. The models rather describe a working principle, and the movement can be derived thereof. Given a specific model, the challenge lies therefore in predicting the precise dynamics in the first place.
\smallskip
 
The models are specified by periodic potentials and reaction rates. Simply put, a potential describes the motor's interaction with the periodic filament when being in a certain chemical state, and the reaction rates describe how likely the motor's chemical state changes---we provide more details in Section~\ref{subsec:model-of-molecular-motor}. Many works investigate Fokker-Planck equations associated to the model. For instance, Hastings, Kinderlehrer and Mcleod study their stationary solution and find sufficient conditions for the occurence of transport in terms of distributions of minima of the potentials and a suitable choice of reaction rates~\cite[Theorem~2.1]{HastingsKinderlehrerMcleod08}, and~\cite[Theorems~3.1,~3.2]{HastingsKinderlehrerMcLeod2008}. Wang, Peskin and Elston provide numerical results for such models~\cite{WangPeskinElston2003}. 
Perthame, Souganidis and Mirrahimi use homogenization techniques to characterize transport in terms of a cell problem~\cite{PerthameSouganidis09a, PerthameSouganidis2009Asymmetric, Mirrahimi2013}. Our work relates closest to their results, and we comment further on their work in Section~\ref{MM:sec:larger-context-and-aim}.
\smallskip

We propose to analyse the underlying stochastic models by means of large deviation theory. While molecular motors are naturally stochastic on the microscale, they move with a nearly deterministic velocity on the macroscale. In the mathematical models, this stability is reflected in the fact that the large-deviation results do not depend on the microscopic details of the dynamics. 
\paragraph{Overview of this chapter.}
We first illustrate in Section~\ref{subsec:model-of-molecular-motor} our general results on a specific example of a molecular-motor model. We sketch, without diving into details, how large deviation theory can be utilized to study stochastic models of molecular motors. The example also provides a picture to interpret the general results that follow.
We continue with outlining in Section~\ref{MM:sec:larger-context-and-aim} why we choose to analyse switching Markov processes. Basically, we want to separate the general arguments of large devation theory from the specific arguments depending on the molecular-motor models. 
The main results are presented in Section~\ref{section:MM:results}. We specify switching Markov processes in a periodic setting in Section~\ref{section:switching-MP}, formulate a general large-deviation theorem for the spatial components in Section~\ref{subsection:results:LDP_switching_MP}, and provide an action-integral representation of the rate functions in Section~\ref{subsection:results:action_integral_rep}.
Then we specialize the large deviation theorems to continuous and discrete models of molecular motors (Section~\ref{subsection:results:LDP_in_MM}). Finally, we give an exact formula for the macroscopic velocity in terms of Hamiltonians that are derived from the microscopic dynamics. The detailed-balance condition in molecular motors implies symmetry of the Hamiltonians and thereby of the large-deviation fluctuations. We show that as a consequence, breaking detailed balance is necessary for transport (Section~\ref{MM:subsec:det-bal-implies-symmetric-H}). While this particular conclusion is not new, it comes from a large-deviation perspective. 
We close with the proof sections and an outline on the literature on principal-eigenvalue problems.
	\section{Example---large deviations for molecular motors}
\begin{figure}[h!]
		\labellist
		\pinlabel $t$ at 1800 90
		\pinlabel  $x$ at 760 90
		\pinlabel  $\psi(x,1)$ at -50 230
		\pinlabel  $\psi(x,2)$ at -50 500
		\pinlabel  \small $1.$ at 290 130
		\pinlabel  \small $2.$ at 310 360
		\pinlabel  \small $3.$ at 410 540
		\pinlabel  \small $4.$ at 490 520
		\pinlabel  \small $5.$ at 570 400
		\pinlabel  \small $6.$ at 570 290
		\pinlabel  \small $1.$ at 1300 80
		\pinlabel  \small $2.$ at 1190 150
		\pinlabel  \small $3.$ at 1260 200
		\pinlabel  \small $4.$ at 1100 430 
		\pinlabel  \small $5.$ at 1510 285
		\pinlabel  \small $6.$ at 1330 670
		\pinlabel {\color{red_one}{$X^{\varepsilon}(t)$}} at 1000 640
		\endlabellist
		\centering
		\includegraphics[scale=.18]{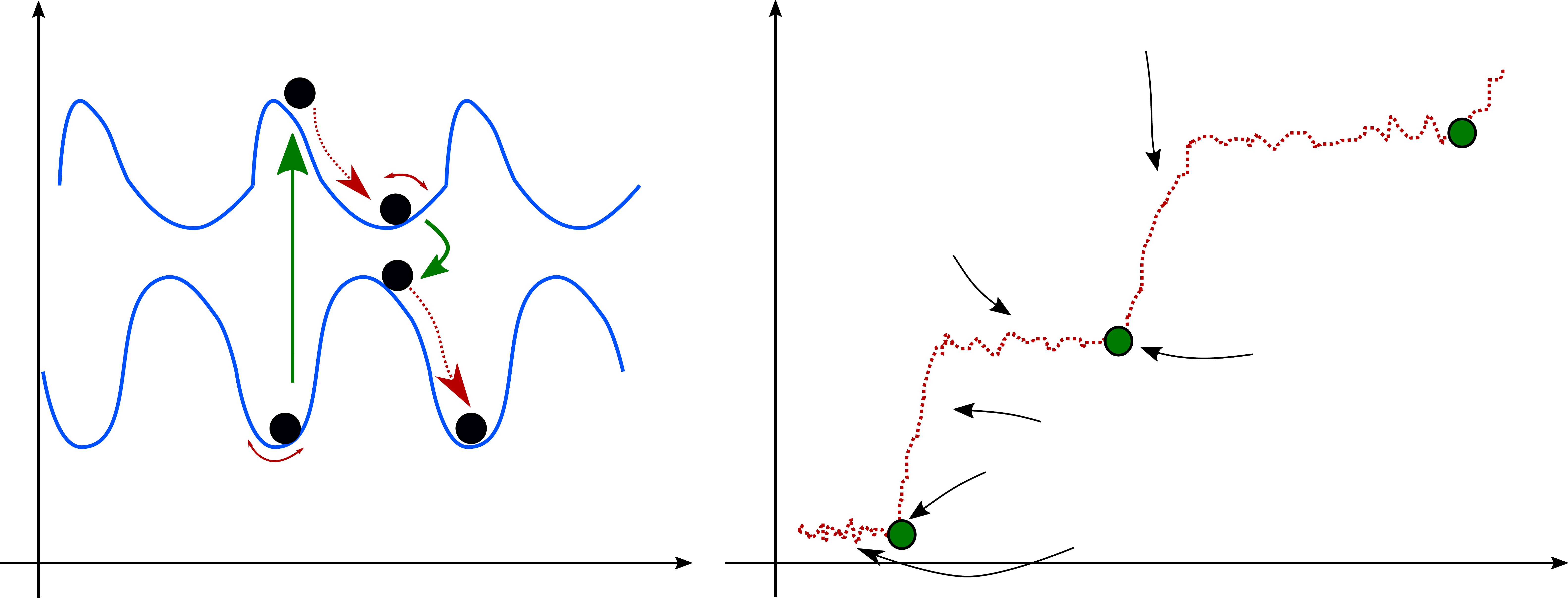}
		\caption{A typical time evolution of $(X^\varepsilon_t, I^\varepsilon_t)$ satisfying~\eqref{eq:intro:example_SDE} and~\eqref{eq:intro:example_evolution_jump_rates}. In the left diagram, the black bullet represents a particle that moves according to~\eqref{eq:intro:example_SDE}. A red arrow corresponds to the value of the spatial component~$X^\varepsilon_t$, and a green arrow indicates a switch of the configurational component~$I^\varepsilon_t$, which changes the potential landscape in which the particle is diffusing. On the right diagram, the spatial component's evolution is shown in a $x$-$t$-diagram, where the red dots represent the values of~$X^\varepsilon_t$. The green bullets indicate a jump of the configurational component~$I^\varepsilon_t$. One forward power stroke consists of the following typical phases: 1. diffusive motion of~$X^\varepsilon$ near minimum; 2. configurational change of~$I^\varepsilon$; 3. flow of~$X^\varepsilon$ towards new minimum.}
		\label{MM:fig:evolution-and-tx-diagram}
	\end{figure}
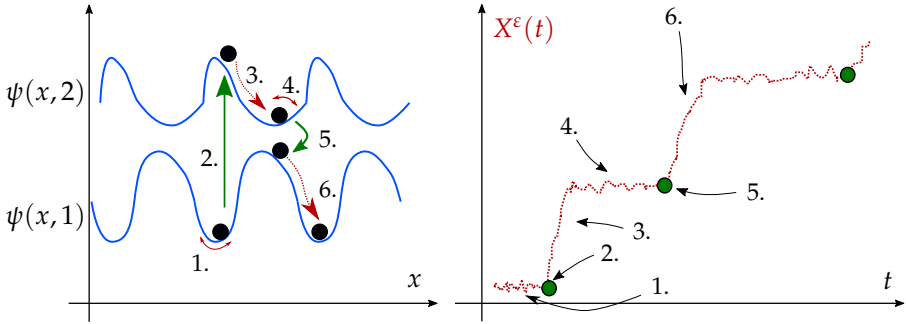

	\sectionmark{Large deviations for molecular motors}
	\label{subsec:model-of-molecular-motor}
	In this example, we consider a two-component Markov process~$(X^\varepsilon,I^\varepsilon)$
	with values in $\mathbb{T} \times \{1,2\}$, where $\mathbb{T} = \mathbb{R} / \mathbb{Z}$ is the one-dimensional flat torus, $\varepsilon = 1/n$ a small parameter, and $n$ an integer. We fix the initial condition; that means for some~$(x_0,i_0) \in \mathbb{T} \times \{1,2\}$, we have~$\left(X^\varepsilon(0),I^\varepsilon(0)\right) =(x_0,i_0)$. Let~$\psi(\cdot,1)$ and~$\psi(\cdot,2)$ be smooth functions on the torus, and write $\psi'(\cdot,i)$ for the derivative of $\psi(\cdot,i) \in C^\infty(\mathbb{T})$. The evolution of $(X^\varepsilon, I^\varepsilon)$ is characterized by the stochastic differential equation
	\begin{equation}\label{eq:intro:example_SDE}
	\dd X_t^\varepsilon = -\psi'\left(\frac{1}{\varepsilon}X_t^\varepsilon,I_t^\varepsilon\right)\,\dd t + \sqrt{\varepsilon}\,\dd B_t.
	\end{equation}
	where $B_t$ is a standard Brownian motion. The process~$I^\varepsilon_t$ is a continuous-time Markov chain on~$\{1,2\}$, which evolves with jump rates~$r_{ij}(\cdot)$ such that
	\begin{equation}
	\mathbb{P}
	\left[I^\varepsilon_{t+\Delta t}
	=
		j\,|\,I^\varepsilon_t =i, X^\varepsilon_t =x
		\right]
		=
		\frac{1}{\varepsilon} r_{ij}\left(\frac{x}{\varepsilon}\right)\Delta t
	+
		\mathcal{O}(\Delta t^2),
	\quad
		\text{as}\;\Delta t \to 0.
	\label{eq:intro:example_evolution_jump_rates}
	\end{equation}
	In summary, the \emph{spatial component}~$X^\varepsilon$ is a drift-diffusion process, the \emph{configurational component}~$I^\varepsilon$ is a jump process on~$\{1,2\}$, and the two are coupled through their respective rates. The drift-term in~\eqref{eq:intro:example_SDE} 
	depends on the value of~$I^\varepsilon$, and thereby the role of~$I^\varepsilon$ is to determine the kind of dynamics that~$X^\varepsilon$ is following. Let~$\{T_k\}_{k=1,2\dots}$ be the jump times of~$I^\varepsilon$, and set~$T_0:=0$. Then~$I^\varepsilon(t)$ is constant in the time windows~$[T_k,T_{k+1})$ (denote its value by~$j_k$), in which the spatial component is a drift-diffusion process with drift term~$-\psi'(\cdot,j_k)$. For details about the construction of such switching diffusions, we refer to~\cite[Chapter~2]{yin2010hybrid}.
	Figure~\ref{MM:fig:evolution-and-tx-diagram} depicts a typical realization of~$(X^\varepsilon, I^\varepsilon)$, where the trajectory of the spatial component is lifted from the torus to~$\mathbb{R}$.
\smallskip

Let us give one possible motivation for the specific $\varepsilon$-scaling. One can start from a process~$(X_t, I_t)$ satisfying
\begin{align*}
\dd X_t&= -\psi'(X_t,I_t)\,\dd t + \dd B_t,
\end{align*}
where the jump process $I_t$ on $\{1,2\}$ evolves according to
\begin{equation*}
\text{Prob}\left[I_{t+\Delta t}=j\,|\,I_t=i, X_t=x\right]=r_{ij}(x)\Delta t+\mathcal{O}(\Delta t^2),\quad i\neq j, 1\leq i,j\leq 2.
\end{equation*}
The large-scale behaviour of $(X_t, I_t)$ is studied by considering the rescaled process~$(X^\varepsilon_t, I^\varepsilon_t)$ defined by~$X^\varepsilon_t := \varepsilon X_{t / \varepsilon}$ and~$I^\varepsilon_t := I_{t / \varepsilon}$.
This rescaling corresponds to zooming out of the $x$-$t$ phase space, which is illustrated below in Figure~\ref{MM:fig:zooming-out}. It\^o calculus implies that the process $(X^\varepsilon_t, I^\varepsilon_t)$ satisfies~\eqref{eq:intro:example_SDE} and~\eqref{eq:intro:example_evolution_jump_rates}.
\begin{figure}[h!]
	\labellist
	\pinlabel $t$ at 1000 50
	\pinlabel $x$ at 50 800
	\pinlabel $t$ at 2125 50
	\pinlabel $x$ at 1170 800
	\pinlabel $v=\partial_p\mathcal{H}(0)$ at 1850 280
	\pinlabel $\varepsilon=1$ at 500 800
	\pinlabel $\varepsilon\ll 1$ at 1600 800 
	\pinlabel {\color{red_one}{$X^{\varepsilon}(t)$}} at 400 600
	\pinlabel {\color{red_one}{$X^\varepsilon(t)$}} at 1500 600
	\endlabellist
	\centering
	\includegraphics[scale=.17]{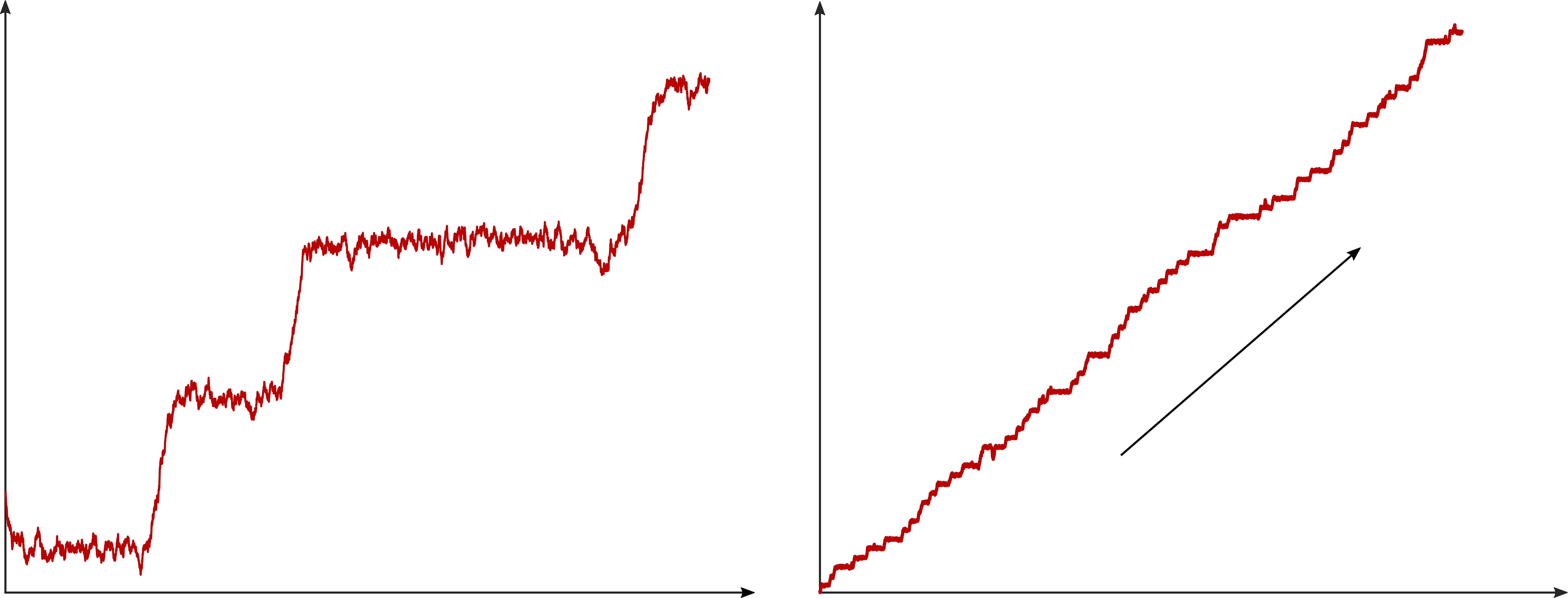}
	\caption{Two typical realizations of the spatial component~$X_t^\varepsilon$ of the two-component process~$(X_t^\varepsilon,I_t^\varepsilon)$ satisfying~\eqref{eq:intro:example_SDE} and~\eqref{eq:intro:example_evolution_jump_rates}. On the left, a realization is depiced for~$\varepsilon$ of order one, and on the right for small~$\varepsilon$. Both graphs depict the lifted trajectory on~$\mathbb{R}$.}
	\label{MM:fig:zooming-out}
\end{figure}
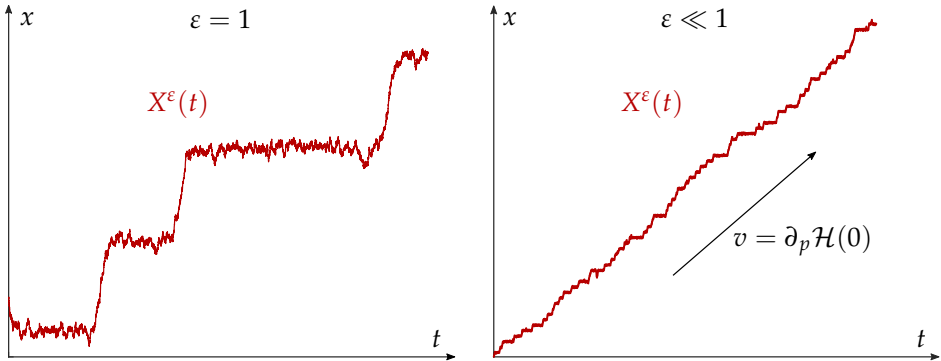
\smallskip

We are interested in the behaviour of the spatial component~$X^\varepsilon$ as $\varepsilon \to 0$. The behaviour of~$X^\varepsilon$ for small~$\varepsilon$ is shown in Figure~\ref{MM:fig:zooming-out}. This figure suggests that for small~$\varepsilon$, the spatial component closely follows a path with a constant velocity. Indeed, when specifying our results of this chapter to the example at hand---the process $(X^\varepsilon, I^\varepsilon)$ defined by~\eqref{eq:intro:example_SDE} and~\eqref{eq:intro:example_evolution_jump_rates}---we find that the spatial component~$X^\varepsilon$ satisfies a pathwise large deviation principle in the limit~$\varepsilon\to 0$. 
\smallskip

To describe this fact more precisely, let~$\mathcal{X}:=C_\mathbb{T}[0,\infty)$ the set of continuous trajectories in~$\mathbb{T}$, equipped with the Skorohod metric, that means the topology of uniform convergence on compact time intervals. The spatial component~$X^\varepsilon$ is a random variable in~$\mathcal{X}$, with a path distribution~$\mathbb{P}(X^\varepsilon\in\cdot)\in\mathcal{P}(\mathcal{X})$. We will show that there exists a rate function~$\mathcal{I}:\mathcal{X}\to[0,\infty]$ with which~$\{X^\varepsilon\}_{\varepsilon>0}$ satisfies a pathwise large deviation principle in the sense of Definition~\ref{def:LDP} from Chapter~1. The gist of this statement is that for any trajectory~$x\in\mathcal{X}$, we have at least intuitively
\begin{equation}\label{MM:eq:intro:example-LDP}
\mathbb{P}\left(X^\varepsilon\approx x\right) \sim e^{-\varepsilon^{-1} \, \mathcal{I}(x)},\quad\varepsilon\to 0.
\end{equation}
The rate function is given by means of a Lagrangian~$\mathcal{L}:\mathbb{R}\to[0,\infty)$,
\begin{equation}\label{MM:eq:intro:example-LDP-RF}
\mathcal{I}(x)=\mathcal{I}_0(x(0)) + \int_0^\infty\mathcal{L}(\partial_t x(t))\,\dd t.
\end{equation}
In there,~$\mathcal{I}_0:\mathbb{T}\to[0,\infty]$ is the rate function of the initial conditions~$X^\varepsilon(0)$, which is given by~$\mathcal{I}_0(x_0)=0$ and~$+\infty$ otherwise---this is because we assume a deterministic initial condition~$X^\varepsilon(0)=x_0$. The Lagrangian is the Legendre dual of a Hamiltonian~$\mathcal{H}:\mathbb{R}\to\mathbb{R}$, that is~$\mathcal{L}(v)=\sup_p [pv-\mathcal{H}(p)]$, and the Hamiltonian is the principal eigenvalue of an associated cell problem described in Lemma~\ref{lemma:LDP_MM:contI:principal_eigenvalue}. We show in Sections~\ref{subsection:results:LDP_switching_MP} and~\ref{section:LDP_in_MM} further below how to obtain the associated cell problem and the principal eigenvalue. 
\smallskip

Here, we focus on how this large-deviation result confirms the claim suggested by Figure~\ref{MM:fig:zooming-out}. The rate function~\eqref{MM:eq:intro:example-LDP-RF} has the following properties:
	\begin{enumerate}[(i)]
	\item $\mathcal{I} : \mathcal{X} \rightarrow [0,\infty]$ is nonnegative.
	\item $\mathcal{I}(x) = 0$ if and only if~$\partial_tx(t)=v$, with~$v=\partial_p\mathcal{H}(0)$.
	\end{enumerate}
These two properties together characterize the unique minimizer of the rate function, and thereby in particular the typical behaviour of~$X^\varepsilon$ for small~$\varepsilon$. Whenever~$\mathcal{I}(x) > 0$ for a path~$x\in\mathcal{X}$, then by~\eqref{MM:eq:intro:example-LDP}, the probability that a realization of~$X^\varepsilon$ is close to~$x$ in the Skorohod metric is exponentially small in $\varepsilon$. More precisely, the large deviation principle implies almost-sure convergence of~$X^\varepsilon$ to the unique minimizer of the rate function (Theorem~\ref{thm:math-formulation-LDP:LDP-implies-as}).
\smallskip

Equipped with the large deviation principle, we can investigate which sets of potentials and rates~$\{\psi_1,\psi_2,r_{12},r_{21}\}$ induce transport, which means observing a non-zero macroscopic velocity~$v=\partial_p\mathcal{H}(0)$. We do not find general sufficient conditions for transport, but can draw some conclusions if the process~$(X^\varepsilon,I^\varepsilon)$ satisfies \emph{detailed balance}. Here, detailed balance is satisfied if~$r_{12}e^{-\psi_1} = C r_{21}e^{-\psi_2}$ for some constant~$C>0$. This condition implies time-reversibility of~$(X^\varepsilon,I^\varepsilon)$ in the sense of Definition~\ref{intro:def:reversibility}---we clarify this connection in Section~\ref{MM:subsec:det-bal-implies-symmetric-H}. There, we will also show that detailed balance implies symmetry of the Hamiltonian, that is~$\mathcal{H}(p)=\mathcal{H}(-p)$. In particular, we find~$v=0$ under detailed balance.
\smallskip


We close this section by describing how the stochastic process~$(X^\varepsilon, I^\varepsilon)$ models the movement of a molecular motor on a polymeric filament. The molecular motor consists of two chains whose heads can attach to the filament. A sequence of chemical reactions provides energy that triggers 
a forward power stroke of a motor head, thereby leading to a spatial displacement of the motor.
The spatial component~$X^\varepsilon_t$ corresponds to the position of the molecular motor on the filament, while a change of the configurational component~$I^\varepsilon_t$ corresponds to the event of a chemical reaction.
The information of how the motor moves forward is encoded in the potentials. The periodicity of the potential reflects the periodic structure of a filament.
One can think of the minima of the potentials as corresponding to the head's binding spots on the filament. Finally, the noise term in~\eqref{eq:intro:example_SDE} models friction arising from collisions of the motor with molecules in the environment. Because of the highly viscous environment, it is common to consider a drift-diffusion process. A justification for this overdamped limit regime is offered for instance by Wang and Elston~\cite{WangElston2007}. When coarse-graining the continuous model to a Markov jump process on the binding spots, we obtain a discrete model. We discuss these jump models in Section~\ref{subsection:results:LDP_in_MM}.
\section{Larger context and aim of this chapter}
\label{MM:sec:larger-context-and-aim}
	One inspiration for the subject of this chapter is a series of papers by Perthame, Souganidis and Mirrahimi~\cite{PerthameSouganidis09a, PerthameSouganidis2009Asymmetric, Mirrahimi2013}. There, the authors start from the Fokker-Planck equations associated with $(X^\varepsilon, I^\varepsilon)$ from \eqref{eq:intro:example_SDE} and \eqref{eq:intro:example_evolution_jump_rates}:
	\begin{align}
	\begin{cases}
		\displaystyle \partial_t \rho_\varepsilon^1 
	=
		\;\varepsilon \frac{1}{2} \partial_{xx} \rho_\varepsilon^1
		+
		\text{div}_x
		\left[\rho_\varepsilon^1\,\psi'_1\left(\frac{x}{\varepsilon}\right)\right]
		+
		\frac{1}{\varepsilon}r_{21}\left(\frac{x}{\varepsilon}\right)\rho_\varepsilon^2
		-
		\frac{1}{\varepsilon}r_{12}\left(\frac{x}{\varepsilon}\right)\rho_\varepsilon^1,\\
		\quad \\
		\displaystyle \partial_t \rho_\varepsilon^2 
	=
		\; \varepsilon \frac{1}{2} \partial_{xx} \rho_\varepsilon^2
		+
		\text{div}_x
		\left[\rho_\varepsilon^2\,\psi'_2\left(\frac{x}{\varepsilon}\right)\right]
		+
		\frac{1}{\varepsilon}r_{12}\left(\frac{x}{\varepsilon}\right)\rho_\varepsilon^1
		-
		\frac{1}{\varepsilon}r_{21}\left(\frac{x}{\varepsilon}\right)\rho_\varepsilon^2.
	\end{cases}
	\label{eq:intro:FP_eq}
	\end{align}
	The functions $\{\psi^i,r_{ij}\}$ are taken to be $1$-periodic and smooth. The system of equations \eqref{eq:intro:FP_eq} describes the evolution of the partial probability densities given in terms of the process by~$\rho^i_\varepsilon(t,\dd x) = \mathbb{P}\left(X^\varepsilon_t \in \dd x, I^\varepsilon_t = i\right)$.
	\smallskip
	
	Perthame and Souganidis define in~\cite{PerthameSouganidis09a} a notion of asymmetry for a given set of functions $\{\psi^i,r_{ij}\}$. This notion is based on migration of density in the stationary Fokker-Planck system ($\partial_t\rho_\varepsilon^i = 0 $ in~\eqref{eq:intro:FP_eq}) on the spatial domain~$(0,1)$ with periodic boundaries.
	The authors consider~$\psi_2=0$ and find a condition under which the densities $\rho^1_\varepsilon$ and $\rho^2_\varepsilon$ converge to a delta mass supported at one end of the interval, which is refered to as the \emph{motor effect} or as \emph{transport}. In all three papers, Perthame, Souganidis and Mirrahimi address the question of what exactly characterizes the class of potentials and rates $\{\psi^i,r_{ij}\}$ that induce transport, and prove convergence statements for the Fokker-Planck system \eqref{eq:intro:FP_eq}. 
	\smallskip
	
	In \cite{PerthameSouganidis09a}, the authors find a sufficient condition for transport in terms of an effective Hamiltonian $\mathcal{H}(p)$ and a total flux $F(p)$, where $p\in\mathbb{R}$. System $\eqref{eq:intro:FP_eq}$ exhibits the motor effect if and only if $\partial_p\mathcal{H} (0)\neq 0$, or equivalently if $F(0)\neq 0$. The effective Hamiltonian is the principal eigenvalue of an associated cell problem, obtained after an exponential change of variables. It is the same principal eigenvalue that appears in the example from above, and we explain in Section~\ref{section:LDP_in_MM} how to obtain it from a large-deviation perspective. Because they consider the stationary system, the information about how fast the density migrates cannot be determined, since that is a question about the dynamics.
	\smallskip
	
	More recently, in~\cite{Mirrahimi2013}, Mirrahimi and Souganidis analysed the system~\eqref{eq:intro:FP_eq} on~$\mathbb{R}^d$, again with~$\psi_2=0$. When taking the limit~$\varepsilon \rightarrow 0$, they find that the sum of partial probabilities converges to a moving delta mass with velocity $v = \partial_p\mathcal{H}(0)$. More precisely, they find~$\rho^1_\varepsilon(t,x)+\rho^2_\varepsilon(t,x)\rightarrow \delta(x-tv)I_0$
	in the sense of measures, where $I_0$ is determined by the initial data.
	This is consistent with the previously found criterion for the motor effect \cite{PerthameSouganidis09a}, $\partial_p\mathcal{H}(0) \neq 0$. Theorem~\ref{thm:results:LDP_cont_MM_Limit_I} further below recovers this result with a stronger form of convergence.
	\smallskip
	
	We point out again that we do not provide any new sufficient conditions for obtaining transport, due to the larger generality of our considerations. We do prove under general conditions that detailed balance leads to a symmetric Hamiltonian (see Theorem~\ref{thm:results:detailed_balance_limit_I} below). This implies that detailed balance has to be broken in order for transport to occur. 
	\smallskip
	
	The methods that Perthame, Souganidis and Mirrahimi apply in~\cite{PerthameSouganidis09a, PerthameSouganidis2009Asymmetric, Mirrahimi2013} are inspired by large deviation theory. However, in their papers, they do not explicitly prove large deviations, but prove convergence statements on the level of Fokker-Planck equations. When proving the associated large deviation principles, as we will do in this chapter, there is a clear distinction between the contributions that come from general large deviation theory on the one hand, and the model-specific contributions on the other hand.
	\smallskip
	
	Our aim is not only to prove the large-deviation results, but also to separate those parts of the argument which are general and come from large deviation theory, from those parts that are specific to the model at hand. We make this explicit by considering so-called Markov processes with random switching, a class of stochastic processes that we introduce in Section~\ref{section:switching-MP}. The process introduced above by \eqref{eq:intro:example_SDE} and \eqref{eq:intro:example_evolution_jump_rates} is an example of such a process, and in particular represents a motivating example for considering this class of processes. In Section~\ref{subsection:results:LDP_switching_MP}, we illustrate by means of example how the argument is then separated into large-deviation parts and model-specific parts.
\section{Main Results}
\label{section:MM:results}
	In this section, we give an overview of our results. We first define in Section~\ref{section:switching-MP} switching Markov processes. We formulate and explain sufficient conditions under which the spatial component of a switching Markov process satsfies a large deviation principle (Theorem~\ref{thm:results:LDP_switching_MP}). Since the rate functions are a priori of intricate form, we cast the rate functions in action-integral form (Theorem~\ref{thm:results:action_integral_representation}). In summary, by working with switching Markov processes we show which properties the large deviation principles for molecular motors depend on.
	\smallskip
	
	We specialize to models of molecular motors in Section~\ref{subsection:results:LDP_in_MM}, where we state large deviation principles for two limit regimes. The Hamiltonians in the action-integral rate functions are principal eigenvalues of certain cell problems. In Section~\ref{MM:subsec:det-bal-implies-symmetric-H}, we work with variational formulas of such principal eigenvalues in order to study the behaviour of molecular motors under the detailed-balance condition. The main challenge is to derive useful formulas for the Hamiltonians that allow us to draw concrete conclusions. We show symmetry of the Hamiltonians under detailed balance. In particular, this implies $v = \partial_p\mathcal{H}(0)=0$, which means transport can only occur if detailed balance is broken. This result about transport is expected and not new, but follows in our case from a more general symmetry of large deviations.
	\subsection{Switching Markov processes in a periodic setting}
	\label{section:switching-MP}
	We introduce switching Markov processes as certain two-component stochastic processes~$(X^\varepsilon, I^\varepsilon)$ taking values in a state space $E_\varepsilon$. The state space of~$I^\varepsilon$ is a finite set~$\{1,\dots,J\}$, while~$X^\varepsilon$ takes values in some compact Polish space~$E_\varepsilon^X$.  
	We are interested in studying processes in a periodic setting. Therefore, we consider the flat $d$-dimensional torus $\mathbb{T}^d := \mathbb{R}^d / (\ell \cdot \mathbb{Z}^d)$, for some fixed length~$\ell\in\mathbb{N}$. We henceforth omit the dependence on~$\ell$. 
	\begin{condition}[Setting]
		Fix~$J\in\mathbb{N}$ and let~$\varepsilon = 1/n > 0$ for an integer~$n$. The state space~$E_\varepsilon$ is a product space~$E_\varepsilon := E_\varepsilon^X \times \{1, \dots, J\}$, where~$E_\varepsilon^X$ be a compact Polish space satisfying the following. There are continuous maps $\iota_\varepsilon : E_\varepsilon^X \rightarrow \mathbb{T}^d$ such that for all $x \in \mathbb{T}^d$ there exist $x_\varepsilon \in E_\varepsilon^X$ with which $\iota_\varepsilon(x_\varepsilon) \rightarrow x$ as $\varepsilon \rightarrow 0$.\qed
		\label{condition:results:general_setting}
	\end{condition}
	This condition means that~$E_\varepsilon^X$ is asymptotically dense in the torus~$\mathbb{T}^d$. The typical example is a finite, discrete and periodic lattice with spacing~$\varepsilon$, so that in the limit of~$\varepsilon$ to zero one obtains the torus. Another example is simply~$E_\varepsilon^X\equiv \mathbb{T}^d$. When it is clear from the context, we omit $\iota_\varepsilon$ in the notation.
	\smallskip
	
	We now define switching Markov processes by specifying their generators from the following ingredients:
	\begin{enumerate}[(1)]
		\item For~$i\in\{1,\dots,J\}$, we have a map~$L_\varepsilon^{i}:\mathcal{D}(L_\varepsilon^{i}) \subseteq C(E_\varepsilon^X) \rightarrow C(E_\varepsilon^X)$ that is the generator of an~$E_\varepsilon^X$-valued Markov process. 
		\item For~$i,j\in\{1,\dots,J\}$, we have a continuous map~$r_{ij}^\varepsilon:E_\varepsilon^X\to[0,\infty)$.
	\end{enumerate}
	With that, define the map~$L_\varepsilon:\mathcal{D}(L_\varepsilon) \subseteq C(E_\varepsilon) \rightarrow C(E_\varepsilon)$ by
	\begin{equation}
	L_\varepsilon f(x,i) := L_\varepsilon^{i} f(\cdot,i) (x) + \sum_{j = 1}^J r_{ij}^\varepsilon(x) \left[ f(x,j) - f(x,i) \right],
	\label{eq:intro:MP_with_switching_L_varepsilon}
	\end{equation}
	where the domain is~$\mathcal{D}(L_\varepsilon) = \{ f \in C(E_\varepsilon) \, : \, f(\cdot,i) \in \mathcal{D}(L_\varepsilon^{i}), i = 1, \dots, J \}$.
	Let~$\mathcal{X}_\varepsilon:=D_{E_\varepsilon}[0,\infty)$ be the set of trajectories in~$E_\varepsilon$ that are right-continuous and have left limits. We equip~$\mathcal{X}_\varepsilon$ with the Skorohod topology~\cite[Section~3.5]{EthierKurtz1986}. For an initial condition~$\mu\in\mathcal{P}(E_\varepsilon)$, we associate a two-component process $(X_t^\varepsilon, I_t^\varepsilon)$ with values in $E_\varepsilon$ to the generator~$L_\varepsilon$ by finding its path distribution~$\mathbb{P}_\mu\in\mathcal{P}(\mathcal{X}_\varepsilon)$. To do so, we assume well-posedness of the associated martingale-problem associated to the pair~$(L_\varepsilon,\mu)$. For the precise statement of the martingale problem, we refer to~\cite[Section~4.3]{EthierKurtz1986}.
	\begin{condition}[Well-posedness]
		Let~$\mu \in \mathcal{P}(E_\varepsilon)$. Then existence and uniqueness holds of the $\mathcal{X}_\varepsilon$-martingale problem for $(L_\varepsilon,\mu)$. Denote the solution to the martingale-problem solution of $L_\varepsilon$ by~$\mathbb{P}_\mu$. The map~$E_\varepsilon \ni z \mapsto \mathbb{P}_{\delta_z} \in \mathcal{P}(\mathcal{X}_\varepsilon)$ is Borel measurable with respect to the weak topology on $\mathcal{P}(\mathcal{X}_\varepsilon)$.\qed
		\label{condition:intro:sol_martinagle_problem}
	\end{condition}
	%
	Condition~\ref{condition:intro:sol_martinagle_problem} is the basic assumption on the processes in~\cite{FengKurtz2006}. A sufficient condition for the measurability is given in~\cite[Theorem~4.4.6]{EthierKurtz1986}. In this chapter, we consider switching Markov processes in the following sense.
	\begin{definition}[Switching Markov processes in a periodic setting]
		Let~$(X^\varepsilon,I^\varepsilon)$ be a two-component Markov proces taking values in~$E_\varepsilon = E^X_\varepsilon \times \{1,\dots,J\}$ satisfying Condition~\ref{condition:results:general_setting}. We call~$(X^\varepsilon,I^\varepsilon)$ a \emph{switching Markov process} if its generator~$L_\varepsilon$ is given by~\eqref{eq:intro:MP_with_switching_L_varepsilon} and satisfies Condition~\ref{condition:intro:sol_martinagle_problem}.\qed
		\label{def:intro:MP_with_switching}
	\end{definition}
	We do not give general conditions on a map~$L_\varepsilon$ that imply Condition~\ref{condition:intro:sol_martinagle_problem}. However, all the examples of stochastic processes modelling molecular motors will satisfy this condition. Further details about existence and regularity properties can be found in the book of Yin and Zhu about switching hybrid diffusions~\cite[Part~I]{yin2010hybrid}.
	\smallskip
	
	We close this section by mentioning that the process of the introductory example satisfying~\eqref{eq:intro:example_SDE} and~\eqref{eq:intro:example_evolution_jump_rates} is a switching Markov process.
	The state space is~$E_\varepsilon=\mathbb{T}\times\{1,2\}$, and its generator is of the form~\eqref{eq:intro:MP_with_switching_L_varepsilon}. The jump rates are given by~$r_{ij}^\varepsilon(x) = r_{ij}(x/\varepsilon) / \varepsilon$, and for~$i=1,2$, we have
	\begin{equation*}
	L_\varepsilon^{i}g(x) := -\psi'(x/\varepsilon,i)\,g'(x) + \varepsilon\,g''(x),\quad x\in\mathbb{T}.
	\end{equation*}
	This is the generator of a drift-diffusion process on $\mathbb{T}$ satisfying
	\begin{equation*}
	\dd Y_t^\varepsilon = -\psi'\left((Y_t^\varepsilon/\varepsilon),i\right)\,\dd t + \sqrt{2\varepsilon}\,\dd B_t.
	\end{equation*}
	A scheme of how to obtain a process on the torus~$\mathbb{T}$ is presented for instance in~\cite[Chapter~3.2]{bensoussan2011asymptotic}.
	\subsection{Large deviation principle for switching Markov processes}
	\label{subsection:results:LDP_switching_MP}
	We consider switching Markov processes~$(X^\varepsilon,I^\varepsilon)$ in the sense of Definition~\ref{def:intro:MP_with_switching}, with generators of the form~\eqref{eq:intro:MP_with_switching_L_varepsilon}. The essence of this section is Theorem~\ref{thm:results:LDP_switching_MP}, which provides general conditions under which a pathwise large deviation principle of the spatial component~$X^\varepsilon$. We alert the reader here that we illustrate the concepts and notations by means of an example below Theorem~\ref{thm:results:LDP_switching_MP}. We state the conditions in terms of the nonlinear generators defined as follows.
	\begin{definition}[Nonlinear generators]
		\label{MM:def:nonlinear-generators-switching-MP}
		Let~$L_\varepsilon$ be the map defined by~\eqref{eq:intro:MP_with_switching_L_varepsilon}. The \emph{nonlinear generator} is the map~$H_\varepsilon : \mathcal{D}(H_\varepsilon) \subseteq C(E_\varepsilon) \rightarrow C(E_\varepsilon)$ defined by
		\begin{equation}
		H_\varepsilon f(x) := \varepsilon \, e^{-f(x)/\varepsilon} L_\varepsilon (e^{f(\cdot)/\varepsilon})(x),
		\label{eq:results:H_varepsilon}
		\end{equation}
		with the domain~$\mathcal{D}(H_\varepsilon) := \{f \in C(E_\varepsilon) \, : \, e^{f(\cdot)/\varepsilon} \in \mathcal{D}(L_\varepsilon)\}$.\qed
	\end{definition}
	We will require the nonlinear generators~$H_\varepsilon$ to converge in the limit~$\varepsilon\to 0$. To formulate this convergence condition, we need to introduce an additional state space~$E'$ for collecting up-scaled variables. The following diagram depicts the relation between the state spaces:
	\begin{equation*}
		\begin{tikzcd}[row sep = tiny]
			&	\mathbb{T}^d \times E^\prime \arrow[dd, "\mathrm{proj}_1"]	\\
		E_\varepsilon \arrow[ur, "(\eta_\varepsilon {,} \eta_\varepsilon^\prime)"] \arrow[dr, "\eta_\varepsilon"']	&	\\
			&	\mathbb{T}^d	
		\end{tikzcd}
	\end{equation*}
	In the diagram,~$\eta_\varepsilon : E_\varepsilon \to \mathbb{T}^d$ is the projection defined by~$\eta_\varepsilon (x,i) := \iota_\varepsilon(x)$, where~$\iota_\varepsilon:E_\varepsilon^X\to\mathbb{T}^d$ is the embedding. The map $\eta_\varepsilon^\prime : E_\varepsilon \to E^\prime$ is continuous. We assume that~$E_\varepsilon$ is asymptotically dense:
	\begin{enumerate}[(C1)]
		\item \label{MM:item:C1} For~$(x,z') \in \mathbb{T}^d \times E'$ there are $y_\varepsilon\in E_\varepsilon$ such that $\eta_\varepsilon(y_\varepsilon) \to x$ and $\eta_\varepsilon'(y_\varepsilon) \to z^\prime$ as $\varepsilon \to 0$.
	\end{enumerate}
	The limit operator of~$H_\varepsilon$ are generally multivalued, that means defined by a subset~$H\subseteq C(\mathbb{T}^d)\times C(\mathbb{T}^d\times E')$. We assume the following convergence condition:
	\begin{enumerate}[(C2)]
		\item \label{MM:item:C2} For any~$(f,g)\in H$, there are functions~$f_\varepsilon\in\mathcal{D}(H_\varepsilon)$ such that
		\[
		\| f\circ\eta_\varepsilon - f_\varepsilon \|_{L^\infty(E_\varepsilon)} \rightarrow 0
		\quad
		\text{ and }
		\quad
		\| g\circ (\eta_\varepsilon,\eta_\varepsilon') - H_\varepsilon f_\varepsilon \|_{L^\infty(E_\varepsilon)} \rightarrow 0.
		\]
	\end{enumerate}
	Frequently, for any~$f$ in the domain of~$H$, the corresponding image functions~$g$ are naturally parametrized by a set of functions on~$E'$. 
	\begin{enumerate}[(C3)]
		\item \label{MM:item:C3} There are a set~$\mathcal{C}\subseteq C(E')$ and functions~$H_{f,\varphi}\in C(\mathbb{T}^d\times E')$ with which
		\begin{equation*}
		H = \left\{\left(f,H_{f,\varphi}\right)\,:\,f\in\mathcal{D}(H),\varphi\in\mathcal{C}\right\}.
		\end{equation*}
	\end{enumerate}
	Below the theorem, we illustrate by an example how to find the multivalued operator starting from the nonlinear generators.
	\begin{theorem}[Large deviation principle for switching processes]
	Let~$(X^\varepsilon,I^\varepsilon)$ be a switching Markov process in the sense of Definition~\ref{def:intro:MP_with_switching}, with nonlinear generators~$H_\varepsilon$ of Definition~\ref{MM:def:nonlinear-generators-switching-MP}. Suppose that there exists a compact metric space~$E^\prime$ satisfying~\ref{MM:item:C1} and a multivalued operator~$H\subseteq C(\mathbb{T}^d)\times C(\mathbb{T}^d\times E')$ with domain~$\mathcal{D}(H)$ satisfying~$C^\infty(\mathbb{T}^d) \subseteq \mathcal{D}(H) \subseteq C^1(\mathbb{T}^d)$ such that:
		\begin{enumerate}[(T1)]
		\item \label{MM:item:T1}
		The operator~$H$ satisfies~\ref{MM:item:C2} and~\ref{MM:item:C3} from above.
		For every $\varphi\in\mathcal{C}$ there is a map~$H_\varphi : \mathbb{R}^d \times E^\prime \to \mathbb{R}$
		such that for all~$f \in \mathcal{D}(H)$,
		\begin{equation*}
		H_{f,\varphi}(x,z^\prime) = H_\varphi(\nabla f(x),z^\prime),\qquad (x,z') \in \mathbb{T}^d\times E'.
		\end{equation*}
		\item \label{MM:item:T2}
		For every $p \in \mathbb{R}^d$, there exists a function $\varphi_p \in \mathcal{C}$ and a constant $\mathcal{H}(p) \in \mathbb{R}$ such that~$H_{\varphi_p}(p,z^\prime)=\mathcal{H}(p)$
		for all $z^\prime \in E^\prime$. 
		\end{enumerate}
		Suppose furthermore that~$\{X^\varepsilon(0)\}_{\varepsilon>0}$ satisfies a large deviation principle in~$\mathbb{T}^d$ with rate function~$\mathcal{I}_0 : \mathbb{T}^d \rightarrow [0,\infty]$.
		Then the family of processes $\{X^\varepsilon\}_{\varepsilon>0}$ satisfies a large deviation principle in $D_{\mathbb{T}^d}[0,\infty)$ with a rate function $\mathcal{I} :D_{\mathbb{T}^d}[0,\infty)\rightarrow[0,\infty]$.
	\label{thm:results:LDP_switching_MP}
	\end{theorem}
	The proof of Theorem~\ref{thm:results:LDP_switching_MP} is given in Section~\ref{section:LDP_switching_MP}. The formula for the rate function is not important here, which is why we give it only in the proofs. While Condition~\ref{MM:item:T1} corresponds to the convergence of nonlinear generators, Condition~\ref{MM:item:T2} usually corresponds to solving a principal-eigenvalue problem. The constant~$\mathcal{H}(p)$ is then uniquely determined as the principal eigenvalue of a certain cell problem. Further below, we give feasible conditions on the map $p\mapsto \mathcal{H}(p)$ under which the rate function admits an action-integral representation (Theorem~\ref{thm:results:action_integral_representation} in Section~\ref{subsection:results:action_integral_rep}). Here, we illustrate by an example how conditions~\ref{MM:item:T1} and~\ref{MM:item:T2} can be obtained starting from the nonlinear generators. Even though the example is not a switching Markov process, it shows the features arising from the mixed scales.
	\paragraph{Example illustrating the general case.}
	Let~$\mathbb{T}$ be the one-dimensional flat torus, $\psi(\cdot) \in C^\infty(\mathbb{T})$, and 
	consider the process~$X_t^\varepsilon$ solving
	\begin{equation*}
	\dd X_{t}^\varepsilon = -\psi'(X_t^\varepsilon / \varepsilon) \, \dd t + \sqrt{\varepsilon} \,\dd B_{t},
	\end{equation*}
	where~$\varepsilon = 1/n$ with some integer $n$. 
	Its generator is given by
	\begin{equation*}
	L_\varepsilon f(x) = -\psi'(x/\varepsilon) f'(x) + \varepsilon\, \frac{1}{2} f''(x),\quad f\in C^2(\mathbb{T}).
	\end{equation*}
	Therefore the nonlinear generators $H_\varepsilon$ are
	\begin{equation*}
	H_\varepsilon f(x)=-\psi'(x/\varepsilon) f'(x)+\frac{1}{2}|f'(x)|^2 +
	\varepsilon\,\frac{1}{2} f''(x).
	\end{equation*}
	The aim is to obtain a limit of $H_\varepsilon$ as $\varepsilon \to 0$. In order to determine the behaviour of $H_\varepsilon f$ for small $\varepsilon$, we have to deal with the problem that the drift-term $\psi'(x/\varepsilon)$ is fastly oscillating as $\varepsilon$ tends to zero. This is solved by considering functions that are of the form~$f_\varepsilon (x) = f(x) + \varepsilon\, \varphi(x/\varepsilon)$. Then we obtain
		\begin{multline*}
		H_\varepsilon f_\varepsilon (x)=-\psi'(x/\varepsilon) \cdot
		\left[f'(x) + \varphi'(x/\varepsilon)\right]+\frac{1}{2} | f'(x) + \varphi'(x/\varepsilon)|^2\\
		+\frac{1}{2} \varphi''(x/\varepsilon)+\frac{\varepsilon}{2} f''(x).
		\end{multline*}
		We want these images to converge in the limit~$\varepsilon \to 0$. The~$\varepsilon \,\frac{1}{2} f''(x)$ term is of order~$\varepsilon$ and therefore not problematic. The remaining terms are in general oscillating in~$\varepsilon$. However, with the right choice of the function~$\varphi$, one can make this term to be independent of the~$(x/\varepsilon)$-variable, and thereby independent of~$\varepsilon$ altogether. In order to see how, we rewrite~$H_\varepsilon f_\varepsilon$ by introducing the fast spatial variable~$y=x/\varepsilon$, with which we find that
		\begin{multline}
		H_\varepsilon f_\varepsilon (x)=e^{-\varphi(y)}\left[
		\frac{1}{2} |f'(x)|^2 -  \psi'(y) f'(x)+ \left(f'(x) - \psi'(y)\right) \partial_y+ \frac{1}{2} \partial_{yy}\right]e^{\varphi(y)}\\
		+\varepsilon\,\frac{1}{2} f''(x)
		\qquad\left(y=\frac x\varepsilon\right).
		\end{multline}
		Hence we aim to find a function $\varphi(\cdot)$ such that the term $e^{-\varphi} [ \cdots ] e^\varphi$ is constant as a function of the $x/\varepsilon$-variable, regarding the $x$-variable as a parameter. This term depends on $x$ only via the derivative of $f$. Hence if we can find such a function~$\varphi$, we can denote the constant by~$\mathcal{H}(\partial_x f(x))$. Then with that choice of the function~$\varphi$, the values of $H_\varepsilon f_\varepsilon$ are given by
		\begin{equation*}
		H_\varepsilon f_\varepsilon (x)
		=
		\mathcal{H}(\partial_x f(x)) + \mathcal{O}(\varepsilon),
		\end{equation*}
		and we find for small $\varepsilon$ that $H_\varepsilon f_\varepsilon(x) \approx \mathcal{H}(\partial_x f(x))$. Making this strategy rigorous can be realized in two steps via~\ref{MM:item:T1} showing convergence of nonlinear generators and~\ref{MM:item:T2} solving a principal-eigenvalue problem, as follows:
		\smallskip
		
		\ref{MM:item:T1}: The images $H_\varepsilon f_\varepsilon(x)$ are given by $H_{f,\varphi}(x,x/\varepsilon) + \varepsilon\, \frac{1}{2} f''(x)$, where
		$$
		H_{f,\varphi}(x,y)
		:=
		-\psi'(y)
		\left[
		f'(x) + \varphi'(y)
		\right]
		+
		\frac{1}{2} |f'(x) + \varphi'(y)|^2
		+
		\frac{1}{2} \varphi''(y).
		$$
		By taking arbitrary $\varphi \in C^2(\mathbb{T})$, we collect all these possible limits of $H_\varepsilon f_\varepsilon$ and summarize them in a multivalued operator 
		$
		H \subseteq C(\mathbb{T}) \times C(\mathbb{T} \times \mathbb{T})
		$
		defined by
		\begin{equation}
		H := \{ (f,H_{f,\varphi}) \, : \, f \in C^2(\mathbb{T}) \text{ and } \varphi \in C^2(\mathbb{T}) \}.
		\label{eq:results:example_H}
		\end{equation}
		The set of upscaled variables is~$E'=\mathbb{T}$, and~\ref{MM:item:C1} is satisfied with~$\eta_\varepsilon'(x)=x/\varepsilon$. The nonlinear generator~$H_\varepsilon$ converges to the limit operator~$H$ as demanded in~\ref{MM:item:C2}: for~$(f,H_{f,\varphi}) \in H$, the functions~$f_\varepsilon (x) = f(x) + \varepsilon\, \varphi(x/\varepsilon)$ satisfy
		\begin{equation*}
		\sup_{x \in \mathbb{T}}
		|
		f(x) - f_\varepsilon(x)
		|
		\xrightarrow{\varepsilon\to 0} 0
		\quad\text{and}\quad
		\sup_{x \in \mathbb{T}}
		|
		H_{f,\varphi}(x,\eta^\prime_\varepsilon(x)) - H_\varepsilon f_\varepsilon (x)
		|
		\xrightarrow{\varepsilon\to 0} 0.
		\end{equation*}
		Condition~\ref{MM:item:C3} is satisfied by construction, with~$\mathcal{C}=C^2(\mathbb{T})$. Finally, the images of the limit operator $H$ are given by
		\begin{align*}
		H_{f, \varphi}(x,y)
		&=
		e^{-\varphi(y)}
		\left[
		V_{f'(x)}(y)
		+
		B_{f'(x)}
		\right]
		e^{\varphi(\cdot)}(y)
		=:
		H_{\varphi}(f'(x),y),
		\end{align*}
		where~$V_p(y) := \frac{1}{2} p^2 - p \, \psi'(y)$ and~$B_p := (p - \psi'(y)) \partial_y + \frac{1}{2} \partial_{yy}$, for $p \in \mathbb{R}$.
		\smallskip
		
		\ref{MM:item:T2}: Fix~$p\in\mathbb{R}$. Finding a function $\varphi(\cdot)$ such that~$H_{\varphi}(p,y)$ becomes constant as a function of~$y$ is equivalent to finding a constant~$\mathcal{H}(p)$ such that on~$\mathbb{T}$,
		\begin{equation*}
		\left[ V_p + B_p \right] e^{\varphi} = \mathcal{H}(p) e^{\varphi}.
		\end{equation*}
		This is a principal-eigenvalue problem, where the constant~$\mathcal{H}(p)$ corresponds to the principal eigenvalue. 
		We come back to principal eigenvalues when considering the results about molecular motor models in Section~\ref{subsection:results:LDP_in_MM} and their proofs in Section~\ref{section:LDP_in_MM}. In Section~\ref{appendix:prinipal_ev}, we further outline to what extend the principal-eigenvalue problems that we encounter in this chapter are solved in the literature.\qed
\smallskip

This example hints at a more general structure comprising all molecular-motor models that we consider in this chapter. The different models are specified by the choice of the $\varepsilon$-scaling, the state space $E^X_\varepsilon$, the spatial dynamics defined by the generators $L^{i}_\varepsilon$, and the reaction rates $r_{ij}^\varepsilon(\cdot)$.
However, the proofs of large deviation principles are independent from these choices; they all follow Theorem~\ref{thm:results:LDP_switching_MP}. The model-specific contribution is only to determine in which setting~\ref{MM:item:T1} and~\ref{MM:item:T2} have to be verified.
	\subsection{Action-integral representation of the rate function}
	\label{subsection:results:action_integral_rep}
	In this section, our main goal is to give a feasible condition under which the rate function of Theorem~\ref{thm:results:LDP_switching_MP} is of action-integral form.
	We say that a rate function~$\mathcal{I}:D_{\mathbb{T}^d}[0,\infty)\to[0,\infty]$ is of action-integral form if there is a convex map~$\mathcal{L}:\mathbb{R}^d\to[0,\infty]$ with which
	\begin{equation*}
	\mathcal{I}(x) =
	\begin{cases}
	\mathcal{I}_0(x(0))+
	\int_0^\infty\mathcal{L}
	\left(\partial_t x(t)\right) \, \dd t 
	&\quad 
	\text{if } x\in \cA\cC([0,\infty); \mathbb{T}^d),\\
	+\infty & 
	\quad
	\text{otherwise}.
	\end{cases}
	\end{equation*}
	\begin{theorem} Consider the setting of Theorem~\ref{thm:results:LDP_switching_MP}. For~$p\in\mathbb{R}^d$, let~$\mathcal{H}(p)$ be the constant obtained in~\ref{MM:item:T2} of Theorem \ref{thm:results:LDP_switching_MP}.
	Suppose the following:
		\begin{itemize}
		\item[(T3)] The map $p \mapsto \mathcal{H}(p)$ is convex and $\mathcal{H}(0) = 0$.
		\end{itemize}
		Then the rate function of Theorem~\ref{thm:results:LDP_switching_MP} is of action-integral form with Lagrangian as the Legendre-Fenchel transform of $\mathcal{H}$, that is~$\mathcal{L}(v)=\sup_{p\in\mathbb{R}^d}\left[p\cdot v - \mathcal{H}(p)\right]$.
	\label{thm:results:action_integral_representation}
	\end{theorem}
	We give the proof in Section~\ref{section:action_integral}. The argument is based on the general strategy of~\cite[Chapter~8]{FengKurtz2006}. The idea of how to obtain such representations is also outlined in Section~\ref{BG:sec:semigroup-flow-HJ-eq} of Chapter~2 in this thesis.
	\subsection{Large deviations for models of molecular motors}
	\label{subsection:results:LDP_in_MM}
	In this section we formulate large deviation theorems for stochastic processes motivated by molecular motors. The proofs are given in Section~\ref{section:LDP_in_MM}. All proofs are based on verifying the conditions of Theorems~\ref{thm:results:LDP_switching_MP} and~\ref{thm:results:action_integral_representation} above. We first define the continuous model---for a motivation, in particular of the $\varepsilon$-scaling, we refer to Section~\ref{subsec:model-of-molecular-motor}.
	\begin{definition}[Continuous model]
		The pair~$(X^\varepsilon, I^\varepsilon)$ is a switching Markov process in~$E_\varepsilon = \mathbb{T}^d \times \{1,\dots,J\}$ with generator~$L_\varepsilon$ defined by
		\begin{multline}
		L_\varepsilon f(x,i) := b^i(x/\varepsilon) \cdot \nabla_x f(\cdot,i) (x) + \frac{\varepsilon}{2} \Delta_x f(\cdot,i) (x)\\
		+ \sum_{j = 1}^J \frac{1}{\varepsilon} \gamma(\varepsilon) r_{ij}(x/\varepsilon)
		\left[
		f(x,j) - f(x,i)
		\right],
		\label{eq:intro:L_varepsilon_cont_MM}
		\end{multline}
		where~$\gamma(\varepsilon)>0$,~$r_{ij}(\cdot) \in C^\infty(\mathbb{T}^d; [0,\infty))$, and~$b^i(\cdot) \in C^\infty(\mathbb{T}^d)$. This is an example of a switching Markov process with generators~$L_\varepsilon^{i}$ defined on the core~$C^2(\mathbb{T}^d)$ by
			\begin{equation*}
			L^{i}_\varepsilon g(x)
			:=
			b^i(x/\varepsilon) \cdot \nabla g(x)
			+
			\varepsilon \,\frac{1}{2} \Delta g(x),
			\end{equation*}
		and rates~$r_{ij}^\varepsilon(x)=(\gamma(\varepsilon)/\varepsilon)r_{ij}(x/\varepsilon)$.
		The domain of~$L_\varepsilon$ is the set given by~$\mathcal{D}(L_\varepsilon) = \{ f(x,i)\,:\, f(\cdot,i) \in \mathcal{D}(L_\varepsilon^{i})\}$.\qed
		\label{def:intro:cont_MM}
	\end{definition}
	\begin{definition}\label{MM:def:irreducible-matrix}
	Let~$J\in\mathbb{N}$. We call a matrix~$A\in\mathbb{R}^{J\times J}$ \emph{irreducible} if there is no decomposition of~$\{1,\dots,J\}$ into two disjoint sets $\mathcal{J}_1$ and $\mathcal{J}_2$ such that $A_{ij} = 0$ whenever $i \in \mathcal{J}_1$ and $j \in \mathcal{J}_2$.\qed
	\end{definition}
	\begin{theorem}[Continuous model, limit I]
	Let $(X^\varepsilon_t, I^\varepsilon_t)$ be the Markov process of Definition \ref{def:intro:cont_MM} with~$\gamma\equiv 1$. Assume that the matrix~$R$ with entries $R_{ij} = \sup_{y \in \mathbb{T}^d}r_{ij}(y)$ is irreducible. Suppose furthermore that
	the family of initial conditions $X^\varepsilon(0)$ satisfies a large deviation principle in $\mathbb{T}^d$ with rate function $\mathcal{I}_0 :\mathbb{T}^d \rightarrow [0,\infty]$.
	
	Then 
	the family of stochastic processes $\{X^\varepsilon\}_{\varepsilon > 0}$ satisfies a large deviation principle in $C_{\mathbb{T}^d}[0,\infty)$ with rate function of action-integral form.
	The Hamiltonian~$\mathcal{H}(p)$ is the principal eigenvalue of an associated cell problem described in~\eqref{eq:LDP_MM:contI:cell_problem} of Lemma~\ref{lemma:LDP_MM:contI:principal_eigenvalue}.
	%
	\label{thm:results:LDP_cont_MM_Limit_I}
	\end{theorem}
	The example of Section~\ref{subsec:model-of-molecular-motor} corresponds to~$d = 1$,~$J = 2$ and~$b^i=-\psi'(\cdot,i)$. The irreducibility condition is imposed to solve the principal-eigenvalue problem that we obtain, and is inspired by what Guido Sweers assumes to solve a coupled system of elliptic PDE's~\cite{Sweers92}. 
	\smallskip
	
	The parameter~$\gamma(\varepsilon)$ models an additional time-scale separation of the components. For large~$\gamma$, the spatial component is effectively driven by potentials averaged over the stationary measure of the fast configurational component. The following theorem shows that if~$\gamma(\varepsilon)\to\infty$, then the large deviation principle is governed by an averaged Hamiltonian.
	\begin{theorem}[Continuous model, limit II]
	Let $(X^\varepsilon_t, I^\varepsilon_t)$ be the Markov process of Definition~\ref{def:intro:cont_MM}, with~$\gamma(\varepsilon) \rightarrow \infty$ as~$\varepsilon \rightarrow 0$. Assume that for every~$y\in\mathbb{T}^d$, the matrix~$R(y)$ with entries~$R(y)_{ij}=r_{ij}(y)$ is irreducible.
	Suppose furthermore that the family of random variables $\{X^\varepsilon(0)\}_{\varepsilon > 0}$ satisfies a large deviation principle in $\mathbb{T}^d$ with rate function $\mathcal{I}_0 :\mathbb{T}^d \rightarrow [0,\infty]$. 
	
	Then~$\{X^\varepsilon\}_{\varepsilon > 0}$ satisfies a large deviation principle in $C_{\mathbb{T}^d}[0,\infty)$ with rate function of action-integral form.
	The Hamiltonian $\overline{\mathcal{H}}(p)$ is the principal eigenvalue of an associated averaged cell problem described in Lemma~\ref{lemma:LDP_MM:contII:principal_eigenvalue}.
	\label{thm:results:LDP_cont_MM_Limit_II}
	\end{theorem}
In the discrete model, the spatial component is not a drift-diffusion process, but a jump process on a discrete periodic lattice. We consider only nearest-neighbor jumps. 
We use the integer $n$ as the scaling parameter. For~$\ell\in\mathbb{N}$, we denote by~$\mathbb{T}_{\ell, n}$ the discrete one-dimensional flat torus of length $\ell$, lattice spacing $1/n$ and with $n \cdot \ell$ points. As a set,~$\mathbb{T}_{\ell, n} \simeq \{0, 1/n, \dots, \ell-1/n\}$ with periodic boundary.
\begin{definition}[Discrete model]
	The pair~$(X^n, I^n)$ be a switching Markov process in $E_n = \mathbb{T}_{\ell, n} \times \{1,\dots,J\}$ with generator~$L_n$ defined on~$\mathcal{D}(L_n) = C(E_n)$,
	\begin{multline}\label{eq:intro:L_N_discr_MM}
	L_nf(x,i) = 	n r_+^i(nx) \left[ f(x + 1/n, i) - f(x,i) \right] + n r_-^i(nx) \left[ f(x - 1/n, i) - f(x,i) \right]	\\
	+ \sum_{j = 1}^J n \gamma(n) r_{ij}(nx) \left[ f(x,j) - f(x,i) \right],
	\end{multline}
	where~$r_{ij}(\cdot) : \mathbb{T}_{\ell, 1} \rightarrow [0,\infty)$,~$r_\pm^i : \mathbb{T}_{\ell, 1} \rightarrow (0,\infty)$, and~$\gamma(n)>0$.
	To connect with our definition of switching Markov processes, the generators $L_n^i$ are
	\begin{align*}
	L^i_n g(x):=n r_+^i(nx) \left[ g(x + 1/n) - g(x) \right]+n r_-^i(nx) \left[ g(x - 1/n) - g(x) \right].\tag*\qed
	\end{align*}
	\label{def:intro:discr_MM}
\end{definition}
The discrete lattice $\mathbb{T}_{\ell, n}$ covers the continuous torus $\mathbb{T}_\ell = \mathbb{R}/({\ell \cdot \mathbb{Z}})$ in the limit~$n\to\infty$.
	\begin{theorem}[Discrete model, limit I]
	Let $(X^n_t,I^n_t)$ be the Markov process from Definition \ref{def:intro:discr_MM}, with $\gamma\equiv 1$. Suppose that
	the matrix~$R$ with entries defined by~$R_{ij} = \sup_{k \in \mathbb{T}_{\ell, 1}}r_{ij}(k)$ is irreducible.
	Suppose furthermore that~$\{X^n(0)\}_{n\in\mathbb{N}}$ satisfies a large deviation principle in $\mathbb{T}_\ell$ with rate function $\mathcal{I}_0 : \mathbb{T}_\ell \rightarrow [0,\infty]$. 
	
	Then~$\{X^n\}_{n \in \mathbb{N}}$ satisfies a large deviation principle in $D_{\mathbb{T}_\ell}[0,\infty)$ with rate function $\mathcal{I} : D_{\mathbb{T}_\ell}[0,\infty) \rightarrow [0,\infty]$ of action-integral form. 
	The Hamiltonian~$\mathcal{H}(p)$ is the principal eigenvalue of a cell problem described in Lemma~\ref{lemma:LDP_MM:discrI:principal_ev}.
	\label{thm:results:LDP_discr_MM_Limit_I}
	\end{theorem}
	If~$\gamma(n)$ is large, the spatial component~$X^n_t$ is driven by the average jump rates that result from averaging over the stationary distribution of the configurational component~$I^n_t$. If~$\gamma(n)\to\infty$, large deviations are characterized by an averaged Hamiltonian.
	\begin{theorem}[Discrete model, limit II]
	Let $(X^n_t,I^n_t)$ be the Markov process from Definition \ref{def:intro:discr_MM}, with $\gamma(n) \rightarrow \infty$ as $n \rightarrow \infty$. In addition to the assumptions of Theorem \ref{thm:results:LDP_discr_MM_Limit_I}, suppose that for each $ k \in \mathbb{T}_{\ell,1} \simeq \{0,1, \dots, \ell-1\}$, there exists a stationary measure $\mu_k \in \mathcal{P}(\{1,\dots,J\})$ for the jump process on $\{1,\dots,J\}$ with frozen jump rates $r_{ij}(k)$.
	Suppose furthermore that~$\{X^n(0)\}_{n \in \mathbb{N}}$ satisfies a large deviation principle in $\mathbb{T}_\ell$ with rate function $\mathcal{I}_0 : \mathbb{T}_\ell \rightarrow [0,\infty]$. 
	
	Then 
	$\{X^n_t|_{t \geq 0}\}_{n \in \mathbb{N}}$ satisfies a large deviation principle in $D_{\mathbb{T}_\ell}[0,\infty)$ with rate function of action-integral form.
	The Hamiltonian is the principal eigenvalue of an averaged cell problem described in Lemma~\ref{lemma:LDP_MM:discrII:principal_ev}.
	\label{thm:results:LDP_discr_MM_Limit_II}
	\end{theorem}
	\subsection{Detailed balance implies symmetric Hamiltonians}
	\label{MM:subsec:det-bal-implies-symmetric-H}
	In this section we show how the large deviation principles can be used to analyse which sets of potentials and rates induce transport on macroscopic scales. To that end, we consider the generator defined in~\eqref{eq:intro:L_varepsilon_cont_MM}, with~$b^i(y)=-\nabla_y\psi_i(y)$. We say that the set of potentials and rates~$\{r_{ij},\psi_i\}$ satisfies \emph{detailed balance} if for all~$i,j\in\{1,\dots,J\}$ and~$y\in\mathbb{T}^d$, we have
	\begin{align}\label{MM:eq:detailed-balance}
	r_{ij}(y) e^{-2\psi_i(y)}=r_{ji}(y) e^{-2\psi_j(y)}.
	\end{align}
	Let us motivate detailed balance. The Fokker-Planck equations of~$(X^\varepsilon,I^\varepsilon)$ are
	\begin{equation*}
	\partial_t \rho_\varepsilon^i = \varepsilon\frac{1}{2}\Delta\rho_\varepsilon^i+\mathrm{div}_x\left[\rho_\varepsilon^i\nabla_y\psi_i(x/\varepsilon)\right] +\frac{1}{\varepsilon} \sum_{j=1}^Jr_{ji}(x/\varepsilon)\rho_\varepsilon^j,\qquad i=1,\dots,J,
	\end{equation*}
	with~$r_{ii}:=-\sum_{j\neq i}r_{ij}$. 
	In general, the stationary measure~$\pi_\varepsilon\in\mathcal{P}(\mathbb{T}^d\times\{1,\dots,J\})$ satisfying~$\partial_t\pi_\varepsilon=0$ is not known explicitly.
	Define the total flux~$J$ by
	\begin{equation*}
	J = \sum_iJ_i\,,\qquad J_i(t,x):=-\varepsilon \frac{1}{2}\nabla_x \rho_\varepsilon^i(t,x) - \rho_\varepsilon^i(t,x) \nabla_y\psi_i(x/\varepsilon).
	\end{equation*}
	Since~$0=\partial_t\sum_i\pi_\varepsilon^i=-\mathrm{div}_x\sum_iJ_i$, the total flux is constant if the system is stationary. Detailed balance is achieved if in addition, 1) each~$J_i$ is constant and 2) the flux between any two configurations~$i$ and~$j$ is balanced. If~$J_i$ is constant, then the stationary component~$\pi_\varepsilon^i$ is a Boltzmann distribution, that means there are constants~$C_i$ such that~$\pi_\varepsilon^i(\dd x)=C_ie^{-2\psi_i(x/\varepsilon)}\dd x$. The constants are inessential and can be absorbed into the potentials (since constant shifts of the potentials do not affect the dynamics), with which we find the stationary measure
	\begin{equation*}
	\pi_\varepsilon(\dd x\dd i) = Z^{-1}\,e^{-2\psi_i(x/\varepsilon)}\,\dd x\dd i,\quad Z = \sum_i\int e^{-2\psi_i(x/\varepsilon)}\,\dd x.
	\end{equation*}
	Regarding 2), we find from~$\partial_t\pi_\varepsilon=0$ by summing over~$i$ that
	\begin{equation*}
	0=\sum_j \left(r_{ji}\pi_\varepsilon^j-r_{ij}\pi_\varepsilon^i\right).
	\end{equation*}
	The detailed-balance condition~\eqref{MM:eq:detailed-balance} requires that each term in the summation vanishes. 
	This motivates the notion of \emph{detailed} balance---if the system is stationary, then the flow from~$i$ to~$j$ is equal to the flow from~$j$ to~$i$.
	\smallskip
	
	Detailed balance implies time-reversibility of the process~$(X^\varepsilon,I^\varepsilon)$ in the sense of Definition~\ref{intro:def:reversibility}. This follows from the fact that symmetry of the generator is equivalent to time-reversibility~\cite[Chapter~II, Proposition~5.3]{Liggett2004}. A calculation shows that
	the generator~$L_\varepsilon$ is symmetric with respect to the stationary measure~$\pi_\varepsilon$; for all~$f,g\in\mathcal{D}(L_\varepsilon)$, we have~$\langle L_\varepsilon f,g\rangle_{\pi_\varepsilon} = \langle f,L_\varepsilon g\rangle_{\pi_\varepsilon}$.
	\begin{theorem}[Detailed balance implies a symmetric Hamiltonian]
	Let~$(X^\varepsilon_t, I^\varepsilon_t)$ be the stochastic process of Definition~\ref{def:intro:cont_MM} with~$b^i=-\nabla_y\psi_i$, where the~$\psi_i$ are smooth potentials. Suppose that the assumptions of Theorem~\ref{thm:results:LDP_cont_MM_Limit_I} and the detailed-balance condition~\eqref{MM:eq:detailed-balance} are satisfied. Then the Hamiltonian satisfies~$\mathcal{H}(p) = \mathcal{H}(-p)$ for all~$p\in\mathbb{R}^d$.
	\label{thm:results:detailed_balance_limit_I}
	\end{theorem}
	\begin{theorem}[Separation of time scales implies a symmetric Hamiltonian]
		Let the stochastic process~$(X^\varepsilon_t, I^\varepsilon_t)$ of Definition~\ref{def:intro:cont_MM} with~$b^i=-\nabla\psi^i$ satisfy the assumptions from Theorem \ref{thm:results:LDP_cont_MM_Limit_II}. Suppose in addition that the rates $r_{ij}(\cdot)$ are constant on $\mathbb{T}^d$. Then $\overline{\mathcal{H}}(p) = \overline{\mathcal{H}}(-p)$, where $\overline{\mathcal{H}}(p)$ is the Hamiltonian from Theorem \ref{thm:results:LDP_cont_MM_Limit_II}.
	\label{thm:results:symmetry_limit_II}
	\end{theorem}
	In both situations, the macroscopic velocity given by $v = \partial_p\mathcal{H}(0)$ vanishes due to the symmetry of the Hamiltonians. Theorem~\ref{thm:results:symmetry_limit_II} confirms the numerical results of Weng, Peskin and Elston~\cite[Section~4.3]{WangPeskinElston2003}. Since the proofs of Theorems~\ref{thm:results:detailed_balance_limit_I} and~\ref{thm:results:symmetry_limit_II} are solely based on a suitable formula for~$\mathcal{H}(p)$, we give them here---the formulas of~$\mathcal{H}(p)$ are proven in Section~\ref{subsection:detailed_balance}. Since the derivation of these formulas is similar, we only give the argument for the more involved case of Theorem~\ref{thm:results:detailed_balance_limit_I}.
	\begin{proof}[Proof of Theorem~\ref{thm:results:detailed_balance_limit_I}]
	We prove in Proposition~\ref{prop:results:detailed_balance_limit_I} that under the detailed-balance condition, the principal eigenvalue~$\mathcal{H}(p)$ is given by 
	$$
		\mathcal{H}(p)
	=
		\sup_{\mu \in \mathbf{P}} \left[ K_p(\mu) - \mathcal{R}(\mu) \right],
	$$
	where $\mathbf{P} \subset \mathcal{P}(E^\prime)$ is a subset of probability measures on $E^\prime = \mathbb{T}^d \times \{1, \dots, J\}$ specified in Proposition~\ref{prop:results:detailed_balance_limit_I}, $\mathcal{R}(\mu)$ is the relative Fisher information specified in~\eqref{MM:eq:relative-Fisher-information}, and $K_p(\mu)$ is given by
	\begin{multline*}
		K_p(\mu)= 
		\inf_{\phi}\bigg\{
		\sum_{i = 1}^J \int_{\mathbb{T}^d}
		\left(
		\frac{1}{2}|\nabla \phi_i(x) + p|^2 
	-
		\sum_{j = 1}^J r_{ij}(x)
		\right) \,
		\dd \mu_i(x)
		\\+
		\sum_{i, j = 1}^J
		\int_{\mathbb{T}^d} \pi_{ij}(x) \sqrt{\overline{\mu}_i(x) \overline{\mu}_j(x)}
		e^{\psi_j(x) + \psi_i(x)} \cosh{(\phi(x,j) - \phi(x,i))}
		\,
		\dd x
		\bigg\},
	\end{multline*}
	where $\pi_{ij}(x) = r_{ij}(x) e^{-2\psi_i(x)}$, the infimum is taken over vectors of functions $\phi_i = \phi(\cdot,i) \in C^2(\mathbb{T}^d)$, and $\dd \mu_i(x) = \overline{\mu}_i(x) \dd x$.
	\smallskip
	
	Let~$\mu \in \mathbf{P}$. We show that~$K_p(\mu) = K_{-p}(\mu)$, which implies~$\mathcal{H}(p) = \mathcal{H}(-p)$. Since~$\cosh(\cdot)$ is symmetric, the sum in which the $\cosh(\cdot)$ terms appear is invariant under transforming as $\phi \to (-\phi)$, in the sense that for
	$$
		C(\phi):=
		\sum_{i, j = 1}^J \int_{\mathbb{T}^d} \pi_{ij}(x)
		\sqrt{\overline{\mu}_i(x) \overline{\mu}_j(x)}
		e^{\psi_j(x) + \psi_i(x)} \cosh{(\phi(x,j) - \phi(x,i))}
		\,
		\dd x,
	$$
	we have $C(\phi) = C(-\phi)$. Hence the bijective transformation $\phi \to (-\phi)$ implies the claimed symmetry~$K_p(\mu) = K_{-p}(\mu)$.
	\end{proof}
	\begin{proof}[Sketch of proof of Theorem~\ref{thm:results:symmetry_limit_II}]
	Under the detailed-balance condition, one can prove that the principal eigenvalue $\overline{\mathcal{H}}(p)$ is given by
	 $$
	 	\overline{\mathcal{H}}(p)
	 =
	 	\sup_{\mu \in \mathbf{P}}
	 	\left[
	 		K_p(\mu) - \mathcal{R}(\mu)
	 	\right],\quad
	 	K_p(\mu) = \inf_{\varphi} \frac{1}{2} \int_{\mathbb{T}^d} | \nabla \varphi + p|^2 \, \dd \mu,
	 $$
	 where $\mathbf{P} \subset \mathcal{P}(\mathbb{T}^d)$ is a subset of the probability measures on $\mathbb{T}^d$,
	 \[
	 \mathbf{P}
	 =
	 \left\{ \mu \in \mathcal{P}(\mathbb{T}^d) \, : \, \mu \ll \dd x,
	 \, \dd \mu = \overline{\mu} \dd x, \,
	 \text{ and } \nabla 
	 \left(
	 \log \overline{\mu} 
	 \right) 
	 \in L^2_{\mu}(\mathbb{T}^d) \right\}.
	 \]
	 The map~$\mathcal{R}$ is the relative Fisher information; with the stationary measure~$\nu$ of the jump process on~$\{1,\dots,J\}$ with rates~$r_{ij}$, we have
	 \begin{equation*}
	 \mathcal{R}(\mu) = \frac{1}{8}\int_{\mathbb{T}^d}\left|\nabla\log\left(\frac{\overline{\mu}}{e^{-2\overline{\psi}}}\right)\right|^2\,\dd\mu,\quad \overline{\psi}(x)=\sum_i \nu_i\,\psi_i(x).
	 \end{equation*}
	 We have~$K_p(\mu) = K_{-p}(\mu)$, since the bijective transformation $\varphi \to (-\varphi)$ leaves the infimum invariant. This implies~$\overline{\mathcal{H}}(p) = \overline{\mathcal{H}}(-p)$.
	\end{proof}
	With a similar analysis, we can study the behaviour of molecular motors under external forces. Let $(X^\varepsilon_t, I^\varepsilon_t)$ be the stochastic process from Theorem \ref{thm:results:LDP_cont_MM_Limit_I} in dimension $d = 1$ with drift~$b^i(y)=F - \psi'(y,i)$, where~$F$ is a constant modelling an external force and~$\psi \in C^\infty(\mathbb{T})$ is a smooth periodic potential. The process~$(X^\varepsilon_t, I^\varepsilon_t)$ is $\mathbb{T} \times \{1,\dots,J\}$-valued and satisfies
	$$
		\dd X^\varepsilon_t = (F - \psi'((X^\varepsilon_t / \varepsilon), I^\varepsilon_t))\, \dd t + \sqrt{\varepsilon} \,\dd B_t,
	$$
	where $I^\varepsilon_t$ a jump process on $\{1,\dots,J\}$ with jump rates~$\frac{1}{\varepsilon} r_{ij}\left(x/\varepsilon\right)$.
	In this case, the Hamiltonian is given by
	\begin{multline*}
		\mathcal{H}(p)
	=
		\frac{1}{2} |p + F|^2 - \frac{1}{2} F^2	\\
	+
		\inf_{\varphi} \sup_{(y,i)} 
		\bigg[
		\frac{1}{2} \varphi''(y,i) + (p + F - \psi'(y,i)) \varphi'(y,i)
		- \psi'(y,i) p
		\\+
		\sum_{i = 1}^J r_{ij}(y) \left( e^{\varphi(y,j) - \varphi(y,i)} - 1 \right)
		\bigg].
	\end{multline*}
	Using detailed balance, this Hamiltonian is symmetric around $(-F)$; one can deal with the variational terms similar as above. Since~$\mathcal{H}(0) = 0$ and~$\mathcal{H}(p)$ is strictly convex, this means that the model predicts a positive force-velocity feedback, since~$F > 0$ implies~$\partial_p\mathcal{H}(0) > 0$, and~$F < 0$ implies~$\partial_p\mathcal{H}(0) < 0$. Establishing a similar result for systems not satisfying detailed balance would be interesting.
	\section{Proof of large deviations of switching processes}
	\label{section:LDP_switching_MP}
	In this section, we prove~Theorem~\ref{thm:results:LDP_switching_MP} (large deviations of switching processes) and Theorem~\ref{thm:results:action_integral_representation} (action-integral representation of the rate function). To do so, we exploit the connection of pathwise large deviations to Hamilton-Jacobi equations~\cite{FengKurtz2006}. In Section~\ref{appendix:LDP_via_CP}, we adapt of~\cite[Theorem~7.18]{FengKurtz2006} to our compact setting. Then we prove Theorem~\ref{thm:results:LDP_switching_MP} in Section~\ref{MM:sec:proof-LDP-switching-MP}, and Theorem~\ref{thm:results:action_integral_representation} in Section~\ref{section:action_integral}.
	\subsection{Pathwise large deviations via comparison principle}
	\label{appendix:LDP_via_CP}
	In the following definitions,~$E$ and~$E^\prime$ are compact metric spaces. In the examples of our note, the space $E$ corresponds to $\mathbb{T}^d$, and $E^\prime$ to the space of upscaled variables. In this section,~$\text{BUSC}(E),\text{BLSC}(E)$ denote the sets of bounded and upper (lower) semicontinuous functions on~$E$, and $\text{BLSC}(E)$ for the bounded and lower semicontinuous functions on~$E$. We adapt~\cite[Definition~7.1]{FengKurtz2006} to the compact setting.
	\begin{definition}[Viscosity solutions]
		Let $H \subseteq C(E) \times C(E\times E^\prime)$ be a multivalued operator with domain~$\mathcal{D}(H) \subseteq C(E)$. Let~$h\in C(E)$ and~$\tau>0$.
		\begin{itemize}
			\item[i)] $u_1\in \text{BUSC}(E)$ is a viscosity subsolution of $(1 - \tau H) u = h$ if for all $(f,g)\in H$ there exists a point $(x,z^\prime)\in E\times E^\prime$ such that
			\[
			(u_1-f)(x)=\sup(u_1-f)\quad\text{and}\quad u_1(x)-\tau g(x,z^\prime)-h(x)\leq 0.
			\]
			\item[ii)] $u_2\in \text{BLSC}(E)$ is a viscosity supersolution of $(1 - \tau H) u = h$ if for all $(f,g)\in H$ there exists a point $(x,z^\prime)\in E\times E^\prime$ such that
			\[
			(f-u_2)(x)=\sup(f-u_2)\quad\text{and}\quad u_2(x)-\tau g(x, z^\prime)-h(x)\geq 0.
			\]
			\item[iii)] $u_1\in\text{BUSC}(E)$ is a strong viscosity subsolution of $(1 - \tau H) u = h$ if for all $(f,g) \in H$ and $x\in E$, whenever
			\[
			(u_1-f)(x)=\sup(u_1-f),
			\]
			then there exists a $z^\prime\in E^\prime$ such that
			\[
			u_1(x)-\tau g(x,z^\prime)-h(x)\leq 0.
			\]
			Similarly for strong supersolutions.
		\end{itemize}
		A function $u\in C(E)$ is called a viscosity solution of $(1 - \tau H) u = h$ if it is both a viscosity sub- and supersolution.
		\label{def:appendix:viscosity_solutions_multivalued_op}
	\end{definition}
	\begin{definition}[Comparison Principle]\quad\\
	We say that the \emph{comparison principle} holds for viscosity sub- and supersolutions of~$(1-\tau H)u=h$ if for any viscosity subsolution~$u_1$ and viscosity supersolution~$u_2$, we have~$u_1\leq u_2$ on~$E$. 
		\label{def:appendix:CP_single_valued_operator}
	\end{definition}	
	In the following adaptation of~\cite[Theorem 7.18]{FengKurtz2006}, the compact Polish spaces $E_n$, $E$ and $E^\prime$ are related with continuous embeddings $\eta_n$ and $\eta_n^\prime$ by
	$$
	\begin{tikzcd}[row sep = tiny]
	&	E \times E^\prime \arrow[dd, "\mathrm{proj}_1"]	\\
	E_n \arrow[ur, "(\eta_n {,} \eta_n^\prime)"] \arrow[dr, "\eta_n"']	&	\\
	&	E	
	\end{tikzcd}
	$$
	such that for any $x \in E$, there exist $x_n \in E_n$ such that $\eta_n(x_n) \to x$ as $n \to \infty$.
	\begin{theorem}
		Let $L_n$ be the generator of an $E_n$-valued process $Y^n$, and let $H_n$ be the nonlinear generators defined by $H_n f = \frac{1}{n} e^{-nf} L_n e^{nf}$. Let the compact Polish spaces $E_n, E$ and $E'$ be related as in the above diagram. In addition, suppose:
		\begin{enumerate}[(i)]
			\item(Condition 7.9 of \cite{FengKurtz2006} on the state spaces)
			There exists an index set $Q$ and approximating state spaces $A_n^q \subseteq E_n$, $q \in Q$, such that the following holds:
			\begin{enumerate}[(a)]
				\item For $q_1, q_2 \in Q$, there exists $q_3 \in Q$ such that $A_n^{q_1} \cup A_n^{q_2} \subseteq A_n^{q_3}$.
				\item For each $x \in E$, there exists $q \in Q$ and $y_n \in A^q_n$ such that $\eta_n(y_n) \to x$ as $n \to \infty$.
				\item For each $q \in Q$, there exist compact sets $K_1^q \subseteq E$ and $K_2^q \subseteq E \times E^\prime$ such that
				\[
				\sup_{y \in A^q_n} \inf_{x \in K_1^q} d_E(\eta_n(y), x) \xrightarrow{n\to\infty} 0,
				\]
				and
				\[
				\sup_{y \in A^q_n} \inf_{(x,z) \in K_2^q} 
				\left[ d_E(\eta_n(y),x)) + d_{E'}(\eta_n'(y),z) \right] \xrightarrow{n\to\infty} 0.
				\]
				\item For each compact $K \subseteq E$, there exists $q \in Q$ such that
				$
				K \subseteq \liminf \eta_n(A_n^q).
				$
			\end{enumerate}
			\item 
			(Convergence Condition 7.11 of \cite{FengKurtz2006}) There exist $H_\dagger, H_{\ddagger} \subseteq C(E) \times C(E \times E^\prime)$ which are the limit of the $H_n$'s in the following sense:
			\begin{enumerate}[(a)]
				\item
				For each $(f,g) \in H_\dagger$, there exist $f_n \in \mathcal{D}(H_n)$ such that
				\[
				\sup_n \left(
				\sup_{x \in E_n} |f_n (x)| + \sup_{x \in E_n} |H_n f_n (x)|
				\right)
				< \infty,
				\]
				and for each $q\in Q$,~$\lim_{n \to \infty} \sup_{y \in A_n^q}|f_n(y) - f(\eta_n(y))| = 0.$
				%
				%
				Furthermore, for each $q \in Q$ and every sequence $y_n \in A_n^q$ such that $\eta_n (y_n) \to x \in E$ and $\eta'_n(y_n) \to z^\prime \in E^\prime$, we have~$\limsup_{n \to \infty} H_n f_n (y_n)
				\leq
				g(x,z^\prime)$.
				\item 
				For each $(f,g) \in H_\ddagger$, there exist $f_n \in \mathcal{D}(H_n)$ (not necessarily the same as above in (a)) such that
				\[
				\sup_n \left(
				\sup_{x \in E_n} |f_n (x)| + \sup_{x \in E_n} |H_n f_n (x)|
				\right)
				< \infty,
				\]
				and for each $q\in Q$,~$\lim_{n \to \infty} \sup_{y \in A_n^q}|f_n(y) - f(\eta_n(y))| = 0$.
				%
				Furthermore, for each $q \in Q$ and every sequence $y_n \in E_n$ such that $\eta_n (y_n) \to x \in E$ and $\eta'_n(y_n) \to z^\prime \in E^\prime$, we have~$\liminf_{n \to \infty} H_n f_n (y_n)
				\geq
				g(x,z^\prime)$.
				%
			\end{enumerate}
			\item
			(Comparison principle) For each $h \in C(E)$ and $\tau > 0$, the comparison principle holds for viscosity subsolutions of $(1 - \tau H_\dagger) u = h$ and viscosity supersolutions of $(1 - \tau H_\ddagger) u = h$.
		\end{enumerate}
		Let $X^n_t := \eta_n (Y^n_t)$ be the corresponding $E$-valued process.
		If $\{X^n(0)\}_{n \in \mathbb{N}}$ satisfies a large deviation principle in $E$ with rate function $\mathcal{I}_0 : E \to [0,\infty]$, then
		%
		$\{X^n\}_{n \in \mathbb{N}}$ 
		satisfies the large deviation principle with rate function $\mathcal{I} : C_E[0,\infty) \to [0,\infty]$ given as in~\eqref{eq:RF-with-semigroup-approach} and~\eqref{BG:eq:RF-1d-marginals} of Chapter~2.
	\label{thm:appendix:LDP_via_CP:Jin_LDP_thm}
	\end{theorem}
The formula for the rate function involves a limiting semigroup~$V(t)$, which we discuss in Chapter~2. We do not repeat its formula here, since we will not work with it.
	\subsection{Proof of large deviation principle}
	\label{MM:sec:proof-LDP-switching-MP}
	Here we prove Theorem~\ref{thm:results:LDP_switching_MP} by verifying the conditions of Theorem~\ref{thm:appendix:LDP_via_CP:Jin_LDP_thm}, which are convergence of nonlinear generators (Proposition~\ref{prop:LDP_switching_MP:convergence_condition_sufficient}) and the comparison principle (Proposition~\ref{prop:LDP_switching_MP:comparison_principle}). The rest of this section below the proof of Theorem~\ref{thm:results:LDP_switching_MP} is devoted to proving the propositions. We point out that the main challenge is to prove the comparison principle using~\ref{MM:item:T1} and~\ref{MM:item:T2}.
	\begin{proposition}
	In the setting of Theorem~\ref{thm:results:LDP_switching_MP}, Condition~(i) of Theorem~\ref{thm:appendix:LDP_via_CP:Jin_LDP_thm} is satisfied. Let $H \subseteq C^1(\mathbb{T}^d) \times C(\mathbb{T}^d \times E^\prime)$ be a multivalued operator satisfying~\ref{MM:item:T1}. Then~$H$ satisfies the convergence condition (ii) of Theorem \ref{thm:appendix:LDP_via_CP:Jin_LDP_thm}.
	\label{prop:LDP_switching_MP:convergence_condition_sufficient}
	\end{proposition}	
	\begin{proposition}
		In the setting of~Theorem \ref{thm:results:LDP_switching_MP}, let $H \subseteq C^1(\mathbb{T}^d) \times C(\mathbb{T}^d \times E^\prime)$ be a multivalued operator satisfying conditions~\ref{MM:item:T1} and~\ref{MM:item:T2}. Then for $\tau > 0$ and $h \in C(\mathbb{T}^d)$, the comparison principle holds for viscosity sub- and supersolutions of
		$
			(1 - \tau H) u =h.
		$
	\label{prop:LDP_switching_MP:comparison_principle}
	\end{proposition}
	\begin{proof}[Proof of Theorem \ref{thm:results:LDP_switching_MP}]
	By Proposition~\ref{prop:LDP_switching_MP:convergence_condition_sufficient}, conditions (i) and (ii) of Theorem~\ref{thm:appendix:LDP_via_CP:Jin_LDP_thm} hold with the single operator $H = H_\dagger = H_\ddagger$. By Proposition~\ref{prop:LDP_switching_MP:comparison_principle}, the comparison principle is satisfied for $(1 - \tau H)u = h$, and hence condition (iii) of Theorem~\ref{thm:appendix:LDP_via_CP:Jin_LDP_thm} holds with a single operator $H = H_\dagger = H_\ddagger$. Therefore the conditions of Theorem~\ref{thm:appendix:LDP_via_CP:Jin_LDP_thm} are satisfied, and the large deviation principle follows.
	\end{proof}
	\begin{proof}[Proof of Proposition \ref{prop:LDP_switching_MP:convergence_condition_sufficient}]
	Recall that with $E_\varepsilon = E_\varepsilon^X \times \{1,\dots,J\}$ and $\iota_\varepsilon : E^X_\varepsilon \to \mathbb{T}^d$ from Condition~\ref{condition:results:general_setting}, the state spaces are related as in the following diagram, in which $\eta_\varepsilon : E_\varepsilon \to \mathbb{T}^d$ is defined by $\eta_\varepsilon(x,i) = \iota_\varepsilon(x)$ and $\eta_\varepsilon' : E_\varepsilon \to E'$ is a continuous map,
	$$
		\begin{tikzcd}[row sep = tiny]
			&	\mathbb{T}^d \times E^\prime \arrow[dd, "\mathrm{proj}_1"]	\\
		E_\varepsilon \arrow[ur, "(\eta_\varepsilon {,} \eta_\varepsilon^\prime)"] \arrow[dr, "\eta_\varepsilon"']	&	\\
			&	\mathbb{T}^d	
		\end{tikzcd}
	$$
	In the notation of Theorem \ref{thm:appendix:LDP_via_CP:Jin_LDP_thm}, we have $E = \mathbb{T}^d$. For verifying the general condition (i) of Theorem \ref{thm:appendix:LDP_via_CP:Jin_LDP_thm} on the approximating state spaces $A_\varepsilon^q$, we take the singleton $Q = \{q\}$ and set $A_\varepsilon^q := E_\varepsilon$. Then part (a) holds, and parts (b) and~(d) are a consequence of Condition~\ref{condition:results:general_setting} on $E_\varepsilon$, which says that for any $x \in \mathbb{T}^d$, there exist $ x_\varepsilon \in E_\varepsilon^X$ such that $\iota_\varepsilon(x_\varepsilon) \to x$. Part (c) follows by taking the compact sets $K_1^q := \mathbb{T}^d$ and $K_2^q := \mathbb{T}^d \times E^\prime$.
	\smallskip 
	
	 We verify the convergence Condition (ii) of Theorem \ref{thm:appendix:LDP_via_CP:Jin_LDP_thm}. By~\ref{MM:item:T1}, part~\ref{MM:item:C2}, there exist $f_\varepsilon \in \mathcal{D}(H_\varepsilon)$ such that
		$$
			\| f\circ\eta_\varepsilon - f_\varepsilon \|_{L^\infty(E_\varepsilon)} \xrightarrow{\varepsilon\to 0} 0
			\quad\text{and}\quad
			\| H_{f,\varphi}\circ (\eta_\varepsilon, \eta_\varepsilon') - H_\varepsilon f_\varepsilon \|_{L^\infty(E_\varepsilon)} \xrightarrow{\varepsilon\to 0} 0.
		$$
	With these $f_\varepsilon$, both conditions (a) and (b) are simultaneously satisfied for the operator $H = H_\dagger = H_\ddagger$, where condition~\ref{MM:item:C1} guarantees that for any point $(x,z') \in \mathbb{T}^d \times E'$ there exist $y_\varepsilon \in E_\varepsilon$ such that both $\eta_\varepsilon(y_\varepsilon) \to x$ and $\eta_\varepsilon'(y_\varepsilon) \to z'$. The boundedness 
	$$
		\sup_{\varepsilon > 0} \left( \sup_{y \in E_\varepsilon}|f_\varepsilon(y)| + \sup_{y \in E_\varepsilon}|H_\varepsilon f_\varepsilon(y)|\right) < \infty
	$$
	follows the uniform-convergence condition~\ref{MM:item:C2} and compactness of~$E_\varepsilon$.
	\end{proof}
	For proving Proposition \ref{prop:LDP_switching_MP:comparison_principle}, we use two operators~$H_1,H_2$ that are derived from a multivalued limit $H$. Define~$H_1,H_2$ by
	\begin{equation*}
	H_1f(x)
	:=
	\inf_{\varphi} \sup_{z^\prime \in E^\prime}
	H_{f, \varphi}(x,z^\prime)
	\quad\text{and}\quad
	H_2f(x)
	:=
	\sup_{\varphi} \inf_{z^\prime \in E^\prime}
	H_{f, \varphi}(x,z^\prime),
	\end{equation*}
	with equal domains~$\mathcal{D}(H_1) = \mathcal{D}(H_2) := \mathcal{D}(H)$. 
	Since the images of~$H$ are of the form $H_{f,\varphi}(x,z^\prime) = H_{\varphi}(\nabla f(x),z^\prime)$, the operators $H_1$ and $H_2$ are as well of the form $H_1 f(x) = \mathcal{H}_1(\nabla f(x))$ and $H_2 f(x) = \mathcal{H}_2(\nabla f(x))$, with two maps $\mathcal{H}_1,\mathcal{H}_2 : \mathbb{R}^d \rightarrow \mathbb{R}$.
	%
	We prove Proposition \ref{prop:LDP_switching_MP:comparison_principle} with the following Lemmas.
	\begin{lemma}[Local operators admit strong solutions]
	Let 
	$
		H \subseteq C^1(\mathbb{T}^d) \times C(\mathbb{T}^d\times E^\prime) 
	$
	be a multivalued limit operator
	satisfying (T1) from Theorem \ref{thm:results:LDP_switching_MP}. Then for any $\tau > 0$ and $h \in C(\mathbb{T}^d)$, viscosity solutions of
	$
		(1 - \tau H) u =h
	$
	coincide with strong viscosity solutions in the sense of Definition \ref{def:appendix:viscosity_solutions_multivalued_op}.
	\label{lemma:LDP_switching_MP:local_op_strong_sol}
	\end{lemma}
	\begin{lemma}[$H_1$ and $H_2$ are viscosity extensions]
		Let $H$ be a multivalued operator satisfying~\ref{MM:item:T1} and~\ref{MM:item:T2} of Theorem~\ref{thm:results:LDP_switching_MP}. For all $h \in C(\mathbb{T}^d)$ and $\tau > 0$, strong viscosity subsolutions $u_1$ of
		$
			(1 - \tau H) u =h
		$
		are strong viscosity subsolutions of
		$
			(1 - \tau H_1) u =h,
		$
		and strong viscosity supersolutions $u_2$ of
		$
			(1 - \tau H) u = h
		$
		are strong viscosity supersolutions of
		$
			(1 - \tau H_2) u =h.
		$
	\label{lemma:LDP_switching_MP:H_to_H1-H2}
	\end{lemma}
	\begin{lemma}[$H_1$ and $H_2$ are ordered]
		Let $H$ be a multivalued operator satisfying~\ref{MM:item:T1} and~\ref{MM:item:T2} of Theorem~\ref{thm:results:LDP_switching_MP}. Then~$\mathcal{H}_1(p)\leq \mathcal{H}_2(p)$ for all~$p\in\mathbb{R}^d$..
	\label{lemma:LDP_switching_MP:H1_leq_H2}
	\end{lemma}
	\begin{proof}[Proof of Proposition \ref{prop:LDP_switching_MP:comparison_principle}]
		Let $u_1$ be a subsolution and $u_2$ be a supersolution of the equation $(1 - \tau H)u = h$. By Lemma~\ref{lemma:LDP_switching_MP:local_op_strong_sol}, $u_1$ is a strong subsolution and $u_2$ a strong supersolution of $(1 - \tau H) u = h$, respectively. By Lemma~\ref{lemma:LDP_switching_MP:H_to_H1-H2}, $u_1$ is a strong subsolution of $(1 - \tau H_1) u = h$, and $u_2$ is a strong supersolution of $H_2$.
		\smallskip 
		
		With that, we establish below the inequality
		\begin{equation}
			\max_{\mathbb{T}^d} (u_1 - u_2)
		\leq
			\tau \left[ \mathcal{H}_1(p_\delta) - \mathcal{H}_2(p_\delta) \right]
		+
			h(x_\delta) - h(x_\delta^\prime),
		\label{eq:LDP_switching_MP:CP_visc_ineq}
		\end{equation}
		with some $x_\delta, x_\delta^\prime \in \mathbb{T}^d$ such that $\text{dist}(x_\delta,x_\delta^\prime) \rightarrow 0$ as $\delta \rightarrow 0$, and certain $p_\delta \in \mathbb{R}^d$. Then using that $h\in C(\mathbb{T}^d)$ is uniformly continuous since~$\mathbb{T}^d$ is compact, and that $\mathcal{H}_1(p_\delta) \leq \mathcal{H}_2(p_\delta)$ by Lemma~\ref{lemma:LDP_switching_MP:H1_leq_H2}, we can further estimate as
		$$
			\max_{\mathbb{T}^d} (u_1 - u_2)
		\leq
			h(x_\delta) - h(x_\delta^\prime)
		\leq
			\omega_h(\text{dist}(x_\delta,x_\delta^\prime)),
		$$
		where $\omega_h : [0,\infty) \rightarrow [0,\infty)$ is a modulus of continuity satisfying $\omega_h(r_\delta) \rightarrow 0$ for $r_\delta \rightarrow 0$. Then $(u_1 - u_2) \leq 0$ follows by taking the limit $\delta \rightarrow 0$.
		\smallskip
		
		We are left with proving \eqref{eq:LDP_switching_MP:CP_visc_ineq}. The line of argument is similar to the one outlined at the end of Section~\ref{subsec:LDP-via-visc-sol} of Chapter~2. Define $\Phi_\delta : \mathbb{T}^d \times \mathbb{T}^d \rightarrow \mathbb{R}$ by
		$$
			\Phi_\delta(x,x^\prime)
		:=
			u_1(x) - u_2(x^\prime) - \frac{\Psi(x,x^\prime)}{2 \delta},
		$$
		where
		\begin{equation}
			\Psi(x,x^\prime)
		:=
			\sum_{j=1}^d \sin^2\left( \pi(x_j - x_j^\prime) \right),
		\qquad
			\text{ for all } x,x^\prime \in \mathbb{T}^d.
		\label{eq:LDP_switching_MP:dist_function_E}
		\end{equation}
	Then $\Psi \geq 0$, and $\Psi(x,x^\prime) = 0$ holds if and only if $x = x^\prime$, and
	\begin{equation}
		\nabla_1 \left[ \Psi(\cdot,x^\prime)\right](x)
	=
		- \nabla_2 \left[ \Psi(x,\cdot) \right](x^\prime)
	\qquad
		\text{ for all } x,x^\prime \in \mathbb{T}^d.
	\label{eq:LDP_switching_MP:dist_function_asymm_deriv}
	\end{equation}
	By boundedness and upper semicontinuity of $u_1$ and $(-u_2)$, and compactness of $\mathbb{T}^d \times \mathbb{T}^d$, for each $\delta > 0$ there exists a pair $(x_\delta,x_\delta^\prime) \in \mathbb{T}^d \times \mathbb{T}^d$ such that
	$$
		\Phi_\delta(x_\delta,x_\delta^\prime)
	=
		\max_{x,x^\prime} \Phi_\delta(x,x^\prime).
	$$
	Using~$\Phi_\delta(x_\delta, x_\delta) \leq \Phi(x_\delta,x_\delta^\prime)$ and boundedness of~$u_2$, we obtain
	$$
		\Psi(x_\delta,x_\delta^\prime)
	\leq
		2\delta \left( u_2(x_\delta) - u_2(x_\delta^\prime) \right)
	\leq	
		4 \delta \|u_2\|_{L^\infty(\mathbb{T}^d)} = \mathcal{O}(\delta),
	$$
	hence~$\Psi(x_\delta,x_\delta^\prime) \rightarrow 0$ as~$\delta \rightarrow 0$.
	\smallskip
	
	In order to use the sub- and supersolution properties of $u_1$ and $u_2$, introduce the smooth test functions $f^\delta_1$ and $f^\delta_2$ as
	$$
		f_1^\delta(x)
	:=
		u_2(x_\delta^\prime) + \frac{\Psi(x,x_\delta^\prime)}{2\delta}
	\quad
		\text{ and }
	\quad
		f_2^\delta(x^\prime)
	:=
		u_1(x_\delta) - \frac{\Psi(x_\delta,x^\prime)}{2\delta},
	$$
	Then $f_1^\delta, f_2^\delta \in C^\infty(\mathbb{T}^d) \subseteq \mathcal{D}(H)$ are both in the domain of $H$, and hence in the domain of $H_1$ and $H_2$, respectively. Furthermore, $(u_1 - f_1)$ has a maximum at $x = x_\delta$, and $(f_2 - u_2)$ has a maximum at $x^\prime = x_\delta^\prime$, by definition of $(x_\delta,x_\delta^\prime)$ and $\Phi_\delta$. Since $u_1$ is a strong subsolution of $(1 - \tau H_1) u = h$, 
	$$
		u_1(x_\delta) - \tau H_1 f_1^\delta(x_\delta) - h(x_\delta)
	\leq
		0,
	$$
	and since $u_2$ is a strong supersolution of $(1 - \tau H_2) u = h$, 
	$$
		u_2(x_\delta^\prime) - \tau H_2 f_2^\delta(x_\delta^\prime) - h(x_\delta^\prime)
	\geq
		0.	
	$$
	Thereby, we can estimate $\max (u_1 - u_2)$ as
	\begin{align*}
		\max_{\mathbb{T}^d}(u_1 - u_2)
	&\leq
		\Phi_\delta(x_\delta,x_\delta^\prime)	\\
	&\leq
		u_1(x_\delta) - u_2(x_\delta^\prime)	\\
	&\leq
		\tau \left[ H_1 f_1^\delta(x_\delta)) - H_2 f_2^\delta(x_\delta^\prime)\right]
	+
		h(x_\delta) - h(x_\delta^\prime)	\\
	&=
		\tau \left[ \mathcal{H}_1 (\nabla f_1^\delta(x_\delta)) - \mathcal{H}_2 (\nabla f_2^\delta(x_\delta^\prime))\right]
	+
		h(x_\delta) - h(x_\delta^\prime).
	\end{align*}
	By \eqref{eq:LDP_switching_MP:dist_function_asymm_deriv}, 
	$
		\nabla f_1^\delta(x_\delta)
	=
		\nabla f_2^\delta(x_\delta^\prime)
	=:
		p_\delta \in \mathbb{R}^d,
	$
	which establishes \eqref{eq:LDP_switching_MP:CP_visc_ineq}, and thereby finishes the proof.
	\end{proof}
	The rest of the section, we prove Lemmas \ref{lemma:LDP_switching_MP:local_op_strong_sol}, \ref{lemma:LDP_switching_MP:H_to_H1-H2} and \ref{lemma:LDP_switching_MP:H1_leq_H2}. Regarding Lemma \ref{lemma:LDP_switching_MP:local_op_strong_sol}, a proof for single valued operators is given in~\cite[Lemma 9.9]{FengKurtz2006}.
	\begin{proof}[Proof of Lemma \ref{lemma:LDP_switching_MP:local_op_strong_sol}]
	Let $\tau>0$, $h\in C(\mathbb{T}^d)$. We verify that subsolutions are strong subsolutions. For a subsolution $u_1\in\text{BUSC}(\mathbb{T}^d)$ of $(1-\tau H)u=h$ and $(f,H_{f,\varphi}) \in H$, let $x\in \mathbb{T}^d$ be such that~$(u_1-f)(x)=\sup(u_1-f)$. The function~$\tilde{f}(x^\prime)=\Psi(x^\prime,x)$ with $\Psi(x^\prime,x)$ from~\eqref{eq:LDP_switching_MP:dist_function_E} is smooth and therefore in the domain $\mathcal{D}(H)$. Then $x$ is the unique maximal point of $(u_1-(f+\tilde{f}))$,
	\[
	(u_1-(f+\tilde{f}))(x)=\sup_{\mathbb{T}^d}(u_1-(f+\tilde{f})).
	\]
	Since~$u_1$ is a subsolution, there exists at least one element $(x,z^\prime)\in \mathbb{T}^d\times E^\prime$ such that the subsolution inequality with test function~$f+\tilde{f}$ holds. Since $x$ is the only point maximising $u_1-(f+\tilde{f})$, the viscosity-subsolution inequality with test function $(f+\tilde{f})$ holds for the point $x\in \mathbb{T}^d$ and some point $z^\prime\in E^\prime$:
	\[
	u_1(x)-\tau H_{f+\tilde{f},\varphi}(x,z^\prime)-h(x)\leq 0.
	\]
	Since~$\nabla\tilde{f}(x)=0$ and~$H$ depends only on gradients by~\ref{MM:item:T1}, we obtain
	\[
	H_{f+\tilde{f},\varphi}(x,z^\prime)
	=
	H_{\varphi}\left((\nabla f+\nabla \tilde{f})(x),z^\prime\right)
	=
	H_{\varphi}(\nabla f(x),z^\prime)
	=
	H_{f,\varphi}(x,z^\prime).
	\]
	Hence the same point~$(x,z')$ satisfies
	\[
	u_1(x)-\tau H_{f,\varphi}(x,z^\prime)-h(x)\leq 0.
	\]
	Thus~$u_1$ is a strong subsolution. The argument is similar for the supersolution case, where one can use $(-\tilde{f})$.
	\smallskip
	
	Vice versa, when given a strong sub- or supersolution $u_1$ or $u_2$, for every $f\in\mathcal{D}(H)$, $(u_1-f)$ and $(f-u_2)$ attain their suprema at some $x_1,x_2\in \mathbb{T}^d$ due to the continuity assumptions on the domain of $H$, the half continuity properties of $u_1$ and $u_2$,  and compactness of $\mathbb{T}^d$. By the strong solution properties, the sub- and supersolution inequalities follow.
	\end{proof}
	\begin{proof}[Proof of Lemma \ref{lemma:LDP_switching_MP:H_to_H1-H2}]
	Let $u_1\in\text{BUSC}(\mathbb{T}^d)$ be a strong subsolution of~$(1-\tau H)u=h$, that is for any $(f,H_{f, \varphi})\in\tilde{H}$, if~$(u_1-f)(x)=\sup(u_1-f)$ for a point~$x\in \mathbb{T}^d$, then there exists~$z^\prime\in E^\prime$ such that
	$$
		u_1(x)-\tau H_{f, \varphi}(x,z^\prime)-h(x)\leq 0.
	$$
	Let $f\in\mathcal{D}(H_1)=\mathcal{D}(H)$ and $x\in \mathbb{T}^d$ be such that $(u_1-f)(x)=\sup(u_1-f)$. For any $\varphi$ there exists a point $z^\prime\in E^\prime$ such that the above subsolution inequality holds. Therefore for all~$x$,
	$$
		u_1(x)-h(x)
	\leq 
		\tau\sup_{z^\prime\in E^\prime} H_{f, \varphi}(x,z^\prime).
	$$
	Since the point $x\in \mathbb{T}^d$ is independent of $\varphi$, we obtain
	\[
	u_1(x)-H_1f(x)-h(x)\overset{\text{def}}{=}
	u_1(x)-\tau\inf_{\varphi}\sup_{z^\prime \in E^\prime} H_{f, \varphi}(x,z^\prime)-h(x)\leq 0.
	\]
	The argument is similar for supersolutions.
	\end{proof}
	\begin{proof}[Proof of Lemma \ref{lemma:LDP_switching_MP:H1_leq_H2}]
	By~\ref{MM:item:T2},~\ref{MM:item:C3}, for every $p\in\mathbb{R}^d$ there exists a $\varphi_p\in C(E^\prime)$ such that for all $z^\prime\in E^\prime$,
	\[
	H_{\varphi_p}(p,z^\prime)=\mathcal{H}(p).
	\]
	Thus
	\[
	\sup_{z^\prime \in E^\prime} H_{\varphi_p}(p,z^\prime)
	=
	\mathcal{H}(p)
	=
	\inf_{z^\prime\in E^\prime}H_{\varphi_p}(p, z^\prime).
	\]
	Taking the infimum and supremum over~$\varphi$, we find
	\begin{align*}
		\mathcal{H}_1(p)
	&=
		\inf_{\varphi}\sup_{z^\prime}H_{\varphi}(p,z^\prime)\\
	&\leq
		\sup_{z^\prime}H_{\varphi_p}(p,z^\prime)
		=
		\mathcal{H}(p)
		=\inf_{z^\prime}H_{\varphi_p}(p,z^\prime)\\
	&\leq
		\sup_{\varphi}\inf_{z^\prime}H_{\varphi}(p,z^\prime)
		=
		\mathcal{H}_2(p),
	\end{align*}
	which finishes the proof.
	\end{proof}
	\subsection{Proof of action-integral representation}
	\label{section:action_integral}
	In this section, we first prove Theorem~\ref{thm:results:action_integral_representation} by means of Proposition~\ref{prop:action_integral:bold_H} below. The rest of the section is then devoted to proving Proposition~\ref{prop:action_integral:bold_H}.
	\begin{proposition}
	Under the same assumptions of Theorems~\ref{thm:results:LDP_switching_MP} and~\ref{thm:results:action_integral_representation}, 
	define the operator $\mathbf{H} : \mathcal{D}(\mathbf{H}) \subseteq C^1(\mathbb{T}^d) \to C(\mathbb{T}^d)$ on the domain $\mathcal{D}(\mathbf{H}) = \mathcal{D}(H)$ by setting $\mathbf{H} f(x) := \mathcal{H}(\nabla f(x))$. Then:
	\begin{enumerate}[(i)]
	\item
	The Legendre-Fenchel transform $\mathcal{L}(v) := \sup_{p \in \mathbb{R}^d} (p\cdot v - \mathcal{H}(p))$ and the operator $\mathbf{H}$ satisfy Conditions 8.9, 8.10 and 8.11 of \cite{FengKurtz2006}.
	\item
	For all $\tau > 0$ and $h \in C(\mathbb{T}^d)$, the comparison principle holds for
	$
		(1 - \tau \mathbf{H}) u = h.
	$
	\item
	For all $\tau > 0$ and $h \in C(\mathbb{T}^d)$, viscosity solutions of $(1 - \tau H) u = h$ are also viscosity solutions of $(1 - \tau \mathbf{H}) u = h$.
	\end{enumerate}
	\label{prop:action_integral:bold_H}
	\end{proposition}	
	\begin{proof}[Proof of Theorem \ref{thm:results:action_integral_representation}]
	Let $V(t) : C(\mathbb{T}^d) \to C(\mathbb{T}^d)$ be the semigroup
	$$
		V(t)
	=
		\lim_{m \to \infty}
		\left[
			\left( 1 - \frac{t}{m} H \right)^{-1}
		\right]^m,
	$$
	where the resolvant $(1 - \tau H)^{-1}$ is defined by means of unique viscosity solutions to the equation $(1 - \tau H) u = h$, and the limit is made precise in Theorem 6.13, (d), of \cite{FengKurtz2006}. Furthermore, let $V_{\mathrm{NS}}(t) : C(\mathbb{T}^d) \to C(\mathbb{T}^d)$ be the Nisio semigroup with cost function $\mathcal{L}$, that is $V_{\mathrm{NS}}(t)$ is defined as
	$$
		V_{\text{NS}}(t) f(x)
	=
		\sup_{ \substack{ \gamma \in \mathrm{AC}_{\mathbb{T}^d}[0,\infty) \\ \gamma(0) = x}}
		\left[
			f(\gamma(t)) - \int_0^t \mathcal{L}(\dot{\gamma}(s)) \, \dd s
		\right],
	$$
	where $\mathrm{AC}_{\mathbb{T}^d}[0,\infty)$ denotes the set of absolutely continuous paths in the torus. In Definition 8.1 and Equation 8.10 in \cite{FengKurtz2006}, relaxed controls are considered in order to cover a general class of possible cost functions. Since the Legendre-Fenchel transform $\mathcal{L}(v)$ is convex, it follows that $V_{\mathrm{NS}}(t)$ equals the semigroup given in 8.10 of \cite{FengKurtz2006} by using that $\lambda_s = \delta_{\dot{x}(s)}$ is an admissible control, and by applying Jensen's inequality, an argument that is given for example in Theorem 10.22 in \cite{FengKurtz2006}. Below we prove that $V(t) = V_{\mathrm{NS}}(t)$; by Theorem 8.14 in \cite{FengKurtz2006}, if $V(t) = V_{\mathrm{NS}}(t)$ on $C(\mathbb{T}^d)$, then the rate function of our Theorem~\ref{thm:results:LDP_switching_MP} satisfies the control representation 8.18 of \cite{FengKurtz2006}. The action-integral representation follows by again applying Jensen's inequality.
	\smallskip
	
	It remains to prove that $V(t) =  V_{\mathrm{NS}}(t)$. By (i) and (ii) of Proposition \ref{prop:action_integral:bold_H}, the Conditions of~\cite[Theorem 8.27]{FengKurtz2006} are satisfied, so that we have $V_{\mathrm{NS}}(t) = \mathbf{V}(t)$, where $\mathbf{V}(t)$ is generated by means of unique viscosity solutions to the equation $(1- \tau \mathbf{H}) u = h$ as shown in~\cite[Theorem 8.27]{FengKurtz2006}, that is
	$$
		\mathbf{V}(t)
	=
		\lim_{m \to \infty}
		\left[
			\left( 1 - \frac{t}{m} \mathbf{H} \right)^{-1}
		\right]^m.
	$$
	Part (iii) of Proposition \ref{prop:action_integral:bold_H} implies by Corollary 8.29 of \cite{FengKurtz2006} that $V(t) = \mathbf{V}(t)$.
	\end{proof}
	\begin{proof}[Proof of (i) in Proposition \ref{prop:action_integral:bold_H}]
	We first show that the following conditions imply Conditions 8.9, 8.10 and 8.11 of \cite{FengKurtz2006}, which are formulated in order to cover a more general and non-compact setting.
	\begin{enumerate}[(i)]
	\item The function $\mathcal{L}:\mathbb{R}^d\rightarrow[0,\infty]$ is lower semicontinuous and for every $C \geq 0$, the level set
	$
		\{v\in \mathbb{R}^d\,:\,\mathcal{L}(v)\leq C\}
	$
	is relatively compact in $\mathbb{R}^d$.
	\item For all $f\in\mathcal{D}(H)$ there exists a right continuous, nondecreasing function $\psi_f:[0,\infty)\rightarrow[0,\infty)$ such that for all $(x_0,v)\in \mathbb{T}^d \times \mathbb{R}^d$,
	\[
	|\nabla f(x_0)\cdot v|\leq \psi_f(\mathcal{L}(v))\qquad \text{and} \qquad
	\lim_{r\rightarrow\infty}\frac{\psi_f(r)}{r}=0.
	\]
	\item For each $x_0\in E$ and every $f\in\mathcal{D}(\mathbf{H})$, there exists an absolutely continuous path $x : [0,\infty) \to \mathbb{T}^d$ such that
	\begin{equation}
		\int_0^t \mathcal{H}(\nabla f (x(s))) \, ds
	=
		\int_0^t \left[
		\nabla f(x(s)) \cdot \dot{x}(s) - \mathcal{L}(\dot{x}(s))
		\right] \, ds.
	\label{eq:action_integral:optimal_path_for_H}
	\end{equation}
	\end{enumerate}
	Then regarding Condition 8.9 of \cite{FengKurtz2006}, the operator $A f(x,v) := \nabla f(x) \cdot v$ on the domain $\mathcal{D}(A) = \mathcal{D}(H)$ satisfies (1). For (2), we can take $\Gamma = \mathbb{T}^d \times \mathbb{R}^d$, and for $x_0 \in \mathbb{T}^d$, take the pair $(x,\lambda)$ with $x(t) = x_0$ and $\lambda(dv \times dt) = \delta_{0} (dv) \times dt$. Part (3) is a consequence of the above Item~(i). Part (4) follows since $\mathbb{T}^d$ is compact. Part (5) is implied by the above Item~(ii). Condition~8.10 is implied by Condition~8.11 and the fact that $\mathbf{H}1 = 0$, see Remark 8.12 (e) in \cite{FengKurtz2006}. Finally, Condition~8.11 is implied by the above Item~(iii), with the control $\lambda(dv \times dt) = \delta_{\dot{x}(t)}(dv) \times dt$.
	\smallskip

	We turn to verifying (i), (ii) and (iii). Since $\mathcal{H}(0) = 0$, we have $\mathcal{L} \geq 0$. The Legendre-transform $\mathcal{L}$ is convex, and lower semicontinuous since the map $\mathcal{H}(p)$ is convex and finite-valued, hence in particular continuous. For $C \geq 0$, we prove that the set $\{v\in\mathbb{R}^d\,:\,\mathcal{L}(v)\leq C\}$ is bounded, and hence is relatively compact. For any $p \in \mathbb{R}^d$ and $v \in \mathbb{R}^d$, we have $p \cdot  v\leq \mathcal{L}(v) + \mathcal{H}(p)$. Thereby, if $\mathcal{L}(v)\leq C$, then
	$
		|v|
	=
		\sup_{|p|=1} p \cdot v 
	\leq
		\sup_{|p|=1} \left[
		\mathcal{L}(v) + \mathcal{H}(p)
		\right]
	\leq
		C + C_1,
	$
	where $C_1$ exists due to continuity of $\mathcal{H}$. Then for $R := C + C_1$,
	$
		\{ v \, : \, \mathcal{L}(v) \leq C \} \subseteq 
		\{ v \, : \, |v| \leq R\},
	$
	thus $\{\mathcal{L}\leq C\}$ is a bounded subset in $\mathbb{R}^d$.
	\smallskip
	
	Item~(ii) can be proven as in~\cite[Lemma 10.21]{FengKurtz2006}. We finish the proof by verifying (iii). This is shown for instance in~\cite[Lemma 3.2.3]{Kraaij2016} under the assumption of continuous differentiability of $\mathcal{H}(p)$, by solving a differential equation with a globally bounded vectorfield. Here, we verify (iii) under the milder assumption of convexity of $\mathcal{H}(p)$ by solving a suitable subdifferential equation. For $p_0 \in \mathbb{R}^d$, define the subdifferential $\partial \mathcal{H}(p_0)$ at $p_0$ as the set
	$$
	\partial \mathcal{H}(p_0)
	:=
	\{\xi \in \mathbb{R}^d \; | \; \forall p \in \mathbb{R}^d\;: 
	\mathcal{H}(p)
	\geq
	\mathcal{H}(p_0) + \langle \xi, p - p_0\rangle\}.
	$$
	We shall solve for any $f \in C^1(\mathbb{T}^d)$ the subdifferential equation~$\dot{x} \in \partial \mathcal{H}(\nabla f(x))$. This means we show that for any initial condition~$x_0\in \mathbb{T}^d$, there exists an absolutely continuous path $x:[0,\infty)\rightarrow \mathbb{T}^d$ satisfying both~$x(0)=x_0$ and~$\dot{x}(t) \in \partial \mathcal{H}(\nabla f(x(t)))$ almost everywhere on $[0,\infty)$. Then \eqref{eq:action_integral:optimal_path_for_H} follows by noting that
	$
	\mathcal{H}(\nabla f(y))\geq \nabla f(y)\cdot v-\mathcal{L}(v)
	$
	for all $y \in \mathbb{T}^d$ and $v \in \mathbb{R}^d$, by convex duality. In particular, 
	$
	\mathcal{H}(\nabla f(x(s)))\geq \nabla f(x(s))\cdot \dot{x}(s)-\mathcal{L}(\dot{x}(s)),
	$
	and integrating gives one inequality in \eqref{eq:action_integral:optimal_path_for_H}. Regarding the other inequality, since $\dot{x}\in\partial \mathcal{H}(\nabla f(x))$, we know that for almost every $t \in [0,\infty)$ and for all $p\in\mathbb{R}^d$, we have
	$
		\mathcal{H}(p)
	\geq 
		\mathcal{H}(\nabla f(x(t)))+\dot{x}(t)\cdot(p-\nabla f(x(t))).
	$
	Therefore, a.e. on $[0,\infty)$,
	\begin{align*}
		\mathcal{H}(\nabla f(x(t))) &\leq \nabla f(x(t))\cdot\dot{x}(t)-\sup_{p\in\mathbb{R}^d}\left[p\cdot \dot{x}(t)-\mathcal{H}(p)\right] \\&= \nabla f(x(t))\cdot\dot{x}(t)-\mathcal{L}(\dot{x}(t)),
	\end{align*}
	and integrating gives the other inequality.
	\smallskip
	
	For solving the subdifferential equation, define $F: \mathbb{R}^d \to 2^{ \mathbb{R}^d}$ by $F(x):=\partial \mathcal{H}(\nabla f(x))$, where the function $f \in C^1(\mathbb{T}^d)$ is regarded as a periodic function on $\mathbb{R}^d$. We apply Lemma 5.1 in \cite{De92} for solving $\dot{x}\in F(x)$. The conditions of Lemma 5.1 in the case of $\mathbb{R}^d$ are satisfied if the following holds: $\sup_{x \in \mathbb{R}^d} \|F(x)\|_{\text{sup}}$ is finite, for all $x\in \mathbb{R}^d$, the set $F(x)$ is non-empty, closed and convex, and the map $x\mapsto F(x)$ is upper semicontinuous.
	\smallskip
	
	For $\xi \in F(x)$, note that for all $p\in\mathbb{R}^d$ 
	$
	\xi \cdot (p-\nabla f(x))\leq \mathcal{H}(p)-\mathcal{H}(\nabla f(x))
	$.
	Therefore, by shifting $p=p^\prime+\nabla f(x)$, we obtain for all $p^\prime\in\mathbb{R}^d$ that
	$
	\xi \cdot p^\prime\leq \mathcal{H}(p^\prime+\nabla f(x))-\mathcal{H}(\nabla f(x))
	$.
	By continuous differentiability and periodicity of~$f$, and continuity of~$\mathcal{H}$, the right-hand side is bounded in $x$, and we obtain
	\begin{align*}
		\sup_{x \in \mathbb{R}^d} \sup_{\xi \in F(x)}|\xi|
	&=
		\sup_{x \in \mathbb{R}^d} \sup_{\xi \in F(x)} \sup_{|p^\prime| = 1} \xi \cdot p^\prime\\
	&\leq 
		\sup_{x \in \mathbb{R}^d} \sup_{\xi \in F(x)} \sup_{|p^\prime|=1}\left[\mathcal{H}(p^\prime+\nabla f(x))-\mathcal{H}(\nabla f(x))\right] < \infty.
	\end{align*}
	For any $x \in \mathbb{R}^d$, the set $F(x)$ is non-empty, since the subdifferential of a proper convex function $\mathcal{H}(\cdot)$ is nonempty at points where $\mathcal{H}(\cdot)$ is finite and continuous~\cite{rockafellar1966characterization}. Furthermore, $F(x)$ is convex and closed, which follows from the properties of a subdifferential set.
	\smallskip
	
	Regarding upper semicontinuity, recall the definition from \cite{De92}: the map $F:\mathbb{R}^d \to 2^{\mathbb{R}^d} \setminus \{\emptyset\}$ is upper semicontinuous if for all closed sets $A\subseteq\mathbb{R}^d$, the set $F^{-1}(A)\subseteq \mathbb{R}^d$ is closed, where
	$
	F^{-1}(A)=\{x\in \mathbb{R}^d\;|\;F(x)\cap A\neq\emptyset\}.
	$
	Let $A\subseteq\mathbb{R}^d$ be closed and $x_n\rightarrow x$ in $\mathbb{R}^d$, with $x_n\in F^{-1}(A)$. That means for all $n\in\mathbb{N}$ that the sets
	$
	\partial \mathcal{H}(\nabla f(x_n))\cap A
	$
	are non-empty, and consequently, there exists a sequence $\xi_n \in F(x_n)\cap A$. We proved above that the set $F(y)\cap A$ is uniformly bounded in $y \in \mathbb{R}^d$. Hence the sequence~$\xi_n$ is bounded, and passing to a subsequence if necessary, it converges to some~$\xi$.
	By definition of $F(x_n)$, for all $p\in\mathbb{R}^d$,
	\begin{align*}
		\xi_n(p-\nabla f(x_n))
	\leq 
		\mathcal{H}(p)-\mathcal{H}(\nabla f(x_n)).
	\end{align*}
	Passing to the limit, we obtain that for all $p \in \mathbb{R}^d$,
	$$
		\xi (p-\nabla f(x))
	\leq 
		\mathcal{H}(p)-\mathcal{H}(\nabla f(x)).
	$$
	This implies by definition that $\xi\in\partial \mathcal{H}(\nabla f(x))$. Since $\xi_n\in A$ and $A$ is closed, we have $\xi\in A$. Hence $x\in F^{-1}(A)$, and $F^{-1}(A)$ is indeed closed.
	\end{proof}
	\begin{proof}[Proof of (ii) in Proposition \ref{prop:action_integral:bold_H}]
	The comparison principle for the operator $\mathbf{H}$ follows from the fact that~$\mathbf{H}f=\mathcal{H}(\nabla f)$ depends on~$x$ only via gradients. Indeed, for subsolutions $u_1$ and supersolutions $u_2$ of $(1-\tau \mathbf{H})u =h $, we have
	$
	\max(u_1-u_2)
	\leq 
	\tau [\mathcal{H}
	\left(
	\nabla f_1(x_\delta)
	\right)
	-
	\mathcal{H}\left(\nabla f_2(x_\delta^\prime)\right)]
	+
	h(x_\delta)-h(x_\delta^\prime),
	$
	with test functions $f_1,f_2\in\mathcal{D}(H)$ satisfying
	$
	\nabla f_1(x_\delta)
	=
	\nabla f_2(x_\delta^\prime),
	$
	and~$\text{dist}(x_\delta,x_\delta^\prime)\rightarrow 0$ as~$\delta \to 0$. Therefore $\mathcal{H}\left(\nabla f_1(x_\delta)\right)-\mathcal{H}\left(\nabla f_2(x_\delta^\prime)\right) = 0$, and~$\max(u_1-u_2)\leq 0$ follows by taking the limit $\delta\rightarrow 0$.
	\end{proof} 
	\begin{proof}[Proof of (iii) in Proposition \ref{prop:action_integral:bold_H}]
	Let $u \in C(\mathbb{T}^d)$ be a viscosity solution of the equation~$(1 - \tau H) u = h$. By Lemmas~\ref{lemma:LDP_switching_MP:local_op_strong_sol} and~\ref{lemma:LDP_switching_MP:H_to_H1-H2}, $u$ is a strong viscosity subsolution of $(1 - \tau H_1) u = h$ and a strong viscosity supersolution of $(1- \tau H_2) u = h$. In the proof of Lemma~\ref{lemma:LDP_switching_MP:H1_leq_H2} we obtained $\mathcal{H}_1 \leq \mathcal{H} \leq \mathcal{H}_2$, which in particular implies the inequalities 
	$
		- H_1 \geq -\mathbf{H} \geq -H_2.
	$
	With that, we find that $u$ is both a strong viscosity sub- and supersolution of $(1 - \tau \mathbf{H}) u = h$.
	\end{proof}
	\section{Proofs of large deviations for molecular motors}
	\sectionmark{LDP for molecular motors}
	\label{section:LDP_in_MM}
	In this section, we prove the theorems of Section \ref{subsection:results:LDP_in_MM} about the stochastic processes motivated by molecular-motor systems. 
	The proofs regarding the continuous model (Theorems~\ref{thm:results:LDP_cont_MM_Limit_I} and~\ref{thm:results:LDP_cont_MM_Limit_II}) are collected in Section~\ref{subsection:cont_MM}, and regarding the discrete model (Theorems~\ref{thm:results:LDP_cont_MM_Limit_I} and~\ref{thm:results:LDP_cont_MM_Limit_II}) in Section~\ref{subsection:discr_MM}. In each proof we verify the conditions of the general theorems for switching Markov processes (Theorems~\ref{thm:results:LDP_switching_MP} and~\ref{thm:results:action_integral_representation}). Finally, we prove in Section~\ref{subsection:detailed_balance} the representation of Hamiltonians~$\mathcal{H}(p)$ that we use to prove symmetry under the detailed balance condition.
	\subsection{Proof for the continuous models}
	\label{subsection:cont_MM}
	In this section, we consider the stochastic process $(X^\varepsilon_t, I^\varepsilon_t)$ of Defintion~\ref{def:intro:cont_MM} and prove Theorems~\ref{thm:results:LDP_cont_MM_Limit_I} and \ref{thm:results:LDP_cont_MM_Limit_II}. The generator~$L_\varepsilon$ is given by
	\begin{multline*}
	L_\varepsilon f(x,i)= \varepsilon \frac{1}{2} \Delta_x f(\cdot,i) (x) +b^i(x/\varepsilon) \cdot \nabla_x f(\cdot,i) (x)\\
	+ \sum_{j = 1}^J \frac{1}{\varepsilon} \gamma(\varepsilon) r_{ij}(x/\varepsilon)
	\left[
	f(x,j) - f(x,i)
	\right],
	\end{multline*}
	with state space $E_\varepsilon = \mathbb{T}^d \times \{1,\dots, J\} = \{(x,i)\}$, drifts $b^i \in C^\infty(\mathbb{T}^d)$, jump rates $r_{ij} \in C^\infty(\mathbb{T}^d;[0,\infty))$, and~$\gamma(\varepsilon)>0$. We frequently write $f(x,i) = f^i(x)$. The nonlinear generators defined by~$H_\varepsilon f = \varepsilon e^{-f / \varepsilon} L_\varepsilon e^{f(\cdot) / \varepsilon}$ are given by
	\begin{multline}\label{eq:LDP_MM:cont:H_varepsilon}
		H_\varepsilon f(x,i)= \frac{\varepsilon}{2} \Delta_x f^i(x)+\frac{1}{2}|\nabla_x f^i(x)|^2+b^i\left(x/\varepsilon\right)\nabla_x f^i(x) 
		\\+ \gamma(\varepsilon) \sum_{j=1}^J r_{ij}\left(x/\varepsilon\right)
		\left[
		e^{ \left(f(x,j)-f(x,i)\right) / \varepsilon}-1
		\right].
	\end{multline}
	\subsubsection{Proof of Theorem \ref{thm:results:LDP_cont_MM_Limit_I}}
	\label{subsubsection:LDP_cont_MM_Limit_I}	
	\begin{proof}[Verification of~\ref{MM:item:T1} of Theorem \ref{thm:results:LDP_switching_MP}]
	With~$f_\varepsilon(x,i)=f(x)+\varepsilon\,\varphi\left(x/\varepsilon, i\right))$, we find
	\begin{multline*}
		H_\varepsilon f_\varepsilon(x,i)
	=
		\frac{\varepsilon}{2}\Delta f(x)
	+
		\frac{1}{2}\Delta_y\varphi^i\left(x/\varepsilon\right)
	+
		\frac{1}{2}\big|\nabla f(x)+\nabla_y\varphi^i\left(x/\varepsilon\right)\big|^2
	\\+
		b^i\left(x/\varepsilon\right)\left(\nabla f(x)+\nabla_y\varphi^i\left(x/\varepsilon\right)\right)	\\
	+
		\sum_{j = 1}^J r_{ij}\left(x/\varepsilon\right)
		\left[
		e^{\varphi\left(x/\varepsilon, j\right)-\varphi\left(x/\varepsilon, i\right)}-1
		\right],
	\end{multline*}
	where $\nabla_y$ and $\Delta_y$ denote the gradient and Laplacian with respect to the variable $y = x/\varepsilon$. The only term of order $\varepsilon$ that remains is $\varepsilon\,\Delta f(x)/2$. This suggests to take the remainder terms as the definition of the multivalued operator $H$. In the notation of Theorem \ref{thm:results:LDP_switching_MP}, we choose $E^\prime = \mathbb{T}^d \times\{1,\dots,J\}$ as the state space of the macroscopic variables, and define
	\begin{align}\label{MM:eq:limit-H-proof-cont-MM-model}
		H:=
		\left\{
		(f, H_{f,\varphi}) \, : \,
		f \in C^2(\mathbb{T}^d), \;
		H_{f, \varphi} \in C(\mathbb{T}^d \times E^\prime) \text{ and }
		\varphi \in C^2(E^\prime)
		\right\}.
	\end{align}
	In $H$, the image functions $H_{f,\varphi} : \mathbb{T}^d \times E^\prime \to \mathbb{R}$ are defined by
	\begin{multline}
		H_{f,\varphi}(x,y,i)
	:=
		\frac{1}{2}\Delta_y\varphi^i(y) + \frac{1}{2} \big|\nabla f(x)+\nabla_y\varphi^i(y) \big|^2 +b^i(y)(\nabla f(x) + \nabla_y\varphi^i(y))
	\\+
		\sum_{j = 1}^J r_{ij}(y)\left[e^{\varphi(y, j)-\varphi(y, i)}-1\right],
	\label{eq:LDP_MM:contI:limit_op_H}
	\end{multline}
	where we write $\varphi = (\varphi^1, \dots, \varphi^J)$ via the identification $C^2(E^\prime) \simeq (C^2(\mathbb{T}^d))^J$. Define the maps $\eta_\varepsilon^\prime:E_\varepsilon \to E^\prime$ by $\eta_\varepsilon^\prime(x,i) := (x/\varepsilon,i)$, and recall that $\eta_\varepsilon(x,i) := x$ are projections $\eta_\varepsilon:E_\varepsilon\to \mathbb{T}^d$.
	\smallskip
	
	We now verify~\ref{MM:item:C1},~\ref{MM:item:C2} and~\ref{MM:item:C3} of~\ref{MM:item:T1}. For~\ref{MM:item:C1}, for any $(x,y,i) \in \mathbb{T}^d \times E^\prime$, we search for elements $(y_\varepsilon,i_\varepsilon) \in \mathbb{T}^d \times \{1,\dots,J\}$ such that both $\eta_\varepsilon(y_\varepsilon,i_\varepsilon) \to x$ and $\eta_\varepsilon^\prime(y_\varepsilon,i_\varepsilon) \to (y,i)$ as $\varepsilon \to 0$. The point $y_\varepsilon := x + \varepsilon(y - x)$ satisfies $y_\varepsilon \to x$ and $y_\varepsilon / \varepsilon = y$, since $x / \varepsilon = x$ in $\mathbb{T}^d$. Therefore,~\ref{MM:item:C1} holds with $y_\varepsilon = x + \varepsilon(y-x)$ and $i_\varepsilon = i$. Regarding~\ref{MM:item:C2}, let $(f, H_{f,\varphi}) \in H$. Then the function~$f_\varepsilon$ defined by $f_\varepsilon(x,i):=f(x)+\varepsilon\,\varphi\left(x/\varepsilon, i\right)$ satisfies
	$$
		\|f\circ\eta_\varepsilon-f_\varepsilon\|_{L^\infty(E_\varepsilon)}
	=
		\sup_{(x,i)\in E_\varepsilon}|f(x)-f_\varepsilon(x,i)|
	=
		\varepsilon\cdot \|\varphi \|_{L^\infty(E_\varepsilon)}\xrightarrow{\varepsilon\rightarrow 0}0,
	$$
	and
	\begin{align*}
		\|H_{f,\varphi}\circ\eta_\varepsilon^\prime-H_\varepsilon f_\varepsilon\|_{L^\infty(E_\varepsilon)}
	&=
		\sup_{(x,i)\in E_\varepsilon}|H_{f,\varphi}(x,x/\varepsilon,i)-H_\varepsilon f_\varepsilon(x,i)|\\
	&= 
		\frac{\varepsilon}{2}\;\sup_{(x,i)\in E_\varepsilon}|\;\Delta f(x)|\leq \varepsilon \frac{1}{2}\sup|\Delta f|
		\xrightarrow{\varepsilon\rightarrow 0}0.
	\end{align*}
	Item~\ref{MM:item:C3}, the fact that the images $H_{f,\varphi}$ depend on $x$ only via the gradients of $f$, can be recognized in~\eqref{eq:LDP_MM:contI:limit_op_H}.
	\end{proof}
	\begin{proof}[Verification of (T2) of Theorem \ref{thm:results:LDP_switching_MP}]
	Let $f$ be a function in $\mathcal{D}(H)=C^2(\mathbb{T}^d)$ and $x\in \mathbb{T}^d$. We establish the existence of a vector function $\varphi=(\varphi^1,\dots,\varphi^J)\in (C^2(\mathbb{T}^d))^J$ such that for all $(y,i) \in E^\prime = \mathbb{T}^d \times \{1,\dots,J\}$ and some constant $\mathcal{H}(\nabla f(x)) \in \mathbb{R}$, we have
	$$
		H_\varphi(\nabla f(x),y,i)
	=
		\mathcal{H}(\nabla f(x)).
	$$
	For the flat torus $E = \mathbb{T}^d$, this means that for fixed $\nabla f(x)=p\in\mathbb{R}^d$, we search for a vector function $\varphi_p$ such that $\tilde{H}_{\varphi_p}(p,y,i) = \mathcal{H}(p)$ becomes independent of the variables $(y,i)\in E^\prime$. We can find this vector function by solving a principal eigenvalue problem. We prove Item~\ref{MM:item:T2} with the following Lemma.
	\begin{lemma}
	Let $E^\prime = \mathbb{T}^d \times \{1, \dots, J\}$ and $H$ be the limit operator~\eqref{MM:eq:limit-H-proof-cont-MM-model}. Then:
	\begin{enumerate}[(a)]
	\item For $f \in \mathcal{D}(H)$, the limiting images $H_{\varphi}(\nabla f(x),y,i)$ are of the form 
	\[
		H_\varphi(\nabla f(x) , y, i)
	=
		e^{-\varphi(y,i)} \left[ (B_p + V_p + R)e^{\varphi} \right] (y,i),
	\]
	with $ p = \nabla f(x) \in\mathbb{R}^d$, and operators $B_p, V_p, R : C^2(E^\prime) \rightarrow C(E^\prime)$ defined as
	\begin{align*}
		(B_p h)(y,i) 
	&:=
		\frac{1}{2} \Delta_y h(y,i) + \left( p +b^i(y)\right)\cdot \nabla_y h(y,i) \\
		(V_p h)(y,i)
	&:=
		\left(\frac{1}{2} p^2 + p\cdot b^i(y)\right) h(y,i), \\
		(R\,h)(y,i)
	&:=
		\sum_{j = 1}^J r_{ij}(y) \left[h(y,j) - h(y,i)\right].
	\end{align*}
	\item 
	For any $p \in \mathbb{R}^d$, there exists an eigenfunction $g_p = (g_p^1,\dots, g_p^J) \in (C^2(\mathbb{T}^d))^J$ with strictly positive component functions, $g^i_p > 0 $ on $\mathbb{T}^d$ for $ i =1,\dots, J$, and an eigenvalue $\mathcal{H}(p) \in \mathbb{R}$ such that
	\begin{equation}
	\left[B_p + V_p + R\right] g_p 
	=
	\mathcal{H}(p)\,g_p.
	\label{eq:LDP_MM:contI:cell_problem}
	\end{equation}
	\end{enumerate}
	\label{lemma:LDP_MM:contI:principal_eigenvalue}
	\end{lemma}
	Now~(T2) follows by (a) and (b), since with $\varphi_p := \log g_p$,
	\begin{align*}
		H_{\varphi_p}(p , y, i)
		&\overset{(a)}{=}
		e^{-\varphi_p(y,i)} \left[ B_p + V_p + R \right] e^{\varphi_p(y,i)}
		\\&=
		\frac{1}{g_p(y,i)} \left[B_p + V_p + R\right] g_p (y,i)
		\overset{(b)}{=}
		\mathcal{H}(p).
	\end{align*}
	\emph{Proof of Lemma \ref{lemma:LDP_MM:contI:principal_eigenvalue}.} Writing $p = \nabla f(x)$, Item~(a) follows directly by regrouping the terms in~\eqref{eq:LDP_MM:contI:limit_op_H}. Regarding Item~(b),~$\left[B_p + V_p +R\right] g_p = \mathcal{H}(p) g_p$ is a system of weakly coupled nonlinear elliptic PDEs on the flat torus. They are weakly coupled in the sense that the component functions $g_p^i$ are only coupled in the lowest order terms by means of the operator $R$, while the operators $B_p$ and $V_p$ act solely on the diagonal. When cast in matrix form, the eigenvalue problem to solve reads as follows: for $D_p + R$, with a diagonal matrix~$D_p$ 
	and a coupling matrix $R$ with entries $R_{ij} = r_{ij}$ ($i \neq j$) and $R_{ii} = - \sum_{j\neq i}r_{ij}$ on the diagonal,
	\begin{align*}
		D_p
	=
		\begin{pmatrix}
    	B_p^1 + V^1_p	&	&	0	\\
    		& \ddots	&	\\
    	0	&	& B_p^J + V^J_p     
  		\end{pmatrix},
  		\qquad
  		R
	=
		\begin{pmatrix}
		R_{11}	&	&	\geq 0	\\
    		& \ddots	&	\\
    	\geq 0	&		&	R_{JJ}
		\end{pmatrix},
	\end{align*}
	find a strictly positive vector function $g_p > 0$ such that $\left[ D_p + R \right] g_p = \mathcal{H}(p) g_p$. 
	Guido Sweers showed how to obtain the principal eigenvalue for such kind of coupled systems for bounded sets~$\Omega \subseteq \mathbb{R}^d$ under Dirichlet boundary conditions~\cite{Sweers92}, but the line of argument applies to the periodic setting as well---we summarize the result in Proposition \ref{proposition:appendix:PrEv:fully_coupled_system}. Under our irreducibility assumption on~$R$, there exists a $\lambda(p)$ and $g_p > 0$ such that $\left[-D_p - R \right] g_p = \lambda(p) g_p$. Thereby, $\left[D_p + R \right] g_p = \mathcal{H}(p) g_p$ follows with the same eigenfunction $g_p > 0$ and the principal eigenvalue $\mathcal{H}(p) = - \lambda(p)$. This finishes the verification of~\ref{MM:item:T2}.
	\end{proof}
	\begin{proof}[Verification of (T3) of Theorem \ref{thm:results:action_integral_representation}]
	We prove that the principal eigenvalue $\mathcal{H}(p)$ of Lemma \ref{lemma:LDP_MM:contI:principal_eigenvalue} is convex in $p\in\mathbb{R}^d$ and satisfies $\mathcal{H}(0)=0$.
	%
	To that end, we use an explicit variational representation formula for the principal eigenvalue. By Proposition \ref{proposition:appendix:PrEv:fully_coupled_system}, the eigenvalue $\mathcal{H}(p) = - \lambda(p)$ admits the representation
	\begin{align*}
		\mathcal{H}(p)
	&=
		- \sup_{g > 0} \inf_{z^\prime \in E^\prime} 
		\left\{ 
		\frac{1}{g(z^\prime)}\left[
		(-B_p - V_p - R)g
		\right](z^\prime) 
		\right\} \\
	&=
		\inf_{g > 0} \sup_{z^\prime \in E^\prime}
		\left\{
		\frac{1}{g(z^\prime)}
		\left[
		(B_p + V_p + R)g
		\right](z^\prime)
		\right\} \\
	&=
		\inf_{\varphi} \sup_{z^\prime \in E^\prime}
		\left\{
		e^{-\varphi(z^\prime)}
		\left[(B_p + V_p + R)e^{\varphi}
		\right] (z^\prime)
		\right\}
	=:
		\inf_{\varphi} \sup_{z^\prime \in E^\prime} F(p,\varphi)(z^\prime).
	\end{align*}
	The map $F$ is given by
	\begin{multline*}
		F(p,\varphi)(y,i)
	=
		\frac{1}{2}\Delta\varphi^i(y)
		+
		\frac{1}{2} |\nabla\varphi^i(y)+p|^2
		+
		b^i(y)(\nabla\varphi^i(y)+p)
		\\+
		\sum_{j = 1}^J r_{ij}(y)
		\left[ e^{\varphi^j(y)-\varphi^i(y)}-1 \right],
	\end{multline*}
	and hence is jointly convex in $p$ and $\varphi$. For the eigenfunction $\varphi = \varphi_p$, equality holds in the sense that for any $z \in E^\prime$, we have $\mathcal{H}(p) = F(p,\varphi_p)(z)$. Therefore, we obtain for $\tau\in[0,1]$ and any $p_1, p_2 \in \mathbb{R}^d$ with corresponding eigenfunctions $g_1 = e^{\varphi_1}$ and $g_2 = e^{\varphi_2}$ that
	\begin{align*}
		\mathcal{H}(\tau p_1 + (1-\tau) p_2)
	&=
		\inf_\varphi\sup_{E^\prime}F\left(\tau p_1 + (1-\tau) p_2, \varphi\right)\\
	&\leq 
		\sup_{E^\prime}F\left(\tau p_1 + (1-\tau) p_2, \tau \varphi_1+(1-\tau)\varphi_2\right)\\
	&\leq 
		\sup_{E^\prime}\left[\tau F(p_1,\varphi_1) + (1-\tau) F(p_2,\varphi_2)\right]\\
	&\leq 
		\tau \sup_{E^\prime}F(p_1,\varphi_1) + (1-\tau) \sup_{E^\prime} F(p_2,\varphi_2) \\
	&= 
		\tau \mathcal{H}(p_1) + (1-\tau) \mathcal{H}(p_2).
	\end{align*}
	Regarding the claim $\mathcal{H}(0) = 0$, we choose the constant function $\varphi = (1,\dots,1)$ in the variational representation of $\mathcal{H}(p)$. Thereby, we obtain the estimate $\mathcal{H}(0) \leq 0$. For the opposite inequality, we show that for any $\varphi \in C^2(E^\prime)$ 
	$$
		\lambda(\varphi) := \sup_{z^\prime \in E^\prime}
	\left\{
	e^{-\varphi(z^\prime)}
	\left[
	(B_0 + V_0 + R) e^{\varphi}
	\right](z^\prime)
	\right\}
		\geq 0,
	$$
	which then implies $\mathcal{H}(0) = \inf_{\varphi} \lambda(\varphi) \geq 0$. Let $\varphi\in C^2(E^\prime)$; the continuous function $\varphi$ on the compact set $E^\prime$ admits a global minimum $z_m = (y_m, i_m) \in E^\prime$. Thereby, noting that $V_0 \equiv 0$, we find
	\begin{multline*}
		\lambda(\varphi) \geq e^{-\varphi(z_m)}(B_0 + R) e^{\varphi(z_m)}
	=
		\underbrace{\frac{1}{2} \Delta_y \varphi(y_m,i_m)}_{ \displaystyle \geq 0} 
	+
		\frac{1}{2}|\underbrace{\nabla_y \varphi(y_m,i_m)}_{\displaystyle = 0}|^2
	\\+
		b^{i_m}(y_m)\cdot \underbrace{\nabla_y \varphi(y_m,i_m)}_{\displaystyle = 0}
	+
		\sum_{j\neq i} r_{ij}(y_m) 
		\underbrace{\left[ e^{\varphi(y_m,j) - \varphi(y_m,i_m)} - 1 \right]}_{\displaystyle \geq 0}
	\geq
		0.
	\end{multline*}
	%
	This finishes the verification of (T3), and thereby the proof of Theorem \ref{thm:results:LDP_cont_MM_Limit_I}.
	\end{proof}
	\subsubsection{Proof of Theorem \ref{thm:results:LDP_cont_MM_Limit_II}}
	\label{subsubsection:LDP_cont_MM_Limit_II}	
	In this section, we consider the process $(X^\varepsilon_t, I^\varepsilon_t)$ from Definition \ref{def:intro:cont_MM} in the limit regime $\gamma(\varepsilon) \to \infty$ as $\varepsilon \to 0$. As above in the proof of Theorem \ref{thm:results:LDP_cont_MM_Limit_I}, we start with the nonlinear generator $H_\varepsilon$ from \eqref{eq:LDP_MM:cont:H_varepsilon}, and verify Conditions (T1), (T2) and (T3) of Theorems \ref{thm:results:LDP_switching_MP} and \ref{thm:results:action_integral_representation}.
	\begin{proof}[Verification of (T1) of Theorem \ref{thm:results:LDP_switching_MP}]
	We can not make the same Ansatz as in the proof of Theorem \ref{thm:results:LDP_cont_MM_Limit_I}, since the reaction terms with $\gamma(\varepsilon)$ diverge whenever the exponent remains of order one. We have three different scales: order $1$ via the variable $x$, of order $1/\varepsilon$ via $(x/\varepsilon)$, and of order $\gamma(\varepsilon)/\varepsilon \gg 1/\varepsilon$ in the variable $i$. Therefore, we choose functions $f_\varepsilon(x,i)$ of the form
	$$
		f_\varepsilon(x,i)
	=
		f(x)
	+
		\varepsilon \, \varphi \left(
		\frac{x}{\varepsilon}
		\right)
	+
		\frac{\varepsilon}{\gamma(\varepsilon)} \, \xi\left(
		\frac{x}{\varepsilon},i
		\right).
	$$
	We abbreviate in the following $y=x/\varepsilon$. Then computing $H_\varepsilon f_\varepsilon$ results in
	\begin{multline*}
		H_\varepsilon f_\varepsilon (x,i)
	=
		\frac{\varepsilon}{2}\Delta f(x)
		+
		\frac{1}{2}\left[\Delta\varphi(y)+\frac{1}{\gamma(\varepsilon)}\Delta\xi^i(y)\right]
		\\+
	\frac{1}{2}
	\big|
	\nabla f(x)+\nabla\varphi(y)+\frac{1}{\gamma(\varepsilon)}\nabla\xi^i(y)
	\big|^2
	+
	b^i(y)\left(\nabla f(x)+\nabla\varphi(y)
	+
	\frac{1}{\gamma(\varepsilon)} \nabla\xi^i(y)\right) 
	\\+
	\gamma(\varepsilon)\sum_{j = 1}^J r_{ij}(y)\left[e^{(\xi(y,j)-\xi(y,i))/\gamma}-1\right].
	\end{multline*}
	The $1/\gamma$ terms vanish as $\gamma\rightarrow\infty$. The last term satisfies
	$$
		\gamma \sum_{j = 1}^J r_{ij}(y) \left[
		e^{(\xi^j-\xi^i)/\gamma}-1
		\right]
	=
		\sum_{j = 1}^Jr_{ij}(y) \left[
		\xi^j(y)-\xi^i(y)
		\right]
	+
		\mathcal{O}(\gamma^{-2}).
	$$
	Therefore, we choose again $E^\prime := \mathbb{T}^d \times \{1,\dots, J\}$ as the state space of the macroscopic variables, and use the following limit operator $H$,
	\begin{align}\label{MM:eq:multi-H-proof-cont-averaged}
		H:=
		\left\{
		(f, H_{f,\varphi,\xi} \, : \,
		f \in C^2(\mathbb{T}^d) \text{ and } 
		H_{f, \varphi,\xi} \in C(\mathbb{T}^d \times E^\prime)
		\right\},
	\end{align}
	with functions $\varphi$ and $\xi$ in the sets $\varphi \in C^2(\mathbb{T}^d)$ and $\xi = (\xi^1, \dots, \xi^J) \in C^2(E^\prime) \simeq (C^2(\mathbb{T}^d))^J$. The image functions $H_{f,\varphi,\xi} : \mathbb{T}^d \times \mathbb{T}^d \times \{1,\dots,J\} \to \mathbb{R}$ are
	\begin{multline}
		H_{f,\varphi,\xi}(x,y,i)
	:=
		\frac{1}{2}\Delta_y\varphi(y)
		+
		\frac{1}{2} |\nabla f(x)+\nabla_y\varphi(y)|^2
		+
		b^i(y)
		\left(\nabla f(x)+\nabla_y\varphi(y)\right) \\
	+
		\sum_{j = 1}^J r_{ij}(y)\left[\xi(y,j)-\xi(y,i)\right].
	\label{eq:LDP_MM:contII:limit_op_H}
	\end{multline}
	Then $H$ satisfies (T1), which is shown by the same line of argument as above in the proof of Theorem \ref{thm:results:LDP_cont_MM_Limit_I}, with the same maps $\eta_\varepsilon$ and $\eta_\varepsilon^\prime$ as there. The image functions depend only on gradients,~$H_{f,\varphi,\xi}(x,y,i)=H_{\varphi,\xi}(\nabla f(x),y,i)$.
	\end{proof}
	\begin{proof}[Verification of (T2) of Theorem \ref{thm:results:LDP_switching_MP}]
	For any $p \in \mathbb{R}^d$, we establish the existence of functions $\varphi_p \in C^2(\mathbb{T}^d)$ and $\xi \in C^2(E^\prime)$ such that $H_{\varphi, \xi}(p,\cdot)$ becomes constant on $E^\prime = \mathbb{T}^d \times \{1,\dots,J\}$. To that end, we find a constant $\mathcal{H}(p) \in \mathbb{R}$ and $\varphi_p$ and $\xi_p$ such that for all $(y,i) \in E^\prime$, we have
	$$
	H_{\varphi_p,\xi_p}(p,y,i) = \mathcal{H}(p).
	$$
	We reduce the problem to finding a principal eigenvalue.
	\begin{lemma}
	Let $E^\prime = \mathbb{T}^d \times \{1,\dots,J\}$ and let~$H$ be the operator~\eqref{MM:eq:multi-H-proof-cont-averaged}. Then:
	\begin{enumerate}[(a)]
	\item For $f \in \mathcal{D}(H)$, the images $H_{\varphi,\xi}$ are given by
	\[
		\tilde{H}_{\varphi,\xi}(p,y,i)
	=
		e^{-\varphi(y)}\left[ (B^i_{p} + V^i_{p}) e^{\varphi}\right](y)
	+
		\sum_{j = 1}^J r_{ij}(y)\left[\xi(y,j)-\xi(y,i)\right],
	\]
	where $p = \nabla f(x) \in \mathbb{R}^d$, $B_p^i = \frac{1}{2}\Delta + (p + b^i(y))\cdot\nabla$ and multiplication operator~$V_p^i(y) = p^2/2 + p\cdot b^i(y)$.
	\item For any $\varphi$ and $y\in \mathbb{T}^d$, there exists a function $\xi(y,\cdot)$ on $\{1,\dots,J\}$ such that $\xi \in C^2(E^\prime)$ and for all $i=1,\dots,J$,
	\[
		e^{-\varphi}\left[ (B^i_{p} + V^i_{p}) e^{\varphi}\right](y)
	+
		\sum_{j = 1}^J r_{ij}(y)\left[\xi(y,j)-\xi(y,i)\right]
	=
		e^{-\varphi(y)} \left[ B_p + V_p \right] e^{\varphi(y)},
	\]
	where $B_p = \frac{1}{2}\Delta + (p + \overline{b}(y))\cdot\nabla$, $V_p(y) = \frac{p^2}{2} + p\cdot\overline{b}(y)$. In the operators, $\overline{b}(y) := \sum_{i = 1}^J \mu_y(i)b^i(y)$ is the average drift with respect to the stationary measure $\mu_y \in \mathcal{P}(\{1, \dots, J\})$ of the jump process with frozen jump rates~$r_{ij}(y)$.
	\item There exists a strictly positive eigenfunction $g_p$ and an eigenvalue $\mathcal{H}(p) \in \mathbb{R}$ such that
	\begin{equation}
	\left[ B_p + V_p \right] g_p = \mathcal{H}(p) g_p.
	\label{eq:LDP_MM:contII:PrEv_eq}
	\end{equation}
	\end{enumerate}
	\label{lemma:LDP_MM:contII:principal_eigenvalue}
	\end{lemma}
	By (a), (b) and (c), taking $\varphi_p = \log g_p$ and the corresponding $\xi(y,i)$, we obtain~(T2) via
	\begin{align*}
	H_{\varphi_p,\xi}(p,y,i)
	&\overset{(a)}{=}
	e^{-\varphi_p(y)}\left[ B^i_{p} + V^i_{p}\right] e^{\varphi_p(y)}
	+
	\sum_{j\in\mathcal{J}}r_{ij}(y)\left[\xi(y,j)-\xi(y,i)\right]
	\\&\overset{(b)}{=}
	e^{-\varphi_p(y)} \left[ (B_p + V_p)e^{\varphi} \right] (y)
	\overset{(c)}{=}
	\mathcal{H}(p).
	\end{align*}
	\emph{Proof of Lemma \ref{lemma:LDP_MM:contII:principal_eigenvalue}.} Regarding (a), writing $\xi(y,i)=\xi_y(i)$ and $p = \nabla f(x)\in\mathbb{R}^d$, for all $(y,i)\in E^\prime$ we find
	\begin{align*}
		H_{\varphi,\xi}(p,y,i)
	&=
		\underbrace{\frac{1}{2}\Delta\varphi + \frac{1}{2} \big|p + \nabla\varphi \big|^2 + b^i(p+\nabla\varphi)}_{\displaystyle=e^{-\varphi}(B_{p,i}+V_{p,i})e^{\varphi}} 
	+
	\underbrace{\sum_{j = 1}^J r_{ij}(y)
	\left[
	\xi(y,j)-\xi(y,i)
	\right]}_{\displaystyle =:R_y \xi(y,\cdot)(i)},
	\end{align*}			
	with a generator $R_y$ of a jump process with frozen jump rates~$r_{ij}(y)$.
	\smallskip
	
	For (b), let $\varphi \in C^2(\mathbb{T}^d)$ and $y\in \mathbb{T}^d$. We wish to find a function~$\xi_y(\cdot) = \xi(y,\cdot) \in C(\{1,\dots,J\})$ such that
	$$
		e^{-\varphi}\left[B_{p,i}+V_{p,i}\right]e^{\varphi} + R_y \xi_y(i)
	$$
	becomes constant in $i = 1,\dots, J$. By the Fredholm alternative, for any vector~$h\in C(\{1,\dots,J\})$, the equation $R_y\xi_y = h$ has a solution $\xi_y(\cdot)\in C(\{1,\dots,J\})$ if and only if $h \perp \text{ker}(R_y^\ast)$. Since $R_y$ is the generator of a jump process on the finite discrete set $\{1, \dots, J\}$ with rates $r_{ij}(y)$, the null space $\text{ker}(R_y^\ast)$ is one-dimensional and spanned by the unique stationary measure $\mu_y\in\mathcal{P}(\{1,\dots,J\})$, which exists by our irreducibility assumption of Theorem~\ref{thm:results:LDP_cont_MM_Limit_II} (e.g.~\cite[Theorem~17.51]{klenke2013probability}).
	Hence
	$
		e^{-\varphi}\left[B_{p,i}+V_{p,i}\right]e^{\varphi} + R_y \xi_y(i)
	=
		h(p,y)
	$
	is independent of $i \in \{1,\dots,J\}$ iff
	\[
	\sum_{i = 1}^J \mu_y(i)\left[(h(p,y)-e^{-\varphi}\left[B_{p,i}+V_{p,i}\right]e^{\varphi}\right]=0.
	\]
	This solvability condition leads to
	\begin{align*}
		\sum_{i = 1}^J \mu_y(i)
		\left[(h(p,y)
		-
		e^{-\varphi}
		\left( B_{p,i}+V_{p,i}\right)
		e^{\varphi}
		\right]
	&=
		h(p,y)-
		e^{-\varphi(y)}\left(
		B_p + V_p
		\right)
		e^{\varphi(y)}
	=
		0.
	\end{align*}
	Hence for $h(p,y) := e^{-\varphi(y)}\left[ B_p + V_p \right]e^{\varphi(y)}$, there exists~$\xi(y,i)$ solving the equation~$R_y \xi(y,\cdot) = h$. Furthermore, since the stationary measure is an eigenvector of a one-dimensional eigenspace, and the rates $r_{ij}(\cdot)$ are smooth by assumption, the eigenfunctions $\xi_y$ depend smoothly on $y$ as well, and (b) follows.
	\smallskip
	
	Regarding (c) in Lemma \ref{lemma:LDP_MM:contII:principal_eigenvalue}, note that \eqref{eq:LDP_MM:contII:PrEv_eq} is a principal eigenvalue problem for a second-order uniformly elliptic operator. By Proposition \ref{proposition:appendix:PrEv:elliptic_op}, the principal eigenvalue problem $\left[ - B_p - V_p \right] g _p = \lambda(p) g_p$ has a solution $g_p > 0$, with eigenvalue $\lambda(p) \in \mathbb{R}$. The same function $g_p$ and the eigenvalue $\mathcal{H}(p) = -\lambda(p)$ solve~\eqref{eq:LDP_MM:contII:PrEv_eq}.
	\end{proof}
	\begin{proof}[Verification of (T3) of Theorem \ref{thm:results:action_integral_representation}]
	The principal eigenvalue $\mathcal{H}(p)$ is of the form
	$$
		\mathcal{H}(p)
	=
		\inf_\varphi\sup_{y\in \mathbb{T}^d} F\left(p,\varphi\right)(y),
	$$
	with $F$ jointly convex in $p$ and $\varphi$. Convexity of~$\mathcal{H}(p)$ and~$\mathcal{H}(0) = 0$ follow as above in the proof of Theorem~\ref{thm:results:LDP_cont_MM_Limit_I}.
	\end{proof}
	\subsection{Proof for the discrete models}
	\label{subsection:discr_MM}
	In this section, we prove the large-deviation Theorems \ref{thm:results:LDP_discr_MM_Limit_I} and \ref{thm:results:LDP_discr_MM_Limit_II} of the stochastic process $(X^n_t, I^n_t)$ from Defintion~\ref{def:intro:discr_MM}. We alert the reader that we use~$n$ as a scaling parameter instead of~$\varepsilon$. The generator $L_n$ in~\eqref{eq:intro:L_N_discr_MM} is
	\begin{multline*}
		L_n f(x,i)
	=
		n r_+^i(nx) \left[ f(x + 1/n, i) - f(x,i) \right]
	+
		n r_-^i(nx) \left[ f(x - 1/n, i) - f(x,i) \right]	\\
	+
		\sum_{j = 1}^J n \gamma(n) r_{ij}(nx) \left[ f(x,j) - f(x,i) \right],
	\end{multline*}
	with~$\gamma(n)>0$, the state space $E_n = \mathbb{T}_{\ell,n} \times \{1,\dots, J\} = \{(x,i)\}$, and $\mathbb{T}_{\ell,n}$ the discrete one-dimensional torus with lattice spacing $1/n$ and of length $\ell$. As in the continuous case, we verify~\ref{MM:item:T1},~\ref{MM:item:T2} and~(T3) of the large-deviation Theorems \ref{thm:results:LDP_switching_MP} and \ref{thm:results:action_integral_representation}. We start from the nonlinear generators~$H_n f = n ^{-nf} L_n e^{nf}$,
	\begin{multline}
		H_n f(x,i)
	=
		r^i_{+}(nx)\left[e^{n\left(f(x+1/n,i)-f(x,i)\right)}-1\right] + r^i_{-}(nx)\left[e^{n\left(f(x-1/n,i)-f(x,i)\right)}-1\right] \\
	+
		\gamma(n)\sum_{j = 1}^J r_{ij}(nx)\left[e^{n\left(f(x,j)-f(x,i)\right)}-1\right].
	\label{eq:LDP_MM:discr:H_N}
	\end{multline}
	\subsubsection{Proof of Theorem \ref{thm:results:LDP_discr_MM_Limit_I}}
	\label{subsubsection:LDP_discr_MM_Limit_I}	
	\begin{proof}[Verification of (T1) of Theorem \ref{thm:results:LDP_switching_MP}] 
	We have~$\gamma\equiv 1$.
	Choose functions of the form $f_n(x,i) = f(x) + \frac{1}{n}\varphi(nx,i)$, where $(x,i) \in \mathbb{T}_{\ell, n} \times \{1,\dots,J\}$ and the~$\varphi(\cdot,i) \in C_{\ell\text{-per}}(\mathbb{Z}) \simeq C(\mathbb{T}_{\ell,1})$ are $\ell$-periodic functions. Then we obtain
	\begin{multline*}
		H_nf_n(x,i)
	=
		r^i_{+}(nx)\left[e^{n\left(f(x+\frac{1}{n})-f(x)\right)}e^{\varphi^i(nx+1)-\varphi^i(nx)}-1\right]
	\\+
		r^i_{-}(nx)\left[e^{n\left(f(x-\frac{1}{n})-f(x)\right)}e^{\varphi^i(nx-1)-\varphi^i(nx)}-1\right] \\
	+
		\sum_{j = 1}^J r_{ij}(nx)\left[e^{\left(\varphi(nx,j)-\varphi(nx,i)\right)}-1\right].
	\end{multline*}
	Then $H_n f_n$ depends on the variables~$x\in\mathbb{T}_{\ell, n}$,~$nx\in\mathbb{T}_{\ell,1}$ and~$i \in \{1,\dots,J\}$. Therefore choose $E^\prime = \mathbb{T}_{\ell,1} \times \{1,\dots,J\}$ for the macroscopic variables, and set
	\begin{multline*}
		H:=
		\bigg\{
		(f, H_{f,\varphi)} \, : \,
		f \in C^1(\mathbb{T}_\ell) \text{ and } 
		H_{f, \varphi} \in C(\mathbb{T}_{\ell} \times E^\prime),\\
		\varphi = (\varphi^1, \dots, \varphi^J) \in C(E^\prime) \simeq (C(\mathbb{T}_{\ell, 1}))^J
		\bigg\}.
	\end{multline*}
	The image funcitons $H_{f,\varphi} : \mathbb{T}_\ell \times E^\prime \to \mathbb{R}$ are defined as
	\begin{multline}
		H_{f,\varphi}(x,y,i)
	:=
		r^i_{+}(y)
		\left[
		e^{\partial_x f(x)}e^{\varphi^i(y+1)-\varphi^i(y)}-1
		\right]
	\\+
		r^i_{-}(y)
		\left[
		e^{-\partial_x f(x)}e^{\varphi^i(y-1)-\varphi^i(y)}-1
		\right]	\\
	+
		\sum_{j = 1}^J r_{ij}(y)\left[e^{\varphi(y,j)-\varphi(y,i)}-1\right].
	\label{eq:LDP_MM:discrI:limit_op_H}
	\end{multline}
	Then with the embedding $\eta_n': E_n \to E^\prime, (x,i) \mapsto \eta_n'(x,i):=(nx,i)$, and the projection $\eta_n (x,i) = x$, (C1) is satisfied. Regarding (C2), for $(f,H_{f,\varphi}) \in H$, the function $f_n(x,i):= f(x)+\frac{1}{n}\varphi(nx,i)$ satisfies $f_n\rightarrow f$ uniformly in $(x,i)\in E_n$ with respect to $\eta_n$, using that $\sup_{E^\prime}\varphi<\infty$. Regarding the images, we note that
	\begin{multline*}
		\sup_{x,i}
		\big|
		H_{f,\varphi}(x,nx,i) - H_n f_n (x,i)
		\big|
	\\\leq
		\sup_{x,i}\bigg|
		r_+^i(nx) \left[ e^{\partial_x f(x)} - e^{n\left(f(x+1/n)-f(x)\right)}\right] 
		e^{\varphi(nx+1,i) - \varphi(nx,i)}
		\bigg| \\
	+
		\sup_{x,i}\bigg|
		r_-^i(nx) \left[ e^{-\partial_x f(x)} - e^{n\left(f(x-1/n)-f(x)\right)}\right] 
		e^{\varphi(nx-1,i) - \varphi(nx,i)}
		\bigg|,
	\end{multline*}
	which converges to zero as $n$ goes to infinity, since $\sup_{E^\prime}\varphi < \infty$ and we have uniformly bounded jump rates $r_\pm^i$. Furthermore, the images depend on $x$ only via the derivatives of $f$: $H_{f,\varphi}(x,y,i) = H_\varphi(\partial_x f(x),y,i)$. Hence (C3) is satisfied, and this finishes the verification of (T1).
	\end{proof}
	\begin{proof}[Verification of (T2) of Theorem \ref{thm:results:LDP_switching_MP}]
	For $p \in \mathbb{R}$, we want to find a function $\varphi_p$ such that the images $H_{\varphi}(p,y,i)$ become constant in $(y,i)$. As in the continuous case, this can be achieved by solving a principal eigenvalue problem. Here in the discrete case, instead of elliptic partial differential equations, we encounter principal eigenvalues of irreducible M-matrices.
	\begin{lemma}\label{lemma:LDP_MM:discrI:principal_ev}
	Let $E^\prime = \mathbb{T}_{\ell,1} \times \{1,\dots,J\}$ and $H \subseteq C^1(\mathbb{T}_\ell) \times C(\mathbb{T}_\ell \times E^\prime)$ be the multivalued operator from \eqref{eq:LDP_MM:discrI:limit_op_H}, and let $p \in \mathbb{R}$. Then:
	\begin{enumerate}[(a)]
	\item Writing $g(y,i) = g^i(y) :=  e^{\varphi^i(y)}$, the images $H_\varphi(p,y,i)$ are of the form
	\[
		\tilde{H}_\varphi(p,y,i)
	=
		\frac{1}{g(y,i)}
		\left[
		B_p
	+
		R
		\right]g (y,i),
	\]
	where 
	\[
		B_p g(y,i) 
	: =
		r^i_+(y)\left[e^pg^i(y+1)-g^i(y)\right]
	+
		r^i_-(y)\left[e^{-p}g^i(y-1)-g^i(y)\right],
	\]
	and
	\[
		R g (y,i)
	:=
		\sum_{j = 1}^J r_{ij}(y) 
		\left[
		g(y,j)-g(y,i)
		\right].
	\]
	\item There exist strictly positive vectors $g^i = \left( g^i(0),\dots,g^i(\ell-1)\right) \in \mathbb{R}^\ell$, $g^i(y) > 0$ for all $i = 1,\dots, J$ and $y = 0, \dots, \ell-1$, and an eigenvalue $\mathcal{H}(p) \in \mathbb{R}$ such that
	\[
		[B_p + R] g(y,i)
	=
		\mathcal{H}(p) g(y,i).
	\]
	\end{enumerate}
	\end{lemma}
	By (a) and (b), choosing $\varphi(y,i) := \log g(y,i)$, we obtain
	$$
	H_\varphi(p,y,i)
	\overset{(a)}{=}
	\frac{1}{g(y,i)}
	\left[
	B_p
	+
	R
	\right]g (y,i)
	\overset{(b)}{=}
	\mathcal{H}(p).
	$$
	\emph{Proof of Lemma \ref{lemma:LDP_MM:discrI:principal_ev}.} Part (a) follows from rewriting the images $H_\varphi(p,y,i)$. Regarding (b), when cast in matrix form, the eigenvalue problem reads
	\begin{align*}
	\left[
	\begin{pmatrix}
		B_p^1	&	&	0	\\
			&	\ddots	&	\\
		0	&	&	B_p^J		
	\end{pmatrix}
	+
	\begin{pmatrix}
		R_{11}	&	&	\geq 0	\\
    		& \ddots	&	\\
    	\geq 0	&		&	R_{JJ}
	\end{pmatrix}
	\right]
	\begin{pmatrix}
		g^1	\\
		\vdots	\\
		g^J
	\end{pmatrix}
	=
	\mathcal{H}(p)
	\begin{pmatrix}
		g^1	\\
		\vdots	\\
		g^J
	\end{pmatrix},
	\end{align*}
	where each $g^i$ is a vector, $g^i = \left( g^i(0), \dots, g^i(\ell-1) \right) \in \mathbb{R}^\ell$, and the square matrices $B^i_p \in \mathbb{R}^{\ell \times \ell}$ are similar to a discretized Laplacian with periodic boundaries.
	More precisely, the matrix $B_p^i$ has entries $-(r_+^i(y) + r_-^i(y))$ on the diagonal that are flanked by $r_+^i(y)e^p$ to the right and $r_-^i(y) e^{-p}$ to the left next entries.
	%
	Each~$R_{ii}$ is a diagonal matrix with~$(R_{ii})_{kk}=-\sum_{j\neq i}r_{ij}(k)$, where~$k=0,1,\dots,\ell-1$. 
	The remaining block matrices in~$R$ are non-negative and mix the different component vectors $g^i$ and $g^j$.
	\smallskip
	
	By the irreducibility assumption in Theorem~\ref{thm:results:LDP_discr_MM_Limit_I} on $R$. Since all off-diagonal terms in $B_p + R$ are non-negative, the off-diagonal elements form an irreducible matrix. Therefore, $M_p := -B_p - R$ is an irreducible M-matrix (Definition~\ref{def:appendix:PrEv:irreducible_M_matrix} further below), and by Proposition~\ref{proposition:appendix:PrEv:M_matrix}, it admits a principal eigenvalue $\lambda(p)$ with strictly positive eigenvector~$g_p$, that is~$M_p g_p = \lambda (p) g_p$. 
	Consequently, we find~$\left[ B_p + R \right] g_p = \mathcal{H}(p) g_p$ with the same eigenvector~$g_p$ and principal eigenvalue~$\mathcal{H}(p) = -\lambda(p)$. This finishes the verification of~\ref{MM:item:T2}.
	\end{proof}
	\begin{proof}[Verification of (T3) of Theorem \ref{thm:results:action_integral_representation}]
	By Proposition~\ref{proposition:appendix:PrEv:M_matrix}, the eigenvalue satisfies
	\begin{multline*}
		\mathcal{H}(p)
	=
		- \sup_{g>0} \inf_{y,i} 
		\left[
		\frac{1}{g(y,i)} \left( -B_p - R \right)g(y,i)
		\right]
	\\=
		\inf_{\varphi} \sup_{y,i}
		\bigg\{
		r^i_{+}(y)
		\left[
		e^{p}e^{\varphi^i(y+1)-\varphi^i(y)}-1
		\right]
	+
		r^i_{-}(y)
		\left[
		e^{-p}e^{\varphi^i(y-1)-\varphi^i(y)}-1
		\right] 
	\\+
		\sum_{j = 1}^J r_{ij}(y)
		\left[
		e^{\varphi(y,j)-\varphi(y,i)}-1
		\right]
		\bigg\}.
	\end{multline*}
	Hence the eigenvalue is of the form 
	\[
		\mathcal{H}(p) = \inf_{\varphi} \sup_{y,i} F(p, \varphi)(y,i),
	\]
	with $F(p,\varphi)$ jointly convex in $p$ and $\varphi$, and convexity of $\mathcal{H}(p)$ follows as demonstrated in the proof of Theorem~\ref{thm:results:LDP_cont_MM_Limit_I}. Choosing the constant vector $\varphi = (1,\dots, 1)$ in the variational representation, we obtain $\mathcal{H}(0) \leq 0$. Conversely, any $\varphi$ admits a global minimum $(y_m,i_m)$. We have the estimate~$ F(0,\varphi)(y_m,i_m) \geq 0$. Therefore,~$\sup F(0,\varphi)\geq0$ for any $\varphi$, and~$\mathcal{H}(0) \geq 0$ follows.
%
	\end{proof}
	\subsubsection{Proof of Theorem \ref{thm:results:LDP_discr_MM_Limit_II}}
	\label{subsubsection:LDP_discr_MM_Limit_II}	
	\begin{proof}[Verification of (T1) of Theorem \ref{thm:results:LDP_switching_MP}]
	We have~$\gamma(n)\to\infty$ as~$n\to\infty$.
	With functions of the form $f_n(x,i) = f(x) + \frac{1}{n} \varphi(nx) + \frac{1}{n\gamma(n)} \xi(nx,i)$, with functions~$\varphi$ and~$\xi(\cdot,i)$ in $C(\mathbb{T}_{\ell,1})$, we obtain
	\begin{multline*}
		H_nf_n(x,i)
	=
		r^i_{+}(nx)\left[e^{n\left(f(x+\frac{1}{n})-f(x)\right)}e^{\varphi(nx+1)-\varphi(nx)} e^{\left( \xi(nx + 1,i) - \xi(nx,i) \right) / \gamma(n)}-1\right] \\
	+
		r^i_{-}(nx)\left[e^{n\left(f(x-\frac{1}{n})-f(x)\right)}e^{\varphi(nx-1)-\varphi(nx)} e^{\left( \xi(nx - 1,i) - \xi(nx,i) \right) / \gamma(n)} -1\right] \\
	+
		\gamma(n)\sum_{j = 1}^J r_{ij}(nx)\left[e^{\left(\xi(nx,j)-\xi(nx,i)\right) / \gamma(n)}-1\right].
	\end{multline*}
	Take $E^\prime := \mathbb{T}_{\ell,1} \times \{1,\dots,J\}$ and set
	\begin{multline*}
		H:=
		\bigg\{
		(f, H_{f,\varphi,\xi}) \, : \,
		f \in C^1(\mathbb{T}_\ell) \text{ and } 
		H_{f, \varphi, \xi} \in C(\mathbb{T}_\ell \times E^\prime),
		\\ \varphi \in C(\mathbb{T}_{\ell,1}), 
		\; \xi = (\xi^1, \dots, \xi^J) \in C(E^\prime) \simeq (C(\mathbb{T}_{\ell,1}))^J
		\bigg\},
	\end{multline*}
	with image functions $H_{f,\varphi} : \mathbb{T}_\ell \times E^\prime \to \mathbb{R}$ defined by
	\begin{multline}
		H_{f,\varphi}(x,y,i)
	:=
		r^i_{+}(y)\left[e^{\partial_x f(x)}e^{\varphi(y+1)-\varphi(y)}-1\right]+r^i_{-}(y)\left[e^{-\partial_x f(x)}e^{\varphi(y-1)-\varphi(y)}-1\right]\\
	+
		\sum_{j = 1}^J r_{ij}(y)\left[\xi(y,j) - \xi(y,i)\right].
	\label{eq:LDP_MM:discrII:limit_op_H}
	\end{multline}
	Then with the embedding $\eta_n': E_n \to E^\prime, (x,i)\mapsto\eta_n^\prime(x,i):=(nx,i)$, and the projection $\eta_n (x,i) = x$, Item~(C1) is satisfied. Regarding Item~(C2), consider a pair~$(f,H_{f,\varphi,\xi}) \in H$. The function $f_n(x,i):=f(x)+\frac{1}{n}\varphi(nx) + \frac{1}{n\gamma(n)} \xi(nx,i)$ satisfies~$f_n\rightarrow f$ uniformly in $(x,i)$ with respect to $\eta_n$. For the convergence of images, use 
	\[
		\sup_{x,i}
		\bigg|
		\left(\xi(nx + 1,i) - \xi(nx,i)\right)/\gamma(n)
		\bigg|
	\leq
		\frac{1}{\gamma(n)}
		\sup_{y,i}
		\bigg|
		\xi(y + 1,i) - \xi(y,i)
		\bigg|
		\xrightarrow{n\to\infty} 0.
	\]
	Expanding the exponential terms in $H_n f_n$ and using the same uniform bounds lead to the claimed convergence. Finally, (C3) is satisfied, since the images \eqref{eq:LDP_MM:discrII:limit_op_H} depend on $x$ only via derivatives of $f$.
	\end{proof}
	\begin{proof}[Verification of (T2) of Theorem \ref{thm:results:LDP_switching_MP}]	
	For any $p \in \mathbb{R}$, we wish to obtain functions $\varphi \in C(\mathbb{T}_{\ell,1})$ and $\xi \in C(E^\prime)$ such that the images $H_{\varphi,\xi}(p,y,i)$ are constant in $(y,i)$. We reduce that to a principal-eigenvalue problem.
	\begin{lemma}\label{lemma:LDP_MM:discrII:principal_ev}
	Let $E^\prime = \mathbb{T}_{\ell,1} \times \{1,\dots,J\}$, $p\in\mathbb{R}$, and let $H \subseteq C^1(\mathbb{T}_\ell) \times C(\mathbb{T}_\ell \times E^\prime)$ be the multivalued operator from \eqref{eq:LDP_MM:discrII:limit_op_H}. Then:
	\begin{enumerate}[(a)]
	\item The images $H_{\varphi,\xi}(p,y,i)$ are of the form
	\[
		\tilde{H}_{\varphi,\xi}(p,y,i)
	=
		\frac{1}{g(y)}
		B^i_p g(y)
	+
		R \xi(y,i),
	\]
	where~$g(y) :=  e^{\varphi(y)}$,
	\[
		B^i_p g(y) 
	: =
		r^i_+(y)\left[e^pg(y+1)-g(y)\right]
	+
		r^i_-(y)\left[e^{-p}g(y-1)-g(y)\right]
	\]
	and
	\[
	R \xi (y,i)
	:=
		\sum_{j = 1}^J r_{ij}(y) 
		\left[
		\xi(y,j)-\xi(y,i)
		\right].
	\]
	\item For any $g(y) = e^{\varphi(y)}$ and $y\in \mathbb{T}_{\ell,1} \simeq \{0,1,\dots,\ell-1\}$, there exists a function $\xi_p(y,\cdot) \in C(\{1,\dots,J\})$ such that for all $i \in \{1,\dots,J\}$,
	\[
		\frac{1}{g(y)} B^i_{p} g(y)
	+
		R\xi(y,i)
	=
		\frac{1}{g(y)} B_p g(y),
	\]
	with 
	\[
		B_p g(y) 
	: =
		\overline{r}_+(y)\left[e^pg(y+1)-g(y)\right]
	+
		\overline{r}_-(y)\left[e^{-p}g(y-1)-g(y)\right],
	\]
	where 
	$\overline{r}_{\pm}(y) = \sum_{i = 1}^J \mu_y(i) r_{\pm}^i(y)$ are the average jump rates with respect to the stationary measure $\mu_y \in \mathcal{P}(\{1,\dots,J\})$ of the jump process with rates~$r_{ij}(y)$.
	\item There exists a strictly positive eigenvector $g_p = \left( g_p(0), \dots, g_p(\ell-1) \right) \in \mathbb{R}^\ell$, $g_p(y) > 0$ for all $y = 0, \dots, \ell-1$, and a corresponding principal eigenvalue $\mathcal{H}(p) \in \mathbb{R}$ such that
	$$
	B_p g_p = \mathcal{H}(p) g_p.
	$$
	\end{enumerate}	
	\end{lemma}
	With $\varphi_p := \log g_p$ and the corresponding function~$\xi_p(y,i)$ of~(b), we find
	$$
	H_{\varphi_p,\xi_p}(p,y,i)
	\overset{(a)}{=}
	\frac{1}{g(y)}
	B^i_p g(y)
	+
	R \xi_p(y,i)
	\overset{(b)}{=}
	\frac{1}{g_p(y)} B_p g_p(y)
	\overset{(c)}{=}
	\mathcal{H}(p).
	$$
	\emph{Proof of Lemma \ref{lemma:LDP_MM:discrII:principal_ev}.} Part (a) follows from rewriting the images in terms of $g(y) = \log\varphi(y)$. For part (b), the argument is similar to the one given in the proof of Theorem \ref{thm:results:LDP_cont_MM_Limit_II}. By the Fredholm alternative, for every $y \in \mathbb{T}_{\ell,1} \simeq \{0,\dots,\ell-1\}$, the equation
	$$
		R\xi(y,i) = \frac{1}{g(y)} \left[B^i_p - B_p\right]g (y)
	$$
	has a solution $\xi(y,\cdot) \in C(\{1,\dots,J\})$ if and only if for the stationary measure $\mu_y \in \mathcal{P}(\{1,\dots,J\})$ satisfying $R^\ast \mu_y = 0$, we have
	$$
		\bigg\langle \mu_y 
	,
		\frac{1}{g(y)} \left[ B^i_p - B_p \right]g (y)
		\bigg\rangle
	=
		0,
	$$
	where the pairing corresponds to a sum over the $i \in \{1,\dots,J\}$. Writing out that condition leads exactly to the average operator $B_p$ as given in (b).
	\smallskip
	
	For part (c), we note that $B_p g_p = \mathcal{H}(p) g_p$ is a matrix eigenvalue problem. The matrix $B_p \in \mathbb{R}^{\ell \times \ell}$ has nonzero entries similar to a discretized Laplacian with periodic boundaries:
	\begin{align*}
		B_p 
	=
		\begin{pmatrix}
			-(\overline{r}_+(0) + \overline{r}_-(0))	&	\overline{r}_+(0) e^p	&	\dots	&	\overline{r}_-(0)e^{-p} \\
			\overline{r}_-(1) e^{-p}	&	-(\overline{r}_+(1) + \overline{r}_-(1))	&	\overline{r}_+(1)e^p	&	\dots \\
			0	&	\overline{r}_-(2)e^{-p}	&	-(\overline{r}_+(2) + \overline{r}_-(2))	&	\dots \\
			\vdots	&	\vdots	&	\ddots	&	\vdots\\	
		\end{pmatrix}.
	\end{align*}
	By the positivity assumptions on the rates $r_\pm^i$ in Definition \ref{def:intro:discr_MM}, the average rates $\overline{r}_\pm$ are positive. Thereby, $M_p := -B_p$ is an irreducible $M$-matrix, so that by Proposition \ref{proposition:appendix:PrEv:M_matrix}, there exists a strictly positive eigenvector $g_p > 0$ and a principal eigenvalue $\lambda(p) \in \mathbb{R}$ such that
	$
	M_p g_p = \lambda(p) g_p.
	$	
	That implies $B_p g_p = \mathcal{H}(p) g_p$ with the same eigenvector $g_p$ and principal eigenvalue $\mathcal{H}(p) = -\lambda(p)$. This finishes the proof of Lemma \ref{lemma:LDP_MM:discrII:principal_ev}, and thereby the verification of (T2).
	\end{proof}
	\begin{proof}[Verification of (T3) of Theorem \ref{thm:results:action_integral_representation}]
	We prove the claimed properties of~$\mathcal{H}(p)$ by means of a variational representation. By Proposition~\ref{proposition:appendix:PrEv:M_matrix}, we have
	\begin{align*}
		\mathcal{H}(p)
	&=
		- \sup_{g>0} \inf_{y} 
		\left[
		\frac{1}{g(y)} (-B_p) g(y)
		\right] \\
	&=
		\inf_{\varphi} \sup_{y}
		\bigg\{
		\overline{r}_+(y)\left[e^p e^{\varphi(y+1)-\varphi(y)} - 1\right]
	+
		\overline{r}_-(y)\left[e^{-p} e^{\varphi(y-1)-\varphi(y)} - 1\right]
		\bigg\}.
	\end{align*}
	The representation is of the form 
	\[
		\mathcal{H}(p) = \inf_{\varphi} \sup_{y} F(p, \varphi(y)),
	\]
	with a joint convex $F$. With that, convexity and~$\mathcal{H}(0)=0$ follow as above.
	\end{proof}
	\subsection{Detailed balance implies symmetry of the Hamiltonian}
	\label{subsection:detailed_balance}
	In Theorem~\ref{thm:results:detailed_balance_limit_I}, we proved that detailed-balance implies symmetric Hamiltonians. The proof was based on a suitable variational representation of the Hamiltonian. In this section, we show how to obtain this representation.
	\smallskip
	
	To that end, we recall the setting. We work with~$E^\prime = \mathbb{T}^d \times \{1, \dots, J\}$, and denote by~$\mathcal{P}(E^\prime)$ the set of probability measures on $E^\prime$. The Hamiltonian~$\mathcal{H}(p)$ is the principal eigenvalue of the cell problem~\eqref{eq:LDP_MM:contI:cell_problem} described in Lemma~\ref{lemma:LDP_MM:contI:principal_eigenvalue}. Here, we start from the fact that this Hamiltonian satisfies
	\begin{equation}
	\mathcal{H}(p)
	=
	\sup_{\mu\in\mathcal{P}(E^\prime)}
	\left[
	\int_{E^\prime} V_p(z) \,\dd \mu(z) - I_p(\mu)
	\right].
	\label{eq:results:LDP_MM:DV_var_rep_H(p)}
	\end{equation}
	In this formula, we have the continuous map
	\begin{equation}\label{MM:eq:function-V-in-Hamiltonian}
	V_p (x,i) := \frac{1}{2} p^2 - p\cdot \nabla\psi^i(x),
	\end{equation}
	and the Donsker-Varadhan functional
	\begin{equation}
	I_p(\mu)=-\inf_{u > 0}\int_{E^\prime} \frac{L_p u}{u}\,\dd \mu,
	\label{eq:results:LDP_MM:DV_functional}
	\end{equation}
	where the infimum is over strictly positive~$u\in C^2(E')$ and the operator~$L_p$ is
	\begin{equation}
	L_p u(x,i) 
	:=
	\frac{1}{2} \Delta_x u(x,i) + (p-\nabla \psi_i(x))\cdot \nabla_x u(x,i)
	+
	\sum_{j = 1}^J r_{ij}(x) \left[u(x,j) - u(x,i)\right].
	\label{eq:results:LDP_MM:L_p_in_DV_functional}
	\end{equation}
	The variational representation~\eqref{eq:results:LDP_MM:DV_var_rep_H(p)} is a special case of Donsker's and Varadhan's results on principal eigenvalues~\cite{DonskerVaradhan75, DonskerVaradhan76}.
	Under their general conditions, the infimum is taken over functions that are in the domain of the infinitesimal generator of the semigroup generated by~$L_p$. Pinsky showed that the infimum can be taken over~$C^2$ functions if the coefficient functions appearing in the operators are smooth~\cite{pinsky1985evaluating, Pi07}.
	\smallskip
	
	Given in the form~\eqref{eq:results:LDP_MM:DV_var_rep_H(p)}, it is not clear why~$\mathcal{H}(p)$ should be symmetric under the detailed-balance condition. We perform a suitable shift in the infimum of the Donsker-Varadhan functional~\eqref{eq:results:LDP_MM:DV_functional} to obtain a suitable representation. Let us first briefly describe this transformation in an informal way. Representing in~\eqref{eq:results:LDP_MM:DV_functional} the strictly positive functions as~$u=\exp(\varphi)$, we find
	\begin{equation*}
	I_p(\mu)=-\inf_{\varphi} \sum_i \int\left[
	\frac{1}{2}\Delta\varphi_i + \frac{1}{2}|\nabla\varphi_i|^2
	+(p-\nabla\psi_i) \nabla\varphi_i+\sum_j r_{ij}\left(e^{\varphi_j-\varphi_j}-1\right)\right]\dd\mu_i.
	\end{equation*}
	Suppose that~$\dd\mu_i=\overline{\mu}_i\,\dd x$ with strictly positive~$\overline{\mu}_i$, where~$\dd x$ is the Lebesgue measure on the torus. Then shifting in the infimum as~$\varphi_i\to\varphi_i+\psi_i+\frac{1}{2}\log\overline{\mu}_i$, we find by calculation that
	\begin{equation}
	\label{MM:eq:DV-functional-after-shift}
	I_p(\mu)=\mathcal{R}(\mu)+\int_{E^\prime} V_p \, \dd\mu - K_p(\mu),
	\end{equation}
	where~$\mathcal{R}(\mu)$ is the Fisher information given by
	\begin{equation}
	\label{MM:eq:relative-Fisher-information}
	\mathcal{R}(\mu) := \frac{1}{8}\sum_i\int_{\mathbb{T}^d}\left|\nabla \left(\log\frac{\overline{\mu}_i}{e^{-2\psi_i}}\right)\right|^2\,\dd\mu_i,
	\end{equation}
	and~$K_p(\mu)$ is given by
	\begin{multline}
	K_p(\mu) = \inf_{\phi}\bigg\{\sum_{i = 1}^J \int_{\mathbb{T}^d}\left(
	\frac{1}{2}|\nabla \phi_i(x) + p|^2 
	-\sum_{j = 1}^J r_{ij}(x)\right) \,\dd \mu_i(x)\\
	+\sum_{i, j = 1}^J \int_{\mathbb{T}^d}r_{ij}(x) e^{-2\psi_i(x)}\sqrt{\overline{\mu}_i(x) \overline{\mu}_j(x)}
	e^{\psi_j(x) + \psi_i(x)} e^{\phi(x,j) - \phi(x,i)}
	\,\dd x\bigg\}.
	\label{eq:results:LDP_MM:K_p(mu)_I}
	\end{multline}
	Plugging formula~\eqref{MM:eq:DV-functional-after-shift} into the variational representation~\eqref{eq:results:LDP_MM:DV_var_rep_H(p)} leads to the desired representation of the Hamiltonian. The transformation we used corresponds to shifting by~$(1/2)\log(\overline{\mu}_i/\pi_i)$, where~$\pi_i=e^{-2\psi_i}$ is the stationary measure (up to multiplicative constant). This transformation is actually reminiscent of a \emph{symmetrization} discussed in Touchette's notes~\cite[Eq.~(36)]{Touchette2018}. 
	Finally, when formulating the detailed-balance condition with additional constants in~\eqref{MM:eq:detailed-balance} (meaning not shifting the potentials by constants to renormalized), one can include these constants in the shift to arrive at the same conclusions.
	\smallskip
	
	To make the strategy displayed above rigorous, we prove that we can work with measures~$\mu$ having the required regularity properties. The central idea is to exploit the fact that~$I_p(\mu)$ is finite since~$\mathcal{H}(p)$ is finite. By a result of Stroock~\cite[Theorem 7.44]{stroock2012introduction}, finiteness of the Donsker-Varadhan functional implies certain regularity properties in case the generator is reversible. Since the generator~$L_p$ is not reversible, we instead bound~$I_p$ by a suitable Donsker-Varadhan functional~$I_\mathrm{rev}$ corresonding to a reversible process, and can then apply~\cite[Theorem 7.44]{stroock2012introduction}. This strategy appears in the proof of the following proposition. The formula we use in the proof of Theorem~\ref{thm:results:detailed_balance_limit_I} of Section~\ref{subsection:results:LDP_in_MM} is given in~(c).
	\begin{proposition}
	The Hamiltonian~$\mathcal{H}(p)$ given by~\eqref{eq:results:LDP_MM:DV_var_rep_H(p)} satisfies the following:
	\begin{enumerate}[(a)]
	\item
	The supremum in \eqref{eq:results:LDP_MM:DV_var_rep_H(p)} can be taken over a smaller set $\mathbf{P}$ of measures, that is
	$$
		\mathcal{H}(p)
	=
		\sup_{\mu\in\mathbf{P}}
		\left[
		\int_{E^\prime} V_p \,\dd\mu - I_p(\mu)
		\right],
	$$
	where
	$ \mathbf{P} \subset \mathcal{P}(E^\prime)$ are the probability measures $\mu = (\mu_1,\dots,\mu_J)$ such that:
	 \begin{enumerate}[(P1)]
	 	\item 
	 	Each $\mu_i$ is absolutely continuous with respect to the uniform measure on $\mathbb{T}^d$.
	 	\item 
	 	For each~$i$, we have~$\nabla (\log \overline{\mu}_i) \in L^2_{\mu_i}(\mathbb{T}^d)$, where $\dd\mu_i(x) = \overline{\mu}_i(x) \dd x$.
	 \end{enumerate}
	\item
	We have
	\begin{equation}
		\mathcal{H}(p)
	=
		\sup_{\mu \in \mathbf{P}}
		\left[
		K_p(\mu) - \mathcal{R}(\mu)
		\right],
	\label{eq:results:LDP_MM:new_representation_H(p)}
	\end{equation}
	with the maps~$\mathcal{R}$ and~$K_p$ given by~\eqref{MM:eq:relative-Fisher-information} and~\eqref{eq:results:LDP_MM:K_p(mu)_I} above. In~$K_p(\mu)$, the infimum can be taken over vectors of functions $\phi_i = \phi(\cdot,i)$ such that $\nabla \phi_i \in L^{2}_{\mu_i}(\mathbb{T}^d)$.
	\item
	Under the detailed balance condition,
	\begin{multline}
		K_p(\mu) 
	= 
		\inf_{\phi}\bigg\{
		\sum_{i = 1}^J \int_{\mathbb{T}^d}
		\left(
		\frac{1}{2}|\nabla \phi^i(x) + p|^2 
	-
		\sum_{j = 1}^J r_{ij}(x)
		\right) \,
		\dd \mu_i(x)	\\
	+
		\sum_{i, j = 1}^J \int_{\mathbb{T}^d} r_{ij}(x)e^{-2\psi_i(x)}
		\sqrt{\overline{\mu}_i(x) \overline{\mu}_j(x)}
		e^{\psi_j(x) + \psi_i(x)} \cosh{(\phi(x,j) - \phi(x,i))}
		\,\dd x
		\bigg\}.
	\label{eq:results:LDP_MM:K_p(mu)_I:det_bal}
	\end{multline}
	\end{enumerate}
	\label{prop:results:detailed_balance_limit_I}
	\end{proposition}
	The representation~\eqref{eq:results:LDP_MM:K_p(mu)_I:det_bal} follows from~\eqref{eq:results:LDP_MM:K_p(mu)_I} by rewriting the sums appearing therein as~$\sum_{ij} a_{ij} = \frac{1}{2} \sum_{ij} (a_{ij} + a_{ji})$, where
	$$
		a_{ij} = \int_{\mathbb{T}^d}
		r_{ij}e^{-2\psi_i}
		\sqrt{\overline{\mu}^i(x) \overline{\mu}^j(x)}
		e^{\psi^j(x) + \psi^i(x)} e^{\phi(x,j) - \phi(x,i)}
		\,\dd x.
	$$
	This leads to the $\cosh(\cdot)$ terms in \eqref{eq:results:LDP_MM:K_p(mu)_I:det_bal}, and proves (c). We now give the proof of (a) and (b) of Proposition \ref{prop:results:detailed_balance_limit_I}.
	\begin{proof}[Proof of (a) in Proposition \ref{prop:results:detailed_balance_limit_I}]
%
Since $\mathcal{H}(p)$ is finite for any $p$ and $V_p(\cdot)$ is bounded, the supremum can be taken over measures $\mu$ such that $I_p(\mu)$ is finite. We show that finiteness of $I_p(\mu)$ implies that $\mu$ must satisfy (P1) and (P2). To that end, define the map $L_{\text{rev}} : \mathcal{D}(L_{\mathrm{rev}}) \subseteq C(E^\prime) \to C(E^\prime)$ by $\mathcal{D}(L_{\mathrm{rev}}) := C^2(E^\prime)$ and
	$$
		L_{\mathrm{rev}} f(x,i)
	=
		\frac{1}{2} \Delta_x f(x,i) - \nabla \psi_i(x) \cdot \nabla_x f(x,i)
	+
		\overline{\gamma} \sum_{j \neq i} s_{ij}(x)
		\left[
		f(x,j) - f(x,i)
		\right],
	$$
	with jump rates $s_{ij}$ defined as
	$
	s_{ij} \equiv 1
	$
	and
	$
	s_{ji} \equiv e^{2\psi_j - 2\psi_i},
	$
	for $i \leq j$, and with 
	$
	\overline{\gamma} := \sup_{\mathbb{T}^d} \left( r_{ij} / s_{ij}\right) < \infty,
	$
	where $r_{ij}(\cdot)$ are the jump rates appearing in~$L_p$. Furthermore, define $I_{L_{\text{rev}}} : \mathcal{P}(E^\prime) \to [0,\infty]$ by 
	$$
		I_{L_\mathrm{rev}}(\mu)
	:=
		-\inf_{\varphi \in C^2(E^\prime)}
		\int_{E^\prime} e^{-\varphi} L_{\mathrm{rev}} (e^{\varphi}) \, \dd\mu.
	$$
	We prove two statements. First, if $I_{L_{\mathrm{rev}}}(\mu)$ is finite, the measure $\mu$ satisfies (P1) and (P2). Second, if $I_p(\mu)$ is finite, then $I_{L_{\mathrm{rev}}}(\mu)$ is finite.
	Since 
	$
		s_{ij} e^{-2\psi_i}
	=
		s_{ji} e^{-2\psi_j},
	$
	the operator $L_\mathrm{rev}$ admits a reversible measure $\nu_{\text{rev}}$ in $\mathcal{P}(E^\prime)$ given by
	$$
		\nu_{\mathrm{rev}}(A_1, \dots, A_J)
	=
		\frac{1}{\mathcal{Z}} \sum_{i = 1}^J \nu_{\mathrm{rev}}^i(A_i),
	\quad
		\text{ where }
		\dd\nu_{\mathrm{rev}}^i = e^{-2 \psi_i} \dd x
		\text{ and }
		\mathcal{Z} = \sum_i \nu_\mathrm{rev}^i(\mathbb{T}^d).
	$$
	The measure $\nu_{\mathrm{rev}}$ is reversible for $L_{\mathrm{rev}}$ in the sense that for all $f,g \in \mathcal{D}(L_\text{rev})$,
	$$
		\langle L_\mathrm{rev} f, g \rangle_{\nu_\mathrm{rev}}
	=	
		\langle f, L_\mathrm{rev} g \rangle_{\nu_\mathrm{rev}},
	\quad
		\text{ where }
		\langle f, h \rangle_{\nu_\mathrm{rev}}
	=
		\frac{1}{\mathcal{Z}} \sum_i \int_{\mathbb{T}^d} f^i(x)h^i(x) \, \dd\nu_\mathrm{rev}^i(x).
	$$
	Hence by Stroock's result~\cite[Theorem 7.44]{stroock2012introduction},
	\begin{align*}
		I_{L_{\text{rev}}}(\mu)
	&=
		\begin{cases}
		\displaystyle
		-\langle f_\mu, L_\text{rev} f_\mu \rangle_{\nu_\text{rev}} ,& f_\mu = \sqrt{g_\mu} \in D^{1/2} := \mathcal{D}\left( \sqrt{-L_\text{rev}} \right) \text{ and } g_\mu = \frac{\dd\mu}{\dd \nu_\text{rev}}, \\
		\displaystyle +\infty ,& \text{otherwise},
		\end{cases}
	\end{align*}
	where $\dd\mu / \dd\nu_\text{rev}$ is the Radon-Nikodym derivative. This statement entails that if~$I_{L_\mathrm{rev}}(\mu)$ is finite, then~$\mu\ll\nu_{\mathrm{rev}}$. Then~$I_{L_\mathrm{rev}}(\mu)$ is explicitly given by
	\begin{multline}
		I_{L_{\text{rev}}}(\mu)
	=
		-\langle f, L_\text{rev} f \rangle_{\nu_\text{rev}}
	\\=
		\frac{1}{\mathcal{Z}}\sum_{i=1}^J
		\left[
		\int_{\mathbb{T}^d} |\nabla f^i(x)|^2 \, d\nu_\text{rev}^i(x)
	+
		\overline{\gamma} \sum_{j = 1}^J \int_{\mathbb{T}^d} s_{ij}(x)|f^j(x) - f^i(x)|^2 \, d\nu_\text{rev}^i(x)
		\right],
	\label{eq:LDP_MM:det_bal:I_L_rev_explicit}
	\end{multline}
	where we write $f^i = (\dd\mu^i/\dd{\nu^i_{\text{rev}}})^{1/2}$.	Furthermore, $\mu^i$ is absolutely continuous with respect to $\nu^i = e^{-2\psi^i} \dd x$. Since~$e^{-2\psi^i} \dd x \ll \dd x$,
	we find that~$\mu^i$ is absolutely continuous with respect to the volume measure on $\mathbb{T}^d$. Hence (P1) holds true.
	\smallskip
	
	We prove that finiteness of $I_{L_\mathrm{rev}}(\mu)$ implies (P2) by showing that the integral
	$
		\int_{\mathbb{T}^d} |\nabla (\log \overline{\mu}^i)|^2 \, d\mu^i
	$
	is finite. Let $g^i_\mu := \dd\mu^i / \dd\nu^i_{\mathrm{rev}}$ be the density of $\mu^i$ with respect to $\nu_\mathrm{rev}^i$. Then the densities $\overline{\mu}^i = \dd \mu^i/ \dd x$ satisfy~$g_\mu^i = \overline{\mu}^i e^{2\psi^i}$, because
	$$
		\overline{\mu}^i  
	=
		\frac{\dd\mu^i}{\dd\nu^i_\text{rev}} \frac{\dd\nu_\text{rev}^i}{\dd x}
	=
		\frac{\dd\mu^i}{\dd\nu^i_\text{rev}} e^{-2\psi^i}.
	$$
	Let $f_\mu^i := \sqrt{g_\mu^i}$. If $I_{L_\mathrm{rev}}(\mu)$ is finite, then by~\eqref{eq:LDP_MM:det_bal:I_L_rev_explicit}, $\int_{\mathbb{T}^d}|\nabla f_\mu^i|^2 \dd\nu_\mathrm{rev}^i$ is finite for every $i = 1,\dots, J$. Hence with the estimate
	\begin{align*}
		\int_{\mathbb{T}^d} |\nabla f_\mu^i|^2 \dd\nu_\text{rev}^i
	&\geq
		\int_{\mathbb{T}^d} |\nabla f_\mu^i|^2 \mathbf{1}_{\{ \overline{\mu}^i > 0\}} \dd\nu_\text{rev}^i
	=
		\frac{1}{4}\int_{\mathbb{T}^d} \frac{|\nabla g_\mu^i|^2}{g_\mu^i} \mathbf{1}_{\{ \overline{\mu}^i > 0\}} \dd\nu^i_\text{rev}	\\
	&=
		\frac{1}{4}\int_{\mathbb{T}^d} \frac{| e^{2\psi^i} \nabla \overline{\mu}^i + 2 \overline{\mu}^i \nabla \psi^i e^{2\psi^i}|^2 }{\overline{\mu}^i} e^{-4\psi^i} \mathbf{1}_{\{ \overline{\mu}^i > 0\}} \dd x\\
	&=
		\frac{1}{4} \int_{\mathbb{T}^d} | \nabla (\log \overline{\mu}^i) + 2 \nabla \psi^i |^2 \mathbf{1}_{\{ \overline{\mu}^i > 0\}} \dd\mu^i	\\
	&\geq
		\frac 18 \int_{\mathbb{T}^d} |\nabla(\log \overline{\mu}^i)|^2 \mathbf{1}_{\{ \overline{\mu}^i > 0\}} \dd\mu^i
		-
		\int_{\mathbb{T}^d} |\nabla \psi^i|^2 \mathbf{1}_{\{ \overline{\mu}^i > 0\}} \dd\mu^i,
	\end{align*}
	we find~$\nabla (\log \overline{\mu}^i) \in L^2_{\mu^i}(\mathbb{T}^d)$.	We are left with proving that if $I_p(\mu)$ is finite, then $I_{L_{\mathrm{rev}}}(\mu)$ is finite. Estimating $r_{ij} / s_{ij}$ from above by $\overline{\gamma} = \sup_{\mathbb{T}^d} (r_{ij} / s_{ij})$, we find
	\begin{multline*}
		I_p(\mu)
		\geq
		\sup_\varphi \sum_i \int_{\mathbb{T}^d}
		-
		\bigg[
		\frac{1}{2} \Delta \varphi^i(x) + \frac{1}{2} |\nabla \varphi^i(x)|^2 + (p - \nabla\psi^i(x)) \nabla \varphi^i(x) 
	\\+ 
		\overline{\gamma} \sum_{j\neq i} s_{ij}(x)
		(e^{\varphi(x,j) - \varphi(x,i)} -1)
		\bigg] \, \dd\mu^i - s_0(\mu),
	\end{multline*}
	where 
	$
		s_0(\mu)
	=
		\sum_{ij} \int_{\mathbb{T}^d} \left[ \overline{\gamma} \, s_{ij}(x) - r_{ij}(x) \right] \dd\mu^i
	$
	is finite. For $p = 0$, this means that $I_0 (\mu) \geq I_{L_\mathrm{rev}}(\mu) - s_0(\mu)$ holds for all $\mu \in \mathcal{P}(E^\prime)$. In particular, if $I_0 (\mu)$ is finite, then $I_{L_\mathrm{rev}}(\mu)$ is finite. For $p \neq 0$, the additional $p$-term is dealt with by Young's inequality applied as $-p\cdot \nabla \phi^i \geq - p^2/(2 \varepsilon) - \frac{\displaystyle \varepsilon}{2} |\nabla\phi^i|^2$. Thereby,
	\begin{multline*}
		I_p(\mu)
	\geq
		\sup_\varphi \sum_i \int_{\mathbb{T}^d}
		-
		\bigg[
		\frac{1}{2} \Delta \varphi^i(x) + \frac{1 + \varepsilon}{2} |\nabla \varphi^i(x)|^2 + - \nabla\psi^i(x) \nabla \varphi^i(x) 
	\\+ 
		\overline{\gamma} \sum_{j\neq i} s_{ij}(x)
		(e^{\varphi(x,j) - \varphi(x,i)} -1)
		\bigg] \, \dd\mu^i -\frac{p^2}{2\varepsilon}- s_0(\mu)
	\\=
		\frac{1}{\lambda}
		\sup_\varphi \sum_i \int_{\mathbb{T}^d}
		-
		\bigg[
		\frac{1}{2} \Delta \varphi^i(x) + \frac{1}{2} |\nabla \varphi^i(x)|^2 + - \nabla\psi^i(x) \nabla \varphi^i(x) 
	\\+ 
		\lambda \overline{\gamma} \sum_{j\neq i} s_{ij}(x)
		(e^{(\varphi(x,j) - \varphi(x,i)) / \lambda} -1)
		\bigg] \, \dd\mu^i
		-\frac{p^2}{2\varepsilon}- s_0(\mu),		
	\end{multline*}
	where the last equality follows by rescaling $\varphi \rightarrow \varphi / \lambda$, with $\lambda = 1 + \varepsilon > 1$. Therefore, apart from the factor $1/\lambda$ in the exponential term and the multiplicative factor $\lambda \overline{\gamma}$, we obtain the same estimate as above in the $p = 0$ case. Denoting the supremum term in the last line by $I^\lambda_{L_\text{rev}}$, we found the estimate
	\[
		I_p(\mu)
	\geq
		\frac{1}{\lambda} I^\lambda_{L_\text{rev}}(\mu) - s_p(\mu),
	\]
	where
	$
		s_p(\mu) = \frac{p^2}{2\varepsilon} + s_0(\mu)
	$
	is finite. We now show that~$I_{L_\text{rev}}(\mu) = \infty$ implies $I^\lambda_{L_{\text{rev}}}(\mu) = \infty$, which proves that finiteness of $I_{L_\text{rev}}^\lambda (\mu)$ implies finiteness of $I_{L_{\text{rev}}}(\mu)$.
	If $I_{L_\text{rev}}(\mu) = \infty$, then by definition, there exist functions $\varphi_n$ such that
	\begin{multline*}
		a(\varphi_n)
	:=
		-\sum_{i=1}^J \int_{\mathbb{T}^d}
		\bigg[
		\frac{1}{2} \Delta\varphi_n^i + \frac{1}{2} |\nabla \varphi_n^i|^2
		-\nabla\psi^i \nabla \varphi_n^i
	\\+ 
		\overline{\gamma} \sum_{j\neq i} s_{ij}
		\left(
		e^{\varphi_n(x,j) - \varphi_n(x,i)} - 1
		\right)
		\bigg]
		\, \dd\mu^i(x)
	\xrightarrow{ n \rightarrow \infty} \infty.
	\end{multline*}
	We aim to prove that $a(\varphi_n) \leq I^\lambda_{L_\text{rev}}(\mu)$ holds for all $n$. To that end, write
	\begin{multline*}
		a^\lambda(\varphi_n)
	:=
		-
		\sum_i \int_{\mathbb{T}^d}
		\bigg[
		\frac{1}{2} \Delta \varphi_n^i + \frac{1}{2} |\nabla \varphi_n^i|^2 + - \nabla\psi^i \nabla \varphi_n^i 
		\\+ \lambda \overline{\gamma} \sum_{j\neq i} s_{ij}
		(e^{(\varphi_n(x,j) - \varphi_n(x,i)) / \lambda} -1)
		\bigg] \, \dd\mu^i
	\end{multline*}
	for the according evaluation of $\varphi_n$ in $I^\lambda_{L_\text{rev}}(\mu)$. we have 
	$
		I^\lambda_{L_\text{rev}}(\mu)
	\geq
		a^\lambda(\varphi_n),
	$
	and show that $a^\lambda(\varphi^n) \rightarrow \infty$. 
	The only difference between $a(\varphi_n)$ and $a^\lambda(\varphi_n)$ lies in the $\lambda$-factors that appear in the exponential terms. Since~$e^x \geq e^x \mathbf{1}_{\{x \geq 0\}}$,
	\begin{multline*}
		\overline{a}_n
	:=
		-\sum_{i=1}^J \int_{\mathbb{T}^d}
		\bigg[
		\frac{1}{2} \Delta\varphi_n^i + \frac{1}{2} |\nabla \varphi_n^i|^2
		-\nabla\psi^i \nabla \varphi_n^i
	\\+ 
		\overline{\gamma} \sum_{j\neq i} s_{ij}
		\left(
		e^{\varphi_n(x,j) - \varphi_n(x,i)} \mathbf{1}_{\{\varphi_n(x,j) - \varphi_n(x,i) \geq 0\}} - 1
		\right)
		\bigg]
		\, \dd\mu^i(x)
	\end{multline*}
	diverges as $n \rightarrow \infty$ (we have $a(\varphi_n) \leq \overline{a}_n$). Define this analogously for $a^\lambda(\varphi_n)$,
	\begin{multline*}
		\overline{a}^\lambda_n
	:=
		-
		\sum_i \int_{\mathbb{T}^d}
		\bigg[
		\frac{1}{2} \Delta \varphi_n^i + \frac{1}{2} |\nabla \varphi_n^i|^2 + - \nabla\psi^i \nabla \varphi_n^i 
	\\+ 
		\lambda \overline{\gamma} \sum_{j\neq i} s_{ij}
		\left(
		e^{(\varphi_n(x,j) - \varphi_n(x,i)) / \lambda} \mathbf{1}_{\{ \varphi_n(x,j) - \varphi_n(x,i) \geq 0 \}} -1
		\right)
		\bigg] \, \dd\mu^i.
	\end{multline*}
	Since~$a^\lambda(\varphi_n)\geq
		\overline{a}^\lambda_n - \sum_{ij} \int_E \lambda \overline{\gamma} s_{ij}(x) \, d\mu^i(x)$, proving that $\overline{a}_n^\lambda \rightarrow \infty$ as $n \to \infty$ is sufficient for obtaining $a^\lambda(\varphi_n) \rightarrow \infty$. Finally, the fact that $\overline{a}^\lambda_n$ diverges as $n \to \infty$ follows by noting that $\overline{a}_n\leq \overline{a}^\lambda_n$, which can be seen via
	\begin{multline*}
		\overline{a}_n - \overline{a}^\lambda_n
	=
		-\sum_i \int_{\mathbb{T}^d} \overline{\gamma} \sum_{j \neq i} s_{ij}
		\left( e^{\varphi_n^j - \varphi_n^i} \mathbf{1}_{\{ \varphi_n^j - \varphi_n^i \geq 0 \}} - 1 \right) d\mu^i
	\\+
		\lambda \overline{\gamma} \sum_i \int_{\mathbb{T}^d} s_{ij} 
		\left(
		e^{\varphi_n^j - \varphi_n^i} \mathbf{1}_{\{ \varphi_n^j - \varphi_n^i \geq 0 \}} - 1
		\right) \dd\mu^i	\\
	=
		\overline{\gamma} 
		(1 - \lambda)
		\sum_{ij} \int_{\mathbb{T}^d} s_{ij} d\mu^i
	+
		\overline{\gamma} \sum_{ij} \int_{\mathbb{T}^d} s_{ij} 
		\left(
		e^{(\varphi_n^j - \varphi_n^i) / \lambda} - e^{\varphi_n^j - \varphi_n^i}
		\right)
		\mathbf{1}_{\{ \varphi_n^j - \varphi_n^i \geq 0 \}}
		\dd\mu^i,
	\end{multline*}
	which is bounded above by zero since $\lambda = 1 + \varepsilon > 1$ and $e^{x / \lambda} \leq e^{x}$ for $x \geq 0$. This finishes the proof of part (a) of Proposition \ref{prop:results:detailed_balance_limit_I}. 
	\end{proof}
	\begin{proof}[Proof of (b) of Proposition \ref{prop:results:detailed_balance_limit_I}]
	It is sufficient to show that for any $\mu \in \mathbf{P}$, the Donsker-Varadhan functional~$I_p(\mu)$ satisfies~\eqref{MM:eq:DV-functional-after-shift}.
	Integration by parts gives
	\begin{multline*}
		I_p(\mu)
	=
		-\inf_{\varphi} \sum_i \int_{\mathbb{T}^d}	
		\bigg[
		-\frac{1}{2} \nabla\varphi^i \nabla(\log \overline{\mu}^i)
		+ \frac{1}{2}|\nabla\varphi^i|^2
		+(p-\nabla\psi^i) \nabla\varphi^i
		\\+ \sum_j r_{ij}
		\left(
		e^{\varphi^j - \varphi^i} - 1
		\right)
		\bigg]
		\dd \mu^i,
	\end{multline*}
	where $\dd \mu^i = \overline{\mu}^i \dd x$. By a density argument, the infimum can be taken over functions in~$L^{1,2}_{\mu^i}(\mathbb{T}^d)$. 
	Now shifting in the infimum as~$\varphi_i\to\varphi_i+\frac{1}{2}\log(\overline{\mu}_i) + \psi^i$, we find after some algebra that
%
	\begin{multline*}
	I_p(\mu) = -\inf_\varphi \sum_i \int_{\mathbb{T}^d}
	\bigg[
	\frac{1}{2} |\nabla\varphi^i + p|^2
	- \frac{1}{2} | (p - \nabla\psi^i) - \frac{1}{2} \nabla \log\overline{\mu}^i|^2
	\\+
	\sum_j r_{ij}
	\left(
	\sqrt{\frac{\overline{\mu}^j}{\overline{\mu}^i}} e^{\psi^j - \psi^i} e^{\varphi^j - \varphi^i}
	- 1\right) \bigg]\,\dd\mu^i.	
	\end{multline*}
	The term containing the square roots and logarithms are not singular since they are integrated against $\dd\mu^i$, so that the integration is over the set $\{\overline{\mu}^i > 0\}$. Now writing out the terms and reorganizing them leads to the claimed equality.
	\end{proof}
\section{Principal eigenvalues and their variational representations}
	\sectionmark{Principal eigenvalues}
	\label{appendix:prinipal_ev}
	In this section, we collect some results about the principal eigenvalue problems that we encounter in this chapter.
	\begin{definition}[Irreducible M-matrix]
		A matrix $P \in \mathbb{R}^{d\times d}$ is an \emph{irreducible M-matrix} if $P = s\mathbf{1} - R$, with some $s\in \mathbb{R}$ and an irreducible matrix $R \geq 0$ with non-negative elements.
		\label{def:appendix:PrEv:irreducible_M_matrix}
	\end{definition}
	The eigenvalue problems are the following:
	\begin{enumerate}[(E1)]
	\item
	For an irreducible M-matrix~$P \in \mathbb{R}^{d \times d}$, find a real eigenvalue $\lambda$ and a corresponding eigenvector $v > 0$ that has strictly positive components $v_j > 0$, such that $P v = \lambda v$. The eigenvalue problems arising for the discrete models (Lemmas \ref{lemma:LDP_MM:discrI:principal_ev} and \ref{lemma:LDP_MM:discrII:principal_ev}) are of that type.
	\item
	For a second-order uniformly elliptic operator given by
	\begin{equation}
		P
	= 
		-\sum_{k \ell} a_{k \ell}(\cdot) \frac{\partial^2}{\partial x^k \partial x^\ell} + \sum_k b_k(\cdot) \frac{\partial}{\partial x^k} + c(\cdot),
	\label{eq:appendix:uniform_elliptic_op}
	\end{equation}
	with smooth coefficients $a_{k \ell}, b_k, c \in C^\infty(\mathbb{T}^d)$, find a real eigenvalue $\lambda$ and a corresponding strictly positive eigenfunction $u$ such that $P u = \lambda u$. This corresponds to the eigenvalue problem in Lemma \ref{lemma:LDP_MM:contII:principal_eigenvalue}, with $a_{k \ell} = 1$.
	\item
	For a coupled system of second-order elliptic operators on $\mathbb{T}^d$, find a real eigenvalue $\lambda$ and a vector of strictly positive functions $u = (u^1, \dotsm u^J)$, $u^i > 0$ on $\mathbb{T}^d$, such that
	\begin{equation}
	\left[
	\begin{pmatrix}
		L^{(1)}	&	&	0	\\
			&	\ddots	&	\\
		0	&	&	L^{(J)}		
	\end{pmatrix}
	-
	\begin{pmatrix}
		R_{11}	&	&	\geq 0	\\
    		& \ddots	&	\\
    	\geq 0	&		&	R_{JJ}
	\end{pmatrix}
	\right]
	\begin{pmatrix}
		u^1	\\
		\vdots	\\
		u^J
	\end{pmatrix}
	=
	\lambda
	\begin{pmatrix}
		u^1	\\
		\vdots	\\
		u^J
	\end{pmatrix},
	\label{eq:appendix:PrEv:eigenproblem:coupled_system_of_elliptic_op}
	\end{equation}
	where $L : C^2(\mathbb{T}^d)^J \rightarrow C(\mathbb{T}^d)^J$ is a $J\times J $ diagonal matrix of uniformly elliptic operators,
	\begin{align}
	L
	=
		\begin{pmatrix}
    	L^{(1)}	&	&	0	\\
    		& \ddots	&	\\
    	0	&	& L^{(J)}       
  		\end{pmatrix},\;L^{(i)} = -\sum_{k \ell}^J a_{k \ell}^{(i)}(\cdot) \frac{\partial^2}{\partial x^k \partial x^\ell} + \sum_k^J b^{(i)}_k(\cdot)\frac{\partial}{\partial x^k} + c^{(i)}(\cdot),
  	\label{eq:appendix:elliptic_operators_in_coupled_system}
	\end{align}
	with
	$
		a_{k \ell}^{(i)}(\cdot), b_k^{(i)}(\cdot), c^{(i)}(\cdot) \in C^\infty(\mathbb{T}^d),
	$
	and $R$ is a $J\times J$ matrix with non-negative functions on the off-diagonal,
	\begin{align*}
	R
	=
	\begin{pmatrix}
		R_{11}	&	&	\geq 0	\\
    		& \ddots	&	\\
    	\geq 0	&		&	R_{JJ}
	\end{pmatrix}, 
	\qquad
	R_{ij} \geq 0 \text{  for all  } i \neq j.
	\end{align*}
	Coupled systems of this type appear in Lemma \ref{lemma:LDP_MM:contI:principal_eigenvalue}.
	\end{enumerate}
	The principal-eigenvalue problems (E1), (E2) and (E3) can be solved by means of the Krein-Rutman Theorem. We recall the setting of the Theorem.
	\begin{definition}[{Ordered Banach space $\left(X,\geq\right)$~ \cite[Appendix 4]{dautray2000spectral}}]
	For a real Banach space $X$, a closed set $K\subseteq X$ with nonempty interior is called a \emph{cone} if i) $0 \in K$, ii) whenever $v,w\in K$ then $av + bw \in K$ for all reals $a,b \geq 0$, iii) if $v\in K$ and $(-v)\in K$, then $v = 0$, and iv) $X = K - K$. For given $v,w\in X$, we write $v \geq w$ if $ v- w \in K$, and denote the elements $v$ in $K$ as $v\geq 0$ the elements in the interior $\mathring{K}$ as $v>0$. Further, $K^\ast \subseteq X^\ast$ is called a \emph{dual cone} if for all $\ell \in K^\ast$, $\langle \ell, v \rangle \geq 0$ whenever $v\ \geq 0$. We write $(X,\geq)$ for an ordered Banach space $X$, where the order $\geq$ is defined by means of a cone $K$.
	\end{definition}
	For an ordered Banach space $(X,\geq)$ and an operator $P:\mathcal{D}(P)\subseteq X \rightarrow X$, we want to find a strictly positive eigenvector~$u>0$ with an associated eigenvalue~$\lambda \in \mathbb{R}$ such that
	\begin{equation}
		P u = \lambda u.
		\tag{PrEv}
		\label{eq:solving_principal_ev:basic_ev_equation}
	\end{equation}
	The problems (E1), (E2) and (E3) are of this type, in the following settings:
	\begin{enumerate}[(E1)]
	\item
	$X = \mathbb{R}^d$, with cone $K = \{v\in\mathbb{R}^d : v_j\geq 0,\,j =1,\dots,d\}$, and corresponding interior $\mathring{K} = \{ v \in K : v_j > 0,\,j=1,\dots,d\}$. The operator $P: \mathbb{R}^d \to \mathbb{R}^d$ is an irreducible M-matrix.
	\item
	$X = C(\mathbb{T}^d)$, with cone $K = \{ f \in X: f \geq 0\}$ and corresponding interior $\mathring{K} = \{f\in X : f > 0\}$. The operator $P : C^2(\mathbb{T}^d) \subseteq C(\mathbb{T}^d) \to C(\mathbb{T}^d)$is~\eqref{eq:appendix:uniform_elliptic_op}.
	\item
 	$X = C(E^\prime)$, with $E^\prime = \mathbb{T}^d \times \{1, \dots, J\}$ and cone $K = \{f \in X: f(x,i) \geq 0,\, x \in \mathbb{T}^d\, i = 1,\dots,J\}$, and interior $\mathring{K} = \{f\in K : f(\cdot,i) > 0,\,i = 1,\dots,J\}$. We identify $C(E^\prime)$ with $C(\mathbb{T}^d)^J$ via $f(x,i) = f^i(x)$, $f = (f^1, \dots, f^J)$. 

	\end{enumerate}
	An operator $B : \mathcal{D}(B) \subseteq X \rightarrow X$ is is called positive if $f\geq 0$ implies $Bf \geq 0$, and is called strongly positive if $f\geq 0$ and $f \neq 0$ imply $Bf > 0$.
	\begin{theorem}[Krein-Rutman, Appendix 4 in \cite{dautray2000spectral}]
	Let $(X,\geq)$ be an ordered Banach space and $T :X \rightarrow X$ be a linear bounded operator. If $T$ is also compact and strongly positive, then there exist unique $g>0$ and $g^\ast > 0$ such that
	\[
	T g = r(T) g, \; \|g\|_X =1,
	\qquad
	\text{and}
	\qquad
	T^\ast g^\ast = r(T) g^\ast, \; \|g^\ast\|_{X^\ast} = 1,
	\]
	with $T^\ast$ the dual operator to $T$, and $\langle g^\ast , f\rangle >0$ whenever $f\geq 0$ and $f\neq 0$. Here, $r(T) = r(T^\ast) $ is the spectral radius of $T$.
	\label{thm:appendix:solving_principal_ev:KR-theorem}
	\end{theorem}
	\begin{theorem}[Positive and compact resolvant implies existence of a principal eigenvalue]
	If for some $\alpha \in \mathbb{R}$, $P_\alpha := P + \alpha \mathbf{1}$ is such that $T_\alpha := P_\alpha^{-1}$ exists as a linear bounded operator $T_\alpha :X \rightarrow X$ that is compact and strongly positive, then \eqref{eq:solving_principal_ev:basic_ev_equation} holds with $\lambda = \frac{1}{r(T_\alpha)} - \alpha$ and eigenfunction $u = T_\alpha g$, where $g$ satisfies $T_\alpha g = r(T_\alpha) g$. Furthermore, $\lambda \in \mathbb{R}$ is the unique eigenvalue with a strictly positive eigenvector.
	\label{thm:appendix:solving_principal_ev:apply_KR}
	\end{theorem}
	\begin{proof}[Proof of Theorem \ref{thm:appendix:solving_principal_ev:apply_KR}]
	By the Krein-Rutman Theorem \ref{thm:appendix:solving_principal_ev:KR-theorem}, there exists a $g>0$ such that $T_\alpha g = r(T_\alpha) g$. By strong positivity of $T_\alpha$, we have $u:= T_\alpha g > 0$, and in particular $r(T_\alpha) >0$. By definition of $T_\alpha$ as the solution operator $h\mapsto f$ of $P_\alpha f =h$, the vector $u \in \mathring{K}$ satisfies $P_\alpha u = \frac{1}{r(T_\alpha)} u$, and \eqref{eq:solving_principal_ev:basic_ev_equation} follows with principal eigenvalue $\lambda = \frac{1}{r(T_\alpha)} - \alpha$ and strictly positive eigenfunction $u > 0$. Regarding uniqueness of the eigenvalue $\lambda$, note that every solution to \eqref{eq:solving_principal_ev:basic_ev_equation} defines an eigenfunction for $T_\alpha$, by shifting with $\alpha$. Thus two independent solutions to \eqref{eq:solving_principal_ev:basic_ev_equation} would correspond to two independent solutions to $T_\alpha g = r(T_\alpha) g$, contradicting the uniqueness (after normalization) of $g>0$ in the Krein-Rutman Theorem. 	
	\end{proof}
	Theorem \ref{thm:appendix:solving_principal_ev:apply_KR} applies to the eigenvalue problems (E1), (E2) and (E3).
	\begin{proposition}
	In the setting (E1), if $P \in \mathbb{R}^{d \times d}$ is an irreducible M-matrix, then there exists an eigenvector $u>0$ and a unique principal eigenvalue $\lambda \in \mathbb{R}$ such that \eqref{eq:solving_principal_ev:basic_ev_equation} holds. The principal eigenvalue~$\lambda$ is given by
	$$
		\lambda = \sup_{w > 0} \inf_{i \in \{1, \dots, J\}} \frac{P w(i)}{w(i)}.
	$$
	\label{proposition:appendix:PrEv:M_matrix}
	\end{proposition}
	\begin{proposition}
	In the setting (E2), let $P$ be given by~\eqref{eq:appendix:uniform_elliptic_op}. Then there exists a strictly positive eigenfunction $u \in C^\infty(\mathbb{T}^d)$ and a unique principal eigenvalue $\lambda \in \mathbb{R}$ satisfying \eqref{eq:solving_principal_ev:basic_ev_equation}. The principal eigenvalue~$\lambda$ is given by
	$$
		\lambda
	=
		\sup_{g > 0} \inf_{x\in \mathbb{T}^d} \left[\frac{P g(x)}{g(x)}\right]
	=
		\inf_{\mu \in \mathcal{P}(\mathbb{T}^d)} 
		\sup_{g > 0}
	\left[
		\int_{\mathbb{T}^d}
		\frac{Pg}{g} \,d\mu
	\right].
	$$
	\label{proposition:appendix:PrEv:elliptic_op}
	\end{proposition}
	\begin{proposition}
	In the setting (E3), let~$L$ be given by~\eqref{eq:appendix:PrEv:eigenproblem:coupled_system_of_elliptic_op} and~\eqref{eq:appendix:elliptic_operators_in_coupled_system}. Suppose that the matrix $\overline{R}$ with entries $\overline{R}_{ij} := \sup_{y \in \mathbb{T}^d} R_{ij}(y)$ is irreducible.
	Then for the operator $P:= L - R$, there exists a unique principal eigenvalue $\lambda \in \mathbb{R}$ and a strictly positive eigenvector $u \in \left(C^\infty(\mathbb{T}^d)\right)^J$, $u^i(\cdot) > 0$ for all $i = 1, \dots, J$, solving \eqref{eq:solving_principal_ev:basic_ev_equation}. Furthermore, the principal eigenvalue is given by
	$$
		\lambda
	=
		\sup_{g > 0} \inf_{z \in E^\prime} \left[\frac{P g(z)}{g(z)}\right]
	=
		\inf_{\mu \in \mathcal{P}(E^\prime)} 
		\sup_{g > 0}
	\left[
		\int_{E^\prime}
		\frac{Pg}{g} \,d\mu
	\right].
	$$
	\label{proposition:appendix:PrEv:fully_coupled_system}
	\end{proposition}
	The principal eigenvalue problem on closed manifolds, such as $\mathbb{T}^d$, is solved for instance by Padilla \cite{padilla1997principal}. Donsker and Varadhan's variational representations \cite{DonskerVaradhan75, DonskerVaradhan76} apply to the case of compact metric spaces without boundary. A proof of how to obtain the principal eigenvalue for coupled systems of equations is given by Sweers \cite{Sweers92} and Kifer \cite{kifer1992principal}. Sweers considers a Dirichlet boundary problem, but his results transfer to the compact setting without boundary. Kifer gives an independent proof for the case of a compact manifold, in Lemma 2.1 and Proposition 2.2 in \cite{kifer1992principal}.
\chapter{Large Deviations of Empirical Measures}
\label{chapter:LDP-of-empirical-measures}
\section{Introduction}
In this chapter, we are motivated by the task of sampling from a distribution~$\pi$ with density with respect to Lebesgue measure given by
\begin{align*}
    \dd \pi (y) = C^{-1}e^{-U(y)}\,\dd y, \quad C = \int_E e^{-U(y)}\,\dd y,
\end{align*}
for some potential function $U: E \to \bR$ and state space $E$. The most common approach is to use Markov chain Monte Carlo (MCMC) methods, which are now essential tools in areas such as computational statistics, molecular dynamics and machine learning~\cite{RobertCasella2004,AsmussenGlynn07,AndrieuEtAl03}.
\smallskip

The idea behind MCMC is to construct a Markov process $Y_t$ with $\pi$ as the invariant measure and use the corresponding empirical measure to obtain approximations. For example, under ergodicity, for any observable $f \in L^1 (\pi)$ we have almost surely
\begin{equation*}
\lim_{t\to\infty} \frac{1}{t}\int_0^t f(Y_s)\,\dd s = \int_E f(y)\,\pi(\dd y).
\end{equation*}
Therefore, for $t>0$ large, $\frac{1}{t} \int _0 ^t f(Y_s)\,\dd s$ can be used to approximate the expected value $\int_E f (y) \pi (\dd y)$. 
Although many standard MCMC constructions, such as the Metropolis-Hastings algorithm \cite{Metropolis1953}, can be used to sample from essentially any target distribution $\pi$, most suffer from slow convergence to the invariant distribution or heavy computational costs per iteration. Designing new, efficient dynamics has therefore become an important research direction within applied probability.
\smallskip

Over the last decade, piecewise-deterministic Markov processes (PDMPs) have emerged as a new tool for the numerical simulation of probability distributions. An introduction to these processes is offered by Davis' monograph~\cite{davis1984piecewise}. The two main examples of such processes used in MCMCs are the Bouncy Particle Sampler and the Zig-Zag Sampler \cite{BouchardCoteVollmerDoucet2017,bierkens2019zig}, after similar ideas appeared first in~\cite{PetersDeWith2012} and~\cite{Monmarche2016}. The idea of using PDMPs extends the ubiquitous discrete time MCMC methodology towards a new continuous time approach, having several advantageous aspects. First, by construction PDMPs are irreversible Markov processes, which typically results in a smaller asymptotic variance as compared to reversible methods. For instance, Duncan, Lelièvre and Pavliotis demonstrate variance reduction for irreversible Langevin samplers~\cite{DuncanLelievrePavliotis2016}. We refer to~\cite{Andrieu2019a} for a recent study of this effect, and to~\cite{bierkens2019zig,Fearnhead2016a} for details of the computational aspects of PDMP trajectories on a computer.
\smallskip

In order to employ this new PDMP methodology a solid understanding of the mathematical properties of these methods is necessary. Whereas the theoretical properties of PDMPs have been an active research area in recent years, our understanding of the performance of the corresponding MCMC methods is still incomplete. In particular knowledge of the speed of convergence of time averages is essential in choosing the most suitable sampling technology for a particular problem and in tuning the parameters of the chosen method. In the spirit of recent work on empirical measure large deviations in the MCMC context \cite{dupuis2012infinite, rey2015irreversible}, we propose the use of large deviation results for studying and comparing the performance of PDMPs. 
\smallskip

In summary, the main contributions we develop in this chapter are:
\begin{itemize}
    \item A semigroup approach to establish the large deviation principles for the empirical measures of a class of Markov processes satisfying assumptions aimed at position-velocity PDMPs.
    \item The large deviation principle for empirical measures of the zig-zag process in both a compact and non-compact setting.
    \item A derivation of an explicit form of the rate function associated with the zig-zag process.
    \item Evaluation of the zig-zag rate function as a function of the additional switching rate $\gamma$, providing an answer to a key question about the switching rate.
\end{itemize}
Donsker and Varadhan studied large deviations for empirical measures in a series of papers~\cite{donsker1975asymptoticI, donsker1975asymptoticII, DonskerVaradhan76}.
In the simulation context it is well-known that for rare-event simulation, sample-path large deviations play an important r\^ole in evaluating and designing efficient algorithms ; see \cite{AsmussenGlynn07, Bucklew04, BudhirajaDupuis2019} and references therein. 
In contrast, empirical measure large deviations are much less explored as a tool for analysing Monte Carlo methods. Standard measures for analysing the efficiency of methods based on ergodic Markov processes include the spectral gap of the associated semigroup and the asymptotic variance for given observables, see for example \cite{Rosenthal2003, BedardRosenthal08,DiaconisHolmesNeal00,FrankeEtAl10,FrigessiEtAl93,HwangHwangSheu05, MengersenTweedie96, RobertsRosenthal04}. 
However these measures are not necessarily appropriate for studying the rate of convergence, as they only link indirectly to the empirical measure, the quantity of interest in Monte Carlo methods. Empirical measure large deviations on the other hand connect explicitly to the relevant properties, such as the transient behaviour of the underlying process. In a similar spirit, \cite{BirrellRB19} recently used concentration inequalities to obtain non-asymptotic performance guarantees for PDMPs.
\smallskip

The first results using empirical measure large deviations for the analysis of MCMC methods were obtained in~\cite{plattner2011infinite, dupuis2012infinite}. Therein empirical measure large deviations, specifically the associated rate function, was proposed as a tool for analysing parallel tempering, one of the computational workhorses of the physical sciences, leading to a new type of simulation method (infinite swapping). In the subsequent work~\cite{DDN2018} empirical measure large deviations were again used, combined with associated stochastic control problems, to analyse the convergence properties of these algorithms. Similarly, in~\cite{rey2015irreversible} Rey-Bellet and Spiliopoulos use empirical measure large deviations to analyse the performance of certain irreversible MCMC samplers. 
\smallskip

The work by Donsker and Varadhan is the starting point for many results and application of empirical measure large deviations and their work has been extended in numerous directions, see e.g.\ \cite{DemboZeitouni1998, FengKurtz2006, BudhirajaDupuis2019} for an overview and further references. However, naively applying the existing theory to PDMPs does not work since the transition probabilities are not sufficiently regular: for every $t > 0$ there is a positive probability that the process has not switched by time $t$, resulting in an atomic component of the Markov transition kernel. As a first step towards using empirical measure large deviations for analysing the performance of PDMPs we must therefore establish the relevant large deviations principles. 
\smallskip

In this chapter, our focus is to establish general large deviation results aimed at PDMPs and then to specialize to the zig-zag process. In the process of proving the necessary large deviation results we consider a general class of Markov processes that can have position-velocity PDMPs, such as the bouncy particle and the zig-zag samplers, as special cases. In particular, this class includes processes that are not of diffusion type and irreversible processes. When specialising to the zig-zag process, we derive an explicit form of the rate function, going beyond the variational form typical for results of Donsker-Varadhan-type. To the best of our knowledge this is the first instance where an explicit form of the rate function has been obtained for irreversible processes that do not have a drift-diffusion character. 
\smallskip

A key question for using the zig-zag process for MCMC is whether or not it is advantageous for convergence to use the minimal (canonical) switching rates, or if one should allow for additional switches according to a fixed refreshment rate $\gamma >0$. Our analysis of the rate function associated with the zig-zag process allows us to give a partial answer to this question: in Section~\ref{sec:rate} we establish that the rate function is decreasing as a function of the additional rate $\gamma$, establishing that from a large deviations perspective it is optimal to use the smallest possible rates, i.e.\ set $\gamma =0$. This goes in the opposite direction of the conclusion drawn from a spectral analysis (see \cite[Section 7.3]{BierkensVerduynLunel2019}), which shows at least a small benefit of increasing gamma beyond zero. This highlights the different nature of convergence of empirical averages by studying large deviations or asymptotic variance (e.g.~\cite{Andrieu2019a,BierkensDuncan2016}) and convergence to equilibrium, using e.g. the spectral gap to describe rate of convergence; We refert to~\cite{Rosenthal2003} for more on this phenomenon. Our conclusion is in line with the earlier observation that having more irreversibility increases the rate function~\cite{rey2015irreversible}: one can view increasing $\gamma$ as decreasing the extent of irreversibility inherent to the process. In that sense,~$\gamma=0$ corresponds to "maximal irreversibility" of the zig-zag process. The fact that the spectral gap can not always detect the benefits of irreversibility is best illustrated with the following example, which can also be found in~\cite[Example~2.9]{rey2015irreversible}.
\paragraph{Example.}
	Consider the diffusion~$\dd X_t = v\,\dd t + \dd B_t$ on the one-dimensional flat torus~$\mathbb{T}$, where~$v$ is a parameter. We focus on the behavior of its empirical measure~$\eta_T$ defined on Borel subsets~$A\subseteq\mathbb{T}$ by
	\begin{equation*}
	\eta_T (A) = \frac{1}{T}\int_0^T\mathbf{1}_A\left(X_t\right)\,\dd t.
	\end{equation*} 
	For a set~$A$,~$\eta_T(A)$ measures the fraction of time that the process~$X_t$ spends in~$A$. As~$T\to\infty$, the empirical measure converges to the uniform measure~$\dd x$ on~$\mathbb{T}$. We are interested in how the convergence rate depends on the drift~$v$.
	The spectrum of its generator~$L_v=v\nabla + (1/2)\Delta$ is
	\begin{equation*}
	\sigma(L_v) =\left\{-n^2 + inv\,:\,n\in\mathbb{Z}\right\}.
	\end{equation*}
	Hence the spectral gap is~$-1$, and is in particular independent of~$v$. Therefore, the spectral gap does not provide us with any information about how the rate of convergence changes with~$v$. However, for a measure~$\dd\mu(x) = u(x)\, \dd x $ with a smooth and positive density~$u$, the Donsker-Varadhan rate function for the empirical measure is
	\begin{equation*}
	\mathcal{I}_v(\mu) = \frac{1}{8}\int_\mathbb{T}\left|\nabla \log u\right|^2\,\dd\mu + \frac{1}{2} \,v^2 \left(1-\frac{1}{\int_\mathbb{T}\frac{1}{u(x)}\dd x}\right).
	\end{equation*}
	The family~$\{\eta_T\}_{T>0}$ satisfies a large deviation principle with this rate function in the limit~$T\to\infty$. Informally, this means
	\begin{equation*}
	\mathbb{P}\left(\eta_T \approx \mu\right) \sim e^{-T\,\mathcal{I}_v(\mu)},\quad T\to\infty.
	\end{equation*}
	In conclusion, for higher values of~$v$, the empirical measure converges faster to the uniform measure, since rate function increases with increasing~$v$. The limit of~$\eta_T$ is independent of~$v$; we always find~$\mathcal{I}_v(\dd x)=0$.\qed
\smallskip

Evaluation of the large deviation rate function for empirical measures, beyond the variational form given by Donsker and Varadhan, is typically a challenging task. For the diffusion setting, including both reversible and irreversible processes, see \cite{dupuis2018large} and the references therein. In \cite{dupuis2015large} the authors consider reversible jump Markov processes and use stochatic control and weak convergence arguments to derive an explicit form of the rate function. Lastly, in the MCMC context, \cite{rey2015irreversible} consider diffusion processes on a compact manifold where the drift can be decomposed into sufficiently smooth reversible and irreversible parts. The rate function can then be expressed in terms of the rate function of a related reversible diffusion and the solution of an elliptic PDE associated with the irreversible component of the drift.
\smallskip

The proofs of the large deviation results are based on the general Hamilton-Jacobi approach to empirical measures developed by Feng and Kurtz in~\cite[Chapter 12]{FengKurtz2006}. We describe this approach, in the context of this paper, in more detail in Section~\ref{sec:aux}. 
\smallskip

The remainder of this chapter is organised as follows. In Section \ref{sec:prelim} we give the necessary preliminaries: notation and relevant definitions, background on the zig-zag process and empirical measure large deviations. In particular we recall well-known large deviation results for empirical measures by Donsker and Varadhan. The main results are then presented in Section \ref{sec:main}. The section is split into the main assumptions and general large deviation statements (Section \ref{sec:gen_LDP}), large deviation results for the zig-zag process (Section \ref{sec:LDP_zigzag}) and an explicit expression of the rate function associated with the zig-zag process (Section \ref{sec:rate}). All proofs are deferred to Section \ref{sec:proofs}. 

\section{Preliminaries}
\label{sec:prelim}

\subsection{Notation and definitions}
\label{sec:notation}
Throughout this chapter, $E$ will denote a complete separable metric space (Polish space) and $\calB (E)$ the relevant $\sigma$-algebra on $E$; unless otherwise stated this is taken to be the Borel $\sigma$-algebra. $C(E)$ and $C_b(E)$ are the spaces of functions $f : E \to \bR$ that are continuous and bounded continuous, respectively. The space of continuous and right-continuous functions from $[0, \infty)$ to $E$ is denoted by $C_E[0, \infty)$ and $D_E [0,\infty)$, respectively.  A sequence of functions $\{ f_n \}_n$ on $E$ converges {\it boundedly and uniformly on compacts} to a function $f$ if and only if $\sup _n \norm{f_n} < \infty$ and for each compact $K \subseteq E$,
	\begin{align*}
		\lim _{n\to \infty} \sup_{x \in K} | f_n(x) - f(x) | =0.
	\end{align*}
	This is denoted as $f = buc-\lim_{n\to \infty} f_n$.
	
	For a Markov process $Y = \{ Y_t: \ t \geq 0\}$, we denote by $S = \{ S(t): t \geq 0 \}$ the associated Markov semigroup. A semigroup $S(t)$ acting on $C(E)$ is {\it Feller continuous} if, for any $t$, $S(t) : C(E) \to C(E)$, {\it strongly continuous} if $S(t) f \to f$ as $t \to 0$ for any $f \in C(E)$ and {\it buc-continuous} if $buc-\lim _{t\to 0} S(t) f = f$ for $f \in C_b(E)$. 
	
	For an operator $L$, $\calD (L)$ denotes the domain of $L$. For functions in $\calD (L)$, $\calD ^+ (L)$ denotes those that are strictly positive and $\calD^{++}(L)$ those that are positive and uniformly bounded from below by a positive constant. For a given $L$ we use $B$ to denote the extended generator associated with $L$.


We use $\mathcal{P}(E)$ to denote the space of probability measures on $E$, and $\calP _c (E)$ is the subset of probability measures with compact support. Throughout the paper we equip $\calP (E)$ with the topology of weak convergence: $\rho _n \to \rho$ in this topology if
\begin{align*}
	\int _E f(x) \rho _n (dx) \to \int _E f(x) \rho (dx), \ n\to \infty, \ \ \forall f \in C_b (E).
\end{align*}
 A special case that will be considered several times is $\mathcal{P}(D_E [0,\infty))$, which is also equipped with the weak topology. For a process $\{ Y(t), t \geq 0 \}$ taking values in $E$ and $y \in E$, we denote by $\mathbb{P}_y \in \mathcal{P}(D_E [0,\infty) )$ the distribution of the process $Y(t)|_{t \geq 0}$ starting at~$y \in E$.
 
The set of positive Borel measures on $E$ is denoted by $\calM (E)$ and the set of finite Borel measures on $E$ are denoted by $\calM _f (E) \subset \calM (E)$. We let $\calL (E)$ denote the following subset of $\calM (E \times [0,\infty))$:
\begin{align*}
	\calL (E) = \{ z \in \calM (E \times [0, \infty)): \ z(E \times [0,t]) = t, \ t \geq 0 \}. 
\end{align*}
The set $\calL (E)$ is endowed with the topology of weak convergence on bounded time intervals: for $\{ \rho_n \} \subset \calL (E)$, $\rho_n \to \rho$ if for all $f \in C_b(E \times [0, \infty))$ and all $t \geq 0$,
\begin{align*}
	\int _{E \times [0,t]} f(x,s) d\rho_n(x,s) \to \int _{E \times [0,t]} f(x,s) d\rho(x,s).
\end{align*}
Then $\calL (E)$ is the set of Borel-measures on $E \times [0, \infty)$  of the form 
\begin{align*}
	d\rho (x,t) = \mu _t (dx) dt,
\end{align*}
for probability measures $\mu_t \in \calP (E)$. That is , for every $\rho \in \calL (E)$, there exists a measurable path $s \mapsto \mu_s \in \calP (E)$ such that 
\begin{align*}
	\rho (A \times [0,t] ) = \int _0 ^t \mu _s (A)ds, \ \textrm{ for any } A \in \calB (E), \ t >0.
\end{align*}


\subsection{Large deviations for empirical measures}
\label{sec:intro_LDP}
Consider a Markov process $Y = \{ Y_t: t \geq 0 \}$ taking values in a Polish space~$E$, with associated generator $L : \mathcal{D}(B) \subseteq C_b(E) \to C_b(E)$ and semigroup $S(t)$. 
The {\it empirical measure} $\eta_t$ associated with $Y_t$ is the stochastic process with values in $\mathcal{P}(E)$ defined by
\begin{equation*}
\eta_t (A) = \frac{1}{t} \int_0^t \boldsymbol{1}_A(Y_s)ds,\quad A \in \mathcal{B}(E).
\end{equation*}
%
%
Empirical measures play an important role in, for example, the settings of MCMC methods and steady-state simulations, via the pairing of measures and observables: For a probability measure $\mu \in \mathcal{P}(E)$ and a function $V \in C_b(E)$, we write
\begin{equation*}
\mu(V) = \int_E V(y) d\mu(y)
\end{equation*}
for the pairing of measures and observables. For the empirical measure $\eta_t$, this pairing corresponds to time averages,
\begin{equation}
\label{eq:etaV}
\eta_t(V) = \frac{1}{t}\int_0^t V(Y_s) ds.
\end{equation}
%
%
If there is an invariant measure $\pi \in \mathcal{P}(E)$ associated with the generator $L$, ergodicity of the process $Y_t$ will ensure the convergence $\eta _t \to \pi$ as $t\to \infty$, w.p.\ 1 in $\calP (E)$, from which it follows that for any $V \in C_b(E)$,
\begin{align*}
	\eta_t(V) \to \pi(V)	\quad \text{ as } t \to \infty, \; \mathbb{P}-a.s.
\end{align*}
Thus, time averages such as \eqref{eq:etaV} are precisely what is used to form approximations in Monte Carlo methods and there is a direct link between the performance of such simulation methods and the properties of the empirical measure.

The theory of large deviations for empirical measures is concerned with deviations of $\eta _t$ from $\pi$ as $t$ grows large. Recall that the gist of the so-called large deviations principle is that for any $\rho \in \mathcal{P} (E)$, for large $t$
\begin{align*}
	\bP_y (\eta _t \approx \rho) \sim \exp \left\{ -t \, \mathcal{I}(\rho) \right\},
\end{align*}
where the function $\mathcal{I} : \calP (E) \to [0,\infty]$ is the rate function associated with the process. This formula is just a short notation for Definition~\ref{def:LDP}; that means~$\mathcal{I}$ has compact sublevel-sets, and for any measurable subset $A \subseteq \calP (E)$, we have
\begin{align*}
- \inf _{\mu \in \mathrm{int}(A)} I (\mu) &\leq \liminf _{t\to \infty} \frac{1}{t} \log \bP \left(\eta _t \in A \right) \\
&\leq \limsup _{t\to \infty} \frac{1}{t} \log \bP \left( \eta _t \in A \right) \leq -\inf _{\mu \in \mathrm{clos}(A)} I(\mu),
\end{align*}
where $\mathrm{int}(A)$ and $\mathrm{clos}(A)$ are the interior and closure of the set $A$.
\smallskip

Under relatively mild conditions on the dynamics of the process $Y$ the rate function will be strictly convex and satisfy $\mathcal{I}(\mu) =0$ if and only if $\mu = \pi$. Thus, the rate function characterises the exponential rate of decay of probabilities of sets not including the invariant distribution $\pi$. Moreover the rate function can be used to characterise {\it how} events may occur - for sets $A$ that do not include $\pi$, the minimisers of $\mathcal{I}$ over $A$ represent the behaviour $\eta _t$ is most likely to exhibit if $A$  occurs. 

For empirical measures of Markov processes, the rate function associated with an LDP can often be expressed using a variational form, obtained by Donsker and Varadhan~\cite{donsker1975asymptoticI}, involving the generator $L$ of the underlying process. For the compact setting, they proved the following result.

\begin{theorem}[{\cite[Theorem 3]{donsker1975asymptoticI}}]\label{thm:proving_LDP:DV_I}
Take $E$ to be a compact, complete separable metric space. Let $S(t)$ be a Markov semigroup acting on $C(E)$ equipped with the supremum norm, and let $L$ be the generator associated to $S(t)$. Assume the following:
\begin{enumerate}[label =(DV.\arabic*)]
\item\label{item:thm_DV_I:Feller} The semigroup is Feller continuous and strongly continuous. 
\item\label{item:thm_DV_I:reference_measure} There exists a probability measure $\lambda \in \mathcal{P}(E)$ such that for each $t > 0$ and $x\in E$, the transition probabilities $P(t,x,dy)$ are absolutely continuous with respect to $\lambda$, that is 
\begin{equation*}
P(t,x,dy) = p(t,x,y) \lambda(dy),
\end{equation*}
for some $p$ with $0< a(t) \leq p(t,x,y) \leq A(t) < \infty$.
\end{enumerate}
Then the associated sequence $\{\eta_t\}_{t > 0}$ satisfies a large deviation principle in $\mathcal{P}(E)$, with rate function $\mathcal{I}: \mathcal{P}(E) \to [0,\infty]$ given by
\begin{equation}\label{eq:proving_LDP:DV_rate_function}
\mathcal{I}(\mu) = -\inf_{u \in \mathcal{D}^+(L)}\int_E \frac{Lu}{u} d\mu.
\end{equation}
\end{theorem}
%
The theorem applies in particular to drift-diffusions taking values in a compact space. Roughly speaking, for such processes, with reasonable coefficients, the Feller-continuity is satisfied and the diffusive part ensures absolute continuity with respect to a volume measure $dx$. In \cite{rey2015irreversible} Rey-Bellet and Spiliopoulos use this result to study performance of specific irreversible MCMC methods based on drift-diffusions; their Assumption~(H) allows for an application of Theorem \ref{thm:proving_LDP:DV_I}. 

%
Condition~\ref{item:thm_DV_I:reference_measure} is a reasonable transitivity assumption for processes that involve a diffusive term. However, this condition excludes many interesting examples, such as continuous-time jump processes, see e.g.\ \cite{dupuis2015large}. 
The issues highlighted therein are present also for the zig-zag process on $\bR \times \{ \pm 1\}$: in a sense, the absence of a diffusive operator excludes the possibility of finding a suitable reference measure.
%
%

If the process $Y_t$ is reversible with respect to the reference measure, that is $p(t,x,y) = p(t,y,x)$, then the rate function takes a more explicit form, see e.g.\ Theorem~5 in~\cite{donsker1975asymptoticI}. However, our interests are explicitly in irreversible processes, such as the zig-zag process, and therefore such representations are not available. 
%
%

In conclusion, while Theorem \ref{thm:proving_LDP:DV_I} can be a starting point for many drift-diffusion processes, it is not a sufficient tool for many other interesting processes, including the position-velocity PDMPs. In order to use large deviation results to study performance of such MCMC algorithms we must first overcome this obstacle and establish the relevant large deviations principles.



In~\cite{dupuis2018large}, Dupuis and Lipshutz consider large deviations of empirical measures of $\mathbb{R}^d$-valued drift-diffusions. Their Condition~2.2 corresponds to a type of stability criterion in terms of a Lyapunov function. A transitivity property similar to Condition~\ref{item:thm_DV_I:reference_measure} of Theorem~\ref{thm:proving_LDP:DV_I} is satisfied due to the diffusive part, and they prove a different, explicit representation of the rate function, assuming only standard regularity conditions on the coefficients. In particular, this representation holds for irreversible drift-diffusions. 

\subsection{The zig-zag process}
\label{sec:zig-zag}

In this section we will discuss very concisely the zig-zag process. As discussed in the introduction the zig-zag process is an example of a piecewise deterministic Markov process \cite{davis1984piecewise}. As the name indicates, a piecewise deterministic Markov process is a Markov process with deterministic trajectories, in between event times at which the process makes a discontinuous change. 

For the one-dimensional zig-zag process, the state space is either $E = \R \times \{\pm 1\}$ or $E = \mathbb{T} \times \{\pm 1\}$ and a typical state is denoted in this paper by $(x,v)$. Here $x$ represents a position and $v$ a velocity. Starting from $(x,v)$ at time $t= T_0 := 0$, the dynamics of a Markov process $(X_t, V_t)$ are given, until the first (random) event time $T_1 > 0$,  by \[ (X_t, V_t) = (x +t v, v), \quad 0 \leq t < T_1.\] In other words, the position changes according to the constant velocity $v$, which itself does not change in between event times. The random time $T_1$ at which the first event happens is distributed according to
\[ \P_{x,v}(T_1 \geq t) = \exp \left( -\int_0^t \lambda(X_s, V_s) \, d s \right) = \exp \left( -\int_0^t \lambda(x + vs, v) \, d s \right),\]
where $\lambda : E \rightarrow [0,\infty)$ is the \emph{event rate}, which is in the case of the zig-zag process also known as the \emph{switching rate}, which we will discuss in more detail below. At an event time $T$ the velocity changes sign and the position remains unchanged: 
\[ V_{T_1} = -V_{T_1-} \quad \text{and} \quad X_{T_1} = X_{T_1-}.\]
From the time $T_1$ onward, the process repeats the dynamics described above: for $i = 1, 2, \dots$
\begin{align*} & X_t = X_{T_{i-1}}  +(t- T_{i-1}) V_{T_{i-1}}, \quad V_t = V_{T_{i-1}}, \quad T_{i-1}\leq t < T_{i}, \\
& \P(T_i \geq t \mid T_{i-1}, X_{T_{i-1}}, V_{T_{i-1}}) = \exp \left( -\int_{T_{i-1}}^t \lambda(X_{T_{i-1}} + s V_{T_{i-1}}, V_{T_{i-1}}) \, d s \right), \\
& X_{T_i} = X_{T_i-}, \quad V_{T_i} = -V_{T_i-}.
\end{align*}
The switching rate $\lambda : E \rightarrow \R$ is assumed to be continuous. If $\lambda$ satisfies
\begin{equation} \label{eq:switching-intensity-condition-1} \lambda(x, 1) - \lambda(x,-1) = U'(x), \end{equation} for a continuously differentiable function $U$, then the measure defined by
\begin{equation*}
\pi(\dd x, \dd v) = \exp(-U(x)) \, \dd x \otimes \unif_{\pm 1}(\dd v)
\end{equation*}
is a stationary measure for $(X_t, V_t)$.  An equivalent condition to~\eqref{eq:switching-intensity-condition-1} is that for some continuous non-negative function $\gamma(x)$, we have
\begin{equation}
\label{eq:switching-intensity-condition-2}    \lambda(x,v) = \max (0, v U'(x)) + \gamma(x).
\end{equation}
Here $(\max 0, vU'(x))$ is called the \emph{canonical switching intensity}, and $\gamma$ is called the \emph{excess switching intensity} or \emph{refreshment rate}. As a rule of thumb, high values of~$\gamma$ lead to many switches of the velocity.
We study the dependence of the empirical measure of the process~$(X_t, V_t)$ on~$\gamma$ in Section~\ref{sec:rate}.
\smallskip

The zig-zag process can be extended in a natural way to a multi-dimensional process in $\R^d \times \{\pm 1\}^d$ (\cite{bierkens2019zig, bierkens2019ergodicity}). Since we focus in this paper on properties of the one-dimensional process we will not discuss this extension here. The ergodic properties of the zig-zag process are essential in order to establish a large deviation principle for the empirical measure. Under mild conditions it can be shown that the zig-zag process is exponentially ergodic, which is proven in~\cite{BierkensRoberts2017} for the one-dimensional case and in~\cite{bierkens2019ergodicity} for the multi-dimensional zig-zag process. Finally, by \cite[Theorem 26.14]{Davis1993}, the extended generator of the zig-zag process is given by
\[ B f(x,v) = v \partial_x f(x,v) + \lambda(x,v) [ f(x,-v) - f(x,v)], \quad (x,v) \in E,\]
with 
\[ \mathcal D(B) = \{ f : E \rightarrow \R : f(\cdot,v)  \, \text{is absolutely continuous for } v = \pm 1\}.\] 

\section
{Large deviations for empirical measures of PDMPs}
\label{sec:main}

In this section we present our main results: we establish a large deviations principle for the empirical measure of a Markov process under fairly general assumptions which include in particular examples of position-velocity PDMPs such as the zig-zag process. After obtaining these general results we focus for concreteness on the zig-zag process, for which we verify the stated assumptions. We also give an explicit characterisation of the corresponding rate function, a necessary step towards using the LDP for analysing the performance and properties of approximations based on the zig-zag process. To streamline the presentation we split the analysis according to whether we consider a compact or non-compact state space~$E$. 
\smallskip

To facilitate the proof of the LDP for the empirical measures, we first formulate in Section~\ref{sec:gen_LDP} two more general large deviations results (compact and non-compact setting) for empirical measures arising from certain continous-time stochastic processes. We then show that the zig-zag process is a special case in this class of processes in Section \ref{sec:LDP_zigzag}. It is worth to emphasise that we do not aim for greatest generality in the large deviations results Theorems \ref{thm:proving_LDP:LDP_compact} and \ref{thm:proving_LDP:LDP_non_compact}. Rather, we settle for conditions that make the general conditions of Lemma \ref{lemma:FengKurtz} more transparent and concrete whilst still allowing us to prove the large deviations principle for the empirical measures of the zig-zag process.

\subsection{Results aimed at position-velocity PDMPs}
\label{sec:gen_LDP}

Before we specialize to the zig-zag process, we consider the setting described in Section~\ref{sec:notation} to PDMPs: $Y$ is a Markov process taking values in a locally compact complete separable metric space $E$, 
with associated semigroup $S(t)$ and 
infinitesimal generator $L$. We also make use of the extended generator $B$; see \cite{Davis1993,EthierKurtz1986} and Section~\ref{sec:zig-zag}. Typically, $E=\mathbb{R}^d\times\mathcal{S}$ where $\mathbb R^d$ is the state space for a position variable $X_t$ and $\mathcal{S}$ is a compact set that models the state space of the velocity variables $V_t$. For the zig-zag process, $\mathcal{S}=\{\pm 1\}^d$, and for the Bouncy Particle Sampler $\mathcal S$ can be taken to be the $(d-1)$-dimensional unit sphere. Note that for $d =1 $ these two choices coincide.
\smallskip

The following are the assumptions we will impose in order to establish an LDP for the empirical measures of the process $Y$. Not all conditions are required at the same time: we impose conditions~\ref{item:thm_LDP_compact:Feller}, \ref{item:thm_LDP_compact:tight} and~\ref{item:thm_LDP_compact:principal_eigenvalue} for the compact case and \ref{item:thm_LDP_compact:Feller}, \ref{item:thm_LDP_compact:tight}, \ref{item:thm_LDP_non_compact:Lyapunov} and~\ref{item:thm_LDP_non_compact:mixing} for the non-compact case.
\begin{enumerate}[label = (A.\arabic*)]
\item\label{item:thm_LDP_compact:Feller} The semigroup $S(t)$ is a
Feller semigroup.
\item\label{item:thm_LDP_compact:tight} For any compact set $K\subseteq E$, the set of measures $\{\mathbb{P}_y: y \in K\}$ is tight in $\mathcal{P}(D_E[0,\infty))$.
\item\label{item:thm_LDP_compact:principal_eigenvalue} For any function $V \in C(E)$, there exists a function $u \in \mathcal{D}^+(L)$ and a real eigenvalue $\beta \in \mathbb{R}$ such that pointwise on $E$,
\begin{equation*}
(V+L) u = \beta u.
\end{equation*}
\item\label{item:thm_LDP_non_compact:Lyapunov} There exist two non-negative functions $g_1,g_2 \in C(E;[0,\infty))$ such that: 
	\begin{enumerate}
		\item\label{item:thm_LDP_non_compact:Lyapunov1} 
		For any $\ell \geq 0$, the sublevel-sets $\{g_i\leq \ell\}$ are compact and $g_i(y)\to\infty$ as $|y|\to\infty$,
		\item\label{item:thm_LDP_non_compact:Lyapunov2} 
		$g_1(y)/g_2(y) \to 0$ as $|y|\to\infty$,
		\item\label{item:thm_LDP_non_compact:Lyapunov3} $e^{g_i} \in \mathcal D(B)$, for any $c\in\mathbb{R}$ the superlevel-sets $\{y\in E\,:\,e^{-g_i(y)} B(e^{g_i})(y)\geq c\}$ are compact, and $e^{-g_1(y)} B(e^{g_1})(y)\to-\infty$ as $|y|\to\infty$,
	\end{enumerate}
	where we recall that $B$ is the extended generator of $Y$. We write $|y_n|\to\infty$ if $d(y_n,z)\to\infty$ for all points~$z\in E$.
	\item\label{item:thm_LDP_non_compact:mixing} For any two compactly supported probability measures $\nu_1,\nu_2 \in \mathcal{P}_c(E)$, there exist constants $T,M > 0$ and measures $\rho_1,\rho_2 \in \mathcal{P}([0,T])$ such that for all Borel sets $A \subseteq E$,
	\begin{equation}\label{eq:proving_LDP:LDP_non_compact:mixing}
	\int_0^T\int_E P(t,y,A)\, d\nu_1(y)d\rho_1(t) \leq M \int_0^T\int_E P(t,z,A)\, d\nu_2(z)d\rho_2(t),
	\end{equation}
	where $P(t,y,dy')$ denotes the transition probabilities associated to $Y_t$.
\end{enumerate}
Conditions \ref{item:thm_LDP_compact:Feller}-\ref{item:thm_LDP_compact:principal_eigenvalue} are enough to prove Theorem \ref{thm:proving_LDP:LDP_non_compact}, the large deviations principle in a compact setting. In this setting conditions \ref{item:thm_LDP_compact:Feller} and \ref{item:thm_LDP_compact:tight} replace Condition~\ref{item:thm_DV_I:Feller} of Theorem \ref{thm:proving_LDP:DV_I};  Condition \ref{item:thm_LDP_compact:tight} can also be weakened to $\mathbb{P}_{y_n} \to \mathbb{P}_y$ in $\mathcal{P}(D_E[0,\infty))$ whenever $y_n \to y$. Together Conditions \ref{item:thm_LDP_compact:Feller} and \ref{item:thm_LDP_compact:tight} imply strong continuity of the semigroup $S$ (see e.g. \ \cite[Remark~11.22]{FengKurtz2006}). 
\smallskip

As pointed out in Section \ref{sec:intro_LDP}, the processes we have in mind do not satisfy a transitivity condition similar to Condition \ref{item:thm_DV_I:reference_measure} of Theorem \ref{thm:proving_LDP:DV_I}. In the compact setting this can be replaced by condition~\ref{item:thm_LDP_compact:principal_eigenvalue}, which corresponds to a principal-eigenvalue problem for the operator $L+V$. In compact settings, such eigenvalue problems can usually be solved if the coefficients of the generator are regular enough. In Section \ref{sec:LDP_zigzag} we show that this is the case for the zig-zag process taking values in the compact torus. 
\smallskip

In the non-compact setting, the eigenvalue problem \ref{item:thm_LDP_compact:principal_eigenvalue} is replaced by conditions \ref{item:thm_LDP_non_compact:Lyapunov} and \ref{item:thm_LDP_non_compact:mixing}.
%
Condition~\ref{item:thm_LDP_non_compact:Lyapunov} is closely related to the stability conditions assumed in~\cite{donsker1976asymptotic} and~\cite{dupuis2018large}. Because $e^{g_1}$ is unbounded, formally we have to use the extended generator $B$ instead of the infinitesimal generator $L$ to formulate Condition~\ref{item:thm_LDP_non_compact:Lyapunov3}. The same problem occurs in Condition~2.2 of~\cite{dupuis2018large}: for a diffusion process $Y_t$ in $\mathbb{R}^d$ satisfying
$
dY_t = -Y_t dt + dW_t,
$
the second-order differential operator 
\begin{equation*}
Bf (x) = \frac{1}{2} \Delta f(x) - x\nabla f(x)
\end{equation*}
acting on $C^2(\mathbb{R}^d)$ is well-defined and is equal to the infinitesimal generator $L$ of the process when restricted to $C^2_b(\mathbb{R}^d)$. With $g_1(x) = \delta |x|^2/2$, the function $e^{-g_1(x)}Be^{g_1}(x)$ goes to minus infinity for $\delta$ small enough. A second Lyapunov function is $g_2(x) = \sqrt{1+|x|^2}$. In the context of the zig-zag process, since $E$ is of the form $E=\mathbb{R}^d\times\{\pm 1\}^d$, using continuous functions that grow to infinity when fixing the velocity variable is sufficient for obtaining compact level sets.
\smallskip

Condition~\ref{item:thm_LDP_non_compact:mixing} plays the role of a transitivity assumption in the non-compact case. While it is feasible to solve a principal-eigenvalue problem for a compact state space, this is much more difficult in the non-compact setting. It would require deriving not only the eigenvalue itself, but also the corresponding eigenfunction on a non-compact space, for which general existence results are not available. 
In this setting the transitivity condition \ref{item:thm_DV_I:reference_measure} of Theorem \ref{thm:proving_LDP:DV_I} is instead partly replaced by the mixing property \ref{item:thm_LDP_non_compact:mixing}. It is a weakened version of~\cite[Condition B.8]{FengKurtz2006},which is based on~\cite[Condition~$\tilde{U}$, page~113]{DeuschelStroock1989}. It is weaker in that it requires the transition probabilities to be comparable only for compactly supported initial conditions~$\nu_1,\nu_2 \in \mathcal{P}_c(E)$. This weakening is crucial for the results in this paper, because the stronger condition fails to be true for the zig-zag process if for instance $\nu_1 = \mathcal{N}(0,1) \otimes \unif_{\pm 1}$ and $\nu_2 = \delta_y$. In that example, while the left-hand side of \eqref{eq:proving_LDP:LDP_non_compact:mixing} is in this case positive for any Borel set $A \subseteq E$, the right-hand side can become zero. This is because the zig-zag process has finite speed propagation, so that for arbitrary $T > 0$, if $\mathrm{dist}(A,y) > T$, then the probability of transitioning from $y$ into $A$ is zero. However for compactly supported measures the condition is satisfied. We verify the conditions of Theorem \ref{thm:proving_LDP:LDP_non_compact} for the zig-zag process in Section \ref{sec:LDP_zigzag}. 
\smallskip

We are now ready to state the two general large deviations results of this paper, which in Section \ref{sec:LDP_zigzag} will be used to derive the large deviations principle for the empirical measures of the zig-zag process. We start with the compact setting.
\begin{theorem}\label{thm:proving_LDP:LDP_compact}
Let $E$ be compact, $S(t)$ a Markov semigroup acting on $C(E)$ equipped with the supremum norm, and $Y_t$ the corresponding Markov process. Let $L$ be the infinitesimal generator of $Y_t$, and assume that $Y_t$ solves the associated martingale problem. Suppose Assumptions \ref{item:thm_LDP_compact:Feller}, \ref{item:thm_LDP_compact:tight} and \ref{item:thm_LDP_compact:principal_eigenvalue} hold. Then 
the empirical measures $\{\eta_t\}_{t > 0}$ associated to $Y_t$ satisfy a large deviations principle in $\mathcal{P}(E)$ with rate function $\mathcal{I}:\mathcal{P}(E) \to [0,\infty]$ given by~\eqref{eq:proving_LDP:DV_rate_function}.
\end{theorem}

Theorem \ref{thm:proving_LDP:LDP_compact} remains valid when replacing the eigenvalue-problem condition~\ref{item:thm_LDP_compact:principal_eigenvalue} by the mixing condition~\ref{item:thm_LDP_non_compact:mixing}. This is because the latter is a weaker condition sufficient for verifying the inequality \eqref{eq:ineqH1H2}, upon which the proof of the theorem hinges.

The next theorem gives the corresponding large deviations result for the non-compact setting; this is the result we use for proving the large deviations principle for the zig-zag process on $\bR \times \{\pm 1\}$ (Theorem \ref{thm:LDP-zigzag:non-compact-1d}).
\begin{theorem}\label{thm:proving_LDP:LDP_non_compact}
Let 
$S(t)$ be a Markov semigroup acting on $C_b(E)$ and $Y_t$ the corresponding Markov process. Let $L$ be the infinitesimal generator of $Y_t$ and assume that $Y_t$ solves the associated martingale problem. Assume \ref{item:thm_LDP_compact:Feller}, \ref{item:thm_LDP_compact:tight},
\ref{item:thm_LDP_non_compact:Lyapunov} and \ref{item:thm_LDP_non_compact:mixing}.
Then, if $Y_0 \in K$ for some compact set $K$, the empirical measures $\{\eta_t\}_{t > 0}$ associated to the Markov process $Y_t$ satisfy a large deviations principle in $\mathcal{P}(E)$, with rate function $\mathcal{I} : \mathcal{P}(E) \to [0,\infty]$ given by
\begin{equation*}
\mathcal{I}(\mu) = -\inf_{u \in \mathcal{D}^{++}(L)} \int_E \frac{Lu}{u} d\mu.
\end{equation*}
%
\end{theorem}

The proofs of Theorems \ref{thm:proving_LDP:LDP_compact} and \ref{thm:proving_LDP:LDP_non_compact} are given in Sections ~\ref{section:proof-of-LDP:compact} and \ref{section:proof-of-LDP:non-compact}, respectively.

%
\subsection{The empirical measures of the zig-zag process}
\label{sec:LDP_zigzag}
Having established the general large deviations results Theorems \ref{thm:proving_LDP:LDP_compact} and \ref{thm:proving_LDP:LDP_non_compact}, we now specialize to the zig-zag process. Throughout the section, $Y_t$ is used to denote the zig-zag process, $Y_t = (X_t, V_t)$ with $X_t$ and $V_t$ as in Section \ref{sec:zig-zag}. However the state space $E$ will change as we split the large deviations statements for the empirical measures of $Y$ into compact (torus) and non-compact ($\bR$) settings. Although \ref{thm:proving_LDP:LDP_non_compact} holds for arbitrary dimension $d \geq 1$, for the zig-zag process we limit ourselves to verifying the conditions for the case $d = 1$. Extending these results to $d > 1$ is substantially more difficult and remains a topic of further research. While conditions~\ref{item:thm_LDP_compact:Feller}, \ref{item:thm_LDP_compact:tight}, \ref{item:thm_LDP_non_compact:Lyapunov} hold true, the main challenge is verifying~\ref{item:thm_LDP_non_compact:mixing}.


We begin by considering the compact state space  $\mathbb{T} \times\{\pm 1\}$. In this case the infinitesemal generator $L$ of the semigroup $S(t)$ is has domain $\mathcal{D}(L) = C^1(\bT \times \{ \pm 1\}) = \{f\in C(\bT \times \{ \pm 1\}): f(\cdot,\pm 1) \in C^1(\mathbb{T})\}$,
and takes the form
\begin{equation}
\label{eq:gen_compact}
Lf (x,v) = v \partial_x f(x,v) + \lambda(x,v) \left[f(x,-v) - f(x,v)\right],
\end{equation}
with $\lambda$ given by \eqref{eq:switching-intensity-condition-2}. The LDP for the empirical measures associated with $Y$ and this state space is given in Theorem \ref{thm:LDP-zigzag:compact}. We prove this result in Section~\ref{section:zig-zag:LDP:compact} by verifying the conditions of Theorem \ref{thm:proving_LDP:LDP_compact}, the large deviations principle for processes taking values in a compact state space.


\begin{theorem}\label{thm:LDP-zigzag:compact}
Suppose that $U \in C^2(\mathbb{T})$. Then the family of empirical measures $\{\eta_t\}_{t>0}$ of the zig-zag process taking values in $\bT \times \{\pm 1\}$ satisfies a large deviations principle in the limit $t \to \infty$, with rate function $\mathcal{I}:\mathcal{P}(\bT \times \{\pm 1\}) \to [0,\infty]$ given by
\begin{equation*}
\mathcal{I}(\mu) = -\inf_{u \in \mathcal{D}^+(L)} \int_{\bT \times \{\pm 1\}} \frac{Lu}{u} d\mu.
\end{equation*}
%
%
\end{theorem}
We now move to the setting of a non-compact state space. Specifically, we consider the zig-zag process $Y_t =(X_t,V_t) $ taking values in $\bR \times \{ \pm 1 \}$. As before, $L$ is the generator of this process, i.e.\ $L : \mathcal{D}(L) \subseteq C_b(\bR \times \{ \pm 1\}) \to C_b(\bR \times \{ \pm 1\})$ is a densely defined linear operator, on the set of functions $\{ f(\cdot, \pm 1) \in C_b ^1(\bR) \}$ we have the representation
\begin{align*}
	Lf(x,v) = v \partial _x f(x,v) + \lambda(x,v) \left[ f(x, -v) - f(x,v) \right], \ \ f \in \calD (L),
\end{align*}
with $\lambda(x,v)=\max(0,vU'(x))+\gamma(x)$. To prove the large deviations principle in this non-compact setting we need additional assumptions on the potential function $U$ determining the jump rates.
\begin{enumerate}[label =($B$.\arabic*)]
	\item\label{item:thm-LDP-zig-zag:Ui} $U(x) \to \infty$ as $|x| \to \infty$ and $U'(x) \to \pm \infty$ as $x \to \pm \infty,$
	\item\label{item:thm-LDP-zig-zag:Uii} $U'(x)/U(x) \to 0$ as $|x| \to \infty,$
	\item\label{item:thm-LDP-zig-zag:Uiii} $U''(x)/U'(x) \to 0$ as $|x| \to \infty.$
\end{enumerate}
Furthermore, we will assume that there exists a second potential $V \in C^2(\mathbb{R})$ such that:
\begin{enumerate}[label =($C$.\arabic*)]
	\item\label{item:thm-LDP-zig-zag:Upsi} $V(x) \to \infty$ as $|x| \to \infty$ and $V'(x) \to \pm \infty$ as $x \to \pm \infty$,
	\item\label{item:thm-LDP-zig-zag:Upsii} $V(x)/U(x) \to 0$, $U'(x)/V(x) \to 0 $ and $V'(x)/U'(x) \to 0$ as $|x| \to \infty$,
	\item\label{item:thm-LDP-zig-zag:Upsiii} $U''(x)/V'(x) \to 0$ as $|x| \to \infty$.
\end{enumerate}
In Section~\ref{section:zig-zag:LDP:non_compact}, we prove the following Theorem.
\begin{theorem}\label{thm:LDP-zigzag:non-compact-1d}
Assume that $U\in C^3(\mathbb{R})$ satisfies \ref{item:thm-LDP-zig-zag:Ui} - \ref{item:thm-LDP-zig-zag:Uiii}, that there is a function $V \in C^2(\mathbb{R})$ satisfying \ref{item:thm-LDP-zig-zag:Upsi} - \ref{item:thm-LDP-zig-zag:Upsiii} and the function $\gamma$ in \eqref{eq:switching-intensity-condition-2} is uniformly bounded by some $\bar \gamma$.
%
%
Suppose that the initial condition $Y_0$ belongs to a compact set $K \subseteq \bR \times \{ \pm 1\}$. Then the empirical measures $\{\eta_t\}_{t > 0}$ of $Y$ satisfies a large deviations principle on $\mathcal{P}(\bR \times \{ \pm 1\})$ with speed $t$ and rate function $\mathcal{I}:\mathcal{P}(\bR \times \{ \pm 1\}) \to [0,\infty]$ given by
\begin{equation*}
\mathcal{I}(\mu) = -\inf_{u \in \mathcal{D}^{++}(L)} \int_{\bR \times \{ \pm 1\}} \frac{Lu}{u} d\mu.
\end{equation*}
%
\end{theorem}
Some comments on the additional assumptions \ref{item:thm-LDP-zig-zag:Ui} - \ref{item:thm-LDP-zig-zag:Uiii} and \ref{item:thm-LDP-zig-zag:Upsi} - \ref{item:thm-LDP-zig-zag:Upsiii} are in place. The condition $U \in C^3(\mathbb{R})$ is imposed to allow for an application of Theorem 4 of \cite{bierkens2019ergodicity}, which is used to verify that \ref{item:thm_LDP_non_compact:mixing} holds. 
The auxiliary potential $V$ is used to find a second Lyapunov function for $L$ that grows slower than $U$ at infinity; roughly speaking, $V$ behaves asymptotically in-between the potential $U$ and its derivative $U'$ as $|x|$ grows. As an example, in the Gaussian case, $U(x) = x^2/2$ satisfies Conditions \ref{item:thm-LDP-zig-zag:Ui} - \ref{item:thm-LDP-zig-zag:Uiii}, and for any $0<\kappa<1$, the potential $V(x) = |x|^{1+\kappa}/(1+\kappa)$ satisfies \ref{item:thm-LDP-zig-zag:Upsi} - \ref{item:thm-LDP-zig-zag:Upsiii}. In general, any potential $U$ growing at infinity as $(1+|x|^2)^{\beta/2}$ with $\beta > 1$ satisfies the conditions, with $0<\kappa < 1$ such that $\beta - \kappa > 1$ and auxiliary potential $V(x) \sim (1+|x|^2)^{(\beta - \kappa)/2}$.

\subsection{Explicit expression for the rate function}
\label{sec:rate}
In Theorems \ref{thm:LDP-zigzag:compact} and \ref{thm:LDP-zigzag:non-compact-1d} we establish the LDP for the empirical measures of the zig-zag process taking values in $\bT \times \{ \pm 1\}$ and $\bR \times \{ \pm 1 \}$, respectively. In those results the rate function is given on the variational form of the results by Donsker and Varadhan, see Section \ref{sec:intro_LDP}. This form follows from the more general large deviations results in Section \ref{sec:gen_LDP} and are not specific to the zig-zag process. Here, we use the properties of the latter to derive a more explicit form of the rate function for the case $E = \mathbb{T} \times \{\pm 1\}$, taking a first step towards using it as a tool for analysing the corresponding simulation algorithms. 


We assume throughout that $E = \mathbb{T} \times \{\pm 1\}$ and the switching rate $\lambda(x,v)$ satisfies $\lambda(x,v) > 0$ for all $(x,v) \in E$. This does not include the canonical rates $\lambda(x,v) = \max(0, v U'(x))_+$; however at the end of this section we present a formal expression for this case. 

Define a reference measure $\nu_0$ on $E$ by $\nu_0(dx, dv) = \Leb(dx) \otimes \unif_{\pm 1}(dv)$.
For any function $f : E \rightarrow \R$ we write $f^+(x) := f(x,+1)$ and $f^-(x) := f(x,-1)$.
Recall the $\arcsinh$ function,
\[ \arcsinh(\xi) = \log \left( \xi + \sqrt{\xi^2 +1} \right), \quad \xi \in \R.\]
The proofs of the following results are given in Section~\ref{sec:proofs-rate-function}.

\begin{proposition}
\label{prop:expression-ratefunction}
Suppose $\mu(dx, dv) = \rho(x,v) \nu_0(dx, dv)$ for a continuously differentiable function $\rho : E \rightarrow [0,\infty)$.
If $\frac{d \rho^+}{dx}(x)= \frac{d \rho^-}{dx}(x)$ and $\rho^+$, $\rho^-$ are strictly positive for all $x \in \mathbb{T}$ then the Donsker-Varadhan functional is given by
\begin{align}
\label{eq:rate-function}
\nonumber  \mathcal I(\mu)& = \int_{\mathbb{T}} \bigg\{\tfrac 1 2 \rho' \log \left( \frac{\lambda^+ \rho^+}{\lambda^- \rho^-} \right)+ \rho'  \arcsinh \left( \frac{ \rho'}{2 \sqrt{ \lambda^+ \lambda^- \rho^+ \rho^-}}\right) \\
 & \quad \quad \quad- \sqrt{ 4 \lambda^+ \lambda^- \rho^+ \rho^- + (\rho')^2} + \lambda^+ \rho^+ + \lambda^- \rho^- \bigg\} \ \dd x.
 \end{align}
If $\rho^+ \geq 0$ and $\rho^- \geq 0$ are constant, then
\begin{equation}
\label{eq:rate-function-constant-density}
\nonumber  \mathcal I(\mu) = \int_{\mathbb{T}} \left( \sqrt{\lambda^+ \rho^+} - \sqrt{\lambda^- \rho^-} \right)^2 \ \dd x.
 \end{equation}
If $\frac{d \rho^+(x)}{dx}(x) \neq \frac{d \rho^-(x)}{dx}$ for some $x \in \mathbb{T}$ then $\mathcal I(\mu) = \infty$.
\end{proposition}


Note that if $\mu(dx,dv) =  \rho(x,v) \nu_0(dx, dv)$ and $\frac{d \rho^+}{dx}= \frac{d \rho^-}{dx}$ on $\mathbb{T}$, then for some constant $c \in \R$ and a probability density function $\rho$ on $\mathbb{T}$ we have $\rho^+(x,v) = \rho(x) + c$ and $\rho^-(x,v) = \rho(x) - c$. 
A useful application of the rate function~$\mathcal I(\mu)$ is in estimating deviations of ergodic averages, which typically requires the computation of 
\[
\inf_{\mu \in \mathcal P} \left(\mathcal I(\mu) - \int_E V \ \dd \mu \right).
\] 
The rigorous statement is the Laplace principle~\cite[Definition~1.6, Theorem~1.8]{BudhirajaDupuis2019}.
If the function $V$ does not depend on $v$, then by the following result we can safely assume $c = 0$ and thus restrict the minimization problem to minimization over probability densities on $\mathbb{T}$.
\begin{proposition}
\label{prop:no-v-dependence}
Let $\rho \in C^1(\mathbb{T})$ be a strictly positive probability density function on $\mathbb{T}$. Let $|k| := \inf_{x} \rho(x)$. Consider the one-parameter family of probability measures $(\mu_c)_{c \in (-k,+k)} \in \mathcal P_{eq}$ with probability density functions $\rho_c : E \rightarrow (0, \infty)$ given by
\[ \rho_c(x,+1) = \rho(x) + c, \quad \rho_c(x,-1) = \rho(x) - c, \quad c \in (-k, +k).\]
Then $c \mapsto \mathcal I(\mu_c)$ is minimized at $c = 0$. Furthermore, for $\mu = \mu_0$,  
\begin{equation}
\label{eq:ratefunction-simplified}
\begin{aligned}
  \mathcal I(\mu) 
 &  = \int_{\mathbb{T}} \bigg\{ \tfrac 1 2 \rho' \log \left( \frac{ \lambda^+}{\lambda^-} \right) + \rho' \arcsinh \left( \frac{\rho'}{2 \rho \sqrt{ \lambda^+ \lambda^-}} \right) \\
 & \quad \quad \quad  - \sqrt{ 4 \lambda^+ \lambda^- \rho^2 + (\rho')^2} + (\lambda^+ + \lambda^-) \rho \bigg\} \ \dd x.
 \end{aligned}
\end{equation}
\end{proposition}

We will specialize to the case in which $\rho = \rho^+ = \rho^-$ and use the representation $\rho = \exp(-W)$, where $W \in C^1(\mathbb{T})$. We  then find 
\begin{equation}
\label{eq:in-W}
\begin{aligned}
\mathcal I(\mu)
& = \int_{\mathbb{T}} \bigg\{ -\tfrac 1 2  W' \log \left( \frac{\lambda^+}{\lambda^-} \right) + W' \arcsinh \left( \frac{ W'}{2 \sqrt{\lambda^+ \lambda^-}} \right) \\
& \quad \quad \quad - \sqrt{4 \lambda^+ \lambda^- + (W')^2} + \lambda^+ + \lambda^- \bigg\} \exp(-W) \ \dd x.
\end{aligned}
\end{equation}
Let 
\begin{equation} \label{eq:switching-rate} \lambda^+(x) = \gamma + \max(0, U'(x)) \quad \mbox{and} \quad \lambda^-(x) = \gamma + \max(0, - U'(x)),\end{equation}
where $\gamma > 0$ is constant, so that $\lambda^{\pm}$ satisfy~\eqref{eq:switching-intensity-condition-2} and hence the measure with $\nu_0$-density $\exp(-U(x))$ is invariant. We call $\gamma$ the \emph{excessive switching intensity} or \emph{refreshment rate}.

We can now investigate the dependence of the rate function $\mathcal I$, through the expression~\eqref{eq:in-W}, on $\gamma$. 
The derivative of the integrand of~\eqref{eq:in-W}  with respect to $\gamma$ can be computed to be
\begin{equation} \label{eq:gamma-derivative} \left( \frac{4 \lambda^+ \lambda^- +(\lambda^+ - \lambda^-)W' - (\lambda^+ + \lambda^-)\sqrt{4 \lambda^- \lambda^+ + (W')^2}}{2 \lambda^+ \lambda^-}\right)\exp(-W),\end{equation}
which is non-positive, and zero only if $W' = \lambda^+ - \lambda^- = U'$ (which can be seen by maximizing with respect to $W'$).
It follows that $\mathcal I(\mu)$ is strictly decreasing as a function of $\gamma$ for $\mu$ not equal to the stationary measure. In other words, for a smaller refreshment rate $\gamma$, the rate function increases. Hence the convergence of empirical averages to equilibrium is faster for smaller~$\gamma$.

Suppose that $\nu_0 (\{ x \in \mathbb{T} : U'(x) = 0\}) = 0$, i.e. the set of points where the derivative of $U$ vanishes is $\nu_0$-negligible. 
In the formal limit~$\gamma \downarrow 0$ in~\eqref{eq:switching-rate}, we obtain the following expression for the rate function:
\begin{equation}
 \label{eq:limit-functional}
 \begin{aligned}
 \mathcal I(\mu) & = \begin{cases} \int_{\mathbb{T}} \left\{ |W'| \left( \log\left( \frac{ W'}{U'} \right) - 1 \right) + |U'| \right\}  \exp(-W) \ \dd x  \quad & \mbox{if $\operatorname{sign}(W') \equiv \operatorname{sign}(U')$}, \\
              \infty \quad & \mbox{otherwise}.
             \end{cases}
  \end{aligned}
\end{equation}


\section{Proofs}
\label{sec:proofs}

\subsection{General large-deviation Theorems \ref{thm:proving_LDP:LDP_compact} and \ref{thm:proving_LDP:LDP_non_compact}}
\label{sec:aux}
In this section we give the proofs of Theorems \ref{thm:proving_LDP:LDP_compact} and \ref{thm:proving_LDP:LDP_non_compact}, which are used to obtain the large deviations principle for the empirical measures of the zig-zag process. The case of a compact state space is treated in Section \ref{section:proof-of-LDP:compact} and the non-compact case in Section \ref{section:proof-of-LDP:non-compact}. Before we embark on these proofs we outline the overall strategy; a more detailed description can be found in the book by Feng and Kurtz~\cite[Chapter~12]{FengKurtz2006}. 

\smallskip

Consider the empirical measure
\begin{equation*}
\eta_t(\cdot) = \frac{1}{t} \int_0^t \delta_{Y_s} (\cdot) \dd s.
\end{equation*}
With a change of variable $s \mapsto ts$ in the integral we can express this as
\begin{align*}
	\eta_t (\cdot) = \int _0 ^1 \delta _{Y_{st} } (\cdot) \dd s,
\end{align*}
the empirical measure for the sped-up process (we can think of $t>1$) $Y^t _s = Y_{st}$ over the time interval $[0,1]$; in fact we will use $t=n \in \bN _+$ below. We can consider the empirical measure of this time-scaled process $Y_t$ on time intervals of lengths other than unity: for $\tau >0$ define $\eta _t ^{\tau}$ as
\begin{align*}
	\eta _t ^{\tau} (\cdot) = \int _0 ^{\tau} \delta _{Y_{st}} \dd s.
\end{align*}
This empirical measure is viewed as an element of $\calL (E)$, the set of Borel measures on $E \times [0, \infty)$ of the form $\dd\rho (x,s) = \mu_s (\dd x)\dd s$, $\mu _s \in \calP (E)$ (see Section \ref{sec:notation}). Any such $\rho \in \calL (E)$ defines a continuous path $t \mapsto \rho _t = \rho (\cdot \times [0,t]) \in \calM _f (E)$ and for $t=1$ this is a probability measure. 

The strategy for proving the large deviations principle for  $\{ \eta _t \}$ is to first show that $\{ \eta _t ^{\tau} \}$ satisfies a large deviations principle in $\calL (E)$. We can then use the fact that projections are continuous maps on $\calL (E)$ (Lemma~\ref{lemma:proofs_LDP:continuous_projection}) and an application of the contraction principle to obtain the sought-after large deviations principle on $\calP (E)$. This is summarised in the following proposition.
\begin{proposition}
\label{prop:LDP_strategy}
	Suppose that the family $\{\eta^\tau_t\}_{t > 0}$ satisfies a large deviations principle in $\mathcal{L}(E)$ with rate function $\mathcal{J}:\mathcal{L}(E) \to [0,\infty]$ given by
\begin{equation*}
\mathcal{J}(\rho) = \int_0^\infty \mathcal{I}(\mu_s) \dd s, \quad \text{ for } \, \rho_t = \int_0^t \mu_s \dd s,
\end{equation*}
where $\mathcal{I} : \mathcal{P}(E) \to [0,\infty]$ is the rate function appearing in the Donsker-Varadhan results,
\begin{equation*}
\mathcal{I}(\mu) = -\inf_{u \in \mathcal{D}^{++}(L)} \int_E\frac{Lu}{u} \dd\mu.
\end{equation*}
Then $\eta_t $ satisfies a large deviations principle in $\mathcal{P}(E)$ with rate function $\mathcal{I}$.
\end{proposition}
\begin{proof}
	Because $\{\eta^\tau_t\}_{t > 0}$ satisfies a large deviations principle on $\mathcal{L}(E)$ and the projection $\pi_1 : \mathcal{L}(E) \to \mathcal{P}(E)$ given by 
	\begin{align*}
		\pi _1 (\rho) = \rho _1,
	\end{align*}
	is continuous (Lemma~\ref{lemma:proofs_LDP:continuous_projection} below), by the contraction principle the sequence evaluated at~$\tau=1$,~$\{\eta_t ^1\}_{t > 0}$, satisfies a large deviation principle on~$\mathcal{P}(E)$ with rate function~$\tilde{\mathcal{I}}: \mathcal{P}(E) \to [0,\infty]$ given by
\begin{align*}
\tilde{\mathcal{I}} (\nu) &= \inf\left\{
\mathcal{J}(\rho) = \int_0^\infty \mathcal{I}(\mu_s) \, \dd s
:
\rho_t = \int_0^t \mu_s \, \dd s \in \mathcal{L}(E), \rho_1 = \nu
\right\}.
\end{align*} 
It remains to show that $\tilde{\calI} (\nu) = \calI  (\nu)$ for every $\nu \in \calP (E)$. First, in the integral defining $\calJ (\rho)$, the integrand is always positive after time $t=1$. It is therefore enough to consider only integrating to time $t=1$ in the infimum, as we are free to chose the the form of $\rho$ after that time. Thus,
\begin{align*}
	\tilde{\calI} (\nu) = \inf \left\{ \int _0 ^1 \calI  (\mu _t) \dd t : \rho _t = \int _0 ^t \mu_s \dd s \in \calL (E), \ \rho_1 = \nu \right\}.
\end{align*}

For a fixed $\nu \in \calP (E)$, take any $\rho_t = \int _0 ^t \mu _s ds \in \calL (E)$ such that $\rho _1 = \nu$. The rate function $\calI $ is convex on $\calP (E)$ and by Jensen's inequality we have
\begin{align*}
	\calI  (\nu) &= \calI  \left( \rho _1 \right) \\
	&=  \calI  \left( \int _0 ^1 \mu _s \dd s \right) \\
	&\leq \int _0 ^1  \calI  \left( \mu _s \right)\dd s.
\end{align*}
Taking the infimum over all such $\rho \in \calL(E)$ yields the inequality
\begin{align*}
	\calI  (\nu) \leq \tilde{\calI} (\nu),
\end{align*}
The constant path $\mu _s = \nu$ gives equality and we have that $\tilde{\calI} = \calI $ as functionals on $\calP(E)$.
\end{proof}
\begin{lemma}[Projection is continuous]\label{lemma:proofs_LDP:continuous_projection}
Let $\mathcal{L}(E)$ be the above space with the topology of weak convergence on bounded time intervals. Let $\mathcal{P}(E)$ be equipped with the weak topology. Then the projection $\pi_1 : \mathcal{L}(E) \to \mathcal{P}(E)$ defined by $\pi_1 (\rho) := \rho_1$ is a continuous map.
\end{lemma}
\begin{proof}[Proof of Lemma~\ref{lemma:proofs_LDP:continuous_projection}]
Let $\rho^n \to \rho$ in $\mathcal{L}(E)$. We need to prove that for any bounded and continuous function $g$ on $E$, we have
\begin{equation*}
\int_E g(u) \dd\rho^n_1(u) \to \int_E g(u) \dd\rho_1.
\end{equation*}
Since
\begin{equation*}
\int_E g(u) \dd\rho^n_1(u) = \int_{E \times [0,1]} g(u) \dd\rho^n(u,s),
\end{equation*}
and~$\varphi(u,s) = g(u)$ is continuous and bounded on~$E \times [0,\infty)$, this is implied by~$\rho^n \to \rho$.
\end{proof}
Armed with Proposition \ref{prop:LDP_strategy}, one way to prove Theorems \ref{thm:proving_LDP:LDP_compact} and \ref{thm:proving_LDP:LDP_non_compact} is to prove the large deviations principle for the empirical measures of the associated sped-up versions of underlying processes and apply the proposition. This is the approach we take and we rely on the following result from \cite{FengKurtz2006} for proving the large deviations principles on $\calL (E)$.
\begin{lemma}[Theorem~12.7 of~\cite{FengKurtz2006}]
\label{lemma:FengKurtz}
Suppose that the following conditions hold:
\begin{enumerate}[label=(FK.\arabic*)]
	\item The martingale problem for $L$ is well-posed. \label{cond:FK1}
	\item The semigroup $S$ is Feller-continuous. \label{cond:FK2}
	\item The semigroup $S$ is $buc$-continuous. \label{cond:FK3}
	\item There is an index set $\calQ$ and a family of subsets of $E$, $\{ \tilde{K} ^q _n \subset E : q \in \calQ\}$, such that for $q_1, q_2 \in \calQ$, there exists $q_3 \in \calQ$ with $\tilde{K} ^{q_1} _n \cup \tilde{K} ^{q_2} _n \subset \tilde{K} ^{q_3} _n$, and for every $y \in E$, there exists $q \in \calQ$ such that $\lim _{n\to \infty} d (y, \tilde{K} ^{q} _n) =0$. Moreover, for each $q \in \calQ$, $T>0$ and $a >0$, there exists a $\hat q (q,a,T) \in \calQ$ satisfying 
	\begin{align*}
		\limsup _{n \to \infty} \sup _{y \in \tilde{K} ^q _n} \frac{1}{n} \log \bP_{y} \left( Y_t \notin \tilde{K} ^{\hat q (q,a,T)} _n, \ \textrm{some } t \leq nT \right) \leq -a.
	\end{align*}
	\label{cond:FK4}
	\item There exists an upper semicontinuous function $\Psi$ on $E$, $\{ \varphi _n \} \subset \calD ^{++} (B_0)$, and $q_0 \in \calQ$ such that $\Psi$ is bounded above, $\{ y\in E: \Psi(y) \geq c \}$ is compact for each $c \in \bR$, $0 < \inf _{y \in K _n ^{q_0}} \varphi _n (y) < 2 \inf _{y\in E} \varphi _n(y)$, $\inf_{n, y \in E} \varphi_n (y) >0$,
	\begin{align*}
		\lim _{n\to \infty} \frac{1}{n} \log \norm{\varphi _n} =0, \ \ \sup_{n,y} \frac{L \varphi _n (y) }{ \varphi_n (y)} < \infty,
	\end{align*} 
	and for each $q \in \calQ$,
	\begin{align*}
		\lim _{n \to \infty} \sup _{y \in \tilde{K} ^q _n} \left( \frac{L \varphi _n (y) }{ \varphi_n (y)} - \Psi (y) \right) \leq 0, \ \ q\in \calQ.
	\end{align*}
	In addition, for each $n$ and $\beta \in (-\infty, 1]$,
	\begin{align*}
		\lim_{t\to 0} \norm{S(t) \varphi_n^{\beta} - \varphi_n^{\beta}} =0.
	\end{align*}
	\label{cond:FK5}
	\item For each $a >0$ there exists compact $K$ and $q \in \calQ$ such that
	\begin{align*}
		\limsup _{n \to \infty} \frac{1}{n} \log \bP \left( Y^n _0 \notin K \cap \tilde{K} ^{q} _n  \right) \leq -a.
	\end{align*}
	\label{cond:FK6}
	\item Take $\calC \subset C_b (E)$ separating and define, with $\Psi$ as in \ref{cond:FK5}, 
	\begin{align*}
	H_1 ^{\beta, \Psi}& = \inf_{0 < \kappa \leq 1} \inf_{f \in \mathcal{D}^{++}(L)} \sup_{y \in E} \left[ \beta (y) \cdot p + (1-\kappa) \frac{Lf(y)}{f(y)} + \kappa \Psi(y)\right],\\
	H_2 ^{\beta, \Psi} &= \sup_{\kappa > 0} \sup_{f \in \mathcal{D}^{++}(L)} \inf_{y \in E} \left[ \beta (y) \cdot p + (1+\kappa) \frac{Lf(y)}{f(y)} - \kappa \Psi(y)\right].
	\end{align*}
	It holds that $H_1 ^{\beta, \Psi} \leq H_2 ^{\beta, \Psi}$, 
	for $\beta \in \calC ^d$, $d=1,2.\dots$.  
	\label{cond:FK7}
\end{enumerate}
Then $\{ \eta ^{\tau} _n \} _n$ satisfies the large deviations principle in $C_E[0, \infty)$ with rate function
\begin{align*}
\label{eq:rateFK}
	\hat \calJ (\rho) = \int _0 ^\infty  I^{\Psi} (\rho_s) \dd s, \ \rho \in \calL (E), 
\end{align*}
where 
\begin{align*}
   I^{\Psi} (\mu) = -\min\left[\inf _{u \in \calD^{++} (L)}  \int _{E} \frac{Lu}{u} \dd\mu , \int _E \Psi \dd\mu\right]. 
\end{align*}
\end{lemma}
\subsubsection{Proof of Theorem~\ref{thm:proving_LDP:LDP_compact}---compact state-space $E$}
\label{section:proof-of-LDP:compact}
As outlined in the previous section, we can prove Theorem \ref{thm:proving_LDP:LDP_compact} by first verifying the conditions of Lemma \ref{lemma:FengKurtz} under the given assumptions and then apply Proposition \ref{prop:LDP_strategy}.
\begin{proof}[Proof of Theorem \ref{thm:proving_LDP:LDP_compact}]
First, Conditions \ref{cond:FK1}-\ref{cond:FK3} follow from the assumption of Feller continuity \ref{item:thm_LDP_compact:Feller} and tightness; see e.g.\ Remark 11.22 in \cite{FengKurtz2006}.


Next, Conditions \ref{cond:FK4} and \ref{cond:FK5} always hold for compact $E$: take $\varphi_n \equiv 1$, $\Psi \equiv 0$, $Q = \{q\}$ (singleton), and $K^q_n :=E$ for every $n \in \mathbb{N}$.  For this choice both conditions are met---it is only for non-compact spaces $E$ that these conditions become non-trivial (see the proof of Theorem \ref{thm:proving_LDP:LDP_non_compact}). Condition \ref{cond:FK6} is trivially true for compact $E$.

Remains to verify the inequality $H_1 ^{\beta} \leq H_2 ^{\beta}$. Take $d \geq 1$ and $\beta \in \calC ^d$. With the choice $\Psi \equiv 0$ the definitions of $H_i ^{\beta}: \bR ^d \to \bR$, $i=1,2$, become
	\begin{align*}
	H_1 ^{\beta}& =  \inf_{f \in \mathcal{D}^{++}(L)} \sup_{y \in E} \left[ \beta (y) \cdot p + \frac{Lf(y)}{f(y)} \right],\\
	H_2 ^{\beta} &= \sup_{f \in \mathcal{D}^{++}(L)} \inf_{y \in E} \left[ \beta (y) \cdot p + \frac{Lf(y)}{f(y)} \right].
	\end{align*}
We now show that the required inequality follows from Assumption \ref{item:thm_LDP_compact:principal_eigenvalue}, solvability of the principal eigenvalue problem. 

For any $\beta \in \calC ^d$ and $p \in \bR ^d$, define the map $V_p (y): E \to \bR$ as
\begin{align*}
	V_p (y) = \beta (y) \cdot p
\end{align*}
This is a continuous function on $E$ and for every $p$ there exists a function $f_p \in \calD ^+ (L)$ and real eigenvalue $\lambda_p$ such that
\begin{align*}
	(L + \beta \cdot p) f_p = \lambda _p f_p.
\end{align*}
It follows that, for any $p \in \mathbb{R}^d$, we have 
	\begin{align*}
		\lambda _p = \sup_{y \in E}
	\left[
	\frac{Lf_p(y)}{f_p(y)} + \beta(y) \cdot p
	\right] = \inf_{y \in E}
	\left[
	\frac{Lf_p(y)}{f_p(y)} + \beta(y) \cdot p
	\right],
	\end{align*}
	which leads to the upper bound
	\begin{align*}
	H^\beta_1(p) &=
	\inf_{f \in \mathcal{D}^{+}(L)} \sup_{y \in E} \left[
	\frac{Lf(y)}{f(y)} + \beta(y) \cdot p
	\right]	\\ &\leq 
	\sup_{y \in E}
	\left[
	\frac{Lf_p(y)}{f_p(y)} + \beta(y) \cdot p
	\right]	\\ &=
	\inf_{y \in E}
	\left[
	\frac{Lf_p(y)}{f_p(y)} + \beta(y) \cdot p
	\right]	\\	&\leq
	\sup_{f \in \mathcal{D}^{+}(L)} \inf_{y \in E} \left[
	\frac{Lf(y)}{f(y)} + \beta(y) \cdot p
	\right] \\ &= H^\beta_2(p).
	\end{align*}
This shows that Condition \ref{cond:FK7} of Lemma \ref{lemma:FengKurtz} follows from \ref{item:thm_LDP_compact:principal_eigenvalue}. 
As a result, in the setting of compact $E$, Assumptions  \ref{item:thm_LDP_compact:Feller} - \ref{item:thm_LDP_compact:principal_eigenvalue} ensure that Lemma \ref{lemma:FengKurtz} is applicable. This gives the large deviations principle for the empirical measures associated with sped-up versions of the process $Y$ and Proposition \ref{prop:LDP_strategy} transfers this to the empirical measures of the original process. This concludes the proof of the large deviations principle. The form of the rate function is trivially seen to be equal to the prescribed form because of the choice of $\Psi \equiv 0$.
\end{proof}
\subsubsection{Proof of Theorem~\ref{thm:proving_LDP:LDP_non_compact}---non-compact state space $E$}
\label{section:proof-of-LDP:non-compact}
We prove large deviations of the family of measures~$\{\eta_t ^\tau\}_{t > 0}$ introduced at the beginning of Section \ref{sec:aux} by verifying the assumptions of Lemma \ref{lemma:FengKurtz}. Proposition \ref{prop:LDP_strategy} then implies the large deviations principle of the empirical measures $\{\eta_t\}_{t > 0}$ with the prescribed rate function. Whereas the conditions of Lemma \ref{lemma:FengKurtz} where straightforward to verify in the compact setting of Theorem \ref{thm:proving_LDP:LDP_compact}, the non-compact case requires more work. Specifically, because we can no longer assume that there is a solution to the principal eigenvalue problem---such an assumption would not allow us to prove the large deviations principle for the zig-zag process---and the state space is no longer compact, \ref{cond:FK4}-\ref{cond:FK7} are more difficult to verify. A crucial component of the proof of Theorem \ref{thm:proving_LDP:LDP_non_compact} is an inequality that is connected to the necessary comparison principle. To streamline the proof we now state this inequality as a separate result.
\smallskip

For any $V \in C_b(E)$ and $\Psi : E \to \bR$, define $H_1 ^{\Psi}, H_2 ^{\Psi} \in \bR$ by
\begin{equation}
	\begin{split}
	H_1 ^{\Psi} &= \inf_{0 < \kappa \leq 1} \inf_{f \in \mathcal{D}^{++}(L)} \sup_{y \in E} \left[ V(y) + (1-\kappa) \frac{Lf(y)}{f(y)} + \kappa \Psi(y)\right], \label{eq:H1_H2} \\
	H_2 ^{\Psi} &= \sup_{\kappa > 0} \sup_{f \in \mathcal{D}^{++}(L)} \inf_{y \in E} \left[ V(y) + (1+\kappa) \frac{Lf(y)}{f(y)} - \kappa \Psi(y)\right].
\end{split}
\end{equation}
\begin{proposition}\label{prop:non_compact:H1_leq_H2}
	Take any $V \in C_b (E)$ and suppose~\ref{item:thm_LDP_non_compact:mixing} holds and that for any $c\in\mathbb{R}$, the superlevel-set $\{\Psi \geq c\}$ is compact. Then
	\begin{align}
	\label{eq:ineqH1H2}
		H_1 ^{\Psi} \leq H _2 ^{\Psi}.
	\end{align}
\end{proposition}
We first complete the proof of Theorem \ref{thm:proving_LDP:LDP_non_compact}.
\begin{proof}[Proof of Theorem \ref{thm:proving_LDP:LDP_non_compact}]
The proof amounts to showing that Conditions \ref{cond:FK1}-\ref{cond:FK7} of Lemma \ref{lemma:FengKurtz} hold. We start with the ones that are straightforward to obtain from the assumptions of the theorem. 

Conditions \ref{cond:FK1}-\ref{cond:FK3} follow from \ref{item:thm_LDP_compact:Feller} and \ref{item:thm_LDP_compact:tight}. For condition \ref{cond:FK6} the existence of such a compact set follows immediately from the assumption that the initial value $Y(0)$ belongs to a compact set $K \subseteq E$.

We now show that Conditions \ref{cond:FK4} and \ref{cond:FK5} follow from \ref{item:thm_LDP_non_compact:Lyapunov}, the existence of Lyapunov functions $g_1$ and $g_2$ with certain growth properties. We start with \ref{cond:FK4} and define the family of compact sets $K^q_n \subseteq E$ by
\begin{equation*}
K^q_n = \{y \in E\,:\, g_2(y) \leq qn\},\quad q,n \in \mathbb{N}
\end{equation*}
For any $q_1, q_2$ and with $q_3 = \max (q_1, q_2)$, it then holds that
\begin{align*}
	K^{q_1}_n \cup K^{q_2}_n \subseteq K^{q_3}_n, \ \ \forall n \in \bN.
\end{align*}
Because $g_2 (y)$ is finite for any $y\in E$, there exists $q, N \in \bN$ such that $n \geq N$ implies that $y \in K_n ^{q}$. In particular, $\textrm{dist} (y, K_n ^q) = 0$. For the last part of Condition \ref{cond:FK4}, take $q \in \bN$ and $T, a > 0$. It remains to find a $\tilde q$ such that
\begin{align*}
	\limsup _{n \to \infty} \sup _{y \in K  ^{q} _n} \frac{1}{n} \log \bP _y \left( Y_t \notin K ^{\tilde q} _n, \ \textrm{some } t \leq nT \right) \leq -a.
\end{align*}
By Lemma 4.20 in \cite{FengKurtz2006}, for any open neighbourhood $\calO$ of $K_n ^q$,
\begin{align}
\label{eq:openN}
	\mathbb{P}\left(Y_t \notin \mathcal{O}, \,\mathrm{some}\, t \leq nT \,|\, Y_0 \in K^q_n\right) \leq \mathbb{P}\left(Y_0 \in K^q_n\right) e^{-\beta _q + nT\gamma(\mathcal{O})},
\end{align}
where the constants $\beta_q$ and $\gamma(\mathcal{O})$ are given by
\begin{align*}
\beta_q = \inf_{E \setminus \mathcal{O}} g_2 - \sup_{K^q_n} g_2,
\end{align*}
and
\begin{align*}
\gamma(\mathcal{O}) = \max \left(\sup_{\mathcal{O}}e^{-g_2} B e^{g_2}, 0\right).
\end{align*}
By the growth condition for $g_2$ (part (a) of \ref{item:thm_LDP_non_compact:Lyapunov}) for any $\tilde q > \hat q > q$ large enough, there exists an open set $\calO$ such that
\begin{align*}
	K_n ^q \subseteq K_n ^{\hat q} \subseteq \calO \subseteq K^{\tilde q} _n.
\end{align*} 
By definition, $g_2 \leq nq$ on $K ^q _n$ and because $K ^{\hat q} _n \subseteq \calO$, we have $g_2 \geq \hat{q} n$ on $E \setminus \calO$. Combined with the upper bound $\gamma (\calO) \leq \gamma (E)$ this gives, starting from \eqref{eq:openN},
\begin{align*}
	\frac{1}{n} \log \bP \left( Y_t \notin K ^{\tilde q} _n, \ \textrm{some } t\leq nT | Y_0 \in K_n ^q \right) &\leq T \gamma(\calO) - \frac{1}{n} \beta _q \\
	&\leq T \gamma(E) +q - \hat{q}.
\end{align*}
The last part of Condition \ref{cond:FK4} is now straightforward to obtain. First, take $\hat q = \hat q (q, a, T)$ large enough that the right-hand side of the last display is bounded by $-a$:
\begin{align*}
	T \gamma(E) +q - \hat{q} \leq -a.
\end{align*}
Next, choose $\tilde q = \tilde q (q, a, T)$ large enough that there is an open set $\calO$ such that $K ^{\hat q} _n \subseteq \calO \subseteq K^{\tilde q} _n$. The asymptotic statement then follows, which concludes the verification of condition \ref{cond:FK4} of Lemma \ref{lemma:FengKurtz}.
\smallskip

To show that condition \ref{cond:FK5} is fulfilled we generalize the arguments used in Example 11.24 in \cite{FengKurtz2006}. The functions $\varphi_n$ are constructed from the Lyapunov functions $g_1$ and $g_2$. First, define
\begin{equation}\label{eq:proving_LDP:r_n}
r_n := \sup \left\{g_1(y) \,: \, y \in E,\, g_1(y) g_2(y) \leq n^2\right\}.
\end{equation}
Then $r_n \to \infty$ and $r_n/n \to 0$ as $n \to \infty$, since~\ref{item:thm_LDP_non_compact:Lyapunov2} and the condition in the set imply
\begin{equation*}
\frac{r_n}{n} = \frac{g_1(y_n)}{n} \leq \sqrt{\frac{g_1(y_n)}{g_2(y_n)}} \to 0.
\end{equation*}
Furthermore, for each $q$ there exists $n_q$ such that $n \geq n_q$ implies
\begin{equation*}
K^q_n \subseteq \left\{y\,:\, g_1(y) \leq r_n\right\}.
\end{equation*}
For a smooth, non-decreasing and concave function $\rho : [0,\infty) \to [0,2]$ satisfying $\rho(r) = r$ for $0 \leq r \leq 1$ and $\rho(r) = 2$ for $r \geq 3$, define the functions $\varphi_n$ by cutting off $g_1$:
\begin{equation}\label{eq:proving_LDP:varphi_n}
\varphi_n(y) := e^{r_n} \rho(e^{-r_n} e^{g_1(y)}).
\end{equation}
We have $\varphi_n = e^{g_1}$ on the compact sets $K^q_n$. Setting $\Psi = e^{-g_1} B e^{g_1}$, we therefore obtain
\begin{equation*}
\frac{L \varphi_n(y)}{\varphi_n(y)} = \Psi(y),\quad y \in K^q_n,\, n \geq n_q.
\end{equation*}
The fact that $r_n/n\to 0$ as $n\to\infty$ implies $n^{-1}\log\|\varphi_n\| \to 0$ as $n\to\infty$. For proving that
\begin{equation*}
    \sup_{n,y}\frac{L\varphi_n(y)}{\varphi_n(y)}<\infty,
\end{equation*}
it is sufficient to show that for any positive function $u:E \to (0,\infty)$ in the domain of $B$ and for any $y_0 \in E$, we have
\begin{equation}\label{eq:proof-LDP-noncompact:ineq-pos-max-pr}
\frac{B (\rho(u))(y_0)}{\rho(u)(y_0)} \leq \frac{\max\left(B u(y_0),0\right)}{u(y_0)}.
\end{equation}
Then with $u=e^{-r_n}e^{g_1}$ and noting that $L\varphi_n=B\varphi_n$, by linearity we obtain
\begin{equation*}
    \frac{L\varphi_n}{\varphi_n} \leq \frac{\max\left(B e^{g_1},0\right)}{e^{g_1}},
\end{equation*}
and the result follows since $\Psi(y)=e^{-g_1(y)} \left(B e^{g_1}\right)(y)\to-\infty$ as $|y|\to\infty$. Hence we are left with verifying~\eqref{eq:proof-LDP-noncompact:ineq-pos-max-pr}.

If $u(y_0) \in (3,\infty)$, then $\rho(u)(y_0)$ is maximal. Hence by the positive maximum principle, $B\rho(u)(y_0) \leq 0$, and the inequality follows. If $u(y_0) \in (0,3]$, then $a_0 := \rho'(u(y_0)) \in [0,1]$, the region where $\rho$ goes from slope one to slope zero. Consider the function $f_0 := \rho(u) - a_0 u$. Since $g_0(s) := \rho(s) - a_0 s$ is maximal for $s_0$ satisfying $\rho'(s_0) = a_0$, we obtain that $y_0$ is an optimizer, that is $f_0(y_0) = \sup_y f(y)$. Furthermore, $g_0(s_0) \geq 0$, so by the positive maximum principle $Bf_0(y_0) \leq 0$. By linearity of $B$ and since $r \leq \max(r,0)$, we obtain the inequality $B \rho(u)(y_0) \leq a_0 \cdot \max\left(Bu(y_0),0\right)$. Hence
\begin{equation*}
\frac{B \rho(u)(y_0)}{\rho(u)(y_0)} \leq \frac{a_0}{\rho(u)(y_0)}\max(Bu(y_0),0) \leq \frac{1}{u(y_0)} \max(Bu(y_0),0),
\end{equation*}
using that $0 \leq g_0(s_0) = \rho(u)(y_0)-a_0 u(y_0)$. This finishes the verification of~\eqref{eq:proof-LDP-noncompact:ineq-pos-max-pr}.
\smallskip

It remains to show that condition \ref{cond:FK7} is fulfilled. However, this is precisely the conclusion of Proposition \ref{prop:non_compact:H1_leq_H2} - the function $V(y) = \beta (y)\cdot p$ where $\beta$ is as in condition (vii) is an element of $C_b (E)$, and by~\ref{item:thm_LDP_non_compact:Lyapunov3}, the function $\Psi$ has compact superlevel-sets.
\smallskip

We have shown that under the assumptions of the theorem, all conditions of Lemma \ref{lemma:FengKurtz} are fulfilled. The large deviations principle for the empirical measures of the sped-up versions thus holds and Proposition \ref{prop:LDP_strategy} then gives the large deviations principle for the empirical measures $\{ \eta _t \}$ associated with~$Y$. 
\smallskip

We are left with showing that the rate function $I ^{\Psi}$ of Proposition \ref{prop:LDP_strategy} satisfies
\begin{align*}
	I ^{\Psi} (\mu) = -\inf _{u \in \calD^{++} (L)}  \int _{E} \frac{Lu}{u} \dd\mu.
\end{align*}
Below, we prove that
\begin{equation*}
    \limsup_{n\to\infty}\int_{E}\frac{L\varphi_n}{\varphi_n} \dd\mu \leq \int_E \Psi \dd\mu.
\end{equation*}
Then 
\begin{equation*}
    \inf_{u \in \calD^{++} (L)}\int_E\frac{Lu}{u} \dd\mu \leq \inf_n \int_E \frac{L\varphi_n}{\varphi_n} \dd\mu \leq \int_E \Psi \dd\mu,
\end{equation*}
and hence the rate function is given by
\begin{align*}
	I^{\Psi} (\mu) = -\min \left[ \inf _{u \in \calD^{++} (L)}  \int _{E} \frac{Lu}{u} \dd\mu,  \int _E \Psi d\mu \right] = -\inf _{u \in \calD^{++} (L)}  \int _{E} \frac{Lu}{u} \dd\mu.
\end{align*}
To see that the functions $\varphi_n$ satisfy the limsup inequality, note that $\Psi$ has compact super-level sets and $\Psi(y)\to-\infty$ as $|y|\to\infty$. Since the compact sets $K_n^q$ exhaust $E$ in the sense that $E = \cup_n K^q_n$ and $K^q_n\subseteq K^q_{n+1}$, there exists a constant $C > 0$ such that
\begin{equation*}
    f_n = -\frac{L\varphi_n}{\varphi_n} + C \geq 0 .
\end{equation*}
Pointwise, we have $f=-\Psi + C = \liminf_n f_n$. Therefore, by Fatou's lemma
\begin{equation*}
    \liminf_{n\to\infty} \int_E \left[-\frac{L\varphi_n}{\varphi_n} + C\right]\,\dd\mu \geq \int_E\left[-\Psi + C\right]\,\dd\mu,
\end{equation*}
and the required limsup inequality follows from reorganizing.
\end{proof}

We now prove the important Proposition \ref{prop:non_compact:H1_leq_H2}. The proof is essentially a combination of different arguments from Chapter 11 and Appendix B of \cite{FengKurtz2006} (see especially Lemmas 11.12, 11.37, B.9-B.11 for full details). We present the proof as to make the presentation self-contained and give a succinct derivation of the results for the setting we consider. The main novelty compared to the arguments in \cite{FengKurtz2006} is that we work with measures $\nu \in \calP _c (E) $ rather than imposing the condition $\int _E \Psi \nu > -\infty$, and we must verify that we can indeed modify the latter.

\begin{proof}[Proof of Proposition \ref{prop:non_compact:H1_leq_H2}]
The strategy is to find two constants, depending on $V$, $c_V^\ast$ and $c_V ^{\ast\ast}$ such that $c_V ^\ast \geq c_V ^{\ast \ast} $ and
\begin{align}
\label{eq:ineqH}
	H_1 ^{\Psi} \leq c_V ^{\ast\ast}, \, \, \textrm{and }\,  H_2 ^{\psi} \geq c_V ^\ast,
\end{align} 
To achieve this we study the following quantity: for $\nu \in \calP_c (E)$, define
\begin{align*}
	c_V (\nu) = \limsup _{t \to \infty} \frac{1}{t} \log \bE \left[ \mathrm{exp}\left\{ \int _0 ^t V(Y(s))\dd s \right\}  \right].
\end{align*}
It can be shown - see e.g. Lemma B.9 in \cite{FengKurtz2006} - that under \ref{item:thm_LDP_non_compact:mixing}, $c_V (\nu)$ exists for each $\nu \in \calP _c (E)$ and the necessary inequalities for $H_i ^{\Psi}$ can be derived for
\begin{align*}
	c_V^{\ast} = \inf_{\nu \in \calP _c (E)} c_V(\nu) \quad \text{ and } \quad
	c_V^{\ast\ast} = \sup_{\nu \in \calP _c (E)} c_V(\nu).
\end{align*}
Cleary $c_V ^\ast \leq c_V ^{\ast \ast}$. However it can be shown, again using~\ref{item:thm_LDP_non_compact:mixing}, that the two quantities are in fact equal, that $c_V (\nu)$ is independent of $\nu$ on $\calP _c (E)$. If we can prove \eqref{eq:ineqH} this would then yield the claim.
We start with the upper bound
\begin{align*}
	H_1 ^{\Psi} \leq c_V ^{\ast\ast}.
\end{align*}
An argument similar to what will follow is also used in \cite{DonskerVaradhan75}, in the proof of their Lemma~2.

Because $\Psi$ has compact superlevel-sets $\{ \Psi \geq c \}$, $c\in \bR$, and $\Psi (y) \to - \infty$ as $\norm{y} \to \infty$, it can be shown using the arguments of Lemma B.11 of \cite{FengKurtz2006} that
\begin{align*}
	\sup_{\nu \in \mathcal{P}_c(E)} \inf_{f \in \mathcal{D}^{++}(L)} \int_E \left(V + \frac{Lf}{f}\right)d\nu \leq c_V^{\ast\ast}.
\end{align*}
It therefore suffices to show that 
\begin{align}
\label{eq:ineqH1}
H_1 ^{\Psi} \leq \sup_{\nu \in \mathcal{P}_c(E)} \inf_{f \in \mathcal{D}^{++}(L)} \int_E \left(V + \frac{Lf}{f}\right)d\nu.
\end{align}
For any finite collection of functions $f_1, \dots , f_m$ in $\calD ^{++} (L)$ and scalars $\alpha_i \geq 0$, $i=1,\dots , m$, $\sum \alpha_i =1$, we have 
\begin{align*}
	H_1 ^{\Psi} \leq \inf_{0 < \kappa \leq 1} \inf_{\alpha_i} \inf_{f_1,\dots,f_m} \sup_{y \in E} \left[
	V(y) + (1-\kappa)\sum_{i=1}^m \alpha_i \frac{Lf_i(y)}{f_i(y)} + \kappa \Psi(y),
	\right]
\end{align*}
which follows as in Lemma~11.35 of~\cite{FengKurtz2006}; define for $t>0$
\begin{align*}
	h_t = \frac{1}{t} \int _0 ^t S(\tau) \prod _{i=1} ^m f_i ^{\alpha _i} d \tau. 
\end{align*}
Then as we let $t\to 0$,
\begin{align*}
\lim_{t \to 0} h_t = \prod_i f_i^{\alpha_i},
\end{align*}
and we have the upper bound
\begin{align*}
 \lim_{t \to 0}Lh_t \leq \left(\prod_i f_i^{\alpha_i}\right) \sum_i \alpha_i \frac{Lf_i}{f_i},
\end{align*}
where the convergence is uniform. Now specializing in the definition of~$H_1^\Psi$ to with these type of functions~$h_t$ with $f_1,\dots,f_m$ and according~$\alpha_i$, and taking the limit $t\to 0$ gives the above estimate for~$H_1^\Psi$.

By Lemma~11.37 of~\cite{FengKurtz2006}, we can select a sequence of functions $\{ f_i \} $ from $\calD ^{++} (L)$ such that for any $\mu \in \calP (E)$,
	\begin{equation*}
	\inf_{i} \int_E \frac{Lf_i}{f_i} d\mu = \inf_{f \in \mathcal{D}^{++}(L)} \int_E \frac{Lf}{f} d\mu.
	\end{equation*}
	Specialising to these functions, for any $m \in \bN$ and $\kappa > 0$ we have the upper bound
	\begin{equation*}
	H_1 ^{\Psi} \leq \inf_{\alpha_i} \sup_{y \in E} \left[
	V(y) + (1-\kappa)\sum_{i=1}^m \alpha_i \frac{Lf_i(y)}{f_i(y)} + \kappa \Psi(y).
	\right]
	\end{equation*}
	The functions $V + Lf_i/f_i$ are bounded, but a priori there is no guarantee that the supremum is attained in a given compact set. However, because $\Psi(y) \to -\infty$ as $\norm{y} \to \infty$, for any $m \in \bN$ and $\kappa > 0$, there exists a constant $\ell = \ell(m,\kappa) > 0$ such that the supremum is attained in the compact set $K_{\ell} = \{\Psi \geq -\ell\}$. Therefore, if we define $\mathcal{K}_\ell = \{\nu \in \mathcal{P}(E): \nu(K_\ell) = 1\}$, then
	\begin{align*}
	H_1 ^{\Psi} &\leq \inf_{\alpha_i} \sup_{y \in K_{\ell}} \left[
	V(y) + (1-\kappa)\sum_{i=1}^m \alpha_i \frac{Lf_i(y)}{f_i(y)} + \kappa \Psi(y)
	\right]	\\
	&= \inf_{\alpha_i} \sup_{\nu \in \mathcal{K}_\ell} \left[
	\int_E\left(V(y) + (1-\kappa)\sum_{i=1}^m \alpha_i \frac{Lf_i(y)}{f_i(y)} + \kappa \Psi(y) \right) d\nu(y)
	\right].
	\end{align*}
	For any $\ell$, we have $\mathcal{K}_\ell \subseteq \mathcal{P}_c(E)$, so that
	\begin{equation*}
	H_1 ^{\Psi} \leq \inf_{\alpha_i} \sup_{\nu \in \mathcal{P}_c(E)} \left[
	\int_E\left(V + (1-\kappa)\sum_{i=1}^m \alpha_i \frac{Lf_i}{f_i} + \kappa \Psi \right) d\nu
	\right].
	\end{equation*}
	For any $m \in \bN$, the set $\{ \alpha_i: \alpha _i \geq 0, \sum _{i=1} ^m \alpha _i =1 \}$ is compact and the infimum and supremum in the last display can be exchanged by Sion's Theorem. This yields
	\begin{align*}
	H_1 ^{\Psi} &\leq  \sup_{\nu \in \mathcal{P}_c(E)} \inf_{\alpha_i} \left[
	\int_E\left(V + (1-\kappa)\sum_{i=1}^m \alpha_i \frac{Lf_i}{f_i} + \kappa \Psi \right) d\nu
	\right]	\\
	&= \sup_{\nu \in \mathcal{P}_c(E)} \left[
	\int_E\left(V + (1-\kappa)\min_{i\leq m} \frac{Lf_i}{f_i} + \kappa \Psi \right) d\nu
	\right],
	\end{align*}
	where we have used that $\inf_{\alpha_i}\sum \alpha_i x_i = \min_i x_i$ for non-negative $x$ and $f_i \in \calD ^{++} (L)$. Taking the infimum over $\kappa$ and the limit $m \to \infty$,
	\begin{align*}
	H_1 ^{\Psi} &\leq \lim_{m \to \infty} \inf_{0 < \kappa \leq 1} \sup_{\nu \in \mathcal{P}_c(E)} \left[
	\int_E V d\nu + (1-\kappa) \min_{i \leq m} \int_E \frac{Lf_i}{f_i} d\nu + \kappa \int_E \Psi d\nu
	\right]	\\
	& = \lim_{m \to \infty} \sup_{\nu \in \mathcal{P}_c(E)} \left[
	\int_E V d\nu + \min\left\{\min_{i \leq m} \int_E \frac{Lf_i}{f_i} d\nu, \int_E \Psi d\nu \right\}
	\right].
	\end{align*}
	The limit and supremum can be shown to commute similarly to the last part of the proof of Lemma 11.12 in \cite{FengKurtz2006}, leading to
	
	\begin{align*}
	H_1 ^{\Psi} &\leq  \sup_{\nu \in \mathcal{P}_c(E)} \left[
	\int_E V d\nu + \min\left\{\inf_{f \in \mathcal{D}^{++}(L)} \int_E \frac{Lf}{f} d\nu, \int_E \Psi d\nu\right\}
	\right] \\
	&\leq \sup_{\nu \in \mathcal{P}_c(E)} \left[
	\int_E V d\nu + \inf_{f \in \mathcal{D}^{++}(L)} \int_E \frac{Lf}{f} d\nu \right] .
	\end{align*}
%
	This completes the proof of the upper bound for $H^ {\Psi} _1$.
	
	Next, we move to the lower bound for $H_2 ^{\Psi}$. Take $\lambda < c_V^\ast$. We prove that for any $\varepsilon > 0$, we have $H_2 ^{\Psi} \geq \lambda - \varepsilon$. To this end, 
	we define the new semigroup $ \{ T(t) \} $ by
	\begin{align*}
	(T(t)f)(y) = \bE \left[ f(Y_t) e^{\int_0 ^t V(Y_s)ds} | Y(0)=y \right],
	\end{align*}
	set
	\begin{align*}
	    R_{\lambda} ^t f = \int _0 ^t e^{-\lambda s} T(s) g ds,
	\end{align*}
	and take $\Gamma$ to be the collection of functions $f_{\gamma}$ of the form
	\begin{align*}
	    f_{\gamma} = \int _0 ^\infty R_{\lambda} ^t 1 \gamma (dt), \ \ \gamma \in \calP ([0,\infty)).
	\end{align*}
	Then $\Gamma \subseteq \mathcal{D}^{++}(L)$ and for any $f \in \Gamma$ we have the uniform lower bound
	\begin{equation}
	\label{eq:uniformLower}
	V(y) + (1+\kappa) \frac{Lf (y)}{f(y)}\geq
	-(1-2\kappa)\|V\|, \ \ y \in E.
	\end{equation}
    Because $\Gamma \subseteq \calD ^{++} (L)$, for any $\kappa >0$ we have the lower bound
	\begin{equation*}
	H_2 ^{\Psi} \geq \sup_{f \in \Gamma} \inf_{y\in E} \left[ V(y) + (1+\kappa)\frac{Lf(y)}{f(y)} - \kappa \Psi (y) \right].
	\end{equation*}
	Due to the uniform lower bound \eqref{eq:uniformLower} and the fact that $\Psi(y) \to -\infty$ as $\norm{y} \to \infty$, for any $\kappa$ there exists an $\ell = \ell(\kappa)$ such that the infimum over $E$ is attained in the compact set $K_\ell = \{ \Psi \geq -\ell\}$. Therefore,

	\begin{align*}
	H_2 ^{\Psi} &\geq \sup_{f \in \Gamma} \inf_{y\in K_\ell} \left[ V(y) + (1+\kappa)\frac{Lf(y)}{f(y)} - \kappa \Psi(y) \right] \\
	&= \sup_{f \in \Gamma} \inf_{\nu \in \mathcal{K}_\ell} \left[ \int_E \left( V + (1+\kappa)\frac{Lf}{f} - \kappa \Psi\right)\,d\nu \right]	\\
	&= \sup_{f \in \Gamma} \inf_{\nu \in \mathcal{K}_\ell} \frac{1}{\int_E f d\nu} \left[ -\kappa \int_E (V + \Psi) f d\nu  + \int_E (1+\kappa) (V+L)f d\nu \right],
	\end{align*}
	where $\mathcal{K}_\ell = \{\nu \in \mathcal{P}(E): \nu(K_\ell) = 1\}$. The second equality follows from the fact that $\inf_\nu \int_E (a/b) d\nu = \inf_\nu (\int_E a d\nu)/ (\int_E b d\nu)$ for $b > 0$. By compactness of $\mathcal{K}_\ell$ and the fact that both $\mathcal{K}_\ell$ and $\Gamma$ are convex, the infimum and supremum are exchangable by Sion's Theorem. This gives the lower bound
	\begin{align*}
	H_2 ^{\Psi}	&\geq \inf_{\nu \in \mathcal{K}_\ell} \sup_{f \in \Gamma} \frac{1}{\int_E f d\nu} \left[ -\kappa \int_E (V + \Psi) f d\nu  + \int_E (1+\kappa) (V+L)f d\nu \right]	\\
	&\geq \inf_{\nu \in \mathcal{P}_c(E)} \sup_{f \in \Gamma}  \frac{1}{\int_E f d\nu} \left[ -\kappa \int_E (V + \Psi) f d\nu  + \int_E (1+\kappa) (V+L)f d\nu \right],
	\end{align*}
	The second estimate follows since $\mathcal{K}_\ell \subseteq \mathcal{P}_c(E)$ for any $\ell$. The rest of the proof follows arguments similar to those used in~\cite{FengKurtz2006}: taking the limit $\kappa \to 0$, and moving it inside the infimum and supremum, we obtain the lower bound
	\begin{align*}
	H_2 ^{\Psi} & \geq \inf_{\nu \in \mathcal{P}_c(E)} \sup_{f \in \Gamma} \left[ \frac{1}{\int_E f d\nu}\int_E (V+L)f d\nu \right].
	\end{align*}
	Therefore, for any $\varepsilon$, there exists a $\nu_\varepsilon \in \mathcal{P}_c(E)$ such that
	\begin{align*}
	H_2 ^{\Psi} \geq \sup_{f \in \Gamma} \left[ \frac{1}{\int_E f d\nu_\varepsilon}\int_E (V+L)f d\nu_\varepsilon \right] - \varepsilon.
	\end{align*}
	%
	There exist functions $f_t \in \Gamma$ satisfying
	\begin{align*}
	\frac{\int_E (V+L) f_t d\nu}{\int_E f_t d\nu} = \lambda + \frac{ \int_E e^{-\lambda t} T(t) 1 d\nu - 1}{\int_E f_t d\nu},
	\end{align*}
	for any $\nu \in \mathcal{P}_c(E)$. Specialising to such $f_t$, we obtain
	\begin{align*}
	H_2 ^{\Psi} \geq \lambda + \frac{\int_E e^{-\lambda t} T(t)1 d\nu_\varepsilon - 1}{\int_E f_t d\nu_\varepsilon} - \varepsilon.
	\end{align*}
	Since $\limsup_{t \to \infty} \int_E e^{-\lambda t} T(t) 1 d\nu_\varepsilon = \infty$, the second term is positive for $t$ large enough, giving the bound
	\begin{equation*}
	H_2 ^{\Psi} \geq \lambda - \varepsilon.
	\end{equation*}
	This completes the proof of the lower bound for $H_2 ^{\Psi}$, and thereby the lemma.
\end{proof}

\subsection{Proofs for the empirical measure of the zig-zag process}
In this section, we prove the large deviations theorems for the empirical measures of the zig-zag process.
\subsubsection{Proof of Theorem~\ref{thm:LDP-zigzag:compact}---compact case}
\label{section:zig-zag:LDP:compact}
For the proof of Theorem~\ref{thm:LDP-zigzag:compact}, recall that the zig-zag generator takes the form
\begin{equation*}
    Lf(x,v)=v\partial_x f(x,v) + \lambda(x,v)\left[f(x,-v)-f(x,v)\right], \quad (x,v)\in E:= \mathbb{T}\times\{\pm 1\}.
\end{equation*}
It is enough to show that assumptions \ref{item:thm_LDP_compact:Feller}-\ref{item:thm_LDP_compact:principal_eigenvalue} hold for the zig-zag process on $E$, the result then follows form Theorem \ref{thm:proving_LDP:LDP_compact}.


We first verify that $L$ is a closed operator that generates the zig-zag process. Note that $L$ is a restriction of the extended generator (see Section~\ref{sec:zig-zag}).
We verify that $L$ is a closed operator. Let $\{ f_n \}$ be a sequence in $\mathcal{D}(L)$ such that $f_n \to f$ and $Lf_n \to g$, some $f,g$, both uniformly on~$E$. Then
\begin{equation*}
\lim_{n \to \infty} v \partial_x f_n(x,v) = g(x,v) - \lambda(x,v) \left[ f(x,-v) - f(x,v) \right].
\end{equation*}
We can represent $f_n (x,v)$ as 
\begin{equation*} f_n(x,v) = f_n(0,v) + v \int_0^x v \partial_x f_n(\xi,v) \, d\xi,
\end{equation*}
and from the dominated convergence theorem we obtain that
\begin{equation*}
f(x,v) = f(0,v) + v \int_0^x \left[ g(\xi,v) - \lambda(\xi,v) \left( f(\xi,-v) - f(\xi,v) \right) \right] \, d\xi.
\end{equation*}
In particular, $f \in \mathcal{D}(L)$, and $Lf = g$ follows from taking derivative $\partial_x$ and multiplying by $v$. 
	
The Feller-continuity property \ref{item:thm_LDP_compact:Feller} of the zig-zag semigroup $S$ is proven in Proposition~4 of~\cite{BierkensRoberts2017}. Since $\mathbb{T}$ is compact, this also follows from the boundedness of the continuous rates $\lambda$, see \cite[Theorem 27.6]{Davis1993}.

It remains to verify assumption \ref{item:thm_LDP_compact:principal_eigenvalue}, the principal-eigenvalue problem. Take $V \in C(E)$. We will show that for any constant $\gamma > \sup_E V$, as a map from $C(E)$ to $\mathcal{D}(L) \subseteq C(E)$, the resolvent
	\begin{align}
	\label{eq:resolvent}
	R_\gamma = \left(
	\gamma - (V+L)
	\right)^{-1},
	\end{align}
	is compact and strongly positive; here strongly positive means that if $f \geq 0$ and $f \neq 0$, then $R_\gamma f > 0$ on $E$. Given strong positivity and compactness, by the Krein-Rutman theorem there exists a strictly positive function $g \in C(E)$ and a real eigenvalue $\beta >0 $ such that
	\begin{equation*}
	\left(
	\gamma - (V+L)
	\right)^{-1} g = \beta g.
	\end{equation*}
 The resolvent maps into the domain of $L$, so that $g \in \mathcal{D}(L)$. An application of $\gamma - (V+L)$ in the eigenvalue equation gives
	\begin{equation*}
	(V+L) g = \left(\gamma - \frac{1}{\beta}\right) g.
	\end{equation*}
	This is precisely \ref{item:thm_LDP_compact:principal_eigenvalue} with function $g$ and eigenvalue $(\gamma - 1/\beta)$. 
	
	We are left with verifying that the resolvents defined by \eqref{eq:resolvent} are strongly positive and compact. For strong positivity, because $V$ is continuous on $E$, it is sufficient to prove strong positivity of $(\gamma - L)^{-1}$; see  \cite[Proposition C-III-3.3]{arendt1986one}. The resolvent $(\gamma - L)^{-1}$ exists for any $\gamma > 0$, and is given by
	\begin{equation}
	\label{eq:resolvent2}
	(\gamma - L)^{-1} f = \int_0^\infty e^{-\gamma t} S(t)f \, dt.
	\end{equation}
	The semigroup associated to the zig-zag process is irreducible in the following sense: for any $f \in C(E)$ such that $f \geq 0$ and $ f \neq 0$, 
	\begin{align*}
	\cup_{t \geq 0}\left\{
	z \in E \,:\, S(t) f(z) > 0
	\right\} = E.
	\end{align*}
	Combined with \eqref{eq:resolvent2} this implies strong positivity of $(\gamma - L) ^{-1}$; see \cite[Definition C-III-3.1]{arendt1986one}. 
	\smallskip
	
	For compactness of $R_\gamma$, let $A \subseteq C(E)$ be bounded. We show that the image $B := R_\gamma(A) \subseteq C(E)$ is bounded and equi-continuous. Compactness of the resolvent then follows from an application of the Arzelà-Ascoli theorem. To show boundedness, by dissipativity of $L$ we obtain, for any $g \in B$,
	\begin{align*}
	(\gamma -\|V\|_E) \|g \| &\leq
	\|(\gamma - (V+L))g\| \\
	& \leq \sup_{f \in A} \|f\| \\
	&< \infty.
	\end{align*}
	Hence $B$ is bounded by $C_A / (\gamma - \|V\|_E)$, where $ C_A := \sup _{f \in A} \|f\|$, and we end the proof by showing that $B$ is equi-continuous. For any $g \in B$ we have $(\gamma - (V + L))g = f$ for some $f \in A$, which implies
	\begin{equation*}
	v \partial_x g(x,v) = f(x,v) + V(x,v) g(x,v) + \gamma g(x,v) - \lambda(x,v) (g(x,-v) - g(x,v).
	\end{equation*}
	By boundedness of the functions $\lambda(x,v)$ and $V$ on $E$ and the sets $A$ and $B$,
	\begin{equation*}
	\sup_{g \in B} \|\partial_x g\| \leq C \sup_{g \in B}\|g\| + \sup_{f \in A}\|f\| < \infty.
	\end{equation*}
	Hence functions in $B$ have uniformly bounded derivatives, and as a consequence, $B$ is equi-continuous. It follows that $R_{\gamma}$ in \eqref{eq:resolvent} is compact and strongly continuous. This finishes the verification of \ref{item:thm_LDP_compact:principal_eigenvalue} and we have shown that assumptions \ref{item:thm_LDP_compact:Feller}-\ref{item:thm_LDP_compact:principal_eigenvalue} hold for the zig-zag process on the compact state space $\bT \times \{ \pm 1\}$. An application of Theorem \ref{thm:proving_LDP:LDP_compact} then proves the claimed large deviations principle.
\subsubsection{Proof of Theorem~\ref{thm:LDP-zigzag:non-compact-1d}---non-compact case}
\label{section:zig-zag:LDP:non_compact}
For notational simplicity we take $E = \bR \times \{ \pm 1\}$.
Similar to the proof of Theorem \ref{thm:LDP-zigzag:compact}, the strategy is to verify the conditions of the more general large deviations result Theorem \ref{thm:proving_LDP:LDP_non_compact}, which covers the non-compact setting. That is, it suffices to verify \ref{item:thm_LDP_compact:Feller}, \ref{item:thm_LDP_compact:tight}, \ref{item:thm_LDP_non_compact:Lyapunov} and \ref{item:thm_LDP_non_compact:mixing}.
\smallskip

Condition \ref{item:thm_LDP_compact:Feller}, Feller-continuity of the Markov semigroup, is proven in Proposition 4 of~\cite{BierkensRoberts2017}. 

Next, we use Theorem 7.2 of \cite{EthierKurtz1986} to verify \ref{item:thm_LDP_compact:tight}. 
Define the metric $d$ on $E$ as 
	\begin{align*}
	d((x,v),(y,v')) = |x-y|_\mathbb{R} + |v-v'|,
	\end{align*} 
	and for any path $\gamma \in D_E[0,\infty)$ set
	\begin{equation*}
	w'(\gamma,\delta,T) = \inf_{\{t_i\}} \max_i \sup_{s,t \in [t_i,t_{i+1})} d(\gamma(s),\gamma(t)),
	\end{equation*}
	where the infimum is taken over finite partitions $\{ t_i \}$ of $[0,T]$ such that $\min_i |t_{i+1}-t_i| > \delta$. Theorem 7.2 of \cite{EthierKurtz1986} states that tightness of $\{ \bP _y : \ y \in K \}$ is equivalent to the following two conditions: 
\begin{enumerate}[label =(\arabic*)]
	\item \label{item:proof-LDP-zigzag:tightness1} For any $\varepsilon > 0$ and rational $t > 0$, there exists a compact set $K_{\varepsilon,t} \subseteq E$ such that
	\begin{equation*}
	\inf_{y \in K} \mathbb{P}_y\left[Y_t \in K_{\varepsilon,t} \right] \geq 1 - \varepsilon.
	\end{equation*}
	\item \label{item:proof-LDP-zigzag:tightness2} 
	For any $\varepsilon > 0$ and $T > 0$, there exists a $\delta > 0$ such that
	\begin{equation*}
	\sup_{y \in K} \mathbb{P}_y\left[w'\left(Y,\delta,T\right) \geq \varepsilon\right] \leq \varepsilon.
	\end{equation*}
\end{enumerate}
The spatial component $X_t$ of the zig-zag process propagates with finite speed. This implies that there exists a compact set $K_t \subseteq \mathbb{R}$ such that if $y \in K$, then
\begin{equation*}
\mathbb{P}_y \left[X_t \in K_t \right] = 1.
\end{equation*}
For any $\varepsilon > 0$ and $t>0$, taking $K_{\varepsilon,t} = K_t \times \{\pm 1\}$ gives~\ref{item:proof-LDP-zigzag:tightness1}. 
\smallskip

For part \ref{item:proof-LDP-zigzag:tightness2}, let $\varepsilon > 0$ and $T > 0$. For any realization $Y(\omega)$ of the zig-zag process on the time interval $[0,T]$, if the sojourn times $\tau_i$ satisfy $\min_i \tau_i > 2\delta$, then $w'(Y(\omega),\delta,T) \leq 2\delta$. In particular, for $\delta$ small enough, $w'(Y(\omega),\delta,T) < \varepsilon$. The probability of having at least one sojourn time that is less than $2\delta$ can be estimated uniformly over starting points $y\in K$. Let $K(T)$ denote the set of points that the zig-zag can reach in the time interval $[0,T]$ when starting in the set $K$ and set $\lambda_K = \sup_{y \in K(T)}\lambda(y)$, a uniform upper bound on the jump rates $\lambda(x,v)$. An estimate for the probability of at least one sojourn time that is less than $2\delta$ is then given by
\begin{equation*}
\sup_{y\in K} \mathbb{P}_y\left[\min_i \tau_i \leq 2\delta\right] \leq 1 - e^{-\lambda_K 2\delta}.
\end{equation*}
For any $y\in K$ we obtain the bound
\begin{align*}
\mathbb{P}_y\left[w'(Y,\delta,T) \geq \varepsilon\right] &= \mathbb{P}_y\left[\{w'(Y,\delta,T)\geq \varepsilon\} \cap \{\min_i \tau_i > 2\delta\}\right] \\
& \quad + \mathbb{P}_y\left[\{w'(Y,\delta,T)\geq \varepsilon\} \cap \{\min_i \tau_i \leq 2\delta\}\right]\\
&\leq 0 + \mathbb{P}_y\left[\{\min_i \tau_i \leq 2\delta\}\right] \\
&\leq 1 - e^{-\lambda_K 2\delta}.
\end{align*}
It follows that, as $\delta \to 0$,
\begin{equation*}
\sup_{y \in K} \mathbb{P}_y\left[\{\min_i \tau_i \leq 2\delta\}\right]  \leq 1 - e^{-\lambda_K 2\delta} \to 0,
\end{equation*}
and \ref{item:proof-LDP-zigzag:tightness2} follows from taking $\delta$ small enough that $1 - e^{-\lambda_K 2\delta} < \epsilon$.

We now move to verifying Condition \ref{item:thm_LDP_non_compact:Lyapunov}, by explicitly defining two Lyapunov functions $g_1,g_2 : E \to \mathbb{R}$ satisfying the condition. For brevity, we carry out the calculations for the case of $\gamma (x) \equiv 0$ in the switching rate $\lambda$ (see \eqref{eq:switching-intensity-condition-2}). Then we can use the following functions: for $\alpha_1,\alpha_2 \in (0,1)$ and $\beta > 0$, let
\begin{align*}
g_1(x,v) &= \alpha_1 V(x) + \beta v U'(x),\\
g_2(x,v) &= \alpha_2 U(x) + \beta v U'(x). 
\end{align*}
For non-constant $\gamma$ that is uniformly bounded by some $\bar \gamma$, the following functions can instead be used:
\begin{align*}
g_1(x,v) &= \alpha_1 V(x) + \phi(v U'(x)),\\
g_2(x,v) &= \alpha_2 U(x) + \phi(v U'(x)),
\end{align*}
where $\phi(s) = \beta \frac{1}{2} \text{sign}(s) \log(\bar{\gamma} + |s|)$ and $\beta \in (0,1)$. For example, for the choice $\beta = 1/2$ calculations analogous to the ones below hold.

We now return to the case $\gamma \equiv 0$ and take $g_1, g_2$ accordingly. Without loss of generality we can assume $g_1, g_2 \geq 0$: we can take $\beta$ small enough and if necessary add a constant to ensure that this holds. We show that for suitable $\alpha_i$ small enough, the functions $g_1,g_2$ satisfy \ref{item:thm_LDP_non_compact:Lyapunov}. For two real-valued functions $f$ and $g$, we write~$f \sim g$ as $x \to \infty$ to say that they asymptotically equivalent in the limit $x \to \infty$, that means~$(f(x)/g(x)))\to 1$ as~$x\to\infty$.
\smallskip

By \ref{item:thm-LDP-zig-zag:Upsi}, $V(x) \to \infty$, and by \ref{item:thm-LDP-zig-zag:Upsii}, $U'(x)/V(x) \to 0$. It follows that  $g_1$ grows to infinity as $|x| \to \infty$. Moreover, $g_2$ grows to infinity by the assumption \ref{item:thm-LDP-zig-zag:Ui} on $U$; since $g_1$ and $g_2$ are continuous, this settles part (a) of \ref{item:thm_LDP_non_compact:Lyapunov}. 

Part (b) of \ref{item:thm_LDP_non_compact:Lyapunov} requires that that $g_2$ grows faster than $g_1$ at infinity. This follows from Assumption \ref{item:thm-LDP-zig-zag:Upsii} on the potentials $U$ and $V$: both dominate the derivative $U'$, and $U$ grows faster than $V$.
\smallskip

To show that \ref{item:thm_LDP_non_compact:Lyapunov} holds for the zig-zag process, we show that both $g_1$ and $g_2$ satisfy
\begin{equation*}
e^{-g_i(x,v)} (B e^{g_i})(x,v) \to -\infty \quad |x|\to\infty.
\end{equation*}
Then since $e^{-g_i(x,v)} (B e^{g_i})(x,v)$ is continuous, the compactness of superlevel-sets follows.
\smallskip

By the definition of $g_1, g_2$ and the extended generator $B$, 
%
\begin{align*}
e^{-g_1(x,v)} (B e^{g_1})(x,v) &= \alpha_1 v V'(x) + \beta U''(x) + \max(vU'(x),0) \left[e^{-2v\beta U'(x)} - 1\right],	
\end{align*}
and
\begin{align*}
e^{-g_2(x,v)} (B e^{g_2})(x,v) &= \alpha_2 v U'(x) + \beta U''(x) + \max(vU'(x),0) \left[e^{-2v\beta U'(x)} - 1\right].
\end{align*}
We first verify the condition for $g_2$. For $v = +1$, we have
\begin{equation*}
e^{-g_2(x,+1)} (B e^{g_2})(x,+1) = \alpha_2 U'(x) + \beta U''(x) + \max(U'(x),0) \left[e^{-2 \beta U'(x)} - 1\right].
\end{equation*}
For $x \to + \infty$, we have $U'(x) \to + \infty$ by~\ref{item:thm-LDP-zig-zag:Ui}, so that 
\begin{align*}
e^{-g_2(x,+1)} (B e^{g_2})(x,+1) &= U'(x) \left[\alpha_2 -1 + \beta \frac{U''(x)}{U'(x)} + e^{-2\beta U'(x)}\right] \\
&\sim U'(x) (\alpha_2 - 1) \to -\infty,\quad x \to + \infty,
\end{align*}
since $U''/U' \to 0$ by~\ref{item:thm-LDP-zig-zag:Uiii} and $\alpha_2 < 1$. 
\smallskip

For $x \to -\infty$, we have $U'(x) \to -\infty$ by~\ref{item:thm-LDP-zig-zag:Ui}, in particular $U'(x) < 0$ for large $x$. Hence,
\begin{align*}
e^{-g_2(x,+1)} (B e^{g_2})(x,+1) &= U'(x) \left[\alpha_2 + \beta \frac{U''(x)}{U'(x)}\right] \\
&\sim U'(x) \alpha_2 \to -\infty,\quad x \to - \infty.
\end{align*}
For $v = -1$, the argument is analogous and we omit the details; this concludes the treatment of $g_2$.
\smallskip

We now consider $g_1$. For $v = +1$, 
\begin{align*}
e^{-g_1(x,+1)} (B e^{g_1})(x,+1) &= \alpha_1 V'(x) + \beta U''(x) + \max(U'(x),0) \left[ e^{-2v\beta U'(x)} - 1\right].
\end{align*}
In the limit $x \to + \infty$, $U'(x) \to +\infty$ and $V'(x)/U'(x) \to 0$ by~\ref{item:thm-LDP-zig-zag:Upsii}. It follows that
\begin{align*}
e^{-g_1(x,+1)} (B e^{g_1})(x,+1) &= \alpha_1 V'(x) + \beta U''(x) + U'(x) \left[ e^{-2v\beta U'(x)} - 1\right]\\
&= U'(x) \left[\alpha_1 \frac{V'(x)}{U'(x)} + \beta \frac{U''(x)}{U'(x)} + e^{-2\beta U'(x)} - 1\right]	\\
&\sim  - U'(x)  \to - \infty, \quad x \to + \infty.
\end{align*}
For $x \to -\infty$, similar to the computations for $g_2$,
\begin{align*}
e^{-g_1(x,+1)} (B e^{g_1})(x,+1) &= \alpha_1 V'(x) + \beta U''(x)\\
&= V'(x) \left[\alpha_1 + \beta \frac{U''(x)}{V'(x)}\right] \to - \infty, \quad x \to + \infty,
\end{align*}
since $U''/V' \to 0$ by~\ref{item:thm-LDP-zig-zag:Upsiii} and $V' \to - \infty$ by~\ref{item:thm-LDP-zig-zag:Upsi}. The case $v = -1$ can be handled using similar arguments. 

The preceding computations conclude the verification of Condition~\ref{item:thm_LDP_non_compact:Lyapunov}. We are left with verifying the mixing property \ref{item:thm_LDP_non_compact:mixing}.

Let $\nu_1,\nu_2 \in \mathcal{P}_c(E)$. Then there exists a compact set $K\subseteq E$ with $\nu_1(K) = \nu_2(K) = 1$. To show that \ref{item:thm_LDP_non_compact:mixing} holds, we must find $T,M > 0$ and $\rho_1,\rho_2 \in \mathcal{P}([0,T])$ such that for all $A \in \mathcal{B}(E)$,
\begin{equation*}
\int_0^T\int_E P(t,y,A) \,d\nu_1(y)d\rho_1(t) \leq M \int_0^T\int_E P(t,y,A) \,d\nu_2(y)d\rho_2(t).
\end{equation*}
By Fubini's theorem, it is sufficient to prove that for any points $y_1 \in \text{supp}(\nu_1)$ and $y_2 \in \text{supp}(\nu_2)$,
\begin{equation}\label{eq:proof-LDP-zigzag:sufficient-Fubini}
\int_0^T P(t,y_1,A)\,d\rho_1(t) \leq M \int_0^T P(t,y_2,A)\,d\rho_2(t),
\end{equation}
with $\rho_1,\rho_2,T,M$ independent of $y_1,y_2$. To that end, let $K\subseteq E$ be a compact set containing the support of both $\nu_1$ and $\nu_2$. Without loss of generality, we can take $K$ of the form $K_\mathbb{R} \times \{\pm 1\}$, where $K_\mathbb{R}$ is a closed interval. For $t_1 > 0$, let $K(t_1)$ be the set of points that the zig-zag process with speed one can reach in the time interval $[0,t_1]$ when starting in $K$:
\begin{equation*}
K(t_1) = \left\{y \in E\,:\, \text{dist}_E(y,K) \leq t_1\right\}.
\end{equation*}
We prove the inequality \eqref{eq:proof-LDP-zigzag:sufficient-Fubini} for arbitrary points $y_1,y_2 \in K$, using the following two steps; in what follows we set $\mu = \Leb \otimes \unif_{\pm 1}$.
\begin{enumerate}[label=(\roman*)]
	\item \label{item:proof-zigzag:step_one} For any $t_1 > 0$ and with $\rho_1$ the uniform distribution over $[0,t_1]$, there is a positive constant $C_{K,t_1}$ depending only on $K$ and $t_1$ such that for any $T > t_1$, we have
	\begin{equation*}
	\int_0^T P(t,y_1,A)\,d\rho_1(t) \leq C_{K,t_1} \cdot \mu\left(A \cap K(t_1)\right), \quad \text{ for all } A \in \mathcal{B}(E),
	\end{equation*}
	with $\mu$ as the reference measure on $\mathcal{B}(E)$.
	\item \label{item:proof-zigzag:step_two} There exist positive constants $T > 0$ and $C'_{K,T}$ such that with $\rho_2$ the uniform distribution over $[0,T]$, we have
	\begin{equation*}
	\int_0^T P(t,y_2,A)\,d\rho_2(t) \geq C'_{K,T}\cdot \mu\left(A \cap K(t_1)\right), \quad \text{ for all } A \in \mathcal{B}(E).
	\end{equation*}
\end{enumerate}
Suppose \ref{item:proof-zigzag:step_one} and \ref{item:proof-zigzag:step_two} hold. Then the estimate \eqref{eq:proof-LDP-zigzag:sufficient-Fubini} also holds, with $M = C_{K,t_1}/C'_{K,T}$, some $t_1 < T$.
\smallskip

To verify \ref{item:proof-zigzag:step_one}, note that the measure
\begin{equation*}
\mu_{y_1}(A) = \int_0^T P(t,y_1,A)\,d\rho_1(t)
\end{equation*}
is absolutely continuous with respect to $\mu= \Leb \otimes \unif_{\pm 1}$ and its density is uniformly bounded in $K$. Now~\ref{item:proof-zigzag:step_one} follows since $P(t,K,A) = 0$ whenever $A \cap K(t_1) = \emptyset$ and $t \leq t_1$.
\smallskip

Next, we use Lemma 8 of \cite{bierkens2019ergodicity} to show \ref{item:proof-zigzag:step_two}. To that end, recall that a tuple $(y,y')$ in $E \times E$ is called reachable if there exists an admissible path from $y$ to $y'$. By Theorem 4 of \cite{bierkens2019ergodicity}, any two points are reachable as long as the potential $U$ has at least one non-degenerate local minimum (which is trivially satisfied on $\mathbb R$ under our assumptions) and satisfies $U\in C^3(\mathbb{R})$. 
\smallskip

By Lemma 8 in \cite{bierkens2019ergodicity}, for any two points $y_a = (x_a,v_a)$ and $y_b = (x_b,v_b)$ in $E$, there are open neighborhoods $U_{y_a}$ of $x_a$ and $U_{y_b}$ of $x_b$, a time interval $(t_0,t_0 + \varepsilon]$ and a constant $c > 0$ such that for all $x_a'\in U_{y_a}$ and $t \in (t_0,t_0 + \varepsilon]$,
\begin{equation*}
P\left(t,(x_a',v_a),A_\mathbb{R} \times \{v_b\}\right) \geq c \, \Leb(A \cap U_{z_b}),\quad \text{ for all }A_\mathbb{R} \in \mathcal{B}(\mathbb{R}).
\end{equation*}
The spatial part of $K \times K(t_1)$ can be covered by open squares associated to all pairs of start and final points $y_a$ and $y_b$, with $y_a \in K$ and $ y_b \in K(t_1)$, as
\begin{equation*}
K_\mathbb{R} \times K(t_1)_\mathbb{R} \subseteq \bigcup_{(y_a,y_b)} U_{y_a} \times U_{y_b},
\end{equation*}
where each $U_{y}$ is an open interval in $\mathbb{R}$. By compactness, there exists a finite subcover by open squares $U_{y_a^i} \times U_{y_b^i}$ corresponding to pairs $(y_a^i,y_b^i)$,
\begin{equation*}
K_\mathbb{R} \times K(t_1)_\mathbb{R} \subseteq \bigcup_{i = 1}^N U_{y_a^i} \times U_{y_b^i}.
\end{equation*}
Thereby, the set $K \times K(t_1) \subseteq E \times E$ is covered as
\begin{equation*}
K \times K(t_1) \subseteq \bigcup_{i = 1}^N \left[ \left(U_{y_a^i}\times \{\pm 1\}\right) \times \left(U_{y_b^i}\times \{\pm 1\}\right)\right].
\end{equation*}
Hence for each $z = (x,v) \in K$, there are finitely many open sets $U_{y_b^i}$ covering $K(t_1)_\mathbb{R}$, with corresponding constants $c_i,t_i,\varepsilon_i$ such that for all $t \in (t_i,t_i+\varepsilon_i]$,
\begin{equation}\label{eq:proof-LDP-zigzag:bound-P-below}
P(t,y,A_\mathbb{R} \times \{v_b^i\}) \geq c_i \Leb(A_\mathbb{R} \cap U_{z_b^i}), \quad \text{ for all } A_\mathbb{R} \in \mathcal{B}(\mathbb{R}).
\end{equation}
For any $A = A_\mathbb{R} \times A_\pm \in \mathcal{B}(E)$, write
\begin{align*}
A^+ &:= A \cap (\mathbb{R} \times \{+1\}),\\
A^- &:= A \cap (\mathbb{R} \times \{-1\}).
\end{align*}
Then with $T > 0$ large enough for all intervals $(t_i,t_i+\varepsilon_i]$ to be contained in $[0,T]$, taking $\rho_2 = \unif([0,T])$, for any $z \in K$ it holds that
\begin{align*}
\int_0^T P(t,y,A)\,d\rho_2(t) &\geq \int_0^T P(t,y,A) \sum_i \mathbf{1}_{(t_i,t_i+\varepsilon_i]}(t)\,d\rho_2(t)	\\
&= \frac{1}{T}\sum_i \int_{t_i}^{t_i+\varepsilon_i} \left[P(t,y,A^+) + P(t,y,A^-)\right]\, dt.
\end{align*}
In each time interval $(t_i,t_i+\varepsilon_i]$, at least one transition probability is bounded from below as in~\eqref{eq:proof-LDP-zigzag:bound-P-below}, while the other one can be bounded from below by zero. Thereby,
\begin{align*}
\int_0^T P(t,y,A)\,d\rho_2(t) &\geq \frac{1}{T} \sum_i \varepsilon_i \cdot c_i\cdot \Leb(A_\mathbb{R} \cap U_{y_b^i})	\\
&\geq \frac{1}{T} \min_i(\varepsilon_i c_i) \sum_i \mu\left[A \cap (U_{y_b^i}\times\{\pm 1\})\right]	\\
&\geq \frac{1}{T} \min_i(\varepsilon_ic_i) \cdot \mu \left[A \cap \bigcup_i \left(U_{y_b^i}\times\{\pm 1\}\right) \right]	\\
&\geq \frac{1}{T} \min_i(\varepsilon_ic_i) \cdot \mu\left[A \cap K(t_1)\right],
\end{align*}
where the last inequality follows from $K(t_1)$ being covered by the $U_{y_b^i}\times\{\pm 1\}$. Hence \ref{item:proof-zigzag:step_two} follows with $C_{K,T}' = \min_i(\varepsilon_i c_i)/T$.
\smallskip

This finishes the verification of Condition \ref{item:thm_LDP_non_compact:mixing}, and thereby the proof of Theorem \ref{thm:LDP-zigzag:non-compact-1d}.
\subsection{Derivation of the explicit form of the rate function}
\label{sec:proofs-rate-function}

Here we prove the results described in Section~\ref{sec:rate}. Recall that the state space is now taken as $E = \mathbb{T} \times \{ \pm 1\}$.

Suppose $\mu$ is absolutely continuous with respect to $\nu_0$ and write $\frac{d \mu}{d \nu_0}(x,v) = \rho(x,v)$ for the Radon-Nikodym density of $\mu$ with respect to $\nu_0$, where $\rho$ is assumed to be absolutely continuous.
Define a mapping $H : \mathcal D^+(L) \rightarrow \R$ by
\begin{equation}
\label{eq:H-functional}
 H (u) := \int_E \frac{ L u}{u}  \ d \mu = \sum_{v\in \{-1,+1\}} \int_{\mathbb{T}} \frac{Lu}{u} (x,v) \rho(x,v) \ d x .
\end{equation}
We compute
\begin{multline}\label{eq:first-steps}
	H(u) =
	\int_{\mathbb{T}} \left\{\frac{d \log u^+}{d x} + \lambda^+ \left( \frac{u^-}{u^+} - 1 \right)  \right\} \rho^+ \ d x 
	\\+
	\int_{\mathbb{T}} \left\{ - \frac{d \log u^-}{d x} + \lambda^- \left( \frac{u^+}{u^-} - 1 \right) \right\} \rho^- \ d x
	\\= 
	\int_{\mathbb{T}} \left\{ -\log u^+ \frac{d \rho^+}{dx} + \lambda^+ \rho^+ \left( \frac{u^-}{u^+} - 1 \right) \right\} \ d x 
	\\+
	\int_{\mathbb{T}} \left \{ \log u^- \frac{d \rho^-}{d x} + \lambda^- \rho^-\left( \frac{ u^+}{u^-} - 1 \right) \right\} \ d x.
\end{multline}
\begin{lemma}
\label{lem:to-infinity}
Suppose $\rho \in C(E)$ is absolutely continuous and satisfies
\[ \nu_0  \left\{\dfrac{d\rho^+}{d x} \neq \dfrac{d \rho^-}{dx} \right\} > 0.\]
Then $\inf_{u \in \mathcal D^+(L)} H(u)=-\infty$.
\end{lemma}
\begin{proof}
Let $u^+_t(x) = u^-_t(x) = \exp \left( - t \left\{ \dfrac{d \rho^+}{d x} - \dfrac{d \rho^-}{d x} \right\}\right)$.
From~\eqref{eq:first-steps} it follows that 
\[ H(u_t) = - t \int_{\mathbb{T}} \left( \dfrac{d \rho^+}{d x} - \dfrac{d \rho^-}{dx} \right)^2 \ d x.\] Now let $t \rightarrow \infty$.
\end{proof}

\begin{lemma}
\label{lem:reparametrization}
Suppose $\rho \in C(E)$ is absolutely continuous and $\frac{d \rho^+}{dx} = \frac{d \rho^-}{dx}$ for all $x \in \mathbb{T}$. Then $\mathcal I$ admits the representation
\begin{equation} \label{eq:I-in-eta} \mathcal I(\mu) = - \inf_{\eta \in C(\mathbb{T})} \int_{\mathbb{T}} \left\{ - \rho' \eta + \lambda^+ \rho^+ (\exp(-\eta) - 1) + \lambda^- \rho^- (\exp(\eta) - 1) \right\} \ d x.\end{equation}
\end{lemma}

\begin{proof}
Write $\rho' : =\frac {d \rho^+}{dx}$, and note that by our assumption $\rho' = \frac{d \rho^-}{dx}$.
By~\eqref{eq:first-steps} we may write
\[ H(u) = \int_{\mathbb{T}} \left\{ - \log \left( u^+/u^- \right) \rho' + \lambda^+\rho^+ (u^-/u^+ -1) + \lambda^- \rho^- (u^+/u^- - 1) \right\} \ d x.\]
We see that only the ratio $u^+/u^-$ determines the value of $H(u)$. To any choice of $u \in \mathcal D^+(L)$ we may associate $\eta = \log u^+ - \log u^- \in C^1(\mathbb{T})$, and correspondingly, to any $\eta \in C^1(\mathbb{T})$ we can associate $u \in \mathcal D^+(L)$ by letting
\[ u^+(x) = \exp(\tfrac 1 2 \eta(x)), \quad u^-(x) = \exp(-\tfrac 1 2 \eta(x)), \quad x \in \mathbb{T}.\]
By the continuous dependence of $H$ on $\eta$, and the fact that $C^1(\mathbb{T})$ is dense in $C(\mathbb{T})$, we obtain the stated representation of $I(\mu)$.
\end{proof}

\begin{lemma}
\label{lem:pointwise-minimization}
Suppose $\rho \in C^1(E)$ and $\frac{ d \rho^+}{d x} = \frac{d \rho^-}{dx}$ for all $x \in \mathbb{T}$. Furthermore suppose $\lambda^- \lambda^+ \rho^- \rho^+ > 0$ on $\mathbb{T}$, and $\lambda^{\pm}$ are continuous. Then $\mathcal I$ is given by~\eqref{eq:rate-function}.
\end{lemma}

\begin{proof}
Differentiating the integrand in~\eqref{eq:I-in-eta} pointwise with respect to $\eta$ gives the first order condition
\[ -\rho' - \lambda^+ \rho^+ \exp(-\eta) + \lambda^- \rho^- \exp(\eta) = 0,\]
which is solved uniquely by 
\[ \eta = \tfrac 1 2 \log \left( \frac{\lambda^+ \rho^+}{\lambda^- \rho^-} \right) + \arcsinh \left( \frac{ \rho'}{2 \sqrt{ \lambda^+ \lambda^- \rho^+ \rho^-}}\right),\]
as long as $\lambda^- \lambda^+ \rho^- \rho^+ \neq 0$. Furthermore $\eta \in C(\mathbb{T})$ by the conditions on $\lambda$ and $\rho$.
The second order derivative with respect to $\eta$ is given by
\[ \lambda^+ \rho^+ \exp(-\eta) + \lambda^- \rho^- \exp(\eta) \geq 0,\]
which shows that the critical value of $\eta$ corresponds to a pointwise global minimum of the integrand. 
\end{proof}
\emph{Proof of Proposition~\ref{prop:expression-ratefunction}}
The result for unequal derivatives is a consequence of Lemma~\ref{lem:to-infinity}, and the expression in case of equality follows from Lemma~\ref{lem:pointwise-minimization}.

\begin{proof}[Proof of Proposition~\ref{prop:no-v-dependence}]
We inspect the dependence of the various terms in the integrand of the expression~\eqref{eq:rate-function}
$I(\mu_c)$ on $c$. For the first term, interchanging integral and derivative,
\begin{align*}
& \frac{d}{dc}  \int_{\mathbb{T}} \rho' \log \left( \frac{\lambda^+ (\rho + c)}{\lambda^- (\rho-c)}\right) \ d x  = \int_{\mathbb{T}} \rho' \left(\frac{1}{\rho + c} + \frac 1 {\rho-c} \right) \ d x \\
 & =  \int_{\mathbb{T}} \frac{d}{dx} \left( \log(\rho +c) + \log(\rho -c) \right) \ d x = 0.
\end{align*}
The following terms (i.e. the $\arcsinh$ and the square root) in the expression for $I(\mu_c)$ are decreasing with respect to the value of $\rho^+ \rho^- = \rho^2 -c^2$. It follows that the integrands are minimized at $c = 0$. Finally, we have that
\[ \int_{\mathbb{T}}  (\lambda^+ \rho^+ + \lambda^- \rho^- )  \ d x = \int_{\mathbb{T}} \left\{  (  \lambda^+ + \lambda^-) \rho + c (\lambda^+ - \lambda^-)  \right\} \ d x.\]
The linear term in $c$ vanishes since $\int_{\mathbb{T}} \{ \lambda^+ - \lambda^- \}  \ d x= \int_{\mathbb{T}} U' \ dx = 0$. It follows that $c = 0$ minimizes $c \mapsto I(\mu_c)$. The stated expression for $I(\mu_0)$ is obtained after a manipulation of~\eqref{eq:rate-function}.
\end{proof}
\chapter{Large Deviations in Stochastic Slow-Fast Systems}
\label{chapter:LDP-in-slow-fast-systems}
\chaptermark{LDP in Stochastic Slow-Fast Systems}
\section{Stochastic slow-fast systems---two time scales}
Our focus in this chapter lies on stochastic systems with two time scales. In various stochastic problems arising for instance in atmospheric models~\cite{BouchetNardiniTangarife2013,BouchetGrafkeTangarifeVandenEijnden2016},
hydrodynamic limits~\cite{KipnisLandim1998}, genetic networks~\cite{CruduDebusscheRadulescu2009}, and statistical physics~\cite{MatosPerez1991,ColletDaiPra2012}, we can identify slow and fast components in the system. The distinction of slow and fast components in a system is based on the observation of different time scales: the fast components converge to their equilibrium state at a much shorter time scale at which the slow components have hardly evolved. As a consequence, the slow components evolve approximately under the averaged effect of the fast components. This observed separation of time scales motivates the term \emph{slow-fast system}. In a stochastic framework, slow-fast systems are frequently modelled by Markov processes that consist of two components, where one component models the slow variables and the other component the fast variables.
\smallskip

A valid approximation of the slow components by averaging over the fast components is also known as the \emph{averaging principle}~\cite{FreidlinWentzell1998}. Establishing an averaging principle in stochastic models has the benefit of rigorously reducing the complexity and leading to simpler models. In stochastic models, this transition from the full system to an approximation via averaging can be justified by the law of large numbers, as demonstrated for instance by Ball, Kan, Kurtz, Popovic and Rempala for the example of reaction networks~\cite{BallKurtzPopovicRempala2006,KangKurtz2013,KangKurtzPopovic2014}. Freidlin and Koralov proved averaging in quasi-linear parabolic PDEs~\cite{FreidlinKoralov2012}, and recently investigated averaging in a slow-fast system whose fast process admits multiple stationary measures~\cite{FreidlinKoralov2020}.
\smallskip

However, the approximation via an averaged evolution is only valid in the limit of infinite time-scale separation. In order to estimate the approximation error, many efforts have therefore concentrated on establishing finer asymptotic results. An example of such an asymptotic result is a pathwise large deviation principle of the slow component as the time-scale separation tends to infinity. Verifying such a large deviation principle is interesting for various reasons. If a large deviation principle is verified for the slow component, then we know the error of the average-approximation to vanish exponentially fast as a function of the time-scale separation. Furthermore, the exponential convergence rate is explicitly known as the so-called large-deviation rate function, which in many situations can be calculated. In the context of multiscale diffusions, Dupuis, Spiliopoulos and Wang show how the large-deviation rate function can be used to design Monte-Carlo  for estimating rare-event probabilities~\cite{DupuisSpiliopoulosWang2012}. The rate function is also the crucial ingredient for characterizing the rare-event behaviour of the system~\cite[Theorem~1.4]{BudhirajaDupuis2019}.  Vanden-Eijnden introduced numerical methods for systems with multiple time scales that do not require to derive the limiting effective equations~\cite{Vanden-Eijnden2003}, and further investigated numerical schemes with Fatkullin~\cite{FatkullinVanden-Eijnden2004} and Weinan and Liu~\cite{WeinanLiuVanden-Eijnden2005,WeinanLiuVanden-Eijnden2007}. We also refer to the monographs of Berglund and Gentz~\cite{BerglundGentz2005}, Kuehn~\cite[Chapter~15]{Kuehn2015} and Pavliotis and Stuart~\cite{PavliotisStuart2008} for more background on both stochastic and deterministic multiscale systems.
\smallskip

Establishing the large deviation principles in slow-fast systems is usually a delicate enterprise, and there has been vivid activity during the last decades to embark on that journey. In their monograph on random perturbations~\cite[Chapter~7]{FreidlinWentzell1998}, Freidlin and Wentzell prove large deviations for several examples where a process is perturbed by a fast process. The two processes are assumed to be weakly coupled in the sense that either the fast process evolves independently of the slow process, or the fast process has a deterministic diffusion coefficient, or the slow process is deterministic. Lipster and Veretennikov also consider a slow diffusion process whose coefficients are coupled to an independent fast diffusion process~\cite{Liptser1996,Veretennikov2000}, and Veretennikov allows for a weak coupling between diffusions in~\cite{Veretennikov1999}, similar to Freidlin and Wentzell. For coupled diffusions where the fast diffusion coefficient is indepedent of the slow process, Feng and Kurtz offer a proof based on Hamilton-Jacobi theory~\cite[Section~11.6, Lemma~11.60~(1)]{FengKurtz2006}. Kifer studies ODE's coupled to fast diffusions~\cite{Kifer1992,Kifer2009}. Bouchet, Grafke, Tangarife and Vanden-Eijnden~\cite{BouchetGrafkeTangarifeVandenEijnden2016} complement this study by specifying the ODEs to concrete examples in order to calculate the Hamiltonians, thereby obtaining a more explicit rate function.
\smallskip

Puhalskii studies fully coupled slow-fast diffusions~\cite{Puhalskii2016} by building up on Lipster's method of considering the joint distribution of the slow process and the empirical measure of the fast process~\cite{Liptser1996}. Spiliopoulos proves large deviations on the path level, and provides importance-sampling schemes for coupled diffusions~\cite{Spiliopoulos2013}. Feng, Fouque and Kumar prove large deviations for the time marginals of a slow diffusion process coupled to fast diffusions~\cite{FengFouqueKumar2012}, and Ghili provides a generalization of their results~\cite{Ghilli2018}.
\smallskip

More recently, fast components modelled by jump processes attracted more interest in the large-deviation context. He and Yin couple diffusions to fast jump processes~\cite{HeYin2014}, considering different scaling regimes of the time-scale separation. Similar in spirit to Puhalskii's paper, Huang, Mandjes and Spreij apply Lipster's idea from~\cite{Liptser1996} to prove large deviations of a slow diffusion process coupled to an independent fast jump process, by first proving joint large deviations of the slow process and the empirical measure process, and then using the contraction principle~\cite{HuangMandjesSpreij2016}. Bressloff and Faugeras start from large-deviation results and derive action-integral representations via contraction~\cite{BressloffFaugeras2017}. Budhiraja, Dupuis and Ganguly~\cite{BudhirajaDupuisGanguly2018} proof process-level large deviations of a slow diffusion process and fast jumps, with fully-coupled components. The rate functions are characterized via an optimal control problem, involving the empirical measure of the fast variable. Popovic and Kumar~\cite{KumarPopovic2017} tackle the general case where both slow and fast components are mixed jump-diffusion processes. They show that under the assumption of well-posedness of a certain Hamilon-Jacobi equation, the one-dimensional time marginals satisfy large deviations.
\smallskip

Despite the enormous interest and the huge literature on the topic of slow-fast systems, the important class of physical models of mean-field interacting particles described by jump processes on a finite state space has not been treated so far in the context of slow-fast systems. These Markovian jump models are frequently consulted as approximations to physical models describing certain non-equilibrium phenomena, such as spin dynamics. An overview involving different spin models is offered for instance by Martinelli~\cite{Martinelli1999}. A typical example is the Glauber dynamics in Ising-models and Potts-models describing ferromagnets. Other fields of applications include communication networks~\cite{AntunesFrickerRobert2006}, game theory with models involving a large number of agents~\cite{GomesMohrSouza2010}, and chemical reactions~\cite{MielkePattersonPeletierRenger2017}.
\smallskip

There is very recent activity in the study of Markovian mean-field jump processes from a large-deviation perspective. Dupuis, Ramanan and Wu prove large deviations of the empirical densities~\cite{DupuisRamananWu2016}, the clue being to allow for more than one jump simultaneously. In another paper with Fischer, they investigate the stability of the nonlinear limit evolution equation of the particle system by constructing Lyapunov functions from relative entropies~\cite{BudhirajaDupuisFischerRamanan2015}. Renger proves large deviations of density-flux pairs of non-interacting particles exploiting Girsanov transformations~\cite{Renger2017}, which Kraaij extended to include weak interactions~\cite{Kr17}. Bertini, Chetrite, Faggionato, and Gabrielli consider a mean-field system with deterministic time-periodic rates~\cite{BertiniChetriteFaggionatoGabrielli2018}, and prove large deviations in the large number of particles limit. Budhiraja and Wu also consider moderate deviations~\cite{BudhirajaWu2017}.
\smallskip

In our work, we contribute to these very recent studies by proving
dynamic large deviation principles in mean-field interacting particles coupled to fast external processes. In general, the main methods used in the literature to prove large deviations in slow-fast systems are the weak-convergence method (\cite{BudhirajaDupuis2019}), classical techniques based on Girsanov transformations, and the method based on convergence of nonlinear generators and Hamilton-Jacobi Theory~\cite{FengKurtz2006}. Despite the interest in mean-field systems, there are few results illuminating the large-deviation behaviour of mean-field particles from a Hamilton-Jacobi point of view. A system of interacting diffusions is considered in~\cite[Chapter~13]{FengKurtz2006}. Feng, Mikami and Zimmer extend the methods therein to prove the comparison principle for equations involving Hamiltonians that arise in this context~\cite{FengMikamiZimmer2019}.
Moreover, proofs about large deviations in coupled systems, like slow-fast systems, assume well-posedness of the comparison principle rather than verifying the comparison principle. The main novelties presented here are the following:
\begin{itemize}[]
	\item We provide a general set of conditions under which we prove pathwise large deviations of slow components in slow-fast systems via Hamilton-Jacobi equations. The conditions allow for irreversible fast processes.
	\item We find Lagrangian rate functions. Next to the standard characterization of the Lagrangian in terms of the dual of a principal eigenvalue, we establish a characterization in terms of a double-optimization.
	\item As our main example, we treat density-flux large deviations of mean-field interacting particles on a finite state space coupled to fast drift-diffusion processes on a compact periodic space. This example requires arguments that are different from those currently available in the literature. We derive an averaging principle from the large deviation principle.
	\item The large-deviation results apply to small-diffusion processes coupled to fast jump processes. This solves a challenge pointed out by Budhiraja, Dupuis and Ganguly in~\cite{BudhirajaDupuisGanguly2018}, which is the fact that in slow-fast systems, classical results about comparison principles are not applicable due to the Hamiltonians having poor regularity properties.
\end{itemize} 
\paragraph{Overview of this chapter.} 
In Section~\ref{SF:sec:toy-ex}, we treat two toy examples of stochastic slow-fast systems. The first toy example 
shows
how fast variables affect the large-deviation behaviour of the slow variables. The second toy example illustrates the characterization of the Lagrangians in terms of a double-optimization, and it's connection to averaging principles. The examples provide the necessary background to have in picture in mind for the general results that follow.
\smallskip

In Section~\ref{SF:sec:main-results}, we state our main results: a general large-deviation theorem for slow components in a slow-fast system (Theorem~\ref{SF:thm:LDP_general}), an action-integral form of the rate functions (Theorem~\ref{SF:thm:rate-function}), a large-deviation theorem for mean-field interacting particles coupled to a fast diffusion process (Theorem~\ref{SF:thm:LDP-mean-field}), and an application of the large deviation principle to the averaging principle for mean-field systems (Theorem~\ref{SF:thm:mean-field-averaging}).
In Section~\ref{SF:sec:assumptions:general-LDP-theorem} we collect the assumptions for the general large-deviation results. The remaining sections contain the proofs.
\section{Two toy examples}\label{SF:sec:toy-ex}
\subsection{Fastly-varying diffusion}
We start from the large deviation principle of a small-diffusion process. Then we couple this process to a fast process to illustrate a slow-fast system.
\paragraph{Single component.}
For fixed $n \in \mathbb{N}$ and a fixed positive constant~$\sigma>0$, consider the stochastic process $Y^n_t \in \mathbb{R}$ solving
\begin{equation*}
	\dd Y^n_t = \frac{1}{\sqrt{n}} \sigma \cdot \dd B_t,
	\quad Y^n(0) = 0. 
\end{equation*}
We call the constant~$\sigma$ the diffusion coefficient, and $B_t$ denotes Brownian motion in $\mathbb{R}$. By Schilder's Theorem~(\cite[Theorem~5.2.3]{DemboZeitouni1998}), the process $Y^n_t$ satisfies a pathwise large deviation principle in the small-diffusion limit $n \to \infty$. The rate function is given by (see also Section~\ref{BG:sec:semigroup-flow-HJ-eq} in Chapter~2)
\begin{equation*}
	J(y) = \int_0^\infty \frac{1}{2 a} |\partial_t y(t)|^2 \, \dd t,
	\quad a := \sigma^2.
\end{equation*}
This corresponds to the Hamiltonian~$\mathcal{H}(p)=ap^2/2$.
\paragraph{Slow-fast system.}
To transition from the single-component process to a slow-fast system, we pass from the constant diffusion coefficient to a switching diffusion coefficient. The switching times depend on a jump process~$Z^n_t$ flipping between~$\pm 1$. The slow-fast system is defined as a two-component process~$(X_t^n,Z_t^n)$, where the first component evolves according to
\begin{equation*}
	\dd X_t^n = \frac{1}{\sqrt{n}}\sigma\left(Z_t^n\right) \dd B_t,\quad X^n(0)=0.
\end{equation*}
The jump process evolves independently of the diffusion~$X_t^n$. In between the jump times of~$Z_t^n$, the dynamics of~$X_t^n$ is just the dynamics of the process~$Y_t^n$ from above, where depending on the value of the second component $Z^n_t$, the diffusion constant is either given by $\sigma_- := \sigma(-1)$ or by $\sigma_+ := \sigma(+1)$. To model the time-scale separation, we take~$Z_t^n$ to evolve with jump rates given by~$r_n(-1,+1)=n\cdot r_-$ and~$r_n(+1,-1)=n\cdot r_{+}$ for some fixed~$r_\pm > 0$. In this two-component process,~$X_t^n$ is the slow and~$Z_t^n$ is the fast component.
\smallskip

The typical questions we ask about such a slow-fast system are: does the slow component satisfy a pathwise large deviation principle? How does the dependence on~$Z_t^n$ affect the large-deviation behaviour? In this toy example, both questions can be answered explicitly. The first question is answered by the fact that the slow component~$X_t^n$ indeed satisfies pathwise large deviations with some rate function~$J_\pm$. 
The second question is answered by describing this rate function. As in the single-component version, the rate function is of action-integral form. The Lagrangian is the Legendre-Fenchel transform $\mathcal{L}_\pm(v)$ of a Hamiltonian $\mathcal{H}_\pm : \mathbb{R} \to \mathbb{R}$. However, due to the fast jump process, the Hamiltonian is no longer quadratic, but given by
\begin{multline}
	\label{SF:eq:intro:Hamiltonian-vastly-varying-diffusion}
	\mathcal{H}_\pm(p) =
	\frac{p^2}{4}(a_{+} + a_{-}) - \frac{1}{2}(r_{+} + r_{-})\\
	+ \frac{1}{2}\sqrt{(r_{+} + r_{-})^2 + p^2(a_{+} - a_{-})(r_{-}-r_{+}) + \frac{1}{4}p^4(a_{+} - a_{-})^2},
\end{multline}
where~$a_\pm := \sigma_\pm^2$.
This Hamiltonian is the principal eigenvalue of the matrix
\begin{equation*}
	M(p) = \frac{p^2}{2}
	\begin{pmatrix}
		a_-	&	0\\
		0	&	a_+
	\end{pmatrix}
	+
	\begin{pmatrix}
		-r_-	&	r_-\\
		r_+	&	-r_+
	\end{pmatrix}.
\end{equation*} 
\begin{figure}[h!]
	\labellist
	\pinlabel $v$ at 2400 -80
	\pinlabel $v$ at 4700 -80
	\pinlabel $\mathcal{L}(v)/\mathcal{L}_\pm(v)$ at 3700 850
	\pinlabel $\sigma_+^2$ at 3000 1500
	\pinlabel {\color{dark_blue}{$\mathcal{L}(v)$}} at 600 1000
	\pinlabel {\color{red_one}{$\mathcal{L}_\pm(v)$}} at 2300 580
	\endlabellist
	\begin{center}
		\includegraphics[scale=.07]{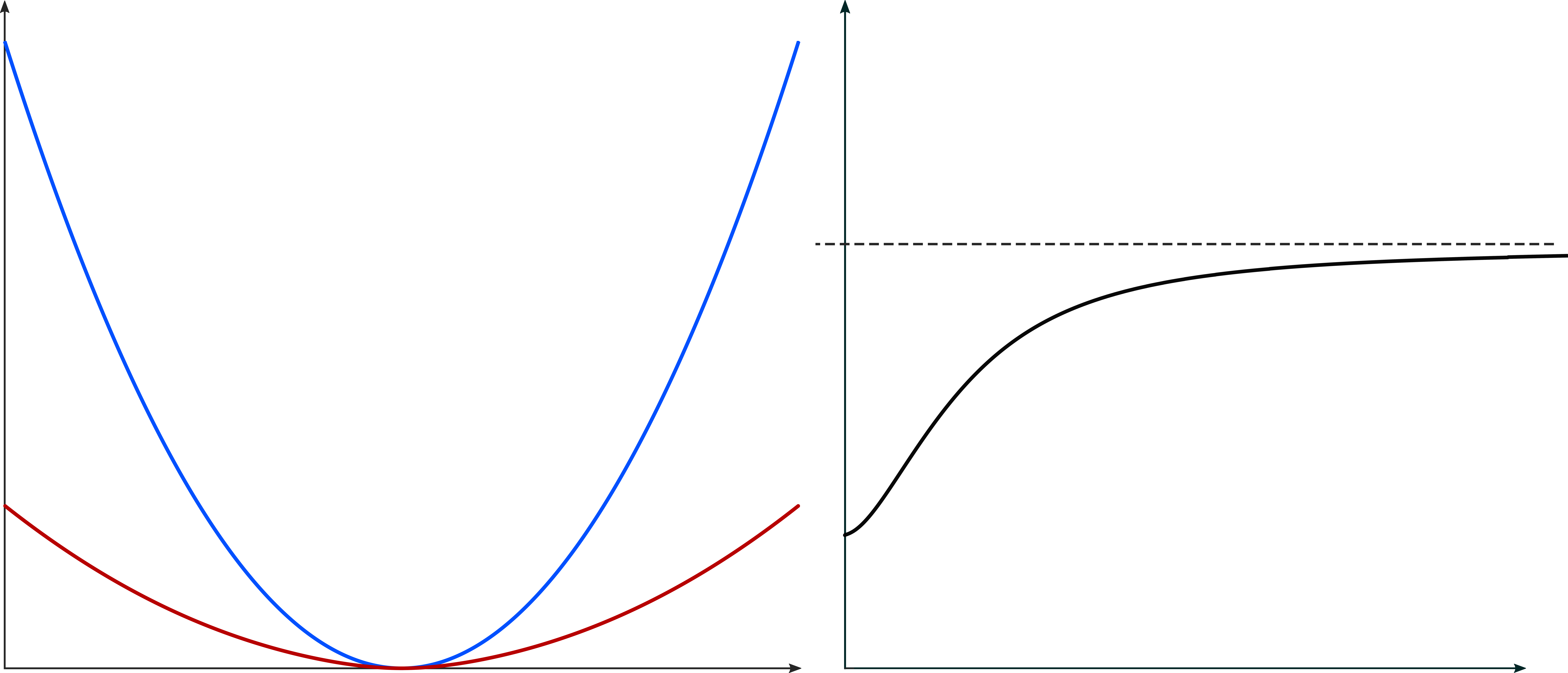}
	\end{center}
	\caption{In blue the quadratic Lagrangian with a constant~$\sigma>0$. In red the Lagrangian of with flipping between two values~$\sigma_\pm$ with~$\sigma_-=\sigma$ and~$\sigma_+>\sigma_-$. In the right graph the quotient~$\mathcal{L}/\mathcal{L}_\pm$.}
	\label{SF:fig:toy-ex:Lagrangians}
\end{figure}
\smallskip

If the diffusion coefficients are equal, then $a_+=a_-$, and we recover the above quadratic Hamiltonian for the small diffusion limit. If the diffusion coefficient is flipping between two values~$\sigma_\pm$ with~$\sigma_+>\sigma_-$, then asymptotically, the exit probabilities
\begin{equation*}
	\mathbb{P}\left[X^n(T)\geq C\right] \sim e^{-n T\mathcal{L}_\pm(C/T)},\quad n\to\infty,
\end{equation*}
are larger. This is illustrated in Figure~\ref{SF:fig:toy-ex:Lagrangians}, which depicts the fact that the Lagrangian~$\mathcal{L}_\pm$ is smaller with flipping than without flipping, provided $\sigma=\sigma_-$ and~$\sigma_+>\sigma_-$. Indeed, one can verify for instance the estimate
\begin{equation*}
\sigma_+^2\frac{1}{1+\sigma_+^2/v^2} \leq \left[\mathcal{L}(v)/\mathcal{L}_\pm(v)\right] \leq \sigma_+^2.
\end{equation*}
In this very specific example, the fast variable does not influence the law of large number limit: in both cases, uncoupled or coupled, the deterministic limit of the small-diffusion process is the path constant equal to zero. While in this simplified toy example, we can find an explicit formula for the Hamiltonian, this is no longer the case in more complicated systems. However, also in more involved systems the large-deviation behaviour is goverened by a principal-eigenvalue problem derived from the slow-fast system---see Section~\ref{SF:sec:proof-general-LDP-Theorem}. We close this example by formally deriving~$\mathcal{}H_\pm(p)$ from the generators~$L_n$ of~$(X^n,Z^n)$,
\begin{equation*}
L_nf(x,z)=\frac{1}{2n}\sigma(z)^2\Delta_x f(x,z) + n\cdot r(z,z')(f(x,z')-f(x,z)),\quad z':=(-1)\cdot z.
\end{equation*}
The nonlinear generators~$H_nf=n^{-1}e^{-nf}L_ne^{nf}$ are given by
\begin{equation*}
H_nf(x,z)=\frac{1}{2n}\sigma(z)^2 \Delta_x f + \frac{1}{2}\sigma(z)^2|\nabla_xf|^2 + r(z,z')\left[e^{n(f(x,z')-f(x,z))}-1\right].
\end{equation*}
Consider functions of the form~$g_n(x,z)=g(x)+n^{-1}\varphi(z)$ to take the scale separation into account. Then
\begin{align*}
H_ng_n(x,z)\xrightarrow{n\to\infty} H_{g,\varphi}(x,z):=\frac{1}{2}\sigma(z)^2|\nabla g(x)|^2 +r(z,z') \left[e^{\varphi(z')-\varphi(z)}-1\right].
\end{align*}
We want this limit to be independent of the fast variable~$z$. We fix~$x$, and thereby also~$p=\nabla g(x)$. By Perron-Frobenius type results, there exists a unique eigenvalue~$\lambda(p)\in\mathbb{R}$ and a vector~$\varphi=(\varphi(-),\varphi(+))$ such that
\begin{equation*}
\frac{1}{2}\sigma(z)^2p^2 +r(z,z') \left[e^{\varphi(z')-\varphi(z)}-1\right]=\lambda(p).
\end{equation*}
This eigenvalue~$\lambda(p)$ is precisely the Hamiltonian~$\mathcal{H}_\pm(p)=\lambda(p)$, and its explicit formula is~\eqref{SF:eq:intro:Hamiltonian-vastly-varying-diffusion} from above.
\subsection{Fastly-alternating drifts}
We slightly vary the previous toy example to illustrate how the Lagrangians in form of a double-optimization provide us with a convenient tool to connect the large-deviation results to the averaging principle.
\paragraph{Single component.} For fixed~$u\in\mathbb{R}$ and~$n\in\mathbb{N}$, consider
\begin{equation*}
\dd Y_t^n = u \,\dd t + \frac{1}{\sqrt{n}}\,\dd B_t,\quad Y^n(0)=0.
\end{equation*}
Then~$\{Y^n\}_{n=1,2,\dots}\in C_\mathbb{R}[0,\infty)$ (equipped with the Skorohod topology) satisfies a large deviation principle by the Freidlin-Wentzell Theorem, with rate function
\begin{equation*}
J(y) = \int_0^\infty \frac{1}{2}|\partial_t y (t) - u|^2\,\dd t.
\end{equation*}
In particular, the limit of~$Y^n$ as~$n\to\infty$ is the path with constant velocity~$u$.
\paragraph{Slow-fast system.} We consider the same setting as above, but let the velocity flip fastly between two values~$u_+,u_-\in\mathbb{R}$. As in the previous toy example, this means we introduce a jump process~$Z^n$ on~$\{-1,+1\}$ and consider the two-component process~$(X^n,Z^n)$, where~$X^n$ follows the dynamics
\begin{equation*}
\dd X_t^n = u(Z_t^n)\,\dd t + \frac{1}{\sqrt{n}}\,\dd B_t,\quad X^n(0)=0,
\end{equation*}
with~$u(\pm 1)=u_\pm$. The jump process~$Z^n$ evolves independently of~$X^n$, with the same jump rates of order~$n$ as in the first toy example. Let~$\pi=\pi_-\delta_-+\pi_+\delta_+$ be the stationary measure of~$Z^n$. Since~$Z^n$ equilibriates fastly, we expect~$X^n$ to converge to the path~$\overline{x}$ with constant average velocity~$\overline{u}=u_-\pi_-+u_+\pi_+$. This statement is an example of the averaging principle.
\smallskip

Let us see how to obtain the averaging principle from a large-deviation perspective. Here, we argue heuristically:
\begin{enumerate}[(a)]
	\item The jump process~$Z^n$ equilibriates exponentially fast at~$\pi$ with speed~$n$,
	\begin{equation*}
	\mathbb{P}\left(\int_0^1 \mathbf{1}_{Z^n(t)}(\cdot)\,\dd t \approx \mu(\cdot)\right) \sim \exp\{- n \mathcal{I}(\mu)\},\quad n\to\infty,
	\end{equation*}
	with the Donsker-Varadhan functional~$\mathcal{I}:\mathcal{P}(\{-1,+1\})\to[0,\infty)$
	\begin{equation*}
	\mathcal{I}(\mu) = \sup_{\xi\in\mathbb{R}^2} \left[r_-\mu_-(1-e^{\xi_+-\xi_-}) + r_+\mu_+(1-e^{\xi_--\xi_+})\right],
	\end{equation*}
	and~$\mathcal{I}(\mu)=0$ if and only if~$\mu=\pi$.
	\item Suppose that~$Z^n$ remains stationary at law~$\pi$. Then Freidlin-Wentzell large deviations suggest that for large~$n$,
	\begin{multline*}
	\mathbb{P}\left(X^n\approx x\,|\,\mathrm{law}(Z^n)=\pi\right) \\
	\sim \exp\left\{-n\int_0^\infty\left[\frac{1}{2}|\partial_tx(t)-u_-|^2 \pi_- + \frac{1}{2}|\partial_tx(t)-u_+|^2 \pi_+\right]\,\dd t\right\}.
	\end{multline*} 
\end{enumerate}
When taking the limit~$n\to\infty$, both the convergence of~$Z^n$ to equilibrium and the convergence of~$X^n$ to a path with constant velocity are competing at the same scale or order~$n$. Therefore both effects contribute to the probabilistic cost when computing the rate function of the slow component deviating from its most likely path. To observe a velocity~$v$ of the process~$X^n(t)$ in a small time-window~$[t,t+\Delta t)$, three events can contribute: the fast process is distributed as~$\mu$ (instead of~$\pi$); while~$Z^n$ is in state~$-1$, the slow component's velocity is~$v_-$ (instead of~$u_-$); while~$Z^n$ is in state~$+1$, the slow component's velocity is~$v_+$ (instead of~$u_+$). The only condition for observing~$v$ is $v=v_-\mu_- + v_+\mu_+$. Since the winner takes it all, the rate~$\mathcal{L}(v)$ is obtained by optimizing over~$\mu,v_-,v_+$, and taking into account the cost of each event:
\begin{equation}\label{SF:toy-ex:Lagrangian-fast-flip-drift}
\mathcal{L}(v)=\inf_{\substack{\mu,v_-,v_+\\v_-\mu_-+v_+\mu_+=v}}\frac{1}{2}|v_--u_-|^2\mu_- + \frac{1}{2}|v_+-u_+|^2\mu_+ + \mathcal{I}(\mu).
\end{equation}
When specializing Theorem~\ref{SF:thm:rate-function} to this example, we find that~$X^n$ satisfies pathwise large deviations with rate function~$J$ with this Lagrangian,
\begin{equation*}
J(x) = \int_0^\infty \mathcal{L}(\partial_t x(t))\,\dd t.
\end{equation*}
This is an instant of a general principle expressed by Theorem~\ref{SF:thm:rate-function}, where we prove this double optimization to hold under the same conditions under which we prove large deviations. 
\smallskip

We close this example by deriving the averaging principle from this rate function. As a consequence of the large deviation principle,~$X^n\to \overline{x}$ almost surely, where~$J(\overline{x})=0$ (Theorem~\ref{thm:math-formulation-LDP:LDP-implies-as}). Hence~$\mathcal{L}(\partial_t\overline{x}(t))=0$ with the Lagrangian~\eqref{SF:toy-ex:Lagrangian-fast-flip-drift}. Since all three terms in the Lagrangian are non-negative, each term must vanish. Thus~$(\mu,v_-,v_+)=(\pi,u_-,u_+)$ is the optimizer in this case, and we must have~$\partial_t\overline{x}(t) = u_-\pi_-+u_+\pi_+$.
\section{Main results}\label{SF:sec:main-results}
\subsection{Preliminaries}\label{SF:sec:preliminaries}
For a Polish space~$E$, we denote by $C(E)$ and $C_b(E)$ the spaces of continuous and bounded continuous functions respectively. If $E \subseteq \bR^d$ then we denote by $C_c^\infty(E)$ the space of smooth functions that vanish outside a compact set in $E$. We denote by $C_{cc}^\infty(E)$ the set of smooth functions that are constant outside of a compact set, and by $\cP(E)$ the space of probability measures on $E$. We equip $\cP(E)$ with the weak topology, that is, the one induced by convergence of integrals against bounded continuous functions.
\smallskip

We consider Markov processes defined via solutions to the martingale problem of a linear operator~$A : \cD(A) \subseteq C_b(E) \rightarrow C_b(E)$. We write~$\mathcal{X}:=D_E[0,\infty)$ for the Skorokhod space of trajectories that are right-continuous and have left limits, equiped with its usual topology~\cite[Section~3.5]{EthierKurtz1986}.
\begin{definition}
	Let $\mu \in \cP(E)$. We say that a measure $\PR \in \cP(\mathcal{X})$ solves \textit{the martingale problem} for $(A,\mu)$ if for all $f \in \cD(A)$ the process
	\begin{equation*}
		M_f(t) := f(X(t)) - f(X(0)) - \int_0^t Af(X(s)) \dd s
	\end{equation*}
	is a martingale with respect to the filtration $t \mapsto \cF_t := \left\{X(s) \, | \, s \leq t\right\}$, and if the projection of $\PR$ on the time $0$ coordinate equals $\mu$. We say that~$\PR \in \cP(\mathcal{X})$ solves the martingale problem for $A$ if it solves the martingale problem for $(A,\mu)$ for some starting measure~$\mu$. We say that the process $\{X(t)\}_{t \geq 0}$ on $\mathcal{X}$ solves the martingale problem for $A$ if its distribution solves the martingale problem. A martingale problem is \textit{well-posed} if there exists a unique solution to the martingale problem for each starting measure.\qed
\end{definition}
\subsection{Setting of slow-fast systems}\label{SF:sec:general-setting}
Here we introduce stochastic \emph{slow-fast systems} as certain two-component Markov processes~$(Y_t^n,Z_t^n)$, where the first component corresponds to the slow process, the second to the fast process. To incorporate the feature of being \emph{slow-fast}, we include a scaling parameter $r_n$ that introduces a separation of time-scales for the two processes.
\smallskip

We first fix the state space of a slow-fast system. To focus only on the features that arise due to the coupling of slow and fast variables, we will assume that the fast process $Z_t^n$ takes values in a compact Polish space $F$. This compactness assumption, as well as the fact that $F$ does not depend on $n$, can both be relaxed at the cost of more but non-trivial technicalities~(e.g.~\cite{Ghilli2018}). Furthermore, for each $n = 1, 2, \dots$, the slow process $X_t^n$ takes values in Polish spaces $E_n$ such that  $\eta_n(E_n) \subseteq E \subseteq \mathbb{R}^d$, where $\eta_n : E_n \to E$ is a continuous embedding and $E$ is a Polish space as well. We assume that $E$ is contained in the $\mathbb{R}^d$-closure of its $\mathbb{R}^d$-interior, which ensures that gradients of functions on $E$ are determined by the values of the function in $E$. The setting of the state spaces is summarized in the following basic condition.
\begin{condition}[Basic condition on the state spaces $E_n$ and $F$]
	\label{SF:condition:compact_setting:state-spaces}
	The state space $F$ is a compact Polish space. The state spaces $E_n$ are Polish spaces that are asymptotically dense in $E \subseteq \mathbb{R}^d$ with respect to continuous embeddings $\eta_n : E_n \to E$; that means for any $x \in E$, there exist $x_n \in E_n$ such that $\eta_n(x_n) \to x$ as $n \to \infty$. Furthermore, suppose that for each compact $K \subseteq E$ the set $\eta_n^{-1}(K)$ is compact in $E_n$ and that there exists a compact set $\widehat{K} \subseteq E$ such that 
	\begin{equation*}
		K \subseteq \liminf_n \eta^{-1}_n(\widehat{K}).
	\end{equation*}
	The last condition means that for every compact	$K \subseteq E$ there is a compact set $\widehat{K} \subseteq E$ such that for all $x \in K$ there is an increasing map $k : \bN\rightarrow \bN$ and $x_{k(n)} \in \eta_{k(n)}^{-1}(\widehat{K})$ such that $\lim_n \eta_{k(n)}(x_{k(n)}) = x$.\qed
\end{condition}
We consider two-component processes~$(Y_t^n,Z_t^n)$ defined by generators that decompose into slow and fast parts in the following sense.
\begin{definition}[Generator of slow-fast system]\label{SF:def:generator-slow-fast}
	We say that a sequence of linear operators $A_n : \mathcal{D}(A_n) \subseteq C_b(E_n \times F) \to C_b(E_n \times F)$ corresponds to a slow-fast system if $A_n$ is given by
	\begin{equation}\label{SF:eq:setting:generator-slow-fast-system}
	A_n f(y,z) := A^\mathrm{slow}_{n,z} f(\cdot,z) (y) +
	r_n \cdot A^\mathrm{fast}_{n,y} f(y,\cdot)(z),
	\end{equation}
	where~$r_n$ is a sequence of positive real numbers such that $r_n \rightarrow \infty$ and
	\begin{enumerate}[label = (\roman*)]
		\item for each $z \in F$ and $n = 1,2,\dots$, there is a generator 
		\begin{equation*}
			A^\mathrm{slow}_{n,z} : \mathcal{D}(A^\mathrm{slow}_{n}) \subseteq C_b(E_n) \to C_b(E_n)
		\end{equation*}
		of an $E_n$-valued Markov process $Y^n_t$. The domain of $A^{\mathrm{slow}}_{n,z}$ is independent of $z$, denoted by $\mathcal{D}(A^\mathrm{slow}_n)$. For $f \in \cD(A_n)$, we have $f(\cdot,z) \in \cD(A^\mathrm{slow}_n)$.
		\item For each $y \in E_n$, there is a generator
		\begin{equation*}
			A^\mathrm{fast}_{n,y} :  \mathcal{D}(A^\mathrm{fast}) \subseteq C(F) \to C(F)
		\end{equation*}
		of a Markov process on $F$. The domain is independent of $n$ and $y$, denoted by $\mathcal{D}(A^\mathrm{fast})$. For $f \in \cD(A_n)$, we have $f(y,\cdot) \in \mathcal{D}(A^\mathrm{fast})$.\qed
	\end{enumerate}
\end{definition}
The sequence of positive numbers~$r_n$ models the time-scale separation of the two processes. The fast component runs at a time scale of order~$r_n$ compared to the slow component. In the law of large number limit, the time separation tends to infinity.
\smallskip

For a sequence of slow-fast systems constructed form operators~ $A_n$ defined as above, we make the following well-posedness assumption regarding solvability of the associated martingale problem.
\begin{condition}[Well-posedness of martingale problem]
	\label{SF:condition:compact_setting:well-posedness-martingale-problem}
	Consider a slow-fast system constructed from operators $A_n$ as in Definition~\ref{SF:def:generator-slow-fast}.
	For each $n \in \mathbb{N}$ and each initial distribution $ \mu \in \mathcal{P}(E_n \times F)$, existence and uniqueness hold for the $(A_n,\mu)$-martingale problem on the Skorohod-space $D_{E_n \times F}[0,\infty)$. Denote the Markov process solving the martingale problem by~$(Y_n(t),Z_n(t))$. The mapping $(y,z) \mapsto P^n_{y,z}$ of $E_n \times F$ into $\mathcal{P}(D_{E_n \times F}[0,\infty))$ is continuous with respect to the weak topology on $\mathcal{P}(D_{E_n \times F}[0,\infty))$, where~$P^n_{y,z}$ is the distribution of the Markov process~$(Y_n(t),Z_n(t))$ starting at~$(y,z)$.\qed
\end{condition}
The assumption that $(y,z) \mapsto P^n_{y,z}$ is continuous is equivalent to the fact that the semigroup of the process~$(Y_n,Z_n)$ is Feller-continuous and that the set $\{P_{x,z}:(x,z)\in K\}$ is tight for any compact $K\subseteq E_n\times F$~(\cite[Remark~11.22]{FengKurtz2006}. We assume continuity in order to apply the general large-deviation results established by Kraaij in~\cite{Kraaij2019ExpResolv}.
\subsection{General large-deviation theorem}
\label{SF:sec:result:general-LDP}
We formulate our main result for general slow-fast systems by two theorems. All assumptions stated in the theorems are collected in Section~\ref{SF:sec:assumptions:general-LDP-theorem}.
For the large-deviation result, we consider a sequence of slow-fast systems~$(Y_t^n,Z_t^n)$ with values in~$E_n\times F$ satisfying Condition~\ref{SF:condition:compact_setting:state-spaces}. The following theorem establishes conditions under which slow components~$Y_t^n$ satisfy a pathwise large deviation principle. Recall the state spaces~$E_n$ imbedded into~$E$ by a continuous map~$\eta_n$.
\begin{theorem}[Large-deviation principle of slow component]\label{SF:thm:LDP_general}
	Let ~$(Y_t^n,Z_t^n)$ be a slow-fast system satisfying the well-posedness Condition~\ref{SF:condition:compact_setting:well-posedness-martingale-problem}. In addition, suppose that Assumptions~\ref{SF:assumption:convergence-slow-nonlinear-generators},~\ref{SF:assumption:convergence-fast-nonlinear-generators},~\ref{SF:assumption:principal-eigenvalue-problem},~\ref{SF:assumption:regularity-V} and~\ref{SF:assumption:regularity-I} are satisfied.
	Let $X_n := \eta_n(Y_n)$ and suppose the large deviation principle holds for~$X_n(0)$ on~$E$ with speed~$r_n$ and rate function~$J_0$.
	Then the process~$X_n$ satisfies a large deviation principle in~$D_E[0,\infty)$ with speed~$r_n$ and with rate function~$J$ given in~\eqref{SF:eq:Richard-general-LDP-thm-rate-function} in Section~\ref{SF:sec:proof-general-LDP-Theorem}. 
\end{theorem} 
We prove Theorem~\ref{SF:thm:LDP_general} in Section~\ref{SF:sec:proof-general-LDP-Theorem}. The rate function is only implicitly characterized by the limit of nonlinear semigroups associated to the slow-fast system, and is therefore not satisfying. This is why we establish two alternative representations of the rate function. These alternative representations establish the rate function as a time-integral over a Lagrangian, which is why we call it action-integral representation.
\begin{theorem}[Action-integral representation]\label{SF:thm:rate-function}
	In addition to the assumptions of Theorem~\ref{SF:thm:LDP_general}, suppose that also Assumption~\ref{SF:assumption:Hamiltonian_vector_field} is satisfied. Then there exists a map~$\mathcal{L}:E\times\mathbb{R}^d\to[0,\infty]$ such that the rate function~$J$ from Theorem~\ref{SF:thm:LDP_general} is
	\begin{equation}\label{SF:eq:action-integral-rate-function}
		J(\gamma) = \begin{cases}
			J_0(\gamma(0)) + \int_0^\infty \cL(\gamma(s),\dot{\gamma}(s)) \dd s & \text{if } \gamma \in \cA\cC([0,\infty);E), \\
			\infty & \text{otherwise}.
		\end{cases}
	\end{equation}
	The map~$\mathcal{L}:E\times\mathbb{R}^d\to[0,\infty]$ admits the two representations specified below in~\eqref{SF:eq:Lagrangian-is-Legendre-dual} and~\eqref{SF:eq:results-Lagr-opt-proc}.
\end{theorem}
The proof is given in Section~\ref{SF:sec:proof-action-integral-form-RF}. The map~$\mathcal{L}$ is called the \emph{Lagrangian}. We characterize the Lagrangian in two ways, as the Legendre-Fenchel transform of a principal eigenvalue (Eq.~\eqref{SF:eq:variational_Hamiltonian} below), and as an optimization problem (Eq.~\eqref{SF:eq:results-Lagr-opt-proc} below). Both characterizations involve a Hamiltonian that we call \emph{slow Hamiltonian}, and the so-called Donsker-Varadhan rate functional. We first describe these two ingredients, and then give the representations of the Lagrangian. 
\smallskip

For fixed~$z\in F$, consider the process~$\widehat{X}_z^n(t)$ with generator $A_{n,z}^\mathrm{slow}$. Intuitively, this means considering the dynamics when freezing the fast process to~$z$. By our assumptions, the process~$\widehat{X}_z^n(t)$ satisfies pathwise large deviations,
\begin{equation*}
	\mathbb{P}\left[\widehat{X}_z^n\approx \gamma(\cdot)\right]\sim \exp\left\{-n\int_0^\infty \widehat{\mathcal{L}}(\gamma,\partial_t\gamma,z)\,\dd t\right\},\quad n\to\infty,
\end{equation*}
where the \emph{slow Lagrangian}~$\widehat{\mathcal{L}}$ is the Legendre-Fenchel transform of a function that we call the \emph{slow Hamiltonian}~$V_{x,p}(z)$, as 
\begin{equation*}
	\widehat{\mathcal{L}}(x,v,z) = \sup_{p}\left[p\cdot v - V_{x,p}(z)\right].
\end{equation*}
The slow Hamiltonian~$V_{x,p}(z)$ is derived from the slow generator. For~$p\in\mathbb{R}^d$ and a function~$f$ with~$\nabla f(x)=p$, it satifies
\begin{equation}\label{SF:eq:main-result:slow-Hamiltonian}
V_{x,p}(z) = \lim_{n\to\infty} \frac{1}{r_n}e^{-r_nf(x)}A_{n,z}^\mathrm{slow}e^{r_nf(x)}.
\end{equation}
We call~$V_{x,p}(z)$ the slow Hamiltonian since it arises from the slow dynamics when completely decoupling the slow and fast processes, which effectively means to consider the slow dynamics only.
\smallskip

Vice versa, consider the process~$\widehat{Z}_x^n(t)$ with generator~$r_n\cdot A_{x}^\mathrm{fast}$. As before, this intuitively means to consider the slow process frozen to~$x$, and to follow the fast dynamics independently of the slow dynamics. The operator~$A_x^\mathrm{fast}$ arises from our assumptions as the limit of the fast generators~$A_{x,n}^\mathrm{fast}$. Under appropriate ergodicity assumptions, the fast process converges exponentially fast to equilibrium with speed~$r_n$ in the sense that for a distribution~$\nu\in\mathcal{P}(F)$,
\begin{equation*}
	\mathbb{P}\left[\int_0^1 \mathbf{1}_{\widehat{Z}_x^n(t)}(\cdot)\dd t\approx \nu\right]\sim \exp\{-r_n \cdot \mathcal{I}(x,\nu)\},\quad n\to\infty.
\end{equation*} 
The map~$\mathcal{I}(x,\cdot):\mathcal{P}(F)\to[0,\infty]$ is the Donsker-Varadhan functional. In terms of the limiting fast generator, it is given by
\begin{equation}\label{SF:eq:main-results:DV-functional}
\mathcal{I}(x,\nu) = -\inf_{\substack{u \in \mathcal{D}(A_x^\mathrm{fast})\\ u > 0}}\int_F \frac{A_x^\mathrm{fast}u}{u}\,\dd\nu.
\end{equation}
We give more background on this type of convergence in Chapter~4, where we prove convergence to equilibrium for piecewise-deterministic processes.
\smallskip

The Lagrangian of Theorem~\ref{SF:thm:rate-function} admits two representations in terms of the slow Hamiltonian~$V_{x,p}(z)$ and the Donsker-Varadhan functional~$\mathcal{I}(x,\nu)$.
\paragraph{Legendre dual of principal eigenvalue.} 
For~$(x,p)\in E\times\mathbb{R}^d$, let~$\mathcal{H}(x,p)$ be the principal eigenvalue of the operator $V_{x,p}(\cdot)+A^\mathrm{fast}_x$, meaning there exists a strictly positive function~$u:F\to(0,\infty)$ such that~$(V_{x,p}(z)+A^\mathrm{fast}_x)u(z)=\mathcal{H}(x,p) u(z)$. This Hamiltonian~$\mathcal{H}(x,p)$ admits the variational representation
\begin{equation}
\label{SF:eq:variational_Hamiltonian}
\cH(x,p) = \sup_{\nu \in \cP(F)} \left\{\int V_{x,p}(z) \, \nu(\dd z) - \cI(x,\nu) \right\}.
\end{equation} 
The Lagrangian is the Legendre dual
\begin{equation}\label{SF:eq:Lagrangian-is-Legendre-dual}
\mathcal{L}(x,v) = \sup_{p\in\mathbb{R}^d}\ip{p}{v} -\mathcal{H}(x,p).
\end{equation}
\paragraph{Optimization over velocities.}
The Lagrangian~$\mathcal{L}$ defined by~\eqref{SF:eq:Lagrangian-is-Legendre-dual} satisfies
\begin{multline}\label{SF:eq:results-Lagr-opt-proc}
	\cL(x,v) = \inf \left\{ \int_F \widehat{\cL}(x,w(z),z) \, \nu(\dd z) + \cI(x,\nu) \, \middle| \, \nu \in \cP(F), \right. \\
	\left. w : F \rightarrow \bR^d \text{ $\nu$-integrable and } \int_F w(z) \, \nu(\dd z) = v  \right\}.
\end{multline}
We close this section by sketching how~\eqref{SF:eq:results-Lagr-opt-proc} follows from~\eqref{SF:eq:variational_Hamiltonian} and~\eqref{SF:eq:Lagrangian-is-Legendre-dual}. Starting from the latter, we find by exchanging infimum and supremum that
\begin{align*}
\mathcal{L}(x,v) &=\sup_{p} \inf_\nu \left[\ip{p}{v}-\int_F V_{x,p}\,\dd \nu + \mathcal{I}(x,\nu)\right]\\
&= \inf_\nu \left[\sup_p \int_F \left(\ip{p}{w}-V_{x,p}\right)\,\dd \nu + \mathcal{I}(x,\nu)\right],
\end{align*}
for any~$w$ averaging to~$v$.
Passing the supremum inside the integral gives~\eqref{SF:eq:results-Lagr-opt-proc}.
\subsection{Mean-field coupled to fast diffusion}
\label{SF:sec:mean-field-fast-diffusion}
In this section, we provide a large-deviation result for mean-field interacting jump processes coupled to a fast diffusion process. Concretely, we take the simultaneous limit of infinitely many particles and inifinte time-scale separation, and are interested in the large deviations of the empirical density-flux pairs of the mean-field system. For formulating the large-deviation result  by Theorem~\ref{SF:thm:LDP-mean-field} below, we first introduce the processes~$X_t^n$ (Eq.~\eqref{SF:eq:def:density-flux} below) and~$Z_t^n$ (Eq.~\eqref{SF:eq:intro-mean-field:fast-generator} below) independently from one another, and then consider the coupling. We start with describing the mean-field system.
\paragraph{The slow process: mean-field system.}
The mean-field system of~$n$ particles is described by~$n$ weakly-interacting jump processes on a finite state space. That means every jump process~$Y_{n,i}$ takes values in~$\{1,\dots,q\}$, for~$i=1,\dots,n$. We collect the states of all particles in a vector
\begin{equation*}
	Y_n(t) := (Y_{n,1}(t),\dots,Y_{n,n}(t)) \in \{1,\dots,q\}^n.
\end{equation*}
Each jump process can jump over an edge~$(a,b)$; for the set of directed edges in~$\{1,\dots,q\}$, we write $\Gamma = \left\{(a,b) \in \{1,\dots,q\}^2 \, \middle| \, a \neq b \right\}$. We assume that only one particle can jump at a time. The time-evolution of~$Y_n$ is specified by jump rates~$r(a,b)$ attached to each bond~$(a,b)\in \Gamma$. To incorporate the assumption of \emph{weak interactions}, the jump rates are assumped to depend on the configuration of the particles only via their distribution. More specifically, consider the \emph{empirical density} 
\begin{equation*}
	\mu_n(Y_n(t)) := \frac{1}{n}\sum_{i=1}^n\delta_{Y_{n,i}(t)}\in\mathcal{P}(\{1,\dots,q\}).
\end{equation*}
Then the transitions of~$Y_n$ as a jump process on~$\{1,\dots,q\}^n$ are determined by a family of rates~$\{r_n(a,b,\mu):(a,b)\in\Gamma,\mu\in\mathcal{P}(\{1,\dots,q\}\}$. For fixed~$\mu$, the scalar~$r_n(a,b,\mu)\geq 0$ is the rate at which transitions from~$a$ to~$b$ occur when the empirical density is in configuration~$\mu$. Put differently, if the particles are in configuration~$Y_n(t)$, then the jump~$Y_{n,i}(t)\to b$ of the i'th particle occurs at rate
\begin{equation*}
	r_n\left(Y_{n,i}(t),b,\mu_n(Y_n(t)\right).
\end{equation*}
Next to the empirical density, we keep track of the number of jumps that occured over each bond. To that end, let~$t \mapsto W_{n,i}(t) \in \bN^\Gamma$ be the process counting the number of times the i'th particle jumps over each bond,
\begin{equation*}
	W_{n,i}(t)(a,b) := \#\left\{0\leq s \leq t \, \middle| \, \left(X_{n,i}(s-), X_{n,i}(s)\right) = (a,b) \right\}.
\end{equation*}
We regard~$W_{n,i}(t)$ as a random vector taking values in~$\mathbb{N}^\Gamma$. The average fluxes over all bonds are captured by the \emph{empirical flux}~$W_{n}$ defined as
\begin{equation*}
	W_n(t) :=\frac{1}{n}\sum_{i=1}^nW_{n,i}(t).
\end{equation*}
The slow process we are interested in is the pair of empirical density and flux,
\begin{equation}\label{SF:eq:def:density-flux}
	X_t^n := \left(\mu_n(Y_n(t)),W_n(t)\right) \in E:= \cP(\{1,\dots,q\}) \times [0,\infty)^\Gamma,
\end{equation}
which we will refer to as the \emph{density-flux process}.
We write~$x$ for the variables in~$E$, which are pairs~$x=(\mu,w)$ of configurations~$\mu\in\mathcal{P}(\{1,\dots,q\})$ and average fluxes~$w\in[0,\infty)^\Gamma$. We identify the probability measures~$\mathcal{P}(\{1,\dots,q\})$ with the simplex in~$\mathbb{R}^q$, 
\begin{equation*}
	\{\mu\in\mathbb{R}^q\,:\, \sum_{i=1}^q \mu_i = 1,\,\mu_i\geq 0\},
\end{equation*}
equipped with the Euclidean topology inherited from~$\mathbb{R}^d$, so that convergence in the simplex coincides with weak convergence in~$\mathcal{P}(\{1,\dots,q\})$. We also identify the~$n$-atomic measures $P_n:=\{(1/n)\sum_{i=1}^n\delta_{q_i}\,:\,q_i\in\{1,\dots,q\}\}$ with the simplex intersected with~$(1/n)\mathbb{Z}^q$. We sometimes write~$\mu_i=\mu(i)$.
\smallskip

Finally, we describe the generator of~$X_t^n$. If~$\mu_n(Y_n(t))=:\mu$, then the transition of a particle from~$a$ to~$b$ occurs at rate~$r_n(a,b,\mu)$. The number of particles in state~$a$ is $n\cdot \mu(a)$. Hence the rate at which the configuration~$\mu$ transitions to the configuration $\mu + (\delta_b-\delta_a)/n$ is given by~$n\cdot \mu(a)\cdot r_n(a,b,\mu)$. Therefore, the generator~$A_n^\mathrm{slow}:C_b(E)\to C_b(E)$ of the jump process~$X_t^n$ is
\begin{equation}\label{SF:eq:intro-mean-field:slow-generator}
	A_n^\mathrm{slow} f(x) = \sum_{a,b;a\neq b} n\cdot \mu(a)\cdot r_n(a,b,\mu)\left[f(x_{a\to b}^n)-f(x)\right],
\end{equation}
where for a state~$x=(\mu,w)\in E$, we denote by~$x_{a\to b}^n$ the state after the jump. Since after the jump, exactly one particle has changed its state from~$a$ to~$b$,
\begin{equation*}
	x_{a\to b}^n = \left(\mu + \frac{1}{n}(\delta_b-\delta_a), w+\frac{1}{n}\delta_{(a,b)}\right).
\end{equation*}
\paragraph{The fast process: drift-diffusion.} The process $Z_t^n$ is a drift-diffusion process on the flat torus~$F=\mathbb{T}^m$, some~$m\in\mathbb{N}$. While our arguments that concern~$Z_t^n$ also hold true on a closed, smooth, compact, connected manifold, we do not consider this generalization in order to avoid geometric discussions. The generator of~$Z_t^n$ is a second-order uniformly-elliptic differential operator given by
\begin{equation}\label{SF:eq:intro-mean-field:fast-generator}
	A_n^\mathrm{fast} f(z) = \sum_{i=1}^m b_n^i(z)\partial_i f(z) + \sum_{ij=1}^ma_n^{ij}(z)\partial_i\partial_j f(z),
\end{equation} 
where~$a_n(z) = \sigma_n(z)\sigma_n(z)^T$ are symmetric positive-definite matrices and~$b_n(z)$ are vector fields. The domain~$\mathcal{D}(A^\mathrm{fast})$ of~$A_n^\mathrm{fast}$ is independent of~$n$ and is dense in~$C(F)$. In one dimension,~$\mathcal{D}(A^\mathrm{fast})=C^2(F)$, while for any dimension~$m\geq 2$, the domain is larger. On functions~$f\in C^2(F)$ however, the action of the generator is always given by~\eqref{SF:eq:intro-mean-field:fast-generator}. For details on the construction of the process from the operator, we refer to Ikeda's and Watanabe's monograph~\cite[Theorem~IV.6.1]{IkedaWatanabe2014} and the discussion thereafter. 
\paragraph{The coupled slow-fast system.}
We described two processes above. First, the density-flux process~$X_t^n$ given in terms of jump rates~$r_n(a,b,\mu)$, whose generator~$A_n^\mathrm{slow}$ is a pure jump process on a finite subset of~$E$. Secondly, the drift-diffusion process~$Z_t^n$ defined in terms of drifts~$b_n^i(z)$, diffusion-coefficient matrices~$\sigma_n(z)$, and the generator~$A_n^\mathrm{fast}$. In order to obtain a coupled system, we consider coefficients depending on both slow and fast variables:
\begin{enumerate}[label=(\roman*)]
	\item The jump rates are in addition~$z$-dependent,~$r_n=r_n(a,b,\mu,z)$.
	\item The drifts and diffusion-coefficients are in addition~$x$-dependent, meaning $b_n^i=b_n^i(x,z)$ and~$\sigma_n^{ij}=\sigma_n^{ij}(x,z)$.
\end{enumerate}
The pair~$(X_t^n,Z_t^n)$ we want to obtain is an example of coupling a jump process to a drift-diffusion process. The following regularity condition is imposed in order to ensure that we obtain a Feller-continuous process~$(X_t^n,Z_t^n)$ solving the martingale problem~\cite[Theorem~2.1, Section~2.5 and Theorem~2.18]{YinZhu2009}.
\begin{condition}[Regularity]\label{SF:condition:mean-field:reg-coefficients}
For each~$i,j\in 1$,~$n=1,2,\dots$, we have:
\begin{enumerate}[label=(\arabic*)] 
	\item For each~$x\in E$,~$a_n^{ij}(x,\cdot)\in C^2(F)$ and~$b_n^i(x,\cdot)\in C^1(F)$. 
	\item There is a constant~$C>0$ such that $\ip{a_n(x,z)\xi}{\xi}\geq C|\xi|^2$ for all~$\xi\in T_z F$ and for all~$(x,z)\in E\times F$.
	\item For each~$(a,b)\in\Gamma$, the jump rates~$r_n(a,b,\mu,z)$ depend continuously on~$(\mu,z)$, and $r_n(a,b,\mu,\cdot)\in C^1(F)$ for each~$\mu\in E$. \qed
\end{enumerate}
\end{condition}
Accordingly, we consider the operators~$A_{n,z}^\mathrm{slow}$ and~$A_{n,x}^\mathrm{fast}$ by replacing the coeffients in~\eqref{SF:eq:intro-mean-field:slow-generator} and~\eqref{SF:eq:intro-mean-field:fast-generator},
\begin{align}
A_{n,z}^\mathrm{slow} g(x) &:= \sum_{a,b;a\neq b} n\cdot \mu(a)\cdot r_n(a,b,\mu,z)\left[g(x_{a\to b})-g(x)\right],\\
A_{n,x}^\mathrm{fast} h(z) &:= \sum_i b_n^i(x,z)\partial_i h(z) + \sum_{ij}a_n^{ij}(x,z)\partial_i\partial_j h(z).
\end{align}
Furthermore, we let the diffusion process run on the time-scale of order~$n$. The generator~$A_n$ of the couple~$(X_t^n,Z_t^n)$ is
\begin{equation}\label{SF:eq:intro-mean-field:generator-slow-fast}
	A_n f(x,z) := A_{n,z}^\mathrm{slow}f(\cdot,z)(x) + n\cdot  A_{n,x}^\mathrm{fast}f(x,\cdot)(z).
\end{equation}
We obtained a two-component process~$(X^n,Z^n)\in D_{E\times F}[0,\infty)$ with generator~$A_n$. The diffusion process~$Z^n$ is running at a time-scale or order~$n$ faster compared to the density-flux process~$X^n$. Therefore we refer to~$Z^n$ as the fast process and to~$X^n$ as the slow process.
\paragraph{Large deviations of the slow component.} We take the limit~$n\to\infty$ and ask the following questions: does the density-flux process~$X_t^n$ satisfy a large deviation principle in~$D_E[0,\infty)$ under the influence of the fast diffusion process~$Z^n$? How exactly does the fast process affect the large-deviation fluctuations of the particle system? We answer these questions by Theorem~\ref{SF:thm:LDP-mean-field}.
\smallskip

A large-deviation result can only be expected if the jump rates of the particle system and the coefficients of the diffusion process converge as~$n\to\infty$. We work under the following convergence assumptions.
\begin{assumption}[Convergence of rates]\label{SF:mean-field:conv-rates}
	There is a kernel~$r=r(a,b,\mu,z)$ such that for each edge~$(a,b)\in\Gamma$,
	\begin{equation*}
	\lim_{n\to\infty}\sup_{\mu\in \mathcal{P}_n}\sup_{z\in F}\left| r_n(a,b,\mu,z)-r(a,b,\mu,z)\right| = 0.
	\end{equation*}
	There are constants~$0<r_\mathrm{min}\leq r_{\mathrm{max}}<\infty$ such that for all edges~$(a,b)\in\Gamma$ satisfying~$\sup_{\mu,z}r(a,b,\mu,z)>0$, we have
	\begin{equation*}
		r_\mathrm{min}\leq \inf_{\mu,z}r(a,b,\mu,z) \leq \sup_{\mu,z}r(a,b,\mu,z) \leq r_\mathrm{max}.
	\end{equation*}
\end{assumption} 
\begin{assumption}[Convergence of coefficients]\label{SF:mean-field:conv-coeff}
	For each~$i,j$, there are functions~$b^i$ and~$\sigma^{ij}$ on~$E\times F$ such that whenever~$x_n=(\mu_n,w_n)\to (\mu,w)$, then
	\begin{equation*}
		\|b^i(\mu,\cdot)-b_n^i(x_n,\cdot)\|_{F}\to 0\quad\text{and}\quad \|\sigma^{ij}(\mu,\cdot)-\sigma_n^{ij}(x_n,\cdot)\|_{F}\to 0,
	\end{equation*}
	where~$\|g\|_F=\sup_F|g|$. The maps~$\mu\mapsto\sigma^{ij}(\mu,\cdot)$ are continuous as functions from~$\mathcal{P}(\{1,\dots,q\})$ to~$C(F)$ equiped with the uniform norm.\qed
\end{assumption}
\begin{theorem}[Large deviations of the density-flux process] \label{SF:thm:LDP-mean-field}
	Let~$(X^n,Z^n)$ be the Markov process with generator~\eqref{SF:eq:intro-mean-field:generator-slow-fast}.
	Suppose that Assumptions~\ref{SF:mean-field:conv-rates} and~\ref{SF:mean-field:conv-coeff} hold true and that $X^n(0)$ satisfies a large deviation principle with good rate function $J_0:E\to[0,\infty]$ on $E = \cP(\{1,\dots,q\}) \times [0,\infty)^\Gamma$.
	\smallskip
	
	Then $\{X^n\}_{n = 1,2\dots}$ satisfies a large deviation principle on $D_{E}[0,\infty)$ with good rate function $J$ given by
	\begin{equation*}
		J(x) = 
		\begin{cases}
			J_0(\gamma(0)) + \int_0^\infty \cL(\gamma(t),\partial_t \gamma(t)) \dd t & \text{if } \gamma \in \cA\cC([0,\infty);E), \\
			\infty & \text{otherwise},
		\end{cases}
	\end{equation*}
	where the Lagrangian $\cL : E \times \bR^d \rightarrow [0,\infty]$ satisfies the two representations shown below.\qed
\end{theorem}
As in the general large-deviation result (Theorem~\ref{SF:thm:rate-function}), the Lagrangian admits two characterizations. To fix notation, a path~$\gamma:[0,\infty)\to E$ is a time-dependent pair~$\gamma_t=(\mu_t,w_t)\in \mathcal{P}(\{1,\dots,q\})\times[0,\infty)^\Gamma$, where we identify the probability measures with the simplex in~$\mathbb{R}^q$. The set~$E$ is a subset of~$\mathbb{R}^d$ with dimension~$d=q+|\Gamma|$. 
\smallskip

We use the terminology from Section~\ref{SF:sec:result:general-LDP} to formulate the Lagrangian in terms of the following ingredients. The slow Hamiltonian; for~$(x,p)\in E\times \mathbb{R}^d$,
\begin{equation*}
	V_{x,p}(z) = \sum_{ab}\mu_ar(a,b,\mu,z)\left[e^{p_b-p_a + p_{ab}}-1\right].
\end{equation*}
The Donsker-Varadhan functional; for~$x\in E$,
\begin{equation*}
	\mathcal{I}(\mu,\pi) = -\inf_{\substack{u>0\\u\in C^2(F)}} \int_F \frac{A_\mu^\mathrm{fast}u}{u}\,\dd\pi,
\end{equation*}
where~$A_\mu^\mathrm{fast}u(z) := \sum_i b^n(\mu,z)\partial_i u(z) + \sum_{ij} a^{ij}(\mu,z)\partial_i u(z)\partial_j u(z)$. The relative entropy function $S(a|b)$,
\begin{equation*}
	S(a \, | \, b) := \begin{cases}
		b & \text{if } a = 0, \\
		a \log \left(a/b\right) - (a-b) & \text{if } a \neq 0, b \neq 0, \\
		+\infty & \text{if } a \neq 0, b = 0.
	\end{cases}
\end{equation*}
\paragraph{Dual of principal eigenvalue.} For~$(x,v)\in E\times\mathbb{R}^d$, the Lagrangian~$\mathcal{L}(x,v)$ is the Legendre dual~$\mathcal{L}(x,v)=\sup_{p\in\mathbb{R}^d}\ip{p}{v}-\mathcal{H}(x,p)$, where the Hamiltonian is the principal eigenvalue of ~$(V_{x,p}+A_x^\mathrm{fast})$. The Hamiltonian satisfies the variational formula~\eqref{SF:eq:variational_Hamiltonian} from Section~\ref{SF:sec:result:general-LDP}.\qed
\paragraph{Optimizing over velocities.}
For a path~$\gamma:[0,\infty)\to E$,~$\gamma=(\mu,w)$, the Lagrangian~$\mathcal{L}$ is finite only if $\partial_t\mu_a = \sum_b \partial_t (w_{ba}-w_{ab})$. If this is the case, then
\begin{equation*}
	\cL\left(\gamma,\partial_t\gamma\right) = \inf_{\pi \in \cP(F)} \inf_{u \in \Phi(\partial_tw,\pi)}  \left\{ \sum_{(a,b) \in \Gamma} \int_F S(u_{ab}(z)\,|\, \mu_a r(a,b,\mu,z)) \pi(\dd z) + \cI(\mu,\pi) \right\},
\end{equation*}
where $\Phi(\partial_tw,\pi)$ is the set of measurable functions $u_{ab}(z)$ for $z \in F$ and $(a,b) \in \Gamma$ such that~$\int u_{ab}(z) \pi(\dd z) = \partial_t w_{ab}$.\qed
\subsection{Averaging principles}
\label{SF:sec:averaging-principles}
We discuss the consequences of the pathwise large-deviation theorems.
\paragraph{Mean-field system.}
We consider the coupled system~$(X^n,Z^n)$ introduced in Section~\ref{SF:sec:mean-field-fast-diffusion}. The pair~$(X^n,Z^n)$ corresponds to the mean-field interaction particles coupled to fast diffusions. The density-flux pair~$X^n=(\rho^n,w^n)$ of the particle-system is a stochastic process in~$D_E[0,\infty)$, where the state space is given by~$E=\mathcal{P}(\{1,\dots,q\})\times[0,\infty)^\Gamma$. 
\smallskip

If the particles are not coupled to the fast diffusion process, then in the limit of large numbers, the evolution of the particle density~$\rho^n$ is characterzied as the solution to a nonlinear ODE, which may be regarded as the finite-dimensional analogue of the McKean-Vlasov equation. We formulate the result in terms of freezing the diffusion process to a value~$z\in F$ and the transition-rate matrix~$R(\rho,z)$ of a jump process with rates~$r(a,b,\rho,z)$, that is
\begin{equation*}
R_{ab}(\rho,z)=r(a,b,\rho,z)\;\; (a\neq b)\qquad\text{and}\qquad R_{aa}(\rho,z)= -\sum_{b\neq a} r(a,b,\rho,z).
\end{equation*}
We write~$\rho R$ for the vector with components~$(\rho R)_a=\sum_{b}\rho_bR_{ba}$.
\begin{proposition}[Law of Large Number limit of mean-field interacting particles]
	Let~$X^n=(\rho^n,w^n)$ be the density-flux process from~$\eqref{SF:eq:def:density-flux}$ with jump rates given by~$r_n(a,b,\cdot,z)$ for some fized~$z\in F$. If the initial density~$\rho^n(0)$ converges in probability to a distribution~$\mu\in\mathcal{P}(\{1,\dots,q\})$, then~$\rho^n$ converges uniformly on compact time intervals to a solution of
	\begin{equation*}
	\partial_t\rho = \rho R(\rho,z),\quad \rho(0)=\mu,
	\end{equation*}
	where~$R(\rho,z)$ is the transition-rate matrix of a jump process with rates~$r(a,b,\rho,z)$.
\end{proposition}
Budhiraja, Dupuis, Fischer and Ramanan proof of this statement~\cite[Theorem~2.2]{BudhirajaDupuisFischerRamanan2015} based on a classical convergence Theorem by Kurtz~\cite{Kurtz1970}. Under a Lipschitz condition on the limiting rates, the limit is unique.
\smallskip

Under the influence of the fast diffusion, we prove the limiting evolution to be altered according to the averaging principle.
\begin{theorem}[Averaging Principle]\label{SF:thm:mean-field-averaging}
	Let~$(X^n,Z^n)$ be the slow-fast system from Theorem~\ref{SF:thm:LDP-mean-field}, with~$X^n=(\rho^n,w^n)$ and initial condition~$\rho^n(0)\to \mu\in\mathcal{P}(\{1,\dots,q\}$ weakly as~$n\to\infty$. Let~$\pi_\nu\in\mathcal{P}(F)$ be the unique measure satisfying~$\mathcal{I}(\nu,\pi_\nu)=0$.
	
	Then~$\rho^n$ converges a.s. with respect to the Skorohod topology to a solution of
	\begin{equation}\label{SF:averaging-principles:mean-field:limit-evol}
	\partial_t\rho = \rho \widehat{R}(\rho),\quad \rho(0) = \mu.
	\end{equation}
	The transition-rate matrix~$\widehat{R}(\rho)$ is an averaged matrix, componentwise given by
	\begin{equation}\label{SF:eq:averaging-principle:mean-field:R-hat}
	\widehat{R}_{ab}(\rho) = \int_F R_{ab}(\rho,z)\,\pi_\rho(\dd z).
	\end{equation}
\end{theorem}
If~$\widehat{R}$ is Lischitz continuous, then the averaged McKean-Vlasov equation~\eqref{SF:averaging-principles:mean-field:limit-evol} has a unique solution. In that case, the minimizer of the rate function is unique, and the pathwise large deviation principle of Theorem~\ref{SF:thm:LDP-mean-field} implies 
that~$\rho^n$ converges to the solution.
In general, if the rate function has multiple minimizers, then the large deviation principle does not contain enough information to determine the limit.
\begin{proof}[Proof of Theorem~\ref{SF:thm:mean-field-averaging}] 
	We show that any density~$\rho$ of the minimizer~$x=(\rho,w)$ of the rate function~$J$ solves~\eqref{SF:averaging-principles:mean-field:limit-evol}. If~$J(\rho,w)=0$, then~$\mathcal{L}(x(t),\partial_t x(t))=0$ for a.e.~$t>0$, where the Lagrangian~$\mathcal{L}$ is given by
	\begin{equation*}
	\cL\left(x,\partial_tx\right) = \inf_{\pi \in \cP(F)} \inf_{u \in \Phi(\partial_tw,\pi)}  \left\{ \sum_{(a,b) \in \Gamma} \int_F S(u_{ab}(z)\,|\, \rho_a r(a,b,\rho,z)) \pi(\dd z) + \cI(\rho,\pi) \right\},
	\end{equation*}
	and by finiteness of the Lagrangian,
	\begin{equation}\label{SF:eq:proof:averaging-principle:mean-field:rho}
	\partial_t\rho_a = \sum_b \partial_t (w_{ba}-w_{ab}).
	\end{equation}
	In the formula for the Lagrangian, all terms inside the infimum are non-negative. Hence for any~$u\in\Phi(\partial_tw,\pi)$ and~$\pi\in\mathcal{P}(F)$, the expression is only zero if we have~$\mathcal{I}(\rho,\pi)=0$. Therefore~$\pi=\pi_\rho$, and
	\begin{equation*}
	0=\mathcal{L}(x,\partial_tx) = \inf_{u\in\Phi(\partial_tw,\pi_\rho)}\sum_{ab} \int_F S\left(u_{ab}(z)\,|\rho_a r(a,b,\rho,z)\right)\,\pi_\rho(\dd z).
	\end{equation*}
	Since~$S(r,s)=0$ if and only if~$r=s$, any optimizer~$u_{ab}(\cdot)$ satisfies
	\begin{equation}\label{SF:eq:proof:averaging-principle:mean-field:uab}
	u_{ab}(z)=\rho_ar(a,b,\rho,z)\qquad \pi_\rho\,\mathrm{a}.\mathrm{e},
	\end{equation}
	and by definition of the set~$\Phi(\partial_tw,\pi_\rho)$,
	\begin{equation}\label{SF:eq:proof:averaging-principle:mean-field:wab}
	\partial_tw_{ab}=\int_F u_{ab}(z)\,\pi_\rho(\dd z).
	\end{equation}
	Combining these equalities, we find
	\begin{align*}
	\partial_t \rho_a 
	&\overset{\eqref{SF:eq:proof:averaging-principle:mean-field:rho}}{=} \sum_{b\neq a} \partial_t (w_{ba}-w_{ab})\\
	&\overset{\eqref{SF:eq:proof:averaging-principle:mean-field:wab}}{=} \sum_{b\neq a} \int_F (u_{ba}-u_{ab})\,\pi_\rho(\dd z)\\
	&\overset{\eqref{SF:eq:proof:averaging-principle:mean-field:uab}}{=} \sum_{b\neq a} \left[\rho_b\int_F r(b,a,\rho,z)\pi_\rho(\dd z) -\rho_a \int_F r(a,b,\rho,z)\pi_\rho(\dd z)\right] \\
	&\overset{\eqref{SF:eq:averaging-principle:mean-field:R-hat}}{=}
	\sum_{b\neq a}\rho_b \widehat{R}_{ba}(\rho) + \rho_a\widehat{R}_{aa}(\rho)\myeqdef (\rho\widehat{R}(\rho))_a,
	\end{align*}
	which finishes the proof.
\end{proof}
The conclusion of Theorem~\ref{SF:thm:mean-field-averaging} remains true when replacing the fast diffusion process by a fast jump process on a finite state space~$\{1,\dots,m\}$. Then the limiting averaged matrix is simply obtained from the equilibrium measure~$\pi=(\pi_1,\dots,\pi_m)$ of the fast jump process as
\begin{equation*}
\widehat{R}(\rho) = \sum_{i=1}^m R(\rho,i)\pi_i.
\end{equation*}
\section{Assumptions of general large-deviation theorem}
\label{SF:sec:assumptions:general-LDP-theorem}
Here we collect the assumptions underlying Theorems~\ref{SF:thm:LDP_general} and~\ref{SF:thm:rate-function}. We pose all assumptions in terms of the slow-fast generators~$A_n$ given by~\eqref{SF:eq:setting:generator-slow-fast-system}, that is
\begin{equation*}
A_n f(y,z) = A^\mathrm{slow}_{n,z} f(\cdot,z) (y) +
r_n \cdot A^\mathrm{fast}_{n,y} f(y,\cdot)(z).
\end{equation*}
The assumptions cluster in three groups, where each group corresponds to one step in the large-deviation proof:
\begin{enumerate}[label=(\roman*)]
	\item Convergence of nonlinear generators.
	\item Comparison principle of a limiting Hamilton-Jacobi equation.
	\item Action-integral form of the rate function.
\end{enumerate}
We make these steps precise when explaining the strategy of proof in Section~\ref{SF:sec:proof-general-LDP-Theorem:strategy-of-proof}. Below, Assumptions~\ref{SF:assumption:convergence-slow-nonlinear-generators} and~\ref{SF:assumption:convergence-fast-nonlinear-generators} correspond to the convergence of nonlinear generators, Assumptions~\ref{SF:assumption:principal-eigenvalue-problem},~\ref{SF:assumption:regularity-V} and~\ref{SF:assumption:regularity-I} are used to prove the comparison principle, and finally Assumption~\ref{SF:assumption:Hamiltonian_vector_field} is made in order to obtain the action-integral form of the rate function.
\smallskip

Recall the setting from Condition~\ref{SF:condition:compact_setting:state-spaces}: the slow-fast process~$(X^n,Z^n$) takes values in the product space~$E_n\times F$. The spaces~$E_n$ are continuously embedded into a Polish space~$E\subseteq{R}^d$ with a map~$\eta_n:E_n\to E$, and the space~$F$ are compact Polish spaces.
\smallskip

We first state the two assumptions concerning the convergence of non-linear transforms~$H_n$ of the generator~$A_n$, defined by~$H_nf=(1/n)e^{-nf}A_ne^{nf}$. We assume the slow and fast parts to converge independently.
\begin{assumption}[Convergence of slow non-linear generators]	\label{SF:assumption:convergence-slow-nonlinear-generators}
	Let $D_0$ be a linear space $C_c^\infty(E) \subseteq D_0 \subseteq C_b^1(E)$ satisfying the following:
	\begin{enumerate}[label=(\roman*)]
		\item For any $n$, $z \in F$ and $f \in D_0$, $z \in F$ we have $e^{r_n f} \in \cD(A_{n,z}^{slow})$ and
		\begin{equation*}
		\sup_{n}\sup_{x \in E_n,z \in F} \left|\frac{1}{r_n} e^{-r_nf(x,z)} A_{n,z}^\mathrm{slow} e^{r_nf(\cdot,z)(x)} \right| < \infty;
		\end{equation*}
		\item There exist continuous functions $V_{x,p} : F \to \mathbb{R}$, where $x \in E$  and $p \in \mathbb{R}^d$, such that for any $f \in D_0$ and all compact sets $K \subseteq E$
		\begin{align*}
		\sup_{ x \in \eta_n^{-1}(K), z \in F} \left|
		\frac{1}{r_n} e^{-r_nf(x)} \left(A_{n,z}^\mathrm{slow} e^{r_nf(\cdot)}\right)(x) - V_{\eta_n(x),\nabla f(\eta_n(x))}(z)
		\right| \to 0.
		\tag*\qed
		\end{align*}
	\end{enumerate}
\end{assumption}
We refer to the function~$V_{x,p}(z)$ as the \emph{slow Hamiltonian}. 
\begin{assumption}[Convergence of fast non-linear generators]
	\label{SF:assumption:convergence-fast-nonlinear-generators}
	For every~$x\in E$, there exists an operator~$A_x^\mathrm{fast}$ with the following properties. For any~$g\in \mathcal{D}(A^\mathrm{fast})$,~$x_n \in E_n$ and~$g_n \in\mathcal{D}(A^\mathrm{fast})$ such that~$\eta_n(x_n)\to x$ and~$g_n\to g$ uniformly on~$F$, we have
	\begin{equation*}
	\|A^\mathrm{fast}_xg-A^\mathrm{fast}_{n,\eta_n(x_n)}g_n\|_F\to 0.
	\end{equation*}
	Furthermore, for any~$\phi\in \mathcal{D}(A^\mathrm{fast})$,
	\begin{align*}
	\sup_{n}\sup_{y\in E_n}\|e^{-\phi(z)} \left(A^\mathrm{fast}_{n,y} e^{-\phi}\right)(z)\|_F < \infty.
	\tag*\qed
	\end{align*}
\end{assumption}
With the above two convergence assumptions, we will obtain a limit operator~$H$ defined in terms of a graph~$H\subseteq C_b(E)\times C_b(E\times F)$. The precise definition of~$H$ is given in Definition~\ref{SF:def:multi-valued-limit:general}. The next three assumptions are imposed in order to prove the comparison principle of the Hamilton-Jacobi equation with this Hamiltonian~$H$.
\smallskip

For the first assumption, consider the Hamiltonian~$\mathcal{H}(x,p)$ defined by~\eqref{SF:eq:variational_Hamiltonian}; given the slow Hamiltonian~$V_{x,p}$ and the limit generator~$A^\mathrm{fast}$,
\begin{equation}\label{SF:eq:sec-assumption:H}
\cH(x,p) = \sup_{\pi \in \cP(F)} \left\{\int V_{x,p}(z) \, \pi(\dd z) - \cI(x,\pi) \right\},
\end{equation}
where~$\mathcal{I}$ is the Donsker-Varadhan rate functional~\eqref{SF:eq:main-results:DV-functional}.
\begin{assumption}[Approximative solution to a principal-eigenvalue problem] 
	\label{SF:assumption:principal-eigenvalue-problem}
	Let~$(x,p)\in E\times\mathbb{R}^d$. The limit~$V_{x,p}$ from Assumption~\ref{SF:assumption:convergence-slow-nonlinear-generators} and the operators~$A^\mathrm{fast}_x$ from Assumption~\ref{SF:assumption:convergence-fast-nonlinear-generators} satisfy the following.
	\begin{enumerate}[label=($\cE$\arabic*)]
		\item \label{SF:item:assumption:PI:domain} For any $\Phi$ such that $e^\Phi\in\cD(A^\mathrm{fast})$, we have $e^{(1-\varepsilon)\Phi} \in \cD(A^\mathrm{fast})$ for any~$0<\varepsilon<1$.
		\item \label{SF:item:assumption:PI:solvePI} For every $\delta > 0$ there exists a strictly positive function $u_\delta \in \cD(A^\mathrm{fast})$ on~$F$ satisfying
		\begin{align*}
		\sup_{z \in F} \left| 	\left( V_{x,p}(z) + A^\mathrm{fast}_x \right) u_\delta(z) - \cH(x,p) u_\delta(z) \right| \leq \delta.
		\tag*\qed
		\end{align*}
	\end{enumerate}
\end{assumption}
The second part of Assumption~\ref{SF:assumption:principal-eigenvalue-problem} is satisfied if the principal-eigenvalue problem for the operator~$V_{x,p}+A_x^\mathrm{fast}$ is well-posed, where~$V_{x,p}$ acts via multiplication. By principal-eigenvalue problem, we mean the existence of a strictly positive function~$u$ in the domain of~$A_x^\mathrm{fast}$ and an eigenvalue~$\lambda\in\mathbb{R}$ such that $(V_{x,p}+A_x^\mathrm{fast})u=\lambda u$ holds pointwise on~$F$. If the principal-eigenvalue problem is well-posed, then~$\lambda=\mathcal{H}(x,p)$ by a result of Donsker and Varadhan~\cite{DonskerVaradhan75}.
\smallskip

We impose the next two assumptions on~$V$ and~$\mathcal{I}$ in order to verify the comparison principle for Hamilton-Jacobi equations involving the above Hamiltonian~$\mathcal{H}(x,p)$. The assumptions are derived from~\cite{KraaijSchlottke2019} (Chapter~6), where we prove the comparison principle for Hamiltonians of the type~\eqref{SF:eq:sec-assumption:H}.
\begin{assumption}[Regularity of the slow Hamiltonian]\label{SF:assumption:regularity-V}
	The slow Hamiltonian from Assumption~\ref{SF:assumption:convergence-slow-nonlinear-generators}, that is the map~$V:E\times\mathbb{R}^d\times F\to\mathbb{R}$, satisfies:
	\begin{enumerate}[label=($V$\arabic*)]
		\item \label{SF:item:assumption:slow_regularity:continuity} For every $(x,p)$ we have $V_{x,p} \in C(F)$, and the map $(x,p) \mapsto V_{x,p}$ is continuous on $C(F)$ for the supremum norm.
		\item \label{SF:item:assumption:slow_regularity:convexity} For any $x \in E$ and $z \in F$, the map $p \mapsto V_{x,p}(z)$ is convex. Furthermore, we have $V_{x,0}(z) = 0$ for all $x,z$.
		\item \label{SF:item:assumption:slow_regularity:compact_containment} There exists a continuous containment function $\Upsilon : E \to [0,\infty)$ in the sense of Definition~\ref{def:results:compact-containment}.
		\item \label{SF:item:assumption:slow_regularity:continuity_estimate} The function $\Lambda(x,p,\nu) := \int V_{x,p}(z) \, \nu(\dd z)$ on~$E\times\mathbb{R}^d\times\mathcal{P}(F)$ satisfies the continuity estimate~\cite[Definition~4.14]{KraaijSchlottke2019}. A definition of the continuity estimate is given in~\ref{def:results:continuity_estimate}.
		\item \label{SF:item:assumption:slow_regularity:controlled_growth} 
		For every compact set $K \subseteq E$, there exist constants $M, C_1, C_2 \geq 0$  such that for all $x \in K$, $p \in \mathbb{R}^d$ and all $z_1,z_2\in F$,
		\begin{align*}
		V_{x,p}(z_1) \leq  \max\left\{M,C_1 V_{x,p}(z_2) + C_2\right\}.
		\tag*\qed
		\end{align*}
	\end{enumerate} 
\end{assumption}
The conditions~\ref{SF:item:assumption:slow_regularity:continuity},~\ref{SF:item:assumption:slow_regularity:convexity} follow from the convergence in Assumption~\ref{SF:assumption:convergence-slow-nonlinear-generators}. We state them nevertheless to clearify the connection to~\cite{KraaijSchlottke2019}. 
\begin{assumption}[Regularity of the Donsker-Varadhan functional]\label{SF:assumption:regularity-I}
	The functional $\mathcal{I}:E\times\cP(F) \to [0,\infty]$ from~\eqref{SF:eq:main-results:DV-functional} satisfies the following.
	\begin{enumerate}[label=($\mathcal{I}$\arabic*)]
		\item \label{SF:item:assumption:I:lsc} The map $(x,\nu) \mapsto \mathcal{I}(x,\nu)$ is lower semi-continuous on $E \times \cP(F)$.
		\item \label{SF:item:assumption:I:zero-measure} For any $x\in E$, there exists a point $\nu_x\in\cP(F)$ such that $\mathcal{I}(x,\nu_x) = 0$. 
		\item \label{SF:item:assumption:I:compact-sublevelsets} For any $x \in E$, compact set $K \subseteq E$ and $C \geq 0$ the set $\left\{\nu \in \cP(F) \, \middle| \cI(x,\nu) \leq C\right\}$ is compact and $\cup_{x\in K}\left\{\nu \in \cP(F) \, \middle| \, \mathcal{I}(x,\nu) \leq C\right\}$ is relatively compact. 
		\item \label{SF:item:assumption:I:finiteness} For any converging sequence $x_n \to x$ in $E$ and sequence $\nu_n \in \cP(F)$, if there is an $M > 0$ such that $\mathcal{I}(x_n,\nu_n) \leq M < \infty$ for all $n \in \mathbb{N}$, then there exists a neighborhood $U_x$ of $x$ and a constant $M' > 0$ such that for any $y \in U_x$ and $n \in \mathbb{N}$,
		\begin{equation*}
		\mathcal{I}(y,\nu_n) \leq M' < \infty.
		\end{equation*}
		\item \label{SF:item:assumption:I:equi-cont} For every compact set $K \subseteq E$ and each $M \geq 0$ the collection of functions $\{\cI(\cdot,\nu)\}_{\nu \in \cP(F)_M}$ with
		\begin{equation*}
		\cP(F)_{M} := \left\{\nu \in \cP(F) \, \middle| \, \forall \, x \in K: \, \mathcal{I}(x,\nu) \leq M \right\}
		\end{equation*}
		is equicontinuous. That is: for all $\varepsilon > 0$, there is a $\delta > 0$ such that for all $\nu \in \cP(F)_M$ and $x,y \in K$ satisfying $d(x,y) \leq \delta$, we have the estimate~$|\mathcal{I}(x,\nu) - \mathcal{I}(y,\nu)| \leq \varepsilon$.\qed
	\end{enumerate}
\end{assumption}
Condition~\ref{SF:item:assumption:I:lsc} follows if the map~$x\mapsto A_x^\mathrm{fast}\phi$ is continuous as a function from~$E$ to~$C(F)$ equiped with the supremum norm. Conditions~\ref{SF:item:assumption:I:zero-measure} and~\ref{SF:item:assumption:I:compact-sublevelsets} are always satisfied by the compactness assumption on~$F$. Again, we state these conditions to make the connection to~\cite{KraaijSchlottke2019} as clear as possible.
\smallskip

Assumptions~\ref{SF:assumption:convergence-slow-nonlinear-generators},~\ref{SF:assumption:convergence-fast-nonlinear-generators},~\ref{SF:assumption:principal-eigenvalue-problem},~\ref{SF:assumption:regularity-V} and~\ref{SF:assumption:regularity-I} suffice for the proof of Theorem~\ref{SF:thm:LDP_general}, which establishes pathwise large deviations. We need one additional assumption to prove the action-integral representation of the rate function. To that end, we denote for a convex function $\Phi : \bR^d \rightarrow (-\infty,\infty]$ its subdifferential by
\begin{equation*}
\partial_p \Phi(p_0)
:= \left\{
\xi \in \mathbb{R}^d \,:\, \Phi(p) \geq \Phi(p_0) + \xi \cdot (p-p_0) \quad (\forall p \in \mathbb{R}^d)
\right\}.
\end{equation*}
The Bouligand tangent cone to $E$ in $\bR^d$ at $x$ is
\begin{equation*}
T_E(x) := \left\{z \in \bR^d \, \middle| \, \liminf_{\lambda \downarrow 0} \frac{d(y + \lambda z, E)}{\lambda} = 0\right\}.
\end{equation*}
\begin{assumption}[] \label{SF:assumption:Hamiltonian_vector_field}
	The slow Hamiltonian $V : E \times \bR^d \times F \rightarrow \bR$ from Assumption~\ref{SF:assumption:convergence-slow-nonlinear-generators} satisfies $\partial_p V_{x,p}(z) \subseteq T_E(x)$ for all $p$, $x$ and $z$.\qed
\end{assumption}
In~\cite{KraaijSchlottke2019}, this assumption is made on the full Hamiltonian~$\mathcal{H}(x,p)$ instead of the slow Hamiltonian~$V_{x,p}(z)$. We will show this property to bootstrap from the slow to the full Hamiltonian.
\section{Proof of large deviations of the slow process}
\label{SF:sec:proof-general-LDP-Theorem}
\subsection{Strategy of the proof}
\label{SF:sec:proof-general-LDP-Theorem:strategy-of-proof}
We outline the large-deviation proof for the slow component of a slow-fast system from Theorem~\ref{SF:thm:LDP_general}. For the generator~$A_n$ of a slow-fast system~$(Y^n,Z^n)$, define the operator~$H_n$ on~$\mathcal{D}(H_n) := \{f\in C_b(E_n\times F)\,:\,e^{r_nf}\in \mathcal{D}(A_n)\}$ by
\begin{equation}\label{SF:eq:nonlinear-generator-slow-fast}
H_n f(y,z):= \frac{1}{r_n}e^{-r_nf(y,z)} \left(A_n e^{r_nf(\cdot)}\right)(y,z)
\end{equation}
We call this operator the \emph{nonlinear generator} of~$(X^n,Z^n)$. To prove large deviations, we exploit the semigroup-convergence method built by Jin Feng and Thomas Kurtz~\cite{FengKurtz2006}. In a nutshell, the large-deviation proof boils down to two steps:
\begin{enumerate}[label=(\roman*)]
	\item Convergence of nonlinear generators to a limit operator.
	\item Verifying the comparison principle for the limit operator.
\end{enumerate}
The definition of the comparison principle is given in~\ref{definition:viscosity_solutions}. Here, we first give precise version of the above steps by Theorem~\ref{SF:thm:Richard-abstract-LDP} below, which is Kraaij's result from~\cite{Kraaij2019ExpResolv} taylored to our setting. After that, we give the proof of Theorem~\ref{SF:thm:LDP_general}. We use the following convergence concepts.
\begin{definition}[LIM-convergence]
	\label{SF:def:LIM-convergence}
	Let $f_n \in C_b(E_n \times F)$ and $f \in C_b(E \times F)$. We say that $\LIM f_n = f$ if 
	\begin{enumerate}[label=(\roman*)]
		\item $\sup_n \vn{f_n} < \infty$,
		\item for all compact sets $K \subseteq E$,
		\begin{align*}
		\lim_{n \rightarrow \infty} \sup_{(y,z) \in \eta_n^{-1}(K) \times F} \left|f_n(y,z) - f(\eta_n(y),z) \right| = 0.
		\tag*\qed
		\end{align*}
	\end{enumerate}
\end{definition}
\begin{definition}[Extended-LIM]\label{SF:def:extended-LIM}
	Let $B_n \subseteq C_b(E_n \times F) \times C_b(E_n\times F)$. The set $\mathrm{ex}-\LIM B_n$ is defined as
	\begin{multline*}
	\mathrm{ex}-\LIM B_n \\
	:= \left\{(f,g) \in C_b(E \times F)^2 \, \middle| \, \exists \, (f_n,g_n) \in B_n: \, \LIM f_n = f, \LIM g_n = g \right\}.
	\tag*\qed
	\end{multline*}
\end{definition}
\begin{definition}[Exponential compact containment condition]
	\label{SF:def:exp-comp-containment}
	Consider the context of Conditions~\ref{SF:condition:compact_setting:state-spaces} and~\ref{SF:condition:compact_setting:well-posedness-martingale-problem}. The sequence of processes $(Y^n,Z^n)$ satisfies the \emph{exponential compact containment condition} at speed $r_n$ if for each compact set $K \subseteq E$, $T >0$ and $a > 0$ there is a compact set $\widehat{K} = \widehat{K}(K,T,a) \subseteq E$ such that 
	\begin{equation*}
	\limsup_{n \rightarrow \infty} \sup_{(y,z) \in \eta_n^{-1}(K) \times F} \frac{1}{r_n} \log P_{y,z}\left[Y_n(t) \notin \eta_n^{-1}(\widehat{K}) \text{ for some } t \in [0,T] \right] \leq - a.
	\end{equation*}	
	\qed
\end{definition}
The following simplified version of~\cite[Theorem~7.10]{Kraaij2019ExpResolv} is sufficient for our purposes.
\begin{theorem}[{Adaptation of~\cite[Theorem~7.10]{Kraaij2019ExpResolv} to our context}] \label{SF:thm:Richard-abstract-LDP}
	Consider a sequence of slow-fast processes~$(Y^n,Z^n)$ in the setting of Conditions~\ref{SF:condition:compact_setting:state-spaces} and~\ref{SF:condition:compact_setting:well-posedness-martingale-problem}, and let $X_n := \eta_n(Y_n)$. Suppose the following conditions hold true:
	\begin{enumerate}[label=(\roman*)]
		\item The exponential compact containment condition (Definition~\ref{SF:def:exp-comp-containment}) is satisfied.
		\item \label{item:LDP_abstract_convergenceHn} There is an operator $H \subseteq C_b(E) \times C_b(E \times F)$ such that $H \subseteq \mathrm{ex}-\LIM H_n$.
		\item \label{item:LDP_abstract_comparison} For all $h \in C_b(E)$ and $\lambda > 0$ the comparison principle holds for $f - \lambda Hf = h$.
	\end{enumerate}
	Suppose furthermore the large deviation principle to hold for $X_n(0) = \eta_n(Y_n(0))$ with speed $r_n$ and good rate function $J_0$.
	\smallskip
	
	Then the processes $X_n = \eta_n(Y_n)$ satisfy a large deviation principle on $D_E[0,\infty)$ with speed $r_n$ and a good rate function~$J$ given by~\eqref{SF:eq:Richard-general-LDP-thm-rate-function}.\qed
\end{theorem}
The rate function~$J$ is implicitly characterized by means of a nonlinear semigroup~$V(t):C_b(E)\to C_b(E)$, as
\begin{equation} \label{SF:eq:Richard-general-LDP-thm-rate-function}
J(\gamma) = J_0(\gamma(0)) + \sup_{k \geq 1} \sup_{\substack{0 = t_0 < t_1 < \dots, t_k \\ t_i \in \Delta_\gamma^c}} \sum_{i=1}^{k} J_{t_i - t_{i-1}}(\gamma(t_i) \, | \, \gamma(t_{i-1})),
\end{equation}
where $\Delta_\gamma^c$ is the set of continuity points of $\gamma$ and the conditional rate functions~$J_t$ are given by
\begin{equation*}
J_t(y \, | \, x) = \sup_{f \in C_b(E)} \left\{f(y) - V(t)f(x) \right\}.
\end{equation*}
The semigroup~$V(t)$ is the limit of nonlinear semigroups~$V_n(t)$ of~$(Y^n,Z^n)$. We do not need the precise statement about the convergence~$V_n\to V$ here. For details, we refer to~\cite{Kraaij2019ExpResolv}.
\subsection{Proof of Theorem~\ref{SF:thm:LDP_general}}
The proof of Theorem~\ref{SF:thm:LDP_general} hinges on the verification of the conditions of Theorem~\ref{SF:thm:Richard-abstract-LDP} above. In the setting of Theorem~\ref{SF:thm:LDP_general}, we consider a sequence of slow-fast systems~$(Y^n,Z^n)$ satisfying Condition~\ref{SF:condition:compact_setting:well-posedness-martingale-problem} (the well-posedness condition on the martingale problem). Their generators~$A_n$ are given by~\eqref{SF:eq:setting:generator-slow-fast-system}, that is
\begin{equation*}
A_n f(y,z) = A^\mathrm{slow}_{n,z} f(\cdot,z) (y) +
r_n \cdot A^\mathrm{fast}_{n,y} f(y,\cdot)(z),
\end{equation*}
Applying the program outlined by Theorem~\ref{SF:thm:Richard-abstract-LDP} requires to establish a limit operator~$H$ of the nonlinear generators~$H_n$ defined above in~\eqref{SF:eq:nonlinear-generator-slow-fast}. The following operator~$H$ defined in terms of a graph~$H\subseteq C_b(E)\times C_b(E\times F)$ serves us as the limit.
\begin{definition}[Multi-valued limit operator]\label{SF:def:multi-valued-limit:general}
	For $f \in C_{cc}^\infty(E)$, $x \in E$ and a function $\phi\in C_b(F)$ such that $e^\phi \in \cD(A^\mathrm{fast}_x)$, set 
	\begin{equation*}
		H_{f,\phi}(x,z) := V_{x,\nabla f(x)}(z) + e^{-\phi(z)} A^\mathrm{fast}_x e^{\phi}(z),
	\end{equation*} 
	and let~$H$ be the graph
	\begin{equation*}
		H := \left\{(f,H_{f,\phi}) \, \middle| \, f \in C_{cc}^\infty(E), \phi: e^{\phi} \in \cD(A^\mathrm{fast}) \right\}.
	\end{equation*}
	The operator $H \subseteq C_b(E) \times C_b(E \times F)$ with $\cD(H) = C_{cc}^\infty(E)$ is multi-valued. \qed
\end{definition}
We prove the following three Lemma's in the subsequent sections under the Assumptions of Theorem~\ref{SF:thm:LDP_general}.
\begin{lemma}[Exponential compact containment]\label{SF:lemma:proof-general-LDP:exp-comp-cont}
	The sequence of slow-fast systems~$(Y^n,Z^n)$ satisfies the exponential compact containment condition.
\end{lemma}
\begin{lemma}[Convergence of nonlinear generators]\label{SF:lemma:proof-general-LDP:conv-to-H}
	Let~$H$ be the operator from Definition~\ref{SF:def:multi-valued-limit:general}. The nonlinear generators~$H_n$ from~\eqref{SF:eq:nonlinear-generator-slow-fast} satisfy~$H\subseteq\mathrm{ex}-\LIM H_n$ in the sense of Definition~\ref{SF:def:extended-LIM}.
\end{lemma}
\begin{lemma}[Comparison principle]\label{SF:lemma:proof-general-LDP:CP-for-H}
	Let~$H$ be operator from Definition~\ref{SF:def:multi-valued-limit:general}. Then for any $h \in C_b(E)$ and $\lambda > 0$, the comparison principle holds for $f - \lambda Hf = h$.
\end{lemma}
\begin{proof}[Proof of Theorem~\ref{SF:thm:LDP_general}]
	By virtue of the above three Lemma's, the conditions of the general large-deviation result from Theorem~\ref{SF:thm:Richard-abstract-LDP} above are satisfied.
\end{proof}
While the verification of exponential compact containment and convergence of nonlinear operators are standard, the proof of the comparison principle takes up the bulk of the argument. We prove Lemma's~\ref{SF:lemma:proof-general-LDP:exp-comp-cont} and~\ref{SF:lemma:proof-general-LDP:conv-to-H} here, and prove Lemma~\ref{SF:lemma:proof-general-LDP:CP-for-H} in Section~\ref{SF:sec:proof-of-CP:general} below.
\begin{proof}[Proof of Lemma~\ref{SF:lemma:proof-general-LDP:conv-to-H}]
	We have to show that for all $(f,g) \in H$ there are functions $f_n \in \cD(H_n)$ satisfying $\LIM f_n = f$ and $\LIM H_n f_n = g$, with the~$\LIM$ convergence from Definition~\ref{SF:def:LIM-convergence}. Recall the slow-fast generator~$A_n$, 
	\begin{equation*}
		A_n f(y,z) = A^\mathrm{slow}_{n,z} f(\cdot,z) (y) +
		r_n \cdot A^\mathrm{fast}_{n,y} f(y,\cdot)(z),
	\end{equation*}
	and the nonlinear generators~$H_nf=(1/r_n)e^{-r_nf}A_ne^{r_nf}$, which amounts to
	\begin{equation*}
		H_n f(y,z) = \frac{1}{n}e^{-nf(y,z)}A^\mathrm{slow}_{n,z}\left( e^{nf(\cdot,z)}\right)(y) + r_n\cdot \frac{1}{n} e^{-nf(x,z)} \left(A^\mathrm{fast}_{n,y}e^{nf(x,\cdot)}\right)(z).
	\end{equation*}
	Now let~$(f,H_{f,\phi})\in H$ be arbitrary. Set~$f_n(y,z) := f(\eta_n(y)) + r_n^{-1} \phi(z)$. Then~$f_n\in\mathcal{D}(H_n)$. We are left with proving $\LIM f_n = f$ and $\LIM H_n f_n = H_{f,\phi}$.
	\smallskip
	
	Since both $f$ and $\phi$ are bounded, $\vn{f_n - f} = r_n^{-1}\vn{\phi} \rightarrow 0$, and $\LIM f_n = f$ follows. The images $H_n f_n$ are given by
	\begin{align*}
		H_n f_n(y,z)
		=
		\frac{1}{r_n} e^{-r_nf(\eta_n(y))} A_{n,z}^\mathrm{slow} e^{r_nf(\eta_n(y))} + e^{-\phi(z)} \left(A^{\mathrm{fast}}_{\eta_n(y)} e^{\phi}\right)(z).
	\end{align*}
	The convergence assumptions on the slow generators and fast generators (Assumptions~\ref{SF:assumption:convergence-slow-nonlinear-generators} and~\ref{SF:assumption:convergence-fast-nonlinear-generators}) imply $\LIM H_nf_n=H_{f,\phi}$.
\end{proof}
\begin{proof}[Proof of Lemma~\ref{SF:lemma:proof-general-LDP:exp-comp-cont}]
	By~\ref{SF:item:assumption:slow_regularity:compact_containment} of Assumption~\ref{SF:assumption:regularity-V} on the slow Hamiltonians~$V$, there are a compact containment function~$\Upsilon$ and a constant~$c_\Upsilon>0$ satisfying 
	\begin{equation*}
	\sup_{x,z} V_{x,\nabla \Upsilon(x)}(z) \leq c_\Upsilon.	
	\end{equation*}
	Choose $\beta > 0$ such that $T\cdot c_\Upsilon + 1 - \beta \leq -a$. By continuity of~$\Upsilon$ is, there is a constant $c>0$ such that
	\begin{equation*}
		K \subseteq \left\{x \, \middle| \, \Upsilon(x) \leq c \right\}
	\end{equation*}
	Next, let $G := \left\{x \, \middle| \, \Upsilon(x) < c + \beta \right\}$, which is an open set. Let $\widehat{K}$ be the closure of~$G$. Then~$\widehat{K}$ is compact since~$\Upsilon$ is a compact containment function.
	\smallskip
	
	Let $f(x) := \iota \circ \Upsilon $, where $\iota$ is some smooth increasing function satisfying
	\begin{equation*}
		\iota(r) = \begin{cases}
			r & \text{if } r \leq \beta +c, \\
			\beta+ c + 1 & \text{if } r \geq \beta + c + 2.
		\end{cases}
	\end{equation*}
	Then $\iota \circ \Upsilon = \Upsilon$ on $\widehat{K}$, and~$f$ is constant outside of a compact set. Set $f_n = f \circ \eta_n$, $g_n := H_n f_n$ and $g := \LIM g_n$. The function~$g$ exists due to Assumption~\ref{SF:assumption:convergence-slow-nonlinear-generators} on the slow generators. Then $g(x,z) = V_{x,\nabla \Upsilon(x)}(z)$ if $x \in \widehat{K}$. Therefore, we have $\sup_{x \in \widehat{K}, z \in F} g(x,z) \leq c_\Upsilon$.

	Let $\tau$ be the stopping time $\tau := \inf \left\{t \geq 0 \, \middle| \, Y_n(t) \notin \eta_n^{-1}(G) \right\}$ and let
	\begin{equation*}
		M_n(t) := \exp \left\{r_n \left( f_n(Y_n(t)) - f_n(Y_n(t)) - \int_0^t g_n(Y_n(t),Z_n(t)) \dd s \right) \right\}.
	\end{equation*}
	By construction $M_n$ is a martingale. By the optional stopping theorem, $t \mapsto M_n(t \wedge \tau)$ is a martingale as well. We obtain
	\begin{align*}
		& \PR_{y,z}\left[Y_n(t) \notin \widehat{K} \text{ for some } t \in [0,T]\right] \\
		& \leq \PR_{y,z}\left[Y_n(t) \notin \eta_n^{-1}(G) \times F \text{ for some } t \in [0,T]\right] \\
		& = \bE_{y,z}\left[\bONE_{\{Y_n(t) \notin \eta_n^{-1}(G) \text{ for some } t \in [0,T]\}} M_n(t \wedge \tau) M_n(t\wedge \tau)^{-1} \right] \\
		& \leq \exp\left\{- r_n \left(\inf_{y_1 \in G^c} \Upsilon(\eta_n(y_1)) - \Upsilon(\eta_n(y))   \right. \right. \\
		& \hspace{4cm} \left. \left. - T \sup_{y_2 \in \eta_n^{-1}(G), z_2 \in F} g_n(y_2,z_2) \right) \right\} \\
		& \hspace{2.5cm} \times  \bE_{y,z}\left[\bONE_{\{Y_n(t) \notin \eta_n^{-1}(G) \text{ for some } t \in [0,T]\}} M_n(t \wedge \tau) \right] 
	\end{align*}
	Since $\LIM f_n = f$ and $\LIM g_n = g$, the term in the exponential is bounded by $ r_n\left(c_\Upsilon T - \beta \right) \leq -r_n a$ for sufficiently large $n$. The final expectation is bounded by $1$ due to the martingale property of $M_n(t \wedge \tau)$. We conclude that
	\begin{align*}
		\limsup_n \sup_{y \in \eta_n^{-1}(K), z \in F} \frac{1}{r_n} \log \PR_{y,z}\left[Y_n(t) \notin \eta^{-1}(\widehat{K}) \text{ for some } t \in [0,T]\right] \leq -a,
	\end{align*}
	which finishes the proof.
\end{proof}
\subsection{Proof of the comparison principle}
\label{SF:sec:proof-of-CP:general}
In this section, we prove Lemma~\ref{SF:lemma:proof-general-LDP:CP-for-H}; the comparison principle for the Hamilton-Jacobi equation $f - \lambda Hf = h$ for the multi-valued limit operator~$H$ from Definition~\ref{SF:def:multi-valued-limit:general} introduce above. We recall the definition of viscosity solutions and the comparison principle in the appendix. 
\smallskip

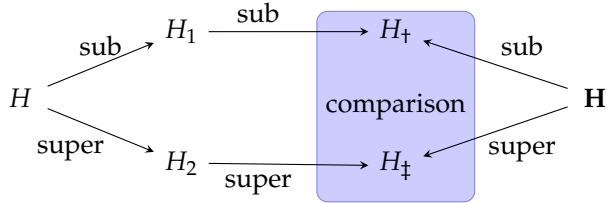
\begin{figure}[h!]
	\begin{center}
		\begin{tikzpicture}
		\matrix (m) [matrix of math nodes,row sep=1em,column sep=4em,minimum width=2em]
		{
			{ } & H_1 &[7mm] H_\dagger &[5mm] { } \\
			H & { } & { } & \bfH \\
			{ }  & H_2 & H_\ddagger & { } \\};
		\path[-stealth]
		(m-2-1) edge node [above] {sub} (m-1-2)
		(m-2-1) edge node [below] {super \qquad { }} (m-3-2)
		(m-1-2) edge node [above] {sub \qquad { }} (m-1-3)
		(m-3-2) edge node [below] {super \qquad { }} (m-3-3)
		(m-2-4) edge node [above] {\qquad sub} (m-1-3)
		(m-2-4) edge node [below] {\qquad super} (m-3-3);
		
		\begin{pgfonlayer}{background}
		\node at (m-2-3) [rectangle,draw=blue!50,fill=blue!20,rounded corners, minimum width=1cm, minimum height=2.5cm]  {comparison};
		\end{pgfonlayer}
		\end{tikzpicture}
	\end{center}
	\caption{An arrow connecting an operator $A$ with operator $B$ with subscript 'sub' means that viscosity subsolutions of $f - \lambda A f = h$ are also viscosity subsolutions of $f - \lambda B f = h$. Similarly for arrows with a subscript 'super'. The box around the operators $H_\dagger$ and $H_\ddagger$ indicates that the comparison principle holds for subsolutions of $f - \lambda H_\dagger f = h$ and supersolutions of $f - \lambda H_\ddagger f = h$.}
	\label{SF:fig:CP-diagram-in-proof-of-CP}
\end{figure}

A key role is played by the principal eigenvalue~$\mathcal{H}(x,p)$ from~\eqref{SF:eq:sec-assumption:H},
\begin{equation*}
\cH(x,p) = \sup_{\pi \in \cP(F)} \left\{\int V_{x,p}(z) \, \pi(\dd z) - \cI(x,\pi) \right\},\quad (x,p)\in E\times\mathbb{R}^d,
\end{equation*}
where the maps~$V$ and~$\mathcal{I}$ satisfy Assumptions~\ref{SF:assumption:regularity-V} and~\ref{SF:assumption:regularity-I}. We associate the following differential operator~$\mathbf{H}$ to this Hamiltonian.
\begin{definition}[Principal-eigenvalue Hamiltonian]\label{SF:def:principal-eigenvalue-Hamiltonian}
	The operator~$\mathbf{H}$ acting on the domain~$\cD(\bfH) = C_{cc}^\infty(E)$ is given by $\bfH f(x) := \cH(x, \nabla f(x))$.\qed
\end{definition}
The prove of Lemma~\ref{SF:lemma:proof-general-LDP:CP-for-H} hinges on being able to reduce the comparison principle of~$H$ to the comparison principle of~$\mathbf{H}$. To that end, we introduce four auxiliary operators and establish the diagram shown in Figure~\ref{SF:fig:CP-diagram-in-proof-of-CP}. 
\smallskip

The theoretical treatment of the Hamitlon-Jacobi equation of~$\mathbf{H}$ is carried out in~\cite{KraaijSchlottke2019}. In there, the comparison principle for the Hamilton-Jacobi equation associated with~$\mathbf{H}$ is proven~\cite[Theorem~3.4]{KraaijSchlottke2019} under a generalization of Assumptions~\ref{SF:assumption:regularity-V} and~\ref{SF:assumption:regularity-I}. The proof establishes the top-right and bottom-right arrows connecting~$\mathbf{H}$ with the auxiliary operators~$H_\dagger,H_\ddagger$. Here, we make the connection the Hamilton-Jacobi equation with the limit operator~$H$ by establishing the remaining arrows in Figure~\ref{SF:fig:CP-diagram-in-proof-of-CP}.
\begin{proof}[Proof of Lemma~\ref{SF:lemma:proof-general-LDP:CP-for-H}]
	Fix~$h\in C_b(E)$ and~$\lambda>0$. Let~$u_1$ be a viscosity subsolution and~$u_2$ be a viscosity supersolution to~$(1-\lambda H)f=h$. By Figure~\ref{SF:fig:CP-diagram-in-proof-of-CP}, the function~$u_1$ is a viscosity subsolution to~$(1-\lambda H_\dagger)f=h$ and~$u_2$ is a viscosity supersolution to~$(1-\lambda H_\ddagger)f=h$. Hence by~\cite[Theorem~3.4]{KraaijSchlottke2019},~$u_1\leq u_2$, which finishes the proof.
\end{proof}
The rest of this section is devoted to proving Figure~\ref{SF:fig:CP-diagram-in-proof-of-CP}.
\paragraph{Definition of auxiliary operators.}
We introduce the auxiliary operators $H_\dagger,H_\ddagger$ and $H_1,H_2$ appearing in Figure~\ref{SF:fig:CP-diagram-in-proof-of-CP}. The new Hamiltonians serve as natural upper and lower bounds for~$\bfH$ and~$H$, respectively. These new Hamiltonians are defined in terms of the containment function~$\Upsilon$ from Assumption~\ref{SF:assumption:regularity-V}, which allows us to restrict the analysis to compact sets. The definitions use the constant $C_\Upsilon := \sup_{x,z} V_{x,\nabla \Upsilon(x)}(z)$. Denote by~$C_l^\infty(E)$ the set of smooth functions on~$E$ that have a lower bound and by~$C_u^\infty(E)$ the set of smooth functions on~$E$ that have an upper bound.
\begin{definition}[$H_\dagger$ and~$H_\ddagger$]
	For $f \in C_l^\infty(E)$  and $\varepsilon \in (0,1)$, set 
	\begin{equation*}
		f^\varepsilon_\dagger := (1-\varepsilon) f + \varepsilon \Upsilon \quad\text{and}\quad
		H_{\dagger,f}^\varepsilon(x) := (1-\varepsilon) \bfH f(x) + \varepsilon C_\Upsilon.
	\end{equation*}
	Then~$H_\dagger$ is the the graph defined by
	\begin{equation*}
		H_\dagger := \left\{(f^\varepsilon_\dagger,H_{\dagger,f}^\varepsilon) \, \middle| \, f \in C_l^\infty(E), \varepsilon \in (0,1) \right\}.
	\end{equation*} 
	For $f \in C_u^\infty(E)$ and $\varepsilon \in (0,1)$, set 
	\begin{equation*}
		f^\varepsilon_\ddagger := (1+\varepsilon) f - \varepsilon \Upsilon \quad\text{and}\quad
		H_{\ddagger,f}^\varepsilon(x) := (1+\varepsilon) \bfH f(x) - \varepsilon C_\Upsilon.
	\end{equation*}
	Then~$H_\ddagger$ is the graph defined by
	\begin{align*}
		H_\ddagger := \left\{(f^\varepsilon_\ddagger,H_{\ddagger,f}^\varepsilon) \, \middle| \, f \in C_u^\infty(E), \varepsilon \in (0,1) \right\}.
		\tag*\qed
	\end{align*}
\end{definition}
\begin{definition}[$H_1$ and~$H_2$]
	For $f \in C_l^\infty(E)$ , $\varepsilon \in (0,1)$ and $\phi$ such that $e^\phi \in \cD(A^\mathrm{fast})$, set 
	\begin{gather*}
		f^\varepsilon_1 := (1-\varepsilon) f + \varepsilon \Upsilon, \\
		H^\varepsilon_{1,f,\phi}(x,z) :=
		(1-\varepsilon) \left( V_{x,\nabla f(x)}(z) + e^{-\phi(z)} A^\mathrm{fast}_x e^{\phi}(z)\right) + \varepsilon C_\Upsilon.
	\end{gather*}
	Then~$H_1$ is the graph defined by
	\begin{equation*}
		H_1 := \left\{(f^\varepsilon_1,H^\varepsilon_{1,f,\phi}) \, \middle| \, f \in C_l^\infty(E), \varepsilon \in (0,1), \phi: \, e^\phi \in \cD(A^\mathrm{fast})  \right\}.
	\end{equation*} 
	For $f \in C_u^\infty(E)$, $\varepsilon \in (0,1)$ and $\phi$ such that $e^\phi \in \cD(A^\mathrm{fast})$, set 
	\begin{gather*}
		f^\varepsilon_2 := (1+\varepsilon) f - \varepsilon \Upsilon, \\
		H^\varepsilon_{2,f,\phi}(x,z) :=
		(1+\varepsilon) \left( V_{x,\nabla f(x)}(z)  + e^{-\phi(z)} A^\mathrm{fast}_x e^{\phi}(z) \right) - \varepsilon C_\Upsilon.
	\end{gather*}
	Then~$H_2$ is the graph defined by
	\begin{align*}
		H_2 := \left\{(f^\varepsilon_2,H^\varepsilon_{2,f,\phi}) \, \middle| \, f \in C_u^\infty(E), \varepsilon \in (0,1), \phi: \, e^\phi \in \cD(A^\mathrm{fast})  \right\}.
		\tag*\qed
	\end{align*}
\end{definition}
\paragraph{Arrows based on the solution of an eigenvalue problem.}
\begin{lemma}\label{SF:lemma:viscosity_solutions_arrows_based_on_eigenvalue}
	Fix $\lambda > 0$ and $h \in C_b(E)$. 
	\begin{enumerate}[(a)]
		\item Every subsolution to $f - \lambda H_1 f = h$ is also a subsolution to $f - \lambda H_\dagger f = h$.
		\item Every supersolution to $f - \lambda H_2 f = h$ is also a supersolution to $f - \lambda H_\ddagger f = h$.
	\end{enumerate}
\end{lemma}
For the proof of this lemma, we need an auxiliary lemma.
\begin{lemma} \label{SF:lemma:strong_viscosity_solutions}
	Fix $\lambda > 0$ and $h \in C_b(E)$. 
	\begin{enumerate}[(a)]
		\item Let $u$ be a subsolution to $f - \lambda H_1 f = h$, then for all $(f,g) \in H_1$ and $x_0 \in E$ such that
		\begin{equation*}
			u_1(x_0) - f(x_0) = \sup_x u_1(x) - f(x)
		\end{equation*}
		we have
		\begin{equation*}
			u_1(x_0) - \lambda g(x_0,z)  \leq h(x_0).
		\end{equation*}
		\item Let $u_2$ be a supersolution to $f - \lambda H_2 f = h$, then for all $(f,g) \in H_2$ and $x_0 \in E$ such that
		\begin{equation*}
			u_2(x_0) - f(x_0) = \inf_x u_2(x) - f(x)
		\end{equation*}
		we have
		\begin{equation*}
			u_2(x_0) - \lambda g(x_0,z)  \geq h(x_0).
		\end{equation*}
	\end{enumerate}
\end{lemma}
The following proof is inspired on \cite[Lemma 9.9]{FengKurtz2006}.
\begin{proof}[Proof of Lemma~\ref{SF:lemma:strong_viscosity_solutions}]
	We only prove (a). Let $u$ be a viscosity subsolution to $f - \lambda H_1 f = h$ and consider $(f,g) \in H_1$. By definition $f$ has compact sublevel-sets. Thus, instead of working with a sequence $x_n$ along which a maximum is attained, we can work with a single point $x_0$. This gives us the existence of a point $(x_0,z) \in E \times F$ satisfying
	\begin{gather*}
		u(x_0) - f(x_0) = \sup_x u_1(x) - f(x), \\
		u(x_0) - \lambda g(x_0,z)  \leq h(x_0),
	\end{gather*}
	rather than having the second inequality for all $x_0$ satisfying $u(x_0) - f(x_0) = \sup_x u(x) - f(x)$.
	
	Now let $x_0$ be such that $u(x_0) - f(x_0) = \sup_x u(x) - f(x)$. Pick a function $\hat{f} \in C_{cc}^\infty(E)$ satisfying $\hat{f}(x_0)  = 0$ and $\hat{f}(x) > 0$ for $x \neq x_0$. Define the function $f_0 = f+ \hat{f}$, and let $g_0$ be the corresponding image, $(f_0,g_0) \in H_1$. Since $\nabla f_0(x_0) = \nabla f(x_0)$ and $g(x_0,z)$ and $g_0(x_0,z)$ only depend on $f$ and $f_0$ via their derivatives at $x_0$, we obtain $g_0(x_0,z) = g(x_0,z)$.
	By construction $x_0$ is the unique point satisfying $u(x_0) - f_0(x_0) = \sup_x u(x) - f_0(x)$. By the sub-solution property, we find
	\begin{equation*}
		u(x_0) - \lambda g(x_0,z) = u(x_0) - \lambda g_0(x_0,z) \leq h_0(x_0),
	\end{equation*}
	establishing the claim.
\end{proof}

\begin{proof}[Proof of Lemma~\ref{SF:lemma:viscosity_solutions_arrows_based_on_eigenvalue}]
	We only prove the subsolution statement. To that end, fix $\lambda > 0$ and $h \in C_b(E)$, and let $u$ be a subsolution of $f - \lambda H_1 f = h$. We prove it is also a subsolution of $f - \lambda H_\dagger f = h$. Let $f^\varepsilon_1 \in \cD(H_1)$ and let $x_0$ be such that
	\begin{equation*}
		u(x_0) - f^\varepsilon_1(x_0) = \sup_x u(x) - f_1^\varepsilon(x).
	\end{equation*}
	For each $\delta > 0$ we find by Assumption~\ref{SF:assumption:principal-eigenvalue-problem} a function $e^{\phi_\delta} \in \cD(A^\mathrm{fast})$ such that
	\begin{equation*}
		\cH(x,p) \geq V_{x_0,\nabla f(x_0)}(z) - e^{-\phi_\delta(z)}\left(A_{x_0}^\mathrm{fast} e^{\phi_\delta}\right)(z) - \delta
	\end{equation*}
	for all $z \in F$. Since
	\begin{equation*}
		\left(f^\varepsilon_1, (1-\varepsilon) \left( V_{x,\nabla f(x)}(z) + e^{-\phi(z)} A^\mathrm{fast}_x e^{\phi(z)}\right) + \varepsilon C_\Upsilon \right) \in H_1,
	\end{equation*}
	we find by the subsolution property of $u$ and Lemma~\ref{SF:lemma:strong_viscosity_solutions} that for all $z$
	\begin{align*}
		h(x_0) & \geq u(x_0) - \lambda \left((1-\varepsilon) \left( V_{x_0,\nabla f(x_0)}(z) + e^{-\phi(z)} A^\mathrm{fast}_{x_0} e^{\phi(z)}\right) + \varepsilon C_\Upsilon\right) \\
		& \geq	u(x_0) - \lambda \left((1-\varepsilon)  \cH(x_0,\nabla f(x_0)) + \varepsilon C_\Upsilon\right) - \lambda(1-\varepsilon)\delta.
	\end{align*}
	Sending $\delta \rightarrow 0$ establishes that $u$ is a subsolution for $f - \lambda H_\dagger f = h$.
\end{proof}
\paragraph{Arrows based on compact containment.}
\begin{lemma} \label{SF:lemma:viscosity_solutions_compactify1}
	Fix $\lambda > 0$ and $h \in C_b(E)$. 
	\begin{enumerate}[(a)]
		\item Every subsolution to $f - \lambda H f = h$ is also a subsolution to $f - \lambda H_1 f = h$.
		\item Every supersolution to $f - \lambda H f = h$ is also a supersolution to $f - \lambda H_2 f = h$.
	\end{enumerate}
\end{lemma}
\begin{lemma} \label{SF:lemma:viscosity_solutions_compactify2}
	Fix $\lambda > 0$ and $h \in C_b(E)$. 
	\begin{enumerate}[(a)]
		\item Every subsolution to $f - \lambda \bfH f = h$ is also a subsolution to $f - \lambda H_\dagger f = h$.
		\item Every supersolution to $f - \lambda \bfH f = h$ is also a supersolution to $f - \lambda H_\ddagger f = h$.
	\end{enumerate}
\end{lemma}
Lemma~\ref{SF:lemma:viscosity_solutions_compactify2} has been proven in~\cite[Lemma~6.3]{KraaijSchlottke2019}.
\begin{proof}[Proof of Lemma~\ref{SF:lemma:viscosity_solutions_compactify1}]
	The proof is similar to the proof of~\cite[Lemma~6.3~(a)]{KraaijSchlottke2019}. We only prove~(a).
	\smallskip
	
	Fix $\lambda > 0$ and $h \in C_b(E)$. Let $u$ be a subsolution to $f - \lambda H f = h$. We prove it is also a subsolution to $f - \lambda H_1 f = h$. Fix $\varepsilon \in (0,1)$, $\phi$ such that $e^\phi \in \mathcal{D}(A^\mathrm{fast})$, and $f \in C_{l}^\infty(E)$, so that $(f^\varepsilon_1,H^\varepsilon_{1,f,\phi}) \in H_1$. We will prove that there are $(x_n,z_n)$ such that
	\begin{gather} 
		\lim_n u(x_n) - f^\varepsilon_1(x_n) = \sup_x u(x) - f^\varepsilon_1(x),\label{SF:eqn:proof_lemma_conditions_for_subsolution_first} \\
		\limsup_n u(x_n) - \lambda H^\varepsilon_{1,f,\phi}(x_n,z_n) - h(x_n) \leq 0. \label{SF:eqn:proof_lemma_conditions_for_subsolution_second}
	\end{gather}
	As $u -(1-\varepsilon)f$ is bounded from above and $\varepsilon \Upsilon$ has compact sublevel-sets, the sequence $x_n$ along which the first limit is attained can be assumed to lie in the compact set $K := \left\{x \, | \, \Upsilon(x) \leq \inf_x \varepsilon^{-1} \left(u(x) - (1-\varepsilon)f(x) \right)\right\}$. We use the constant~$M := \inf_x \varepsilon^{-1} \left(u(x) - (1-\varepsilon)f(x) \right)$. Let $\gamma : \bR \rightarrow \bR$ be a smooth increasing function such that
	\begin{equation*}
		\gamma(r) = \begin{cases}
			r & \text{if } r \leq M, \\
			M + 1 & \text{if } r \geq M+2.
		\end{cases}
	\end{equation*}
	Denote by $f_\varepsilon$ the function on $E$ defined by 
	\begin{equation*}
		f_\varepsilon(x) := \gamma\left((1-\varepsilon)f(x) + \varepsilon \Upsilon(x) \right).
	\end{equation*}
	By construction $f_\varepsilon$ is smooth and constant outside of a compact set and thus lies in $\cD(H) = C_{cc}^\infty(E)$. As $e^\phi \in \mathcal{D}(A^\mathrm{fast})$ we have by Assumption~\ref{SF:assumption:principal-eigenvalue-problem} that also $e^{(1-\varepsilon)\phi} \in \mathcal{D}(A^\mathrm{fast})$. We conclude that $(f_\varepsilon, H_{f,(1-\varepsilon)\phi}) \in H$. 
	
	As $u$ is a viscosity subsolution for $f - \lambda Hf = h$ there exist $x_n \in K \subseteq E$ (by our choice of $K$) and $z_n \in F$ with
	\begin{gather}
		\lim_n u(x_n) - f_\varepsilon(x_n) = \sup_x u(x) - f_\varepsilon(x), \label{SF:eqn:visc_subsol_sup} \\
		\limsup_n u(x_n) - \lambda H_{f,(1-\varepsilon)\phi}(x_n,z_n) - h(x_n) \leq 0. \label{SF:eqn:visc_subsol_upperbound}
	\end{gather}
	As $f_\varepsilon$ equals $f$ on $K$, we have from \eqref{SF:eqn:visc_subsol_sup} that also
	\begin{equation*}
		\lim_n u(x_n) - f(x_n) = \sup_x u(x) - f(x),
	\end{equation*}
	establishing \eqref{SF:eqn:proof_lemma_conditions_for_subsolution_first}. Convexity of $p \mapsto V_{x,p}$ and $\psi \mapsto e^{-\psi(z)}\left(A_x^\mathrm{fast} e^\psi\right)(z)$ yields for arbitrary $(x,z)$ the elementary estimate
	\begin{align*}
		H_{f,(1-\varepsilon)\phi}(x,z) = & V_{x,\nabla f_\varepsilon}(z) + e^{-(1-\varepsilon)\phi(z)} \left(A^\mathrm{fast}_x e^{(1-\varepsilon)\phi}\right)(z) \\
		& \leq (1-\varepsilon) V_{x,\nabla f(x)}(z) + \varepsilon V_{x,\nabla \Upsilon(x)} + (1-\varepsilon)e^{-\phi(z)} \left(A^\mathrm{fast}_x e^{\phi}\right)(z) \\
		& = H^\varepsilon_{1,f,\phi}(x,z).
	\end{align*} 
	Combining this inequality with \eqref{SF:eqn:visc_subsol_upperbound} yields
	\begin{multline*}
		\limsup_n u(x_n) - \lambda H^\varepsilon_{1,f,\phi}(x,z) - h(x_n) \\
		\leq \limsup_n u(x_n) - \lambda H_{f,(1-\varepsilon)\phi}(x_n,z_n) - h(x_n) \leq 0,
	\end{multline*}
	establishing \eqref{SF:eqn:proof_lemma_conditions_for_subsolution_second}. This concludes the proof.
\end{proof}
\section{Proof of action-integral representation}\label{SF:sec:proof-action-integral-form-RF}
\subsection{Structure of proof}
In this section, we outline the structure of proof of Theorem~\ref{SF:thm:rate-function}. To that end, recall the Hamiltonian~$\mathcal{H}:E\times\mathbb{R}^d\to\mathbb{R}$ from~\eqref{SF:eq:variational_Hamiltonian},
\begin{equation*}
	\cH(x,p) = \sup_{\nu \in \cP(F)} \left\{\int V_{x,p}(z) \, \nu(\dd z) - \cI(x,\nu) \right\},
\end{equation*}
and the Lagrangian~$\mathcal{L}(x,v)$ defined as the Legendre dual,
\begin{equation}\label{SF:prof-RF:Lagrangian}
	\mathcal{L}(x,v)=\sup_{p\in\mathbb{R}^d}\left[\ip{p}{v}-\mathcal{H}(x,p)\right].
\end{equation}
Theorem~\ref{SF:thm:rate-function} consists of two claims:
\begin{enumerate}[label=(\Roman*)]
	\item \label{SF:item:proof-action:I} The rate function~$J$ obtained in Theorem~\ref{SF:thm:LDP_general} satisfies the action-integral representation~\eqref{SF:eq:action-integral-rate-function} with the Lagrangian~\eqref{SF:prof-RF:Lagrangian}.
	\item \label{SF:item:proof-action:Lagr} The Lagrangian satisfies the formula~\eqref{SF:eq:results-Lagr-opt-proc}.
\end{enumerate}
We prove the statement~\ref{SF:item:proof-action:Lagr} in Section~\ref{SF:sec:proof-alternative-form-of-L}. To outline the proof of~\ref{SF:item:proof-action:I}, recall that the rate function~$J$ obtained in Theorem~\ref{SF:thm:LDP_general} given  by~\eqref{SF:eq:Richard-general-LDP-thm-rate-function} is characterized in terms of a semigroup~$V(t)$, as
\begin{equation*}
J(\gamma) = J_0(\gamma(0)) + \sup_{k \geq 1} \sup_{\{t_i\}} \sum_{i=1}^{k} J_{t_i - t_{i-1}}(\gamma(t_i) \, | \, \gamma(t_{i-1})),
\end{equation*}
where the conditional rate functions~$J_t$ are defined by
\begin{equation*}
	J_t(y \, | \, x) = \sup_{f \in C_b(E)} \left\{f(y) - V(t)f(x) \right\}.
\end{equation*}
We describe the semigroup~$V(t):C_b(E)\to C_b(E)$ further below. Our proof is based on the following definitions. First, we define the variational semigroup~$\mathbf{V}(t)$ arising from an optimal control problem with cost function~$\mathcal{L}$.
\begin{definition}[Variational semigroup~$\mathbf{V}(t)$]
	For~$f\in C_b(E)$ and~$t\geq 0$, define
	\begin{align*}
		\bfV(t) f(x) := \sup_{\substack{\gamma \in \mathcal{A}\mathcal{C}\\ \gamma(0) = x}} h(\gamma(t)) -  \int_0^t \mathcal{L}(\gamma(s),\dot{\gamma}(s)),
	\end{align*}
	where the supremum is over absolutely continuous maps~$\gamma:[0,\infty)\to E$.\qed
\end{definition}
Secondly, we will exploit the fact that the semigroups~$V(t)$ and~$\mathbf{V}(t)$ are determined by means of the existence of unique viscosity solutions to the equations~$(1-\lambda H)f=h$ and~$(1-\lambda\mathbf{H})f=h$, respectively. To do so, we introduce the following resolvents.
\begin{definition}[Resolvent~$R(\lambda)$]
	Let~$H$ be the multi-valued operator from Definition~\ref{SF:def:multi-valued-limit:general},~$\lambda>0$ and~$h\in C_b(E)$. Define~$R(\lambda):C_b(E)\to C_b(E)$ by setting~$R(\lambda)h := u$, the unique viscosity solution to~$(1-\lambda H)u=h$.\qed
\end{definition}
\begin{definition}[Resolvent~$\mathbf{R}(\lambda)$]
	For~$\lambda>0$ and~$h\in C_b(E)$, define the resolvent~$\mathbf{R}(\lambda)h:E\to\mathbb{R}$ by
	\begin{equation}
	\mathbf{R}(\lambda) h(x):= \sup_{\substack{\gamma \in \mathcal{A}\mathcal{C}\\ \gamma(0) = x}} \int_0^\infty \lambda^{-1} e^{-\lambda^{-1}t} \left[h(\gamma(t)) - \int_0^t \mathcal{L}(\gamma(s),\dot{\gamma}(s))\right] \, \dd t,
	\end{equation} 
	where the supremum is over absolutely-continuous maps~$\gamma:[0,\infty)\to E$.\qed
\end{definition}
The semigroup~$V(t)$ is determined via~$R$ as (e.g.\cite[Prop.~6.6]{Kraaij2019ExpResolv})
\begin{equation*}
	V(t)f(x) = \lim_{m\to\infty}\left[R\left(t/m\right)\right]^{m}f.
\end{equation*}
The statement~\ref{SF:item:proof-action:I} is a direct consequence of the following Lemmas. 
\begin{lemma}\label{SF:lemma:R-equals-mathbf-R}
	Let Assumptions~\ref{SF:assumption:regularity-V},~\ref{SF:assumption:regularity-I} and~\ref{SF:assumption:Hamiltonian_vector_field} be satisfied. Then for all~$\lambda>0$, we have~$R(\lambda)=\mathbf{R}(\lambda)$.
\end{lemma}
\begin{lemma}\label{SF:lemma:V-equals-mathbf-V}
	Let Assumptions~\ref{SF:assumption:regularity-V},~\ref{SF:assumption:regularity-I} and~\ref{SF:assumption:Hamiltonian_vector_field} be satisfied. If~$R(\lambda)=\mathbf{R}(\lambda)$ for all~$\lambda>0$, then~$V(t)=\mathbf{V}(t)$.
\end{lemma}
\begin{lemma}\label{SF:lemma:equal-semigroup-implies-RF}
	If~$V(t)=\mathbf{V}(t)$, then~\ref{SF:item:proof-action:I} holds true. 
\end{lemma}
Lemma~\ref{SF:lemma:equal-semigroup-implies-RF} can be proven as shown in~\cite[Theorem~8.14]{FengKurtz2006}, using convexity of~$v\mapsto \mathcal{L}(x,v)$. We show Lemmas~\ref{SF:lemma:R-equals-mathbf-R} and~\ref{SF:lemma:V-equals-mathbf-V} in Section~\ref{SF:thm:proof-of-RF-action}.
\subsection{Proof of Theorem~\ref{SF:thm:rate-function}}
\label{SF:thm:proof-of-RF-action}
For the proof of Lemma~\ref{SF:lemma:R-equals-mathbf-R}, we argue with the diagram in~Figure~\ref{SF:fig:CP-diagram-in-proof-of-CP} that we established in the proof of the comparison principle. We recall it here below.
\begin{figure}[h!]
	\begin{center}
		\begin{tikzpicture}
		\matrix (m) [matrix of math nodes,row sep=1em,column sep=4em,minimum width=2em]
		{
			{ } & H_1 &[7mm] H_\dagger &[5mm] { } \\
			H & { } & { } & \bfH \\
			{ }  & H_2 & H_\ddagger & { } \\};
		\path[-stealth]
		(m-2-1) edge node [above] {sub} (m-1-2)
		(m-2-1) edge node [below] {super \qquad { }} (m-3-2)
		(m-1-2) edge node [above] {sub \qquad { }} (m-1-3)
		(m-3-2) edge node [below] {super \qquad { }} (m-3-3)
		(m-2-4) edge node [above] {\qquad sub} (m-1-3)
		(m-2-4) edge node [below] {\qquad super} (m-3-3);
		
		\begin{pgfonlayer}{background}
		\node at (m-2-3) [rectangle,draw=blue!50,fill=blue!20,rounded corners, minimum width=1cm, minimum height=2.5cm]  {comparison};
		\end{pgfonlayer}
		\end{tikzpicture}
	\end{center}
\end{figure}
\smallskip

\begin{proof}[Proof of Lemma~\ref{SF:lemma:R-equals-mathbf-R}]
	 Figure~\ref{SF:fig:CP-diagram-in-proof-of-CP} shows that if~$u$ is a viscosity solution to~$(1-\tau H)f=h$ and~$v$ is a viscosity solution to~$(1-\tau\mathbf{H})f=h$, then~$u=v$. 
	 \smallskip
	 
	 Let~$\lambda>0$ and~$h\in C_b(E)$. Then by definition,~$R(\lambda)h$ is the viscosity solution to~$(1-\lambda H)f=h$. We prove that the function~$\mathbf{R}(\lambda)h$ is the viscosity solution to~$(1-\lambda\mathbf{H})f=h$. Then by virtue of Figure~\ref{SF:fig:CP-diagram-in-proof-of-CP}, we obtain~$R(\lambda)h=\mathbf{R}(\lambda)h$. Since~$\lambda$ and~$h$ are arbitrary, this establishes Lemma~\ref{SF:lemma:R-equals-mathbf-R}.
	 \smallskip
	 
	 The fact that~$\mathbf{R}(\lambda)h$ is a viscosity solution is established in~\cite[Theorem~3.7]{KraaijSchlottke2019}, under Assumptions~3.12,~3.13 and~3.16 therein. Here, our Assumption~\ref{SF:assumption:regularity-V} corresponds exactly to \cite[Assumption~3.12]{KraaijSchlottke2019}, and our Assumption~\ref{SF:assumption:regularity-I} to~\cite[Assumption~3.13]{KraaijSchlottke2019}. We are left with showing that~\cite[Assumption 3.16]{KraaijSchlottke2019}, follows from our Assumption~\ref{SF:assumption:Hamiltonian_vector_field}. 
	 \smallskip
	 
	 To that end, for a convex function $\Phi$ $\Phi : \bR^d \rightarrow (-\infty,\infty]$, define the subdifferential set by
	 \begin{equation*}
	 	\partial_p \Phi(p_0)
	 	:= \left\{
	 	\xi \in \mathbb{R}^d \,:\, \Phi(p) \geq \Phi(p_0) + \xi \cdot (p-p_0) \quad (\forall p \in \mathbb{R}^d)
	 	\right\}.
	 \end{equation*}
	 Fix $x \in E$ and $p_0 \in \bR^d$. We aim to prove that $\partial_p \cH(x,p_0) \subseteq T_E(x)$. Since the map $p \mapsto \cH(x,p)$ is proper and convex as a supremum over convex functions, the subdifferential $\partial_p \cH(x,p_0)$ is non-empty. 
	 \smallskip
	 
	 Let $\Omega$ be the set that of measures~$\pi$ that optimize 
	 \begin{equation} \label{SF:eqn:Hamiltonian_in_subdifferential_proof}
	 \cH(x,p_0) = \sup_\pi \left\{\int V_{x,p_0}(z) \pi(\dd z) - \cI(x,\pi) \right\}
	 \end{equation}
	 
	 We first aim to relate $\partial_p \cH(x,p_0)$ to $\bigcup_{\pi^* \in \Omega}\partial_p \int V_{x,p_0}(z) \pi^*(\dd z)$. Afterwards, we show that for all $\pi$ we have $\partial_p \int V_{x,p_0}(z) \pi(\dd z) \subseteq T_E(x)$.
	 
	 \smallskip
	 
	 For each fixed $p$ we can restrict our supremum in~\eqref{SF:eqn:Hamiltonian_in_subdifferential_proof} to the compact set of measures $\pi$ such that $\cI(x,\pi) \leq 2\vn{V_{x,p}(\cdot)}$. For various $p$, this set might change and we might end up with a non-compact set. However, to study the subdifferential at $p_0$ we can instead at the map $p \mapsto \cH(x,p)$ with its domain restricted to $[p_0-1,p_0 + 1]$ which leaves the subdifferential set unchanged.
	 \smallskip
	 
	 Set $C_V := \sup_{p \in [p_0-1,p_0+1]} \sup_z \vn{v_{\cdot,p}(z)}_\infty < \infty$, and let~$\Xi$ to be the closure of
	 \begin{equation}
	 \{\pi \in \cP(F) \, \| \, \cI(x,\pi) \leq 2C_V\}.
	 \end{equation}
	 By Assumption~\ref{SF:assumption:regularity-I}, the set $\Xi$ is compact. Therefore, for all $p \in [p_0-1,p_0+1]$, we can restrict the supremum in~\eqref{SF:eqn:Hamiltonian_in_subdifferential_proof} to the compact set $\Xi$.
	 \smallskip
	 
	 Using the definition of $\Omega$ as the set of optimizers and that $\cI$ is lower semicontinuous by Assumption~\ref{SF:assumption:regularity-I}, it follows by \cite[Theorem 4.4.2]{HiriartLemarechal2012} that 
	 \begin{equation*}
	 	\partial_p \cH(x,p_0) = ch \left(\bigcup_{\pi^* \in \Omega} \partial_p \left(\int V_{x,p_0}(z) \pi^*(\dd z) - \cI(x,\pi^*) \right)\right).
	 \end{equation*}
	 Here $ch$ denotes the convex hull. Since $\cI(x,\pi^*)$ does not depend on $p$,
	 \begin{equation*}
	 	\partial_p \cH(x,p_0) = ch \left(\bigcup_{\pi^* \in \Omega} \partial_p \left(\int V_{x,p_0}(z) \pi^*(\dd z)\right)\right).
	 \end{equation*}
	 Since $\partial_p V_{x,p}(z) \subseteq T_E(x)$ for all $p$ and $z$, we find by~\cite[Theorem 3]{Papageorgiou1997} applied with $\varepsilon = 0$ that $\partial_p \cH(x,p_0) \subseteq T_E(x)$. This establishes~\cite[Assumption 3.16]{KraaijSchlottke2019}.
\end{proof} 
\begin{proof}[Proof of Lemma~\ref{SF:lemma:V-equals-mathbf-V}]
	By~\cite[Theorem 7.10]{Kraaij2019ExpResolv} and \cite[Theorem 6.1]{Kraaij2019GenConv}, there is some sequentially strictly dense set $D \subseteq C_b(E)$ such that for $f \in D$
	\begin{equation} \label{SF:eqn:convergence_R_to_V}
	\lim_{m \rightarrow \infty} \vn{R\left(\frac{t}{m}\right)^m f - V(t)f} = 0.
	\end{equation}
	Similarly, we find by \cite[Lemma 8.18]{FengKurtz2006} that for all $f \in C_b(E)$ and $x \in E$
	\begin{equation} \label{SF:eqn:convergence_bfR_to_bfV}
	\lim_{m \rightarrow \infty} \mathbf{R}\left(\frac{t}{m}\right)^m f(x) = \bfV(t)f(x).
	\end{equation}
	Combining~\eqref{SF:eqn:convergence_R_to_V} and~\eqref{SF:eqn:convergence_bfR_to_bfV} we conclude that $V(t)f = \bfV(t)f$ for all $t$ and $f \in D$. Since Figure~\ref{SF:fig:CP-diagram-in-proof-of-CP} implies $R(\lambda)h$ = $\bfR(\lambda)h$ for $h \in C_b(E)$, we conclude from~\eqref{SF:eqn:convergence_R_to_V} and~\eqref{SF:eqn:convergence_bfR_to_bfV} that $V(t)f = \bfV(t)f$ for all $t$ and $f \in D$.
	\smallskip
	
	Since $D$ is sequentially strictly dense by assumption, the equality for all $f \in C_b(E)$ follows if $V(t)$ and $\bfV(t)$ are sequentially continuous. The semigroup~$V(t)$ is seequentially continuous by~\cite[Theorem 7.10]{Kraaij2019ExpResolv} and \cite[Theorem~6.1]{Kraaij2019GenConv}, and~$\mathbf{V}(t)$ is sequentially continuous by~\cite[Lemma 8.22]{FengKurtz2006}. We conclude that $V(t)f = \bfV(t)f$ for all $f \in C_b(E)$ and $t \geq 0$.
\end{proof}
\subsection{Proof of alternative form of Lagrangian}
\label{SF:sec:proof-alternative-form-of-L}
We prove that the Lagrangian~$\mathcal{L}(x,v)$ defined in~\eqref{SF:prof-RF:Lagrangian} as the Legendre dual of~$\mathcal{H}(x,p)$ satisfies~\eqref{SF:eq:results-Lagr-opt-proc}. Recall that the Hamiltonian is
\begin{equation*}
	\mathcal{H}(x,p)=\sup_{\pi\in\mathcal{P}(F)}\left[\int_F V_{x,p}(z)-\mathcal{I}(x,\pi)\right].
\end{equation*}
For the proof, we write for~$\pi\in\mathcal{P}(F)$
\begin{equation*}
	\mathcal{H}_\pi(x,p) := \int_F V_{x,p}(z) \pi(\dd z)\quad\text{and}\quad
	\mathcal{L}_\pi(x,v) := \sup_p  \ip{p}{v} - \int_F V_{x,p}(z) \pi(\dd z).
\end{equation*}
Furtheremore, let $\Phi(v,\pi)$ be the set of measurable functions $w \in L^1(F;\mathbb{R}^d,\pi)$ such that 
\begin{equation*}
	\int w(z) \pi(\dd z) = v.
\end{equation*}
\begin{proposition}\label{SF:proposition:representation_of_Lagrangian}
	Suppose that Assumptions~\ref{SF:assumption:regularity-V} and~\ref{SF:assumption:regularity-I} are satisfied.
	Then
	\begin{equation*}
		\cL(x,v) = \inf_\pi \inf_{w \in \Phi(v,\pi)} \left\{\int \cL_z(x,w(z)) \pi(\dd z) + \mathcal{I}(x,\pi) \right\}.
	\end{equation*}
\end{proposition} 
This rewrite is a consequence of results in convex analysis and follows under much weaker assumptions, namely convexity of $\pi \mapsto \cI(x,\pi)$ and $p \mapsto V_{x,p}(z)$, which is satisfied in our setting. The proof below only uses results from convex analysis in~\cite{Papageorgiou1997,HiriartLemarechal2012,Rockafellar1970}. These results have been stated for~$\bR^d$, which is also the setting to which we restrict ourselves in this chapter. We believe, however, that the result should extend to a more general setting, but we were unable to find their generalizations in the literature on convex analysis.
\begin{lemma} \label{SF:lemma:constant_momentum}
	Fix~$\pi\in\mathcal{P}(F)$ and suppose that $v \in \text{rel. int. dom } \cL_\pi(x,\cdot)$. Then:
	\begin{enumerate}[(a)]
		\item There is a $p^* \in \bR^d$ such that $v \in \partial_p \cH_\pi(x,p^*)$,
		\item There is a $w^* \in \Phi(v,\pi)$ such that $w^*(z) \in \partial_p V_{x,p^*}(z)$ $\pi$ almost surely.
		\item \label{SF:item:equality_of_Lagrangians_in_relint} We have 
		\begin{multline*}
			\sup_{p} \ip{p}{v} - \int V_{x,p}(z) \pi(\dd z) \\
			= \inf_{\substack{w(z) \\ \int w \dd \pi = v}} \sup_{p(z)} \int \ip{p(z)}{w(z)} - V_{x,p(z)}(z) \pi(\dd z).
		\end{multline*}
	\end{enumerate}
\end{lemma}
For the proof of this lemma, we use the notion of the relative interior of a convex set. If $A \subseteq \bR^d$ is a convex set, then $\text{rel. int. } A$ is the interior of $A$ inside the smallest affine hyperplane in $\bR^d$ that contains $A$. For a convex functional $\Phi : E \rightarrow (-\infty,\infty]$ the domain of $\Phi$, denoted by $\text{dom } \Phi$, is the set of points $x \in E$ where $\Phi(x) < \infty$.
\begin{proof}[Proof of Lemma~\ref{SF:lemma:constant_momentum}]
	The Legendre transform of $\cH_\pi$ equals $\cL_\pi$. Since we have $v \in \text{rel. int. dom } \cL_\pi(x,\cdot)$, we have by~\cite[Theorem 23.4]{Rockafellar1970} or~\cite[Theorem E.1.4.2]{HiriartLemarechal2012} that $\partial_v \cL_\pi(x,v)$ is non-empty. Let $p^* \in \partial_v \cL_\pi(x,v)$. Then by~\cite[Theorem 23.5]{Rockafellar1970} or~\cite[Proposition E.1.4.3]{HiriartLemarechal2012}, we obtain that $v \in \partial_p H_\pi(x,p^*)$.
	\smallskip
	
	By~\cite[Theorem~3]{Papageorgiou1997} applied for $\varepsilon = 0$, we find a $\pi$-integrable function~$w^*$ such that $\int w^*(z) \pi(\dd z) = v$ and $w(z) \in \partial_p V_{x,p^*}(z)$ $\pi$ almost surely.	We proceed with the proof of \ref{SF:item:equality_of_Lagrangians_in_relint}. To that end, note that
	\begin{multline*}
		\sup_{p} \ip{p}{v} - \int V_{x,p}(z) \pi(\dd z) \\
		\leq \inf_{\substack{w(z) \\ \int w \dd \pi = v}} \sup_{p(z)} \int \ip{p(z)}{w(z)} - V_{x,p(z)}(z) \pi(\dd z).
	\end{multline*}
	For the other inequality, note that 
	\begin{align*}
		&  \inf_{\substack{w(z) \\ \int w \dd \pi = v}} \sup_{p(z)} \int \ip{p(z)}{w(z)} - V_{x,p(z)}(z) \pi(\dd z) \\
		& \quad \leq \sup_{p(z)} \int \ip{p(z)}{w^*(z)} - V_{x,p(z)}(z) \pi(\dd z) \\
		& \quad = \int \ip{p^*}{w^*(z)} - V_{x,p^*}(z) \pi(\dd z) \\
		& \quad = \ip{p^*}{v} - \int V_{x,p^*}(z) \pi(\dd z) \\
		& \quad \leq \sup_p \ip{p}{v} - \int V_{x,p}(z) \pi(\dd z)
	\end{align*}
	where we used in the third line that $w(z) \in \partial_p V_{x,p^*}(z)$ $\pi$ almost surely by~\cite[Theorem 23.5]{Rockafellar1970} or~\cite[Propposition E.1.4.3]{HiriartLemarechal2012}.
\end{proof}
\begin{proof}[Proof of Proposition~\ref{SF:proposition:representation_of_Lagrangian}]
	We have
	\begin{align}
		\cL(x,v) & = \sup_p \ip{p}{v} - \cH(x,p)  \notag \\
		& = \sup_p \inf_{\pi \in \cP(F)} \ip{p}{v} - \int V_{x,p}(z) \pi(\dd z) + I_x(\pi) \notag \\
		& = \inf_{\pi \in \cP(F)} \sup_p  \ip{p}{v} - \int V_{x,p}(z) \pi(\dd z) + I_x(\pi) \label{SF:eqn:sion_lagrangian}
	\end{align}
	by Sion's minimax lemma, since the map $\pi \mapsto \cI(x,\pi)$ is convex. Fix $\pi \in \cP(F)$, and write
	\begin{align*}
		\widehat{\cL}_\pi(x,v) = \inf_{w: \int w(z) \pi(\dd z) = v} \int \cL_z(x,w(z)) \pi(\dd z).
	\end{align*}
	By \eqref{SF:eqn:sion_lagrangian}, our proposition follows if for all $(x,v)$ and $\pi$ we have
	\begin{equation} \label{SF:eqn:equality_Lagrangians_pi}
	\cL_\pi(x,v) = \widehat{\cL}_\pi(x,v).
	\end{equation}
	Fix $(x,v)$ and $\pi$.
	\smallskip
	
	\textit{Step 1:} We establish $\cL_\pi(x,v) \leq \widehat{\cL}_\pi(x,v)$. For any integrable function $z \mapsto w(z)$ such that $\int w(z) \pi(\dd z) = v$, we have
	\begin{equation*}
		\cL_\pi(x,v) = \sup_p \int \ip{p}{w(z)} - V_{x,p}(z) \pi(\dd z),
	\end{equation*}
	implying that
	\begin{equation*}
	\cL_\pi(x,v) = \inf_{w: \int w(z) \pi(\dd z) = v} \sup_p \int \ip{p}{w(z)} - V_{x,p}(z) \pi(\dd z) \leq \widehat{\cL}_\pi(x,v) 
	\end{equation*}
	by taking the supremum over $p$ inside the integral. We conclude that $\cL_\pi(x,v) \leq \widehat{\cL}_\pi(x,v)$.
	\smallskip
	
	\textit{Step 2:} We now establish that if $v \in \text{rel. int. dom } \cL(x,\cdot)$ then $\cL_\pi(x,v) = \widehat{\cL}_\pi(x,v)$. Indeed, by Lemma \ref{SF:lemma:constant_momentum} there is a measurable function $z \mapsto w(z)$ such that $\int w(z) \pi(\dd z) = v$ and
	\begin{align*}
		\widehat{\cL}_\pi(x,v) & = \inf_{\substack{w(z) \\ \int w \dd \pi = v}} \sup_{p(z)} \int \ip{p(z)}{w(z)} - V_{x,p(z)}(z) \pi(\dd z) \\
		& = \sup_{p}   \ip{p}{v} - \int V_{x,p(z)}(z) \pi(\dd z) \\
		& = \cL_\pi(x,v).
	\end{align*}
	\textit{Step 3:} We now establish $\cL_\pi(x,v) = \widehat{\cL}_\pi(x,v)$. By step 1, we have 
	\begin{equation*}
		\text{rel. int. dom } \widehat{\cL}_\pi(x,\cdot) \subseteq \text{rel. int. dom } \cL_\pi(x,\cdot),
	\end{equation*}
	since the Lagrangians are ordered point-wise. By step 2,
	\begin{equation*}
		\text{if}\;v \in \text{rel. int. dom } \cL_\pi(x,\cdot),\quad\text{then}\quad v \in \text{dom } \cL_\pi(x,\cdot)
	\end{equation*}
	We conclude that
	\begin{equation*}
		\text{rel. int. dom } \widehat{\cL}_\pi(x,\cdot) = \text{rel. int. dom } \cL_\pi(x,\cdot),
	\end{equation*}
	and $\cL_\pi(x,\cdot) = \widehat{\cL}_\pi(x,\cdot)$ holds on this set. We conclude that $\cL_\pi(x,v) = \widehat{\cL}_\pi(x,v)$ by \cite[Corollary 7.3.4]{Rockafellar1970} (this can also be derived from \cite[Proposition B.1.2.6]{HiriartLemarechal2012}).
\end{proof}
\section{Proof of mean-field large deviations}
\label{SF:sec:proof-mean-field}
In this section, we prove Theorem~\ref{SF:thm:LDP-mean-field} by verifying the assumptions of our general large-deviation result (Theorem~\ref{SF:thm:LDP_general}) and the action-integral representation (Theorem~\ref{SF:thm:rate-function}), that is means Assumptions~\ref{SF:assumption:convergence-slow-nonlinear-generators},~\ref{SF:assumption:convergence-fast-nonlinear-generators},~\ref{SF:assumption:principal-eigenvalue-problem},~\ref{SF:assumption:regularity-V},~\ref{SF:assumption:regularity-I} and~\ref{SF:assumption:Hamiltonian_vector_field}.
\smallskip

We recall the setting: The slow-fast process~$(X^n,Z^n)$ takes values in~$E_n\times F$, where we embed~$E_n$ by identity into~$E=\mathcal{P}(\{1,\dots,q\})\times[0,\infty)^\Gamma$, and regard~$(X^n,Z^n)$ as a process on~$E\times F$. The set~$F$ is a finite-dimensional torus~$\mathbb{T}^,$. The generator of the slow-fast system is
\begin{equation*}
A_n f(x,z) := A_{n,z}^\mathrm{slow}f(\cdot,z)(x) + n\cdot  A_{n,x}^\mathrm{fast}f(x,\cdot)(z),
\end{equation*}
with slow and fast generators given by
\begin{align*}
	A_{n,z}^\mathrm{slow} g(x) &:= \sum_{ab,a\neq b} n\cdot \mu(a)\cdot r_n(a,b,\mu,z)\left[g(x_{a\to b})-g(x)\right],\\
	A_{n,x}^\mathrm{fast} h(z) &:= \sum_i b_n^i(x,z)\partial_i h(z) + \sum_{ij}a_n^{ij}(x,z)\partial_i\partial_j h(z).
\end{align*}
\begin{proof}[Verification of Assumption~\ref{SF:assumption:convergence-slow-nonlinear-generators}]
	We have to find the slow Hamiltonian~$V_{x,p}(z)$ such that
	\begin{equation*}
		\frac{1}{n}e^{-nf}A_{n,z}^\mathrm{slow}e^{nf}\xrightarrow{n\to\infty} V_{x,\nabla f(x)}(z)
	\end{equation*}
	as specified in Assumption~\ref{SF:assumption:convergence-slow-nonlinear-generators}. We have
	\begin{align*}
		\frac{1}{n}e^{-nf(x)}A_{n,z}^\mathrm{slow}e^{nf(x)}&= \sum_{ab,a\neq b}\mu(a) r_n(a,b,\mu,z) \left[\exp\{n(f(x_{a\to b}^n)-f(x))\}-1\right].
	\end{align*}
	Suppose that~$x_n=(\mu_n,w_n)\to x$. Then by Taylor expansion,
	\begin{equation*}
		n(f(x_{a\to b}^n)-f(x_n)) \xrightarrow{n\to\infty} \ip{\nabla f(x)}{e_{b}-e_a+e_{ab}},
	\end{equation*}
	for all~$f\in C^2(E)$, uniformly on compacts~$K\subseteq E$. By the convergence assumption on~$r_n$, we obtain the claimed convergence with~$D_0=C_b^2(E)$ and
	\begin{equation*}
		V_{x,p}(z) = \sum_{ab,a\neq b}\mu(a)r(a,b,\mu,z)\left[e^{p_b-p_a+p_{ab}}-1\right].
	\end{equation*}
\end{proof}
\begin{proof}[Verification of Assumption~\ref{SF:assumption:convergence-fast-nonlinear-generators}]
	Let~$x_n\to x$ in~$E$ and~$h_n\to h$ in~$C(F)$. By the convergence assumptions on the coefficients~$b_n^i$ and~$a_n^{ij}$, we obtain
	\begin{equation*}
		A_{n,x_n}^\mathrm{fast}h_n(z) \to \sum_ib^i(x,z)\partial_ih(z)+\sum_{ij}a^{ij}(x,z)\partial_i\partial_j h(z) =:A_x^\mathrm{fast}h(z)
	\end{equation*}
	uniformly over~$z \in F$.
\end{proof}
\begin{proof}[Verification of Assumption~\ref{SF:assumption:principal-eigenvalue-problem}]
	Part~\ref{SF:item:assumption:PI:domain} follows since~$F$ is compact. Let~$x\in E$ and~$p\in\mathbb{R}^d$. We find a strictly positive eigenfunction~$u:F\to(0,\infty)$ and an eigenvalue~$\mathbb{R}$ such that
	\begin{equation}\label{SF:eq:proof-mean-field-eigenvalue-problem}
		(V_{x,p} + A_x^\mathrm{fast})u=\lambda u.
	\end{equation}
	Then as remarked below Assumption~\ref{SF:assumption:principal-eigenvalue-problem}, part~\ref{SF:item:assumption:PI:solvePI} follows. Equation~\eqref{SF:eq:proof-mean-field-eigenvalue-problem} is a principal-eigenvalue problem for an uniformly elliptic operator. Uniform ellipticity follows by Condition~\ref{SF:condition:mean-field:reg-coefficients} on the diffusion coefficients and the uniform convergence in Assumption~\ref{SF:mean-field:conv-coeff}. Hence there exists a unique eigenfunction~$u$ with a real eigenvalue~$\lambda$ (e.g. Sweers~\cite{Sweers92}). By~\cite{DonskerVaradhan75}, this principal eigenvalue satisfies the variational representation
	\begin{equation*}
		\lambda = \sup_{\pi\in\mathcal{P}(F)}\left[\int_F V_{x,p}(z)\,\dd\pi(z) - \mathcal{I}(x,\pi)\right],
	\end{equation*} 
	with the functional
	\begin{equation*}
		\mathcal{I}(x,\pi)=-\inf_{\phi>0}\int_F\frac{A_x^\mathrm{fast}\phi(z)}{\phi(z)}\,\dd\pi(z).
	\end{equation*}
	Hence Assumption~\ref{SF:assumption:principal-eigenvalue-problem} holds with the Hamtilonian~\eqref{SF:eq:sec-assumption:H} as claimed.
\end{proof}
\begin{proof}[Verification of Assumptions~\ref{SF:assumption:regularity-V} and~\ref{SF:assumption:regularity-I}]
	The verification can be found in Proposition~\ref{prop:mean-field-coupled-to-diffusion} of Chapter~\ref{chapter:CP-for-two-scale-H}. The proofs are also found in~\cite[Proposition~8.2 and~8.4]{KraaijSchlottke2019}. Our example at hand is considered in Remark~8.5 in there.
\end{proof}
\begin{proof}[Verification of Assumption~\ref{SF:assumption:Hamiltonian_vector_field}]
	This follows via computation. For instance, consider~$E=\mathcal{P}(\{a,b\})$ (ignoring the flux for the moment), and identify~$E$ with the simplex in~$\mathbb{R}^2$. Fix the external variable~$z$. We have to show~$\partial_pV_{x,p}(z)\subseteq T_E(x)$. Recall that~$T_E(x)$ is the tangent cone at~$x$, that means the vectors at~$x$ pointing inside of~$E$. We compute the vector~$\nabla_p V_{x,p}(z)\in\mathbb{R}^2$,
	\begin{equation*}
	\nabla_pV_{x,p}(z) = \begin{pmatrix}
	\mu_a r(a,b,\mu,z) (-1) e^{p_b-p_a} + \mu_b r(b,a,\mu,z) e^{p_a-p_b}\\
	\mu_a r(a,b,\mu,z) e^{p_b-p_a} + \mu_b r(b,a,\mu,z) (-1)e^{p_a-p_b}
	\end{pmatrix}.
	\end{equation*}
	For~$\mu=(\mu_a,\mu_b)\in E$ with~$\mu_a,\mu_b > 0$, the tangent cone~$T_E(x)$ is spanned by~$(1,-1)^T$. Since~$\nabla_pV_{x,p}(z)$ is orthogonal to~$(1,1)^T$, we indeed find~$\partial_pV_{x,p}(z)\subseteq T_E(x)$ in that case. For~$\mu=(1,0)$, the tangent cone is~$T_E(1,0)=\{\lambda(-1,1)^T\,:\,\lambda \geq 0\}$. We have
	\begin{equation*}
	\nabla_pV_{\mu,p}(z) = \begin{pmatrix}
	r(a,b,\mu,z) (-1) e^{p_b-p_a}\\
	r(a,b,\mu,z) e^{p_b-p_a}
	\end{pmatrix},
	\end{equation*}
	which is parallel to~$(-1,1)^T$, and therefore~$\partial_pV_{\mu,p}(z)\subseteq T_E(x)$. The argument is similar for~$\mu=(0,1)$. The general case (including the fluxes) follows from writing out the definitions.
\end{proof}
\chapter{Comparison Principle for Two-Scale Hamiltonians}
\label{chapter:CP-for-two-scale-H}
\chaptermark{CP for Two-Scale Hamiltonians}
\section{Introduction and aim} 
\label{section:introduction} 
	The main purpose of this chapter is to establish well-posedness for first-order nonlinear partial differential equations of Hamilton-Jacobi-Bellman type on subsets $E$ of $\mathbb{R}^d$,
	\begin{equation}\label{eq:intro:HJ-general}
	u(x)-\lambda \,\mathcal{H}\left[u(x),\nabla u(x)\right] = h(x),\quad x \in E\subseteq \mathbb{R}^d.\tag{\text{HJB}}
	\end{equation}
	In there, $\lambda > 0$ is a scalar and  $h$ is a continuous and bounded function on~$E$. The Hamiltonian $\mathcal{H}:E\times \mathbb{R}^d\to \mathbb{R}$ is given by
	\begin{equation}\label{eq:intro:variational_hamiltonian}
	\mathcal{H}(x,p) = \sup_{\theta \in \Theta}\left[\Lambda(x,p,\theta) - \mathcal{I}(x,\theta)\right],
	\end{equation}
	where $\theta \in \Theta$ plays the role of a control variable. For fixed $\theta$, the function $\Lambda$ can frequently be interpreted as an Hamiltonian itself. We call it the \emph{internal} Hamiltonian. The function $\cI$ can be interpreted as the cost of applying the control $\theta$. This type of Hamiltonians typically arises in two-scale problems such as discussed in the previous chapter.
	
	\smallskip
	
	We will establish existence of viscosity solutions (e.g.~\cite{CIL92}) in the sense of Definition~\ref{definition:viscosity_solutions} via a resolvent defined in terms of a standard discounted control procedure. However, the main problem we overcome in this chapter is to verify a \emph{comparison principle} in order to establish uniqueness of viscosity solutions. The comparison principle for Hamilton-Jacobi equations is a well-studied problem in the literature. The standard assumption that allows one to obtain the comparison principle in the context of optimal control problems (e.g.~\cite{BardiDolcetta1997}) is that either there is a modulus of continuity $\omega$ such that
	\begin{equation}\label{eq:intro:standard-reg-estimate-on-H}
	|\mathcal{H}(x,p)-\mathcal{H}(y,p)| \leq \omega\left(|x-y|(1+|p|)\right),
	\end{equation} 
	or that $\cH$ is uniformly coercive:
	\begin{equation} \label{eqn:uniformly_coercive}
	\lim_{|p| \rightarrow \infty} \inf_x \cH(x,p) = \infty.
	\end{equation}	
	The estimate~\eqref{eq:intro:standard-reg-estimate-on-H} can be translated into conditions for~$\Lambda$ and $\cI$, which include (e.g. \cite[Chapter~III]{BardiDolcetta1997})
	\begin{itemize}
		\item $|\Lambda(x,p,\theta)-\Lambda(y,p,\theta)|\leq \omega_\Lambda(|x-y|(1+|p|))$, uniformly in $\theta$, and
		\item $\cI$ is bounded, continuous and $|\cI(x,\theta) - \cI(y,\theta)| \leq \omega_\cI(|x-y|)$.
	\end{itemize}
	However, such type of estimates are not satisfied for the examples that we are interested in. We make these examples of~\eqref{eq:intro:HJ-general} more precise in Section~\ref{section:intro:motivation}. There we also explain why the standard assumptions are not satisfied and where the challenge of solving~\eqref{eq:intro:HJ-general} is pointed out in the literature. Here, we focus on the motivation for our assumptions. They mainly build up on two observations:
	\begin{enumerate}[label=(\roman*)]
		\item Fix a control variable $\theta_0\in\Theta$ and consider the Hamiltonian $H(x,p):= \Lambda(x,p,\theta_0)$. In all our examples, the comparison principle is satisfied for sub- and supersolutions of~$u(x)-\lambda H(x,\nabla u(x))=h(x)$. 
		\item In all our examples, the cost function $\mathcal{I}(x,\theta)$ satisfies an estimate of the type $|\cI(x,\theta) - \cI(y,\theta)| \leq \omega_{\cI,C}(|x-y|)$ on sublevel sets $\{\mathcal{I} \leq C\}$.
	\end{enumerate}
	Our main idea is to take advantage of viscosity sub- and supersolution inequalities in order to work on sublevel sets of the cost function $\mathcal{I}$. To do so, we assume that $H(x,p)=\Lambda(x,p,\theta_0)$ satisfies a \emph{continuity estimate} uniformly for $\theta_0$ varying in a compact set. This continuity estimate captures the key information that allows to prove the comparison principle for~$H$. In the end, this is what we call the bootstrap principle: given sufficient regularity of $\cI$, one can bootstrap the comparison principle for the internal Hamiltonian $\Lambda$ to obtain a comparison principle for the full Hamiltonian $\mathcal{H}$. In examples, this approach proves to be a crucial improvement over known results.
	\smallskip
	
	In summary, the novelties we present in this chapter are:
	\begin{enumerate}[$\bullet$]
		\item Motivated by examples violating the standard regularity estimate~\eqref{eq:intro:standard-reg-estimate-on-H} on Hamiltonians, we find different conditions under which the comparison principle for~\eqref{eq:intro:HJ-general} is satisfied for variational Hamiltonians $\mathcal{H}$ of the type~\eqref{eq:intro:variational_hamiltonian}. The result is formulated in Theorem~\ref{theorem:comparison_principle_variational}. The main bootstrapping argument is explained in simplified form in Section~\ref{section:bootstrap-argument-in-nutshell} and carried out in Section~\ref{section:comparison_principle}.
		\item A proof of the comparison principle that covers a class of non-coercive Hamiltonians which typically arises in mean-field interacting particle systems that are coupled to external variables. This example has not been treated before, and we make it explicit in Proposition~\ref{prop:mean-field-coupled-to-diffusion} of Section~\ref{section:examples-H}.
		\item A proof of existence of a viscosity solution based on solving subdifferential inclusions in the non-compact setting. The proof relies on continuity of $\mathcal{H}$ and finding a priori estimates on the range of solutions to associated differential inclusions. The result is formulated in Theorem~\ref{theorem:existence_of_viscosity_solution}, and the structure of the proof is explained in Section~\ref{section:strategy_existence_viscosity_solution}.
	\end{enumerate}
	With these results established, we can study large deviation problems with two time-scales from a Hamilton-Jacobi point-of-view in more generality. This is the subject of Chapter~5, where we exploit the semigroup approach to large deviations. We remark that in~\cite[Lemmas~9.3,~9.19,~9.25]{FengKurtz2006}, a different technique is introduced, based on introducing an extra parameter~$\lambda$. We give further comments on that in the discussion section in Chapter~\ref{chapter:discussion}.
%
	\paragraph{Overview of this chapter.} In Section~\ref{section:intro:motivation}, we discuss Hamiltonians violating the standard regularity assumptions. The main results are formulated in Section \ref{section:results}. We proceed with a discussion of the strategy of the proofs in Section \ref{section:strategy}. In Section \ref{section:regularity-of-H-and-L} we establish regularity properties of $\cH$ used in the later proof sections. In Section \ref{section:comparison_principle} we establish the comparison principle. In Section \ref{section:construction-of-viscosity-solutions} we establish that a resolvent operator $R(\lambda)$ in terms of an exponentially discounted control problem gives rise to viscosity solutions of the  Hamilton-Jacobi-Bellman equation~\eqref{eq:intro:HJ-general}. Finally, in Section~\ref{section:verification-for-examples-of-Hamiltonians} we verify the assumptions for examples.
	\section{Examples violating the standard assumptions}
	\label{section:intro:motivation}
	Hamiltonians of the type~\eqref{eq:intro:variational_hamiltonian} arise in a range of fields. In this section, we mention two examples of Hamiltonians arising in the context of stochastic systems with two time scales. We explain why they violate the standard regularity estimates. These examples illustrate the need for an alternative set of assumptions allowing to treat these cases. These Hamiltonians frequently arise in the study of systems with multiple time-scales, e.g. geophysical flows, planetary motion, finance, weather-climate interaction models, molecular dynamics and models in statistical physics---we provide more background in Chapter~5.
    In such systems, one can often recognize a slow and a fast component. Typically, one is interested in the behaviour of the slow component in the limit in which the separation of time scales goes to infinity. As the fast system equilibrates before the slow system has made a significant difference, the limit of such systems can be described by a ordinary or partial differential equation involving only the average behaviour of the fast component.
	\smallskip
	
	However, in applications an infinite separation of time scales is never really achieved. Thus, the slow process still shows fluctuations around its limiting behaviour while the fast process fluctuates around its average. The effective fluctuations arise from the combination of both sources. In this two-scale context, when analysing the fluctuations by means of large-deviation techniques, one obtains Hamiltonians of the type~\eqref{eq:intro:variational_hamiltonian}. We refer to~\cite{KumarPopovic2017} for derivations in this context, and to~\cite{BouchetGrafkeTangarifeVandenEijnden2016} for an extensive explanation in which the authors study ODE's coupled to fast diffusion. In these examples, the internal Hamiltonians~$\Lambda$ capture the fluctuations of the slow component, while the cost function~$\mathcal{I}$ arises from fluctuations of averages of the fast component. The full Hamiltonian~$\mathcal{H}$ takes both contributions into account.
	\smallskip
	
	\textit{Example 1}. In~\cite{BudhirajaDupuisGanguly2018}, the authors study large deviations of a diffusion processes with vanishing noise on $E = \bR^d$ coupled to a fast jump process on a finite discrete set $\{1,\dots,J\}$. They identified the challenge of proving comparison principles for Hamiltonians arising in such two-scale systems, where the Hamiltonians can be casted in the form~\eqref{eq:intro:variational_hamiltonian}. We consider this general setting in Proposition~\ref{prop:diffusion-coupled-to-jumps} in Section~\ref{section:examples-H}. We illustrate the issues arising in a simpler but more concrete form. With $d=1$ and $J=2$, when approaching this problem from the Hamilton-Jacobi perspective, a key step (e.g.~\cite{KumarPopovic2017}) is to solve~\eqref{eq:intro:HJ-general} with $\cH$ consisting of the following ingredients: 
	\begin{enumerate}[label=(\roman*)]
		\item The internal state space is $E=\mathbb{R}^d$.
		\item The set of control variables is $\Theta=\mathcal{P}(\{1,2\})$.
		\item The internal Hamiltonian~$\Lambda$ is given by
		\begin{equation*}
		\Lambda(x,p,\theta) = \frac{1}{2}a(x,1)|p|^2 \theta_1 + \frac{1}{2} a(x,2)|p|^2 \theta_2,
		\end{equation*}
		where $a(x,i) > 0$ and $\theta_i = \theta(\{i\})$.
		\item The cost function~$\mathcal{I}$ is given by
		\begin{equation*}
		\mathcal{I}(x,\theta) = \sup_{w\in\mathbb{R}^2} \left[r_{12}(x)\theta_1 \left(1-e^{w_2-w_1}\right) + r_{21}(x)\theta_2 (1-e^{w_1-w_2})\right],
		\end{equation*}
		where $r_{ij}(x) \geq 0$.
	\end{enumerate}
	In this example, the cost function is unbounded if $r_{ij}(x)$ is unbounded. For instance, consider $\theta_1 = 1$ and $\theta_2 = 0$. Then by choosing $w=(1,0)$ in the supremum,
	\begin{equation*}
	\mathcal{I}(x,\theta) \geq C \, r_{12}(x),
	\end{equation*}
	and thus $\mathcal{I}(x,\theta)$ diverges as $|x|\to\infty$.
	\smallskip
	
	We now turn to another notable problem with two time-scales that motivates our considerations: a system of mean-field interacting particles coupled to fast external variables.
	\smallskip
	
	\textit{Example 2}. In~\cite{BertiniChetriteFaggionatoGabrielli2018}, the authors prove large-deviation principles of mean-field interacting particles that are coupled to fast time-periodic variables. In this setting, the associated Hamilton-Jacobi equations are solved in~\cite{Kr17}. However, when considering a coupling to general fast random variables such as diffusions, then solving the corresponding Hamilton-Jacobi equations remained an open challenge. In full generality, we formulate this case in Proposition~\ref{prop:mean-field-coupled-to-diffusion}. For a corresponding large-deviation analysis, we refer to Chapter~\ref{chapter:LDP-in-slow-fast-systems} (based on~\cite{KraaijSchlottke2020}). Here we illustrate the difficulties that arise by considering the Hamiltonian in a simplified setting:
	\begin{enumerate}[label=(\roman*)]
		\item The internal state space is $E=\mathcal{P}(\{a,b\}) \times [0,\infty)\times[0,\infty)$, embedded in~$\mathbb{R}^4$ by identifying~$\mathcal{P}(\{a,b\})$ with the simplex in~$\mathbb{R}^2$. We denote the variables as $x = (\mu,w)$, with $\mu\in\mathcal{P}(\{a,b\})$ and $w\in[0,\infty)^2$.
		\item The set of control variables is $\Theta=\mathcal{P}(\mathbb{T})$, that is the probability measures on the torus~$\mathbb{T}$.
		\item The internal Hamiltonian~$\Lambda$ is given by
		\begin{multline*}
		\Lambda(x,p,\theta) = \mu_a r_{ab}(\mu,\theta)\left[\exp\left\{p_b-p_a + p_{ab}\right\}-1\right] \\ + \mu_b r_{ba}(\mu,\theta)\left[\exp\left\{p_a-p_b + p_{ba}\right\}-1\right],
		\end{multline*}
		with $p = (p_a,p_b,p_{ab},p_{ba}) \in \mathbb{R}^4$ and $\mu_i := \mu(\{i\})$. The rates $r_{ij}$ are non-negative.
		\item The cost function~$\mathcal{I}:\Theta\to[0,\infty]$ is independent of $x$ and is given by
		\begin{equation*}
		\mathcal{I}(\theta) = \sup_{\substack{u\in C^\infty(\mathbb{T})\\ u > 0}} \int_\mathbb{T} \left(-\frac{u''(y)}{u(y)}\right)\,\dd\theta(y)
		\end{equation*}
	\end{enumerate}
	In this example, the internal Hamiltonian~$\Lambda$ is not uniformly coercive. For instance, take momenta $p$ such that $p_b-p_a+p_{ab}$ is constant. Then if $|p|\to\infty$, we do not necessarily have that $\Lambda(x,p,\theta)\to\infty$. A similar effect occurs when choosing $p_a \rightarrow \infty$ and $\mu_a = 0$. Regarding the cost function, for any singular measure $\delta_z$ with a point $z\in S$ we have $\mathcal{I}(\delta_z) = \infty$. This similarly holds for finite convex combinations of Dirac measures. Since this linear span is dense in $\mathcal{P}(\mathbb{T})$, this implies that $\mathcal{I}$ can not be continuous.
	\section{Main results}
	\label{section:results}
	In this section, we start with preliminaries in Section~\ref{section:preliminaries} which includes the definition of viscosity solutions and that of the comparison principle. 
	\smallskip
	
	We proceed in Section~\ref{section:results:HJ-of-Perron-Frobenius-type} with the main results: a comparison principle for the Hamilton-Jacobi-Bellman equation~\eqref{eq:intro:HJ-general} based on variational Hamiltonians of the form~\eqref{eq:intro:variational_hamiltonian}, and the existence of viscosity solutions. 
	
	\smallskip
	
	In Section \ref{section:assumptions} we collect all the assumptions that are needed for all main results in one place and discuss the applicability of our results. In Section~\ref{section:examples-H}, we verify the assumptions for the examples that motivate the Hamilton-Jacobi equations we discuss in this chapter.
	
	\subsection{Preliminaries} \label{section:preliminaries}
	For a Polish space $\cX$, we denote by $C(\cX)$ and $C_b(\cX)$ the spaces of continuous and bounded continuous functions respectively. If $\cX \subseteq \bR^d$ then we denote by $C_c^\infty(\cX)$ the space of smooth functions that vanish outside a compact set. We denote by $C_{cc}^\infty(\cX)$ the set of smooth functions that are constant outside of a compact set, and by $\cP(\cX)$ the space of probability measures on $\cX$. We equip $\cP(\cX)$ with the weak topology, that is, the one induced by convergence of integrals against bounded continuous functions.
	
	\smallskip
	
	Throughout this chapter, $E$ will be the set on which we base our Hamilton-Jacobi equations. We assume that $E$ is a subset of $\bR^d$ that is a Polish space which is contained in the $\bR^d$ closure of its $\bR^d$ interior. This ensures that gradients of functions are determined by their values on $E$. Note that we do not assume that $E$ is open. We assume that the space of controls $\Theta$ is Polish.
	\smallskip
	
	We next introduce viscosity solutions for the Hamilton-Jacobi equation with Hamiltonians like $\mathcal{H}(x,p)$ of our introduction.
	\begin{definition}[Viscosity solutions and comparison principle] \label{definition:viscosity_solutions}
		Let $A : C_b(E) \to C_b(E)$ be an operator with domain $\mathcal{D}(A)$, $\lambda > 0$ and $h \in C_b(E)$. Consider the Hamilton-Jacobi equation
		\begin{equation}
		f - \lambda A f = h. \label{eqn:differential_equation} 
		\end{equation}
		We say that $u$ is a \textit{(viscosity) subsolution} of equation \eqref{eqn:differential_equation} if $u$ is bounded, upper semi-continuous and if, for every $f \in \cD(A)$ there exists a sequence $x_n \in E$ such that
		\begin{gather*}
		\lim_{n \uparrow \infty} u(x_n) - f(x_n)  = \sup_x u(x) - f(x), \\
		\lim_{n \uparrow \infty} u(x_n) - \lambda A f(x_n) - h(x_n) \leq 0.
		\end{gather*}
		We say that $v$ is a \textit{(viscosity) supersolution} of equation \eqref{eqn:differential_equation} if $v$ is bounded, lower semi-continuous and if, for every $f \in \cD(H)$ there exists a sequence $x_n \in E$ such that
		\begin{gather*}
		\lim_{n \uparrow \infty} v(x_n) - f(x_n)  = \inf_x v(x) - f(x), \\
		\lim_{n \uparrow \infty} v(x_n) - \lambda Af(x_n) - h(x_n) \geq 0.
		\end{gather*}
		We say that $u$ is a \textit{(viscosity) solution} of equation \eqref{eqn:differential_equation} if it is both a subsolution and a supersolution to \eqref{eqn:differential_equation}.
		
		We say that \eqref{eqn:differential_equation} satisfies the \textit{comparison principle} if for every subsolution $u$ and supersolution $v$ to \eqref{eqn:differential_equation}, we have $u \leq v$.
	\end{definition}
	\begin{remark} \label{remark:existence of optimizers}
		Consider the definition of subsolutions. Suppose that the testfunction $f \in \cD(A)$ has compact sublevel sets, then instead of working with a sequence $x_n$, there exists $x_0  \in E$ such that
		\begin{gather*}
		u(x_0) - f(x_0)  = \sup_x u(x) - f(x), \\
		u(x_0) - \lambda A f(x_0) - h(x_0) \leq 0.
		\end{gather*}
		A similar simplification holds in the case of supersolutions. For an explanatory text on the notion of viscosity solutions and fields of applications, we refer to~\cite{CIL92}.
	\end{remark}
	\subsection{Hamilton-Jacobi-Bellman Equations}
	\label{section:results:HJ-of-Perron-Frobenius-type}
	In this Section, we state our main results, which are the comparison principle (Theorem~\ref{theorem:comparison_principle_variational}) and the existence of solutions (Theorem \ref{theorem:existence_of_viscosity_solution}).
	We consider the variational Hamiltonian $\cH : E \times \bR^d \rightarrow  \bR$ given by
	\begin{equation}\label{eq:results:variational_hamiltonian}
	\mathcal{H}(x,p) = \sup_{\theta \in \Theta}\left[\Lambda(x,p,\theta) - \mathcal{I}(x,\theta)\right].
	\end{equation}
	The precise assumptions on the maps $\Lambda$ and $\mathcal{I}$ are formulated in Section~\ref{section:assumptions}. Define the operator $\bfH f(x) := \cH(x,\nabla f(x))$ with domain $\cD(\bfH) = C_{cc}^\infty(E)$. Our first main result is that this operator $\bfH$ constructed out of $\cH$ satisfies the comparison principle.
	\begin{theorem}[Comparison principle]\label{theorem:comparison_principle_variational}
		Suppose that the maps~$\Lambda$ and~$\mathcal{I}$ satisfy Assumptions~\ref{assumption:results:regularity_of_V} and~\ref{assumption:results:regularity_I}, respectively.
		Then for any $h \in C_b(E)$ and $\lambda > 0$, the comparison principle holds for
		\begin{equation}\label{eq:results:HJ-eq}
		f - \lambda \, \bfH f = h.
		\end{equation}
	\end{theorem}
	\begin{remark}[Uniqueness]
		If $u$ and $v$ are two viscosity solutions of~\ref{eq:results:HJ-eq}, then we have $u\leq v$ and $v\leq u$ by the comparison principle, giving uniqueness.
	\end{remark}
	\begin{remark}[Domain]
		The comparison principle holds with any domain that satisfies $C_{cc}^\infty(E)\subseteq \mathcal{D}(\mathbf{H})\subseteq C^1_b(E)$. We state it with $C^\infty_{cc}(E)$ to connect it with the existence result of Theorem~\ref{theorem:existence_of_viscosity_solution}, where we need to work with test functions whose gradients have compact support.
	\end{remark}
	We turn to the existence of a viscosity solution for~\eqref{eq:results:HJ-eq}. As mentioned in the introduction, the viscosity solution is given in terms of an optimization problem with discounted cost. The Legendre dual $\cL : E \times \bR^d \rightarrow [0,\infty]$ of $\cH$, given by
	\begin{equation*}
	\cL(x,v) := \sup_{p\in\mathbb{R}^d} \left[\ip{p}{v} - \cH(x,p)\right],
	\end{equation*}
	plays the role of a running cost. In the following Theorem, $\cA\cC$ is the collection of absolutely continuous paths in $E$. For each $\lambda > 0$, let $R(\lambda)$ be the operator
	\begin{equation}\label{HJB:eq:resolvent}
	R(\lambda) h(x) = \sup_{\substack{\gamma \in \mathcal{A}\mathcal{C}\\ \gamma(0) = x}} \int_0^\infty \lambda^{-1} e^{-\lambda^{-1}t} \left[h(\gamma(t)) - \int_0^t \mathcal{L}(\gamma(s),\dot{\gamma}(s))\right] \, \dd t.
	\end{equation}
	\begin{theorem}[Existence of viscosity solution] \label{theorem:existence_of_viscosity_solution}
		Suppose that~$\Lambda$ and $\mathcal{I}$ satisfy Assumptions~\ref{assumption:results:regularity_of_V} and~\ref{assumption:results:regularity_I}, respectively, and that $\mathcal{H}$ satisfies Assumption~\ref{assumption:Hamiltonian_vector_field}.
		Then the function~$R(\lambda)h$ is the unique viscosity solution to~\eqref{eq:results:HJ-eq}.
	\end{theorem}

	\subsection{Assumptions} \label{section:assumptions}
	In this section, we formulate and comment on the assumptions imposed on the Hamiltonians defined in the previous sections. We first motivate the assumptions that are required for proving the comparison principle, Theorem~\ref{theorem:comparison_principle_variational}.
	
	\smallskip
	
	Usually, proofs of the comparison principle for a subsolution $u$ and a supersolution $v$ for the equation $f - \lambda \bfH f = h$ are reduced to establishing an estimate of the type
	\begin{equation*}
	\liminf_{\varepsilon \downarrow 0} \liminf_{\alpha \rightarrow \infty}\cH(x_{\alpha,\varepsilon},\alpha(x_{\alpha,\varepsilon}-y_{\alpha,\varepsilon})) - \cH(y_{\alpha,\varepsilon},\alpha (x_{\alpha,\varepsilon}-y_{\alpha,\varepsilon})) \leq 0
	\end{equation*}
	where $(x_{\alpha,\varepsilon},y_{\alpha,\varepsilon})$ are elements of $E$ such that
	\begin{multline} \label{eqn:heuristic:doubling_variable}
	u(x_{\alpha,\varepsilon}) - v(x_{\alpha,\varepsilon}) - \frac{\alpha}{2} |x_{\alpha,\varepsilon}-y_{\alpha,\varepsilon}|^2 -  \frac{\varepsilon}{2} (|x_{\alpha,\varepsilon}|^2 + |y_{\alpha,\varepsilon}|^2) \\
	= \sup_{x,y \in E} u(x) - v(y) - \frac{\alpha}{2} |x-y|^2 - \frac{\varepsilon}{2} (|x|^2 + |y|^2).
	\end{multline}
	Equation~\eqref{eqn:heuristic:doubling_variable}, together with the sub- and supersolution property of $u$ and $v$ respectively, has the following consequences:
	\begin{enumerate}[(1)]
		\item \label{item:relative_compact} For all $\varepsilon >0 $, the set $\{x_{\alpha,\varepsilon}, y_{\alpha,\varepsilon} \, | \, \alpha > 0\}$ is relatively compact in $E$;
		\item \label{item:distance} For all $\varepsilon >0 $, we have $|x_{\alpha,\varepsilon} - y_{\alpha,\varepsilon}| + \alpha |x_{\alpha,\varepsilon} - y_{\alpha,\varepsilon}|^2 \rightarrow 0$ as $\alpha \rightarrow \infty$;
		\item \label{item:coercivity} For all $\varepsilon > 0$, we have
		\[
		 \inf_{\alpha} \cH(x_{\alpha,\varepsilon}, \alpha(x_{\alpha,\varepsilon} - y_{\alpha,\varepsilon})) > - \infty
		 \quad\text{and}\quad
		  \sup_{\alpha} \cH(y_{\alpha,\varepsilon}, \alpha(x_{\alpha,\varepsilon} - y_{\alpha,\varepsilon})) < \infty.
		\]
	\end{enumerate}
	In our bootstrap procedure, we aim to lift the comparison principle that holds for the Hamilton-Jacobi equation in terms of $\Lambda$ to that for $\bfH$. Thus, we need to establish an estimate of the type \eqref{eqn:heuristic:doubling_variable} under assumptions of the type \ref{item:relative_compact}, \ref{item:distance} and \ref{item:coercivity} where in addition, we have to vary our control variable~$\theta$. It turns out that it suffices to vary $\theta$ in a compact set in $\Theta$ that depends on $\varepsilon$. In addition, to make sure that we can bootstrap, we have to relax the $\sup$ and $\inf$ in \ref{item:coercivity} to a $\limsup$ and $\liminf$. 
	\smallskip
	
	To establish the comparison principle, the quadratic distance is not special, except for being symmetric and well suited for quadratic Hamiltonians.
	We will work with a general non-negative function~$\Psi$ to penalize the distance between~$x$ and~$y$, and use a function~$\Upsilon$ to penalize points~$x$ and~$y$ far away from the `origin'.
	\begin{definition}[Penalization function]\label{def:results:good_penalization_function}
		We say that $\Psi : E^2 \rightarrow [0,\infty)$ is a \textit{ penalization function} if $\Psi \in C^1(E^2)$ and if $x = y$ if and only if $\Psi(x,y) = 0$.
	\end{definition}
	\begin{definition}[Containment function]\label{def:results:compact-containment}
		We say that a function $\Upsilon : E \rightarrow [0,\infty]$ is a \textit{containment function} for $\Lambda$ if there is a constant $c_\Upsilon$ such that
		\begin{itemize}
			\item For every $c \geq 0$, the set $\{x \, | \, \Upsilon(x) \leq c\}$ is compact;
			\item We have $\sup_\theta\sup_x \Lambda\left(x,\nabla \Upsilon(x),\theta\right) \leq c_\Upsilon$.
		\end{itemize}	
	\end{definition}
	\begin{definition}[Continuity estimate] \label{def:results:continuity_estimate}
		Let  $\Psi$ be a penalization function and let $\cG: E \times \mathbb{R}^d\times\Theta : (x,p,\theta)\mapsto \cG(x,p,\theta)$ be a function. Suppose that for $\varepsilon > 0$ and $\alpha > 0$, we have a collection of variables $(x_{\varepsilon,\alpha},y_{\varepsilon,\alpha})$ in $E^2$ and variables $\theta_{\varepsilon,\alpha}$ in $\Theta$. We say that this collection is \textit{fundamental} for $\cG$ with respect to $\Psi$ if:
		\begin{enumerate}[label = (C\arabic*)]
			\item \label{item:def:continuity_estimate:1} For each $\varepsilon$, there are compact sets $K_\varepsilon \subseteq E$ and $\widehat{K}_\varepsilon\subseteq\Theta$ such that for all $\alpha$ we have $x_{\varepsilon,\alpha},y_{\varepsilon,\alpha} \in K_\varepsilon$ and $\theta_{\varepsilon,\alpha}\in\widehat{K}_\varepsilon$.
			\item \label{item:def:continuity_estimate:2} For each $\varepsilon > 0$, we have limit points $x_{\varepsilon} \in K_\varepsilon$ and $y_{\varepsilon} \in K_\varepsilon$ of $x_{\alpha,\varepsilon}$ and $y_{\alpha,\varepsilon}$ as $\alpha\rightarrow \infty$. For these limit points we have 
			\begin{equation*}
			\lim_{\alpha \rightarrow \infty} \alpha \Psi(x_{\alpha,\varepsilon},y_{\alpha,\varepsilon}) = 0, \qquad 	\Psi(x_{\varepsilon},y_{\varepsilon}) = 0.
			\end{equation*}
			\item \label{item:def:continuity_estimate:3} We have
			\begin{align} 
			& \limsup_{\varepsilon \rightarrow 0} \limsup_{\alpha \rightarrow \infty} \cG\left(y_{\alpha,\varepsilon}, - \alpha (\nabla \Psi(x_{\alpha,\varepsilon},\cdot))(y_{\alpha,\varepsilon}),\theta_{\varepsilon,\alpha}\right) < \infty, \label{eqn:control_on_Gbasic_sup} \\
			& \liminf_{\varepsilon \rightarrow 0} \liminf_{\alpha \rightarrow \infty} \cG\left(x_{\alpha,\varepsilon}, \alpha (\nabla \Psi(\cdot,y_{\alpha,\varepsilon}))(x_{\alpha,\varepsilon}),\theta_{\varepsilon,\alpha}\right) > - \infty. \label{eqn:control_on_Gbasic_inf} 	
			\end{align} \label{itemize:funamental_inequality_control_upper_bound}
			In other words, the operator $\cG$ evaluated in the proper momenta is eventually bounded from above and from below.
		\end{enumerate}
		We say that $\cG$ satisfies the \textit{continuity estimate} if for every fundamental collection of variables we have
		\begin{multline}\label{equation:Xi_negative_liminf}
		\liminf_{\varepsilon \rightarrow 0} \liminf_{\alpha \rightarrow \infty} \cG\left(x_{\alpha,\varepsilon}, \alpha (\nabla \Psi(\cdot,y_{\alpha,\varepsilon}))(x_{\alpha,\varepsilon}),\theta_{\varepsilon,\alpha}\right) \\
		- \cG\left(y_{\alpha,\varepsilon}, - \alpha (\nabla \Psi(x_{\alpha,\varepsilon},\cdot))(y_{\alpha,\varepsilon}),\theta_{\varepsilon,\alpha}\right) \leq 0.
		\end{multline}
	\end{definition}
	The continuity estimate is a sensible notion because it is satisfied in a wide range of examples. Indeed, all our examples satisfy the continuity estimate., and in Section \ref{section:verification_of_continuity_estimate}, we verify the continuity estimate in three different contexts. In the appendix of~\cite{KraaijSchlottke2019}, we state a slightly more general continuity estimate on the basis of multiple penalization functions. For the first reading of the proofs below, the use of this more general setting would be distracting. We want to mention, however, that all arguments below can be carried out on the basis of this more elaborate continuity estimate. Following~\cite{Kr17} a continuity estimate of this more elaborate type can be established in the context of Markov jump processes and their fluxes.
	\smallskip
	
	Our first assumption essentially states that we can solve the comparison principle for the Hamilton-Jacobi equation for $\Lambda$ uniformly over compact sets in $\Theta$. In addition to this assumption, we assume \ref{item:assumption:slow_regularity:controlled_growth}  which states the function $\Lambda$ grows roughly equally fast in $p$ for different control variables.
	\begin{assumption}\label{assumption:results:regularity_of_V}
		The function $\Lambda:E\times\mathbb{R}^d\times\Theta\to\mathbb{R}$ in the Hamiltonian~\eqref{eq:results:variational_hamiltonian} satisfies the following.
		\begin{enumerate}[label=($\Lambda$\arabic*)]
			\item \label{item:assumption:slow_regularity:continuity} The map $\Lambda : E\times\mathbb{R}^d\times\Theta\to\mathbb{R}$ is continuous and for any $(x,p)$, we have boundedness: $\|\Lambda(x,p,\cdot)\|_{\Theta}:= \sup_{\theta\in\Theta}|\Lambda(x,p,\theta)| < \infty$.
			\item \label{item:assumption:slow_regularity:convexity} For any $x\in E$ and $\theta\in\Theta$, the map $p\mapsto \Lambda(x,p,\theta)$ is convex. For $p_0=0$, we have $\Lambda(x,p_0,\theta) = 0$ for all $x\in E$ and all $\theta \in \Theta$.
			\item \label{item:assumption:slow_regularity:compact_containment} There exists a containment function $\Upsilon : E \to [0,\infty)$ in the sense of Definition~\ref{def:results:compact-containment}.
			\item \label{item:assumption:slow_regularity:continuity_estimate} The function $\Lambda$ satisfies the continuity estimate.
			\item \label{item:assumption:slow_regularity:controlled_growth} 
			For every compact set $K \subseteq E$, there exist constants $M, C_1, C_2 \geq 0$  such that for all $x \in K$, $p \in \mathbb{R}^d$ and all $\theta_1,\theta_2\in\Theta$, we have 
			\begin{equation*}
			\Lambda(x,p,\theta_1) \leq \max\left\{M,C_1 \Lambda(x,p,\theta_2) + C_2\right\}.
			\end{equation*}
		\end{enumerate} 
	\end{assumption}
	Our next assumption is on the regularity of the cost functional $\cI$. They are satisfied for continuous and bounded $\cI$ and $\Theta$ a compact space. 
	\begin{assumption}\label{assumption:results:regularity_I}
		The functional $\mathcal{I}:E\times\Theta \to [0,\infty]$ in~\eqref{eq:results:variational_hamiltonian} satisfies the following.
		\begin{enumerate}[label=($\mathcal{I}$\arabic*)]
			\item \label{item:assumption:I:lsc} The map $(x,\theta) \mapsto \mathcal{I}(x,\theta)$ is lower semi-continuous on $E \times \Theta$.
			\item \label{item:assumption:I:zero-measure} For any $x\in E$, there exists a point $\theta_x\in\Theta$ such that $\mathcal{I}(x,\theta_x) = 0$. 
			\item \label{item:assumption:I:compact-sublevelsets} For any $x \in E$, compact set $K \subseteq E$ and $C \geq 0$ the set $\left\{\theta \in \Theta \, \middle| \cI(x,\theta) \leq C\right\}$ is compact and $\cup_{x\in K}\left\{\theta\in\Theta \, \middle| \, \mathcal{I}(x,\theta) \leq C\right\}$ is relatively compact.
			\item \label{item:assumption:I:finiteness} For any converging sequence $x_n \to x$ in $E$ and sequence $\theta_n \in \Theta$, if there is an $M > 0$ such that $\mathcal{I}(x_n,\theta_n) \leq M < \infty$ for all $n \in \mathbb{N}$, then there exists a neighborhood $U_x$ of $x$ and a constant $M' > 0$ such that for any $y \in U_x$ and $n \in \mathbb{N}$,
			\begin{equation*}
			\mathcal{I}(y,\theta_n) \leq M' < \infty.
			\end{equation*}
			\item \label{item:assumption:I:equi-cont} For every compact set $K \subseteq E$ and each $M \geq 0$ the collection of functions $\{\cI(\cdot,\theta)\}_{\theta \in \Theta_M}$ with
			\begin{equation*}
			\Theta_{M} := \left\{\theta\in\Theta \, \middle| \, \forall \, x \in K: \, \mathcal{I}(x,\theta) \leq M \right\}
			\end{equation*}
			is equicontinuous. That is: for all $\varepsilon > 0$, there is a $\delta > 0$ such that for all $\theta \in \Theta_M$ and $x,y \in K$ such that $d(x,y) \leq \delta$ we have $|\mathcal{I}(x,\theta) - \mathcal{I}(y,\theta)| \leq \varepsilon$.
		\end{enumerate}
	\end{assumption}
	\begin{remark}[Gamma-convergence]
		The assumptions on $\mathcal{I}$ imply that for any sequence $x_n\to x$ in $E$, the functionals defined by $\mathcal{I}_n(\theta):=\mathcal{I}(x_n,\theta)$ $\Gamma$-converge to $\mathcal{I}_\infty$ defined by $\mathcal{I}_\infty(\theta) := \mathcal{I}(x,\theta)$. We give a proof in Proposition \ref{prop:Gamma-convergence-of-I} below.
	\end{remark}
	We turn to Theorem \ref{theorem:existence_of_viscosity_solution}. A key ingredient in establishing the existence of a viscosity solution to Hamilton-Jacobi equations is the existence of `optimally' controlled paths. The optimal controls can, for continuously differentiable Hamiltonians, be found from the Hamiltonian flow. In our context, $\cH$ is not continuously differentiable. We will show in Proposition \ref{prop:reg-of-H-and-L:reg-H}, however, that $\cH$ is convex in $p$. We can therefore define the subdifferential set
	\begin{multline} \label{eqn:subdifferential_H}
	\partial_p \mathcal{H}(x_0,p_0) \\
	:= \left\{
	\xi \in \mathbb{R}^d \,:\, \mathcal{H}(x_0,p) \geq \mathcal{H}(x_0,p_0) + \xi \cdot (p-p_0) \quad (\forall p \in \mathbb{R}^d)
	\right\}.
	\end{multline}
	Instead using solutions arising from the differential equation arising from the gradient of $\cH$, we will use solutions to differential inclusions arising from $\partial_p \cH$. As our set $E$ is not necessarily equal to $\bR^d$, but could be, e.g. a domain with corners like $[0,\infty)^d$, we need some conditions to make sure that the solutions to our differential inclusions remain within $E$. Assumption \ref{assumption:Hamiltonian_vector_field} below will make sure that the Hamiltonian vector field points `inside' $E$.	
	\begin{definition} \label{definition:tangent_cone}
		The tangent cone (sometimes also called \textit{Bouligand cotingent cone}) to $E$ in $\bR^d$ at $x$ is
		\begin{equation*}
		T_E(x) := \left\{z \in \bR^d \, \middle| \, \liminf_{\lambda \downarrow 0} \frac{d(y + \lambda z, E)}{\lambda} = 0\right\}.
		\end{equation*}
	\end{definition}
	\begin{assumption} \label{assumption:Hamiltonian_vector_field}
		The map $\cH : E \times \bR^d \rightarrow \bR$ defined in~\eqref{eq:results:variational_hamiltonian} is such that $\partial_p \cH(x,p) \subseteq T_E(x)$ for all $p$.
	\end{assumption}
		Assumption \ref{assumption:Hamiltonian_vector_field} is intuitively implied by the comparison principle for $\bfH$. We therefore expect Assumption \ref{assumption:Hamiltonian_vector_field} to be satisfied in any situation in which Theorem \ref{theorem:comparison_principle_variational} holds. We argue in a simple case why this is to be expected. The main reason is that the comparison principle for $\bfH$ builds upon the \emph{maximum principle}.  
		\smallskip
		
		Let $E = [0,1]$ and $f,g \in C^1_b(E)$. Suppose that $f(0) - g(0) = \sup_x f(x) - g(x)$. As $0$ is a boundary point, we conclude that $f'(0) \leq g'(0)$. If the maximum principle holds, we must have
		\begin{equation*}
		\cH(0,f'(0)) = Hf(0) \leq Hg(0) = \cH(0,g'(0)).
		\end{equation*}
		Hence the map $p \mapsto \cH(0,p)$ is increasing, which just means 
		\begin{equation*}
		\partial_p \cH(x,p)) \subseteq [0,\infty) = T_{[0,1]}(0).
		\end{equation*}
	\subsection{Examples of Hamiltonians}
	\label{section:examples-H}
	The purpose of this section is to demonstrate via different examples that the method introduced is versatile enough to capture a variety of interesting examples. Propositions~\ref{prop:diffusion-coupled-to-jumps} and~\ref{prop:diffusion-coupled-to-diffusion} correspond to the Hamiltonian that one encounters in two-scale systems as studied in~\cite{BudhirajaDupuisGanguly2018,KumarPopovic2017}. The example of Proposition~\ref{prop:mean-field-coupled-to-diffusion} arises in models of mean-field interacting particles that are coupled to fast external variables, whose Hamiltonians can not be treated with standard methods. Recall the form of the Hamiltonian,
	\begin{equation}\label{eq:results:H-example-section}
	\mathcal{H}(x,p) = \sup_{\theta \in \Theta}\left[\Lambda(x,p,\theta) - \mathcal{I}(x,\theta)\right],
	\end{equation} 
	Each definition below corresponds to a specification of the elements involved in~\eqref{eq:results:H-example-section}. 
	All propositions are proven in Section~\ref{section:verification-for-examples-of-Hamiltonians}, by verifying the general Assumptions~\ref{assumption:results:regularity_of_V} and~\ref{assumption:results:regularity_I} on the functions~$\Lambda$ and~$\mathcal{I}$. Let us start with Hamiltonians arising from a diffusion process coupled to a fast jump process.
	\begin{proposition}[Diffusion coupled to jumps]\label{prop:diffusion-coupled-to-jumps}
		Let $E=\mathbb{R}^d$ and $F=\{1,\dots,J\}$ be a finite set. Suppose the following.
		\begin{enumerate}[label=(\roman*)]
			\item The set of control variables is $\Theta:=\mathcal{P}(\{1,\dots,J\})$, that is probability measures over the finite set $F$.
			\item The function $\Lambda$ is given by
			\begin{equation*} 
			\Lambda(x,p,\theta) := \sum_{i\in F}\left[\ip{a(x,i)p}{p}+\ip{b(x,i)}{p}\right]\theta_i,
			\end{equation*}
			where $a:E\times F\to\mathbb{R}^{d\times d}$ and $b:E\times F\to\mathbb{R}^d$, and $\theta_i:=\theta(\{i\})$.
			\item The cost function $\mathcal{I}:E\times\Theta\to[0,\infty)$ is given by
			\begin{equation*}
			\mathcal{I}(x,\theta) := \sup_{w\in\mathbb{R}^J}\sum_{ij}r(i,j,x)\theta_i \left[1-e^{w_j-w_i}\right],
			\end{equation*}
			with non-negative rates $r:F^2\times E\to[0,\infty)$.
		\end{enumerate}
		Suppose that the cost function~$\mathcal{I}$ satisfies the assumptions of Proposition~\ref{prop:verify:DV-for-Jumps} and the function~$\Lambda$ satisfies the assumptions of Proposition~\ref{prop:verify-ex:Lambda_quadratic}. Then Theorems~\ref{theorem:comparison_principle_variational} and~\ref{theorem:existence_of_viscosity_solution} apply to the Hamiltonian~\eqref{eq:results:H-example-section}.
	\end{proposition}
		Under irreducibility conditions on the rates, as assumed below in Proposition~\ref{prop:verify:DV-for-Jumps}, by~\cite{DonskerVaradhan75} the Hamiltonian~$\mathcal{H}(x,p)$ is the principal eigenvalue of the matrix $A_{x,p} \in \mathrm{Mat}_{J \times J}(\mathbb{R})$ given by
		\[
		A_{x,p} = \mathrm{diag}\left[\ip{a(x,1)p}{p}+\ip{b(x,1)}{p}, \dots, \ip{a(x,J)p}{p}+\ip{b(x,J)}{p}\right] + R_x,
		\]
		where $x,p \in \mathbb{R}^d$ and $R_x$ is the matrix with~$(R_x)_{ii} = -\sum_{j \neq i} r_{ij}(x)$ on the diagonal and $(R_x)_{ij} = r_{ij}(x)$ for $i \neq j$.
		\smallskip
		
	The next Hamiltonian arises from a diffusion process coupled to a diffusion.
	\begin{proposition}[Diffusion coupled to diffusion]\label{prop:diffusion-coupled-to-diffusion}
		Let $E=\mathbb{R}^d$ and $F$ be a smooth compact Riemannian manifold without boundary. Suppose the following.
		\begin{enumerate}[label=(\roman*)]
			\item The set of control variables $\Theta$ equals the space $\cP(F)$.
			\item The function $\Lambda$ is given by
			\begin{equation*} 
			\Lambda(x,p,\theta) := \int_F\left[\ip{a(x,z)p}{p}+\ip{b(x,z)}{p}\right]\,d\theta(z),
			\end{equation*}
			with $a:E\times F\to\mathbb{R}^{d\times d}$ and $b:E\times F\to\mathbb{R}^d$.
			\item The cost function $\mathcal{I}:E\times\Theta\to[0,\infty]$ is given by
			\begin{equation*}
			\mathcal{I}(x,\theta) := \sup_{\substack{u\in \mathcal{D}(L_x)\\u>0}}\left[ -\int_F \frac{L_xu}{u}\,d\theta\right],
			\end{equation*}
			where $L_x$ is a second-order elliptic operator locally of the form
			\begin{equation*}
			L_x = \frac{1}{2}\nabla\cdot\left(a_x \nabla\right) + b_x\cdot \nabla,
			\end{equation*}
			on the domain $\mathcal{D}(L_x):=C^2(F)$, with positive-definite matrix $a_x$ and co-vectors $b_x$.
		\end{enumerate}
		Suppose that the cost function~$\mathcal{I}$ satisfies the assumptions of Proposition~\ref{prop:verify:DV-functional-of-drift-diffusion} and the function~$\Lambda$ satisfies the assumptions of Proposition~\ref{prop:verify-ex:Lambda_quadratic}. Then Theorems~\ref{theorem:comparison_principle_variational} and~\ref{theorem:existence_of_viscosity_solution} apply to the Hamiltonian~\eqref{eq:results:H-example-section}.
	\end{proposition}
	In the context of weakly interacting jump processes on a collection of states $\{1,\dots,q\}$ the dynamics of the empirical measures takes place on $\cP(\{1,\dots,q\})$. Transitions occur over the bonds $(a,b) \in E^2$ with $a \neq b$. We denote the set of bonds with $\Gamma$. 
	\begin{definition}[Proper kernel] \label{definition:proper_kernel}
		Let $v : \Gamma \times \cP(\{1,\dots,q\}) \times \Theta \rightarrow \bR^+$. We say that $v$ is a \textit{proper kernel} if $v$ is continuous and if for each $(a,b) \in \Gamma$, the map $(\mu,\theta) \mapsto v(a,b,\mu,\theta)$ is either identically equal to zero or satisfies the following two properties:
		\begin{enumerate}[(a)]
			\item $v(a,b,\mu,\theta) = 0$ if $\mu(a) = 0$ and $v(a,b,\mu,\theta) > 0$ for all $\mu$ with $\mu(a) > 0$. 
			\item There exists a decomposition $v(a,b,\mu,\theta) = v_{\dagger}(a,b,\mu(a)) v_{\ddagger}(a,b,\mu,\theta)$ such that $v_{\dagger}$ is increasing in the third coordinate and such that $v_{\ddagger}(a,b,\cdot,\cdot)$ is continuous and satisfies $v_{\ddagger}(a,b,\mu,\theta) > 0$.
		\end{enumerate}
	\end{definition}
	A typical example of a proper kernel is given by 
	\begin{equation*}
	v(a,b,\mu,\theta) = \mu(a) r(a,b,\theta) e^{ \partial_a V(\mu) - \partial_b V(\mu)}, 
	\end{equation*}
	with $r > 0$ continuous and $V \in C^1_b(\cP(\{1,\dots,q\})$.
	\begin{proposition}[Mean-field coupled to diffusion]\label{prop:mean-field-coupled-to-diffusion}
		Let the space $E$ be given by the embedding of $E:=\mathcal{P}(\{1,\dots,J\})\times[0,\infty)^\Gamma\subseteq \mathbb{R}^d$ and $F$ be a smooth compact Riemannian manifold without boundary. Suppose the following. 
		\begin{enumerate}[label=(\roman*)]
			\item The set of control variables $\Theta$ equals $\mathcal{P}(F)$.
			\item The function $\Lambda$ is given by
			\begin{equation*} 
			\Lambda((\mu,w),p,\theta) = \sum_{(a,b) \in \Gamma} v(a,b,\mu,\theta)\left[\exp\left\{p_b - p _a + p_{(a,b)} \right\} - 1 \right]
			\end{equation*}
			with a proper kernel $v$ in the sense of Definition~\ref{definition:proper_kernel}.
			\item The cost function $\mathcal{I}:E\times\Theta\to[0,\infty]$ is given by
			\begin{equation*}
			\mathcal{I}(x,\theta) := \sup_{\substack{u\in \mathcal{D}(L_x)\\u>0}}\left[ -\int_F \frac{L_xu}{u}\,d\theta\right],
			\end{equation*}
			where $L_x$ is a second-order elliptic operator locally of the form
			\begin{equation*}
			L_x = \frac{1}{2}\nabla\cdot\left(a_x \nabla\right) + b_x\cdot \nabla,
			\end{equation*}
			on the domain $\mathcal{D}(L_x):=C^2(F)$, with positive-definite matrix $a_x$ and co-vectors $b_x$.
		\end{enumerate}
		Suppose that the cost function~$\mathcal{I}$ satisfies the assumptions of Proposition~\ref{prop:verify:DV-functional-of-drift-diffusion} and the function~$\Lambda$ satisfies the assumptions of Proposition~\ref{prop:verify-ex:Lambda_exponential}. Then Theorems~\ref{theorem:comparison_principle_variational} and~\ref{theorem:existence_of_viscosity_solution} apply to the Hamiltonian~\eqref{eq:results:H-example-section}.
	\end{proposition}
	An analogous proposition can be formulated for mean-field particles coupled to jumps as in Proposition~\ref{prop:diffusion-coupled-to-jumps}.
	\section{Strategy of the proofs} \label{section:strategy}
	We comment on the overall strategy of proofs. In Section~\ref{section:bootstrap-argument-in-nutshell}, we explain informally without the details how the bootstrap argument works in a simple setting in which $E$ is taken to be compact.  This allows us to focus on the bootstrapping argument without having to bother with the reduction to compact sets. We proceed with a discussion on the existence of a viscosity solution in Section~\ref{section:strategy_existence_viscosity_solution}.
	\subsection{The bootstrap argument in a nutshell}
	\label{section:bootstrap-argument-in-nutshell}
	In this section, we explain informally the main bootstrapping idea  behind proving the comparison principle with Hamiltonians of the type~\eqref{eq:results:variational_hamiltonian} for
	\begin{equation*}
	u(x)-\mathcal{H}(x,\nabla u(x)) = 0,
	\end{equation*}
	assuming compactness of $E$ and $\Psi(x,y) = \frac{1}{2}|x-y|^2$. In what follows, $u_1$ is a subsolution and $u_2$ is a supersolution. Recall that for smooth functions $f$, if $(u_1-f)$ is maximal at a point $x$, then
	\begin{equation*}
	u_1(x) -\mathcal{H}(x,\nabla f(x)) \leq 0.
	\end{equation*}
	Similarly for the supersolution $u_2$: If $(f-u_2)$ is maximal at a point $y$, then
	\begin{equation*}
	u_2(y)-\mathcal{H}(y,\nabla f(y)) \geq 0.
	\end{equation*}
	We sketch how to prove $u_1\leq u_2$ in several steps. 
	\begin{enumerate}[label=(\roman*)]
		\item By the classical doubling of variables procedure (e.g. \cite{CIL92}), choosing for each $\alpha > 0$ points $x_{\alpha},y_\alpha$ such that
		\begin{equation*}
		u_1(x_\alpha) - u_2(y_\alpha) - \alpha \Psi(x_\alpha,y_\alpha) = \sup_{x,y \in E} u_1(x) - u_2(y) - \alpha \Psi(x,y),
		\end{equation*}
		then by the properties of $\Psi$, we have 
		\begin{equation} \label{eqn:informal_control_Psi}
		\alpha \Psi(x_\alpha,y_\alpha) \rightarrow 0
		\end{equation}
		and the difference $\sup_x u_1(x) - u_2(x)$ can be approximated as
		\begin{equation*}
		\sup(u_1-u_2) \leq \liminf_{\alpha \rightarrow \infty} u_1(x_\alpha) - u_2(y_\alpha).
		\end{equation*}
		Set $p_\alpha := \alpha(x_\alpha - y_\alpha)$. Using the subsolution inequality $u_1(x_\alpha)\leq \mathcal{H}(x_\alpha,p_\alpha)$ and the supersolution inequality $u_2(y_\alpha) \geq \mathcal{H}(y_\alpha,p_\alpha)$, one arrives at the estimate
		\begin{equation*}
		\sup(u_1-u_2) \leq \liminf_{\alpha \rightarrow \infty} \mathcal{H}(x_\alpha,p_\alpha)-\mathcal{H}(y_\alpha,p_\alpha).
		\end{equation*}
		\item Recall that the Hamiltonian is given by
		\begin{equation*}
		\mathcal{H}(x,p) = \sup_{\theta\in\Theta}\left[\Lambda(x,p,\theta) - \mathcal{I}(x,\theta)\right].
		\end{equation*}
		Taking the optimizer $\theta_\alpha$ for $\mathcal{H}(x_\alpha,p_\alpha)$ and estimating the Hamiltonian at $y_\alpha$ with this optimizer, we obtain
		\begin{multline*}
		\sup(u_1-u_2) \leq \\
		\liminf_{\alpha \rightarrow \infty} \left[\Lambda(x_\alpha,p_\alpha,\theta_\alpha) - \Lambda(y_\alpha,p_\alpha,\theta_\alpha)\right] + \left[  \mathcal{I}(y_\alpha,\theta_\alpha)- \mathcal{I}(x_\alpha,\theta_\alpha)\right].
		\end{multline*}
		\item \label{item:informal_continuity_estimate} We assume the \emph{continuity estimate} on $\Lambda$. That means that if we have
		\begin{align}
		& \alpha \Psi(x_\alpha,y_\alpha) \rightarrow 0 \label{eqn:informal_control_cont_estimate1} \\
		& \liminf_{\alpha\to\infty}\Lambda(x_\alpha,p_\alpha,\theta_\alpha) > -\infty, \label{eqn:informal_control_cont_estimate2} \\
		& \limsup_{\alpha\to\infty}\Lambda(y_\alpha,p_\alpha,\theta_\alpha) < \infty, \label{eqn:informal_control_cont_estimate3}
		\end{align}
		and that $\theta_{\alpha}$ are in a compact set, then the difference of $\Lambda$'s is controlled as
		\begin{equation*}
		\liminf_{\alpha\to\infty}\left[\Lambda(x_\alpha,p_\alpha,\theta_\alpha) - \Lambda(y_\alpha,p_\alpha,\theta_\alpha)\right] \leq 0.
		\end{equation*}
		We postpone the verification that $\theta_{\alpha}$ are in a compact set to the next step \ref{item:informal_boundedness_I} below. Part~\eqref{eqn:informal_control_cont_estimate1} is just what we already know~\eqref{eqn:informal_control_Psi}. We show how the other two bounds follow from the sub- and supersolution inequalities. By the subsolution inequality,
		\begin{equation*}
		u_1(x_\alpha) \leq \mathcal{H}(x_\alpha,p_\alpha) = \Lambda(x_\alpha,p_\alpha,\theta_\alpha) - \mathcal{I}(x_\alpha,\theta_\alpha) \leq \Lambda(x_\alpha,p_\alpha,\theta_\alpha),
		\end{equation*}
		and~\eqref{eqn:informal_control_cont_estimate2} follows since $u_1$ is bounded. Letting $\theta_\alpha^0$ be the control variable such that $\mathcal{I}(y_\alpha,\theta_\alpha^0) = 0$, we obtain from the supersolution inequality that
		\begin{equation} \label{eqn:informal_supersolution_bound}
		u_2(y_\alpha) \geq \mathcal{H}(y_\alpha,p_\alpha) \geq \Lambda(y_\alpha,p_\alpha,\theta_\alpha^0),
		\end{equation}
		and therefore $\Lambda(y_\alpha,p_\alpha,\theta_\alpha^0)$ is bounded above. Assuming that
		\begin{equation*}
		\Lambda(y_\alpha,p_\alpha,\theta_\alpha) \leq C_1 \Lambda(y_\alpha,p_\alpha,\theta_\alpha^0) + C_2,
		\end{equation*}
		the bound~\eqref{eqn:informal_control_cont_estimate3} follows. In summary, if indeed $\theta_{\alpha}$ are in a compact set, taking the $\liminf_{\alpha\to\infty}$ in the last estimate on $(u_1-u_2)$, we obtain
		\begin{equation*}
		\sup(u_1-u_2) \leq 0 + \liminf_{\alpha\to\infty}\left[\mathcal{I}(y_\alpha,\theta_\alpha)- \mathcal{I}(x_\alpha,\theta_\alpha)\right].
		\end{equation*}
		\item \label{item:informal_boundedness_I} We assume that if the cost functions are uniformly bounded,
		\begin{equation} \label{eqn:informal_bounds_on_I}
		\mathcal{I}(x_\alpha,\theta_\alpha) \leq M \quad \text{and} \quad \mathcal{I}(y_\alpha,\theta_\alpha) \leq M,
		\end{equation}
		then (1) the control variables $\theta_\alpha$ are in a compact set, implying that we can carry out the argument of step~\ref{item:informal_continuity_estimate} above, and (2) the cost functions are continuous as a function of the internal variables $x$, giving
		\begin{equation*}
		\limsup_{\alpha\to\infty}\left[\mathcal{I}(y_\alpha,\theta_\alpha)- \mathcal{I}(x_\alpha,\theta_\alpha)\right] = 0.
		\end{equation*}
		The required bounds on $\mathcal{I}$ in~\eqref{eqn:informal_bounds_on_I} follow as well from the sub- and supersolution inequalities. From the subsolution inequality, we have
		\begin{equation*}
		u_1(x_\alpha) \leq \mathcal{H}(x_\alpha,p_\alpha) = \Lambda(x_\alpha,p_\alpha,\theta_\alpha) - \mathcal{I}(x_\alpha,\theta_\alpha).
		\end{equation*}
		Thus the bound on $\mathcal{I}(x_\alpha,\theta_\alpha)$ follows if we establish an upper bound on $\Lambda(x_\alpha,p_\alpha,\theta_\alpha)$.	Note that
		\begin{equation*}
		\Lambda(x_\alpha,p_\alpha,\theta_\alpha) \leq C_1 \Lambda(x_\alpha,p_\alpha,\theta_\alpha^0) + C_2
		\end{equation*}
		and
		\begin{equation*}
		\Lambda(x_\alpha,p_\alpha,\theta_\alpha^0) = \Lambda(y_\alpha,p_\alpha,\theta_\alpha^0) + \left[\Lambda(x_\alpha,p_\alpha,\theta_\alpha^0)-\Lambda(y_\alpha,p_\alpha,\theta_\alpha^0)\right].
		\end{equation*}
		We have an upper bound for the first term on the right-hand side by \eqref{eqn:informal_supersolution_bound}. The second term is bounded above by the continuity estimate, which can be carried out as we know that the $\theta_{\alpha}^0$ are in a compact set because they satisfy $\cI(y_\alpha,\theta_{\alpha}^0) = 0$. Since $y_\alpha$ is close to $x_\alpha$ and $\mathcal{I}$ is continuous as a function of $x$ when bounded, the bound on $\mathcal{I}(x_\alpha,\theta_\alpha)$ carries over to $\mathcal{I}(y_\alpha,\theta_\alpha)$.
	\end{enumerate}
	In summary, by using the information contained in the sub- and supersolution inequalities, the continuity estimate of the functions $\Lambda$ bootstraps to a continuity estimate of $\mathcal{H}$, giving the comparison principle.
	\subsection{Proof of the existence of a viscosity solution}\label{section:strategy_existence_viscosity_solution}
	For the existence of a viscosity solution to $f - \lambda \bfH f = h$, we will use the results of Chapter 8 of \cite{FengKurtz2006}. We will briefly discuss the method to obtain this result.
	
	To establish that $R(\lambda)h$ given by~\eqref{HJB:eq:resolvent}
	yields a viscosity solution to $f - \lambda \bfH f = h$, we follow a general strategy, first used in \cite{FengKurtz2006} and summarized in a more general context in~\cite[Proposition~3.4]{Kraaij2019GenConv}. For this strategy, we need to check three properties (see also~Section~\ref{BG:sec:semigroup-flow-HJ-eq} of Chapter~2):
	\begin{enumerate}[(a)]
		\item For all $(f,g) \in \bfH$, we have $f = R(\lambda)(f - \lambda g)$ ;
		\item The operator $R(\lambda)$ is a pseudo-resolvent: for all $h \in C_b(E)$ and $0 < \alpha < \beta$ we have 
		\begin{equation*}
		R(\beta)h = R(\alpha) \left(R(\beta)h - \alpha \frac{R(\beta)h - h}{\beta} \right).
		\end{equation*}
		\item The operator $R(\lambda)$ is contractive. 
	\end{enumerate}
	In other words: if $R(\lambda)$ serves as a classical left-inverse to $\bONE - \lambda \bfH$ and is also a pseudo-resolvent, then it is a viscosity right-inverse of $(\bONE- \lambda \bfH)$.
	Establishing (c) is a matter of writing out the definition. The proof of (a) and (b) stems from two main properties of exponential random variable. Let~$\tau_\lambda$ be the measure on~$\bR^+$ corresponding to the exponential random variable with mean~$\lambda^{-1}$.
	\begin{itemize}
		\item (a) is related to integration by parts: for bounded measurable functions $z$ on $\bR^+$, we have
		\begin{equation*}
		\lambda \int_0^\infty  z(t) \, \tau_\lambda( \dd t) = \int_0^\infty \int_0^t z(s) \, \dd s \, \tau_\lambda(\dd t).
		\end{equation*}
		\item (b) is related to a more involved integral property of exponential random variables. For $0 < \alpha < \beta$, we have
		\begin{multline*}
		\int_0^\infty z(s) \tau_\beta(\dd s) \\
		= \frac{\alpha}{\beta} \int_0^\infty z(s) \tau_\alpha(\dd s) + \left(1 - \frac{\alpha}{\beta}\right) \int_0^\infty \int_0^\infty z(s+u) \, \tau_\beta(\dd u) \, \tau_\alpha(\dd s).
		\end{multline*}
	\end{itemize}
	Establishing (a) and (b) can then be reduced by a careful analysis of optimizers in the definition of $R(\lambda)$, and concatenation or splittings thereof. This was carried out in Chapter 8 of \cite{FengKurtz2006} on the basis of three assumptions, namely \cite[Assumptions 8.9, 8.10 and 8.11]{FengKurtz2006}. We verify these in Section~\ref{section:construction-of-viscosity-solutions}.

	\section{Regularity of the Hamiltonian}
	\label{section:regularity-of-H-and-L}
	
	In this section, we establish continuity, convexity and the existence of a containment function for the Hamiltonian $\cH$ of \ref{eq:results:variational_hamiltonian}. We repeat its definition for convenience:
	\begin{equation} \label{eqn:regularity_section_H_variational_rep}
	\mathcal{H}(x,p) = \sup_{\theta \in \Theta}\left[\Lambda(x,p,\theta) - \mathcal{I}(x,\theta)\right].
	\end{equation}
	\begin{proposition}[Regularity of the Hamiltonian]\label{prop:reg-of-H-and-L:reg-H}
		Let $\mathcal{H} : E \times \mathbb{R}^d\to \mathbb{R}$ be the Hamiltonian as in \eqref{eqn:regularity_section_H_variational_rep}, and suppose that Assumptions~\ref{assumption:results:regularity_of_V} and~\ref{assumption:results:regularity_I} are satisfied. Then:
		\begin{enumerate}[label=(\roman*)]
			\item \label{item:prop:reg-of-H-and-L:convex} For any $x \in E$, the map $p \mapsto \mathcal{H}(x,p)$ is convex and $\mathcal{H}(x,0) = 0$.
			\item \label{item:prop:reg-of-H-and-L:compact-contain} With the containment function $\Upsilon : E \to \mathbb{R}$ of~\ref{item:assumption:slow_regularity:compact_containment}, we have
			\[
			\sup_{x \in E}\mathcal{H}(x,\nabla\Upsilon(x)) \leq C_\Upsilon < \infty.
			\]
		\end{enumerate}
	\end{proposition}
	\begin{proof}[Proof]
		The map $p \mapsto \mathcal{H}(x,p)$ is convex as it is the supremum over convex functions.
		
		\smallskip
		
		For proving $\mathcal{H}(x,0) = 0$, let $x \in E$. Then by~\ref{item:assumption:slow_regularity:convexity} of Assumption~\ref{assumption:results:regularity_of_V}, we have $\Lambda(x,0,\theta) = 0$, and therefore
		\[
		\mathcal{H}(x,0) = - \inf_{\theta\in\Theta} \mathcal{I}(x,\theta) = 0,
		\]
		since $\cI \geq 0$ and $\mathcal{I}(x,\theta_x)=0$ for some $\theta_x$ by~\ref{item:assumption:I:zero-measure} of Assumption~\ref{assumption:results:regularity_I}. Regarding~\ref{item:prop:reg-of-H-and-L:compact-contain}, we note that by~\ref{item:assumption:slow_regularity:compact_containment},
		\begin{align*}
		\mathcal{H}(x,\nabla \Upsilon(x)) \leq \sup_\theta \Lambda(x,\nabla \Upsilon(x),\theta) \leq \sup_{\theta\in\Theta}\sup_{x \in E} \Lambda(x,\nabla \Upsilon(x),\theta) \leq C_\Upsilon.
		\end{align*}
	\end{proof}
	To prove that $\cH$ is continuous, we use Assumption \ref{assumption:results:regularity_I}. What we truly need, however, is that $\cI$ Gamma converges as a function of $x$. We establish this result first.

	\begin{proposition}[Gamma convergence of the cost functions]\label{prop:Gamma-convergence-of-I}
		Let a cost function $\mathcal{I}:E\times\Theta\to[0,\infty]$ satisfy Assumption~\ref{assumption:results:regularity_I}. Then if $x_n\to x$ in $E$, the functionals $\mathcal{I}_n$ defined by
		\begin{equation*}
		\mathcal{I}_n(\theta) := \mathcal{I}(x_n,\theta)
		\end{equation*}
		converge in the $\Gamma$-sense to $\mathcal{I}_\infty(\theta) := \mathcal{I}(x,\theta)$. That is:
		\begin{enumerate}
			\item If $x_n \rightarrow x$ and $\theta_n \rightarrow \theta$, then $\liminf_{n\to\infty} \cI(x_n,\theta_n) \geq \cI(x,\theta)$,
			\item For $x_n \rightarrow x$ and all $\theta \in \Theta$ there are $\theta_n \in \Theta$ such that 
			\[
			\theta_n \rightarrow \theta\quad\text{and}\quad \limsup_{n\to\infty} \cI(x_n,\theta_n) \leq \cI(x,\theta).
			\]
		\end{enumerate}
	\end{proposition}
	\begin{proof}[Proof]
		Let $x_n\to x$. If $\theta_n\to \theta$, then by lower semicontinuity~\ref{item:assumption:I:lsc},
		\begin{equation*}
		\liminf_{n\to\infty}\mathcal{I}(x_n,\theta_n) \geq \mathcal{I}(x,\theta).
		\end{equation*}
		For the $\text{lim-sup}$ bound, let $\theta\in\Theta$. If $\mathcal{I}(x,\theta)=\infty$, there is nothing to prove. Thus suppose that $\mathcal{I}(x,\theta)$ is finite. Then by~\ref{item:assumption:I:finiteness}, there is a neighborhood $U_x$ of $x$ and a constant $M < \infty$ such that for any $y\in U_x$,
		\begin{equation*}
		\mathcal{I}(y,\theta) \leq M.
		\end{equation*}
		Since $x_n\to x$, the $x_n$ are eventually contained in $U_x$. Taking the constant sequence $\theta_n:=\theta$, we thus get that $\mathcal{I}(x_n,\theta_n) \leq M$ for all $n$ large enough. By~\ref{item:assumption:I:equi-cont}, 
		\begin{equation*}
		\lim_{n\to\infty}|\mathcal{I}(x_n,\theta_n)-\mathcal{I}(x,\theta)| \leq 0,
		\end{equation*}
		and the $\text{lim-sup}$ bound follows.
	\end{proof}
	\begin{proposition}[Continuity of the Hamiltonian]\label{prop:reg-of-H-and-L:continuity}
		Let $\mathcal{H} : E \times \mathbb{R}^d\to \mathbb{R}$ be the Hamiltonian defined in~\eqref{eq:results:variational_hamiltonian}, and suppose that Assumptions~\ref{assumption:results:regularity_of_V} and~\ref{assumption:results:regularity_I} are satisfied. Then the map $(x,p) \mapsto \cH(x,p)$ is continuous and the Lagrangian $(x,v) \mapsto \cL(x,v) := \sup_{p} \ip{p}{v} - \mathcal{H}(x,p)$ is lower semi-continuous.
	\end{proposition}
	Before we start with the proof, we give a remark on the generality of its statement and on the assumption that $\Theta$ is Polish.
	\begin{remark}
		The proof of upper semi-continuity of $\cH$ works in general, using continuity properties of $\Lambda$, lower semi-continuity of $(x,\theta) \mapsto I(x,\theta)$ and the compact sublevel sets of $\cI(x,\cdot)$. To establish lower semi-continuity,we need the that the functionals $\cI$ Gamma converge as a function of $x$. This was established in Proposition \ref{prop:Gamma-convergence-of-I}.
	\end{remark}
	\begin{remark}
		In the lemma we use a sequential characterization of upper hemi-continuity. This is inspired by the natural formulation of Gamma convergence in terms of sequences. An extension of our results to spaces $\Theta$ beyond the Polish context should take care of this issue. Without introducing the complicated matter, an extension is possible to Hausdorff $\Theta$ that are k-spaces in which all compact sets are metrizable.
	\end{remark}
	We will use the following technical result to establish upper semi-continuity of $\cH$.
	\begin{lemma}[Lemma 17.30 in \cite{AlBo06}] \label{lemma:upper_semi_continuity_abstract}
		Let $\cX$ and $\cY$ be two Polish spaces. Let $\phi : \cX \rightarrow \cK(\cY)$, where $\cK(\cY)$ is the space of non-empty compact subsets of $\cY$. Suppose that $\phi$ is upper hemi-continuous, that is if $x_n \rightarrow x$ and $y_n \rightarrow y$ and $y_n \in \phi(x_n)$, then $y \in \phi(x)$. 
		
		Let $f : \text{Graph} (\phi) \rightarrow \bR$ be upper semi-continuous. Then the map defined by~$m(x) := \sup_{y \in \phi(x)} f(x,y)$ is upper semi-continuous.
	\end{lemma} 
	\begin{proof}[Proof of Proposition~\ref{prop:reg-of-H-and-L:continuity}]
		We start by establishing upper semi-continuity of $\cH$. We argue on the basis of Lemma \ref{lemma:upper_semi_continuity_abstract}. Recall the representation of $\cH$ of \eqref{eqn:regularity_section_H_variational_rep}. Set $\mathcal{X} = E\times\mathbb{R}^d$ for the $(x,p)$ variables, $\mathcal{Y} = \Theta$, and $f(x,p,\theta) = \Lambda(x,p,\theta) - \cI(x,\theta)$ and note that this function is upper semi-continuous by Assumption~\ref{assumption:results:regularity_I}~\ref{item:assumption:I:lsc} and by Assumption \ref{assumption:results:regularity_of_V}~\ref{item:assumption:slow_regularity:continuity}.
		
		\smallskip
		
		By Assumption \ref{assumption:results:regularity_I}~\ref{item:assumption:I:zero-measure}, we have $\cH(x,p) \geq \Lambda(x,p,\theta_x)$.  Thus, it suffices to restrict the supremum over $\theta \in \Theta$ to $\theta \in \phi(x,p)$ where 
		\begin{equation*}
		\phi(x,p) := \left\{\theta \in \Theta \, \middle| \, \cI(x,\theta) \leq 2 \vn{\Lambda(x,p,\cdot)}_\Theta \right\},
		\end{equation*}
		in the sense that we have
		\begin{equation*}
		\mathcal{H}(x,p) = \sup_{\theta\in \phi(x,p)}\left[\Lambda(x,p,\theta)-\mathcal{I}(x,\theta)\right].
		\end{equation*}
		$\phi(x,p)$ is non-empty as $\theta_x \in \phi(x,p)$ and it is compact due to Assumption \ref{assumption:results:regularity_I}~\ref{item:assumption:I:compact-sublevelsets}. We are left to show that $\phi$ is upper hemi-continuous.
		
		\smallskip
		
		Thus, let $(x_n,p_n,\theta_n) \rightarrow (x,p,\theta)$ with $\theta_n \in \phi(x_n,p_n)$. We establish that $\theta \in \phi(x,p)$. By \ref{item:assumption:I:lsc} and the definition of $\phi$ we find
		\begin{equation*}
		\cI(x,\theta) \leq \liminf_n \cI(x_n,\theta_n) \leq \liminf_n 2\vn{\Lambda(x_n,p_n,\cdot}_\Theta = 2 \vn{\Lambda(x,p,\cdot)}_\Theta
		\end{equation*}
		which implies indeed that $\theta  \in \phi(x,p)$. Thus, upper semi-continuity follows by an application of Lemma \ref{lemma:upper_semi_continuity_abstract}.

		\smallskip
		
		We proceed with proving lower semi-continuity of $\cH$. Suppose that $(x_n,p_n) \rightarrow (x,p)$, we prove that $\liminf_n \cH(x_n,p_n) \geq \cH(x,p)$.

		Let $\theta$ be the measure such that $\cH(x,p) = \Lambda(x,p,\theta) - \cI(x,\theta)$. We have
		\begin{itemize}
			\item By Proposition \ref{prop:Gamma-convergence-of-I} there are $\theta_n$ such that $\theta_n \rightarrow \theta$ and $\limsup_n \cI(x_n,\theta_n) \leq \cI(x,\theta)$. 
			\item $\Lambda(x_n,p_n,\theta_n)$ converges to $\Lambda(x,p,\theta)$ by Assumption \ref{item:assumption:slow_regularity:continuity}.
		\end{itemize}
		Therefore,
		\begin{align*}
		\liminf_{n\to\infty}\mathcal{H}(x_n,p_n)&\geq \liminf_{n\to\infty} \left[\Lambda(x_n,p_n,\theta_n)-\mathcal{I}(x_n,\theta_n)\right]\\
		&\geq \liminf_{n\to\infty}\Lambda(x_n,p_n,\theta_n)-\limsup_{n\to\infty}\mathcal{I}(x_n,\theta_n)\\
		&\geq \Lambda(x,p,\theta)-\mathcal{I}(x,\theta) = \mathcal{H}(x,p),
		\end{align*}
		establishing that $\mathcal{H}$ is lower semi-continuous.
		
		\smallskip
		
		The Lagrangian $\cL$ is obtained as the supremum over continuous functions. This implies $\cL$ is lower semi-continuous.
	\end{proof}
	\section{The comparison principle} \label{section:comparison_principle}
	In this section, we establish the comparison principle for $f - \lambda \bfH f = h$ in the context of Theorem \ref{theorem:comparison_principle_variational}, using the general strategy of Section \ref{section:bootstrap-argument-in-nutshell}. Before being able to use this strategy, we need to restrict our analysis to compact sets in $E$. We will use a classical penalization technique that we will write down in operator form.
	
	\smallskip
	
	We thus introduce two new operators $H_\dagger$ and $H_\ddagger$, which are defined in terms of $\cH$ and the containment function $\Upsilon$ from Assumption \ref{assumption:results:regularity_of_V} \ref{item:assumption:slow_regularity:compact_containment}. We will then show that the comparison principle holds for a pair of Hamilton-Jacobi equations in terms of $H_\dagger$ and $H_\ddagger$. 
	This procedure allows us to clearly separate the reduction to compact sets on one hand, and the proof of the comparison principle on the basis of the bootstrap procedure on the other. Schematically, we will establish the following diagram:
	\begin{center}
		\begin{tikzpicture}
		\matrix (m) [matrix of math nodes,row sep=1em,column sep=4em,minimum width=2em]
		{
			{ } &[7mm] H_\dagger \\
			\bfH & { } \\
			{ }  & H_\ddagger \\};
		\path[-stealth]
		(m-2-1) edge node [above] {sub \qquad { }} (m-1-2)
		(m-2-1) edge node [below] {super \qquad { }} (m-3-2);
		
		\begin{pgfonlayer}{background}
		\node at (m-2-2) [rectangle,draw=blue!50,fill=blue!20,rounded corners, minimum width=1cm, minimum height=2.5cm]  {comparison};
		\end{pgfonlayer}
		\end{tikzpicture}
	\end{center}
	In this diagram, an arrow connecting an operator $A$ with operator $B$ with subscript 'sub' means that viscosity subsolutions of $f - \lambda A f = h$ are also viscosity subsolutions of $f - \lambda B f = h$. Similarly for arrows with a subscript 'super'.
	\smallskip
	
	We introduce the operators $H_\dagger$ and $H_\ddagger$ in Section~\ref{subsection:definition_of_Hamiltonians}. The arrows will be established in Section \ref{subsection:implications_from_compact_containment}. Finally, we will establish the comparison principle for $H_\dagger$ and $H_\ddagger$ in Section~\ref{subsection:proof_of_comparison_principle}, which by the arrows implies the comparison principle for $\bfH$.
	\begin{proof}[Proof of Theorem~\ref{theorem:comparison_principle_variational}]
		Fix $h_1,h_2 \in C_b(E)$ and $\lambda > 0$.
		\smallskip 
		
		Let $u_1,u_2$ be a viscosity sub- and supersolution to $f - \lambda \bfH f = h_1$ and  $f - \lambda \bfH f = h_2$ respectively. By Lemma \ref{lemma:viscosity_solutions_compactify2} proven in Section~\ref{subsection:implications_from_compact_containment}, $u_1$ and $u_2$ are a sub- and supersolution to $f - \lambda H_\dagger f = h_1$ and $f - \lambda H_\ddagger f = h_2$ respectively. Thus $\sup_E u_1 - u_2 \leq \sup_E h_1 - h_2$ by Proposition~\ref{prop:CP} of Section~\ref{subsection:proof_of_comparison_principle}. Specialising to $h_1=h_2$ gives Theorem~\ref{theorem:comparison_principle_variational}.
	\end{proof}
	\subsection{Definition of auxiliary operators} \label{subsection:definition_of_Hamiltonians}
	In this section, we repeat the definition of $\bfH$, and introduce the operators $H_\dagger$ and $H_\ddagger$.
	\begin{definition} \label{definition_effectiveH}
		The operator $\bfH \subseteq C_b^1(E) \times C_b(E)$ has domain $\cD(\bfH) = C_{cc}^\infty(E)$ and satisfies $\bfH f(x) = \cH(x, \dd f(x))$, where $\cH$ is the map
		\begin{equation*}
		\mathcal{H}(x,p) = \sup_{\theta \in \Theta}\left[\Lambda(x,p,\theta) - \mathcal{I}(x,\theta)\right].
		\end{equation*}
	\end{definition}
	We proceed by introducing $H_\dagger$ and $H_\ddagger$. These new Hamiltonians will serve as natural upper and lower bound for $\bfH$. They are defined in terms of the containment function $\Upsilon$, and essentially allow us to restrict our analysis to compact sets.
	\smallskip
	
	For the following definition, recall Assumption~\ref{item:assumption:slow_regularity:compact_containment} and the constant $C_\Upsilon := \sup_{\theta}\sup_x \Lambda(x,\nabla \Upsilon(x),\theta)$ therein. Denote by $C_\ell^\infty(E)$ the set of smooth functions on $E$ that have a lower bound and by $C_u^\infty(E)$ the set of smooth functions on $E$ that have an upper bound.
	\begin{definition}[The operators $H_\dagger$ and $H_\ddagger$]
		For $f \in C_\ell^\infty(E)$ and $\varepsilon \in (0,1)$  set 
		\begin{gather*}
		f^\varepsilon_\dagger := (1-\varepsilon) f + \varepsilon \Upsilon \\
		H_{\dagger,f}^\varepsilon(x) := (1-\varepsilon) \cH(x,\nabla f(x)) + \varepsilon C_\Upsilon.
		\end{gather*}
		and set
		\begin{equation*}
		H_\dagger := \left\{(f^\varepsilon_\dagger,H_{\dagger,f}^\varepsilon) \, \middle| \, f \in C_\ell^\infty(E), \varepsilon \in (0,1) \right\}.
		\end{equation*} 
		For $f \in C_u^\infty(E)$ and $\varepsilon \in (0,1)$  set 
		\begin{gather*}
		f^\varepsilon_\ddagger := (1+\varepsilon) f - \varepsilon \Upsilon \\
		H_{\ddagger,f}^\varepsilon(x) := (1+\varepsilon) \cH(x,\nabla f(x)) - \varepsilon C_\Upsilon.
		\end{gather*}
		and set
		\begin{equation*}
		H_\ddagger := \left\{(f^\varepsilon_\ddagger,H_{\ddagger,f}^\varepsilon) \, \middle| \, f \in C_u^\infty(E), \varepsilon \in (0,1) \right\}.
		\end{equation*} 
	\end{definition}
	\subsection{Implications based on compact containment} \label{subsection:implications_from_compact_containment}
	The operator $\bfH$ is related to $H_\dagger, H_\ddagger$ by the following Lemma.
	\begin{lemma}\label{lemma:viscosity_solutions_compactify2}
		Fix $\lambda > 0$ and $h \in C_b(E)$. 
		\begin{enumerate}[(a)]
			\item Every subsolution to $f - \lambda \bfH f = h$ is also a subsolution to $f - \lambda H_\dagger f = h$.
			\item Every supersolution to $f - \lambda \bfH f = h$ is also a supersolution to~$f-\lambda H_\ddagger f=~h$.
		\end{enumerate}
	\end{lemma}
	We only prove (a) of Lemma~\ref{lemma:viscosity_solutions_compactify2}, as (b) can be carried out analogously.

	\begin{proof}[Proof]
		Fix $\lambda > 0$ and $h \in C_b(E)$. Let $u$ be a subsolution to $f - \lambda \mathbf{H}f = h$. We prove it is also a subsolution to $f - \lambda H_\dagger f = h$.
		\smallskip
		
		Fix $\varepsilon > 0 $ and $f\in C_\ell^\infty(E)$ such that $(f^\varepsilon_\dagger,H^\varepsilon_{\dagger,f,\phi}) \in H_\dagger$. We will prove that there are $x_n\in E$ such that
		\begin{gather}
		\lim_{n\to\infty}\left(u-f_\dagger^\varepsilon\right)(x_n) = \sup_{x\in E}\left(u-f_\dagger^\varepsilon\right),\label{eqn:proof_lemma_conditions_for_subsolution_first}\\
		\limsup_{n\to\infty} \left[u(x_n)-\lambda H_{\dagger,f}^\varepsilon(x_n) - h(x_n)\right]\leq 0.\label{eqn:proof_lemma_conditions_for_subsolution_second}
		\end{gather}
		As the function $\left[u -(1-\varepsilon)f\right]$ is bounded from above and $\varepsilon \Upsilon$ has compact sublevel-sets, the sequence $x_n$ along which the first limit is attained can be assumed to lie in the compact set 
		\begin{equation*}
		K := \left\{x \, | \, \Upsilon(x) \leq \varepsilon^{-1} \sup_x \left(u(x) - (1-\varepsilon)f(x) \right)\right\}.
		\end{equation*}
		Set $M = \varepsilon^{-1} \sup_x \left(u(x) - (1-\varepsilon)f(x) \right)$. Let $\gamma : \bR \rightarrow \bR$ be a smooth increasing function such that
		\begin{equation*}
		\gamma(r) = \begin{cases}
		r & \text{if } r \leq M, \\
		M + 1 & \text{if } r \geq M+2.
		\end{cases}
		\end{equation*}
		Denote by $f_\varepsilon$ the function on $E$ defined by 
		\begin{equation*}
		f_\varepsilon(x) := \gamma\left((1-\varepsilon)f(x) + \varepsilon \Upsilon(x) \right).
		\end{equation*}
		By construction $f_\varepsilon$ is smooth and constant outside of a compact set and thus lies in $\cD(H) = C_{cc}^\infty(E)$. As $u$ is a viscosity subsolution for $f - \lambda Hf = h$ there exists a sequence $x_n \in K \subseteq E$ (by our choice of $K$) with
		\begin{gather}
		\lim_n \left(u-f_\varepsilon\right)(x_n) = \sup_x \left(u-f_\varepsilon\right)(x), \label{eqn:visc_subsol_sup} \\
		\limsup_n \left[u(x_n) - \lambda \mathbf{H} f_\varepsilon(x_n) - h(x_n)\right] \leq 0. \label{eqn:visc_subsol_upperbound}
		\end{gather}
		As $f_\varepsilon$ equals $f_\dagger^\varepsilon$ on $K$, we have from \eqref{eqn:visc_subsol_sup} that also
		\begin{equation*}
		\lim_n \left(u-f_\dagger^\varepsilon\right)(x_n) = \sup_{x\in E}\left(u-f_\dagger^\varepsilon\right),
		\end{equation*}
		establishing~\eqref{eqn:proof_lemma_conditions_for_subsolution_first}. Convexity of $p \mapsto \mathcal{H}(x,p)$ yields for arbitrary points $x\in K$ the estimate
		\begin{align*}
		\mathbf{H} f_\varepsilon(x) &= \mathcal{H}(x,\nabla f_\varepsilon(x)) \\
		& \leq (1-\varepsilon) \mathcal{H}(x,\nabla f(x)) + \varepsilon \mathcal{H}(x,\nabla \Upsilon(x)) \\
		&\leq (1-\varepsilon) \mathcal{H}(x,\nabla f(x)) + \varepsilon C_\Upsilon = H^\varepsilon_{\dagger,f}(x).
		\end{align*} 
		Combining this inequality with \eqref{eqn:visc_subsol_upperbound} yields
		\begin{multline*}
		\limsup_n \left[u(x_n) - \lambda H^\varepsilon_{\dagger,f}(x_n) - h(x_n)\right] \\
		\leq \limsup_n \left[u(x_n) - \lambda \mathbf{H} f_\varepsilon(x_n) - h(x_n)\right] \leq 0,
		\end{multline*}
		establishing \eqref{eqn:proof_lemma_conditions_for_subsolution_second}. This concludes the proof.
	\end{proof}
	\subsection{The comparison principle} \label{subsection:proof_of_comparison_principle}
	In this section, we prove the comparison principle for the operators $H_\dagger$ and $H_\ddagger$.
	\begin{proposition}\label{prop:CP} 
		Fix $\lambda > 0$ and $h_1,h_2 \in C_b(E)$. 	Let $u_1$ be a viscosity subsolution to $f - \lambda H_\dagger f = h_1$ and let $u_2$ be a viscosity supersolution to $f - \lambda H_\ddagger f = h_2$. Then we have $\sup_x u_1(x) - u_2(x) \leq \sup_x h_1(x) - h_2(x)$.
	\end{proposition}
	The proof uses an estimate that was proven in the proof of Proposition~A.11 of \cite{CoKr17} for one penalization function $\Psi$, or in the context of the more general continuity estimate of the Appendix of~\cite{KraaijSchlottke2019}, in the proof of Proposition~4.5 of \cite{Kr17} for two penalization functions $\{\Psi_1,\Psi_2\}$. In both contexts we use the containment function $\Upsilon$ of Assumption~\ref{assumption:results:regularity_of_V}, \ref{item:assumption:slow_regularity:compact_containment}. We start with a key result that allows us to find optimizing points that generalize the argument of Section \ref{section:bootstrap-argument-in-nutshell} to the non compact setting. 
	\smallskip
	
	The result is a copy of Lemma A.11 of \cite{CoKr17}, which is in turn a variant of Lemma 9.2 in \cite{FengKurtz2006} and Proposition 3.7 in \cite{CIL92}. We have included it for completeness.
	
	\begin{lemma}\label{lemma:doubling_lemma}
		Let $u$ be bounded and upper semi-continuous, let $v$ be bounded and lower semi-continuous, let $\Psi : E^2 \rightarrow \bR^+$ be penalization functions and let $\Upsilon$ be a containment function.
		\smallskip
		
		Fix $\varepsilon > 0$. For every $\alpha >0$ there exist $x_{\alpha,\varepsilon},y_{\alpha,\varepsilon} \in E$ such that
		\begin{multline} \label{eqn:existence_optimizers}
		\frac{u(x_{\alpha,\varepsilon})}{1-\varepsilon} - \frac{v(y_{\alpha,\varepsilon})}{1+\varepsilon} - \alpha \Psi(x_{\alpha,\varepsilon},y_{\alpha,\varepsilon}) - \frac{\varepsilon}{1-\varepsilon}\Upsilon(x_{\alpha,\varepsilon}) -\frac{\varepsilon}{1+\varepsilon}\Upsilon(y_{\alpha,\varepsilon}) \\
		= \sup_{x,y \in E} \left\{\frac{u(x)}{1-\varepsilon} - \frac{v(y)}{1+\varepsilon} -  \alpha \Psi(x,y)  - \frac{\varepsilon}{1-\varepsilon}\Upsilon(x) - \frac{\varepsilon}{1+\varepsilon}\Upsilon(y)\right\}.
		\end{multline}
		Additionally, for every $\varepsilon > 0$ we have that
		\begin{enumerate}[(a)]
			\item The set $\{x_{\alpha,\varepsilon}, y_{\alpha,\varepsilon} \, | \,  \alpha > 0\}$ is relatively compact in $E$.
			\item All limit points of $\{(x_{\alpha,\varepsilon},y_{\alpha,\varepsilon})\}_{\alpha > 0}$ as $\alpha \rightarrow \infty$ are of the form $(z,z)$ and for these limit points we have $u(z) - v(z) = \sup_{x \in E} \left\{u(x) - v(x) \right\}$.
			\item We have 
			\[
			\lim_{\alpha \rightarrow \infty}  \alpha \Psi(x_{\alpha,\varepsilon},y_{\alpha,\varepsilon}) = 0.
			\]
		\end{enumerate}
	\end{lemma}
	\begin{proof}[Proof of Proposition~\ref{prop:CP}]
		Fix $\lambda >0$ and $h_1,h_2 \in C_b(E)$. Let $u_1$ be a viscosity subsolution and $u_2$ be a viscosity supersolution of $f - \lambda H_\dagger f = h_1$ and  $f - \lambda H_\ddagger f = h_2$ respectively.  We prove Theorem~\ref{prop:CP} in two steps.
		\smallskip
		
		\underline{\emph{Step 1}}: We prove that for $\varepsilon > 0 $ and $\alpha > 0$, there exist points $x_{\varepsilon,\alpha},y_{\varepsilon,\alpha} \in E$ and momenta $p_{\varepsilon,\alpha}^1,p_{\varepsilon,\alpha}^2 \in \mathbb{R}^d$ such that
		\begin{multline} \label{eqn:estimate_step_1}
		\sup_E(u_1-u_2) \leq \lambda \liminf_{\varepsilon\to 0}\liminf_{\alpha \to \infty} \left[\mathcal{H}(x_{\varepsilon,\alpha},p^1_{\varepsilon,\alpha}) - \mathcal{H}(y_{\varepsilon,\alpha},p^2_{\varepsilon,\alpha})\right] \\ + \sup_{E}(h_1 - h_2).
		\end{multline}
		This step is solely based on the sub- and supersolution properties of $u_1,u_2$, the continuous differentiability of the penalization function $\Psi(x,y)$, the containment function $\Upsilon$, and convexity of $p \mapsto \mathcal{H}(x,p)$.
		\smallskip
		
		\underline{\emph{Step 2}}: Using Assumptions~\ref{assumption:results:regularity_of_V} and~\ref{assumption:results:regularity_I}, we prove that
		\begin{equation*}
		\liminf_{\varepsilon\to 0}\liminf_{\alpha \to \infty} \left[\mathcal{H}(x_{\varepsilon,\alpha},p^1_{\varepsilon,\alpha}) - \mathcal{H}(y_{\varepsilon,\alpha},p^2_{\varepsilon,\alpha})\right] \leq 0.
		\end{equation*}
		\smallskip
		
		\underline{\emph{Proof of Step 1}}: For any $\varepsilon > 0$ and any $\alpha > 0$, define the map $\Phi_{\varepsilon,\alpha}: E \times E \to \mathbb{R}$ by
		\begin{equation*}
		\Phi_{\varepsilon,\alpha}(x,y) := \frac{u_1(x)}{1-\varepsilon} - \frac{u_2(y)}{1+\varepsilon} - \alpha \Psi(x,y) - \frac{\varepsilon}{1-\varepsilon} \Upsilon(x) - \frac{\varepsilon}{1+\varepsilon}\Upsilon(y).
		\end{equation*}
		Let $\varepsilon > 0$. By Lemma \ref{lemma:doubling_lemma}, there is a compact set $K_\varepsilon \subseteq E$ and there exist points $x_{\varepsilon,\alpha},y_{\varepsilon,\alpha} \in K_\varepsilon$ such that
		\begin{equation} \label{eqn:comparison_optimizers}
		\Phi_{\varepsilon,\alpha}(x_{\varepsilon,\alpha},y_{\varepsilon,\alpha}) = \sup_{x,y \in E} \Phi_{\varepsilon,\alpha}(x,y),
		\end{equation}
		and 
		\begin{equation}\label{eq:proof-CP:Psi-xy-converge}
		\lim_{\alpha \to \infty} \alpha \Psi(x_{\varepsilon,\alpha},y_{\varepsilon,\alpha}) = 0.
		\end{equation}
		As in the proof of Proposition~A.11 of~\cite{Kr17}, it follows that
		\begin{equation}\label{eq:proof-CP:general-bound-u1u2}
		\sup_E (u_1 - u_2) \leq \liminf_{\varepsilon \to 0} \liminf_{\alpha \to \infty} \left[ \frac{u_1(x_{\varepsilon,\alpha})}{1-\varepsilon} - \frac{u_2(y_{\varepsilon,\alpha})}{1+\varepsilon}\right].
		\end{equation}
		At this point, we want to use the sub- and supersolution properties of $u_1$ and $u_2$. Define the test functions $\varphi^{\varepsilon,\alpha}_1 \in \cD(H_\dagger), \varphi^{\varepsilon,\alpha}_2 \in \cD(H_\ddagger)$ by
		\begin{multline*}
			\varphi^{\varepsilon,\alpha}_1(x):=(1-\varepsilon) \bigg[\frac{u_2(y_{\varepsilon,\alpha})}{1+\varepsilon} + \alpha \Psi(x,y_{\varepsilon,\alpha}) + \frac{\varepsilon}{1-\varepsilon}\Upsilon(x) + \frac{\varepsilon}{1+\varepsilon}\Upsilon(y_{\varepsilon,\alpha})\\ + (1-\varepsilon)(x-x_{\varepsilon,\alpha})^2\bigg]
		\end{multline*}
	and
	\begin{multline*}
	\varphi^{\varepsilon,\alpha}_2(y):=	(1+\varepsilon)\bigg[\frac{u_1(x_{\varepsilon,\alpha})}{1-\varepsilon} - \alpha \Psi(x_{\varepsilon,\alpha},y) - \frac{\varepsilon}{1-\varepsilon}\Upsilon(x_{\varepsilon,\alpha}) - \frac{\varepsilon}{1+\varepsilon}\Upsilon(y)\\-(1+\varepsilon) (y-y_{\varepsilon,\alpha})^2\bigg].
	\end{multline*}
		Using \eqref{eqn:comparison_optimizers}, we find that $u_1 - \varphi^{\varepsilon,\alpha}_1$ attains its supremum at $x = x_{\varepsilon,\alpha}$, and thus
		\begin{equation*}
		\sup_E (u_1-\varphi^{\varepsilon,\alpha}_1) = (u_1-\varphi^{\varepsilon,\alpha}_1)(x_{\varepsilon,\alpha}).
		\end{equation*}
		Denote $p_{\varepsilon,\alpha}^1 := \alpha \nabla_x \Psi(x_{\varepsilon,\alpha},y_{\varepsilon,\alpha})$. By our addition of the penalization $(x-x_{\varepsilon,\alpha})^2$ to the test function, the point $x_{\varepsilon,\alpha}$ is in fact the unique optimizer, and we obtain from the subsolution inequality that
		\begin{equation}\label{eq:proof-CP:subsol-ineq}
		u_1(x_{\varepsilon,\alpha}) - \lambda \left[ (1-\varepsilon) \mathcal{H}\left(x_{\varepsilon,\alpha}, p_{\varepsilon,\alpha}^1 \right) + \varepsilon C_\Upsilon\right] \leq h_1(x_{\varepsilon,\alpha}).
		\end{equation}	
		With a similar argument for $u_2$ and $\varphi^{\varepsilon,\alpha}_2$, we obtain by the supersolution inequality that
		\begin{equation}\label{eq:proof-CP:supersol-ineq}
		u_2(y_{\varepsilon,\alpha}) - \lambda \left[(1+\varepsilon)\mathcal{H}\left(y_{\varepsilon,\alpha}, p_{\varepsilon,\alpha}^2 \right) - \varepsilon C_\Upsilon\right] \geq h_2(y_{\varepsilon,\alpha}),
		\end{equation}
		where $p_{\varepsilon,\alpha}^2 := -\alpha \nabla_y \Psi(x_{\varepsilon,\alpha},y_{\varepsilon,\alpha})$. With that, estimating further in~\eqref{eq:proof-CP:general-bound-u1u2} leads to
		\begin{multline*}
		\sup_E(u_1-u_2) \leq \liminf_{\varepsilon\to 0}\liminf_{\alpha \to \infty} \bigg[\frac{h_1(x_{\varepsilon,\alpha})}{1-\varepsilon} - \frac{h_2(y_{\varepsilon,\alpha})}{1+\varepsilon} + \frac{\varepsilon}{1-\varepsilon} C_\Upsilon \\ + \frac{\varepsilon}{1+\varepsilon} C_\Upsilon  + \lambda \left[\mathcal{H}(x_{\varepsilon,\alpha},p^1_{\varepsilon,\alpha}) - \mathcal{H}(y_{\varepsilon,\alpha},p^2_{\varepsilon,\alpha})\right]\bigg].
		\end{multline*}
		Thus, \eqref{eqn:estimate_step_1} in Step 1 follows.
		
		\smallskip

		\underline{\emph{Proof of Step 2}}: Recall that $\mathcal{H}(x,p)$ is given by
		\begin{equation*}
		\mathcal{H}(x,p) = \sup_{\theta \in \Theta}\left[\Lambda(x,p,\theta) - \mathcal{I}(x,\theta)\right].
		\end{equation*}
		Since $\Lambda(x_{\varepsilon,\alpha},p^1_{\varepsilon,\alpha},\cdot) : \Theta \to \mathbb{R}$ is bounded and continuous by \ref{item:assumption:slow_regularity:continuity} and the map $\mathcal{I}(x_{\varepsilon,\alpha},\cdot) : \Theta \to [0,\infty]$ has compact sub-level sets in $\Theta$ by~\ref{item:assumption:I:compact-sublevelsets}, there exists an optimizer $\theta_{\varepsilon,\alpha} \in\Theta$ such that
		\begin{equation} \label{eqn:choice_of_optimal_measure}
		\mathcal{H}(x_{\varepsilon,\alpha},p^1_{\varepsilon,\alpha}) = \Lambda(x_{\varepsilon,\alpha},p^1_{\varepsilon,\alpha},\theta_{\varepsilon,\alpha}) - \mathcal{I}(x_{\varepsilon,\alpha},\theta_{\varepsilon,\alpha}).
		\end{equation}
		Choosing the same point in the supremum of the second term $\mathcal{H}(y_{\varepsilon,\alpha},p^2_{\varepsilon,\alpha})$, we obtain for all $\varepsilon > 0$ and $\alpha > 0$ the estimate
		\begin{multline} \label{eqn:basic_decomposition_Hamitlonian_difference}
		\mathcal{H}(x_{\varepsilon,\alpha},p^1_{\varepsilon,\alpha})-
		\mathcal{H}(y_{\varepsilon,\alpha},p^2_{\varepsilon,\alpha})
		\leq 
		\Lambda(x_{\varepsilon,\alpha},p^1_{\varepsilon,\alpha},\theta_{\varepsilon,\alpha})-
		\Lambda(y_{\varepsilon,\alpha},p^2_{\varepsilon,\alpha},\theta_{\varepsilon,\alpha})\\
		+\mathcal{I}(y_{\varepsilon,\alpha},\theta_{\varepsilon,\alpha})-
		\mathcal{I}(x_{\varepsilon,\alpha},\theta_{\varepsilon,\alpha}).
		\end{multline}
		We will establish an upper bound for this difference using the continuity estimate  \ref{item:assumption:slow_regularity:continuity_estimate} and equi-continuity \ref{item:assumption:I:equi-cont}.
		\smallskip
		
		To apply the continuity estimate \ref{item:assumption:slow_regularity:continuity_estimate}, we need to verify \eqref{eqn:control_on_Gbasic_sup} and \eqref{eqn:control_on_Gbasic_inf} (see \eqref{eq:proof-CP:Vx-unif-bound-below} and \eqref{eq:proof-CP:Vy-unif-bound-above} below) for the variables $\theta_{\varepsilon,\alpha}$. In addition, we need to establish that $\theta_{\varepsilon,\theta}$ are contained in a compact set. 
		\smallskip
		
		To apply \ref{item:assumption:I:equi-cont}, we need to control the size of $\cI(x_{\varepsilon,\alpha},\theta_{\varepsilon,\alpha})$ and $\cI(y_{\varepsilon,\alpha},\theta_{\varepsilon,\alpha})$ along subsequences, which by Assumption \ref{item:assumption:I:compact-sublevelsets} implies the above requirement that along these subsequences $\theta_{\varepsilon,\theta}$ are contained in a compact set. To obtain control on the size of $\cI$, we employ an auxiliary argument based on the continuity estimate for the measures $\theta_{\varepsilon,\alpha}^0$, obtained by \ref{item:assumption:I:zero-measure}, satisfying
		\begin{equation}\label{eqn:choice_of_zero_measure}
		\mathcal{I}(y_{\alpha,\varepsilon},\theta_{\varepsilon,\alpha}^0) = 0.
		\end{equation}
		The application of the continuity estimate for $\theta_{\varepsilon,\alpha}^0$ only requires to check \eqref{eqn:control_on_Gbasic_sup} and \eqref{eqn:control_on_Gbasic_inf} as the measures $\theta_{\varepsilon,\alpha}^0$ are contained in a compact set by \eqref{eqn:choice_of_zero_measure} and~\ref{item:assumption:I:compact-sublevelsets}.
		Thus, we will first establish
		\begin{align}
		& \liminf_{\varepsilon \to 0}\liminf_{\alpha \to \infty} \Lambda(x_{\varepsilon,\alpha},p^1_{\varepsilon,\alpha},\theta_{\varepsilon,\alpha}) > - \infty,\label{eq:proof-CP:Vx-unif-bound-below} \\
		& \liminf_{\varepsilon \to 0}\liminf_{\alpha \to \infty} \Lambda(x_{\varepsilon,\alpha},p^1_{\varepsilon,\alpha},\theta_{\varepsilon,\alpha}^0) > - \infty, \label{eq:proof-CP:Vx-unif-bound-below_zero_measure} \\
		& \limsup_{\varepsilon\to 0}\limsup_{\alpha \to \infty}\Lambda(y_{\varepsilon,\alpha},p^2_{\varepsilon,\alpha},\theta_{\varepsilon,\alpha}) < \infty, \label{eq:proof-CP:Vy-unif-bound-above} \\
		& \limsup_{\varepsilon\to 0}\limsup_{\alpha \to \infty}\Lambda(y_{\varepsilon,\alpha},p^2_{\varepsilon,\alpha},\theta_{\varepsilon,\alpha}^0) < \infty. \label{eq:proof-CP:Vy-unif-bound-above_zero_measure}
		\end{align}
		Note that by \ref{item:assumption:slow_regularity:controlled_growth} the bounds in \eqref{eq:proof-CP:Vx-unif-bound-below} and \eqref{eq:proof-CP:Vx-unif-bound-below_zero_measure} are equivalent. Similarly \eqref{eq:proof-CP:Vy-unif-bound-above} and \eqref{eq:proof-CP:Vy-unif-bound-above_zero_measure} are equivalent.
		
		\smallskip
		
		By the subsolution inequality~\eqref{eq:proof-CP:subsol-ineq},
		\begin{align}\label{eq:proof-CP:estimate-via-subsol}
		\frac{1}{\lambda} \inf_E\left(u_1 - h\right) & \leq (1-\varepsilon) \mathcal{H}(x_{\varepsilon,\alpha},p^1_{\varepsilon,\alpha}) + \varepsilon C_{\Upsilon}	\\
		&\leq (1-\varepsilon) \Lambda(x_{\varepsilon,\alpha},p^1_{\varepsilon,\alpha},\theta_{\varepsilon,\alpha}) + \varepsilon C_\Upsilon,\notag
		\end{align}
		and the lower bounds~\eqref{eq:proof-CP:Vx-unif-bound-below} and \eqref{eq:proof-CP:Vx-unif-bound-below_zero_measure} follow.
		
		\smallskip
		
		By the supersolution inequality~\eqref{eq:proof-CP:supersol-ineq}, we can estimate 
		\begin{align*}
		(1+\varepsilon) \Lambda(y_{\varepsilon,\alpha},p^2_{\varepsilon,\alpha},\theta_{\varepsilon,\alpha}^0) &= (1+\varepsilon) \left[\Lambda(y_{\varepsilon,\alpha},p^2_{\varepsilon,\alpha},\theta_{\varepsilon,\alpha}^0) - \mathcal{I}(y_{\varepsilon,\alpha},\theta_{\varepsilon,\alpha}^0)\right]	\\
		&\leq \left((1+\varepsilon) \mathcal{H}\left(y_{\varepsilon,\alpha},p^2_{\varepsilon,\alpha}\right) - \varepsilon C_\Upsilon\right) + \varepsilon C_\Upsilon	\\
		&\leq \frac{1}{\lambda} \sup_E (u_2-h) + \varepsilon C_{\Upsilon} < \infty,
		\end{align*}
		and the upper bounds~\eqref{eq:proof-CP:Vy-unif-bound-above} and \eqref{eq:proof-CP:Vy-unif-bound-above_zero_measure} follow. 
		\smallskip
		
		Since the $\theta_{\varepsilon,\alpha}^0$ are contained in a compact set by~\ref{item:assumption:I:compact-sublevelsets}, we conclude by the continuity estimate~\ref{item:assumption:slow_regularity:continuity_estimate} that
		\begin{equation*}
		\liminf_{\varepsilon\to 0}\liminf_{\alpha\to\infty}\left[\Lambda\left(x_{\varepsilon,\alpha},p^1_{\varepsilon,\alpha},\theta_{\varepsilon,\alpha}^0\right)-\Lambda\left(y_{\varepsilon,\alpha},p^2_{\varepsilon,\alpha},\theta_{\varepsilon,\alpha}^0\right)\right]\leq 0.
		\end{equation*}
		Without loss of generality, we can choose for all small $\varepsilon$ subsequences $(x_{\varepsilon,\alpha},y_{\varepsilon,\alpha})$ (denoted the same) such that also 
		\begin{equation}
		\liminf_{\varepsilon\to 0}\limsup_{\alpha \to \infty} \left[\Lambda\left(x_{\varepsilon,\alpha},p^1_{\varepsilon,\alpha},\theta_{\varepsilon,\alpha}^0\right)-\Lambda\left(y_{\varepsilon,\alpha},p^2_{\varepsilon,\alpha},\theta_{\varepsilon,\alpha}^0\right)\right] \leq 0.\label{eqn:proof_comparison_application_continuity_estimate_zero_measures_strenghtened}
		\end{equation}
		
		We proceed to establish that along this collection of subsequences we have $\limsup_{\alpha \rightarrow \infty} \mathcal{I}(x_{\varepsilon,\alpha},\theta_{\varepsilon,\alpha}) < \infty$. We return to the first inequality of \eqref{eq:proof-CP:estimate-via-subsol}, combined with \eqref{eqn:choice_of_optimal_measure}, to obtain
		\begin{align*}
		\frac{1}{\lambda} \inf_E\left(u_1 - h\right) & \leq (1-\varepsilon) \cH(x_{\varepsilon,\alpha},p^1_{\varepsilon,\alpha}) + \varepsilon C_{\Upsilon} \\
		& = (1-\varepsilon) \left[\Lambda(x_{\varepsilon,\alpha},p^1_{\varepsilon,\alpha},\theta_{\varepsilon,\alpha}) - \mathcal{I}(x_{\varepsilon,\alpha},\theta_{\varepsilon,\alpha}) \right] + \varepsilon C_{\Upsilon}.
		\end{align*}
		We conclude that $\limsup_{\alpha \rightarrow \infty} \mathcal{I}(x_{\varepsilon,\alpha},\theta_{\varepsilon,\alpha}) < \infty$ is implied by
		\begin{equation*}
		\limsup_{\alpha \rightarrow \infty} \Lambda(x_{\varepsilon,\alpha},p^1_{\varepsilon,\alpha},\theta_{\varepsilon,\alpha}) < \infty
		\end{equation*}
		which by \ref{item:assumption:slow_regularity:controlled_growth} is equivalent to
		\begin{equation*}
		\limsup_{\alpha \rightarrow \infty} \Lambda(x_{\varepsilon,\alpha},p^1_{\varepsilon,\alpha},\theta_{\varepsilon,\alpha}^0) < \infty.
		\end{equation*}
		This, however, yields what we want by \eqref{eq:proof-CP:Vy-unif-bound-above_zero_measure} and \eqref{eqn:proof_comparison_application_continuity_estimate_zero_measures_strenghtened}:
		\begin{multline*}
		\limsup_{\alpha \rightarrow \infty}\Lambda(x_{\varepsilon,\alpha},p^1_{\varepsilon,\alpha},\theta_{\varepsilon,\alpha}^0) \leq \limsup_{\alpha \rightarrow \infty} \Lambda(y_{\varepsilon,\alpha},p^2_{\varepsilon,\alpha},\theta_{\varepsilon,\alpha}^0) \\
		\qquad + \limsup_{\alpha \rightarrow \infty} \left[\Lambda\left(x_{\varepsilon,\alpha},p^1_{\varepsilon,\alpha},\theta_{\varepsilon,\alpha}^0\right)-\Lambda\left(y_{\varepsilon,\alpha},p^2_{\varepsilon,\alpha},\theta_{\varepsilon,\alpha}^0\right)\right] < \infty.
		\end{multline*}
		We thus obtain
		\begin{equation*}
		\limsup_{\alpha \rightarrow \infty} \mathcal{I}(x_{\varepsilon,\alpha},\theta_{\varepsilon,\alpha}) < \infty.
		\end{equation*}
		Therefore, by~\ref{item:assumption:I:compact-sublevelsets}, for each $\varepsilon > 0$ the $\theta_{\varepsilon,\alpha}$ are contained in a compact set. With the bounds~\eqref{eq:proof-CP:Vx-unif-bound-below} and~\eqref{eq:proof-CP:Vy-unif-bound-above}, we conclude by the continuity estimate~\ref{item:assumption:slow_regularity:continuity_estimate} that
		\begin{equation}\label{eq:proof-CP:final-bound-by-cont-estimate}
		\liminf_{\varepsilon\to 0}\liminf_{\alpha \to \infty} \left[\Lambda\left(x_{\varepsilon,\alpha},p^1_{\varepsilon,\alpha},\theta_{\varepsilon,\alpha}\right)-\Lambda\left(y_{\varepsilon,\alpha},p^2_{\varepsilon,\alpha},\theta_{\varepsilon,\alpha}\right)\right] \leq 0.
		\end{equation}
		By~\eqref{eq:proof-CP:Psi-xy-converge}, we have along a subsequence $(x_{\varepsilon,\alpha},y_{\varepsilon,\alpha}) \to (z_\varepsilon,z_\varepsilon) \in K_\varepsilon \times K_\varepsilon$ as $\alpha \to \infty$. Therefore by~\ref{item:assumption:I:finiteness} there exists a subsequence of $(x_{\varepsilon,\alpha},y_{\varepsilon,\alpha})$ (denoted the same) and a constant $M_\varepsilon' < \infty$ such that for all $\alpha > 0$ large enough,
		\begin{equation*}
		\mathcal{I}(x_{\varepsilon,\alpha},\theta_{\varepsilon,\alpha}) \leq M_\varepsilon' \quad \text{ and } \quad \mathcal{I}(y_{\varepsilon,\alpha},\theta_{\varepsilon,\alpha}) \leq M_\varepsilon'.
		\end{equation*}
		Hence by~\ref{item:assumption:I:equi-cont}, for any $\varepsilon > 0$,
		\begin{equation}\label{eq:proof-CP:Iy-Ix-leq-0}
		\limsup_{\alpha \to \infty}|\mathcal{I}(y_{\varepsilon,\alpha},\theta_{\varepsilon,\alpha}) - \mathcal{I}(x_{\varepsilon,\alpha},\theta_{\varepsilon,\alpha})| = 0.
		\end{equation}
		Then combining~\eqref{eq:proof-CP:final-bound-by-cont-estimate} with~\eqref{eq:proof-CP:Iy-Ix-leq-0} gives an estimate on \eqref{eqn:basic_decomposition_Hamitlonian_difference} which completes Step 2. 
	\end{proof}
	\section{Construction of viscosity solutions}
	\label{section:construction-of-viscosity-solutions}
	In this Section, we will show that $R(\lambda)h$, for $h \in C_b(E), \lambda > 0$ of Theorem \ref{theorem:existence_of_viscosity_solution} is indeed a viscosity solution to $f - \lambda \bfH f = h$. To do so, we will use the methods of Chapter 8 of \cite{FengKurtz2006} which are based on the strategy laid out in Section \ref{section:strategy_existence_viscosity_solution}.
	\smallskip
	
	In particular, we will verify \cite[Conditions 8.9, 8.10 and 8.11]{FengKurtz2006} which imply by~\cite[Theorem 8.27]{FengKurtz2006} and the comparison principle for $f - \lambda \bfH f = h$ that $R(\lambda)h$ is a viscosity solution to $f - \lambda \bfH f = h$.
	\begin{proof}[Verification of Conditions 8.9, 8.10 and 8.11]
		In the notation of \cite{FengKurtz2006}, we use $U = \bR^d$, $\Gamma = E \times U$, one operator $\bfH = \bfH_\dagger = \bfH_\ddagger$  and $Af(x,u) = \ip{\nabla f(x)}{u}$ for $f \in \mathcal{D}(\mathbf{H}) = C_{cc}^\infty(E)$.
		
		\smallskip
		
		Regarding Condition~8.9, by continuity and convexity of $\cH$ obtained in Propositions \ref{prop:reg-of-H-and-L:reg-H} and \ref{prop:reg-of-H-and-L:continuity}, parts 8.9.1, 8.9.2, 8.9.3 and 8.9.5 can be proven e.g. as in the proof of \cite[Lemma 10.21]{FengKurtz2006} for $\psi = 1$. Part 8.9.4 is a consequence of the existence of a containment function, and follows as shown in the proof of~\cite[Theorem~A.17]{CoKr17}. Since we use the argument further below, we briefly recall it here. We need to show that for any compact set $K \subseteq E$, any finite time $T > 0$ and finite bound $M \geq 0$, there exists a compact set $K' = K'(K,T,M) \subseteq E$ such that for any absolutely continuous path $\gamma :[0,T] \to E$ with $\gamma(0) \in K$, if
		\begin{equation} \label{eqn:control_on_L}
		\int_0^T \mathcal{L}(\gamma(t),\dot{\gamma}(t)) \, dt \leq M,
		\end{equation}
		then $\gamma(t) \in K'$ for any $0\leq t \leq T$.
		\smallskip
		
		For $K\subseteq E$, $T>0$, $M\geq 0$ and $\gamma$ as above, this follows by noting that
		\begin{align}
		\label{eq:action_integral_representation:Lyapunov_bound}
		\Upsilon(\gamma(\tau)) &=
		\Upsilon(\gamma(0)) + \int_0^\tau \nabla\Upsilon(\gamma(t)) \dot{\gamma}(t) \, dt \notag	\\
		&\leq \Upsilon(\gamma(0)) + \int_0^\tau \left[
		\mathcal{L}(\gamma(t),\dot{\gamma}(t))) + \mathcal{H}(x(t),\nabla \Upsilon(\gamma(t)))
		\right] \, dt \notag	\\
		&\leq \sup_K \Upsilon + M + T \sup_{x \in E} \mathcal{H}(x,\nabla \Upsilon(x)) =: C < \infty,
		\end{align}
		for any $0 \leq \tau \leq T$, so that the compact set $K' := \{z \in E \,:\, \Upsilon(z) \leq C\}$ satisfies the claim.
		
		\smallskip
		
		We proceed with the verification of Conditions~8.10 and 8.11 of~\cite{FengKurtz2006}. By Proposition~\ref{prop:reg-of-H-and-L:reg-H}, we have $\cH(x,0) = 0$ and hence $\bfH 1 = 0$. Thus, Condition 8.10 is implied by Condition 8.11 (see Remark 8.12 (e) in~\cite{FengKurtz2006}).
		\smallskip
		
		We establish that Condition 8.11 is satisfied: for any function $f\in \mathcal{D}(\bfH) = C_{cc}^\infty(E)$ and $x_0 \in E$, there exists an absolutely continuous path $x:[0,\infty) \to E$ such that $x(0) = x_0$ and for any $t \geq 0$,
		\begin{equation}
		\label{eq:action_integral_representation:solution_control_problem}
		\int_0^t \mathcal{H}(x(s),\nabla f(x(s)) \, ds =
		\int_0^t \left[
		\dot{x}(s) \cdot \nabla f(x(s)) - \mathcal{L}(x(s),\dot{x}(s))
		\right] \, ds.
		\end{equation}
		To do so, we solve the differential inclusion
		\begin{equation}
		\label{eq:action_integral_representation:subdifferential_eq}
		\dot{x}(t) \in \partial_p \mathcal{H}(x(t),\nabla f(x(t))), \qquad x(0) = x_0,
		\end{equation}
		where the subdifferential of $\cH$ was defined in \eqref{eqn:subdifferential_H} on page \pageref{eqn:subdifferential_H}. 
		\smallskip
		
		Since the addition of a constant to $f$ does not change the gradient, we may assume without loss of generality that $f$ has compact support. A general method to establish existence of differential inclusions~$\dot{x} \in F(x)$ is given by Lemma 5.1 of Deimling~\cite{De92}.
		We use this result for $F(x) := \partial_p \cH(x,\nabla f(x))$. To apply this lemma, we need to verify that:
		\begin{enumerate}[(F1)]
			\item $F$ is upper hemi-continuous and $F(x)$ is non-empty, closed, and convex for all $x \in E$.
			\item $\|F(x)\| \leq c(1 + |x|)$ on $E$, for some $c > 0$.
			\item $F(x) \cap T_E(x) \neq \emptyset$ for all $x \in E$. (The Definition~\ref{definition:tangent_cone} of~$T_E(x)$ is given on page~\pageref{eqn:subdifferential_H} of this thesis).
		\end{enumerate}
		While part (F1) follows from the properties of a subdifferential set and (F3) is a consequence of Assumption \ref{assumption:Hamiltonian_vector_field}, part (F2) is in general not satisfied. To circumvent this problem, we use properties of $\cH$ to establish a-priori bounds on the range of solutions.
		
		\smallskip
		
		\emph{Step 1:} Let $T > 0$, and assume that $x(t)$ solves \eqref{eq:action_integral_representation:subdifferential_eq}. We establish that there is some $M$ such that~\eqref{eqn:control_on_L} is satisfied. By~\eqref{eq:action_integral_representation:subdifferential_eq} we obtain for all $p \in \mathbb{R}^d$,
		\[
		\mathcal{H}(x(t),p) \geq \mathcal{H}(x(t),\nabla f(x(t))) + \dot{x}(t) \cdot (p - \nabla f(x(t))),
		\]
		and as a consequence
		\[
		\dot{x}(t) \nabla f(x(t)) - \mathcal{H}(x(t),\nabla f(x(t))) \geq
		\mathcal{L}(x(t),\dot{x}(t)).
		\]
		Since $f$ has compact support and $\mathcal{H}(y,0) = 0$ for any $y \in E$, we estimate 
		\begin{align*}
		\int_0^T \mathcal{L}(x(t),\dot{x}(t)) \, ds &\leq
		\int_0^T \dot{x}(t) \nabla f(x(t)) \, dt - T\inf_{y \in \mathrm{supp}(f)} \mathcal{H}(y,\nabla f(y)).
		\end{align*}
		By continuity of $\mathcal{H}$ the field $F$ is bounded on compact sets, so the first term can be bounded by
		\[
		\int_0^T \dot{x}(t) \nabla f(x(t)) \, dt \leq
		T \sup_{y \in \mathrm{supp}(f)}\|F(y)\| \sup_{z \in \mathrm{supp}(f)}|\nabla f(z)|.
		\]
		Therefore, for any $T>0$, we obtain that the integral over the Lagrangian is bounded from above by $M = M(T)$, with
		\[
		M := T \sup_{y \in \mathrm{supp}(f)}\|F(y)\| \sup_{z \in \mathrm{supp}(f)}|\nabla f(z)| -
		\inf_{y \in \mathrm{supp}(f)} \mathcal{H}(y,\nabla f(y)).
		\]
		From the first part of the, see the argument concluding after \eqref{eq:action_integral_representation:Lyapunov_bound}, we find that the solution $x(t)$ remains in the compact set
		\begin{equation} \label{eqn:containment_set_existence}
		K' := \left\{
		z \in E \, \middle| \, \Upsilon(z) \leq C 
		\right\}, \quad C := \Upsilon(x_0) + M + T \sup_x \cH(x,\nabla \Upsilon(x)),
		\end{equation}
		for all $t \in [0,T]$.
		
		\smallskip
		
		\emph{Step 2}: We prove that there exists a solution $x(t)$ of  \eqref{eq:action_integral_representation:subdifferential_eq} on $[0,T]$. 

		Using $F$, we define a new multi-valued vector-field $F'(z)$ that equals $F(z) = \partial_p \mathcal{H}(z,\nabla f(z))$ inside $K'$, but equals $\{0\}$ outside a neighborhood of $K$. This can e.g. be achieved by multiplying with a smooth cut-off function $g_{K'} : E \to [0,1]$ that is equal to one on $K'$ and zero outside of a neighborhood of $K'$.
		\smallskip
		
		The field $F'$ satisfies (F1), (F2) and (F3) from above, and hence there exists an absolutely continuous path $y : [0,\infty) \to E$ such that $y(0) = x_0$ and for almost every $t \geq 0$,
		\[
		\dot{y}(t) \in F'(y(t)).
		\]
		By the estimate established in step 1 and the fact that $\Upsilon(\gamma(t)) \leq C$ for any $0 \leq t \leq T$, it follows from the argument as shown above in \eqref{eq:action_integral_representation:Lyapunov_bound} that the solution $y$ stays in $K'$ up to time $T$. Since on $K'$, we have $F' = F$, this implies that setting $x = y|_{[0,T]}$, we obtain a solution $x(t)$ of \eqref{eq:action_integral_representation:subdifferential_eq} on the time interval $[0,T]$.
	\end{proof}
	\section{Verification for examples of Hamiltonians}
	\label{section:verification-for-examples-of-Hamiltonians}

	In this section, we verify the conditions on $\Lambda$ and $\mathcal{I}$ for the example Hamiltonians of Section~\ref{section:examples-H}. Since the conditions on the functions~$\Lambda$ and~$\mathcal{I}$ are independent of each other, we verify these conditions separately. In Section \ref{section:verify-ex:cost-functions}, we consider Assumption \ref{assumption:results:regularity_I} for $\cI$. In Sections \ref{section:verify-ex:functions-Lambda}, we consider Assumption \ref{assumption:results:regularity_of_V} for $\Lambda$. The continuity estimates will be verified separately in Section \ref{section:verification_of_continuity_estimate}.
	\subsection{Verifying assumptions for cost functions $\mathcal{I}$}
	\label{section:verify-ex:cost-functions}
	We verify Assumption~\ref{assumption:results:regularity_I} for two types of cost functions $\mathcal{I}(x,\theta)$, corresponding to the examples of Section~\ref{section:examples-H}.
	\smallskip
	
	We start by considering the case in which the cost function is the large-deviation rate function for the occupation-time measures of jump process taking values in a finite set $\{1,\dots,J\}$ (e.g.~\cite{donsker1975asymptoticI,denHollander2000}). We follow this example in Proposition~\ref{prop:verify:DV-functional-of-drift-diffusion} in which the cost function stems from occupation-time large deviations of a drift-diffusion process on a compact manifold, see e.g. \cite{DonskerVaradhan75,Pi07}. We expect these results to extend also to non-compact spaces, but we feel this is better suited for a separate work.
	\begin{proposition}[Donsker-Varadhan functional for jump processes]
		\label{prop:verify:DV-for-Jumps}
		Consider a finite set $F = \{1,\dots,J\}$ and let $\Theta := \mathcal{P}(\{1,\dots,J\})$ be the set of probability measures on $F$. For $x\in E$, let $L_x : C_b(F) \rightarrow C_b(F)$ be the operator given by
		\begin{equation*}
		L_x f(i) := \sum_{j=1}^Jr(i,j,x)\left[f(j)-f(i)\right],\quad f :\{1,\dots,J\}\to\mathbb{R}.
		\end{equation*}
		Suppose that the rates $r:\{1,\dots,J\}^2\times E \rightarrow \bR^+$ are continuous as a function on $E$ and moreover satisfy the following:
		\begin{enumerate}[label=(\roman*)]
			\item For any $x\in E$, the matrix $R(x)$ with entries $R(x)_{ij} := r(i,j,x)$ for $i\neq j$ and $R(x)_{ii} = -\sum_{j\neq i}r(i,j,x)$ is irreducible.
			\item For each pair $(i,j)$, we either have $r(i,j,\cdot)\equiv 0$ or for each compact set $K\subseteq E$, it holds that
			\begin{equation*}
			r_{K}(i,j) := \inf_{x\in K}r(i,j,x) > 0.
			\end{equation*}
		\end{enumerate}
		Then the Donsker-Varadhan functional $\mathcal{I}:E\times\Theta \rightarrow \bR^+$ defined by
		\begin{equation*}
		\mathcal{I}(x,\theta) := \sup_{w\in\mathbb{R}^J}\sum_{ij}r(i,j,x)\theta_i \left[1-e^{w_j-w_i}\right]
		\end{equation*}
		satisfies Assumption~\ref{assumption:results:regularity_I}.
	\end{proposition}

	\begin{proof}[Proof]
		\underline{\ref{item:assumption:I:lsc}}: For a fixed vector $w\in\mathbb{R}^J$, the map
		\begin{equation*}
		(x,\theta)\mapsto \sum_{ij}r(i,j,x)\theta_i \left[1-e^{w_j-w_i}\right]
		\end{equation*}
		is continuous on $E\times\Theta$. Hence $\mathcal{I}(x,\theta)$ is lower semicontinuous as the supremum over continuous functions.
		\smallskip
		
		\underline{\ref{item:assumption:I:zero-measure}}: Let $x\in E$. First note that for all $\theta$, the choice $w = 0$ implies that $\cI(x,\theta) \geq 0$. By the irreducibility assumption on the rates $r(i,j,x)$, there exists a unique measure $\theta_x\in\Theta$ such that for any $f:\{1,\dots,J\}\to\mathbb{R}$,
		\begin{equation} \label{eqn:example_jump_DV_stationarity}
		\sum_i L_x f(i) \theta_x(i)=0.
		\end{equation}
		We establish $\cI(x,\theta_x) = 0$. Let $w \in \mathbb{R}^J$. By the elementary estimate
		\begin{equation*}
		\left(1-e^{b - a}\right)\leq -(b-a) \quad \text{ for all } \; a,b > 0,
		\end{equation*}
		we obtain
		\begin{align*}
		\sum_{ij}r(i,j,x) \theta_x(i) \left(1-e^{w_j - w_i}\right)
		&\leq  -\sum_{ij}r(i,j,x) \theta_x(i) \left(w_j - w_i \right)\\
		&= -\sum_i (L_x w)(i) \theta_x(i) \overset{\eqref{eqn:example_jump_DV_stationarity}}{=} 0
		\end{align*}
		Since $\mathcal{I} \geq 0$, this implies $\mathcal{I}(x,\theta_x) = 0$.
		\smallskip
		
		\underline{\ref{item:assumption:I:compact-sublevelsets}}: Any closed subset of $\Theta$ is compact.
		\smallskip
		
		\underline{\ref{item:assumption:I:finiteness}}: Let $x_n\to x$ in $E$. It follows that the sequence is contained in some compact set $K \subseteq E$ that contains the $x_n$ and $x$ in its interior. For any $y\in K$,
		\begin{equation*}
		\mathcal{I}(y,\theta) \leq \sum_{ij, i \neq j} r(i,j,y) \theta_i \leq
		\sum_{ij, i\neq j} r(i,j,y) \leq
		\sum_{ij, i \neq j} \bar{r}_{ij}, \quad \bar{r}_{ij} := \sup_{y \in K} r(i,j,y).
		\end{equation*}
		Hence $\mathcal{I}$ is uniformly bounded on $K\times\Theta$, and~\ref{item:assumption:I:finiteness} follows with $U_x$ the interior of $K$.
		\smallskip
		
		\underline{\ref{item:assumption:I:equi-cont}}: Let $d$ be some metric that metrizes the topology of $E$. We will prove that for any compact set $K\subseteq E$ and $\varepsilon > 0$ there is some $\delta > 0$ such that for all $x,y \in K$ with $d(x,y) \leq \delta$ and for all $\theta \in \mathcal{P}(F)$, we have
		\begin{equation} \label{eqn:proof_equi_cont_I}
		|\mathcal{I}(x,\theta) - \mathcal{I}(y,\theta)| \leq \varepsilon.
		\end{equation}
		Let $x,y \in K$. By continuity of the rates the $\mathcal{I}(x,\cdot)$ are uniformly bounded for $x \in K$:
		\begin{equation*}
		0 \leq \mathcal{I}(x,\theta) \leq \sum_{ij, i \neq j} r(i,j,x) \theta_i \leq
		\sum_{ij, i\neq j} r(i,j,x) \leq
		\sum_{ij, i \neq j} \bar{r}_{ij}, \quad \bar{r}_{ij} := \sup_{x \in K} r(i,j,x).
		\end{equation*}
		For any $n \in \mathbb{N}$, there exists $w^n \in \mathbb{R}^J$ such that
		\begin{equation*}
		0 \leq \mathcal{I}(x,\theta) \leq \sum_{ij, i \neq j} r_{ij}(x) \theta_i (1 - e^{w^n_j - w^n_i}) + \frac{1}{n}.
		\end{equation*}
		By reorganizing, we find for all bonds $(a,b)$ the bound 
		\begin{equation*}
		\theta_a e^{w^n_b - w^n_a} \leq
		\frac{1}{r_{K,a,b}} \left[
		\sum_{ij, i \neq j, a \neq i, b \neq j} r(i,j,x)\theta_i + \frac{1}{n}
		\right]
		\leq \frac{1}{r_{K,a,b}} \left[
		\sum_{ij, i \neq j} \bar{r}_{ij} + \frac{1}{n}
		\right].
		\end{equation*}
		Thereby, evaluating in $\mathcal{I}(y,\theta)$ the same vector $w^n$ to estimate the supremum,
		\begin{align*}
		& \mathcal{I}(x,\theta) - \mathcal{I}(y,\theta) \\
		&\leq \frac{1}{n} + \sum_{ab, a\neq b} r(a,b,x) \theta_a (1 - e^{w^n_b - w^n_a}) - \sum_{ab, a\neq b} r(a,b,y) \theta_a (1 - e^{w^n_b - w^n_a})	\\
		&\leq \frac{1}{n} + \sum_{ab, a\neq b} |r(a,b,x) - r(a,b,y)| \theta_a + \sum_{ab, a\neq b} |r(a,b,y) - r(a,b,x)| \theta_a e^{w^n_b - w^n_a}	\\
		&\leq \frac{1}{n} + \sum_{ab, a\neq b}|r(a,b,x) - r(a,b,y)| \left(
		1 + \frac{1}{r_{K,a,b}} \left[\sum_{ij, i \neq j} \bar{r}_{ij} + 1 \right]
		\right) 
		\end{align*}
		We take $n \to \infty$ and use that the rates $x \mapsto r(a,b,x)$ are continuous, and hence uniformly continuous on compact sets, to obtain \eqref{eqn:proof_equi_cont_I}.
	\end{proof}
	\begin{proposition}[Donsker-Varadhan functional for drift-diffusions]
		\label{prop:verify:DV-functional-of-drift-diffusion}
		Let $F$ be a smooth compact Riemannian manifold without boundary and set $\Theta:=\mathcal{P}(F)$, the set of probability measures on $F$. For $x\in E$, let $L_x : C^2(F) \subseteq C_b(F) \rightarrow C_b(F)$ be the second-order elliptic operator that in local coordinates is given by
		\begin{equation*}
		L_x = \frac{1}{2}\nabla\cdot\left(a_x \nabla\right) + b_x\cdot \nabla,
		\end{equation*}
		where $a_x$ is a positive definite matrix and $b_x$ is a vector field having smooth entries $a_x^{ij}$ and $b_x^i$ on $F$. Suppose that for all $i,j$ the maps 
		\begin{equation} \label{eqn:examples_DV_diffusion_continuous_coefficients}
		x \mapsto a_x^{i,j}(\cdot), \qquad x \mapsto b_x^i(\cdot)
		\end{equation}	
		are continuous as functions from $E$ to $C_b(F)$, where we equip $C_b(F)$ with the supremum norm. Then the functional $\mathcal{I}:E\times\Theta\to[0,\infty]$ defined by
		\begin{equation*}
		\mathcal{I}(x,\theta) := \sup_{\substack{u\in \mathcal{D}(L_x)\\u>0}}\left[ -\int_F \frac{L_xu}{u}\,d\theta\right]
		\end{equation*}
		satisfies Assumption~\ref{assumption:results:regularity_I}.
	\end{proposition}
	\begin{proof}[Proof]
		\underline{\ref{item:assumption:I:lsc}}: For any fixed function $u\in\mathcal{D}(L_x)$ that is strictly positive on $F$, the function $(-L_xu/u)$ is continuous on $F$. For any fixed $u$ it follows by~\eqref{eqn:examples_DV_diffusion_continuous_coefficients} and compactness of $F$ that 
		\begin{equation*}
		(x,\theta)\mapsto -\int_F \frac{L_xu}{u}\,d\theta
		\end{equation*}
		is continuous on $E\times\Theta$. As a consequence $\mathcal{I}(x,\theta)$ is lower semicontinuous as the supremum over continuous functions.
		\smallskip
		
		\underline{\ref{item:assumption:I:zero-measure}}: Let $x\in E$. The stationary measure $\theta_x\in\Theta$ satisfying
		\begin{equation}\label{eq:proof:verify-I-diffusion:I2}
		\int_F L_xg(z)\,d\theta_x(z) = 0\quad \text{for all}\;g\in \mathcal{D}(L_x)
		\end{equation}
		is the minimizer of $\mathcal{I}(x,\cdot)$, that is $\mathcal{I}(x,\theta_x) = 0$. This follows by considering the Hille-Yosida approximation $L_x^\varepsilon$ of $L_x$ and using the same argument (using $w = \log u$) as in Proposition~\ref{prop:verify:DV-for-Jumps} for these approximations. For any $u>0$ and for any~$\varepsilon>0$, 
		\begin{align*}
		-\int_F \frac{L_xu}{u}\,d\theta & = -\int_F \frac{L^\varepsilon_xu}{u}\,d\theta + \int_F \frac{(L^\varepsilon_x-L_x)u}{u}\,d\theta\\
		&\leq -\int_F \frac{L^\varepsilon_xu}{u}\,d\theta + \frac{1}{\inf_F u} \|(L_x^\varepsilon-L_x)u\|_F\\
		&\leq -\int_F L^\varepsilon_x \log(u)\,d\theta + o(1)_{\varepsilon\to 0}.
		\end{align*}
		Sending $\varepsilon\to 0$ and then using~\eqref{eq:proof:verify-I-diffusion:I2} gives~\ref{item:assumption:I:zero-measure}.
		\smallskip
		
		\underline{\ref{item:assumption:I:compact-sublevelsets}}: Since $\Theta = \mathcal{P}(F)$ is compact, any closed subset of $\Theta$ is compact. Hence any union of sub-level sets of $\mathcal{I}(x,\cdot)$ is relatively compact in~$\Theta$.
		\smallskip

		\underline{\ref{item:assumption:I:finiteness}}: Let $x_n\to x$ in $E$ and $\theta_n$ be a sequence in $\Theta$, and suppose that $\mathcal{I}(x_n,\theta_n)\leq M$ for some constant $M$ independent of $n$. Let $\dd z$ be the Riemannian measure on $F$. By Pinsky's results in~\cite{pinsky1985evaluating,Pi07}, if $\cI(y,\theta) < \infty$, then the density $\frac{\dd \theta}{\dd z}$ exists. In addition, there are constants $c_1,c_2,c_3,c_4$ depending only on $a_{y},b_{y}$, and not on $\theta$, such that
		\begin{equation} \label{eqn:Pinksy_bootstrap}
		c_1(y)\int_F|\nabla g_\theta |^2\,dz - c_2(y) \leq \mathcal{I}(y,\theta) \leq c_3(y) \int_F|\nabla g_\theta |^2\,dz + c_4(y),
		\end{equation}
		where $g_\theta = (\dd \theta/\dd z)^{1/2}$. In particular, as can be seen by the derivation of~\cite[Eq.~(2.18),~(2.19)]{pinsky1985evaluating}, the constants depend continuously on $y\in E$ by our continuity assumptions on $a_y$ and $b_y$.

		Applying this to our sequences $x_n$ and $\theta_n$, we have
		\begin{equation*}
		\int_F |\nabla g_{\theta_n}|^2\,dz \leq M',
		\end{equation*}
		for a constant $M'$. This implies again by \eqref{eqn:Pinksy_bootstrap} that for any $y$ in some neighborhood of $x$ that
		\begin{equation*}
		\mathcal{I}(y,\theta_{x_n}) \leq C <\infty,
		\end{equation*}
		with a constant independent of $n$.
		\smallskip
		
		\underline{\ref{item:assumption:I:equi-cont}}: Since the coefficients $a_x$ and $b_x$ of the operator $L_x$ depend continuously on $x$, assumption~\ref{item:assumption:I:equi-cont} follows from Theorem~2 of~\cite{Pi07}.
	\end{proof}
	\subsection{Verifying assumptions for functions $\Lambda$}
	\label{section:verify-ex:functions-Lambda}
	We verify Assumption~\ref{assumption:results:regularity_of_V} for three types of functions~$\Lambda$ corresponding to the examples of Section \ref{section:examples-H}. We start with $\Lambda$'s that are given as integrals over quadratic polynomials in $p$.
	\begin{proposition}[Quadratic function $\Lambda$]\label{prop:verify-ex:Lambda_quadratic}
		Let $E=\mathbb{R}^d$ and $\Theta=\mathcal{P}(F)$ for some compact Polish space $F$. Suppose that the function $\Lambda :E\times\mathbb{R}^d\times\Theta\to\mathbb{R}$ is given by
		\begin{equation*} 
		\Lambda(x,p,\theta) = \int_F\ip{a(x,z)p}{p}\,d\theta(z) + \int_F\ip{b(x,z)}{p}\,d\theta(z),
		\end{equation*}
		where $a:E\times F\to\mathbb{R}^{d\times d}$ and $b:E\times F\to\mathbb{R}^d$ are continuous. Suppose that for every compact set $K \subseteq \bR^d$,
		\begin{align*}
		a_{K,min} & := \inf_{x \in K, z \in F, |p|=1} \ip{a(x,z)p}{p} > 0, \\
		a_{K,max} & := \sup_{x \in K, z \in F, |p| = 1} \ip{a(x,z)p}{p} < \infty, \\
		b_{K,max} & := \sup_{x \in K, z \in F, |p|=1} |\ip{b(x,z)}{p}| < \infty.
		\end{align*}
		Furthermore, there exists a constant $L>0$ such that for all $x,y\in E$ and $z\in F$,
		\begin{equation*}
		\|a(x,z)-a(y,z)\| \leq L|x-y|,
		\end{equation*}
		and suppose that the functions $b$ are one-sided Lipschitz continuous. Then Assumption~\ref{assumption:results:regularity_of_V} holds.
	\end{proposition}
	\begin{proof}[Proof]
		\underline{\ref{item:assumption:slow_regularity:continuity}}: Let $(x,p)\in E\times\mathbb{R}^d$. By the boundedness assumptions on $a$ and $b$, 
		\begin{equation*}
		\sup_\theta |\Lambda(x,p,\theta)| \leq a_{\{x\},\text{max}} + b_{\{x\},\text{max}} < \infty,
		\end{equation*}
		and hence the function $\theta\mapsto |\Lambda(x,p,\theta)|$ is bounded on $\mathcal{P}(F)$. Continuity of $\Lambda$ is a consequence of the fact that
		\begin{equation*}
		\Lambda(x,p,\theta) = \int_F V(x,p,z)\,\dd\theta(z)
		\end{equation*}
		is the pairing of a continuous bounded function $V(x,p,\cdot)$ with~$\theta\in\mathcal{P}(F)$.
		\smallskip
		
		\underline{\ref{item:assumption:slow_regularity:convexity}}: Let $x\in E$ and $\theta \in \mathcal{P}(F)$. Convexity of $p\mapsto \Lambda(x,p,\theta)$ follows since $a(x,z)$ is positive definite by assumption. If $p_0 = 0$, then evidently $\Lambda(x,p_0,\theta) = 0$.
		\smallskip
		
		\underline{\ref{item:assumption:slow_regularity:compact_containment}}: We show that the map $\Upsilon : E\to\mathbb{R}$ defined by
		\begin{equation*}
		\Upsilon(x) := \frac{1}{2}\log\left(1 + |x|^2\right)
		\end{equation*}
		is a containment function for $\Lambda$. For any $x\in E$ and $\theta\in\mathcal{P}(F)$, we have
		\begin{align*}
		\Lambda(x,\nabla\Upsilon(x),\theta) &= \int_F \ip{a(x,z)\nabla\Upsilon(x)}{\nabla\Upsilon(x)}\,d\theta(z) + \int_F\ip{b(x,z)}{\nabla\Upsilon(x)}\,d\theta(z)\\
		&\leq a_{\{x\},\text{max}} |\nabla\Upsilon(x)|^2 + b_{\{x\},\text{max}}|\nabla\Upsilon(x)|\\
		&\leq C (1+|x|) \frac{x^2}{(1+x^2)^2} + C(1+|x|) \frac{x}{(1+x^2)},
		\end{align*}
		and the boundedness condition follows with the constant 
		\begin{equation*}
		C_\Upsilon := C \,\sup_x (1+|x|) \left[\frac{x^2}{(1+x^2)^2} + \frac{x}{(1+x^2)} \right] <\infty.
		\end{equation*}
		\smallskip
		
		\underline{\ref{item:assumption:slow_regularity:continuity_estimate}}: By the assumption on $a(x,z)$, the function $\Lambda$ is uniformly coercive in the sense that for any compact set $K\subseteq E$, 
		\begin{equation*}
		\inf_{x\in K, \theta\in\Theta}\Lambda(x,p,\theta) \to \infty \quad \text{ as }\; |p|\to \infty,
		\end{equation*}
		and the continuity estimate follows by Proposition~\ref{proposition:continuity_estimate_coercivity}.
		\smallskip
		
		\underline{\ref{item:assumption:slow_regularity:controlled_growth}}: Let $K\subseteq E$ be compact. We have to show that there exist constants $M, C_1, C_2 \geq 0$  such that for all $x \in K$, $p \in \mathbb{R}^d$ and all $\theta_1,\theta_2 \in \mathcal{P}(F)$, we have
		\begin{equation} \label{eqn:growth_bound_general}
		\Lambda(x,p,\theta_1)
		\leq
		\max \left\{M, C_1 \Lambda(x,p,\theta_2) + C_2 \right\}.
		\end{equation}
		Fix $\theta_1,\theta_2 \in \mathcal{P}(F)$.	We have for $x \in K$
		\begin{equation*}
		\int \ip{a(x,z)p}{p} d\theta_1(z) \leq \frac{a_{K,max}}{a_{K,min}} \int \ip{a(x,z)p}{p} d\theta_2(z)
		\end{equation*}
		In addition, as $a_{K,min} > 0$ and $b_{K,max} < \infty$ we have for any $C > 0$ and sufficiently large $|p|$ that
		\begin{equation*}
		\int \ip{b(x,z)}{p} \,d\theta_1(z) - (C+1)\int \ip{b(x,z)}{p} \,d\theta_2(z) \leq C \int \ip{a(x,z)p}{p} \,d\theta_2(z)
		\end{equation*}
		Thus, for sufficiently large $|p|$ (depending on $C$) we have
		\begin{equation*}
		\Lambda(x,p,\theta_1) \leq (1+C) \Lambda(x,p,\theta_2).
		\end{equation*}
		Fix a $C =: C_1$ and denote the set of `large' $p$ by $S$. The map $(x,p,\theta) \mapsto \Lambda(x,p,\theta)$ is bounded on $K \times \times S^c\times \Theta$. Thus, we can find a constant $C_2$ such that \eqref{eqn:growth_bound_general} holds.
	\end{proof}
	We proceed with an example in which $\Lambda$ depends on $p$ through exponential functions. Let $q \in \mathbb{N}$ be an integer and 
	\begin{equation*}
	\Gamma := \left\{(a,b)\,:\,a,b\in\{1,\dots,q\}, \,a\neq b\right\}
	\end{equation*}
	be the set of oriented edges in $\{1,\dots,q\}$.
	\begin{proposition}[Exponential function $\Lambda$]\label{prop:verify-ex:Lambda_exponential}
		Let $E\subseteq \mathbb{R}^d$ be the embedding of $E = \cP(\{1,\dots,q\}) \times (\bR^+)^{|\Gamma|}$ and $\Theta$ be a topological space. Suppose that $\Lambda$ is given by
		\begin{equation*} 
		\Lambda((\mu,w),p,\theta) = \sum_{(a,b) \in \Gamma} v(a,b,\mu,\theta)\left[\exp\left\{p_b - p _a + p_{(a,b)} \right\} - 1 \right]
		\end{equation*}
		where $v$ is a proper kernel in the sense of Definition~\ref{definition:proper_kernel}. Suppose in addition that there is a constant $C > 0$ such that for all $(a,b) \in \Gamma$ such that $v(a,b, \cdot,\cdot) \neq 0$ we have
		\begin{equation}\label{eq:prop-verify-Lambda-exp:boundedness-kernel}
		\sup_{\mu} \sup_{\theta_1,\theta_2} \frac{v(a,b,\mu,\theta_1)}{v(a,b,\mu,\theta_2)} \leq C.
		\end{equation}
		Then $\Lambda$ satisfies Assumption~\ref{assumption:results:regularity_of_V}.
	\end{proposition}
		Similar to previous proposition, the assumptions on~$\Lambda$ are satisfied if~$\Theta = \mathcal{P}(F)$ for some Polish space $F$, and if~$v(a,b,\mu,\theta) = \mu(a) \int r(a,b,\mu,z) \theta(\dd z)$ and there are constants~$0 < r_{min} \leq r_{max} < \infty$ such that for all~$(a,b) \in \Gamma$ such that~$\sup_{\mu,z} r(a,b,\mu,z) > 0$, we have
		\begin{equation*}
		r_{min} \leq \inf_{z} \inf_{\mu} r(a,b,\mu,z) \leq \sup_{z} \sup_{\mu} r(a,b,\mu,z) \leq r_{max}.
		\end{equation*}
		Regarding~\eqref{eq:prop-verify-Lambda-exp:boundedness-kernel}, for $(a,b) \in \Gamma$ for which $v(a,b,\cdot,\cdot)$ is non-trivial, we have
		\begin{equation*}
		\frac{v(a,b,\mu,\theta_1)}{v(a,b,\mu,\theta_2)} = \frac{\int r(a,b,\mu,z) \theta_1(\dd z)}{\int r(a,b,\mu,z) \theta_2(\dd z)} \leq \frac{r_{max}}{r_{min}}.
		\end{equation*}
	\begin{proof}[Proof of Proposition~\ref{prop:verify-ex:Lambda_exponential}]
		\underline{\ref{item:assumption:slow_regularity:continuity}}: The function $\Lambda$ is continuous as the sum of continuous functions. Boundedness of $\Lambda$ as a function of $\theta$ follows from the boundedness assumption~\eqref{eq:prop-verify-Lambda-exp:boundedness-kernel}.
		\smallskip
		
		\underline{\ref{item:assumption:slow_regularity:convexity}}: Convexity of $\Lambda$ as a function of $p$ follows from the fact that $\Lambda$ is a finite sum of convex functions, and $\Lambda(x,0,\theta)=0$ is evident.
		\smallskip
		
		\underline{\ref{item:assumption:slow_regularity:compact_containment}}: The function $\Upsilon : E\to\mathbb{R}$ defined by
		\begin{equation*}
		\Upsilon(\mu,w) := \sum_{(a,b)\in\Gamma}\log\left[1 + w_{(a,b)}\right]
		\end{equation*}
		is a containment function for $\Lambda$ (an explicit verification is given in~\cite{Kr17}).
		\smallskip
		
		\underline{\ref{item:assumption:slow_regularity:continuity_estimate}}: The continuity estimate is the content of Proposition~\ref{proposition:continuity_estimate_directional_extended} below.
		\smallskip
		
		\underline{\ref{item:assumption:slow_regularity:controlled_growth}}: Note that	
		\begin{align*}
		\Lambda((\mu,w),\theta_1,p) & \leq \sum_{(a,b)\in \Gamma} v(a,b,\mu,\theta_1) e^{p_{a,b} + p_b - p_a} \\
		& \leq C \sum_{(a,b)\in \Gamma} v(a,b,\mu,\theta_2) e^{p_{a,b} + p_b - p_a}  \\
		& \leq C \sum_{(a,b)\in \Gamma} v(a,b,\mu,\theta_2) \left[e^{p_{a,b} + p_b - p_a} - 1 \right] + C_2 .
		\end{align*}
		Thus the estimate holds with $M = 0$, $C_1 = C$ and $C_2 = \sup_{\mu,\theta} \sum_{a,b} v(a,b,\mu,\theta)$.
	\end{proof}
	\subsection{Verifying the continuity estimate} \label{section:verification_of_continuity_estimate}
	
	With the exception of  the verification of the continuity estimate in Assumption \ref{assumption:results:regularity_of_V} the verification in Section \ref{section:verify-ex:functions-Lambda} is straightforward. On the other hand, the continuity estimate is an extension of the comparison principle, and is therefore more complex. We verify the continuity estimate in three contexts, which we hope illustrates that the continuity estimate follows from essentially the same arguments as the standard comparison principle. We will do this for:
	\begin{itemize}
		\item Coercive Hamiltonians
		\item One-sided Lipschitz Hamiltonians
		\item Hamiltonians arising from large deviations of empirical measures.
	\end{itemize}
	This list is not meant to be an exhaustive list, but to illustrate that the continuity estimate is a sensible extension of the comparison principle, which is satisfied in a wide range of contexts. In what follows, $E\subseteq \mathbb{R}^d$ is a Polish subset and $\Theta$ a topological space.
	\begin{proposition}[Coercive $\Lambda$] \label{proposition:continuity_estimate_coercivity}
		Let $\Lambda : E \times \bR^d \times \Theta \rightarrow \bR$ be continuous and uniformly coercive: that is, for any compact $K \subseteq E$ we have
		\begin{equation*}
		\inf_{x \in K, \theta\in\Theta} \Lambda(x,p,\theta) \to \infty \quad \mathrm{as} \; |p| \to \infty.
		\end{equation*}
		Then the continuity estimate holds for $\Lambda$ with respect to any penalization function $\Psi$.
	\end{proposition}
	
	\begin{proof}
		Let $\Psi(x,y) = \tfrac{1}{2}(x-y)^2$. Let $(x_{\alpha,\varepsilon},y_{\alpha,\varepsilon},\theta_{\varepsilon,\alpha})$ be fundamental for $\Lambda$ with respect to $\Psi$. Set $p_{\alpha,\varepsilon} = \alpha(x_{\varepsilon,\alpha} - y_{\varepsilon,\alpha})$. By the upper bound~\eqref{eqn:control_on_Gbasic_sup}, we find that for sufficiently small $\varepsilon > 0$ there is some $\alpha(\varepsilon)$ such that
		\begin{equation*}
		\sup_{\alpha \geq \alpha(\varepsilon)} \Lambda\left(y_{\varepsilon,\alpha}, p_{\varepsilon,\alpha}, \theta_{\varepsilon,\alpha}\right) < \infty.
		\end{equation*}
		As the variables $y_{\alpha,\varepsilon}$ are contained in a compact set by property (C1) of fundamental collections of variables, the uniform coercivity implies that the momenta $p_{\varepsilon,\alpha}$ for $\alpha \geq \alpha(\varepsilon)$ remain in a bounded set. Thus, we can extract a subsequence $\alpha'$ such that $(x_{\varepsilon,\alpha'},y_{\varepsilon,\alpha'},p_{\varepsilon,\alpha'},\theta_{\varepsilon,\alpha'})$ converges to $(x,y,p,\theta)$ with $x = y$ due to property (C2) of fundamental collections of variables. By continuity of $\Lambda$ we find
		\begin{align*}
		& \liminf_{\alpha \rightarrow \infty} \Lambda\left(x_{\varepsilon,\alpha}, p_{\varepsilon,\alpha},\theta_{\varepsilon,\alpha}\right) - \Lambda\left(y_{\alpha,\varepsilon},p_{\varepsilon,\alpha},\theta_{\varepsilon,\alpha}\right) \\
		& \leq \lim_{\alpha'\rightarrow \infty} \Lambda\left(x_{\varepsilon,\alpha'}, p_{\varepsilon,\alpha'},\theta_{\varepsilon,\alpha'}\right) - \Lambda\left(y_{\varepsilon,\alpha'},p_{\varepsilon,\alpha'},\theta_{\varepsilon,\alpha'}\right) = 0
		\end{align*}
		establishing the continuity estimate.
	\end{proof}

	\begin{proposition}[One-sided Lipschitz $\Lambda$]  \label{proposition:continuity_estimate_lipschitz}
		Let $\Lambda : E \times \bR^d \times \Theta\rightarrow \bR$ satisfy
		\begin{equation} \label{eqn:one_sided_Lipschitz_G}
		\Lambda(x,\alpha(x-y),\theta) - \Lambda(y,\alpha(x-y),\theta) \leq  c(\theta) \omega( \alpha (x-y)^2)
		\end{equation}
		for some collection of constants $c(\theta)$ satisfying $\sup_\theta c(\theta) < \infty$ and a function $\omega : \bR^+ \rightarrow \bR^+$ satisfying $\lim_{\delta \downarrow 0} \omega(\delta) = 0$.
		
		Then the continuity estimate holds for $\Lambda$ with respect to $\Psi(x,y) = \tfrac{1}{2}(x-y)^2$.
	\end{proposition}
	
	\begin{proof} [Proof]
		Let $\Psi(x,y) = \tfrac{1}{2}(x-y)^2$. Let $(x_{\alpha,\varepsilon},y_{\alpha,\varepsilon},\theta_{\varepsilon,\alpha})$ be fundamental for $\Lambda$ with respect to $\Psi$. Set $p_{\alpha,\varepsilon} = \alpha(x_{\varepsilon,\alpha} - y_{\varepsilon,\alpha})$. We find
		\begin{align*}
		& \liminf_{\alpha \rightarrow \infty} \Lambda\left(x_{\varepsilon,\alpha}, p_{\varepsilon,\alpha}, \theta_{\varepsilon,\alpha}\right) - \Lambda\left(y_{\alpha,\varepsilon},p_{\varepsilon,\alpha},\theta_{\varepsilon,\alpha}\right) \\
		& \leq \liminf_{\alpha\rightarrow \infty} c(\theta) \omega( \alpha (x-y)^2)
		\end{align*}
		which equals $0$ as $\sup_\theta c(\theta) < \infty$, $\lim_{\delta \downarrow 0} \omega(\delta) = 0$ and property (C1) of a fundamental collection of variables.
	\end{proof}

	For the empirical measure of a collection of independent processes one obtains maps $\Lambda$ that are neither uniformly coercive nor Lipschitz. Also in this context one can establish the continuity estimate. We treat a simple 1d case and then state a more general version for which we refer to \cite{Kr17}.
	
	\begin{proposition} \label{proposition:continuity_estimate_directional}
		Suppose that $E = [-1,1]$ and that $\Lambda(x,p,\theta)$ is given by
		\begin{equation*}
		\Lambda(x,p,\theta) = \frac{1-x}{2} c_+(\theta) \left[e^{2p} -1\right] +  \frac{1+x}{2} c_-(\theta) \left[e^{-2p} -1\right]
		\end{equation*}
		with $c_-,c_+$ non-negative functions of $\theta$.	Then the continuity estimate holds for $\Lambda$ with respect to $\Psi(x,y) = \tfrac{1}{2}(x-y)^2$.
	\end{proposition}

	\begin{proof}[Proof]
		Let $\Psi(x,y) = \tfrac{1}{2}(x-y)^2$. Let $(x_{\alpha,\varepsilon},y_{\alpha,\varepsilon},\theta_{\varepsilon,\alpha})$ be fundamental for $\Lambda$ with respect to $\Psi$. Set $p_{\alpha,\varepsilon} = \alpha(x_{\varepsilon,\alpha} - y_{\varepsilon,\alpha})$. We have
		\begin{align*}
		& \Lambda\left(x_{\varepsilon,\alpha}, p_{\varepsilon,\alpha}, \theta_{\varepsilon,\alpha}\right) - \Lambda\left(y_{\alpha,\varepsilon},p_{\varepsilon,\alpha},\theta_{\varepsilon,\alpha}\right) \\
		& = \frac{y_{\varepsilon,\alpha}-x_{\varepsilon,\alpha}}{2} c_+(\theta_{\varepsilon,\alpha}) \left[e^{2p_{\varepsilon,\alpha}} -1\right] + \frac{x_{\varepsilon,\alpha}-y_{\varepsilon,\alpha}}{2} c_-(\theta_{\varepsilon,\alpha}) \left[e^{-2p_{\varepsilon,\alpha}} -1\right]
		\end{align*}
		Now note that $y_{\varepsilon,\alpha}-x_{\varepsilon,\alpha}$ is positive if and only if $e^{2p_{\varepsilon,\alpha}} -1$ is negative so that the first term is bounded above by $0$. With a similar argument the second term is bounded above by $0$. Thus the continuity estimate is satisfied.
	\end{proof}
	\begin{proposition} \label{proposition:continuity_estimate_directional_extended}
		Suppose $E = \cP(\{1,\dots,q\} \times (\bR^+)^\Gamma$ and suppose that $\Lambda$ is given by
		\begin{equation*} 
		\Lambda((\mu,w),\theta,p) = \sum_{(a,b) \in \Gamma} v(a,b,\mu,\theta)\left[\exp\left\{p_b - p _a + p_{(a,b)} \right\} - 1 \right]
		\end{equation*}
		where $v$ is a proper kernel. Then the continuity estimate holds for $\Lambda$ with respect to penalization functions
		\begin{align*}
		\Psi_1(\mu,\hat{\mu}) & := \frac{1}{2} \sum_{a} ((\hat{\mu}(a) - \mu(a))^+)^2, \\
		\Psi_2(w,\hat{w}) & := \frac{1}{2} \sum_{(a,b) \in \Gamma} (w_{(a,b)} - \hat{w}_{(a,b)})^2.
		\end{align*}
		Here we denote $r^+ = r \vee 0$ for $r \in \bR$.
	\end{proposition}
	
	In this context, one can use coercivity like in Proposition \ref{proposition:continuity_estimate_coercivity} in combination with directional properties used in the proof of Proposition \ref{proposition:continuity_estimate_directional} above.
	The proof of this proposition can be carried out exactly as the proof of~\cite[Theorem~3.8]{Kr17}. Namely at any point, a converging subsequence is constructed, and the variables~$\alpha$ need to be chosen such that we also get convergence of the measures~$\theta_{\varepsilon,\alpha}$ in~$\cP(F)$.
\chapter{Gradient Flow to Non-Gradient-Flow}
\label{chapter:GF-to-NGF}
\section{Introduction}
\subsection{Diffusion in an asymmetric potential landscape}
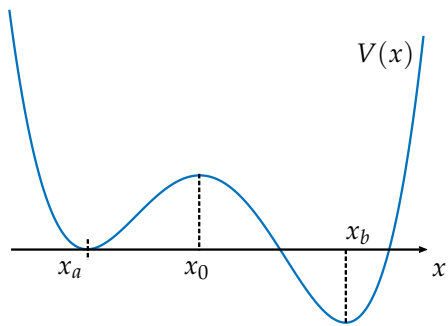
\begin{figure}[h!]
	\labellist
	\pinlabel $x$ at 1600 200
	\pinlabel $x_a$ at 230 200
	\pinlabel $x_0$ at 700 200
	\pinlabel $x_b$ at 1300 330
	\pinlabel $V(x)$ at 1400 1000
	\endlabellist
	\centering
	\includegraphics[scale=.1]{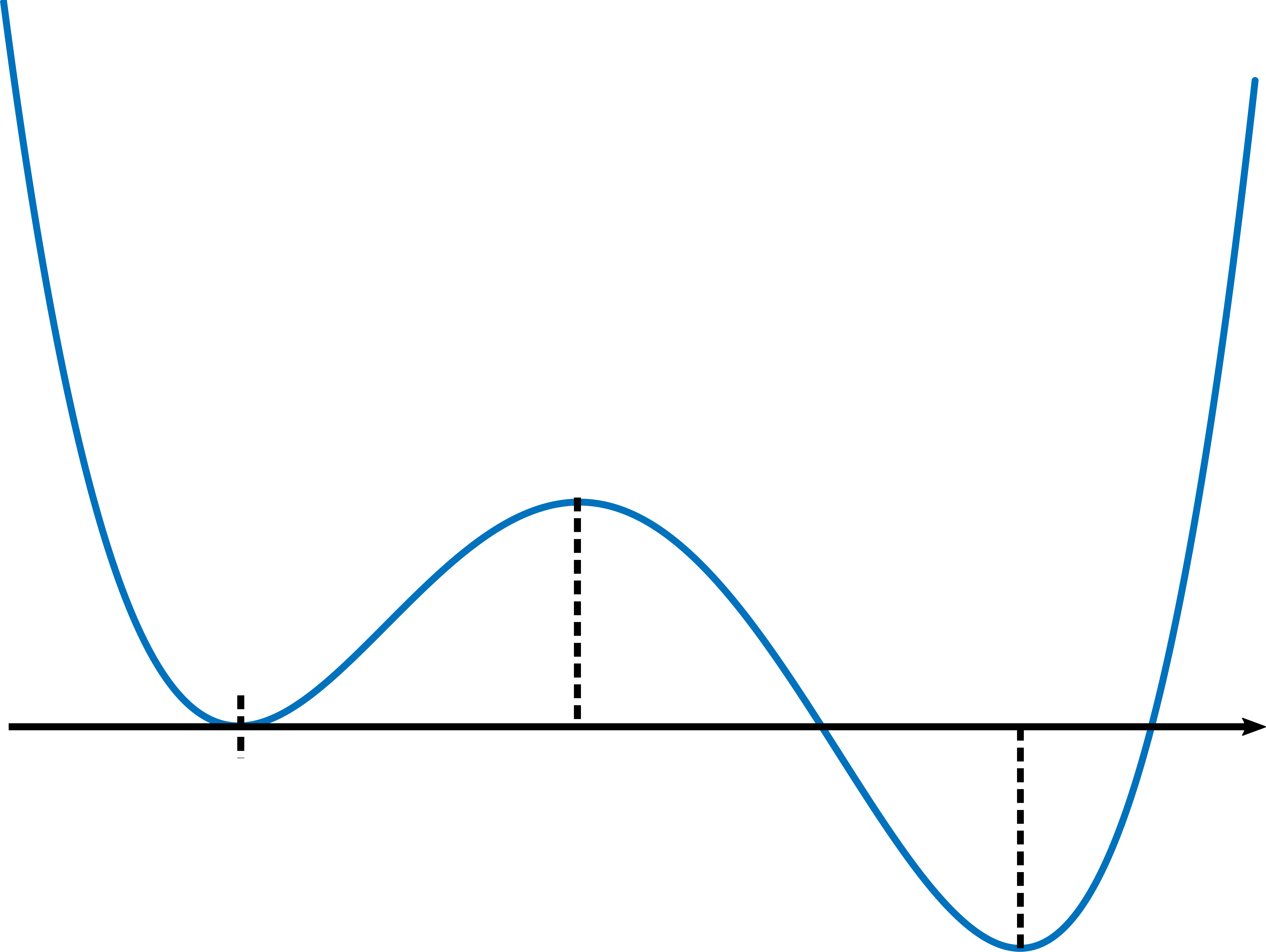}
	\caption{The typical asymmetric potential $V(x)$.}
	\label{GF_NGF:fig:asymmetric-doublewell-potential}
\end{figure}

Our main interest in this chapter is the family of Fokker-Planck equations in one dimension defined by 
\begin{equation}\label{GF_NGF:eq:intro:upscaled-FP}
\partial_t \rho_\varepsilon = \tau_\varepsilon\left[ \varepsilon\, \Delta \rho_\varepsilon + \mathrm{div} \left( \rho_\varepsilon \nabla V\right)\right], 
\quad t\geq 0,\, x\in\mathbb{R}.
\end{equation}
Here, we take an asymmetric double-well potential~$V:\mathbb{R}\to\mathbb{R}$ as depicted in Figure~\ref{GF_NGF:fig:asymmetric-doublewell-potential}.
\smallskip

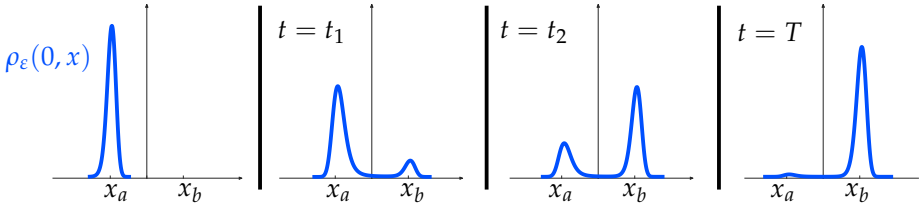
\begin{figure}[h!]
	\labellist
	\pinlabel $x_a$ at 130 0
	\pinlabel $x_b$ at 240 0
	\pinlabel $x_a$ at 460 0
	\pinlabel $x_b$ at 570 0
	\pinlabel $x_a$ at 790 0
	\pinlabel $x_b$ at 900 0
	\pinlabel $x_a$ at 1120 0
	\pinlabel $x_b$ at 1230 0
	\pinlabel $t=t_1$ at 420 240
	\pinlabel $t=t_2$ at 750 240
	\pinlabel $t=T$ at 1100 240
	\pinlabel {\color{dark_blue}{$\rho_\varepsilon(0,x)$}} at 30 200
	\endlabellist
	\centering
	\includegraphics[scale=.25]{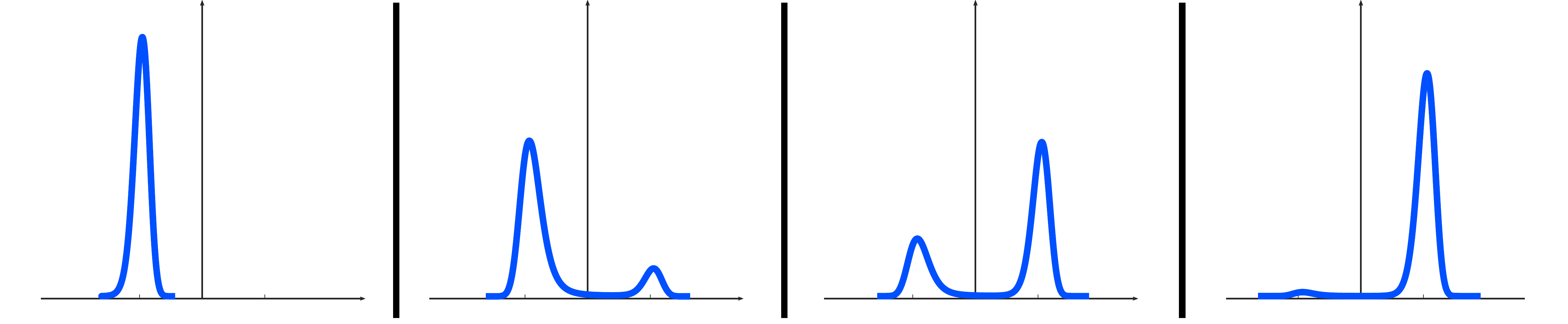}
	\caption{The time-evolution of a solution~$\rho_\varepsilon(t,x)$ to~\eqref{GF_NGF:eq:intro:upscaled-FP} whose initial distribution is supported solely on the left. Time is increasing from left to right. At final time, the solution os close to the equilibrium distribution, which is given by~$\exp\{-V(x)/\varepsilon\}$ up to normalization. The smaller the value of~$\varepsilon$, the sharper the equilibrium distribution concentrates around the global minimum~$x_b$.}
	\label{GF_NGF:fig:intro:evolution-FP}
\end{figure}

A typical solution~$\rho_\varepsilon(t,x)$ is displayed in Figure~\ref{GF_NGF:fig:intro:evolution-FP}, showing a mass transition from left to right. There are two parameters~$\varepsilon>0$ and~$\tau_\varepsilon>0$ appearing in the Fokker-Planck equation. The parameter~$\varepsilon$ controls how fast mass can transition in between the potential's wells. In general, smaller values of~$\varepsilon$ correspond to larger transition times of mass flowing from left to right. The time-scale parameter~$\tau_\varepsilon$ is chosen such that transitions from the local minimum~$x_a$ to the global minimum~$x_b$ happen at rate of order one. Below, we make our choice of~$\tau_\varepsilon$ precise.
\smallskip

We regard the PDE~\eqref{GF_NGF:eq:intro:upscaled-FP} as derived from a stochastic model for metastability introduced by Kramers to study chemical reactions~\cite{Kramers1940}. The stochastic model he considered is the diffusion process $Y^\varepsilon_t = Y^\varepsilon(t)$ in $\mathbb{R}$ solving
\begin{equation*}
\dd Y^\varepsilon_t = -\nabla V(Y^\varepsilon_t)\,\dd t + \sqrt{2\varepsilon}\;\dd B_t,
\end{equation*}
where~$B_t$ denotes the standard Brownian motion. For example, consider a particle initiall starting in the left minimum~$x_a$ and propagating from left to right. This propagation may model an reaction-event in which a molecule's state changes from a low-energy state~$x_a$ via a high-energy state~$x_0$ to another low-energy state~$x_b$. Then the assumption of asymmetry of the potential~$V$ corresponds to modeling a reaction resulting in a molecule's state whose final energy is lower than its initial energy. The energy barrier that the particle has to overcome,~$V(x_0)-V(x_a)$, corresponds to the minimal \emph{activation energy} for the reaction to occur. Kramers discussed various examples of reactions that may be modeled this way~\cite[Paragraph~6]{Kramers1940}. His interest lied in deriving formulas for the average reaction rates from the average transition time of a particle from~$x_a$ to~$x_b$. In the stochastic model, as~$\varepsilon$ decreases, a transition from~$x_a$ to~$x_b$ becomes more unlikely, and hence the average-time for a transition~$x_a\to x_b$ to occur increases. Kramers derived an asymptotic expression for this average-time,
\begin{equation*}
\mathbb{E}_{x_a}\left[T(x_a \to x_b)\right] = \left[1 + o(1)_{\varepsilon\to 0}\right] \frac{2\pi}{\sqrt{V''(x_a) |V''(x_0)|}} \exp\{\varepsilon^{-1}(V(x_0) - V(x_a))\},
\end{equation*}
which is also known as the \emph{Kramers formula}. It shows that the average transition time scales exponentially with respect to the energy barrier~$V(x_0)-V(x_a)$ and the inverse of the diffusion coefficient,~$\varepsilon^{-1}$. For further details and background on this model, we refer to the monographs on metastability of Bovier and den Hollander~\cite{BovierDenHollander2016}, and of Berglund and Gentz~\cite{BerglundGentz2005}.
\smallskip
 
Motivated by Kramers' formula, we define the time-scale parameter~$\tau_\varepsilon$ by
\begin{equation}\label{GF_NGF:intro:eq:def-time-scale-parameter}
\tau_\varepsilon := \frac{2\pi}{\sqrt{V''(x_a) |V''(x_0)|}} \exp^{\varepsilon^{-1}(V(x_0) - V(x_a))},
\end{equation}
in order to be at a time-scale at which jumps from left to right happen at rate of order one.
One way to motivate the PDE~\eqref{GF_NGF:eq:intro:upscaled-FP} from the small-diffusion process is to speed up the process~$Y(t)$ by exactly that time-scale parameter: consider the upscaled process $X^\varepsilon(t):= Y^\varepsilon(\tau_\varepsilon t)$. Then by Itô calculus,~$X^\varepsilon$ satisfies the SDE
\begin{equation}\label{GF_NGF:eq:intro:upscaled-diffusion-process}
\dd X^\varepsilon_t = - \tau_\varepsilon V'(X^\varepsilon_t)\,\dd t + \sqrt{\tau_\varepsilon} \sqrt{2\varepsilon} \; \dd B_t,
\end{equation}
and the equation~\eqref{GF_NGF:eq:intro:upscaled-FP} is the Fokker-Planck equation for  the transition probabilities $\rho_\varepsilon(t,\dd x) = \mathbb{P}\left[X^\varepsilon_t \in \dd x\right]$.
\smallskip

\begin{figure}[h!]
	\labellist
	\pinlabel $x_a$ at 130 0
	\pinlabel $x_b$ at 240 0
	\pinlabel $x_a$ at 460 0
	\pinlabel $x_b$ at 570 0
	\pinlabel $x_a$ at 790 0
	\pinlabel $x_b$ at 900 0
	\pinlabel $x_a$ at 1120 0
	\pinlabel $x_b$ at 1230 0
	\pinlabel $t=t_1$ at 420 240
	\pinlabel $t=t_2$ at 750 240
	\pinlabel $t=T$ at 1080 240
	\pinlabel {\color{red_one}{$\rho_0(0,x)$}} at 30 200
	\endlabellist
	\centering
	\includegraphics[scale=.25]{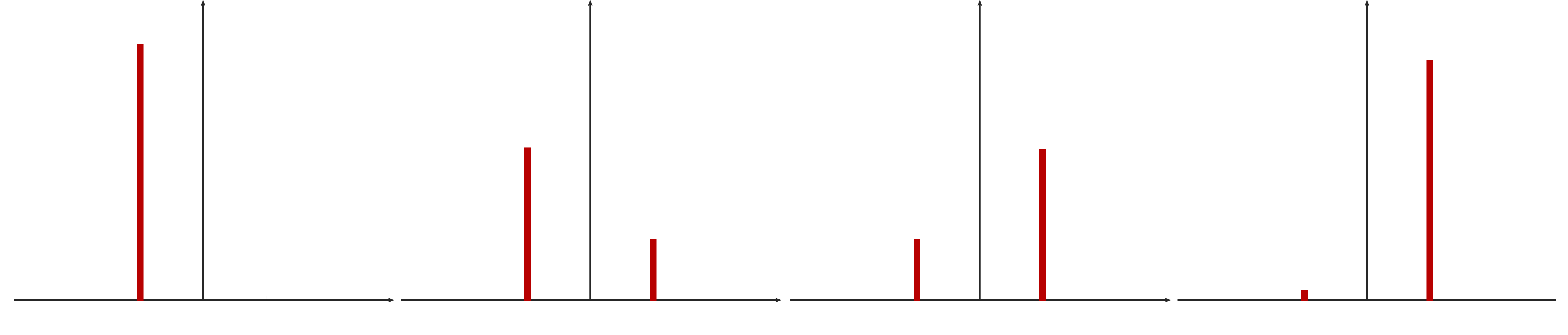}
	\caption{The time evolution of~$\rho_0$, defined as the~$\varepsilon\to 0$ limit of the solution~$\rho_\varepsilon(t,x)$ to~\eqref{GF_NGF:eq:intro:upscaled-FP}. The initial distribution is supported solely on the left. As time progresses, mass is flowing only from left to right, with rate one.}
	\label{GF_NGF:fig:intro:limit-evolution-FP}
\end{figure}
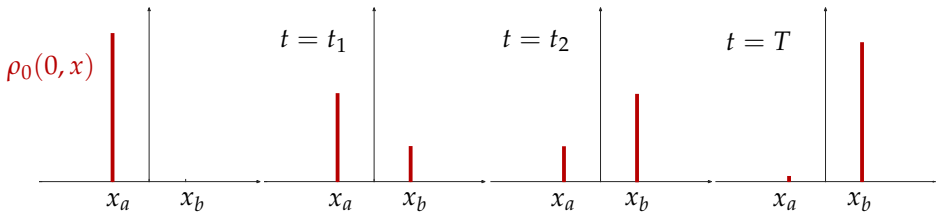

We are interested in the limit~$\varepsilon\to 0$ in the diffusion system~\eqref{GF_NGF:eq:intro:upscaled-FP}. In the limit, we expect the solution $\rho_\varepsilon$ to concentrate at the minima $x_a$ and $x_b$. This is because for small values of~$\varepsilon$, the particle spends most of its time around the minima of the potential. Furthermore, transitions from left to right occur frequently than from right to left due to the lower energy barrier. Since the transition frequency scales exponentially with~$\varepsilon$ and the potential barrier, in the limit~$\varepsilon\to 0$, we expect transitions to occur only from left to right. By our choice of the time-scale~$\tau_\varepsilon$, the limiting dynamics is characterized by mass being transfered at rate one from the local minimum~$x_a$ to the global minimum~$x_b$. In summary, $\rho_\varepsilon \to \rho_0 = z \delta_{x_a} + (1-z) \delta_{x_b}$, with a density~$z=z(t)$ decaying at rate one according to~$\partial_t z = - z$. The time evolution of the limiting density is depicted in Figure~\ref{GF_NGF:fig:intro:limit-evolution-FP}.
\subsection{From diffusion to reaction---a singular limit}
There has been recent interest in studying the limit~$\varepsilon\to 0$ for the case of \emph{symmetric} double-well potentials, that is potentials~$V$ satisfying~$V(x_a)=V(x_b)$. Peletier, Savaré and Veneroni have proved the concentration of solutions~$\rho_\varepsilon$ on the two potential-minima as~$\varepsilon$ tends to zero. The limiting densities are coupled by ODE's and correspond to a reaction-system~\cite{PeletieSavareVeneroni2010}. These results give a rigorous mathematical meaning to Kramers' program of approximating reactions by upscaling diffusions. The authors also included an additional spatial variable undergoing diffusive motion that we will not consider here.
\smallskip

A question left unanswered in~\cite{PeletieSavareVeneroni2010} rests on the fact that the Fokker-Planck equation~\eqref{GF_NGF:eq:intro:upscaled-FP} is the gradient flow of the entropy with respect to the Wasserstein metric---we give the precise definitions of gradient flows further below in Section~\ref{GF_NGF:sec:GF}. It is then natural to ask whether a convergence result such as established in~\cite{PeletieSavareVeneroni2010} can be achieved when working only with the gradient-flow structure rather than with the specific equation itself. Indeed, motivated by that question, Arnrich, Mielke, Peletier, Savaré and Veneroni soon after proved in~\cite{ArnrichMielkePeletierSavareVeneroni2012} the convergence of the corresponding Wasserstein gradient flow to a limit gradient flow. Their result comprises the convergence of the upscaled diffusion to the limiting reaction system as a special case. The proof is no longer based on the linearity of the problem, but exploits information derived solely from boundedness of the players involved in the Wasserstein gradient flow: the entropy, the Fisher information, and the Wasserstein metric. 
\smallskip

A convergence result of a variational structure, such as the Wasserstein gradient-flow, is interesting for multiple reasons. First, as Adams, Dirr, Peletier and Zimmer demonstrated~\cite{AdamsDirrPeletierZimmer2011}, the Wasserstein gradient flow is special since it arises naturally from a micro-macro limit using the theory of large deviations. That constitutes a probabilistic argument for working with the Wasserstein gradient flow rather than with other gradient-flow structures. Secondly, since many evolution equations are of Wasserstein gradient-flow type (e.g.~\cite{AmbrosioGigliSavare2008,BlanchetCalvezCarrillo2008,CarrilloDiFrancescoFigalliLaurentSlepcev2011,CarrilloYoungPilTse2019,CarlenGangbo2004,Gigli2010,GianazzaSavareToscani2009,MatthesMcCannSavare2009,Savare2007,Lisini2009}), arguments based on this variational structure have the potential to be applicable to other but similar systems as well. 
\smallskip

The abovementioned studies concentrated on \emph{symmetric} potentials. Our work presented in this chapter contributes to these studies by tackling the case of \emph{asymmetric} double-well potentials~$V$ such as shown above in Figure~\ref{GF_NGF:fig:asymmetric-doublewell-potential}. In the spirit of~\cite{PeletieSavareVeneroni2010,ArnrichMielkePeletierSavareVeneroni2012}, we establish a Gamma-convergence result for variational structures associated to the family of equations~\eqref{GF_NGF:eq:intro:upscaled-FP} in the limit~$\varepsilon\to 0$.
\smallskip

On the level of a gradient-flow structure of~\eqref{GF_NGF:eq:intro:upscaled-FP}, the \emph{asymmetry} of the potential landscape comes with a couple of challenges that we need to address. We will describe gradient flows and their related convergence concepts in more detail in Section~\ref{GF_NGF:sec:GF}. When taking the limit~$\varepsilon\to 0$, for two reasons we do in fact not expect the Wasserstein gradient flow to converge as in~\cite{ArnrichMielkePeletierSavareVeneroni2012}: first, the \emph{energies diverge in the limit}, and secondly, taking the limit means \emph{passing from reversible to irreversible}. Let us comment further on these two reasons.
\smallskip

First, various convergence concepts for gradient structures have in common that they require control of energies. In the Wasserstein gradient-flow structure of~\eqref{GF_NGF:eq:intro:upscaled-FP}, the energy is given by the relative entropy with respect to the equilibrium distribution. Due to the asymmetry of the potential, these relative entropies diverge in the limit~$\varepsilon\to 0$ (Section~\ref{GF_NGF:sec:why-GF-not-converge}). This is the main reason for which we can not follow the line of argument in~\cite{ArnrichMielkePeletierSavareVeneroni2012}, which exploits boundedness of entropies.
\smallskip

Secondly, Mielke, Peletier and Renger identified that under quite general conditions, gradient-flow structures arising from density large deviations are linked one-to-one to time-reversible stochastic processes~\cite{MielkePeletierRenger2014}. From their result, we infer a general rule of thumb: when passing from reversible stochastic processes to an irreversible stochastic process, then we do not expect the limit of the corresponding gradient-flow structures to be of gradient-flow type anymore. In our problem at hand, in the pre-limit regime the Fokker-Planck equation~\eqref{GF_NGF:eq:intro:upscaled-FP} corresponds to a reversible drift-diffusion process. When taking the limit~$\varepsilon\to 0$, we obtain a jump process with jumps only from left to right, which is an irreversible dynamics. This is why we do not expect the limit of the Wasserstein gradient flow to be a gradient flow anymore. We confirm this reasoning in our context by proving that the limiting \emph{variational structure} we obtain in Theorem~\ref{GF_NGF:thm:intro-main-result} is indeed not a gradient flow (Section~\ref{GF_NGF:sec:why-limit-is-not-GF}), even though the limiting equation can be given a gradient-flow structure.
\subsection{Flux-density functionals}
For the two abovementioned reasons, we can not take the limit of the Wasserstein gradient-flow structure. Therefore, we propose to work instead with a higher-level variational structure. While the Wasserstein gradient-flow structure can be motivated from \emph{density} large deviations, we take our motivation  from so-called \emph{flux-density} large deviations~\cite{BertiniDeSoleGabrielliJonaLasinioLandim2015}. We will introduce the rigorous terms in Section~\ref{GF_NGF:sec:flux-density-functionals}. Here, we give a brief description of the flux-density functionals and its central ingredients in order to formulate our main result. 
\smallskip

We define variational structures by functionals that act on time-dependent measures, where the minimizers of these functionals correspond to the dynamics of a Fokker-Planck equation. The Wasserstein gradient-flow of~\eqref{GF_NGF:eq:intro:upscaled-FP} is described by a map~$\mathcal{A}_\varepsilon$ acting on time-dependent probability measures~$\rho$ such that~$\mathcal{A}_\varepsilon(\rho)\geq 0$ for all~$\rho$. The solution~$\rho_\varepsilon$ to the dynamics of~\eqref{GF_NGF:eq:intro:upscaled-FP} minimizes the functional, that means~$\mathcal{A}_\varepsilon(\rho_\varepsilon)=0$. We describe this functional and its relation to gradient flows in more detail in Section~\ref{GF_NGF:sec:GF}. 
\smallskip

For defining the flux-density funtionals, we reformulate the upscaled Fokker-Planck equation~\eqref{GF_NGF:eq:intro:upscaled-FP} as an upscaled continuity equation,
\begin{equation}\label{GF_NGF:eq:upscaled-CE}
\partial_t\rho_\varepsilon + \mathrm{div} \,j_\varepsilon = 0,\quad t\geq 0,\, x\in\mathbb{R},
\end{equation}
where the function~$j_\varepsilon$ is the so-called~\emph{flux} defined by
\begin{equation}\label{GF_NGF:eq:intro:def-flux}
j_\varepsilon(t,x) := -\tau_\varepsilon\left[\varepsilon\, \nabla \rho_\varepsilon + \rho_\varepsilon\nabla V\right].
\end{equation}
In general, for a density~$\rho$, we write~$J_\varepsilon^\rho := -\tau_\varepsilon\left[\varepsilon\, \nabla \rho + \rho\nabla V\right]$. We denote the set of pairs~$(\rho,j)$ satisfying the continuity equation distributionally as
\begin{equation*}
\mathrm{CE}([0,T];\mathbb{R}) := \{(\rho,j)\,:\, \partial_t\rho + \mathrm{div}\, j = 0\;\text{in}\;\mathcal{D}'((0,T)\times\mathbb{R})\}.
\end{equation*}
The precise conditions on~$(\rho,j)$ are stated in Definition~\ref{GF_NGF:def:continuity-equation}.
\paragraph{Pre-limit functional}
For~$\varepsilon>0$, the map~$\mathcal{I}_\varepsilon:\mathrm{CE}([0,T];\mathbb{R})\to[0,\infty]$ is defined by
\begin{equation}\label{GF_NGF:eq:level-2p5-rate-function}
\mathcal{I}_\varepsilon(\rho,j) := \frac{1}{4} \int_0^T\int_\mathbb{R} \frac{1}{\varepsilon\, \tau_\varepsilon} \frac{1}{\rho(t,x)} \big|j(t,x) -J_\varepsilon^\rho(t,x)\big|^2\,dxdt.
\end{equation}
This formal expression assumes that the measure~$\rho(t,\dd x)$ is absolutely continuous with respect to the Lebesgue measure on~$\mathbb{R}$ and strictly positive. In Definition~\ref{GF_NGF:def:pre-limit-RF}, we give the mathematically rigorous expression of~$\mathcal{I}_\varepsilon$, which is a dual formulation of~\eqref{GF_NGF:eq:level-2p5-rate-function}.
\smallskip

The map~$\mathcal{I}_\varepsilon$ is a functional whose minimizer corresponds to the solution of the upscaled continuity equation~\eqref{GF_NGF:eq:upscaled-CE}, and hence the upscaled Fokker Planck equation~\eqref{GF_NGF:eq:intro:upscaled-FP}. 
The formula for~$\mathcal{I}_\varepsilon$ is motivated from large-deviation theory of flux-density pairs~\cite[Eq.~(1.3)]{BertiniDeSoleGabrielliJonaLasinioLandim2015}.
\smallskip

The flux-density functional leads by contraction to the Wasserstein gradient flow, and in that sense comprises the Wasserstein gradient-flow,
\begin{equation*}
\frac{1}{2} \mathcal{A}_\varepsilon(\rho) = \inf_{\substack{j\\(\rho,j)\in\mathrm{CE}}}\mathcal{I}_\varepsilon(\rho,j),
\end{equation*}
where the infimum is over fluxes~$j$ such that~$j\ll\rho$.
Other examples of such contraction principles from flux-density functionals to density functionals can be found for example in~\cite{Feng1994,Leonard1995,BertiniFaggionatoGabrielli2015}.
\smallskip
\subsection{Main result---$\Gamma$-convergence of flux-density functionals}
In the spirit of~$\Gamma$-convergence of functionals, we would like to obtain a limit of the functionals~$\mathcal{I}_\varepsilon$ as~$\varepsilon\to 0$. Thus the main questions that we ask in this chapter are:
\begin{enumerate}[label=(\roman*)]
	\item \emph{Compactness:} For a family of pairs~$(\rho_\varepsilon',j_\varepsilon')$ depending on~$\varepsilon$, does boundedness of~$\mathcal{I}_\varepsilon(\rho_\varepsilon',j_\varepsilon')$ imply the existence of a subsequence of~$(\rho_\varepsilon',j_\varepsilon')$ that converges in a certain topology~$\mathcal{T}$ on the set~$\mathrm{CE}([0,T];\mathbb{R})$ as~$\varepsilon\to 0$ ?
	\item \emph{Convergence along sequences:} Is there a limit functional~$\mathcal{I}_0$ satisfying
	\begin{equation*}
	(\rho_\varepsilon',j_\varepsilon')\xrightarrow{\mathcal{T}}(\rho,j)\quad\Rightarrow\quad \mathcal{I}_\varepsilon(\rho_\varepsilon',j_\varepsilon') \xrightarrow{\varepsilon\to 0} \mathcal{I}_0(\rho,j)\,?
	\end{equation*}
\end{enumerate}
We answer the first question in Theorem~\ref{GF_NGF:thm:compactness}, which establishes that sequences $(\rho_\varepsilon',j_\varepsilon')$ such that~$\mathcal{I}_\varepsilon(\rho_\varepsilon',j_\varepsilon')$ remains bounded are compact with respect to a certain topology. In Theorem~\ref{GF_NGF:thm:compactness}, we make the additional assumption that the densities~$\rho_\varepsilon'$ have uniformly bounded Radon-Nikodym derivatives with respect to a stationary measure we specify in Definition~\ref{GF_NGF:def:transformed-stationary-measure}.
\smallskip

\begin{figure}[h!]
	\labellist
	\pinlabel \Large \color{black}{$\mathcal{I}_\varepsilon$} at 1300 1100
	\pinlabel \color{dark_blue}{$\text{reversible}$} at -650 1100
	\pinlabel \color{black}{$\text{Stochastic}$} at -100 1150
	\pinlabel \color{black}{$\text{Process}$} at -100 1050
	\pinlabel \color{black}{$(\varepsilon,n)$} at 225 1050
	\pinlabel \Large \color{red_one}{$\mathcal{I}_0$} at 1300 120
	\pinlabel \color{dark_blue}{$\text{Gradient Flow}$} at 1900 1100
	\pinlabel \color{red_one}{$\text{irreversible}$} at -650 150
	\pinlabel \color{black}{$\text{Stochastic}$} at -100 200
	\pinlabel \color{black}{$\text{Process}$} at -100 100
	\pinlabel \color{black}{$(0,n)$} at 225 100
	\pinlabel \color{red_one}{$\text{Non-Gradient-Flow}$} at 2000 150
	\pinlabel \color{black}{$\text{Large deviations}$} at 750 1170
	\pinlabel \color{black}{$n\to\infty$} at 750 1050
	\pinlabel \color{black}{$\text{Large deviations}$} at 750 170
	\pinlabel \color{black}{$n\to\infty$} at 750 50 
	\pinlabel \color{black}{$\varepsilon$} at -100 760
	\pinlabel \Large \color{black}{$\downarrow$} at -100 620
	\pinlabel \color{black}{$0$} at -100 480
	\pinlabel \color{black}{$\varepsilon$} at 1400 760
	\pinlabel \Large \color{black}{$\downarrow$} at 1400 620
	\pinlabel \color{black}{$0$} at 1400 480
	\endlabellist
	\centering
	\includegraphics[scale=.1]{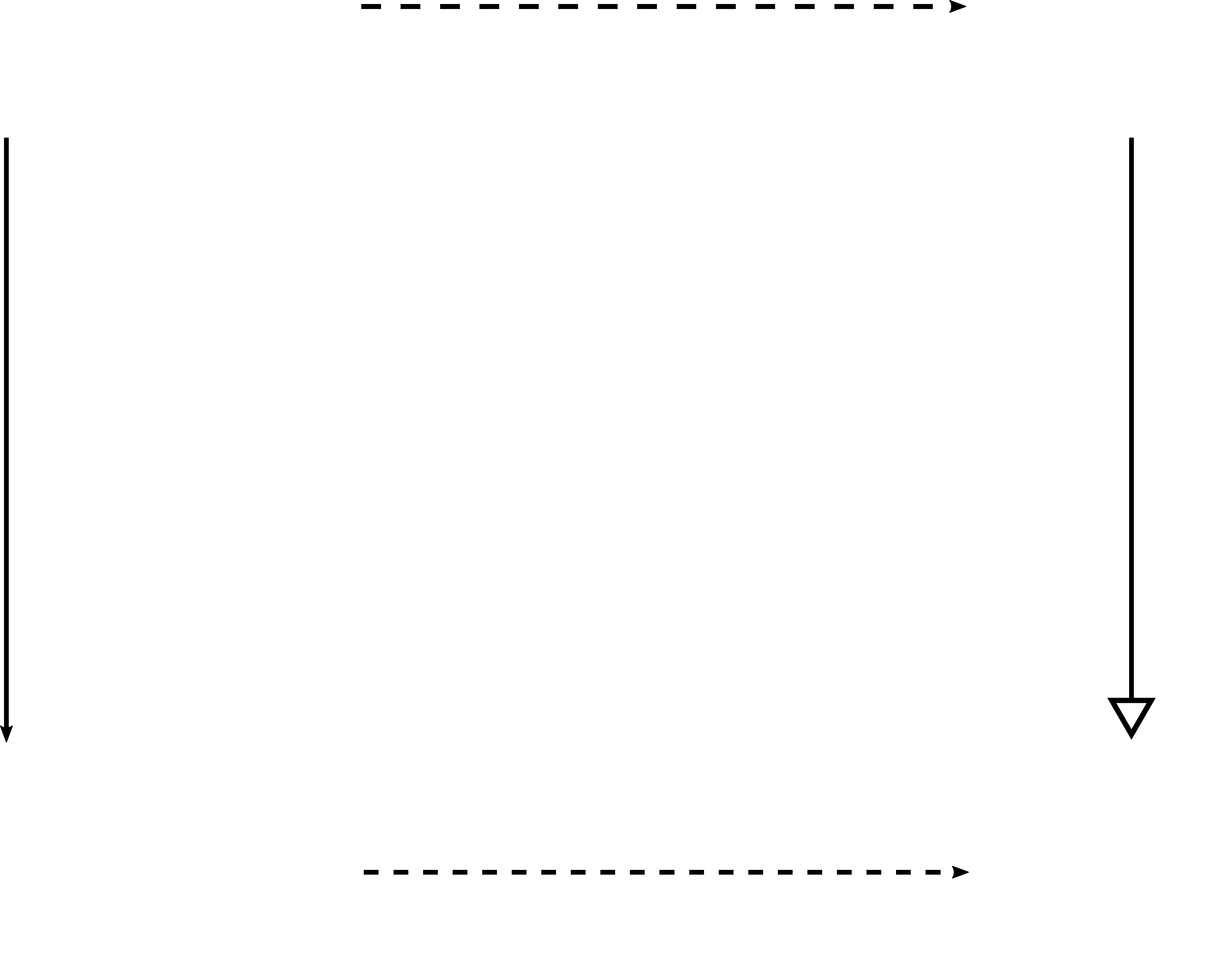}
	\caption{The top row corresponds to the empirical flux-density pairs~\eqref{GF_NGF:eq:empirical-flux-density-pairs} stemming from i.i.d.~copies of the reversible diffusion process~$X^\varepsilon_i(t)$ from~\eqref{GF_NGF:eq:intro:upscaled-diffusion-process}, whose Fokker-Planck equation is~\eqref{GF_NGF:eq:intro:upscaled-FP}. The bottom row corresponds similarly to a jump process defined on two states~$\{-,+\}$, with jumps only from~$-$ to~$+$. The bottom arrow is a rigorous large-deviation principle, and we prove the right arrow by Theorem~\ref{GF_NGF:thm:intro-main-result}. Whether the left and top arrows correspond to theorems is still an open question.}
	\label{GF_NGF:fig:commuting-diagram}
\end{figure}
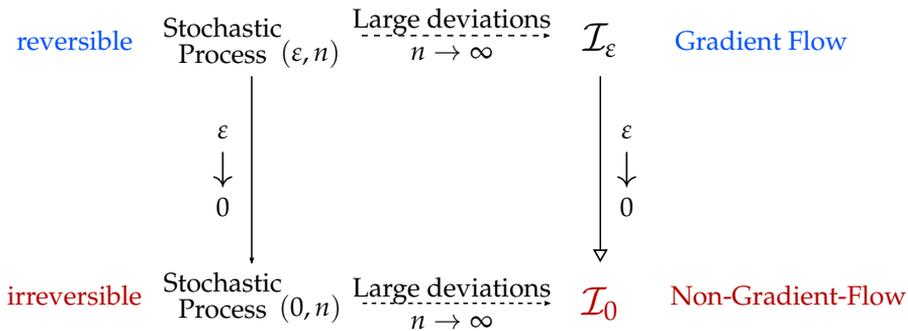

The second question is answered by Theorems~\ref{GF_NGF:thm:lower-bound} (liminf bound) and Theorem~\ref{GF_NGF:thm:upper-bound} (limsup bound), which together establish a limit of~$\mathcal{I}_\varepsilon$ in the sense of~$\Gamma$-convergence. Here, we give a short version that combines these theorems into one statement. We will consider convergence in~$\mathrm{CE}([0,T];\mathbb{R})$ in the distributional sense, meaning convergence against any smooth and compactly supported test function (Definition~\ref{GF_NGF:def:converge-in-CE}). Furthermore, we introduce a variable transformation in Definition~\ref{GF_NGF:def:coordinate-transformation} akin to our problem at hand, and we give the reason for including the transformation when defining~$y_\varepsilon$.
In brief, the purpose of this transformation is to map, in the limit~$\varepsilon\to 0$, the region around the left-minimum~$x_a$ to one point and the region around~$x_b$ to another point. The effect of this transformation for finite~$\varepsilon$ is shown much further below in Figure~\ref{GF_NGF:fig:coordinate-transformation}.
\begin{theorem}
	[Main result]
	\label{GF_NGF:thm:intro-main-result}
	There is a functional~$\mathcal{I}_0$ such that under the assumptions of Theorems~\ref{GF_NGF:thm:lower-bound} and~\ref{GF_NGF:thm:upper-bound}, we have~$\lim_{\varepsilon\to 0}\mathcal{I}_\varepsilon = \mathcal{I}_0$ in the following sense of~$\Gamma$-convergence: for any~$(\rho,j)\in \mathrm{CE}([0,T];\mathbb{R})$ such that~$\mathcal{I}_0(\rho,j)$ is finite, there are~$(\rho_\varepsilon,j_\varepsilon)\in \mathrm{CE}([0,T];\mathbb{R})$ such that
	\begin{equation*}
	(\hat{\rho}_\varepsilon,\hat{\jmath}_\varepsilon) \xrightarrow{\varepsilon\to 0} (\rho,j)\quad\text{and}\quad \mathcal{I}_\varepsilon(\rho_\varepsilon,j_\varepsilon) \xrightarrow{\varepsilon\to 0}\mathcal{I}_0(\rho,j).
	\end{equation*}
\end{theorem}
This Theorem is a first step into proving commutativity of the diagram shown in Figure~\ref{GF_NGF:fig:commuting-diagram}.
\vspace{.3cm}
In there, the stochastic process depending on~$(\varepsilon,n)$ is the so-called empirical flux-density pair~$(\rho_{\varepsilon,n},j_{\varepsilon,n})$ formally defined by
\begin{equation}\label{GF_NGF:eq:empirical-flux-density-pairs}
\rho_{\varepsilon,n} = \frac{1}{n}\sum_{i=1}^n \delta_{X_i^\varepsilon(t)}\quad\text{and}\quad j_{\varepsilon,n}\approx\frac{1}{n}\sum_{i=1}^n\delta_{X_i^\varepsilon(t)} \partial_t X_i^\varepsilon(t),
\end{equation}
where~$X_i^\varepsilon(t)$ are independent copies of the upscaled diffusion process satisfying~\eqref{GF_NGF:eq:intro:upscaled-diffusion-process}. For every fixed~$\varepsilon>0$, this process is time-reversible. In the limit~$\varepsilon\to 0$, we expect to obtain a jump process on two states~$\{-,+\}$ with jumps only from left to right. This limit process is no longer reversible.
\subsection{The limiting flux-density functional}
\label{GF_NGF:sec:intro:limiting-flux-density-functional}
We introduce the limiting functional~$\mathcal{I}_0$ from our main result, Theorem~\ref{GF_NGF:thm:intro-main-result} above, and then give a heuristic argument for why this functional is to be expected. The rate function is given in terms of the function
\begin{equation}\label{GF_NGF:eq:S-fct}
S(a,b):=
\begin{cases}
a\log(a/b) -(a-b),	& 	a,b>0,\\
b,					&	a=0,b>0,\\
+\infty,			&	\text{otherwise}.
\end{cases}
\end{equation}
\paragraph{Limit functional}
The map~$\mathcal{I}_0:\mathrm{CE}([0,T];\mathbb{R})\to[0,\infty]$ is defined by
\begin{equation*}
\mathcal{I}_0(\rho,j) := \int_0^T S(j(t)|z(t))\,dt,
\end{equation*}
whenever~$\rho(t,\dd x) = z(t)\delta_{-1/2}(\dd x) + (1-z(t))\delta_{+1/2}(\dd x)$ and the flux is piecewise constant and given by~$j(t,x) = j(t)\mathbf{1}_{(-1/2,+1/2)}(x)$. Otherwise, we set $\mathcal{I}_0(\rho,j) = +\infty$.
\smallskip

The limiting functional~$\mathcal{I}_0$ is finite only for measures~$\rho$ that are concentrated on the points~$\{\pm 1/2\}$. By continuity equation, the flux is given by $j(t)=-\partial_t z(t)$. If~$\mathcal{I}_0$ vanishes, then~$S(j(t)|z(t))=0$,  which in turn implies that~$j(t)=z(t)$. Hence the minimizer~$\rho$ of the functional~$\mathcal{I}_0$ with density~$z$ satisfies the evolution equation~$\partial_tz=-z$. The evolution of such a~$\rho$ is depicted in Figure~\ref{GF_NGF:fig:intro:limit-evolution-FP}. The fact that the limit concentrates on~$\{\pm 1/2\}$ rather than~$\{x_a,x_b\}$ is a consequence of the coordinate transformation~$y_\varepsilon$ from Definition~\ref{GF_NGF:def:coordinate-transformation}. The choice of the left-point is arbitrary, but fixes the right point. We choose~$y_\varepsilon$ such that in the limit, the distance between the points is equal to one. 
\smallskip

Just like the pre-limit functional, the functional~$\mathcal{I}_0$ is as well motivated from large-deviation theory. More precisely, it is the large-deviation rate function of flux-density pairs of independent jump processes on a set of two points~$\{-,+\}$, with jump rates~$r_{-+}=1$ and~$r_{+-}=0$.  We prove these type of large-deviation results from a Hamilton-Jacobi point-of-view in more generality in Chapter~5. Renger gives a proof based on Girsanov-transformation arguments~\cite{Renger2017}, and Kraaij provides a proof based on Hamilon-Jacobi theory~\cite{Kr17}. Heuristically, we expect the limit of the Fokker-Planck equation~\eqref{GF_NGF:eq:intro:upscaled-FP} to be characterized by exactly that dynamics: transition of mass occurs from left to right at rate one, while transitions from right to left do not occur at all.
\subsection{Overview}
The rest of this chapter is structured as follows. In Section~\ref{GF_NGF:sec:GF}, we introduce gradient-flow structures and their relation to large-deviation theory. This part provides the background to convergence to equilibirum in the Fokker-Planck equation~\eqref{GF_NGF:eq:intro:upscaled-FP} emphasising the role of entropy, and the Wasserstein gradient flow formulation. In Section~\ref{GF_NGF:sec:flux-density-functionals}, we start from the action formulation (Definition~\ref{GF_NGF:def:wassterstein-GF}) to demonstrate in Section~\ref{GF_NGF:sec:why-GF-not-converge} that the Wasserstein functional does not converge for our problem at hand, as opposed to~\cite{ArnrichMielkePeletierSavareVeneroni2012}. After that, we define the flux-density functionals~$\mathcal{I}_\varepsilon$ and~$\mathcal{I}_0$. In Section~\ref{GF_NGF:sec:proof-of-gamma-convergence}, we prove the main Theorem by splitting it in three statements: compactness (Theorem~\ref{GF_NGF:thm:compactness}), the lower bound (Theorem~\ref{GF_NGF:thm:lower-bound}) and the upper bound (Theorem~\ref{GF_NGF:thm:upper-bound}). For the proof of the lower bound, we work under the assumption of uniformly bounded densities.
\section{Gradient flows}
\label{GF_NGF:sec:GF}
Gradient flows are an example of variational structures that appear naturally in partial differential equations modelling dissipative phenomena. We refer to Peletier's lectures on variational modelling~\cite{Peletier2014} for background and physical motivations on gradient flows . Here in this section, we do not cover novel results, but provide the background to gradient-flow structures that underlie the type of Fokker-Planck equations we consider in this chapter. 
\smallskip

By Fokker-Planck equations, we generally refer to a class of partial differential equations that describe the time evolution of distributions of stochastic processes. We illustrate gradient flows without the parameters present in~\eqref{GF_NGF:eq:intro:upscaled-FP}, that is we consider~$\rho=\rho(t,x)$ solving an equation of the type
\begin{equation}\label{GF_NGF:eq:Fokker-Planck-intro}
\partial_t \rho = \Delta \rho + \mathrm{div}\left(\rho \nabla V\right),\quad t\geq 0, \, x\in\mathbb{R},
\end{equation}
which is a basic model for the probability distribution of a particle diffusing in one dimension in a confining potential landscape~$V(x)$. The asymmetric double-well potential as depicted in Figure~\ref{GF_NGF:fig:asymmetric-doublewell-potential} is an example of a confining potential: the particle is trapped by the potential, which effectively forces the particle to diffuse closely around the minima of~$V$.
\smallskip

We think of the solution~$\rho(t,x)$ to~\eqref{GF_NGF:eq:Fokker-Planck-intro} as the probability of observing the particle at time~$t$ being in state~$x$. As illustrated above by Kramers problem, the particle's state~$x$ can represent the value of a system's parameter that is fluctuating in time according to Brownian motion, and does not necessarily correspond to the position of a physical particle diffusing in a viscous fluid. The function~$V(x)$ then usually corresponds to an effective energy landscape. For a background on Fokker-Planck equations and their derivations from a phenomenological point of view, we refer to Risken's monograph on this type of equations~\cite{Risken1996}.
\smallskip

We first discuss in Section~\ref{GF_NGF:sec:GF:conv-to-equil} convergence to equilibrium in~\eqref{GF_NGF:eq:Fokker-Planck-intro}. The main point is to illustrate at the same time both the phenomenon of entropy-dissipation and the dynamics of~\eqref{GF_NGF:eq:Fokker-Planck-intro}. 
Then in Section~\ref{GF_NGF:sec:GF:diss-mech}, we recall the so-called JKO-scheme that Jordan, Kinderlehrer and Otto identified in~\cite{JordanKinderlehrerOtto1998}, to illuminate the fact that~\eqref{GF_NGF:eq:Fokker-Planck-intro} is the a solution to the gradient-flow of the entropy with respect to the Wasserstein distance. 
In Section~\ref{GF_NGF:sec:GF:var-struct}, we choose the formulation of gradient flows that we will use in later sections, by means of certain functionals~$\mathcal{I}$. The minimizers of those functionals are called \emph{curves of maximal slope} that correspond to gradient solutions in the classical case. This way of formulating a gradient flow in metric spaces goes back to Ennio De Giorgi and collaborators~\cite{DeGiorgiMarinoTosques1980}. The Wasserstein gradient-flow of~\eqref{GF_NGF:eq:Fokker-Planck-intro} is one example of such a structure: the gradient-flow dynamics is characterized as the minimizer of a functional that involves the entropy, Fisher information and Wasserstein distance. This variational formulation via a functional is the starting point for the~$\Gamma$-convergence results in~\cite{ArnrichMielkePeletierSavareVeneroni2012}, which is why we recall it in Section~\ref{GF_NGF:sec:GF:var-struct}. We use this formulation to showcase where exactly the line of argument in~\cite{ArnrichMielkePeletierSavareVeneroni2012} is limited to symmetric potentials~$V$.
\smallskip

The presentation draws from different sources: the overview of Markowich and Villani on convergence to equilibrium~\cite{MarkowichVillani2000}, the monograph on gradient flows~\cite{AmbrosioGigliSavare2008} written by Ambrosio, Gigli and Savaré, and the program of deriving gradient flows from large deviations put forward by Adams, Dirr, Mielke, Peletier, Renger and Zimmer~\cite{AdamsDirrPeletierZimmer2011,MielkePeletierRenger2014}.
\subsection{Convergence to equilibrium via dissipation of entropy}
\label{GF_NGF:sec:GF:conv-to-equil}
Boltzmann discovered the celebrated~$\mathrm{H}$-Theorem: according to Boltzmann's equation, an ideal gas of particles evolves in time in such a way that its so-called entropy is monotonically decreasing. As a consequence, after enough time has passed, we find the gas in a state minimizing the entropy. In this state, the distribution of particle's velocities~$v$ is stable, known as the Maxwell-Boltzmann distribution. More generally, we usually refer to a state minimizing the entropy as equilibrium. If Boltzmann's~$\mathrm{H}$-Theorem applies to a system of particles, an initial distribution of particles will eventually converge to equilibrium.
\smallskip

We can also observe such a convergence phenomenon for the solutions to the Fokker-Planck equation~\eqref{GF_NGF:eq:Fokker-Planck-intro}.
Its equilibrium state, defined by satisfying~$\partial_t\gamma=0$, is given by~$\gamma = e^{-V}$, and we will refer to it as the \emph{Boltzmann distribution} or simply \emph{equilibrium}. We will assume that~$\gamma$ has mass one (otherwise, we add a constant to the potential~$V$). For measuring how far a solution~$\rho$ of~\eqref{GF_NGF:eq:Fokker-Planck-intro} is away from equilibirum, it is natural to introduce the density~$u(t,x)$ by~$\rho(t,\dd x)=u(t,x)\gamma(\dd x)$. Then the solution~$\rho$ is in equilibrium if and only if~$u\equiv 1$. The density~$u$ evolves in time according to
\begin{equation*}
\partial_t u = \Delta u - \nabla u \nabla V,
\end{equation*}
which we can infer from~\eqref{GF_NGF:eq:Fokker-Planck-intro}. For fixed time~$t>0$, the \emph{relative entropy} of~$\rho$ with respect to equilibrium is defined as 
\begin{equation}\label{GF_NGF:eq:intro:def-entropy}
\mathrm{Ent}(\rho(t,\cdot)|\gamma) := \int_\mathbb{R} u(t,x)\log u(t,x) \,\dd \gamma(x).
\end{equation}
In notation, we shall suppress the dependence on time.
If the equilibrium distribution~$\gamma$ is clear from the context, we call~$\mathrm{Ent}$ simply the~\emph{entropy}. Here, we will point out the special role that the entropy plays in the study of the convergence to equilibrium. For further details, we refer to Markowich's and Villani's overview~\cite{MarkowichVillani2000}, where the authors connect convergence to equilibrium with various functional inequalities. 
\smallskip

The entropy vanishes if~$u\equiv 1$, and therefore vanishes if~$\rho$ is in equilibrium. By the estimate~$x\log x \geq x-1$, the entropy is non-negative:
\begin{equation*}
\mathrm{Ent}(\rho|\gamma)\myeqdef \int_\mathbb{R}\frac{\dd\rho}{\dd\gamma} \log\left(\frac{\dd\rho}{\dd\gamma}\right)\,\dd\gamma \geq \int_\mathbb{R}\left(\frac{\dd\rho}{\dd\gamma}-1\right)\,\dd\gamma = 0.
\end{equation*}
Hence the equilibrium distribution~$\gamma$ indeed minimizes the entropy. For a solution~$\rho(t,\dd x)$ of the Fokker-Planck equation~\eqref{GF_NGF:eq:Fokker-Planck-intro}, let us see how the entropy evolves in time. A calculation involving integration by parts yields
\begin{equation}\label{GF_NGF:eq:intro:diss-ent-is-fisher}
\frac{\dd}{\dd t}\mathrm{Ent}(\rho|\gamma) = - \mathrm{I}(\rho|\gamma),\quad \text{where}\quad\mathrm{I}(\rho|\gamma):=\int_\mathbb{R} |\nabla (\log u)|^2\,\dd\rho.
\end{equation}
The functional~$\mathrm{I}$ is non-negative and zero only if~$u$ is constant. Hence the entropy decreases in time unless~$\rho$ is in equilibrium. The functional~$\mathrm{I}$ is known as the Fisher information. In Chapter~4, we encounter the Fisher information as the exponential convergence rate of the empirical measure associated to~\eqref{GF_NGF:eq:Fokker-Planck-intro}.
\smallskip

Under suitable assumptions on the potential~$V$, we can be more precise about how fast the entropy decays. We say that the distribution~$\gamma=e^{-V}$ satisfies the logarithmic Sobolev inequality with a constant~$\lambda>0$ if
\begin{equation}\label{GF_NGF:eq:intro:log-sobolev}
\mathrm{Ent}(\rho|\gamma) \leq \frac{1}{2\lambda} \mathrm{I}(\rho|\gamma).
\end{equation}
If that inequality is satisfied, then
\begin{equation*}
\frac{\dd}{\dd t} \mathrm{Ent}(\rho|\gamma) \overset{\eqref{GF_NGF:eq:intro:diss-ent-is-fisher}}{=} -\mathrm{I}(\rho|\gamma) \overset{\eqref{GF_NGF:eq:intro:log-sobolev}}{\leq} -2\lambda \mathrm{Ent}(\rho|\gamma).
\end{equation*}
Hence by Gr\"{o}nwall's inequality, the entropy decays exponentially fast:
\begin{equation}\label{GF_NGF:eq:intro:ent-decays-exp}
\mathrm{Ent}(\rho(t,\cdot)|\gamma) \leq \mathrm{Ent}(\rho(0,\cdot)|\gamma) e^{-2\lambda t}. 
\end{equation}
According to the estimate~\eqref{GF_NGF:eq:intro:ent-decays-exp}, the entropy is being dissipated under the time evolution of~\eqref{GF_NGF:eq:Fokker-Planck-intro} under two conditions. First, the initial distribution~$\rho(0,\cdot)$ must be non-singular with respect to equilibrium in the sense that the relative entropy is finite. Secondly, the logarithmic Sobolev inequality must be satisfied. The latter is satisfied for confining potentials (for instance~\cite[Theorem~1]{MarkowichVillani2000} and the discussion thereafter). Otto and Villani give geometric derivations of Talagrand- and logarithmic Sobolev inequalities in~\cite{OttoVillani2000}.
\smallskip

The above analysis demonstrates that for solutions~$\rho$ of the Fokker-Planck equation, the dissipation of entropy happens exponentially fast under fairly general conditions on the potential~$V$, and that the amount of dissipation is quantified by the Fisher information~\eqref{GF_NGF:eq:intro:diss-ent-is-fisher}. Jordan, Kinderlehrer and Otto revealed in~\cite{JordanKinderlehrerOtto1998} an exciting and deeper geometric connection between the entropy and the Fokker-Planck equation: the solution~$\rho$ flows in the direction of the Wasserstein gradient of the entropy. Their analysis leads to a variational structure that we will call a \emph{Wasserstein gradient flow}, and we shall discuss their insights next.
%

\subsection{Gradient flow---a dissipation mechanism}
\label{GF_NGF:sec:GF:diss-mech} 
As we saw above, solutions to the Fokker-Planck equation~\eqref{GF_NGF:eq:Fokker-Planck-intro} evolve such that entropy decays exponentially fast. Here, we shall discuss the dissipation mechanism that  Jordan, Kinderlehrer and Otto identified in~\cite{JordanKinderlehrerOtto1998}, in which the entropy plays the role of the energy being dissipated. Before we describe this dissipation mechanism for the Fokker-Planck equation~\eqref{GF_NGF:eq:Fokker-Planck-intro}, we illustrate the central ingredients of a dissipation mechanism in a simpler context.
\smallskip

In the one-dimensional Euclidian setting, a gradient-flow is an equation of the type
\begin{equation}\label{GF_NGF:eq:intro:basic-GF-in-R}
\partial_t x = -\nabla\mathrm{E}(x),\qquad x(0)=x_0,
\end{equation}
where~$x:[0,T]\to\mathbb{R}$ is a sufficiently regular path,~$\mathrm{E}:\mathbb{R}\to\mathbb{R}$ is a confining potential that we refer to as an \emph{energy} and~$\nabla \mathrm{E}$ is the gradient of~$\mathrm{E}$, here the derivative. 
By confining we mean that~$E(x)\to\infty$ as~$|x|\to\infty$,~$\mathrm{E}\in C^2(\mathbb{R})$ and that its second derivative is uniformly bounded from below. 
We will write~$x(t)=x_t$ for the evaluation of the path~$x$ at time~$t$.
\smallskip

By definition of the gradient in 1d, in each time step the solution~$x$ to~\eqref{GF_NGF:eq:intro:basic-GF-in-R} follows the direction that dissipates as much energy as possible. As a result, as time tends to infinity, the solution converges to a local minimum of~$\mathrm{E}$. One way to understand this evolution is to start from a time-discretization.
For an infinitesimal time-step~$\tau>0$, the backward Euler approximation to~\eqref{GF_NGF:eq:intro:basic-GF-in-R} 
is
\begin{equation*}
\frac{x(t+\tau)-x(t)}{\tau} + \nabla \mathrm{E}(x(t+\tau)) \approx 0,
\end{equation*}
which motivates the implicit Euler scheme: define the set of points~$\{x_k^\tau\}_{k=0,1,\dots}$ iteratively by~$x_0^\tau := x_0$ and
\begin{equation}\label{GF_NGF:eq:intro:min-movement-1d-R}
x_k^\tau := \mathrm{argmin}_{\substack{x}\in\mathbb{R}}\left[\frac{d(x,x_{k-1}^\tau)^2}{2\tau} + \mathrm{E}(x)\right],\qquad d(x,y) := |x-y|.
\end{equation}
The map~$d$ is just the standard Euclidian metric. We have the following convergence statement of this time-discretization: if both~$\tau\to 0$ and~$k\to\infty$ such that~$k\tau \to t$, then $x_k^\tau \to x(t)$, where~$x$ solves~\eqref{GF_NGF:eq:intro:basic-GF-in-R}.
\smallskip

In the formulation~\eqref{GF_NGF:eq:intro:min-movement-1d-R}, we can recognize a couple of aspects. First, in each time-step the solution~$x$ minimizes not merely the potential, but rather the combination of both the metric~$\mathrm{d}$ and the potential~$\mathrm{E}$. 
Secondly, we can also interpret the precise role of the metric. To that end, consider the step~$k\to k+1$ for a fixed and small value of~$\tau$. In the minimization procedure~\eqref{GF_NGF:eq:intro:min-movement-1d-R}, points far away from the starting point~$x_{k}^\tau$ are heavily punished since the metric is upscaled by~$\tau^{-1}$, whereas points close to~$x_k^\tau$ 
that decrease the value of~$\mathrm{E}$ are favored. In combination, roughly speaking, the faster the metric grows (the map $y\mapsto d(y,x_k^\tau)$), the less the energy will decrease in the step~$k\to k+1$. In this way, the metric determines how much energy is dissipated in each time step. 
\smallskip

We therefore call~\eqref{GF_NGF:eq:intro:min-movement-1d-R} a \emph{dissipation meachanism} underlying the gradient flow equation~\eqref{GF_NGF:eq:intro:basic-GF-in-R}. In geometric terms, the potential and metric together determine the direction of movement while the metric controls the amount of dissipated energy per step. Let us summarize the players of the dissipation mechanism that leads to the gradient flow~\eqref{GF_NGF:eq:intro:basic-GF-in-R}:
\begin{enumerate}[label=(\roman*)]
	\item A state space~$\mathrm{M}$; here~$\mathrm{M}=\mathbb{R}$.
	\item A map~$\mathrm{E}:\mathrm{M}\to\mathbb{R}$; here~$\mathrm{E}$ is a confining potential, which we call energy.
	\item A metric~$d:\mathrm{M}\times \mathrm{M}\to[0,\infty]$; here~$d$ is the standard Euclidian metric.
\end{enumerate}
Jordan, Kinderlehrer and Otto made the remarkable discovery~\cite{JordanKinderlehrerOtto1998} that the Fokker-Planck equation~\eqref{GF_NGF:eq:Fokker-Planck-intro} admits a dissipation mechanism in which the Boltzmann entropy serves as the energy. They identified the corresponding distance as a transport cost that arises in the theory of optimal transport. The scheme they developed is made from the following three ingredients: 
\begin{enumerate}[label=(\roman*)]
	\item $\mathrm{M}:=\mathcal{P}_2(\mathbb{R})$, the set of probability measures with finite second moments.
	\item $\mathrm{E}:=\mathrm{Ent}(\cdot|\gamma):\mathrm{M}\to\mathbb{R}\cup\{+\infty\}$, the entropy defined as in~\eqref{GF_NGF:eq:intro:def-entropy} by
	\begin{equation}\label{GF_NGF:eq:JKO-scheme-relative-entropy}
	\mathrm{Ent}(\mu|\gamma) := \int_\mathbb{R} u\log u \,\dd \gamma,\quad \text{with}\; u(x) := \frac{\dd\mu}{\dd\gamma}(x).
	\end{equation}
	If~$\mu$ is not absolutely continuous with respect to~$\gamma$, then~$\mathrm{Ent}(\mu|\gamma):=+\infty$.
	\item $d:=\mathcal{W}$, the Wasserstein metric defined by
	\begin{equation*}
	\mathcal{W}(\mu,\nu) := \inf_{m\in \Pi(\mu,\nu)}\int_{\mathbb{R}\times\mathbb{R}}|x-y|^2\,m(\dd x\dd y),
	\end{equation*}
	where~$\Pi(\mu,\nu)$ is the set of probability measures on~$\mathbb{R}\times\mathbb{R}$ whose first marginal equals~$\mu$ and whose second marginal equals~$\nu$.
\end{enumerate}
\smallskip

The Wasserstein metric can be interpreted as the minimal cost required to transport a pile of sand distributed as~$\mu$ to a pile of sand distributed as~$\nu$, where the cost of transporting a sand grain from~$x$ to~$y$ is given by~$|x-y|^2$. For a thorough historical and mathematical overview of the topic of optimal transport we refer to Villani's monograph~\cite{Villani2008}.
\smallskip

The main result discovered by Jordan, Kinderlehrer and Otto is the following dissipation mechanism (\cite[Theorem~5.1]{JordanKinderlehrerOtto1998}). For an initial condition~$\rho_0\in \mathcal{P}_2(\mathbb{R})$ and fixed time-step~$\tau>0$, define~$\{\rho_k^\tau\}_{k=0,1,\dots}$ iteratively by
\begin{equation}\label{GF_NGF:eq:intro:diss-mech-for-FP}
\rho_k^\tau := \text{argmin}_{\mu\in \mathcal{P}_2(\mathbb{R})}\left[\frac{\mathcal{W}(\mu,\rho_{k-1}^\tau)^2}{2\tau}+\mathrm{Ent}(\mu|\gamma)\right].
\end{equation}
Define the piecewise-constant path~$\rho^\tau$ by~$\rho^\tau(t):=\rho_k^\tau$ for~$t\in[k\tau,(k+1)\tau)$. Then we have~$\rho^\tau\to \rho$ strongly in~$L^1((0,T)\times\mathbb{R})$ as~$\tau\to 0$, where~$\rho$ is the solution to~\eqref{GF_NGF:eq:Fokker-Planck-intro} with initial condition~$\rho(0,\dd x)=\rho_0(\dd x)$.
\smallskip

This time-discretization scheme, also refered to as the JKO-scheme, is one way of making sense of the one-dimensional gradient flow~\eqref{GF_NGF:eq:intro:basic-GF-in-R} in infinite dimensions (with $\mathrm{M}=\mathcal{P}_2(\mathbb{R})$ instead of~$\mathrm{M}=\mathbb{R}$). On top of the fact that the entropy decays exponentially fast, the JKO-scheme reveals that solutions to the Fokker-Planck equation~\eqref{GF_NGF:eq:Fokker-Planck-intro} flow along the \emph{steepest descent} of the entropy. The Wasserstein metric determines the amount of dissipated entropy per time-stepm just as the Euclidian metric determines the amount of dissipated energy in~\eqref{GF_NGF:eq:intro:min-movement-1d-R}. We say the solution~$\rho$ to the scheme~\eqref{GF_NGF:eq:intro:diss-mech-for-FP} is the solution to the \emph{Wasserstein gradient-flow}. 
Soon after, Otto further attached a precise geometrical meaning to a Wasserstein gradient-flow~\cite{Otto2001}.
\smallskip

Ambrosio, Gigli and Savaré generalize the formulation of gradient flows via a dissipation mechanism, such as~\eqref{GF_NGF:eq:intro:diss-mech-for-FP}, to arbitrary metric spaces~\cite[Chapter~2]{AmbrosioGigliSavare2008}. This generalization of~\eqref{GF_NGF:eq:intro:diss-mech-for-FP} is called a \emph{minimizing movement scheme}, \cite[Definition~2.0.6]{AmbrosioGigliSavare2008}. The main assumptions on the energy functional in order to obtain solutions to a minimizing movement scheme are suitable coercivity, lower-semicontinuity and compactness prtoperties. For the precise set of assumptions, we refer in particular to~\cite[Section~2.2]{AmbrosioGigliSavare2008}. 
%
%
\subsection{Gradient flow---a variational structure}
\label{GF_NGF:sec:GF:var-struct}
In the previous section, we discussed the minimizing movement scheme or JKO-scheme~\eqref{GF_NGF:eq:intro:diss-mech-for-FP}, a dissipation mechanism build up from a triple~$(\mathrm{M},\mathrm{E},d)$. The JKO-scheme represents one way of regarding the solution~$\rho$ to~\eqref{GF_NGF:eq:Fokker-Planck-intro} as a solution to a \emph{gradient flow}, since the limiting solution obtained from~\eqref{GF_NGF:eq:intro:diss-mech-for-FP} flows along the steepest descent of~$\mathrm{E}$. 
\smallskip

While this formulation of a gradient flow in terms of discrete time steps is conceptually enlightening, it is not well suited for passing to limits in gradient flows. For instance, given a family of triples~$(\mathrm{M},\mathrm{E}_\varepsilon,d_\varepsilon)$, under which convergence conditions on the energies~$\mathrm{E}_\varepsilon$ and metrics~$d_\varepsilon$ will solutions to the~$(\mathrm{M},\mathrm{E}_\varepsilon,d_\varepsilon)$-scheme converge to solutions of a limiting scheme $(\mathrm{M},\mathrm{E}_0,d_0)$? To answer questions of that type for gradient flows and to simplify the treatment of convergence of gradient flows, we introduce in this section a different but formally equivalent formulation of a gradient-flow. This formulation is known as the \emph{energy-dissipation principle}, and defines a gradient flow in terms of a functional.
As in the previous section, we first illustrate the formulation on the example of the real-valued gradient flow~\eqref{GF_NGF:eq:intro:basic-GF-in-R}. Then we turn to the Wasserstein gradient flow of the Fokker-Planck equation~\eqref{GF_NGF:eq:Fokker-Planck-intro}.
\subsubsection{Example illustrating the energy-dissipation principle}
Recall that~\eqref{GF_NGF:eq:intro:basic-GF-in-R} is the equation~$\partial_t x = -\nabla\mathrm{E}(x)$ on~$\mathbb{R}$. 
For any~$a,b\in\mathbb{R}$, if~$2ab\leq -a^2-b^2$, then $a=-b$. Hence a path~$x$ solves~\eqref{GF_NGF:eq:intro:basic-GF-in-R} if and only if
\begin{equation}\label{GF_NGF:eq:GF:Young-bound-GF-on-R}
\partial_tx \cdot \nabla E(x) \leq -\frac{1}{2}|\partial_t x|^2 - \frac{1}{2}|\nabla E(x)|^2.
\end{equation}
By the chain rule,~$\partial_t\mathrm{E}(x)=\nabla\mathrm{E}(x)\partial_t x$. Hence performing integration in time, the solution~$x$ to~\eqref{GF_NGF:eq:intro:basic-GF-in-R} satisfies the inequality
\begin{equation*}
E(x_T) + D(x;0,T) \leq \mathrm{E}(x_0),\quad 
\end{equation*}
where we introduced the \emph{dissipation}
\begin{equation}\label{GF_NGF:eq:GF:dissipation-euclidian-example}
D(x;0,T):=\int_0^T\left[\frac{1}{2} |\partial_t x|^2+ \frac{1}{2}|\nabla E(x)|^2\right]\dd t.
\end{equation}
For an absolutely continuous path~$y:[0,T]\to\mathbb{R}$,
\begin{equation*}
\mathcal{I}(y) := \mathrm{E}(y_T)-\mathrm{E}(y_0) + D(y;0,T).
\end{equation*}
The map~$y\mapsto\mathcal{I}(y)$ carries two important features. First, it is non-negative for \emph{any} path~$y$. This follows from the chain rule and the estimate~$2ab \geq -a^2-b^2$,
\begin{equation*}
E(y_T)-E(y_0) = \int_0^T\partial_t y\cdot \nabla E(y)\,\dd t \geq \int_0^T\left[-\frac{1}{2}|\partial_t y|^2 - \frac{1}{2}|\nabla E(y)|^2\right]\dd t.
\end{equation*}
Secondly,~$\mathcal{I}(y)$ vanishes if and only if~$y=x$ is a solution to the gradient flow~$\partial_tx=-\nabla\mathrm{E}(x)$, which follows from the bound~\eqref{GF_NGF:eq:GF:Young-bound-GF-on-R} for~$y$. 
\smallskip

Since~$\mathcal{I}$ is non-negative and zero only for the solution, we can reformulate the solution of a gradient-flow as
\begin{equation*}
\partial_tx=-\nabla\mathrm{E}(x) \quad\Leftrightarrow\quad \mathrm{E}(x_T)+D(x)\leq E(x_0)\quad\Leftrightarrow\quad \mathcal{I}(x) = 0.
\end{equation*}
The map~$\mathcal{I}$ is an example of a gradient-flow structure. The reformulation in terms of a dissipation~$\mathrm{D}(\cdot)$ by integrating the infinitesimal gradient-flow in time is an example of Ennio De Giorgi's \emph{Energy-Dissipation-Principle}. This connection between solutions to gradient flows is also refered to as the \emph{Energy-Dissipation Theorem}, which holds in much more generality than presented here. 
%
\subsubsection{Energy-Dissipation Principle}
Motivated by the previous example, we will consider the following formulation of a gradient flow, which corresponds to~\cite[Definition~1.1]{ArnrichMielkePeletierSavareVeneroni2012}.
\begin{definition}[Gradient flow]\label{GF_NGF:def:GF}
Let~$\mathrm{M}$ be a metric space,~$\mathrm{E}:\mathrm{M}\to\mathbb{R}\cup\{+\infty\}$ be a function and~$\mathrm{D}(\cdot;t_1,t_2):C([0,T];\mathrm{M})\to\mathbb{R}\cup\{+\infty\}$ be a functional defined for all~$0\leq t_1<t_2\leq T$. We call the triple~$(\mathrm{M},\mathrm{E},\mathrm{D})$ a \emph{gradient-flow structure} if for any~$\mu\in C([0,T];\mathrm{M})$ and all~$t_1<t_2$, the inequality
\begin{equation}\label{GF_NGF:eq:def:GF}
\mathrm{E}(\mu_{t_2}) + D(\mu;t_1,t_2) \geq E(\mu_{t_1})
\end{equation}
is satisfied.\qed
\end{definition}
We call~$\mathrm{E}$ the \emph{energy} and~$\mathrm{D}$ the \emph{dissipation}. 
\begin{definition}[Solution to gradient flow]\label{GF_NGF:def:sol-to-GF}
We call a curve~$\mu \in C([0,T];\mathrm{M})$ a solution to the gradient-flow structure~$(\mathrm{M},\mathrm{E},\mathrm{D})$ if~$\mathrm{E}(\mu_0)<\infty$ and
\begin{equation}\label{GF_NGF:eq:def:GF-solution}
\mathrm{E}(\mu_t) + \mathrm{D}(\mu;0,t) = \mathrm{E}(\mu_0)\qquad \text{for all } t\in [0,T].
\end{equation}
\qed
\end{definition} 
A solution to the gradient flow is defined by maximizing the dissipation of energy; hence the equality~\eqref{GF_NGF:eq:def:GF-solution}, which in the classical case corresponds by~\eqref{GF_NGF:eq:GF:Young-bound-GF-on-R} to curves of maximal slope. This equality is the \emph{Energy-Dissipation Principle}.
\smallskip

There are various ways in which the dissipation~$\mathrm{D}$ may depend on the energy~$\mathrm{E}$. In the definition of a gradient flow, the inequality~\eqref{GF_NGF:eq:def:GF} plays the role of replacing the chain rule. In the Euclidian example from above, this inequality holds true as a consequence of two aspects: the dissipation is related to the energy by~\eqref{GF_NGF:eq:GF:dissipation-euclidian-example} and the chain rule applies in~$\mathbb{R}$. On the other hand, the gradient flow~$\partial_tx=-\nabla\mathrm{E}(x)$ is recovered from~\eqref{GF_NGF:eq:def:GF-solution} via the fact that~$2ab\leq -a^2-b^2$ implies $a=-b$. In general, the formulation of a gradient-flow solution via Definition~\ref{GF_NGF:def:sol-to-GF} is equivalent if the dissipation is given via so-called dissipation potentials~$\mathcal{R}$-$\mathcal{R}^\ast$ (~\cite[Theorem~3.3.1]{Mielke2016}).
\smallskip

The Wasserstein gradient flow of~\eqref{GF_NGF:eq:Fokker-Planck-intro} however can be well-motivated from the quadratic structure as in~\eqref{GF_NGF:eq:GF:dissipation-euclidian-example}. We give an example based on generalizing the modulus of the gradient~$|\nabla\mathrm{E}|$ and the velocity~$|\partial_tx|$. To that end, let~$(\mathrm{M},\mathrm{E},\mathrm{D})$ be a triple as in Definition~\ref{GF_NGF:def:GF}. The local slope of the functional~$\mathrm{E}:\mathrm{M}\to\mathbb{R}\cup\{+\infty\}$ is defined by (\cite[Definition~1.2.4]{AmbrosioGigliSavare2008})
\begin{equation*}
|\partial\mathrm{E}|(\mu) := \limsup_{\nu\to\mu}\frac{(\mathrm{E}(\mu)-\mathrm{E}(\nu))_+}{d(\mu,\nu)}.
\end{equation*}
For an absolutely-continuous curve $\mu:[0,T]\to\mathrm{M}$, define its metric velocity as (\cite[Eq.~(1.1.3)]{AmbrosioGigliSavare2008})
\begin{equation*}
|\partial_t\mu|(t) := \lim_{\Delta t\to 0}\frac{d(\mu(t),\mu(t+\Delta t))}{\Delta t}.
\end{equation*}
We consider the dissipation to be given by
\begin{equation}\label{GF_NGF:eq:intro:def-of-dissipation}
\mathrm{D}(\mu,t_1,t_2) := \int_{t_1}^{t_2}\left[\frac{1}{2}|\partial_t\mu|(t)^2 + \frac{1}{2}|\partial\mathrm{E}|(\mu_t)^2\right]\,\dd t.
\end{equation}
Assume that~$|\partial\mathrm{E}|$ is a \emph{strong upper gradient} (\cite[Definition~1.2.1]{AmbrosioGigliSavare2008}), meaning
\begin{equation*}
|\mathrm{E}(\gamma_{t_2})-\mathrm{E}(\gamma_{t_1})| \leq \int_{t_1}^{t_2} |\partial_t \gamma|(t)\cdot |\partial\mathrm{E}|(\gamma_t)\,\dd t
\end{equation*}
holds for every absolutely-continuous curve~$\gamma:[0,T]\to\mathrm{M}$. Then by Young's inequality,
\begin{equation*}
\mathrm{E}(\gamma_{t_1})-\mathrm{E}(\gamma_{t_2}) \leq \int_{t_1}^{t_2} \left[\frac{1}{2}|\partial_t \gamma|(t)^2 + \frac{1}{2}|\partial\mathrm{E}|(\gamma_t)^2\right]\,\dd t 
\overset{\eqref{GF_NGF:eq:intro:def-of-dissipation}}{=} \mathrm{D}(\gamma,t_1,t_2).
\end{equation*}
Hence if the local slope of the energy functional~$\mathrm{E}$ is a strong upper gradient, then the triple~$(\mathrm{M},\mathrm{E},\mathrm{D})$ is indeed a gradient-flow structure in the sense of Definition~\ref{GF_NGF:def:GF}.
\subsubsection{Wasserstein gradient flow}
As the JKO-scheme suggests, we choose the state space~$\mathrm{M}=\mathcal{P}_2(\mathbb{R})$ equiped with the Wasserstein metric and $\mathrm{E}(\mu):=\mathrm{Ent}(\mu|\gamma)$, the relative entropy as defined in~\eqref{GF_NGF:eq:JKO-scheme-relative-entropy} with the equilbirum distribution~$\gamma=e^{-V}$. We take the dissipation~$\mathrm{D}$ defined by~\eqref{GF_NGF:eq:intro:def-of-dissipation}. To complete the description, we give the characterization of the relative entropy's local slope and the Wasserstein velocity~$|\partial_t\mu|(t)$.
\smallskip

The local slope of the relative entropy is the relative Fisher information \cite[Theorem~10.4.7]{AmbrosioGigliSavare2008} given by
\begin{equation*}
\mathrm{I}(\mu|\gamma) = \int_\mathbb{R}|\nabla \log u|^2\,\dd\mu,\qquad u(x) := \frac{\dd\mu}{\dd\gamma}(x).
\end{equation*}
The Wasserstein velocity is characterized by a particular velocity field~$v(t,x)$ satisfying~$\partial_t\mu_t+\mathrm{div}(v_t\mu_t)=0$ in the sense of distributions, as (\cite[Proposition~8.4.5]{AmbrosioGigliSavare2008})
\begin{equation*}
|\partial_t \mu|(t)^2 = \int_\mathbb{R}|v(t,x)|^2\,\mu(t,\dd x).
\end{equation*}
This characterization is closely related to the dynamical formulation of the Wasserstein distance discovered by Benamou and Brenier~\cite{BenamouBrenier2000},
\begin{equation*}
\mathcal{W}(\mu_0,\mu_1)^2 = \inf_{\mu_t}\left\{\int_0^1\int_\mathbb{R}|v(t,x)|^2\,\mu_t(\dd x)dt\,:\,\partial_t\mu_t+\mathrm{div}(v_t\mu_t)=0\right\}\,
\end{equation*}
where~$v(t,x)=v_t(x)$.
With these remarks, we find the dissipation~$\mathrm{D}$
\begin{equation}\label{GF_NGF:eq:wasserstein-dissipation}
\mathrm{D}(\mu;t_{1},t_{2}) = \int_{t_1}^{t_2}\left[\frac{1}{2}\int_\mathbb{R}|v(t,x)|^2\mu(t,\dd x) + \frac{1}{2}\mathrm{I}(\mu_t|\gamma)\right]\,\dd t.
\end{equation}
\begin{definition}[Wasserstein gradient flow]\label{GF_NGF:def:wassterstein-GF}
	Let~$\mathrm{M}:=\mathcal{P}_2(\mathbb{R})$, the set of probability measures on~$\mathbb{R}$ with finite second moments. Let~$\mathcal{A}:C([0,T];\mathrm{M})\to[0,\infty]$ be the map given by
	\begin{equation}\label{GF_NGF:non-eps-Wasserstein-action}
	\mathcal{A}(\rho) := \mathrm{Ent}(\rho_T|\gamma)-\mathrm{Ent}(\rho_0|\gamma) + \mathrm{D}(\rho;0,T),
	\end{equation}
	with the relative entropy~\eqref{GF_NGF:eq:JKO-scheme-relative-entropy} and the dissipation~\eqref{GF_NGF:eq:wasserstein-dissipation}. We call~$\mathcal{A}$ the \emph{Wasserstein action functional} corresponding to the \emph{Wasserstein gradient-flow structure} given by~$(\mathcal{P}_2(\mathbb{R}),\mathrm{Ent},\mathcal{W})$. The curve~$\rho$ satisfying~$\mathcal{A}(\rho)=0$ is called the solution to the Wasserstein gradient flow.\qed
\end{definition}
The minimizer of~$\mathcal{A}$ is also the solution to the JKO-scheme~\cite[Theorem~11.2.1]{AmbrosioGigliSavare2008}, and hence the solution to the Fokker-Planck equation~\eqref{GF_NGF:eq:Fokker-Planck-intro}. The solution~$\rho$ satisfies the energy-dissipation equality (e.g.~\cite[Eq.~(11.2.4)]{AmbrosioGigliSavare2008})
\begin{equation*}
\mathrm{Ent}(\rho_T|\gamma)+\int_0^T\mathrm{I}(\rho_t|\gamma)\,\dd t=\mathrm{Ent}(\rho_0|\gamma),
\end{equation*}
which is also obtained by integrating~\eqref{GF_NGF:eq:intro:diss-ent-is-fisher} in time. 
\subsubsection{The $\mathcal{R}$-$\mathcal{R}^\ast$ formulation}
The quadratic structure of the dissipation~\eqref{GF_NGF:eq:intro:def-of-dissipation} is an example of a more general structure of grandient flows. We introduce this more general structure here. In the~$\mathcal{R}$-$\mathcal{R}^\ast$ formulation, generalized gradient-flow structures in a smooth setting arise from a combination of the following three ingredients:
\begin{enumerate}[label=(\roman*)]
	\item A \emph{state space}~$\mathrm{M}$, which is a set with a sufficiently rich differentiable structure that gives gradients a meaning, such as Riemannian manifolds.
	\item A function~$\mathrm{E}:\mathrm{M}\to\mathbb{R}$, the \emph{energy}.
	\item A function~$\mathcal{R}:T\mathrm{M}\to[0,\infty]$, which we call \emph{dissipation potential}.,
	such that for each state~$x\in \mathrm{M}$:
	\begin{enumerate}[label=($\mathcal{R}$\arabic*)]
		\item $\mathcal{R}(x,\cdot):T_x\mathrm{M}\to[0,\infty]$ is convex and lower semicontinuous.
		\item $\mathcal{R}(x,0) = \min_{v\in T_x\mathrm{M}}\mathcal{R}(x,v) = 0$.
	\end{enumerate}
\end{enumerate}
We denote by~$\mathcal{R}^\ast:T^\ast \mathrm{M}\to[0,\infty]$ the Legendre transform of~$\mathcal{R}$ defined by
\begin{equation*}
\mathcal{R}^\ast(x,\xi) := \sup_{v\in T_q\mathrm{M}}\left[\ip{\xi}{v} - \mathcal{R}(x,v)\right],
\end{equation*}
For an energy~$\mathrm{E}:\mathrm{M}\to\mathbb{R}$, we denote by~$\mathrm{d}\mathrm{E}$ its differential, that is the map
\begin{equation*}
\mathrm{d}\mathrm{E}:\mathrm{M}\to T^\ast \mathrm{M}, \quad x\mapsto \mathrm{d}\mathrm{E}_x \in T_x^\ast \mathrm{M}.
\end{equation*}
A functional~$\mathcal{I}$ acting on trajectories~$C([0,T];\mathrm{M})$ is of \emph{gradient-flow structure} if
\begin{equation}\label{GF_NGF:eq:action-of-R-R_star-functional}
\mathcal{I}(x) = \mathrm{E}(x_T)-\mathrm{E}(x_0) + \int_0^T \left[\mathcal{R}\left(x_t,\partial_t{x}_t\right)+\mathcal{R}^\ast\left(x_t,-\mathrm{d}\mathrm{E}(x_t)\right)\right]dt.
\end{equation}
The dissipation~\eqref{GF_NGF:eq:GF:dissipation-euclidian-example} from above corresponds to the flat space~$\mathbb{R}$, where~$\dd\mathrm{E}$ gets identified with~$\nabla\mathrm{E}$, and with the quadratic dissipation potentials
\begin{equation*}
\mathcal{R}(x,v) = \frac{1}{2}v^2\quad\text{and}\quad \mathcal{R}^\ast(x,\xi) = \frac{1}{2}\xi^2.
\end{equation*}
Typical examples where the~$\mathcal{R}$-$\mathcal{R}^\ast$ dissipation occurs are large-deviation rate functions of jump processes. The limit problem of this chapter is an example: if the potential~$V$ in~\eqref{GF_NGF:eq:intro:upscaled-FP} is symmetric, then our limit variational structure admits an~$\mathcal{R}$-$\mathcal{R}^\ast$ formulation with a $\cosh$-type dissipation potential~$\mathcal{R}$; however, we show there is no such formulation for \emph{asymmetric} potentials (Section~\ref{GF_NGF:sec:why-limit-is-not-GF}). In infinite-dimensional settings such as in our pre-limit problem, a careful definition of the tangent and cotangent spaces~\cite[Section~12.4]{AmbrosioGigliSavare2008} is required to make the above display~\eqref{GF_NGF:eq:action-of-R-R_star-functional} rigorous.
\subsection{Taking limits of gradient-flow structures}
There are many variants of taking limits of gradient flows. Mielke provides several different definitions of types of convergences in~\cite{Mielke2016}. The natural concept of a gradient-flow convergence is to demand both the energies and dissipations to converge separately. Sandier and Serfaty introduced this concept in~\cite{SandierSerfaty2004}, which since then has found applications in a variety of other problems. This convergence concept is also applied in~\cite{ArnrichMielkePeletierSavareVeneroni2012}, and we introduce it next. For~$\rho_n, \rho\in C([0,T],\mathrm{M})$, we say~$\rho_n\to\rho$ if the convergence is uniform in time.
\begin{definition}[EDP convergence]\label{GF_NGF:def:EDP-convergence}
	Let~$(\mathrm{M},\mathrm{E}_\varepsilon,\mathrm{D}_\varepsilon)$ be a family of gradient-flow structures in the sense of Definition~\ref{GF_NGF:def:GF}. We say that~$(\mathrm{M},\mathrm{E}_\varepsilon,\mathrm{D}_\varepsilon)$ converges in the \emph{EDP sense} to a gradient-flow structure~$(\mathrm{M},\mathrm{E}_0,\mathrm{D}_0)$ if: 
	\begin{enumerate}[label=(\roman*)]
		\item $\mathrm{E}_\varepsilon\xrightarrow{\Gamma}\mathrm{E}_0$ in~$\mathrm{M}$;
		\item For each~$t\in[0,T]$, $\mathrm{D}_\varepsilon\xrightarrow{\Gamma}\mathrm{D}_0$ in~$C([0,T];\mathrm{M})$.\qed
	\end{enumerate}
\end{definition}
EDP convergence implies convergence of solutions (e.g.~\cite[Lemma~2.8]{MielkeMontefuscoPeletier2020}).
\begin{proposition}\label{GF_NGF:prop:simple-EDP-convergence}
Assume that a family of gradient-flow structures~$(\mathrm{M},\mathrm{E}_\varepsilon,\mathrm{D}_\varepsilon)$ converges in the EDP sense to a gradient-flow structure~$(\mathrm{M},\mathrm{E}_0,\mathrm{D}_0)$. Let~$\rho_\varepsilon$ be the solutions of~$(\mathrm{M},\mathrm{E}_\varepsilon,\mathrm{D}_\varepsilon)$. Suppose that
\begin{equation*}
\rho_\varepsilon\to\rho_0 \;\text{in}\; C([0,T];\mathrm{M})\quad\text{and}\quad \mathrm{E}_\varepsilon(\rho_\varepsilon(0)) \to \mathrm{E}_0(\rho_0(0)).
\end{equation*}
Then~$\rho_0$ is a solution of~$(\mathrm{M},\mathrm{E}_0,\mathrm{D}_0)$.
\end{proposition}
\begin{proof}
	By the assumption of~$\Gamma$-convergences, exploiting~$\rho_\varepsilon\to \rho$ gives
	\begin{multline*}
	\mathrm{E}_0(\rho_0(T))-\mathrm{E}_0(\rho_0(0)) + \mathrm{D}_0(\rho_0;0,T) \\\leq \liminf_{\varepsilon\to 0}\left[ \mathrm{E}_\varepsilon(\rho_\varepsilon(T))-\mathrm{E}_\varepsilon(\rho_\varepsilon(0)) + \mathrm{D}_\varepsilon(\rho_\varepsilon;0,T)\right] = 0.
	\end{multline*}
	The other inequality is satisfied by Definition~\ref{GF_NGF:def:GF}, and hence~$\rho_0$ is a solution in the sense of Definition~\ref{GF_NGF:def:sol-to-GF}.
\end{proof}
EDP convergence is therefore a suitable limit concept in problems where the energies remain bounded in the limit~$\varepsilon\to 0$. In~\cite{ArnrichMielkePeletierSavareVeneroni2012}, the authors proof boundedness of the entropies for the case of a symmetric potential~$V$. As a consequence of their~$\Gamma$-convergence result, one obtains by Proposition~\ref{GF_NGF:prop:simple-EDP-convergence} convergence of solutions for free. For the asymmetric case we consider in this chapter, boundedness of the entropies is no longer satisfied (Section~\ref{GF_NGF:sec:why-GF-not-converge}).
\section{Flux-density functionals}
\label{GF_NGF:sec:flux-density-functionals}
The main point of this section is to give the rigorous definitions of the flux-density functionals (Definitions~\ref{GF_NGF:def:pre-limit-RF} and~\ref{GF_NGF:def:limit-RF}). We first demonstrate why the Wasserstein gradient-flow functionals do not converge as in~\cite{ArnrichMielkePeletierSavareVeneroni2012}, which motivates our choice of working with flux-density functionals in the first place. 
\subsection{Why the gradient-flow structure does not converge}
\label{GF_NGF:sec:why-GF-not-converge}
By Definition~\ref{GF_NGF:def:wassterstein-GF}, the action functional defining the Wasserstein gradient-flow structure of the Fokker-Planck equation~\eqref{GF_NGF:eq:intro:upscaled-FP} is
\begin{equation}\label{GF_NGF:epsilon-Wasserstein-GF-action}
\mathcal{A}_\varepsilon(\rho) = \mathrm{Ent}(\rho(T)|\gamma_\varepsilon) -\mathrm{Ent}(\rho(0)|\gamma_\varepsilon) + \mathrm{D}_\varepsilon(\rho;0,T),
\end{equation}
where the time-scale parameter~$\tau_\varepsilon$ enters the dissipation~$\mathrm{D}_\varepsilon$ as in~\cite{ArnrichMielkePeletierSavareVeneroni2012}; the dissipation part is not important for our argument however. The equilibrium distribution~$\gamma_\varepsilon$ has mass one and is given by
\begin{equation*}
\gamma_\varepsilon(\dd x) = \mathcal{N}_\varepsilon^{-1}e^{-V(x)/\varepsilon}\,\dd x,\quad \mathcal{N}_\varepsilon := \int_\mathbb{R} e^{-V(x)/\varepsilon}\,\dd x.
\end{equation*}
In the limit~$\varepsilon\to 0$, the equilibrium distribution concentrates solely on the global minimum~$x_b$ of the potential~$V$, that is~$\gamma_\varepsilon$ converges weakly to $\delta_{x_b}$.
For proving EDP convergence (Definition~\ref{GF_NGF:def:EDP-convergence}), we have to verify that the energies and dissipations converge independently from one another.
\begin{proposition}[Diverging entropies]\label{GF_NGF:prop:diverging-entropies}
	Let~$\gamma_\varepsilon$ be the equilibrium distribution to the Fokker-Planck equation~\eqref{GF_NGF:eq:intro:upscaled-FP} and let~$\mu_\varepsilon\in\mathcal{P}(\mathbb{R})$ be a family of probability measures converging weakly to~$\mu_0$. Suppose that~$\mu_0\neq \delta_{x_b}$. Then
	\begin{equation*}
	\mathrm{Ent}(\mu_\varepsilon|\gamma_\varepsilon) \to +\infty.
	\end{equation*}
\end{proposition}
\begin{proof}
	This follows from the fact that~$\gamma_\varepsilon\xrightharpoonup{\ast} \delta_{x_b}$ and~$\mathrm{Ent}(\mu_0|\delta_{x_b})=+\infty$.
\end{proof}
Let us demonstrate why this excludes the possibility of proving EDP convergence of the Wasserstein gradient flow of~\eqref{GF_NGF:eq:intro:upscaled-FP}.
A EDP convergence requires the entropies of the initial conditions to converge to some \emph{finite} limiting energy functional~$\mathrm{E}_0$,
\begin{equation*}
\mathrm{Ent}(\rho_\varepsilon(0)|\gamma_\varepsilon) \to \mathrm{E}_0(\rho_0(0)).
\end{equation*}
By Proposition~\ref{GF_NGF:prop:diverging-entropies}, any initial condition that is not concentrating on~$\{x_b\}$ leads to diverging relative entropies. Therefore, if we insist on finite entropies in the limit, the only initial conditions we could consider are those converging to equilibrium in the limit~$\varepsilon\to 0$, meaning only the initial condition~$\delta_{x_b}$. Since this excludes any dynamics in the limit~$\varepsilon\to 0$, we do not work with~$\mathcal{A}_\varepsilon$.
\smallskip

While we can not use the Wasserstein action~\eqref{GF_NGF:epsilon-Wasserstein-GF-action}, the density-flux functional~$\mathcal{I}_\varepsilon$ from~\eqref{GF_NGF:eq:level-2p5-rate-function} provides a natural way to cope with this divergence of entropies by including them into the dissipation. We sketch this observation. To that end, consider the density-flux functional without the $\varepsilon$-dependent parameters,
\begin{equation*}
	\mathcal{I}(\rho,j) =\frac{1}{4}\int\frac{1}{\rho}|j-j^\rho|^2\,\dd x \dd t,
\end{equation*}
where~$\partial_t\rho +\mathrm{div}\,j=0$ and~$J^\rho=-\nabla\rho + \rho\nabla V$. Expanding the square leads to
\begin{equation*}
	\mathcal{I}(\rho,j) = \frac{1}{4}\left[\int \left(\frac{\dd j}{\dd\rho}\right)^2\,\dd\rho +\int \left(\frac{\dd J^\rho}{\dd\rho}\right)^2\,\dd\rho - 2 \int \frac{\dd j}{\dd\rho}\frac{\dd J^\rho}{\dd\rho}\,\dd\rho\right].
\end{equation*}
Using~$J^\rho = -\rho \nabla \log (\rho/e^{-V})=-\rho\nabla (\log u)$, and hence~$\dd J^\rho/\dd\rho = -\nabla\log(u)$, we find the second term to be the Fisher information in the dissipation~$\mathrm{D}$ from~\eqref{GF_NGF:eq:wasserstein-dissipation}. For the cross term, integrating by parts and using~$\partial_t\rho + \mathrm{div}\,j=0$,
\begin{align*}
	-\int \frac{\dd j}{\dd\rho}\frac{\dd J^\rho}{\dd\rho}\,\dd\rho &= \int(\nabla\log u)\,\dd j\\
	&= - \int \log (u) \,\mathrm{div}\,j \,\dd x\dd t\\
	&= \int \log(u) \,\partial_t\rho = \int \partial_t u \log(u)\,\dd\gamma \dd t \\
	&= \int_0^T \partial_t \left(\int_\mathbb{R} u\log u\,\dd x\right) \,\dd t + 0,
\end{align*}
which leads to the entropy terms. Hence taking the infimum over fluxes~$j$ satisfying the continuity equation and such that~$j\ll \rho$, and using the Benamou-Brenier characterization of the Wasserstein distance, we find back the Wasserstein functional~\eqref{GF_NGF:non-eps-Wasserstein-action},
\begin{equation*}
	\inf_j \mathcal{I}(\rho,j) = \frac{1}{2}\mathcal{A}(\rho).
\end{equation*}
\subsection{Stationary measure and coordinate transformation}
\label{GF_NGF:sec:cooridnate-transformation}
We introduce in Definition~\ref{GF_NGF:def:left-normalized-stationary-measure} stationary measures that are not normalized to mass one on~$\mathbb{R}$. Rather, they are normalized to mass one when integrating from~$-\infty$ to the local maximum~$x_0$ of~$V$. We call these measures \emph{left-normalized stationary measures}. The normalization is chosen in order to capture the transitions from left to right, which are or order one. We use the superscript~$\ell$ to distinguish objects that are derived from this choice of normalization.
\begin{definition}[Left-normalized stationary measure]\label{GF_NGF:def:left-normalized-stationary-measure}
	For a potential $V$ as in Figure~\ref{GF_NGF:fig:asymmetric-doublewell-potential}, define the left-normalized stationary measure $\gamma_\varepsilon^\ell \in \mathcal{M}(\mathbb{R})$ by
	\begin{equation*}
	\gamma_\varepsilon^\ell(A) := Z_\varepsilon^{-1} \int_A e^{-V(x)/\varepsilon}\,dx, \quad Z_\varepsilon^{-1} := \int_{-\infty}^{x_0} e^{-V(x)/\varepsilon}\, dx.
	\end{equation*}
\end{definition}
With this left-normalization, these measures concentrate in the limit~$\varepsilon\to 0$ on the potential wells, that is the set $\{V\leq 0\}$.
\begin{proposition}[Concentration of measure] \label{prop:GF_NGF:left-normalized-stationary-measure}
	For any~$\delta > 0$,
	\begin{equation*}
	\lim_{\varepsilon\to 0}\gamma_\varepsilon^\ell(\{V>\delta\}) = 0.
	\end{equation*}
\end{proposition}
\begin{proof}[Proof of Proposition~\ref{prop:GF_NGF:left-normalized-stationary-measure}]
	Let~$\delta > 0$. For any finite $M > \delta$, we find
	\begin{align*}
	\gamma_\varepsilon^\ell\left(\left\{M > V \geq \delta\right\}\right) &= Z_\varepsilon^{-1} \int_{\left\{M > V \geq \delta\right\}} e^{-V(x)/\varepsilon}\,dx \\
	&\leq Z_\varepsilon^{-1}\, e^{-\delta/\varepsilon}\, \mathcal{L}\left(\left\{M > V \geq \delta\right\}\right) \xrightarrow{\varepsilon\to 0}0,
	\end{align*}
	since $V(x_a) = 0$ and $Z_\varepsilon=[1+o(1)] \sqrt{2\pi\varepsilon/V''(x_a)} e^{-V(x_a)/\varepsilon}$, which follows from Laplace's method (Lemma~\ref{lemma:watson}).
	\smallskip
	
	Since $V(x) \to \infty$ as $|x| \to \infty$, there is a sequence $M_n \to \infty$ such that the sets $A_n := \left\{V \geq M_n\right\}$ are decreasing in $n$ in the sense that $A_{n+1}\subseteq A_n$. By finiteness of the measure~$\gamma_\varepsilon^\ell$, this implies
	\begin{equation*}
	\lim_{n \to\infty} \gamma_\varepsilon^\ell(A_n) = \gamma_\varepsilon^\ell\left(\cap_{n} A_n \right) = 0.
	\end{equation*}
	Let $\kappa > 0$. Since
	\begin{equation*}
	\gamma_\varepsilon^\ell\left(\left\{V \geq \delta\right\}\right) = \gamma_\varepsilon^\ell\left(\left\{V\geq M_n\right\}\right) + \gamma_\varepsilon^\ell\left(\left\{M_n > V \geq \delta\right\}\right),
	\end{equation*}
	when choosing $n$ large enough such that $\gamma_\varepsilon^\ell(A_n) < \kappa/2$, then for all $\varepsilon$ small enough such that $\gamma_\varepsilon^\ell\left(\left\{M_n > V \geq \delta\right\}\right) < \kappa/2$, we obtain the estimate
	\[
	\gamma_\varepsilon^\ell\left(\left\{V \geq \delta\right\}\right) < \kappa.
	\]
	Since~$\kappa$ is arbitrary, the claim follows.
\end{proof}
With these stationary measures at hand, we now motivate the coordinate transformation~$y_\varepsilon$. To that end, we start from the flux-density rate function (\cite[Eq.~(1.3)]{BertiniDeSoleGabrielliJonaLasinioLandim2015}) specialized to the Fokker-Planck equation~\eqref{GF_NGF:eq:intro:upscaled-FP},
\begin{equation*}
\mathcal{I}_\varepsilon(\rho,j) := \frac{1}{4} \int_0^T\int_\mathbb{R} \frac{1}{\varepsilon \tau_\varepsilon} \frac{1}{\rho(t,x)} \big|j(t,x) -j_\varepsilon^\rho(t,x)\big|^2\,\dd x\dd t,
\end{equation*}
where $j_\varepsilon^\rho(t,x) := -\tau_\varepsilon\left[\varepsilon \nabla \rho + \rho \nabla V\right]$. This formula suggests that we should consider measures~$\rho(t,\dd x)$ that have a Lebesgue density. Then~$\rho$ also has a density with respect to the left-normalized stationary measure, and we write $\rho(t,\dd x) = u_\varepsilon^\ell(t,x)\gamma_\varepsilon^\ell(\dd x) = u_\varepsilon^\ell(t,x)g_\varepsilon^\ell(x)\dd x$ with $g_\varepsilon^\ell(x)=Z_\varepsilon^{-1}e^{-V(x)/\varepsilon}$. With that, the flux-density rate function can be written as
\begin{equation*}
\mathcal{I}_\varepsilon(\rho,j) = \frac{1}{4} \int_0^T\int_\mathbb{R} \frac{1}{\varepsilon\tau_\varepsilon} \frac{1}{g_\varepsilon^\ell(x) u_\varepsilon^\ell(t,x)}\big| j(t,x) + \varepsilon\,\tau_\varepsilon\, g_\varepsilon^\ell(x) \partial_x u_\varepsilon^\ell(t,x)\big|^2\, \dd x \dd t.
\end{equation*}
The transformation is chosen such that the parameters~$\varepsilon$ and~$\tau_\varepsilon$ are shifted to the densities and do not appear explicitly in the flux-density functional. This suggests to introduce the coordinate transformation~$y_\varepsilon$ on~$\mathbb{R}$ satisfying
\begin{equation*}
\dd y_\varepsilon(x) = \frac{1}{\varepsilon\, \tau_\varepsilon\, g_\varepsilon^\ell(x)} \,\dd x.
\end{equation*}
Then with $\hat{u}_\varepsilon^\ell(t,y_\varepsilon(x)) := u_\varepsilon^\ell(t,x)$ and $\hat{\jmath}_\varepsilon(t,y_\varepsilon(x)) := j(t,x)$, we obtain
\begin{equation*}
\mathcal{I}_\varepsilon(\rho,j) = \frac{1}{4}\int_0^T\int_\mathbb{R} \frac{1}{\hat{u}_\varepsilon^\ell(t,y)}\big|\hat{\jmath}_\varepsilon(t,y) + \partial_y \hat{u}_\varepsilon^\ell(t,y)\big|^2\,\dd y\dd t.
\end{equation*}
Written in this form, all the parameters are absorbed into the density~$\hat{u}_\varepsilon^\ell$. The coordinate transformation~$y_\varepsilon$ is the almost the same as in~\cite{ArnrichMielkePeletierSavareVeneroni2012}; the only difference is that we use the left-normalized stationary measure, whereas in the symmetric case, one can use the stationary measure normalized to one. 
\begin{definition}[Coordinate transformation $y_\varepsilon$]\label{GF_NGF:def:coordinate-transformation}
	For a potential~$V$ as in Figure~\ref{GF_NGF:fig:asymmetric-doublewell-potential}, the left-normalized stationary measure $\dd\gamma_\varepsilon^\ell = g_\varepsilon^\ell(z)dz$ of Definition~\ref{GF_NGF:def:left-normalized-stationary-measure} and the time-scale $\tau_\varepsilon$ defined by~\eqref{GF_NGF:intro:eq:def-time-scale-parameter}, define the map $y_\varepsilon : \mathbb{R} \to \mathbb{R}$ by
	\begin{align*}
	y_\varepsilon(x) := \frac{1}{\varepsilon\, \tau_\varepsilon} \int_{x_0}^x \frac{1}{g_\varepsilon^\ell(z)}\,\dd z.
	\tag*\qed
	\end{align*}
\end{definition}
\begin{figure}[h!]
	\labellist
	\pinlabel $V(x)$ at 1700 1250
	\pinlabel $x_0$ at 800 600
	\pinlabel $x_a$ at 400 780
	\pinlabel $x_b$ at 1350 750
	\pinlabel $x$ at 1600 750
	\pinlabel $-\frac{1}{2}$ at 400 -60
	\pinlabel $0$ at 800 -60
	\pinlabel $+\frac{1}{2}$ at 1150 -60
	\pinlabel $y_\varepsilon(x)$ at 1750 -50
	\endlabellist
	\centering
	\includegraphics[scale=.1]{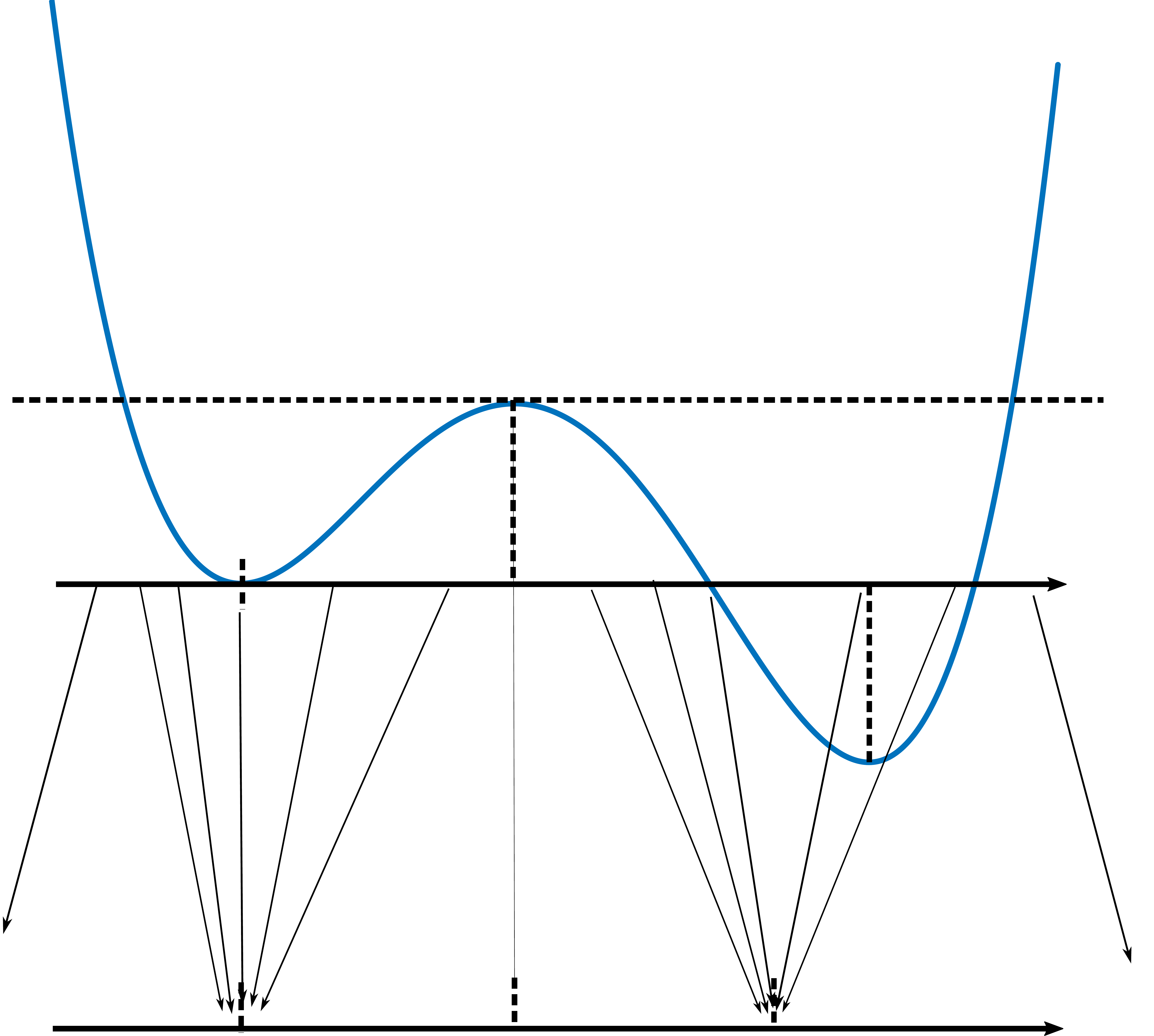}
	\caption{The effect of the coordinate transformation $y_\varepsilon$ of Defintion~\ref{GF_NGF:def:coordinate-transformation}. Points to the left of~$x_0$ s.t.~$V(x)<V(x_0)$ are mapped to~$-1/2$, and similarly, points to the right of~$x_0$ are mapped to~$+1/2$. The smaller the value of~$\varepsilon$, the sharper is the concentration effect. As~$\varepsilon\to 0$, points far to the left of~$x_a$ and far to the right of~$x_b$ are mapped to~$\mp\infty$, respectively.}
	\label{GF_NGF:fig:coordinate-transformation}
\end{figure}
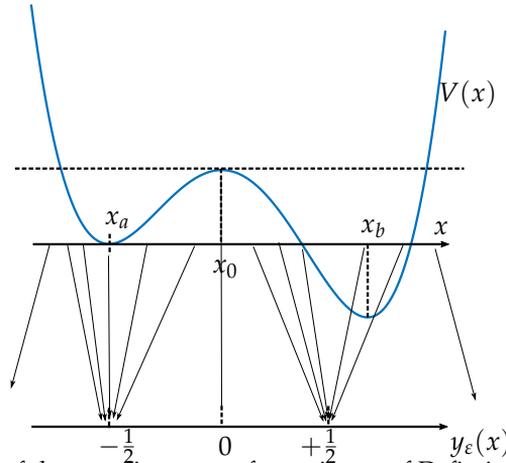
This coordinate transformation indeed maps the minima of the potential~$V$ to~$\pm 1$, as the following Proposition shows.
\begin{proposition}[Coordinate transformation]\label{GF_NGF:prop:coordinate-transformation}
	The map $y_\varepsilon:\mathbb{R}\to \mathbb{R}$ of Definition~\ref{GF_NGF:def:coordinate-transformation} satisfies the following:
	\begin{enumerate}[label=(\roman*)]
		\item The map $y_\varepsilon$ is strictly increasing and bijective.
		\item For any $x<x_0$ such that $V(x) < V(x_0)$, we have $y_\varepsilon(x) \to -\frac{1}{2}$ as $\varepsilon \to 0$.
		\item For any $x>x_0$ such that $V(x) < V(x_0)$, we have $y_\varepsilon(x) \to +\frac{1}{2}$ as $\varepsilon \to 0$.
	\end{enumerate}
\end{proposition}
\begin{proof}[Proof of Proposition~\ref{GF_NGF:prop:coordinate-transformation}]
	Since $y_\varepsilon'(x) > 0$ for any $x \in \mathbb{R}$ and $y_\varepsilon(x) \to \pm \infty$ as $x \to \pm\infty$, the map $y_\varepsilon$ is strictly increasing and bijective. For $x < x_0$ satisfying $V(x) < V(x_0)$, we obtain
	\begin{align*}
	y_\varepsilon(x) &= \frac{1}{\varepsilon \tau_\varepsilon} \cdot Z_\varepsilon \cdot  \int_{x_0}^x e^{V(z)/\varepsilon}\, dz	\\
	&=[1+o(1)] \frac{1}{\varepsilon \tau_\varepsilon} \cdot e^{-V(x_a)/\varepsilon} \sqrt{\frac{2\pi\varepsilon}{V''(x_a)}} \cdot \frac{1}{2}e^{V(x_0)/\varepsilon} \sqrt{\frac{2\pi\varepsilon}{|V''(x_0)|}} (-1) = -\frac{1}{2},
	\end{align*}
	by applying Lemma~\ref{lemma:watson} to $Z_\varepsilon = \int_{-\infty}^{x_0} e^{-V(z)/\varepsilon}\,dz$ and the integral. The factor~$\frac{1}{2}$ stems from the fact that the exponential $e^{V/\varepsilon}$ achieves its maximum at the boundary of the interval~$[x,x_0]$. The argument for the case~$x>x_0$ is similar.
\end{proof}
\begin{definition}[Transformed left-normalized stationary measure]\label{GF_NGF:def:transformed-stationary-measure}
	For the measure~$\gamma_\varepsilon^\ell\in \mathcal{M}(\mathbb{R})$ from Definition~\ref{GF_NGF:def:left-normalized-stationary-measure} and the coordinate transformation~$y_\varepsilon$ from Definition~\ref{GF_NGF:def:coordinate-transformation}, we let~$\hat{\gamma}_\varepsilon^\ell\in\mathcal{M}(\mathbb{R})$ be the push-forward measure
	\begin{equation*}
	\hat{\gamma}_\varepsilon^\ell(A) := (y_\varepsilon)_\#\gamma_\varepsilon^\ell(A) = \gamma_\varepsilon^\ell\left(y_\varepsilon^{-1}(A)\right).
	\end{equation*}
\end{definition}
\begin{proposition}[Concentration of measure]\label{prop:GF_NGF:transformed-stationary-measure}
	For any~$\kappa>0$ small, let~$ U_\kappa$ be the neighborhood $U_\kappa := B_\kappa(-1/2) \cup B_\kappa(+1/2)$ of~$\{\pm1/2\}$. Then
	\begin{equation*}
	\lim_{\varepsilon\to 0} \hat{\gamma}_\varepsilon^\ell(\mathbb{R}\setminus U_\kappa) = 0.
	\end{equation*}
\end{proposition}
\begin{proof}[Proof of Proposition~\ref{prop:GF_NGF:transformed-stationary-measure}]
	Let~$\kappa>0$ and fix~$\delta>0$. By Proposition~\ref{GF_NGF:prop:coordinate-transformation}, if~$\varepsilon>0$ is small enough, then~$\{V<\delta\} \subseteq y_\varepsilon^{-1}(U_\kappa)$. Therefore~$\mathbb{R}\setminus y_\varepsilon^{-1}(U_\kappa) \subseteq \{V>\delta\}$, and we find that
	\begin{equation*}
	\hat{\gamma}_\varepsilon^\ell(\mathbb{R}\setminus U_\kappa) \leq \gamma_\varepsilon^\ell(\{V>\delta\}).
	\end{equation*}
	By Proposition~\ref{prop:GF_NGF:left-normalized-stationary-measure}, the right-hand side vanishes in the limit~$\varepsilon\to 0$.
\end{proof}
\subsection{Definition of flux-density functionals}
The flux-density functionals are defined on pairs of measures~$(\rho,j)$ satisfying the continuity equation $\partial_t\rho + \mathrm{div} j = 0$ in the following sense.
\begin{definition}[Continuity Equation]\label{GF_NGF:def:continuity-equation}
	Fix~$T>0$ and let~$E:=(0,T)\times\mathbb{R}$. We say that a pair~$(\rho(t,(\cdot),j(t,\cdot))$ of time-dependent Borel measures on~$\mathbb{R}$ satisfies the continuity equation if:
	\begin{enumerate}[label=(\roman*)]
		\item For each~$t\in(0,T)$, $\rho(t,\cdot)$ is a probability measure on~$\mathbb{R}$. The map~$t\mapsto \rho(t,\cdot)\in\mathcal{P}(\mathbb{R})$ is measurable with respect to the weak topology on~$\mathcal{P}(\mathbb{R})$.
		\item For each~$t\in(0,T)$, $j(t,\cdot)$ is a finite Borel measure on~$\mathbb{R}$. The map~$t\mapsto j(t,\cdot)\in \mathcal{M}(\mathbb{R})$ is measurable with respect to the weak topology on~$\mathcal{M}(\mathbb{R})$.
		\item The pair solves $\partial_t\rho + \mathrm{div}\, j = 0$ in~$\mathcal{D}'(E)$; that means for any test function $\varphi\in C_c^\infty(E)$, we have
		\begin{equation}
		\int_0^T\int_\mathbb{R} \left[\rho(t,dy)\, \partial_t \varphi(t,y) + j(t,dy)\, \partial_y \varphi(t,y) \right]\,dt = 0.
		\end{equation}
	\end{enumerate}
	We denote by~$\mathrm{CE}(0,T;\mathbb{R})$ the set of all pairs~$(\rho,j)$ satisfying the continuity equation.\qed
\end{definition}
As discussed in the previous section, with the coordinate transformation~$y_\varepsilon$, the flux-density rate function takes the form
\begin{equation*}
\mathcal{I}_\varepsilon(\rho,j) = \frac{1}{4}\int_0^T\int_\mathbb{R} \frac{1}{\hat{u}_\varepsilon^\ell(t,y)}\big|\hat{\jmath}_\varepsilon(t,y) + \partial_y \hat{u}_\varepsilon^\ell(t,y)\big|^2\,dydt.
\end{equation*}
We take the dual formulation for integrals over convex functions (Lemma~\ref{GF_NGF:lemma:appendix:dual-of-convex-functions}), which generalizes
\begin{equation*}
\frac{1}{2}\frac{x^2}{y} = \sup_{\substack{b \in \mathbb{R}}} \left[(-b^2/2)y + bx\right],\quad x\in\mathbb{R},\,y>0.
\end{equation*}
Shifting derivatives to test functions, we arrive at the following Definition.
\begin{definition}[Prel-limit Rate Function]\label{GF_NGF:def:pre-limit-RF}
	For~$E:=(0,T)\times\mathbb{R}$ and~$\varepsilon>0$, define the map $\widehat{\mathcal{I}}_\varepsilon:\mathrm{CE}(0,T;\mathbb{R})\to[0,\infty]$ by
	\begin{equation}\label{def:GF-NGF:pre-limitRF}
	\widehat{\mathcal{I}}_\varepsilon\left(\hat{\rho},\hat{\jmath}\right) := \frac{1}{2} \sup_{b\in C_c^\infty(E)}\int_E \left[\hat{u}_\varepsilon^\ell\left(-\partial_y b - \frac{1}{2}b^2\right) + \hat{\jmath}\cdot b\right]\,dydt,
	\end{equation}
	if~$\hat{\rho}(t,dy)=\hat{u}_\varepsilon^\ell(t,y) \hat{\gamma}_\varepsilon^\ell(dy)$, that is if~$\hat{\rho}(t,\cdot)$ is absolutely continuous with respect to~$\hat{\gamma}_\varepsilon^\ell$, the left-normalized stationary measure from Definition~\ref{GF_NGF:def:transformed-stationary-measure}. Otherwise, we set~$\widehat{\mathcal{I}}_\varepsilon\left(\hat{\rho},\hat{\jmath}\right)$ equal to~$+\infty$.\qed
\end{definition}
\begin{definition}[Limit Rate Function]\label{GF_NGF:def:limit-RF}
	With the funciton~$S$ from~\eqref{GF_NGF:eq:S-fct}, define~$\widehat{\mathcal{I}}_0:\mathrm{CE}(0,T;\mathbb{R})\to[0,\infty]$ by
	\begin{equation}
	\widehat{\mathcal{I}}_0\left(\hat{\rho},\hat{\jmath}\right) := \int_0^T S\left(\hat{\jmath}(t),\hat{z}(t)\right)\,dt,
	\end{equation}
	if~$\hat{\rho}(t,dy) = \hat{z}(t)\delta_{-\frac{1}{2}}(dy) + (1-\hat{z}(t))\delta_{+\frac{1}{2}}(dy)$ and $\hat{\jmath} = \hat{\jmath}(t)\mathbf{1}_{\left(-\frac{1}{2},+\frac{1}{2}\right)}(y)dtdy$, with coefficient~$\hat{z} \in H^1([0,T])$ and~$\hat{\jmath}(t)>0$. Otherwise, we set~$\widehat{\mathcal{I}}_0\left(\hat{\rho},\hat{\jmath}\right)$ equal to~$+\infty$.\qed
\end{definition}
In fact, we have~$\hat{\jmath}(t)=-\partial_t\hat{z}(t)$, which follows from the continuity equation~$\partial_t\hat{\rho}+\partial_y\hat{\jmath}=0$ and the special form of~$\hat{\rho}$ (see Theorem~\ref{GF_NGF:thm:compactness},~\ref{item:compact:j-is-nice}). That means the rate function is finite only if~$\hat{z}(t)$ is decreasing, meaning that mass is flowing only from left to right. We will show furthermore that finiteness of~$\mathcal{I}_0$ implies~$z\in W^{1,1}([0,T])$.
\subsection{Why the limit is not a gradient flow}
\label{GF_NGF:sec:why-limit-is-not-GF}
We give a formal argument. Suppose that~$\mathcal{I}_0$ is a generalized gradient-flow given by~\eqref{GF_NGF:eq:action-of-R-R_star-functional}, with an energy~$\mathrm{E}$ and dissipation potentials~$\mathcal{R}$,~$\mathcal{R}^\ast$. Then locally, we obtain
\begin{equation*}
S(j,z) = \dd \mathrm{E}(\rho) \cdot j + \mathcal{R}(\rho,j) + \mathcal{R}^\ast(\rho,-\dd\mathrm{E}(\rho)).
\end{equation*}
Taking the derivative with respect to~$j$, 
\begin{equation*}
\partial_j S(j,z) = \dd\mathrm{E}(\rho) + \partial_j \mathcal{R}(\rho,j).
\end{equation*}
Since~$\mathcal{R}(\rho,0)=0$ is minimal by definition, 
\begin{equation*}
\partial_j S(0,z) = \dd \mathrm{E}(\rho) + 0.
\end{equation*}
But~$\partial_jS(j,z) = \log(j/z)$ diverges to~$-\infty$ as~$j\to 0$.
\smallskip

Informally, the limit structure is not an \emph{entropy} gradient flow because there is no decay in entropy. As we saw in Section~\ref{GF_NGF:sec:why-GF-not-converge}, the entropies diverge in the limit~$\varepsilon\to 0$. In the limit, the mass of the stationary measure is concentrated on one point, while the dynamics has mass on both points. Hence for any finite time, the relative entropy equals~$+\infty$ and is not decaying.
\section{Proof of Gamma-convergence}
\label{GF_NGF:sec:proof-of-gamma-convergence}
In~$\mathrm{CE}(0,T;\mathbb{R})$, we consider convergence in distribution.
\begin{definition}[Convergence of solutions to continuity equation]\label{GF_NGF:def:converge-in-CE} We say that a
	sequence~$(\hat{\rho}_n,\hat{\jmath}_n)$ in~$\mathrm{CE}(0,T;\mathbb{R})$ converges to~$(\hat{\rho},\hat{\jmath})$ if and only if for any test function~$\varphi\in C_c^\infty((0,T)\times\mathbb{R})$, we have
	\begin{equation*}
	\int_{0}^T\int_\mathbb{R} \left[\hat{\rho}_n \partial_t\varphi + \hat{\jmath}_n \,\mathrm{div}\,\varphi\right]\,dydt \xrightarrow{n\to\infty}
	\int_{0}^T\int_\mathbb{R} \left[\hat{\rho} \partial_t\varphi + \hat{\jmath}\, \mathrm{div}\,\varphi\right]\,dydt.
	\end{equation*}
\end{definition}
\begin{theorem}[Lower Bound]\label{GF_NGF:thm:lower-bound}
	For any sequence~$(\hat{\rho}_\varepsilon,\hat{\jmath}_\varepsilon)\in\mathrm{CE}(0,T;\mathbb{R})$ such that
	\begin{equation*} 
		(\hat{\rho}_\varepsilon,\hat{\jmath}_\varepsilon)\to(\hat{\rho}_0,\hat{\jmath}_0)\in\mathrm{CE}(0,T;\mathbb{R})\quad\text{and}\quad \sup_{(t,y)\in E}|\hat{u}_\varepsilon^\ell(t,y)| \leq C,
	\end{equation*}
	where~$\hat{\rho}_\varepsilon=\hat{u}_\varepsilon^\ell\hat{\gamma}_\varepsilon^\ell$, we have
	\begin{equation*}
	\liminf_{\varepsilon\to 0}\widehat{\mathcal{I}}_\varepsilon(\hat{\rho}_\varepsilon,\hat{\jmath}_\varepsilon) \geq \widehat{\mathcal{I}}_0(\hat{\rho}_0,\hat{\jmath}_0).
	\end{equation*}
\end{theorem}
\begin{theorem}[Upper bound]\label{GF_NGF:thm:upper-bound}
	For any~$(\hat{\rho}_0,\hat{\jmath}_0)\in\mathrm{CE}(0,T;\mathbb{R})$ such that the rate function is finite, $\widehat{\mathcal{I}}_0(\hat{\rho}_0,\hat{\jmath}_0)<\infty$, there exist~$(\hat{\rho}_\varepsilon,\hat{\jmath}_\varepsilon)\in\mathrm{CE}(0,T;\mathbb{R})$ such that
	\[
	(\hat{\rho}_\varepsilon,\hat{\jmath}_\varepsilon)\xrightarrow{\varepsilon\to 0}(\hat{\rho}_0,\hat{\jmath}_0)\quad\text{and}\quad \limsup_{\varepsilon\to 0}\widehat{\mathcal{I}}_\varepsilon(\hat{\rho}_\varepsilon,\hat{\jmath}_\varepsilon) \leq \widehat{\mathcal{I}}_0(\hat{\rho}_0,\hat{\jmath}_0).
	\]
\end{theorem}
\subsection{Proof of compactness}
Recall we denote by~$\mathrm{CE}(0,T;\mathbb{R})$ the set of solutions to the continuity equation in the sense of Definition~\ref{GF_NGF:def:continuity-equation}, and by~$\widehat{\mathcal{I}}_\varepsilon:\mathrm{CE}(0,T;\mathbb{R})\to[0,\infty]$ the rate function from Definition~\ref{GF_NGF:def:pre-limit-RF}. The measures~$\hat{\gamma}_\varepsilon^\ell$ are the transformed left-normalized stationary measures introduced in Definition~\ref{GF_NGF:def:transformed-stationary-measure}.
\begin{theorem}[Sequential Compactness]\label{GF_NGF:thm:compactness}
	For $\varepsilon>0$, let $(\hat{\rho}_\varepsilon,\hat{\jmath}_\varepsilon)\in \mathrm{CE}(0,T;\mathbb{R})$ be pairs such that~$\hat{\rho}_\varepsilon(t,\cdot)$ is absolutely continuous w.r.t.~$\hat{\gamma}_\varepsilon^\ell$, with Radon-Nikodym derivative $\hat{u}_\varepsilon^\ell(t,\cdot)$. Let~$E:=(0,T)\times\mathbb{R}$. Suppose that there exists a constant~$C>0$ such that for all~$\varepsilon>0$,
	\begin{equation*}
	\widehat{\mathcal{I}}_\varepsilon\left(\hat{\rho}_\varepsilon,\hat{\jmath}_\varepsilon\right) \leq C \quad\text{and}\quad 
	\sup_{(t,y)\in E}|\hat{u}_\varepsilon^\ell(t,y)| \leq C.
	\end{equation*}
	Then there exists a pair~$(\hat{\rho}_0,\hat{\jmath}_0)\in\mathrm{CE}(0,T;\mathbb{R})$ and a limiting function~$\hat{u}_0^\ell$ such that:
	\begin{enumerate}[label=(\roman*)]
		\item \label{item:compact:u-converge}$\hat{u}_0^\ell \in L^\infty(E)$ and along a subsequence~$\hat{u}_\varepsilon^\ell\xrightharpoonup{\ast}\hat{u}_0^\ell$ in~$\sigma(L^\infty(E),L^1(E))$.
		\item \label{item:compact:rho-converge}The densities converge weakly: $\hat{\rho}_\varepsilon\xrightharpoonup{\ast}\hat{\rho}_0$ w.r.t.~$\sigma(\mathcal{M}(E),C_b(E))$, where for every $t\in[0,T]$, the measure $\hat{\rho}_0(t,\cdot) \in \mathcal{P}(\mathbb{R})$ is absolutely continuous with respect to $\hat{\gamma}_0 := \delta_{-1/2}+\delta_{+1/2}$.
		\item \label{item:compact:j-converge} Along a subsequence: $\hat{\jmath}_\varepsilon\xrightharpoonup{\ast}\hat{\jmath}_0$ in~$\sigma(Y^\ast,Y)$, where $Y:=L^2\left(0,T;X\right)$ with the Banach space~$X$ obtained by taking the closure of~$C_c^\infty(\mathbb{R})$ under the norm
		\[
		\|f\|_X := \|f\|_{L^2(\mathbb{R})} + \|\partial_y f\|_{L^1(\mathbb{R})}.
		\]
	\end{enumerate}
	Moreover, the limiting objects~$\hat{u}_0^\ell, \hat{\rho}_0$ and~$\hat{\jmath}_0$ satisfy the following regularity properties:
	\begin{enumerate}[label=(\roman*)]
		\setcounter{enumi}{3}
		\item \label{item:compact:rho-is-nice} The density
		$\hat{z}_0(t)$ such that~$\hat{\rho}_0(t,dy)=\hat{z}_0(t)\delta_{-1/2}(dy) + (1-\hat{z}_0(t))\delta_{+1/2}(dy)$ satisfies~$\hat{z}_0 \in H^1([0,T])$.
		\item \label{item:compact:j-is-nice} The limiting flux~$\hat{\jmath}_0$ is unique and is given by
		\begin{equation*}
		\hat{\jmath}_0(dt\,dy) = -\partial_t \hat{z}_0(t) \mathbf{1}_{(-1/2,+1/2)}(y)\,dtdy.
		\end{equation*}
		\item \label{item:compact:u-is-nice} The limiting function~$\hat{u}_0^\ell$ satisfies for a.e.~$t\in(0,T)$
		\begin{equation*}
		\partial_y\hat{u}_0^\ell \in L^2\left((0,T)\times\mathbb{R}\right),\quad \hat{u}_0^\ell(t,-1/2) = \hat{z}_0(t) \quad\text{and}\quad \hat{u}_0^\ell(t,+1/2) = 0.
		\end{equation*}
	\end{enumerate}
\end{theorem}
\begin{proof}[Proof of Theorem~\ref{GF_NGF:thm:compactness}]
	If both~\ref{item:compact:rho-converge} and~\ref{item:compact:j-converge} hold, then a limiting pair~$(\hat{\rho}_0,\hat{\jmath}_0)$ satisfies the continuity equation: for any~$\varphi\in C_c^\infty(E)$, $\partial_t\varphi\in C_c^\infty(E)$ and~$\partial_y\varphi \in Y$, so that
	\[
	\int_E \partial_t\varphi\, \hat{\rho}_\varepsilon \xrightarrow{\varepsilon\to 0}\int_E \partial_t \varphi\,\hat{\rho}_0\quad\text{and}\quad \int_E \partial_y\varphi(t,y)\,\hat{\jmath}_\varepsilon \xrightarrow{\varepsilon\to 0} \int_E \partial_y\varphi(t,y)\,\hat{\jmath}_0,
	\]
	Hence $\partial_t\hat{\rho}_0+\partial_y\hat{\jmath}_0=0$ in~$\mathcal{D}'$ is inherited from $\partial_t\hat{\rho}_\varepsilon+\partial_y\hat{\jmath}_\varepsilon=0$ in~$\mathcal{D}'$.
	\smallskip
	
	For proving~\ref{item:compact:u-converge} and~\ref{item:compact:rho-converge}, we only exploit the uniform boundedness assumption on the densities~$\hat{u}_\varepsilon^\ell$. The remaining properties follow from the boundedness assumption on the rate function.
	\smallskip
	
	\underline{\ref{item:compact:u-converge}}: The family of functions $\{\hat{u}_\varepsilon^\ell\}_{\varepsilon>0}$ is bounded in~$L^\infty(E)$, the topological dual of $L^1(E)$. Hence by the Banach-Alaoglu Theorem, there is a subsequence converging in~$\sigma(L^\infty(E),L^1(E))$.
	\smallskip
	
	\underline{\ref{item:compact:rho-converge}}: For any Borel subset $I\times A \subseteq E$,
	\begin{align*}
	|\hat{\rho}_\varepsilon(I\times A)| = \left|\int_I\int_A \hat{u}_\varepsilon^\ell(t,y)\hat{\gamma}_\varepsilon^\ell(dy)dt\right|&\leq T \cdot \|\hat{u}_\varepsilon^\ell\|_{L^\infty(E)}\cdot \hat{\gamma}_\varepsilon^\ell(A)\\
	&\leq T\cdot C \cdot \hat{\gamma}_\varepsilon^\ell(A).
	\end{align*}
	For~$\kappa>0$, let $U_\kappa := B_\kappa(-1/2)\cup B_\kappa(+1/2)$. Then $\hat{\gamma}_\varepsilon^\ell(U_\kappa)\to 0$ as~$\varepsilon\to 0$ (Proposition~\ref{prop:GF_NGF:transformed-stationary-measure}). Since~$\kappa>0$ is arbitrary, this means that~$\hat{\rho}_\varepsilon$ converges weakly to a measure~$\hat{\rho}_0$ that is supported on the set~$(0,T)\times \{-\frac{1}{2},+\frac{1}{2}\}$.
	\smallskip
	
	\underline{\ref{item:compact:j-converge}}: By the uniform-boundedness assumption of the rate function,
	\[
	\sup_{b\in C_c^\infty(E)}\int_0^T\int_\mathbb{R} \left[\hat{u}_\varepsilon^\ell \left(-\partial_y b - \frac{1}{2}b^2\right) + \hat{\jmath}_\varepsilon \cdot b\right]\,dtdy \leq C < \infty.
	\]
	Therefore, for any~$b\in C_c^\infty(E)$,
	\begin{align*}
	\left|\int_0^T\int_\mathbb{R} \hat{\jmath}_\varepsilon\cdot b(t,y)\,dydt\right| &\leq C + \|\hat{u}_\varepsilon^\ell\|_{\infty} \int_0^T\int_\mathbb{R} \left(|\partial_y b| + \frac{1}{2}b^2\right)\,dydt\\
	&\leq C'\left[1+\int_0^T\left(\|b(t,\cdot)\|_{L^2(\mathbb{R})}^2 + \|\partial_yb(t,\cdot)\|_{L^1(\mathbb{R})}\right)\,dt\right],
	\end{align*}
	where $C'=\max(C,\|\hat{u}_\varepsilon^\ell\|)$.
	We henceforth abbreviate integrals $\int_0^T\int_\mathbb{R} \hat{\jmath}_\varepsilon b\,dtdy$ simply by $\langle \hat{\jmath}_\varepsilon,b\rangle$. Rescaling in the above estimate as~$b\to \lambda b$, we obtain that for any~$\lambda>0$ and~$b\in C_c^\infty(E)$,
	\begin{align*}
	\left|\langle \hat{\jmath}_\varepsilon, b \rangle\right| \leq C'\left[\frac{1}{\lambda} + \lambda \int_0^T\|b(t,\cdot)\|_2^2\,dt + \int_0^T\|\partial_y b(t,\cdot)\|_1 \,dt\right].
	\end{align*}
	Optimizing the right-hand side in~$\lambda$ gives
	\[
	\lambda_{\mathrm{op}} = \left(\int_0^T\|b(t,\cdot)\|_2^2\,dt\right)^{-1/2}.
	\]
	With that optimal~$\lambda=\lambda_\mathrm{op}$, we found the estimate
	\begin{align*}
	\left|\langle \hat{\jmath}_\varepsilon, b \rangle\right| \leq C'\left[\left(\int_0^T\|b(t,\cdot)\|_2^2\,dt\right)^{1/2} + \int_0^T\|\partial_y b(t,\cdot)\|_1 \,dt\right].
	\end{align*}
	With the elementary estimate $(a+b)^2 \leq 2a^2+2b^2$, we arrive at
	\begin{align*}
	\left|\langle \hat{\jmath}_\varepsilon, b \rangle\right|^2 &\leq 2C'^2\left[\int_0^T\|b(t,\cdot)\|_2^2\,dt + \left(\int_0^T\|\partial_y b(t,\cdot)\|_1 \,dt\right)^2\right] \\
	&\leq 2C'^2\left[\int_0^T\|b(t,\cdot)\|_2^2\,dt + \int_0^T\|\partial_y b(t,\cdot)\|_1^2 \,dt\right],
	\end{align*}
	where the second estimate is a consequence of Jensen's inequality. Therefore, for some constant~$C>0$ and for all~$b\in C_c^\infty(E)$,
	\begin{align*}
	\left|\langle \hat{\jmath}_\varepsilon, b \rangle\right| \leq C \|b\|_{Y}.
	\end{align*}
	Hence~$\hat{\jmath}_\varepsilon$ is bounded in~$Y^\ast$, and by the Banach-Alaoglu Theorem, there exists a converging subsequence of~$\hat{\jmath}_\varepsilon$ in the~$\sigma(Y^\ast,Y)$ topology.
	\smallskip
	
	\underline{\ref{item:compact:rho-is-nice}}: The density~$\hat{z}_0(t)$ is measurable as the limit of measurable functions, since
	\[
	\hat{z}_0(t) = \lim_{\varepsilon\to 0} \hat{\rho}_\varepsilon(t,U_-)
	\]
	for $U_-$ a small neighborhood of~$- 1/2$. We now prove the claimed regularity.
	For a test function~$b\in C_c^\infty(E)$, we write $B(t,y) := \int_{-\frac{1}{2}}^y b(t,z)\,dz$, so that in particular $\partial_y B = b$. Since the pair~$(\hat{\rho}_\varepsilon,\hat{\jmath}_\varepsilon)$ satisfies the continuity equation, we obtain for any~$b$ such that~$B$ is compactly supported that 
	\[
	\int_0^T\int_\mathbb{R} \partial_y B\,\hat{\jmath}_\varepsilon = - \int_0^T\int_\mathbb{R} \partial_t B \,\hat{\rho}_\varepsilon.
	\]
	From the boundedness of the rate function, we find for such~$b$ the estimate
	\begin{align*}
	C\geq \widehat{\mathcal{I}}_\varepsilon(\hat{\rho}_\varepsilon,\hat{\jmath}_\varepsilon) &\geq \frac{1}{2}\int_0^T \int_\mathbb{R} \left[ \hat{u}_\varepsilon^\ell\left( - \partial_y b - \frac{1}{2}b^2\right) + \hat{\jmath}_\varepsilon \partial_y B \right]\,dydt	\\
	&= \frac{1}{2} \int_0^T \int_\mathbb{R} \left[ \hat{u}_\varepsilon^\ell\left( - \partial_y b - \frac{1}{2}b^2\right) - \hat{\rho}_\varepsilon \partial_t B \right]\,dydt.
	\end{align*}
	We specialize further to functions~$b$ such that~$B(t,y)=\varphi(t)\psi(y)$, where the function~$\varphi \in C_c^\infty(0,T)$ is arbitrary and~$\psi$ is a fixed function that has compact support and satisfies $\psi(-1/2)=0$ and $\psi(1/2)=1$. Writing~$\hat{z}_+:=1-\hat{z}_0$,
	\begin{align*}
	\left|\int_0^T \varphi'(t) \hat{z}_+(t)\,dt\right| &=
	\left|\int_0^T\int_\mathbb{R} \partial_t B(t,y)\, \hat{\rho}_0(dtdy) \right|	\\
	&\leq
	C + \left|\int_0^T\int_\mathbb{R} \hat{u}_0^\ell \left(\varphi(t) \psi''(y) + \frac{1}{2}\varphi(t)^2 |\psi'(y)|^2\right)\,dydt\right|	\\
	&\leq C + \|\hat{u}_0^\ell\|_\infty \left(\|\psi''\|_\infty \|\varphi\|_{L^1}+\frac{1}{2} \|\psi'\|_2^2 \|\varphi\|_{L^2}^2 \right)	\\
	&\leq C\left(1 + \|\varphi\|_{L^1} + \|\varphi\|_{L^2}^2\right).
	\end{align*}
	Rescaling as~$\varphi\to\lambda\varphi$ and optimizing the resulting estimate in~$\lambda$ in the same fashion as in the proof of~\ref{item:compact:j-converge} above, we arrive at
	\begin{align*}
	\left|\int_0^T\varphi'(t)\hat{z}_+(t)\,dt\right| \leq C\left(\|\varphi\|_{L^2} + \|\varphi\|_{L^1}\right) \leq C' \|\varphi\|_{L^2}, 
	\end{align*}
	where the second estimate uses $\|\varphi\|_{L^1} \leq \sqrt{T}\|\varphi\|_{L^2}$. Since this bound holds for any~$\varphi\in C_c^\infty(0,T)$, we obtain that~$\hat{z}_+\in H^1(0,T)$, and hence also~$\hat{z}_0 \in H^1(0,T)$.
	\smallskip
	
	\underline{\ref{item:compact:j-is-nice}}: Let~$\hat{\jmath}_0$ be such that~$\hat{\jmath}_\varepsilon\xrightharpoonup{\ast}\hat{\jmath}_0$ in~$\sigma(Y^\ast,Y)$. First, we show that~$\hat{\jmath}_0$ is piecewise constant in the sense that
	\[
	\hat{\jmath}_0(dt\,dy) = \left[j^-(t)\mathbf{1}_{(-\infty,-\frac{1}{2})}(y) + j(t)\mathbf{1}_{(-\frac{1}{2},+\frac{1}{2})}(y) +  j^+(t)\mathbf{1}_{(+\frac{1}{2},+\infty)}(y)\right]\,dtdy,
	\] 
	where~$j^-,j,j^+$ are measurable functions. Secondly, we verify that
	\[
	j^-\equiv 0,\quad j\equiv -\partial_t \hat{z}_0,\quad j^+\equiv 0.
	\]
	Combining the two statements proves the claim.
	\smallskip
	
	The limiting density is of the form $\hat{\rho}_0(t,dy) = \hat{z}_-(t)\delta_{-\frac{1}{2}}(dy) + \hat{z}_+(t)\delta_{+\frac{1}{2}}(dy)$. Specializing the continuity equation $\partial_t\hat{\rho}_0 + \partial_y\hat{\jmath}_0=0$ to test functions of the form $b(t,y) = \varphi(t)\psi(y)$, where~$\varphi\in C_c^\infty(0,T)$ and~$\psi\in C_c^\infty((-\infty,-1/2))$, we find
	\begin{align*}
	\int_0^T\left[\int_\mathbb{R}\partial_y\psi(y)\hat{\jmath}_0(t,dy)\right]\varphi(t)\,dt = 0
	\end{align*}
	Therefore $\partial_y\hat{\jmath}_0 = 0$ in $(0,T)\times(-\infty,-1/2)$, since~$\varphi,\psi$ are arbitrary. Repeating the arument on $(-1/2,+1/2)$ and $(+1/2,\infty)$, we find that~$\hat{\jmath}_0$ is piecewise constant as claimed. We are left with verifying that the flux vanishes outside the interval~$(-1/2,+1/2)$ and in the interval is given by~$\partial_t\hat{z}_+=-\partial_t\hat{z}_0$. 
	\smallskip
	
	By boundedness of the rate function, for any~$b\in C_c^\infty(E)$, 
	\begin{align*}
	C \geq \int_0^T\int_\mathbb{R} \left[\hat{u}_\varepsilon^\ell\left(-\partial_y b - \frac{1}{2}b^2\right)+\hat{\jmath}_\varepsilon b\right]\,dydt.
	\end{align*}
	The densities~$\hat{u}_\varepsilon^\ell$ converge along a subsequence in~$\sigma(L^\infty(E),L^1(E))$, the fluxes~$\hat{\jmath}_\varepsilon$ converge in~$\sigma(Y^\ast,Y)$, and any test function~$b$ together with its derivatives is both in $L^1(E)$ and~$Y$. Therefore we can pass to the limit to obtain
	\begin{align*}
	C \geq \int_0^T\int_\mathbb{R} \left[\hat{u}_0^\ell\left(-\partial_y b - \frac{1}{2}b^2\right)+\hat{\jmath}_0 b\right]\,dydt.
	\end{align*}
	Specializing to a sequence of functions~$b_n=\varphi(t)\psi_n(y)$ with functions~$\psi_n$ that are supported in~$(-\infty,-1)$ and satisfy
	\[
	\|\psi_n\|_{L^1} \xrightarrow{n\to\infty}\infty,\quad\|\partial_y\psi_n\|_{L^1}\leq C\quad\text{and}\quad \|\psi_n^2\|_{L^1} \xrightarrow{n\to\infty}0,
	\]
	we find that
	\begin{align*}
	\left|\int_0^T \varphi(t)j^-(t)\,dt\right|\cdot\left| \int_\mathbb{R}\psi_n \,dy\right| &= \left|\int_0^T\int_{-\infty}^{-1}b_n \hat{\jmath}_0\,dydt\right|\\
	&\leq C \left[1 + T\|\varphi\|_\infty\left( \|\hat{u}_\varepsilon^\ell\|_\infty \|\partial_y\psi_n\|_{L^1} + \|\hat{u}_\varepsilon^\ell\|_\infty \|\psi_n^2\|_{L^1}\right)\right]\\
	&\leq C'.
	\end{align*}
	Since~$\|\psi_n\|_{L^1}\to\infty$ and since~$\varphi$ is arbitrary, this implies~$j^-\equiv 0$.	Examples of~$\psi_n$ are smoothend versions of the step functions~$n^{-2/3}\mathbf{1}_{(-n,-1)}$. The argument for the region~$(+1/2,+\infty)$ is similar. Therefore,
	\[
	\hat{\jmath}_0(dt,dy) = j(t)\mathbf{1}_{\left(-\frac{1}{2},+\frac{1}{2}\right)}(y)\,dtdy.
	\]
	Testing the continuity equation $\partial_t\hat{\rho}_0+\partial_y\hat{\jmath}_0=0$ with functions~$b(t,y)=\varphi(t)\psi(y)$ such that~$\psi(-1/2)=0$ and~$\psi(+1/2)=1$, we find
	\begin{align*}
	\int_0^T\hat{z}_+(t) \partial_t\varphi(t)\,dt + \int_0^Tj(t)\varphi(t)\,dt = 0.
	\end{align*}
	Since~$\hat{z}_+$ is in~$H^1$, integration by parts and arbitraryness of the test function~$\varphi$ imply that the flux is given by~$j(t)=\partial_t\hat{z}_+(t)$.
	\smallskip
	
	\underline{\ref{item:compact:u-is-nice}}: As shown above, we have for any~$b\in C_c^\infty(E)$ the bound
	\[
	C\geq \int_E \left[\hat{u}_0^\ell\left(-\partial_y b - \frac{1}{2}b^2\right) + \hat{\jmath}_0 b\right]\,dydt.
	\]
	Using that $\hat{\jmath}_0\in L^2(E)$, we find
	\begin{align*}
	\left|\int_E \hat{u}_0^\ell \, \partial_y b\,dydt\right| &\leq C + \frac{1}{2} \|\hat{u}_0^\ell\|_\infty \|b\|_{L^2(E)}^2 + \|\hat{\jmath}_0\|_{L^2(E)} \|b\|_{L^2(E)} \\
	&\leq C'\|b\|_{L^2(E)},
	\end{align*}
	where the second estimate follows after rescaling $b\to\lambda b$ and optimizing in~$\lambda$. This shows boundedness of the map
	\[
	b\mapsto \langle\hat{u}_0^\ell,\partial_y b\rangle_{L^2(E)}.
	\]
	Since~$L^2(E)$ is self-dual, the fact that~$\partial_y \hat{u}_0^\ell\in L^2(E)$ follows by the Banach-Alaoglu Theorem.
	\smallskip
	
	We now show $\hat{u}_0^\ell(t,-1/2) = \hat{z}_0(t)$ for a.e.~$t\in(0,T)$. The density~$\hat{z}_0(t)$ satisfies
	\[
	\hat{z}_0(t) = \lim_{\varepsilon\to 0}\hat{\rho}_\varepsilon(t,U_\kappa),
	\]
	where~$U_\kappa$ is a small ball of radius~$\kappa$ around~$\{-1/2\}$. In~$U_\kappa$, the densities~$\hat{u}_\varepsilon^\ell$ and~$\hat{u}_0^\ell$ are close in the sense that
	\[
	\int_{U_\kappa}\hat{u}_\varepsilon^\ell\,\hat{\gamma}_\varepsilon^\ell = \int_{U_\kappa}\hat{u}_0^\ell\,\hat{\gamma}_\varepsilon^\ell + \int_{U_\kappa}(\hat{u}_\varepsilon^\ell-\hat{u}_0^\ell)\,\hat{\gamma}_\varepsilon^\ell = \int_{U_\kappa}\hat{u}_0^\ell\,\hat{\gamma}_\varepsilon^\ell + o(1)_{\varepsilon\to 0},
	\] 
	since~$\hat{\gamma}_\varepsilon^\ell$ is left-normalized and concentrates on~$U_\kappa$ in the sense that 
	\[
	\hat{\gamma}_\varepsilon^\ell(U_\kappa) = \frac{1}{Z_\varepsilon^\ell}\int_{U_\kappa} e^{-V(y	)/\varepsilon}\,dy \xrightarrow{\varepsilon\to 0}1,
	\]
	and since~$\hat{u}_\varepsilon^\ell\xrightharpoonup{\ast} \hat{u}_0^\ell$ in~$\sigma(L^\infty(E),L^1(E))$. With~$\hat{\rho}_\varepsilon(t,dy) = \hat{u}_\varepsilon^\ell(t,y)\hat{\gamma}_\varepsilon^\ell(dy)$, we therefore find
	\begin{align*}
	o(1)_{\varepsilon\to 0} + \left(\inf_{U_\kappa}\hat{u}_0^\ell(t,\cdot)\right)\hat{\gamma}_\varepsilon^\ell(U_\kappa) \leq \hat{\rho}_\varepsilon(t,U_\kappa) \leq \left(\sup_{U_\kappa}\hat{u}_0^\ell(t,\cdot)\right)\hat{\gamma}_\varepsilon^\ell(U_\kappa) + o(1)_{\varepsilon\to 0}.
	\end{align*}
	First passing to the limit~$\varepsilon\to 0$ and then taking~$\kappa\to 0$ gives
	\begin{align*}
	\lim_{\kappa\to 0}\inf_{U_\kappa}\hat{u}_0^\ell(t,\cdot) \leq \hat{z}_0(t) \leq \lim_{\kappa \to 0}\sup_{U_\kappa}\hat{u}_0^\ell(t,\cdot).
	\end{align*}
	Hence~$\hat{z}_0(t)$ is bounded from below by the lower-semicontinuous regularization of~$\hat{u}_0^\ell(t,y)$ at~$y=-1/2$, and from above by the upper-semicontinuous regularization. Since~$\hat{u}_0^\ell(t,\cdot)\in H^1(\mathbb{R})$ for almost every~$t\in(0,T)$, the function~$\hat{u}_0^\ell(t,\cdot)$ is continuous, implying that
	\[
	\hat{u}_0^\ell(t,-1/2) = \lim_{\kappa \to 0}\sup_{U_\kappa}\hat{u}_0^\ell(t,\cdot) = \lim_{\kappa\to 0}\inf_{U_\kappa}\hat{u}_0^\ell(t,\cdot).
	\]
	Hence~$\hat{u}_0^\ell(t,-1/2)=\hat{z}_0(t)$ for a.e.~$t\in(0,T)$.
	\smallskip
	
	The fact that~$\hat{u}_0^\ell(t,+1/2)=0$ follows from observing that in the limit~$\varepsilon\to 0$, the left-normalized measure~$\hat{\gamma}_\varepsilon^\ell$ blows up in a neighborhood~$U_\kappa^+$ of~$\{+\frac{1}{2}\}$ while
	\[
	1\geq \hat{\rho}_\varepsilon(t,U_\kappa^+) = \int_{U_\kappa^+}\hat{u}_\varepsilon^\ell(t,y)\,\hat{\gamma}_\varepsilon^\ell(dy)
	\]
	remains bounded. Hence continuity of~$\hat{u}_0^\ell(t,\cdot)$ enforces~$\hat{u}_0^\ell(t,+1/2)=0$. 
\end{proof}
\subsection{Proof of lower bound}
\begin{proof}[Proof of Theorem~\ref{GF_NGF:thm:lower-bound}]
	The limiting measure~$\hat{\rho}_0$ is supported on the set $(0,T)\times\{-1/2,+1/2\}$, that is~$\hat{\rho}_0=\hat{z}_0(t)\delta_{-1/2}+(1-\hat{z}_0(t))\delta_{+1/2}$, and the limiting flux~$\hat{\jmath}_0$ is piecewise constant, given as in~\ref{item:compact:j-is-nice} of Theorem~\ref{GF_NGF:thm:compactness}. 
	By definition of~$\widehat{\mathcal{I}}_\varepsilon$, for any~$b\in C_c^\infty(E)$,
	\begin{align*}
	C\geq \widehat{\mathcal{I}}_\varepsilon(\hat{\rho}_\varepsilon,\hat{\jmath}_\varepsilon) &\geq \frac{1}{2}\int_0^T\int_\mathbb{R} \left[\hat{u}_\varepsilon^\ell\left(-\partial_y b - \frac{1}{2}b^2\right)+\hat{\jmath}_\varepsilon b\right]\,dydt.
	\end{align*}
	As in the proof of Theorem~\ref{GF_NGF:thm:compactness}, we can pass to the limit since~$\hat{u}_\varepsilon^\ell$ and~$\hat{\jmath}_\varepsilon$ converge: $\hat{u}_\varepsilon^\ell\xrightarrow{\ast}\hat{u}_0^\ell$ in $\sigma(L^\infty(E),L^1(E))$ and~$\hat{\jmath}_\varepsilon\xrightarrow{\varepsilon\to 0}\hat{\jmath}_0$ in $\sigma(Y^\ast,Y)$,
	where the limiting flux is given by~$\hat{\jmath}_0(dt,dy)=\hat{\jmath}_0(t)\mathbf{1}_{(-1/2,+1/2)}(y)\,dtdy$ with~$\hat{\jmath}_0(t)=-\partial_t\hat{z}_0(t)$.
	This leads to
	\begin{align*}
	\liminf_{\varepsilon\to 0}\widehat{\mathcal{I}}_\varepsilon(\hat{\rho}_\varepsilon,\hat{\jmath}_\varepsilon) &\geq  \frac{1}{2}\int_0^T\int_\mathbb{R} \left[\hat{u}_0^\ell\left(-\partial_y b - \frac{1}{2}b^2\right)+\hat{\jmath}_0 b\right]\,dydt\\
	&= \frac{1}{2}\int_0^T\int_\mathbb{R}\left[\hat{u}_0^\ell\left(-\frac{1}{2}b^2\right)+b\left(\hat{\jmath}_0 + \partial_y\hat{u}_0^\ell\right)\right]\,dydt,
	\end{align*}
	using that~$\partial_y \hat{u}_0^\ell \in L^2(E)$ by~Theorem~\ref{GF_NGF:thm:compactness} and integration by parts.
	Taken the supremum over smooth functions~$b$ that have compact support in $E_0:=(0,T)\times(-1/2,+1/2)$, we find that
	\begin{align*}
	\liminf_{\varepsilon\to 0}\widehat{\mathcal{I}}_\varepsilon(\hat{\rho}_\varepsilon,\hat{\jmath}_\varepsilon) &\geq
	\frac{1}{2}\sup_{b\in C_c^\infty(E_0)}\int_0^T\int_{-1/2}^{+1/2} \left[\hat{u}_0^\ell\left(- \frac{1}{2}b^2\right)+ b\left(\hat{\jmath}_0(t)+\partial_y\hat{u}_0^\ell\right)\right]\,dydt\\
	&\overset{\ref{GF_NGF:lemma:appendix:dual-of-convex-functions}}{=} \frac{1}{4} \int_0^T\int_{-1/2}^{+1/2}\frac{1}{\hat{u}_0^\ell}\left|\hat{\jmath}_0(t)+\partial_y\hat{u}_0^\ell\right|^2\,dydt, 
	\end{align*}
	the last equality following from Lemma~\ref{GF_NGF:lemma:appendix:dual-of-convex-functions}. For fixed~$t>0$, we have
	\begin{equation*}
	\int_{-1/2}^{+1/2}\frac{1}{\hat{u}_0^\ell}\big|\hat{\jmath}_0(t) + \partial_y\hat{u}_0^\ell\big|^2\,dy \geq \inf_{\substack{\hat{u}=\hat{u}(y)\\ \hat{u}(\pm 1/2) = \hat{u}_0^\ell(t,\pm 1/2)}} \int_{-1/2}^{+1/2}\frac{1}{\hat{u}(y)}\big|\hat{\jmath}_0(t) + \partial_y \hat{u}\big|^2\,dy.
	\end{equation*}
	By Theorem~\ref{GF_NGF:thm:compactness}, the boundary conditions are given by~$\hat{u}_0^\ell(t,+1/2)=0$ and $\hat{u}_0^\ell(t,-1/2)=\hat{z}_0(t)$. With the function $S(a,b):=a\log(a/b) - (a-b)$ from~\eqref{GF_NGF:eq:S-fct}, the infimum is
	\begin{equation*}
	\inf_{\substack{\hat{u}(- 1/2) = \hat{z}_0(t)\\\hat{u}(t,+1/2)= 0}} \int_{-1/2}^{+1/2}\frac{1}{\hat{u}(y)}\big|\hat{\jmath}_0(t) + \partial_y \hat{u}\big|^2\,dy = 4 \, S\left(\hat{\jmath}_0(t), \hat{z}_0(t)\right),
	\end{equation*}
	which we prove in Lemma~\ref{lemma:GF_NGF:variational-problem}.
	Therefore,
	\begin{align*}
	\liminf_{\varepsilon\to 0}\widehat{\mathcal{I}}_\varepsilon(\hat{\rho}_\varepsilon,\hat{\jmath}_\varepsilon) &\geq
	\frac{1}{4} \int_0^T\int_{-1/2}^{+1/2}\frac{1}{\hat{u}_0^\ell}\left|\hat{\jmath}_0(t)+\partial_y\hat{u}_0^\ell\right|^2\,dydt\\
	&\geq \frac{1}{4}\int_0^T\inf_{\substack{\hat{u}(- 1/2) = \hat{z}_0(t)\\\hat{u}(+1/2)= 0}} \int_{-1/2}^{+1/2}\frac{1}{\hat{u}(y)}\big|\hat{\jmath}_0(t) + \partial_y \hat{u}\big|^2\,dy\,dt\\
	&= \int_0^T S\left(\hat{\jmath}_0(t), \hat{z}_0(t)\right)\,dt = \widehat{\mathcal{I}}_0(\hat{\rho}_0,\hat{\jmath}_0),
	\end{align*}
	which finishes the proof of the lower bound.
\end{proof}
\subsection{Proof of upper bound}
We first comment on the idea of proof of Theorem~\ref{GF_NGF:thm:upper-bound}. To that end, let
\[
E:= (0,T)\times\mathbb{R}\quad\text{and}\quad E_0:=(0,T)\times(-1/2,+1/2).
\]
If the limiting rate function is finite, then by definition the pair~$(\hat{\rho}_0,\hat{\jmath}_0)$ is given by
\begin{align}\label{eq:GF_NGF:limsup:rho_0}
	\hat{\rho}_0(t,dy) &= \hat{z}_0(t)\delta_{-1/2}(dy) + (1-\hat{z}_0(t))\delta_{+1/2}(dy),\\
	\label{eq:GF_NGF:limsup:j_0}
	\hat{\jmath}_0(t,dy) &= \hat{\jmath}_0(t) \mathbf{1}_{(-1/2,+1/2)}(y)\,dy,
\end{align}
with~$\hat{z}_0\in H^1(0,T)$ and~$\hat{\jmath}_0(t)=-\partial_t\hat{z}_0(t) \geq 0$. We will first work under the following regularity assumption.
\begin{assumption}\label{assump:GF_NGF:z0-is-ct-and-pos}
	The density~$\hat{z}_0:[0,T]\to[0,1]$ satisfies
\begin{equation}
\partial_t\hat{z}_0\in C([0,T]),\; \inf_{t\in(0,T)}|\partial_t \hat{z}_0(t)|>0 \quad\text{and}\quad
\sup_{t\in(0,T)}|\partial_{tt}\hat{z}_0(t)|<\infty.
\end{equation}
\end{assumption}
The proof of Theorem~\ref{GF_NGF:thm:upper-bound} consists of the following four steps.
\begin{enumerate}
	\item We show that the limiting rate function satisfies
	\begin{equation}\label{eq:GF_NGF:limit-RF-via-density}
		\widehat{\mathcal{I}}_0(\hat{\rho}_0,\hat{\jmath}_0) = 
		\frac{1}{4} \int_{E_0} \hat{b}_0^2\,\hat{u}_0\,dydt,
	\end{equation}
	where~$\hat{u}_0:E_0\to[0,\infty)$ is the function given by
	\begin{equation}\label{eq:GF_NGF:limit-density-u0}
		\hat{u}_0(t,y) = -\left(\hat{\jmath}_0-\hat{z}_0\right)\left(y+\hat{y}_t\right)\left(y-\frac{1}{2}\right),\quad \hat{y}_t := \frac{1}{2} \frac{\hat{\jmath}_0+\hat{z}_0}{\hat{\jmath}_0-\hat{z}_0},
	\end{equation}
	and~$\hat{b}_0:E_0\to\mathbb{R}$ is defined by
	\begin{equation}\label{eq:GF_NGF:limit-function-b0}
		\hat{b}_0(t,y) := \frac{\hat{\jmath}_0(t)+\partial_y\hat{u}_0(t,y)}{\hat{u}_0(t,y)} = \frac{2}{\hat{y}_t + y}.
	\end{equation}
	The second-order polynomial~$\hat{u}_0(t,\cdot)$ is either concave ($\hat{\jmath}_0>\hat{z}_0$), linear ($\hat{\jmath}_0=\hat{z}_0$) or convex ($\hat{\jmath}_0<\hat{z}_0$). These three cases are sketched in Figure~\ref{fig:u_concance_linear_convex}. 
	\begin{figure}[h!]
		\labellist
		\pinlabel $y$ at 950 0
		\pinlabel $\hat{z}_0(t)$ at 900 500
		\pinlabel $\hat{u}_0(t,y)$ at 550 750
		\endlabellist
		\centering
		\includegraphics[scale=.18]{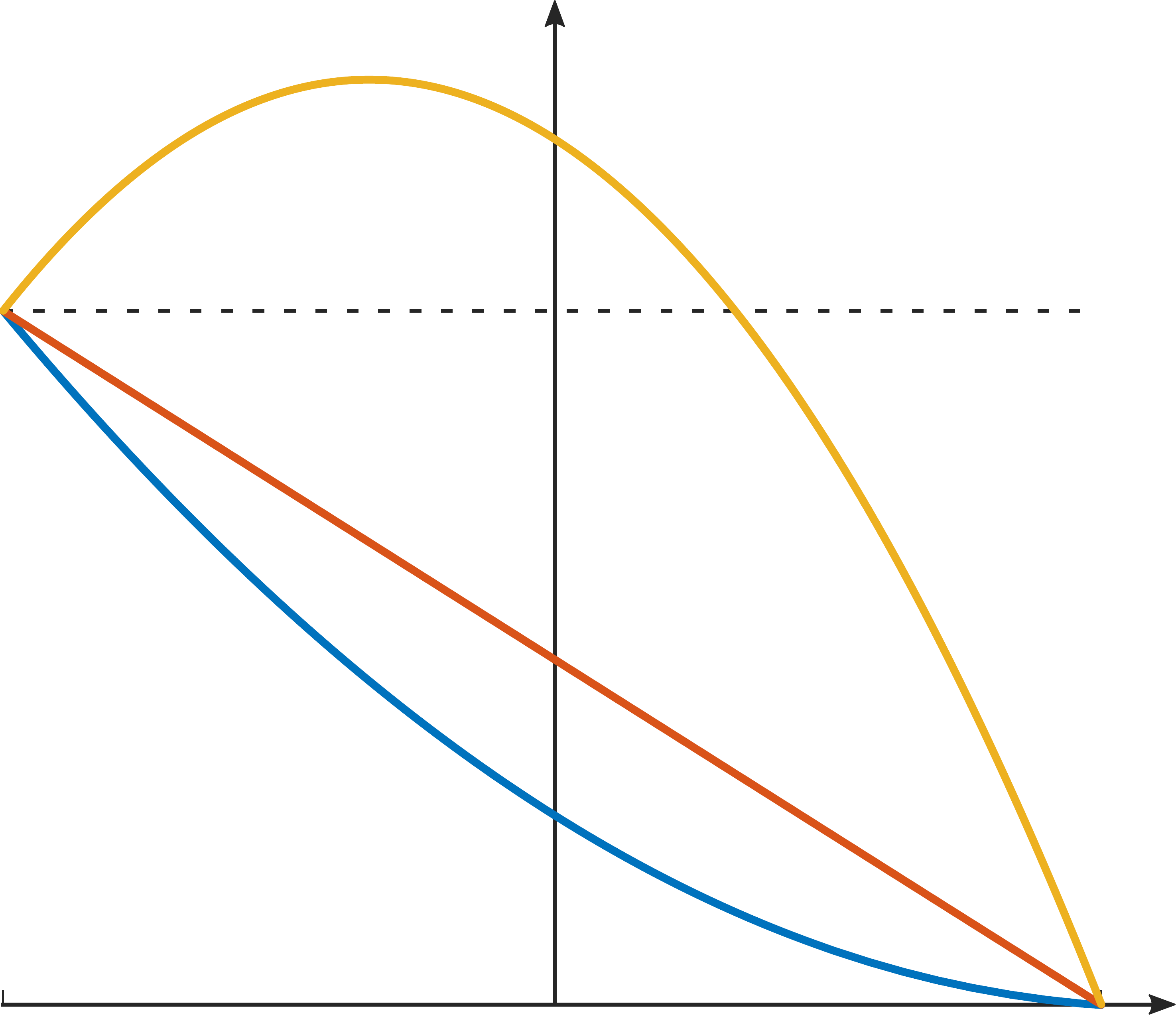}
		\caption{The polynomial $y\mapsto \hat{u}_0(t,y)$ on~$[-1/2,+1/2]$ for the three cases~$\hat{\jmath}_0(t)>\hat{z}_0(t)$ (yellow), $\hat{\jmath}_0(t)=\hat{z}_0(t)$ (red) and~$\hat{\jmath}_0(t)<\hat{z}_0(t)$ (blue). In particular, the function always satisfies $\hat{u}_0(t,-1/2) = \hat{z}_0(t)$ and~$\hat{u}_0(t,+1/2)=0$.}
		\label{fig:u_concance_linear_convex}
	\end{figure}
	\item We define the function~$\hat{u}_\varepsilon^\ell:E\to[0,\infty)$ as the weak solution to the auxiliary PDE 
	\begin{equation}\label{eq:GF_NGF:steps-auxiliary-PDE}
		\hat{g}_\varepsilon^\ell \partial_t \hat{u}_\varepsilon^\ell = \partial_{yy}\hat{u}_\varepsilon^\ell - \partial_y(\hat{b}_0\mathbf{1}_{E_0} \hat{u}_\varepsilon^\ell),
	\end{equation}
	where~$\hat{g}_\varepsilon^\ell:\mathbb{R}\to(0,\infty)$ denotes the Lebesgue density of the left-stationary measure~$\hat{\gamma}_\varepsilon^\ell$ from Definition~\ref{GF_NGF:def:transformed-stationary-measure}, that is~$\hat{\gamma}_\varepsilon^\ell(dy) = \hat{g}_\varepsilon^\ell(y)dy$.
	With that, we define the pair~$(\hat{\rho}_\varepsilon,\hat{\jmath}_\varepsilon)$ by setting 
	\begin{align*}
		\hat{\rho}_\varepsilon(t,dy) := \hat{u}_\varepsilon^\ell(t,y)\hat{\gamma}_\varepsilon^\ell(dy)\quad\text{and}\quad \hat{\jmath}_\varepsilon := -\partial_y\hat{u}_\varepsilon^\ell + \hat{b}_0\mathbf{1}_{E_0} \hat{u}_\varepsilon^\ell.
	\end{align*}
	We choose the initial condition~$\hat{u}_\varepsilon^\ell(0,\cdot)$ such that the measure
	\begin{align*}
		\hat{\rho}_\varepsilon(0,dy):=\hat{u}_\varepsilon^\ell(0,y)\hat{\gamma}_\varepsilon^\ell(dy)
	\end{align*}
	has mass one and converges weakly to~$\hat{\rho}_0(0,dy)$.
	\item We show that the solution~$\hat{u}_\varepsilon^\ell$ to the auxiliary PDE~\eqref{eq:GF_NGF:steps-auxiliary-PDE} is such that
	\begin{equation}\label{eq:GF_NGF:u_eps-converges-to-u-in-E0}
	\hat{u}_\varepsilon^\ell\mathbf{1}_{E_0} \xrightarrow{\varepsilon\to 0} \hat{u}_0 \quad\text{weakly in}\;L^2(E_0),
	\end{equation}
	and that the pair~$(\hat{\rho}_\varepsilon,\hat{\jmath}_\varepsilon)$ converges to~$(\hat{\rho}_0,\hat{\jmath}_0)$ in the sense of Definition~\ref{GF_NGF:def:converge-in-CE}.
	\item We verify that with the choice of~$(\hat{\rho}_\varepsilon,\hat{\jmath}_\varepsilon)$ as above, the rate function is given by
	\begin{equation}\label{eq:GF_NGF:RF-via-density}
	\widehat{\mathcal{I}}_\varepsilon(\hat{\rho}_\varepsilon,\hat{\jmath}_\varepsilon) = \frac{1}{4}\int_{E_0}\hat{b}_0^2 \hat{u}_\varepsilon^\ell\,dydt.
	\end{equation}
\end{enumerate}
With these steps accomplished, the limsup-bound follows as
\begin{align*}
	\lim_{\varepsilon\to 0}\widehat{\mathcal{I}}_\varepsilon(\hat{\rho}_\varepsilon,\hat{\jmath}_\varepsilon) &\overset{\eqref{eq:GF_NGF:RF-via-density}}{=} \lim_{\varepsilon\to 0} \frac{1}{4}\int_{E_0}\hat{b}_0^2 \hat{u}_\varepsilon^\ell\,dydt	\\
	&\overset{\eqref{eq:GF_NGF:u_eps-converges-to-u-in-E0}}{=} \frac{1}{4}\int_{E_0}\hat{b}_0^2\hat{u}_0\,dydt \\
	&\overset{\eqref{eq:GF_NGF:limit-RF-via-density}}{=} \widehat{\mathcal{I}}_0(\hat{\rho}_0,\hat{\jmath}_0).
\end{align*}
We now formulate the Lemmas we need in order to rigorously carry out the abovementioned steps. After that, we give the proof of Theorem~\ref{GF_NGF:thm:upper-bound} and then prove the Lemmas.
\begin{lemma}\label{GF_NGF:lemma:upperbound:energy-dense-set-of-densities}
	The set of densities~$\{\hat{z}\}$ satisfying Assumption~\ref{assump:GF_NGF:z0-is-ct-and-pos} is energy-dense; that means if~$\widehat{\mathcal{I}}_0(\hat{\rho}_0,\hat{\jmath}_0)$ is finite, then there are denities~$\hat{z}_0^\delta$ satisfying Assumption~\ref{assump:GF_NGF:z0-is-ct-and-pos} such that the pair~$(\hat{\rho}_0^\delta,\hat{\jmath}_0^\delta)$ defined via~$\hat{z}_0^\delta$ as in~\eqref{eq:GF_NGF:limsup:rho_0} and~\eqref{eq:GF_NGF:limsup:j_0} satisfies
	\begin{equation*}
	\lim_{\delta\to 0}(\hat{\rho}_0^\delta,\hat{\jmath}_0^\delta) = (\hat{\rho}_0,\hat{\jmath}_0) \quad\text{and}\quad \lim_{\delta\to 0}\widehat{\mathcal{I}}_0(\hat{\rho}_0^\delta,\hat{\jmath}_0^\delta) = \widehat{\mathcal{I}}_0(\hat{\rho}_0,\hat{\jmath}_0).
	\end{equation*}
\end{lemma}
\begin{lemma}[Limiting Rate Function]\label{lemma:GF_NGF:limsup:limiting-RF}
	The rate function~$\widehat{\mathcal{I}}_0$ is given by~\eqref{eq:GF_NGF:limit-RF-via-density}.
\end{lemma}
In the next Lemmas,~$\hat{\gamma}_\varepsilon^\ell(dy) = \hat{g}_\varepsilon^\ell(y)dy$ is the transformed left-normalized stationary measure from Definition~\ref{GF_NGF:def:transformed-stationary-measure}.
\begin{lemma}[Auxiliary PDE]\label{lemma:GF_NGF:limsup:auxiliary-PDE}
	For any~$\varepsilon>0$ and any initial condition~$\hat{u}_\varepsilon^\ell(0,\cdot)\in C(\mathbb{R};[0,\infty))$, under Assumption~\ref{assump:GF_NGF:z0-is-ct-and-pos} there exists a weak solution~$\hat{u}_\varepsilon^\ell$ to the PDE~\eqref{eq:GF_NGF:steps-auxiliary-PDE}; that means there exists a function~$\hat{u}_\varepsilon^\ell:E\to[0,\infty)$ satisfying
	\begin{equation*}
		\hat{u}_\varepsilon^\ell \in L^2(0,T;H^1(\mathbb{R}))\cap C(0,T;L^2(\mathbb{R})) \quad\text{and}\quad \partial_t \hat{u}_\varepsilon^\ell \in L^2(0,T;H^{-1}(\mathbb{R}))
	\end{equation*}
	such that for any~$\varphi\in H^1(\mathbb{R})$,
	\begin{equation}
	\int_E \varphi \, \hat{g}_\varepsilon^\ell \partial_t\hat{u}_\varepsilon^\ell\,dydt = -\int_E \partial_y \hat{u}_\varepsilon^\ell \, \partial_{y}\varphi \,dydt + \int_E \hat{b}_0\mathbf{1}_{E_0}\hat{u}_\varepsilon^\ell \, \partial_y \varphi\,dydt.
	\end{equation}
\end{lemma}
\begin{lemma}[Uniform Energy Estimates]\label{lemma:GF_NGF:unif-energy-est}
	For~$\varepsilon>0$, let~$\hat{u}_\varepsilon^\ell$ be the solution to the auxiliary PDE~\eqref{eq:GF_NGF:steps-auxiliary-PDE} under Assumption~\ref{assump:GF_NGF:z0-is-ct-and-pos} and with initial condition~$\hat{u}_\varepsilon^\ell(0,\cdot)$ such that $\hat{\rho}_\varepsilon(0,dy) = \hat{u}_\varepsilon^\ell(0,y)\hat{\gamma}_\varepsilon^\ell(dy)$ has mass one and~$\|\hat{u}_\varepsilon^\ell(0,\cdot)\|_{L^\infty(\mathbb{R})}$ is uniformly bounded in~$\varepsilon>0$. 
	
	Let~$B:E\to[0,\infty)$ be the function defined by~$B(t,y) := \int_{-1/2}^{y}\hat{b}_0(t,z)\mathbf{1}_{E_0}(t,z)\,dz$, and  define~$\hat{v}_\varepsilon^\ell := e^{-B}\hat{u}_\varepsilon^\ell$ and~$\Omega := (-17,+17)$. 
	
	Then there exists a constant~$C>0$ such that for every~$\varepsilon>0$,
	\begin{align}
	&\int_0^T\int_\mathbb{R} e^{B} |\partial_y\hat{v}_\varepsilon^\ell|^2\,dydt + \sup_{t\in(0,T)}\int_\mathbb{R} e^{B}\hat{g}_\varepsilon^\ell \hat{v}_\varepsilon^\ell(t)^2\,dy \leq C,\label{eq:GF_NGF:limsup:energy-bound-v_eps}\\
	& \int_0^T\|\hat{u}_\varepsilon^\ell(t,\cdot)\|_{H^1(\Omega)}^2\,dt \leq C.\label{eq:GF_NGF:limsup:u_eps-bounded-in-H1-Omega}
	\end{align}
\end{lemma}
The size of~$\Omega$ can be chosen arbitrary, as long as it is finite and contains the inverval~$[-1/2,+1/2]$ in its interior.
\begin{lemma}[Limiting density]\label{lemma:GF_NGF:limsup:limit-density}
	Let~$(\hat{\rho}_0,\hat{\jmath}_0)$ be a pair given by~\eqref{eq:GF_NGF:limsup:rho_0} and~\eqref{eq:GF_NGF:limsup:j_0}. Let~$\hat{u}_\varepsilon^\ell$ be the solution to the PDE~\eqref{eq:GF_NGF:steps-auxiliary-PDE} with initial condition~$\hat{u}_\varepsilon(0,\cdot)\in C^\infty(\mathbb{R})$ such that $\hat{u}_\varepsilon^\ell(0,y)\hat{\gamma}_\varepsilon^\ell(dy)\xrightharpoonup{\ast} \hat{\rho}_0(0,dy)$ as~$\varepsilon\to 0$. Define
	\begin{equation}\label{eq:GF_NGF:limsup:rho_eps-and-j_eps}
		\hat{\rho}_\varepsilon(t,dy) := \hat{u}_\varepsilon^\ell(t,y)\hat{\gamma}_\varepsilon^\ell(dy)\quad\text{and}\quad \hat{\jmath}_\varepsilon := -\partial_y\hat{u}_\varepsilon^\ell + \hat{b}_0\mathbf{1}_{E_0} \hat{u}_\varepsilon^\ell,
	\end{equation}
	and let~$\Omega:=(-17,+17)$. Then we can choose a subsequence of~$(\hat{\rho}_\varepsilon,\hat{\jmath}_\varepsilon)$ (denoted the same) such that:
	\begin{enumerate}[label=(\roman*)]
		\item \label{item:GF_NGF:lemma-limit-density:i}$(\hat{\rho}_\varepsilon,\hat{\jmath}_\varepsilon)\in \mathrm{CE}(0,T;\mathbb{R})$ and~$\widehat{\mathcal{I}}_\varepsilon(\hat{\rho}_\varepsilon,\hat{\jmath}_\varepsilon)$ is given by~\eqref{eq:GF_NGF:RF-via-density}.
		\item \label{item:GF_NGF:lemma-limit-density:rho-j-converge} The pair~$(\hat{\rho}_\varepsilon,\hat{\jmath}_\varepsilon)$ converges to~$(\hat{\rho}_0,\hat{\jmath}_0)$ in the sense of Definition~\ref{GF_NGF:def:converge-in-CE}.
		\item \label{item:GF_NGF:lemma-limit-density:ii} There exists a function~$\hat{u}_0^\ell\in L^2(0,T;H^1(\Omega))$ such that $\hat{u}_\varepsilon^\ell$ converges to $\hat{u}_0^\ell$ weakly in~$L^2(0,T;L^2(\Omega))$, and~$\hat{u}_0^\ell \mathbf{1}_{E_0}=\hat{u}_0$.
	\end{enumerate}
\end{lemma} 
\begin{proof}[Proof of Theorem~\ref{GF_NGF:thm:upper-bound}]
	Let~$(\hat{\rho}_0,\hat{\jmath}_0)$ be given by~\eqref{eq:GF_NGF:limsup:rho_0} and~\eqref{eq:GF_NGF:limsup:j_0}. By Lemma~\ref{GF_NGF:lemma:upperbound:energy-dense-set-of-densities}, we can suppose without loss of generality that~$\hat{z}_0$ satisfies Assumption~\ref{assump:GF_NGF:z0-is-ct-and-pos}.
	\smallskip
	
	By Lemma~\ref{lemma:GF_NGF:limsup:auxiliary-PDE}, we can define the function~$\hat{u}_\varepsilon^\ell:E\to[0,\infty)$ as the weak solution to the PDE~\eqref{eq:GF_NGF:steps-auxiliary-PDE}. We take the initial condition~$\hat{u}_\varepsilon^\ell(0,\cdot)\in C(\mathbb{R})$ to be such that~$\hat{\rho}_\varepsilon(0,dy)\xrightharpoonup{\ast}\hat{\rho}_0(0,dy)$ as~$\varepsilon\to 0$ and define the measures $(\hat{\rho}_\varepsilon,\hat{\jmath}_\varepsilon)$ by~\eqref{eq:GF_NGF:limsup:rho_eps-and-j_eps}.
	\smallskip
	
	By~\ref{item:GF_NGF:lemma-limit-density:rho-j-converge} of Lemma~\ref{lemma:GF_NGF:limsup:limit-density}, $\lim_{\varepsilon\to 0}(\hat{\rho}_\varepsilon,\hat{\jmath}_\varepsilon)=(\hat{\rho}_0,\hat{\jmath}_0)$ in the sense of Definition~\ref{GF_NGF:def:converge-in-CE}.
	We are left with verifying the limsup bound. By~\ref{item:GF_NGF:lemma-limit-density:i} of Lemma~\ref{lemma:GF_NGF:limsup:limit-density},
	\begin{equation*}
	\widehat{\mathcal{I}}_\varepsilon(\hat{\rho}_\varepsilon,\hat{\jmath}_\varepsilon) = \frac{1}{4}\int_{E_0}\hat{b}_0^2 \hat{u}_\varepsilon^\ell\,dydt.	
	\end{equation*}
	As a consequence of Assumption~\ref{assump:GF_NGF:z0-is-ct-and-pos}, the function~$\hat{b}_0^2$ is in~$L^2(E_0)$, and by~\ref{item:GF_NGF:lemma-limit-density:ii} of Lemma~\ref{lemma:GF_NGF:limsup:limit-density}, the function~$\hat{u}_\varepsilon^\ell$ converges to~$\hat{u}_0$ weakly in~$L^2(E_0)$. Therefore
	\begin{equation*}
	\lim_{\varepsilon\to 0} \frac{1}{4}\int_{E_0}\hat{b}_0^2 \hat{u}_\varepsilon^\ell\,dydt = \frac{1}{4}\int_{E_0}\hat{b}_0^2 \hat{u}_0\,dydt,
	\end{equation*}
	and by Lemma~\ref{lemma:GF_NGF:limsup:limiting-RF}, 
	\begin{equation*}
	\frac{1}{4}\int_{E_0}\hat{b}_0^2 \hat{u}_0\,dydt = \widehat{\mathcal{I}}_0(\hat{\rho}_0,\hat{\jmath}_0).
	\end{equation*}
	Thus $\lim_{\varepsilon\to 0}\widehat{\mathcal{I}}_\varepsilon(\hat{\rho}_\varepsilon,\hat{\jmath}_\varepsilon) = \widehat{\mathcal{I}}_0(\hat{\rho}_0,\hat{\jmath}_0)$, which proves the limsup bound.
\end{proof}
\subsubsection{Proofs of the Lemmas}
\begin{proof}[Proof of Lemma~\ref{GF_NGF:lemma:upperbound:energy-dense-set-of-densities}]
	One can follow the same line of argument of~\cite[Theorem~6.1]{ArnrichMielkePeletierSavareVeneroni2012}.
\end{proof}
\begin{proof}[Proof of Lemma~\ref{lemma:GF_NGF:limsup:limiting-RF}] By definition, for a pair~$(\hat{\rho}_0,\hat{\jmath}_0)$ satisfying~\eqref{eq:GF_NGF:limsup:rho_0} and~\eqref{eq:GF_NGF:limsup:j_0},
	\begin{equation*}
	\widehat{\mathcal{I}}_0(\hat{\rho}_0,\hat{\jmath}_0) = \int_0^T S(\hat{\jmath}_0(t)|\hat{z}_0(t))\,dt,	
	\end{equation*}
We show in~\ref{item:GF_NGF:lemma:var-probl:i} of Lemma~\ref{lemma:GF_NGF:variational-problem} that the function~$S$ arises as the solution to the following variational problem:
	\begin{equation*}
		S(\hat{\jmath}_0(t)|\hat{z}_0(t)) = \frac{1}{4} \inf_u \int_{-1/2}^{+1/2} \frac{1}{u(y)}\big|\hat{\jmath}_0(t) + \partial_y u(y)\big|^2\,dy,
	\end{equation*} 
	where the infimum is taken over smooth functions~$u$ satisfying the boundary conditions $u(-1/2)=\hat{z}_0(t)$ and~$u(+1/2)=0$. By~\ref{item:GF_NGF:lemma:var-probl:ii} of Lemma~\ref{lemma:GF_NGF:variational-problem}, the optimizer of this variational problem is the polynomial $y\mapsto \hat{u}(t,y)$ given above in~\eqref{eq:GF_NGF:limit-density-u0}. Hence we find that
	\begin{align*}
		\widehat{\mathcal{I}}_0(\hat{\rho}_0,\hat{\jmath}_0) &\overset{\mathrm{def}}{=} \int_0^T S(\hat{\jmath}_0(t)|\hat{z}_0(t))\,dt\\
		&= \frac{1}{4} \int_0^T\int_{-1/2}^{+1/2} \frac{1}{\hat{u}_0(t,y)}\left|\hat{\jmath}_0(t) + \partial_y \hat{u}_0(t,y)\right|^2\,dydt \overset{\mathrm{def}}{=} \frac{1}{4}\int_{E_0}\hat{b}_0^2 \hat{u}_0\,dydt,
	\end{align*}
	where the last equality follows from the Definition of~$\hat{b}_0$ in~\eqref{eq:GF_NGF:limit-function-b0}.
\end{proof}
\begin{proof}[Proof of Lemma~\ref{lemma:GF_NGF:limsup:auxiliary-PDE}]
	This follows from the fact that~$\hat{b}_0$ is bounded.
\end{proof}
\begin{proof}[Proof of Lemma~\ref{lemma:GF_NGF:unif-energy-est}]
	We first prove the estimate~\eqref{eq:GF_NGF:limsup:energy-bound-v_eps}. As a consequence of Assumption~\ref{assump:GF_NGF:z0-is-ct-and-pos}, the function~$B$ is bounded. We find by calculation that the function $\hat{v}_\varepsilon^\ell = e^{-B}\hat{u}_\varepsilon^\ell$ is the weak solution to
	\begin{equation}\label{eq:GF_NGF:limsup:eq-for-v_eps}
		\hat{g}_\varepsilon^\ell \partial_t \left(e^{B}\hat{v}_\varepsilon^\ell\right) = \partial_y\left(e^{B}\partial_y \hat{v}_\varepsilon^\ell\right).
	\end{equation}
	Multiplying~\eqref{eq:GF_NGF:limsup:eq-for-v_eps} with $\hat{v}_\varepsilon^\ell$ and integrating over~$\mathbb{R}$, (that is specialising the test function to the weak solution~$\hat{v}_\varepsilon^\ell$), we find
	\begin{equation*}
		\frac{1}{2}\frac{d}{dt}\left(\int_{\mathbb{R}}\hat{g}_\varepsilon^\ell e^B \frac{1}{2}(\hat{v}_\varepsilon^\ell)^2\,dy\right) + \frac{1}{2}\int_\mathbb{R}\hat{g}_\varepsilon^\ell (\hat{v}_\varepsilon^\ell)^2 \,\partial_t e^{B}\,dy= -\int_\mathbb{R} e^{B}|\partial_y\hat{v}_\varepsilon^\ell|^2\,dy.
	\end{equation*}
	Integrating over the time inverval~$(0,t)$ for~$t\in (0,T)$,
	\begin{multline}\label{eq:GF_NGF:proof-lemma-energy-dens:after-mult-v}
		\frac{1}{2}\int_\mathbb{R} e^B \hat{g}_\varepsilon^\ell \frac{1}{2} \hat{v}_\varepsilon^\ell(t,y)^2\,dy + \frac{1}{2} \int_0^t\int_\mathbb{R} \left(\partial_t e^B\right) \hat{g}_\varepsilon^\ell (\hat{v}_\varepsilon^\ell)^2\,dydt' + \int_0^t\int_\mathbb{R} e^B |\partial_y\hat{v}_\varepsilon^\ell|^2\,dydt \\
		= \frac{1}{2}\int_\mathbb{R} e^{B(0,y)} \hat{g}_\varepsilon^\ell \frac{1}{2}\hat{v}_\varepsilon^\ell(0,y)^2\,dy
	\end{multline}
	Applying the estimate $ac\geq -(1/2)a^2-(1/2)c^2$ to the second term with
	\begin{equation*}
		a = \sqrt{\hat{g}_\varepsilon^\ell} \hat{v}_\varepsilon^\ell e^{B/2} \partial_t B \quad \text{and}\quad c = \sqrt{\hat{g}_\varepsilon^\ell} \hat{v}_\varepsilon^\ell e^{B/2}
	\end{equation*}
	leads to the estimate
	\begin{multline*}
	\int_0^t\int_\mathbb{R} \left(\partial_t e^B\right) \hat{g}_\varepsilon^\ell (\hat{v}_\varepsilon^\ell)^2\,dydt \\ \geq -\frac{1}{2}\int_0^t\int_\mathbb{R} \left(\partial_t B\right)^2 e^B \hat{g}_\varepsilon^\ell (\hat{v}_\varepsilon^\ell)^2\,dydt'- \frac{1}{2} \int_0^t\int_\mathbb{R}e^{B}\hat{g}_\varepsilon^\ell (\hat{v}_\varepsilon^\ell)^2\,dydt'	\\
	\geq (-t)\left[\frac{1}{2}\sup_{E}|\partial_t B|^2 + \frac{1}{2} \right] \sup_{t'\in(0,T)} \int_\mathbb{R} e^{B(t',y)}\hat{g}_\varepsilon^\ell(y) \hat{v}_\varepsilon^\ell(t',y)^2\,dy.
	\end{multline*}
	With the abbreviations
	\begin{equation*}
		C := \frac{1}{2}\sup_{E}|\partial_t B|^2 + \frac{1}{2} \quad\text{and}\quad C' := \frac{1}{2}\int_\mathbb{R}e^{B(0,y)} \hat{g}_\varepsilon^\ell \frac{1}{2}\hat{v}_\varepsilon^\ell(0,y)^2\,dy,
	\end{equation*}
	we find via the above estimate and by~\eqref{eq:GF_NGF:proof-lemma-energy-dens:after-mult-v} that
	\begin{equation*}
		(1-C\cdot t) \sup_{t'\in(0,T)} \int_\mathbb{R} e^{B(t',y)}\hat{g}_\varepsilon^\ell(y) \hat{v}_\varepsilon^\ell(t',y)^2\,dy + \int_0^t\int_\mathbb{R} e^{B}|\partial_y\hat{v}_\varepsilon^\ell|^2\,dydt \leq C'.
	\end{equation*}
	Now iterating this estimate,~\eqref{eq:GF_NGF:limsup:energy-bound-v_eps} follows.
	\smallskip
	
	We now prove~\eqref{eq:GF_NGF:limsup:u_eps-bounded-in-H1-Omega}, that is boundedness of~$\hat{u}_\varepsilon^\ell$ in~$L^2(0,T;H^1(\Omega))$. Since both~$B$ and~$\partial_y B$ are bounded functions as a consequence of Assumption~\ref{assump:GF_NGF:z0-is-ct-and-pos}, proving boundedness of~$\hat{v}_\varepsilon^\ell=e^{-B}\hat{u}_\varepsilon^\ell$ in~$L^2(0,T;H^1(\Omega))$ is sufficient for proving~\eqref{eq:GF_NGF:limsup:u_eps-bounded-in-H1-Omega}. To prove boundedness of~$\hat{v}_\varepsilon^\ell$ in~$L^2(0,T;H^1(\Omega))$, we will establish the following estimates:
	\begin{enumerate}[label=(\arabic*)]
		\item There exists a constant~$C_1>0$ such that for all~$\varepsilon>0$,
		\begin{equation}\label{eq:GF_NGF:limsup:seminorm-bounded}
		p(\hat{v}_\varepsilon^\ell)\leq C_1,
		\end{equation}
		where $p$ is the seminorm~$p(f):= \left(\int_0^T|f(t,+1/2)|^2\,dt\right)^{-1/2}$.
		\item There exists a constant~$C_2>0$ such that for all~$\varepsilon>0$,
		\begin{equation}\label{eq:GF_NGF:limsup:poincare}
		\int_0^T\|\hat{v}_\varepsilon^\ell(t,\cdot)\|_{L^2(\Omega)}^2\,dt \leq C_2 \left[\int_0^T\|\partial_y\hat{v}_\varepsilon^\ell(t,\cdot)\|_{L^2(\Omega)}^2\,dt + p(\hat{v}_\varepsilon^\ell)\right].
		\end{equation}
	\end{enumerate}
	With these estimates established, boundedness of~$\hat{v}_\varepsilon^\ell$ and hence~\ref{eq:GF_NGF:limsup:u_eps-bounded-in-H1-Omega} follow as
	\begin{align*}
		\int_0^T\|\hat{v}\|_{L^2(\Omega)}^2\,dt + \int_0^T\|\partial_y \hat{v}_\varepsilon^\ell\|_{L^2(\Omega)}^2\,dt &\overset{\eqref{eq:GF_NGF:limsup:poincare}}{\leq}(1+C_2)\int_0^T \|\partial_y \hat{v}_\varepsilon^\ell\|_{L^2(\Omega)}^2\,dt + C_2p(\hat{v}_\varepsilon^\ell)\\
		&\overset{\eqref{eq:GF_NGF:limsup:seminorm-bounded}}{\leq} (1+C_2)\int_0^T \|\partial_y \hat{v}_\varepsilon^\ell\|_{L^2(\Omega)}^2\,dt + C_2 C_1\\
		&\overset{\eqref{eq:GF_NGF:limsup:energy-bound-v_eps}}{\leq} (1+C_2) \cdot \sup_E (e^{-B}) \cdot C + C_2 C_1 < \infty.
	\end{align*}
	The estimate~\eqref{eq:GF_NGF:limsup:poincare} is the Poincaré inequality~\cite[Eq.~1.35]{Temam2012}, which holds since~$p$ is a seminorm that is a norm when restricted to constant functions, that is
	\begin{equation*}
	(p(c)=0,c\in\mathbb{R}) \Rightarrow c= 0.	
	\end{equation*}
	For verifying~\eqref{eq:GF_NGF:limsup:seminorm-bounded}, we prove that
	\begin{equation*}
	\lim_{\varepsilon\to 0}p(\hat{v}_\varepsilon^\ell) = 0.
	\end{equation*}
	Since there exists a constant~$C>0$ such that for all~$\varepsilon>0$,
	\begin{equation*}
		\int_0^T\int_{\Omega} |\partial_y\hat{v}_\varepsilon^\ell|^2\,dydt \leq \sup_E(e^{-B}) \int_0^T\int_\mathbb{R}|\partial_y\hat{v}_\varepsilon^\ell|^2\,dydt \overset{\eqref{eq:GF_NGF:limsup:energy-bound-v_eps}}{\leq} C,
	\end{equation*}
	there are a function~$\alpha_\varepsilon\in L^2(0,T)$ and a constant~$C'>0$ such that for all~$\varepsilon>0$,
	\begin{equation*}
		\|\alpha_\varepsilon\|_{L^2(0,T)}^2\leq C'\quad\text{and}\quad \hat{v}_\varepsilon^\ell(t,y) \geq \hat{v}_\varepsilon^\ell(t,1/2) - \alpha_\varepsilon(t)|y-(1/2)|^{1/2}.
	\end{equation*}
	Let~$\delta>0$ be arbitrary and set~$U_\delta:=(\frac{1}{2}-\delta,\frac{1}{2}+\delta)$. Then since~$\hat{u}_\varepsilon^\ell(t,y)\hat{\gamma}_\varepsilon^\ell(dy)$ has mass one for every~$t>0$ and~$\hat{\gamma}_\varepsilon^\ell(dy)=\hat{g}_\varepsilon^\ell(y)dy$, 
	\begin{align*}
		T &= \int_0^T\int_\mathbb{R} e^B e^{-B}\hat{u}_\varepsilon^\ell \hat{g}_\varepsilon^\ell\,dydt \\
		&\geq \inf_E(e^B)\int_0^T\int_{U_\delta}\left(\hat{v}_\varepsilon^\ell(t,1/2) - \alpha_\varepsilon(t)|y-(1/2)|^{1/2}\right)\hat{g}_\varepsilon^\ell\,dydt\\
		&\geq \inf_E(e^B)\int_0^T\int_{U_\delta}\hat{v}_\varepsilon^\ell(t,1/2)\hat{g}_\varepsilon^\ell\,dydt - \inf_E(e^B)\delta^{1/2}\int_0^T\int_{U_\delta} \alpha_\varepsilon(t)\hat{g}_\varepsilon^\ell\,dydt.
	\end{align*}
	By re-organizing, we deduce the estimate
	\begin{align*}
		\int_0^T\hat{v}_\varepsilon^\ell(t,1/2)\,dt \int_{U_\delta}\hat{g}_\varepsilon^\ell\,dy\leq C \left( 1 + \delta^{1/2} \int_{U_\delta}\hat{g}_\varepsilon^\ell\,dy\right),
	\end{align*}
	where~$C = \max(T \cdot (\inf_E e^B)^{-1}, C' T^{1/2})$, and therefore arrive at
	\begin{equation*}
		\int_0^T\hat{v}_\varepsilon^\ell(t,1/2)\,dt \leq C\delta^{1/2} + C \left(\int_{U_\delta}\hat{g}_\varepsilon^\ell\,dy\right)^{-1}.
	\end{equation*}
	For any~$\delta>0$, we have~$\int_{U_\delta}\hat{g}_\varepsilon^\ell\,dy\to+\infty$ as $\varepsilon\to 0$, so that
	\begin{equation*}
		\limsup_{\varepsilon\to 0}\int_0^T\hat{v}_\varepsilon^\ell(t,1/2)\,dt \leq C \delta^{1/2}.
	\end{equation*}
	Since~$\delta>0$ is arbitrary, this implies~$p(\hat{v}_\varepsilon^\ell)\to 0$ as~$\varepsilon\to 0$.
\end{proof}
\begin{proof}[Proof of Lemma~\ref{lemma:GF_NGF:limsup:limit-density}]
	\ref{item:GF_NGF:lemma-limit-density:i}: The fact that~$(\hat{\rho}_\varepsilon,\hat{\jmath}_\varepsilon)\in\mathrm{CE}(0,T;\mathbb{R})$ follows directly from~\ref{eq:GF_NGF:limsup:rho_eps-and-j_eps} and the definition of~$\hat{u}_\varepsilon^\ell$. The rate function is given by
	\begin{align*}
		\widehat{\mathcal{I}}_\varepsilon(\hat{\rho}_\varepsilon,\hat{\jmath}_\varepsilon) &\overset{\mathrm{def}}{=} \frac{1}{2} \sup_{b\in C_c^\infty(E)}\int_E \left[\hat{u}_\varepsilon^\ell\left(-\partial_y b - \frac{1}{2}b^2\right) + \hat{\jmath}_\varepsilon\cdot b\right]\,dydt \\
		&= \frac{1}{4}\int_E \frac{1}{\hat{u}_\varepsilon^\ell(t,y)}\big|\hat{\jmath}_\varepsilon(t,y) + \partial_y \hat{u}_\varepsilon^\ell(t,y)\big|^2\,dydt,
	\end{align*}
	where the last equality follows from Lemma~\ref{GF_NGF:lemma:appendix:dual-of-convex-functions}. The form of the rate function~\eqref{eq:GF_NGF:RF-via-density} is an immediate consequence of the definition of~$\hat{\jmath}_\varepsilon$ in~\eqref{eq:GF_NGF:limsup:rho_eps-and-j_eps}.
	\smallskip
	
	We prove~\ref{item:GF_NGF:lemma-limit-density:rho-j-converge} and~\ref{item:GF_NGF:lemma-limit-density:ii} via the following steps:
	\begin{enumerate}[label=(\arabic*)]
		\item \label{item:GF_NGF:proof:rho-j-conv:1} We show that there is a pair~$(\hat{\rho}_0^\ell,\hat{\jmath}_0^\ell)\in\mathrm{CE}(0,T;\mathbb{R})$ of the form
		\begin{align*}
			\hat{\rho}_0^\ell(t,dy) &= \hat{z}_0^\ell(t)\delta_{-1/2}(dy) + (1-\hat{z}_0^\ell(t))\delta_{+1/2}(dy),\\
			\hat{\jmath}_0^\ell(t,y)&= -\partial_t\hat{z}_0^\ell(t)\mathbf{1}_{(-1/2,+1/2)}(y),
		\end{align*}
		such that~$(\hat{\rho}_\varepsilon,\hat{\jmath}_\varepsilon)\to(\hat{\rho}_0^\ell,\hat{\jmath}_0^\ell)$ as~$\varepsilon\to 0$. We use the superscript~$\ell$ to distinguish this limit~$(\hat{\rho}_0^\ell,\hat{\jmath}_0^\ell)$ from~$(\hat{\rho}_0,\hat{\jmath}_0)$.
		\item \label{item:GF_NGF:proof:rho-j-conv:2} We show that the density~$\hat{u}_\varepsilon^\ell$ converges to a function~$\hat{u}_0^\ell$ in $L^2(0,T;H^1(\Omega))$ that satisfies for a.e. $t\in(0,T)$ the ODE
		\begin{align*}
			\begin{cases}
				\displaystyle-\partial_y\hat{u}_0^\ell(t,y) + \hat{b}_0(t,y) \hat{u}_0^\ell(t,y) = -\partial_t\hat{z}_0^\ell(t) \quad \text{in}\; \left(-\frac{1}{2},+\frac{1}{2}\right),\\
				\displaystyle\hat{u}_0^\ell(t,-1/2) = \hat{z}_0^\ell(t),\\
				\displaystyle\hat{u}_0^\ell(t,+1/2) = 0.
			\end{cases}
		\end{align*}
		\item \label{item:GF_NGF:proof:rho-j-conv:3} We show that the ODE enforces both~$\hat{z}_0^\ell(t) = \hat{z}_0(t)$ for a.e.~$t\in(0,T)$ and~$\hat{u}_0^\ell\mathbf{1}_{E_0}=\hat{u}_0$ in~$L^2(E_0)$.
	\end{enumerate}
	Then the convergence statement~\ref{item:GF_NGF:lemma-limit-density:rho-j-converge} follows as
	\begin{equation*}
		\lim_{\varepsilon\to 0}\left(\hat{\rho}_\varepsilon,\hat{\jmath}_\varepsilon\right)\overset{\ref{item:GF_NGF:proof:rho-j-conv:1}}{=} (\hat{\rho}_0^\ell,\hat{\jmath}_0^\ell) \overset{\ref{item:GF_NGF:proof:rho-j-conv:3}}{=}\left(\hat{\rho}_0,\hat{\jmath}_0\right),
	\end{equation*}
	and likewise~\ref{item:GF_NGF:lemma-limit-density:ii} as
	\begin{equation*}
	\lim_{\varepsilon\to 0}\hat{u}_\varepsilon^\ell\mathbf{1}_{E_0} \overset{\ref{item:GF_NGF:proof:rho-j-conv:2}}{=} \hat{u}_0^\ell \overset{\ref{item:GF_NGF:proof:rho-j-conv:3}}{=} \hat{u}_0.
	\end{equation*}
	We are left with verifying~\ref{item:GF_NGF:proof:rho-j-conv:1},~\ref{item:GF_NGF:proof:rho-j-conv:2} and~\ref{item:GF_NGF:proof:rho-j-conv:3}.
	\smallskip
	
	\ref{item:GF_NGF:proof:rho-j-conv:1}: For any test function~$\varphi\in C_b(E)$,
	\begin{align*}
		\left|\int_E \varphi \hat{\rho}_\varepsilon\right|^2 &\overset{\mathrm{def}}{=} \left|\int_E\varphi\, e^{B}\hat{g}_\varepsilon^\ell \hat{v}_\varepsilon^\ell\,dydt\right|^2\\
		&\overset{\mathrm{CS}}{\leq} \left(\int_E|e^{B}\hat{g}_\varepsilon^\ell \hat{v}_\varepsilon^\ell(t)^2|\,dydt\right) \left(\int_E |\varphi e^{B}\hat{g}_\varepsilon^\ell| \,dydt\right)\\
		&\leq C \left(\sup_t \int_\mathbb{R} e^{B(t)} \hat{g}_\varepsilon^\ell \hat{v}_\varepsilon^\ell(t)^2\,dy\right) \int_{\text{supp}(\varphi)}\hat{g}_\varepsilon^\ell(y)\,dy\\
		&\overset{\eqref{eq:GF_NGF:limsup:energy-bound-v_eps}}{\leq} C' \int_{\text{supp}(\varphi)}\hat{g}_\varepsilon^\ell(y)\,dy.
	\end{align*}
	Hence for any testfunction with support outside of~$\{\pm 1/2\}$,
	\begin{equation*}
		\int_E \varphi \hat{\rho}_\varepsilon \xrightarrow{\varepsilon\to 0} 0.
	\end{equation*}
	Therefore in the limit~$\varepsilon\to 0$, the family of measures~$\hat{\rho}_\varepsilon$ converges weakly to a measure~$\hat{\rho}_0^\ell$ that is concentrated on~$[0,T]\times\{\pm 1/2\}$.
	\smallskip
	
	The flux~$\hat{\jmath}_\varepsilon$ is given by~$\hat{\jmath}_\varepsilon = e^{-B}\partial_y\hat{v}_\varepsilon^\ell$. Since the function~$B$ is bounded, we find by virtue of the estimate~\eqref{eq:GF_NGF:limsup:energy-bound-v_eps} that~$\hat{\jmath}_\varepsilon$ is bounded in~$L^2(E)$, because
	\begin{equation*}
		\int_E |\hat{\jmath}_\varepsilon|^2\,dydt \leq C \int_E e^B |\partial_y\hat{v}_\varepsilon^\ell|^2\,dydt \overset{\eqref{eq:GF_NGF:limsup:energy-bound-v_eps}}{\leq} C'.
	\end{equation*}
	Hence the flux converges weakly in~$L^2(E)$ along a subsequence (denoted the same) to some~$\hat{\jmath}_0^\ell \in L^2(E)$. This finishes the proof of~$(\hat{\rho}_\varepsilon,\hat{\jmath}_\varepsilon)\to(\hat{\rho}_0^\ell,\hat{\jmath}_0^\ell)$, since weak $L^2$ convergence is stronger than convergence in distribution in the sense of Definition~\ref{GF_NGF:def:converge-in-CE}.
	\smallskip
	
	Combining the above convergence statements of~$\hat{\rho}_\varepsilon$ and~$\hat{\jmath}_\varepsilon$, we find for any test function~$\varphi\in C^\infty_c(E)$,
	\begin{equation*}
		0 \overset{\mathrm{CE}}{=} \int_E \partial_t \varphi \, \hat{\rho}_\varepsilon + \int_E \partial_y\varphi  \, \hat{\jmath}_\varepsilon \xrightarrow{\varepsilon\to 0} \int_E \partial_t \varphi \, \hat{\rho}_0^\ell + \int_E \partial_y\varphi \,\hat{\jmath}_0^\ell.
	\end{equation*}
	Since~$\hat{\rho}_0^\ell$ is concentrated on~$[0,T]\times\{\pm1/2\}$, the limiting flux is piecewise constant with jumps only at~$\{\pm1/2\}$, and due to the fact that~$\hat{\jmath}_0^\ell$ is in~$L^2(E)$, this limiting flux must vanish outside of~$(-1/2,+1/2)$. Therefore, the continuity equation~$0=\partial_t\hat{\rho}_0^\ell + \partial_y\hat{\jmath}_0^\ell$ in the distributional sense implies that the flux is given by
	\begin{equation*}
	\hat{\jmath}_0^\ell(t,y)= -\partial_t\hat{z}_0^\ell(t)\mathbf{1}_{(-1/2,+1/2)}(y).
	\end{equation*}
	\smallskip
	
	\ref{item:GF_NGF:proof:rho-j-conv:2}: By definition, the flux~$\hat{\jmath}_\varepsilon$ is given by
	\begin{equation*}
		\hat{\jmath}_\varepsilon = -e^{B}\partial_y\hat{v}_\varepsilon^\ell = -\partial_y\hat{u}_\varepsilon^\ell + \hat{b}_0\mathbf{1}_{E_0}\hat{u}_\varepsilon^\ell.
	\end{equation*}
	By the estimate~\eqref{eq:GF_NGF:limsup:u_eps-bounded-in-H1-Omega} from Lemma~\ref{lemma:GF_NGF:unif-energy-est}, the function~$\hat{u}_\varepsilon^\ell$ converges to some function~$\hat{u}_0^\ell$ weakly in~$L^2(0,T;L^2(\Omega))$, and as shown above, the flux~$\hat{\jmath}_\varepsilon$ is bounded in~$L^2(E)$. Hence for any test function~$\varphi\in C^\infty_c((0,T)\times\Omega)$,
	\begin{equation*}
		0=\lim_{\varepsilon\to 0}\int_{(0,T)\times \Omega}\left[\varphi\hat{\jmath}_\varepsilon -\hat{u}_\varepsilon^\ell\partial_y\varphi -\hat{b}_0\mathbf{1}_{E_0}\hat{u}_\varepsilon^\ell\right]=
		\int_{(0,T)\times \Omega}\left[\varphi\hat{\jmath}_0^\ell -\hat{u}_0^\ell\partial_y\varphi -\hat{b}_0\mathbf{1}_{E_0}\hat{u}_0^\ell\right].
	\end{equation*} 
	Therefore~$\hat{\jmath}_0^\ell = -\partial_y\hat{u}_0^\ell + \hat{b}_0\mathbf{1}_{E_0}\hat{u}_0^\ell$ weakly in~$L^2((0,T)\times\Omega)$. Since we also found above that~$\hat{\jmath}_0^\ell(t,y) = -\partial_t \hat{z}_0^\ell(t)\mathbf{1}_{(-1/2,+1/2)}(y)$, this means that in~$E_0$, the function~$\hat{u}_0^\ell$ is the weak solution to the ODE
	\begin{equation*}
		-\partial_y\hat{u}_0^\ell(t,y) + \hat{b}_0(t,y)\hat{u}_0^\ell(t,y) = -\partial_t\hat{z}_0^\ell(t).
	\end{equation*}
	We are left with verifying the boundary conditions. We will prove that for any test function~$\psi\in C_c^\infty(0,T)$,  
	\begin{align*}
		0&=\int_0^T\left(\hat{z}_0^\ell(t)-\hat{u}_0^\ell(t,-1/2)\right)\psi (t)\,dt,\\
		0&= \int_0^T\hat{u}_0^\ell(t,+1/2)\psi(t)\,dt.
	\end{align*}
	For~$\delta>0$, let~$U_\delta$ be a small neighborhood around~$(-1/2)$ of length~$2\delta$. Since~$\hat{u}_\varepsilon^\ell$ is uniformly bounded in~$L^2(0,T;H^1(\Omega))$ by~\eqref{eq:GF_NGF:limsup:u_eps-bounded-in-H1-Omega} of Lemma~\ref{lemma:GF_NGF:unif-energy-est}, there is a $C(t)$ such that
	\begin{equation*}
		\hat{u}_\varepsilon^\ell(t,y)\leq \hat{u}_\varepsilon^\ell(t,-1/2) + C(t)\,|y+(1/2)|^{1/2}.
	\end{equation*}
	With that, we can estimate
	\begin{multline*}
		\int_0^T\psi(t)\hat{\rho}_\varepsilon(t,U_\delta)\,dt = \int_0^T\int_{U_\delta} \psi(t)\hat{u}_\varepsilon^\ell(t,y)\hat{g}_\varepsilon^\ell(y)\,dydt\\
		\leq \int_0^T\psi(t)\hat{u}_\varepsilon^\ell(t,-1/2)\,dt \int_{U_\delta}\hat{g}_\varepsilon^\ell(y)\,dy + C \|\psi\|_{L^\infty} \int_{U_\delta} |y+(1/2)|^{1/2}\hat{g}_\varepsilon^\ell(y)\,dy\\
		\leq \int_0^T\psi(t)\hat{u}_\varepsilon^\ell(t,-1/2)\,dt \int_{U_\delta}\hat{g}_\varepsilon^\ell(y)\,dy + C \|\psi\|_{L^\infty} \delta^{1/2} \int_{U_\delta}\hat{g}_\varepsilon^\ell(y)\,dy.
	\end{multline*}
	For each~$\delta>0$,~$\int_{U_\delta}\hat{g}_\varepsilon^\ell(y)dy$ converges to one as~$\varepsilon\to 0$, and
	\begin{equation*}
		\lim_{\varepsilon\to 0}\int_0^T\psi(t)\hat{\rho}_\varepsilon(t,U_\delta)\,dt = \int_0^T\psi(t)\hat{z}_0^\ell(t)\,dt.
	\end{equation*}
	Therefore,
	\begin{equation*}
		\int_0^T\psi(t)\hat{z}_0^\ell(t)\,dt \leq \liminf_{\varepsilon\to 0}\int_0^T\psi(t)\hat{u}_\varepsilon^\ell(t,-1/2)\,dt + C'\delta^{1/2}.
	\end{equation*}
	Noting that~$\delta>0$ is arbitrary and repeating the argument for the reversed inequality, we find that
	\begin{equation*}
		\int_0^T\psi(t)\hat{z}_0^\ell(t)\,dt = \lim_{\varepsilon\to 0}\int_0^T\psi(t)\hat{u}_\varepsilon^\ell(t,-1/2)\,dt,
	\end{equation*}
	and the first boundary conditions follows since~$\hat{u}_\varepsilon^\ell$ converges in~$L^2(\Omega)$ and~$\hat{u}_0^\ell(t,\cdot)$ is continuous. The argument for the second boundary condition is similar, using that~$\hat{g}_\varepsilon^\ell(+1/2)\to\infty$ as~$\varepsilon\to 0$.
	\smallskip
	
	\ref{item:GF_NGF:proof:rho-j-conv:3}: With $B(t,y) = \int_{-1/2}^yb(t,z)\,dz$, the solution in~$E_0$ satisfying the boundary condition~$\hat{u}_0^\ell(t,-1/2)=\hat{z}_0^\ell(t)$ is given by
	\begin{equation*}
		\hat{u}_0^\ell(t,y) = e^{B(t,y)}\left[\hat{z}_0^\ell(t) + \partial_t\hat{z}_0^\ell(t) \int_{-1/2}^ye^{-B(t,z)}\,dz\right].
	\end{equation*}
	A calculation yields that
	\begin{equation*}
		\int_{-1/2}^{+1/2}e^{-B(t,z)}\,dz = -\frac{\hat{z}_0(t)}{\partial_t\hat{z}_0(t)} > 0.
	\end{equation*}
	The boundary condition~$\hat{u}_0^\ell(t,+1/2)=0$ therefore enforces
	\begin{equation}\label{GF_NGF:eq:upper-bound:proof-lemma-limiting-density}
		\partial_t \log\hat{z}_0^\ell(t) = \partial_t\log\hat{z}_0(t).
	\end{equation}
	The convergence assumption on the initial condition~$\hat{\rho}_\varepsilon(0,dy)$ implies~$\hat{z}_0^\ell(0)=\hat{z}_0(0)$. Hence by \eqref{GF_NGF:eq:upper-bound:proof-lemma-limiting-density}, we obtain~$\hat{z}_0^\ell = \hat{z}_0$. Now the fact that~$\hat{u}_0^\ell$ equals~$\hat{u}_0$ on~$E_0$ follows from an explicit calculation. Alternatively, we note that~$\hat{z}_0^\ell = \hat{z}_0$ and~$\partial_t\hat{z}_0^\ell=\partial_t\hat{z}_0$ implies that~$\hat{u}_0$ is a solution to the ODE, and the result follows from uniqueness of solutions.
\end{proof}
\section{Appendix---useful lemmas}
\begin{lemma}[Laplace's method]\label{lemma:watson}
	Let $f : [a,b] \to \mathbb{R}$ be twice differentiable. Suppose that for some $x_i \in (a,b)$, we have $f(x_i) = \inf_{[a,b]} f$. Then
	\begin{equation*}
	\int_a^b e^{-nf(x)}dx = \left[ 1 + o(1)\right] \sqrt{\frac{2\pi}{n f''(x_i)}} e^{-n f(x_i)}, \quad n \to \infty.
	\end{equation*}
	If $x_i = a$ or $x_i=b$, then
	\begin{equation*}
	\int_a^b e^{-nf(x)}dx = \left[ 1 + o(1)\right] \frac{1}{2} \sqrt{\frac{2\pi}{n f''(x_i)}} e^{-n f(x_i)}, \quad n \to \infty.
	\end{equation*}
\end{lemma}
\begin{lemma}[Dual of convex functions]\label{GF_NGF:lemma:appendix:dual-of-convex-functions}
	For $X=[0,T]\times \mathbb{R}^d$ and $f,g: X \to \mathbb{R}$ measurable with $g > 0$, any nonnegative Borel measure $\mu$ satisfies
	\begin{equation*}
	\int_X \frac{1}{2} \frac{|f(x)|^2}{g(x)} \, \dd \mu(x) = \sup_{\substack{b \in C_c^\infty(X)}} \int_X \left[\left(-\frac{b(x)^2}{2}\right)g(x) + b(x)f(x)\right]\,\dd\mu(x),
	\end{equation*}
	with the integral diverging when the supremum is infinity.
\end{lemma}
A proof is given for instance in~\cite[Lemma~3.4]{ArnrichMielkePeletierSavareVeneroni2012}. The representation in there can be further simplified by setting $a=-b^2/2$.
\begin{lemma}[Variational Problem]\label{lemma:GF_NGF:variational-problem}
	Define the function~$S$ by~\eqref{GF_NGF:eq:S-fct}. Then:
	\begin{enumerate}[label=(\roman*)]
		\item \label{item:GF_NGF:lemma:var-probl:i} We have
		\begin{equation}\label{eq:GF_NGF:lemma:var-probl}
		S(j,z) = \frac{1}{4} \inf_u \int_{-1/2}^{+1/2} \frac{1}{u(y)}\big|j + \partial_y u(y)\big|^2\,dy,
		\end{equation}
		where the infimum is taken over smooth functions~$u:[-1/2,+1/2]\to[0,\infty)$ satisfying the boundary conditions $u(-1/2)=z$ and~$u(+1/2)=0$.
		\item \label{item:GF_NGF:lemma:var-probl:ii} The optimizer in~\eqref{eq:GF_NGF:lemma:var-probl} is the polynomial
		\begin{equation}\label{eq:GF_NGF:lemma-var-probl:optimal-u}
			u(y) = -(j-z)(y-y_0)\left(y-\frac{1}{2}\right), \quad y_0 = \frac{1}{2}\frac{z+j}{z-j}.
		\end{equation}
	\end{enumerate}
\end{lemma}
The proof is best carried out by exploiting the fact that the energy is conserved, since the function in the variational problem does not depend explicitly on~$y$.
\chapter{Discussion and Future Questions}
\label{chapter:discussion}
Here we summarize the main results presented in each chapter and discuss interesting future questions related to them.
\subsubsection{Chapter~\ref{chapter:LDP-for-switching-processes}: Large Deviations of Switching Processes.}
\paragraph{Summary.} We consider a class of switching processes~$(X^\varepsilon,I^\varepsilon)$ in a periodic setting and prove pathwise large deviation principles of their spatial components~$X^\varepsilon$ in the limit~$\varepsilon\to 0$. The switching processes are motivated by stochastic models describing the spatial position of molecular motors walking on filaments within a cell, where the parameter~$\varepsilon>0$ corresponds to the ratio of microscopic to macroscopic scales. Our results embed existing results about molecular motors in a large-deviation context. The proofs of large deviation principle for the various models of molecular motors are examples of a general strategy outlined by Theorems~\ref{thm:results:LDP_switching_MP} and~\ref{thm:results:action_integral_representation}. In particular, the large-deviation proofs are independent of the specific choices involved in the models. Our method of proof exploits the connection of large deviations to Hamilton-Jacobi equations~\cite{FengKurtz2006}. Based on this connection, we find a strategy of proof consisting of two steps: first, identifying a multivalued limit operator, and second, solving a principal-eigenvalue problem.
\smallskip

In the models, 
the periodic setting reflects the periodic stucture of the filaments. As a consequence of this periodicity, the motor cannot advance without coupling to the chemically active environment, but switching mechanism between different configurations can generate motion. We derive an exact formula for the motor's large-scale velocity,~$v=\partial_p\mathcal{H}(0)$. This formula, based on the principal eigenvalue~$\mathcal{H}(p)$ of a cell problem, coincides with the findings of Perthame and Souganidis~\cite{PerthameSouganidis09a}.
We work with variational representations of the principal eigenvalues to derive from the large deviation principles the following fact: a non-zero velocity~$v$ can only be achieved if detailed balance is broken (Theorem~\ref{thm:results:detailed_balance_limit_I}).
\paragraph{Discussion and future questions.} Our more concrete conclusions based on working with the Hamiltonians~$\mathcal{H}(p)$ are limited to detailed balance. It would be interesting to investigate the Hamiltonians for systems not satisfying detailed balance. 
Hastings, Kinderlehrer and Mcleod for instance showed that transport occurs if potentials and rates collaborate in a suitable way~\cite[Theorem~2.1]{HastingsKinderlehrerMcleod08}. These conditions should consequently imply a non-trivial velocity~$v$. Another interesting question is the behaviour of the motor under load- or external forces. We showed under detailed balance that with a constant external force~$F$, the Hamiltonians are symmetric around~$-F$, which means that a positive (negative) force leads to a positive (negative) velocity. In general, does the velocity depend monotonically on external forces? Is there a stalling force~$F$ with which the motor's velocity vanishes? We could not find suitable symmetries of~$\mathcal{H}(p)$ for answering these questions.
\smallskip

Another open question is related to the coupled Fokker-Planck equations of the molecular-motor models. 
Chipot, Kinderlehrer and Kowalczyk considered a variational formulation for molecular motors~\cite{ChipotKinderlehrerKowalczyk2003}, similar in spirit to the JKO-scheme of Jordan, Kinderlehrer and Otto for the diffusion equation~\cite{JordanKinderlehrerOtto1998}.
It would be interesting to know whether we can derive such variational formulations from large deviations of empirical densities, in the same manner as 
Adams, Dirr, Peletier and Zimmer derived Wasserstein gradient flows~\cite{AdamsDirrPeletierZimmer2011}.
We do not expect a gradient-flow structure for molecular motors, since molecular motors are modelled by irreversible processes, and reversible processes lead to gradient flows as shown by Mielke, Peletier and Renger~\cite{MielkePeletierRenger2014}. However, once one knows how to derive a meaningful variational formulation in this example, one might be able to obtain variational formulations for similar irreversible processes as well.
\subsubsection{Chapter~\ref{chapter:LDP-of-empirical-measures}: Large Deviations of Empirical Measures.} 
\paragraph{Summary.}
The zig-zag process is an irreversible piecewise-deterministic Markov process designed to have a specific Gibbs-type stationary measure. We prove that its empirical measure satisfies a large deviation principle. Classical results in large deviation theory are not applicable due to the finite-speed and non-diffusive character of the zig-zag process. Therefore we derive suitable conditions based on the semigroup approach to large deviations~\cite{FengKurtz2006}. Our main contribution lies in proving that the Lyapunov functions in~\ref{item:thm_LDP_non_compact:Lyapunov} and the mixing property~\ref{item:thm_LDP_non_compact:mixing} suffice for a proof of large deviations in a non-compact state space.
\smallskip

We cannot characterize the rate functions by the Donsker-Varadhan formula for reversible diffusions, due to the inherent irreversibility of the zig-zag process. We derive an explicit formula of the rate function for the compact case. Based on this characterization, we conclude that the optimal rate of convergence is achieved by setting the refreshment rate~$\gamma$ (in Eq.~\eqref{eq:switching-intensity-condition-2}) to zero.
\paragraph{Discussion and future questions.}
Our conclusions about the zig-zag process are limited to one dimension. It would be very interesting to know whether our results also hold in higher dimensions.
While the Lyapunov functions that we found for the zig-zag process are suitable for arbitrary dimensions, we were not able to verify the mixing property~\ref{item:thm_LDP_non_compact:mixing}.
\smallskip

Another open question is whether we can also explicitly characterize the rate functions in higher dimensions---for the zig-zag process as well as for other PDMPs such as the bouncy particle sampler~\cite{BouchardCoteVollmerDoucet2017}. The idea of using large-deviation rate functions to compare the performance of MCMC algorithms was introduced in~\cite{plattner2011infinite, dupuis2012infinite}. Rey-Bellet and Spiliopoulos showed that adding irreversible drifts to a diffusion process increases the rate functions~\cite[Theorem~2.2]{rey2015irreversible} and decreases the asymptotic variance~\cite[Theorem~2.7]{rey2015irreversible}. Explicit characterizations of the rate functions would be useful to address similar performance questions for PDMPs. Further natural steps are to compare samplers based on drift-diffusion processes and PDMPs, and to investigate how the rate functions scale with the dimension. An answer to the latter question would give interesting insights into how the various algorithms deal with the curse of dimensionality.
%
%
\smallskip

Nicolás García Trillos and Daniel Sanz-Alonso recently demonstrated that samplers based on drift-diffusion processes converge faster to equilibrium when choosing a suitable non-Euclidean metric for the space of position variables~\cite[Theorem~4.1, Proposition~4.3]{TrillosSanz-Alonso2018}. The authors call these processes geometry-informed Langevin diffusions, and their conclusions are based on an analysis of the spectral gap. Their results raise the question of whether a similar effect can be observed from a large-deviation point of view and for geometry-informed PDMPs. For instance, in between jumps of the velocity variables, the zig-zag and bouncy-particle samplers move in straight lines. It would be interesting to explore whether these samplers can benefit from modifying the piecewise-deterministic dynamics to follow geodesics with respect to a non-Euclidean metric.
\subsubsection{Chapter~\ref{chapter:LDP-in-slow-fast-systems}: Large Deviations in Stochastic Slow-Fast Systems.} 
\paragraph{Summary.} We consider two-component stochastic processes whose individual components run at different time scales. Our main results are a proof of large deviation principles in the limit of an infinite time-scale separation and an interpretation of the Lagrangian rate functions we obtain. The analytical challenge in the proof (the comparison principle for an associated Hamilton-Jacobi equation) is solved in Chapter~\ref{chapter:CP-for-two-scale-H}. The results apply in particular to irreversible diffusions as fast processes.
Our main example are mean-field interacting particles coupled to fast diffusion processes, for which we deduce an averaging principle from the large deviation principle. A key ingredient for this argument is a suitable formula for the Lagrangians. 
\paragraph{Discussion and future questions.}
We assumed the fast variables to live in a compact space to focus only on the effects coming from the scale separation. It would be worthwhile to extend the analysis to the non-compact setting in order to cover for instance a fast Ornstein–Uhlenbeck process. Another interesting question we left unanswered is whether one can treat \emph{degenerate} diffusions as in~\cite{BudhirajaDupuisGanguly2018} with our methods---we always worked under uniform ellipticity assumptions. In all these examples, the key problem one has to solve is the comparison principle. We further comment on that in the discussion below.
\subsubsection{Chapter~\ref{chapter:CP-for-two-scale-H}: Comparison Principle for Two-Scale Hamiltonians.} 
\paragraph{Summary.} We prove existence and uniqueness of solutions of a Hamilton-Jacobi equation, where the Hamiltonian is given by an optimization over control variables. The Hamiltonians 
appearing in 
large-deviation problems for slow-fast systems (Chapter~\ref{chapter:LDP-in-slow-fast-systems}) are of this type. We propose a bootstrap procedure to solve the comparison principle, for which we have to assume sufficient regularity of the cost functions. The method applies to non-coercive Hamiltonians arising in mean-field models. Furthermore, it addresses a problem pointed out in~\cite{BudhirajaDupuisGanguly2018}, which is that classical comparison results are not readily applicable due to the poor regularity properties of this type of two-scale Hamiltonians.
\paragraph{Discussion and future questions.}
There are various examples that we cannot treat with our method, but which are important to address. Let us mention two examples. First, if the internal Hamiltonians correspond to degenerate diffusions---we use uniform ellipticity in the proof of Proposition~\ref{prop:verify-ex:Lambda_quadratic}. Under Lipschitz conditions on the diffusion coeffcients, the comparison principle for degenerate diffusions is proven in~\cite[Lemma~9.25]{FengKurtz2006} by means of an auxiliary variable~$\lambda$. It would be interesting to investigate whether one can combine this method of proof to include the case of degenerate diffusions in two-scale Hamiltonians.
\smallskip

Second, we considered Hamiltonians arising from a scale separation in a weakly-coupled regime. That is reflected in the fact that the  cost functionals do not depend on~$p$. But there are interesting problems leading to such as setting.
For instance, in the molecular-motor models, we only discussed potentials and rates $\{\psi_i, r_{ij}\}$ depending on the up-scaled variables, in the sense that~$\psi_i=\psi_i(x/\varepsilon)$ and~$r_{ij}=r_{ij}(x/\varepsilon)$. That assumption leads to a simplification, since then the well-posedness of the principal-eigenvalue problem is sufficient for proving the comparison principle---this is basically the content of~\ref{MM:item:T1} and~\ref{MM:item:T2} of Theorem~\ref{thm:results:LDP_switching_MP}, which state that finding an eigenvalue~$\mathcal{H}(p)$ and an eigenfunction~$\varphi_p$ are sufficient. When we consider instead potentials~$\psi_i=\psi_i(x,x/\varepsilon)$ and rates~$r_{ij}=r_{ij}(x,x/\varepsilon)$, the eigenvalue Hamiltonians are---similar to~\eqref{eq:results:LDP_MM:DV_var_rep_H(p)} of Section~\ref{subsection:detailed_balance}---of the form
\begin{equation*}
\mathcal{H}(x,p)=\sup_{\mu\in\mathcal{P}(E')}\left[\int_{E'}V_{x,p}(z)\,\dd\mu(z) - \mathcal{I}_{x,p}(\mu)\right].
\end{equation*} 
In there, the maps~$V_{x,p}$ and~$\mathcal{I}_{x,p}$ are obtained from~\eqref{MM:eq:function-V-in-Hamiltonian} and~\eqref{eq:results:LDP_MM:DV_functional} by replacing the potentials and rates. It is unknown whether the comparison principle is satisfied for the Hamilton-Jacobi equation with these Hamiltonians. The Hamiltonians we obtain in slow-fast systems in Chapter~\ref{chapter:LDP-in-slow-fast-systems} are simpler in the sense that the Donsker-Varadhan functionals~$\mathcal{I}_{x,p}(\mu)$ in there are indepedent of the momentum variable~$p$. It would be interesting to explore whether the method developed in Chapter~\ref{chapter:CP-for-two-scale-H} can be extended to include this type of Hamiltonians.
\subsubsection{Chapter~\ref{chapter:GF-to-NGF}: Gradient Flow to Non-Gradient-Flow.} 
\paragraph{Summary.}
We study a family of Fokker-Planck equations corresponding to a particle diffusing in an asymmetric double-well potential. The associated gradient-flow structures do not converge in a certain limit due to the relative entropies diverging, which originates from the asymmetry of the potential. We propose to work instead with a different variational formulation based on functionals that include fluxes, and show~$\Gamma$-convergence of these functionals.
Our motivation is taken from the fact that reversible processes give us gradient flows via large deviation theory. Therefore the convergence of gradient-flow structures appears in many contexts, and it is natural to ask which convergence concepts are suitable for treating cases in which the underlying processes become irreversible.
\paragraph{Discussion and future questions.} 
It would be exciting to investigate other cases of gradient-flow structures that are not converging due to the relative entropies diverging. On the level of the functionals, one may regard the inclusion of fluxes as "absorbing" or "including" the relative entropies into the dissipation functional. While we use a special coordinate transformation that is akin to the problem we study, a natural question is whether the techniques we employed in our case also apply to~$\Gamma$-convergence problems for other density-flux functionals.
\qed
\backmatter
\cleardoublepage
\chapter{Summary}
\begin{center}
\Large
\myTitle
\end{center}
\thispagestyle{empty}
In this thesis we study path distributions of stochastic processes by means of large deviation theory.
We focus on processes that are typically time-irreversible, which means that inverting time leads to a different path distribution. Our main motivation comes from the fact that while reversible processes lead via large deviation theory to gradient flows, it is an open question of which variational structures can be obtained from irreversible processes. In this thesis we make a first step to answering this question by deriving large deviation principles for irreversible processes.
\smallskip
 
The stochastic processes we consider
depend on a parameter
characterzing the concrete process at hand; a length-scale separation, the number of particles in a system, the time variable itself, or a parameter modeling separation of time scales. If the parameter is sent off to infinity, the stoachstic process becomes deterministic. That means in the limit, realizations of the processes are with probability one equal to a particular limiting trajectory. We use large deviation theory to show that 
the probability of obtaining an atypical realization of the process vanishes exponentially fast, with a rate depending on the atypical trajectory. Our aim is to express this rate as an integral over time involving a so-called Lagrangian, which provides one way of determining the limiting typical behaviour of the stochastic process. In the reversible case, the connection to gradient flows case is derived using symmetries of the Lagrangians.
\smallskip

The first two chapters introduce the basic concept of large deviation theory applied in the context of stochastic processes. In particular, we illustrate the Feng-Kurtz method of how to rigorously derive Lagrangians of a sequence of stochastic processes starting from the infinitesmial generators of the processes.
\smallskip

In Chapter~3, we contribute to the analysis of stochastic models of molecular motors, which are proteins transporting cargo in living cells. The stochastic processes model the position of a molecular motor walking on a filament. We use the Feng-Kurtz method to prove large deviations principles of the position variable in the large-scale limit. Our results provide one way of analysing the macroscopic behaviour of the molecular motor starting from the microscopic dynamics. For instance, the influence under external forces, or the fact that transport can only occur if time-reversibility is broken.
\smallskip

Chapter~4 is dedicated to the analysis of Markov chain Monte Carlo (MCMC) methods based on piecewise-deterministic Markov processes (PDMP). The idea behind MCMC is to approximate a probability distribution by the occupation time measure of a stochastic process. Under ergodicty assumptions, the occupation time measure converge to the stationary measure of the process, which is designed to be equal to the desired probability distribution. Our results offer a framework for proving that the associated occupation time measures of PDMP's satisfy in fact a large deviation principle. Classical theorems do not apply due to the singular nature of PDMPs, but we show how the nonlinear semigroup approach provides one way of overcoming this difficulty.
We learn from our results that for the zig-zag process, maximal irreversibility corresponds to the optimal rate of convergence to stationarity.
\smallskip

In Chapters~5 and~6, we consider stochastic slow-fast systems. In particular, we are interested in mean-field interacting particles where the interaction rates are fluctuating on a much faster time scale than the particle's evolution. Intuitively, one expects the particle system to evolve under averaged interaction rates, which is refered to the averaging principle. We first prove large deviation principles of the particle densities and fluxes in the simultaneous limit of infinitely many particles and time-scale separation tending to infinity. Then we show that the averaging principle holds as a consequence of the large deviation principle. The techniques are based on Hamilton-Jacobi theory. Chapter~5 contains the large-deviation analysis, while we solve in Chapter~6 more general Hamilton-Jacobi-Bellman equations arising in this context.
\smallskip

Finally, we analyse in Chapter~7 partial differential equations arising in models of chemical reactions. The equations contain parameters modeling the activation energy of certain chemical reactions and the time-scale of reaction events. A crucial role in our analysis is the variational formulation of such PDEs by means of density-flux functionals. In the limit of large-activation energy, we prove~$\Gamma$-convergence of these functionals. On the level of underlying stochastic processes, this convergence result corresponds to passing from reversible to irreversible processes. With this problem we address the question of which variational formulations beyond gradient flows are suitable for studying such limits from reversible to irreversible.
\smallskip
\cleardoublepage
\cleardoublepage
\chapter{Acknowledgments}
First I want to thank you, Mark. This thesis would not have been possible without your support and guidance. I am grateful that you gave me the opportunity to become a mathematician; that you patiently taught me how to approach mathematical questions, how to write papers, how to present science; also that erasers are personal belongings rather than public goods.
I enjoyed in particular working with you in the office; I always learned a lemma, a theorem, and your humor made the discussions enjoyable. I missed that during times of corona.
I also appreciate that you nudged me to organize CASA Days and the Wednesday Morning Sessions---this gentle kick making me become active is exactly what I needed. 
\smallskip

I thank the committee members for accepting the invitation, and for their many comments, remarks and suggestions after reading the first manuscript, which motivated me to read up on many other works. Thank you for your interest in my thesis and in our work.
\smallskip

Frank, Francesca and Federico, I greatly enjoyed the frequent and inspiring meetings we had in Delft and Eindhoven. I learned to appreciate stochastic processes, martingales and Brownian motion during our discussions, and liked the atmosphere that you created in our meetings. 
\smallskip

I am also grateful to my collaborators. Richard, you always made me feel welcome in Delft by chatting about life and mathematics. Thank you for offering me numerous cups of coffee; I liked our discussions about comparison principles,~$f=h$, and many more topics. Joris and Pierre, thank you for introducing me to the world of MCMCs. I learned a lot from your style of writing and from our frequent attempts to show that the zig-zag process is a nice process. Mario, thank you for our discussions in Eindhoven and Bonn.
\smallskip

Jin, thank you for answering all my questions via many emails and for our discussions at the Leiden conference; they were crucial for making progress.
\smallskip

The four years at CASA have been a lot of fun. 
Diane, Enna and Jolijn, you were a great support during my stay at CASA. I want to thank my office mates I had over the years; Saeed, Upanshu, Koondi, Jasper, Anastasiia and Alberto. Thank you Anastasiia and Jasper 
for the fun times during the ODE course and discussions on the whiteboard, and for saving my plants (before and during corona). Arthur (Zigge-zagge); thanks for all the amusing and entertaining talks. Anastasiia, Harshit, Xingang, it was fun to get our chip-predictions right (Brownian motion). Jim, I really enjoyed being part of your measure-theory course, discussing exams and homeworks with you, and going for pizza to celebrate the end of a course; thank you for this time, I learned a lot from you about measure theory. Georg, thank you for discussing compact operators and principal eigenvalues with me; it was also great preparing the workshop on quantum computing with you. Thanks Oliver for frequently dropping by at our office to chat about math. Thank you Jan-Cees for the interesting projects during the ODE course. Alberto, Carlo, Oxana, thank you for your suggestions during the Wednesday Morning Sessions we had so far. Finally, I also want to thank all others who make CASA a welcoming place. 
\smallskip

I am also glad to be part of the random people in Delft; Andrea, Bart, Federico, Francesca, Mario, Martina, Richard, Rik, Sebastiano, Simone, I always enjoyed being at TU Delft and chatting with you over coffee!
\smallskip

Special thanks goes to the Nederlandse Spoorwegen and Julia's. Thanks to your trains and coffees during these four years, I could live in (and in between) Eindhoven, Delft and Leiden.
\smallskip

Ich m{\"o}chte besonders meiner Familie danken, die mich in den vergangenen Jahren begleitet hat und zu meiner Verteidigung anreist: Mama, Papa, Joana, Nielsson, Oma, Max und Olha.
Mama, danke dass du dich so sehr daf{\"u}r eingesetzt hast, mich auf ein Gymnasium zu bringen; ohne all deine Bem{\"u}hungen und Unterst{\"u}tzung h{\"a}tte ich niemals anfangen k{\"o}nnen diese Arbeit zu schreiben.
Papa, danke dass du mir kurz nach meiner Geburt den kleinen Fermat vorgelesen hast, in der Hoffnung ich h{\"a}tte einen Geistesblitz. Zwar blieb dieser bis heute aus, aber der Funke ist {\"u}bergesprungen weil du meine Freude an der Physik und Mathematik immer unterst{\"u}tzt und befeuert hast.
\smallskip

Mia, thank you so much for supporting me during the whole time of my PhD. You were there for me. You discussed math with me, we prepared exams and homeworks for measure theory together, you encouraged me during the process of writing the thesis when I needed you; thank you for all of that! I am happy about all the memories we share from the last years, and when looking ahead with you. Thank you for coming into my life, it is wonderful with you!
\chapter{Curriculum Vitae}
\pagestyle{empty}
Mikola Christoph Schlottke was born on 29-06-1991 in Erlangen, Germany. After finishing high school in 2010 at the Friedrich-Alexander-Gymnasium in Neustadt an der Aisch, he did a Voluntary Year of Social Service at the Bavarian Red Cross as a paramedic. He then started his studies of Physics at the University of Potsdam in Germany, which he completed in 2014 with distinction and a thesis on the tunnel effect under the supervision of prof.dr. Markus Klein. He continued to study Theoretical Physics at the University of Amsterdam, where he graduated in 2016 with a thesis on the Atiyah-Singer-Index Theorem under the supervision of dr. H.B. Posthuma.
\smallskip

In October 2016, he started a PhD project at the Eindhoven University of Technology under the supervision of prof.dr. Mark A. Peletier. The results obtained during this project are presented in this dissertation. The PhD project was part of the TOP-1 project \emph{Large deviations and gradient flows: beyond equilibrium}, which included regular meetings with prof.dr. F.H.J. Redig, dr. Francesca Collet and dr. Federico Sau, and was funded by the NWO grant 613.001.552.
\addcontentsline{toc}{chapter}{Bibliography}
\bibliography{bib_Mikola,ref_Mikola_zigzag,KraaijBib}
\bibliographystyle{alpha}
\end{document}